\documentclass[12pt]{amsbook}

\usepackage{amsmath, amsthm, amssymb}

\usepackage[cal=boondox]{mathalfa}
\usepackage{bbm}

\usepackage[toc,page]{appendix} 

\makeindex

\usepackage{mathrsfs}
\usepackage{calc}
\usepackage{graphicx}
\usepackage[all,cmtip]{xy}
\usepackage{amscd}
\usepackage{lmodern}
\usepackage{amsfonts}
\usepackage{picinpar}
\usepackage{tikz}
\usepackage{tikz-cd}
\usepackage{mathrsfs}
\usetikzlibrary{arrows}
\usepackage{enumerate}
\usepackage{stmaryrd}
\usepackage{array}
\usepackage{extarrows}
\usepackage{hyperref} 
\usepackage{enumerate}
\usepackage{verbatim}
\usepackage{esint}
\usepackage[T1]{fontenc}
\usepackage{xcolor}
\usepackage{enumitem}
\usepackage{hyperref}
\usepackage{mathtools}

\makeatletter
\def\paragraph{\@startsection{paragraph}{4}%
  \z@\z@{-\fontdimen2\font}%
  {\normalfont\bfseries}}
\makeatother

\newfont{\bigbf}{cmbx10 scaled\magstep1}
\normalbaselines
\newtheorem{thm}{Theorem}[chapter]
\newtheorem{assumption}{Assumption}[chapter]
\newtheorem{obs}[thm]{Observation}

\newtheorem{lem}[thm]{Lemma}
\newtheorem{prop}[thm]{Proposition}

\newtheorem{cor}[thm]{Corollary}
\theoremstyle{definition}

\newtheorem{defi}[thm]{Definition}
\newtheorem{definition}[thm]{Definition}

\newtheorem{example}[thm]{Example}
\newtheorem{question}[thm]{Question}

  \newtheorem{introthm}{Theorem}

\theoremstyle{remark}

\newtheorem{rem}[thm]{Remark}
\newtheorem{introrem}[introthm]{Remark}

\numberwithin{section}{chapter}
\numberwithin{subsection}{section}

\numberwithin{figure}{chapter}

\DeclareMathOperator{\C}{\mathbb{C}}
\DeclareMathOperator{\Sh}{Sh}

\DeclareMathOperator{\PSH}{PSH}

\DeclareMathOperator{\GL}{GL}
\DeclareMathOperator{\SL}{SL}
\DeclareMathOperator{\SO}{SO}
\DeclareMathOperator{\Sp}{Sp}

\DeclareMathOperator{\SU}{SU}

\DeclareMathOperator{\PGL}{PGL}
\DeclareMathOperator{\PSL}{PSL}

\DeclareMathOperator{\Th}{Th}
\DeclareMathOperator{\coh}{H}
\DeclareMathOperator{\Opp}{Opp}
\DeclareMathOperator{\vol}{vol}
\DeclareMathOperator{\Ad}{Ad}

\DeclareMathOperator{\Isom}{Isom}
\DeclareMathOperator{\Stab}{Stab}
\DeclareMathOperator{\hull}{hull}
\DeclareMathOperator{\pos}{pos}
\DeclareMathOperator{\supp}{supp}

\def\half{\frac{1}{2}}

\DeclareMathOperator{\acts}{\curvearrowright}

\DeclareMathOperator{\R}{\mathbb{R}}

\DeclareMathOperator{\diag}{diag}

\newcommand\bac[1]{{\color{blue} #1} }

\newcommand\intcur[1]{\llbracket #1 \rrbracket}

\date{\today}

\def\<{\langle}
\def\>{\rangle}

\def\tits{\partial_{Tits}}

\def\eps{\epsilon}
\def\tangle{\angle_{Tits}}
\def\taumod{\tau_{mod}}

\def\simod{\sigma_{mod}}

\def\Fmod{F_{mod}}
\def\amod{a_{mod}}
\def\dF{d_{\mathcal{F}}}
\def\dyl{>0}

\def\Uupt{U(P,t)}

\def\refthmmain{\ref{thm:main}}

\def\bB{{\mathbb{B}}}
\def\bC{{\mathbb{C}}}

\def\bR{{\mathbb{R}}}

\def\P{{\mathbb{P}}}

\def\bH{{\mathbb{H}}}

\def\bZ{{\mathbb{Z}}}

\def\bP{{\mathbb{P}}}
\def\bN{{\mathbb{N}}}

\def\cB{{\mathscr{B}}}
\def\cC{{\mathscr{C}}}
\def\cD{{\mathscr{D}}}
\def\cE{{\mathscr{E}}}
\def\cF{{\mathscr{F}}}
\def\cH{{\mathscr{H}}}

\def\cM{{\mathscr{M}}}

\def\cO{{\mathscr{O}}}
\def\cP{{\mathscr{P}}}
\def\cQ{{\mathscr{Q}}}
\def\cR{{\mathscr{R}}}
\def\cS{{\mathscr{S}}}
\def\cT{{\mathscr{T}}}

\def\cW{{\mathscr{W}}}
\def\cX{{\mathscr{X}}}
\def\cY{{\mathscr{Y}}}

\def\sK{{\mathscr{K}}}

\def\fC{{\mathfrak{C}}}
\def\fm{{\mathfrak{m}}}
\def\fM{{\mathfrak{M}}}
\def\fd{{\mathfrak{d}}}
\def\fS{{\mathfrak{S}}}
\def\fu{{\mathfrak{u}}}
\def\fuwp{{\mathfrak{u}_w^+}}
\def\fuwm{{\mathfrak{u}_w^-}}
\def\Psiwp{{\Psi_w^+}}
\def\Psiwm{{\Psi_w^-}}
\def\Swx{{S_w(x^+)}}
\def\Swxd{{S_{w^\vee}(x^-)}}
\def\Oppxwd{{\operatorname{Opp}(x_{w^\vee})}}
\def\xpair{{(x^+,x^-)}}
\def\cah{{\mathcal{h}}}

\def\tga{{\tilde{\gamma}}}
\def\Cyli{{\mathcal{C}_I}}
\def\Cylii{{\mathcal{C}_{II}}}

\def\CC{{\mathbf{C}}}

\def\PP{{\mathbf{P}}}

\def\RR{{\mathbf{R}}}
\def\ZZ{{\mathbf{Z}}}

\def\ba{\mathbf{a}}

\def\fk{{\mathbf{k}}}
\def\bu{\mathbf{u}}
\def\bv{\mathbf{v}}
\def\bx{\mathbf{x}}

\def\F{\mathrm{F}}
\def\re{\mathrm{e}}

\def\Pc{{P^{\vee}}}

\def\Ga{\Gamma}
\def\ga{\gamma}
\def\La{\Lambda}
\def\la{\lambda}
\def\Om{\Omega}
\def\al{\alpha}

\def\geo{\partial_\infty}
\def\id{\mathrm id}

\def\ol{\overline}

\def\Span{\operatorname{Span}}

\def\inte{\operatorname{int}}

\def\diam{\operatorname{diam}}

\def\cG{{\mathscr{G}}}

\def\cL{{\mathscr{L}}}

\DeclarePairedDelimiter\ceil{\lceil}{\rceil}

\def\mini{\scriptsize}

\def\kap{}

\title[Equidistribution of currents]{Equidistribution of currents\\ under Anosov group actions}

\author[N.-B. Dang]{Nguyen-Bac Dang}
\address[N.-B. Dang]{Institut de Math\'ematiques d'Orsay\\
Universit\'e Paris-Saclay\\
F-91405 Orsay Cedex, France}
\email{nguyen-bac.dang@universite-paris-saclay.fr}
\author[M. Kapovich]{Michael Kapovich}
\address[M.~Kapovich]{Department of Mathematics\\1 Shields Ave.\\
  University of California, Davis\\
  CA 95616, USA}
\email{kapovich@math.ucdavis.edu}
\author[M. Lyubich]{Mikhail Lyubich} 
\address[M. Lyubich]{ Institute for Mathematical Sciences\\Stony Brook University\\
Stony Brook, NY 11794-3660, USA}
\email{mlyubich@math.stonybrook.edu}
\author[S. Zhao]{Shengyuan Zhao}
\address[S. Zhao]{Universit\'e Paul Sabatier\\
Institut de Math\'ematiques de Toulouse\\
118, route de Narbonne\\
F-31062 Toulouse Cedex 9, France}
\email{shengyuan.zhao@math.univ-toulouse.fr}

\begin{document}

\begin{abstract}
We study the behavior of currents on flag-manifolds under actions of Anosov subgroups $\Ga$ of complex semisimple Lie groups $G$. 
Given a current $T$ of bidimension $(k,k)$ on the full flag-manifold $\cF=G/B$, we average it under the $\Ga$-action via a construction analogous to the construction 
of Patterson--Sullivan measures. We show that under certain genericity conditions (in the case of currents of integration over subvarieties), 
the limiting current is a Gibbs current, i.e. is given by integration (with respect to a Gibbs measure) over the flag-limit set of $\Ga$ of  suitable multiples of currents of integration along Schubert varieties based at the limit points of $\Ga$. We prove that the same equidistribution result (but without any genericity assumptions) for currents defined by smooth forms on $\cF$. {We also prove a form of { Axiom A property} for the suspension flow of the $\Ga$-action on $\cF$.  }
\end{abstract}

\maketitle

\tableofcontents

\mainmatter

\chapter*{Introduction}

This book is partly motivated by exploring similarities between discrete subgroups of higher rank complex Lie groups $G$ and holomorphic dynamics in dimensions $\geqslant 2$.
While the analogy between rank $1$ group actions and rational mappings on the Riemann sphere initiated by Sullivan is very fruitful and has led to important results in both fields, we intend here to extend certain aspects of the dictionary to discrete subgroups of Lie groups of higher rank and holomorphic dynamics in higher dimension.

We first explain the holomorphic side of the dictionary. For $n\geqslant 2$, consider a holomorphic map $F \colon \bP^n(\C) \to \bP^n(\C)$,  defined by homogeneous polynomials of degree $d \geqslant 2$ with no common factors. One studies the action of (direct and inverse) iterations of $F$ on points, curves and algebraic subvarieties. 
For a given point $p$ outside an exceptional algebraic subset $E \subsetneq \bP^n(\C)$, Briend and Duval, \cite{briend_duval}, show that the preimages under $F^k$ of $p$ equidistribute to a unique $F$-invariant measure  $\mu$. This results  generalizes the equidistribution of points due to Brolin and Lyubich, \cite{brolin,lyubich_equidistribution}, for rational mappings on the Riemann sphere to higher dimensional varieties.
This  measure $\mu$ is an important invariant of $F$ and is the unique measure of maximal entropy. However,  measures are a special class of currents (of degree $0$) as they control the distribution of points under the (semi)group actions. In higher dimension, there are higher dimensional invariant currents which characterize the distribution of subvarieties under the action. 
  For self-maps of $\bP^2(\C)$, Favre and Jonsson, \cite{favre_jonsson_brolin}, show that the preimages under $F^k$ of an algebraic curve not passing through finitely many exceptional points equidistribute to a unique invariant (up to scaling) current $T$, called the {\em Green current}. 
  This current is directly connected to the measure of maximal entropy of the map $F$ as its self intersection $T\wedge T$ can be defined (via Bedford-Taylor's method \cite{bedford_taylor}) and is equal to $\mu$. Taking the support of the current, one obtains a family of $F$-invariant subsets stratified by the codimension, $J_0= \supp \mu \subset J_1 := \supp T$.
   The first set is called the {\em small Julia set} and  the second one is called the {\em big Julia set}.  This stratification is also known in higher dimensions and yields $d$ invariant subsets.  
While most of the work mentioned above study  a single holomorphic map,  the construction of natural currents and measures invariant by a large class of  subgroups of birational surface maps was carried by  Cantat-Dujardin \cite{cantat_dujardin}. Their approach uses naturally some arguments  from rank $1$-group dynamics, making the Sullivan dictionary more relevant. 

{On the side of group actions,} we will be studying equidistribution of complex subvarieties in flag manifolds $G/P$ under actions of 
{\em Anosov subgroups} $\Gamma<G$ of complex semisimple Lie groups $G$. While some recent work was done regarding equidistribution of points in the setting of Patterson--Sullivan measures (see e.g. \cite{Dey-Kapovich} for a  review of the current literature), the equidistribution for complex submanifolds of positive dimension remained untouched until now. However, even in the setting of equidistribution of {\em points}, we propose a different setup. Namely, instead of using invariant distance functions on the symmetric space of $G$ in the definition of the Patterson--Sullivan measures, we use {\em hyperbolic} length functions $\beta$ on the hyperbolic group $\Gamma$, i.e. length functions coming from geometric $\Ga$-actions on geodesic metric spaces (see Definition \ref{def:hyperbolic length function} and Lemma \ref{lem:regular} for details). 

Given a (hyperbolic) length function $\beta$, one defines the Poincar\'e series \index{Poincar\'e series $\cP(s)$}
$$
\cP(s)=\sum_{\gamma\in \Gamma}e^{-s\beta(\gamma)},
$$
its critical exponent $s_0$ and a family of probability measures on $\Gamma$, 
$$
 \mu(s) = \frac{1}{\cP(s)}\sum_{\gamma\in \Gamma}e^{-s\beta(\gamma)}\delta_{\gamma}, \quad s> s_0. 
$$

This leads to a measure $\mu$ on the Gromov-boundary $\geo\Gamma$ 
of $\Gamma$, by taking a suitable limit as $s\to s_0^+$. Unlike the  Patterson--Sullivan measure, 
the $\Gamma$-orbit of {\em every point} $\zeta$ in 
$\overline\Gamma= \Gamma\cup \geo \Gamma$ equidistributes to $\mu$, as long as $\Gamma$ is {\em nonelementary}. 
More precisely, define 
$$
 \mu_\zeta(s) = \frac{1}{\cP(s)}\sum_{\gamma\in \Gamma}e^{-s\beta(\gamma)}\delta_{\gamma(\zeta)}, \quad s> s_0. 
$$
(Here $\Ga$ acts on itself via the left multiplication.) 
Then, whenever a limiting measure  \index{Gibbs measure $\mu$}
$$
\mu= \lim_{s_j\to s_0+} \mu_\zeta(s_j)
$$ 
on $\geo \Gamma$ exists for some sequence $(s_j)$ and some $\zeta\in \ol\Ga$, it exists for all $\zeta$'s and is independent of $\zeta$. 
See Lemma \ref{lem_Gibbs_independant_basepoint} (for $\zeta\in \Ga$) and Theorem \ref{thm:rank1-equi} (for $\zeta\in \geo \Ga$).  
Every limiting measure appearing this way is called a {\em Gibbs measure} on $\geo \Ga$. From now on, we fix a hyperbolic length function $\beta$,  
a sequence $(s_j)$ and the resulting Gibbs measure $\mu$ on $\geo \Ga$. 

\begin{introrem}
Here is an example of a hyperbolic length function which is especially relevant for our purposes. Let $\cX=G/K$ be a symmetric space of noncompact type, let $\Gamma< G$ be a $P$-Anosov  subgroup with respect to a (self-dual) parabolic subgroup $P< G$. 
Then $\cX$ has a $G$-invariant Finsler metric and associated distance function $d_F$ (see \cite{Kapovich-Leeb-finsler, Dey-Kapovich} for details). Suppose that $x\in \cX$ is a point not fixed by any nontrivial element of $\Gamma$. 
Then, according to \cite{Dey-Kapovich}, for the length function $\beta(\gamma):= 
d_F(x, \gamma x)$ the limit 
$$
\mu= \lim_{s\to s_0+} \mu_\zeta(s)
$$ 
exists and there is no need to work with sequences. See also Theorem \ref{thm:uniqueness} for the existence of hyperbolic length functions 
on arbitrary hyperbolic groups, for which the Gibbs measure is unique. Another example of hyperbolic length functions is the word-lengths defined with respect to finite generating sets of $\Ga$. 
 \end{introrem}

The main result of this book is that the Gibbs measure describes equidistribution of currents (satisfying a suitable genericity condition) 
in flag manifolds $\cF=G/B$ under actions of nonelementary $B$-Anosov subgroups $\Ga<G$, where $G$ is a complex semisimple Lie group and $B< G$ is a Borel subgroup. While similar equidistribution results (with similar proofs) are likely to hold in the case of $P$-Anosov subgroups for general parabolic subgroups $P< G$, we decided to limit ourselves to the case $P=B$ in order to reduce the technicalities.

Let $\Ga< G$ be a nonelementary $B$-Anosov subgroup, $f: \geo \Ga\to \La\subset \cF=G/B$ be an equivariant homeomorphism to the {\em flag-limit set} (see \S~\ref{sec:flag_limit_set}) $\La=\Lambda(\Ga)\subset \cF$ of the group $\Ga$ (the {\em boundary} map or the {\em asymptotic embedding}, see \cite{KLP17}). For each $x\in \cF$ we define 
probability measures 
\begin{equation}
\mu_{x}(s) = \frac{1}{\cP(s)}\sum_{\gamma\in \Gamma}e^{-s\beta(\gamma)}\delta_{\gamma(x)}, \quad s> s_0 
 \end{equation}
 and the Gibbs measure $\mu_{\cF}$ on $\La$, the push-forward of the Gibbs measure on $\geo \Ga$ via the homeomorphism $f$.

 For each point $x\in \cF$ we have a projective subvariety $\Th^{n-1}(x)=\cF\setminus \Opp(x)\subset \cF$,  consisting of all points $y\in \cF$ which are not opposite (antipodal) to $x$ (here and in what follows, $n$ is the dimension of $\cF$). A point $x\in \cF$ is {\em $\La$-generic} if the subset \index{generic point}
 $$
 \La_x= \{\la\in \La: x\in \Th^{n-1}(\la)\}\subset \La
 $$
has zero Gibbs measure. 
 We can now state our equidistribution theorem for points in $\cF$: 
 
 \begin{introthm}\label{thm:main-points}
For every $\Lambda$-generic point $x\in  \cF$, we have  
 $$
\mu_{\cF}= \lim_{s_j\to s_0+} \mu_x(s_j). 
$$ 
 \end{introthm}

In order to appreciate this equidistribution result, consider a discrete cocompact subgroup $\Ga< \PSL(2,\C)$ and 
its diagonal action on the flag-manifold $\cF=S^2\times S^2$ of the group $G= \PSL(2,\C)\times \PSL(2,\C)$. This action 
is ergodic (since the geodesic flow of $\bH^3/\Ga$  is ergodic), hence, almost every $\Ga$-orbit in $\cF$ 
is dense in $\cF$. Nevertheless,  in the Patterson--Sullivan sense, every orbit equidistributes to the Gibbs measure, which is supported on a much smaller subset, namely, the diagonal in $\cF$.

\medskip 
Let us now formulate a generalization of Theorem \ref{thm:main-points} to the case subvarieties $S\subset \cF$ instead of points.  
For each point $x\in \cF$ we define its $m$-dimensional thickening $\Th^m(x)$ as the union of all $m$-dimensional Schubert subvarieties in $\cF$ based at $x$. Schubert subvarieties $X_w\subset \cF$ are indexed by elements of the Weyl group $W$ of $G$ and every Schubert subvariety $X_w$ equals some $\Th_w(x)$, the {\em $w$-thickening} 
of a  point $x$  in $\cF$.  (We refer the reader to Section \ref{sec:Thickenings} for the precise definitions.) 
We say that a $k$-dimensional subvariety $S\subset \cF$ is {\em $\La$-generic} if the subset \index{generic subvariety}
$$
\La_S:= \{\la\in \La: S\cap \Th^{n-k-1}(\la)\ne \emptyset\}= \Th^{n-k-1}(S)\cap \La 
$$
has zero Gibbs measure. (Here  $n$ is the dimension of $\cF$.) 

\begin{introrem}\label{rem:genericity}
In Proposition \ref{prop:genericity} we will prove that a subvariety $S\subset \cF$ is $\La$-generic if and only if $\La_S\ne \La$. In particular, genericity does not depend on the choice of a hyperbolic length function $\beta$. Furthermore, for Zariski dense Anosov subgroups $\Ga< G$ {\em every} subvariety $S\subsetneq \cF$ is $\La$-generic. 
\end{introrem}

Note that if $k=0$, we obtain the $\La$-genericity condition for points in $\cF$. 
Observe also that the sum of dimensions of $\Th^{n-k-1}(x)$ and $S$ is strictly less than the dimension of $\cF$, which means that, 
given $S\subset \cF$, for a.e. $x\in \cF$ the intersection $S\cap \Th^{n-k-1}(x)$ is empty. 
The subvariety $S\subset \cF$ defines the current of integration, $\llbracket S\rrbracket$ in $\cF$ (of the 
bidegree $(n-k,n-k)$). Our main theorem (stated below) describes asymptotics of images of this current under the action of $\Ga$. 
Recall that fundamental classes of Schubert cycles define a basis of homology groups of $\cF$. In particular, given a $k$-dimensional projective subvariety $S\subset \cF$, we have
$$
[S]= \sum_{w\in W, |w|=k} a_w [X_w]\in \coh_{2k}(\cF), a_w\in {\mathbb Z}_+.   
$$ 
Here and in what follows, $|w|$ is the word-length of $w$ with respect to the generating set of $W$ consisting of simple reflections. In order to simplify the notation we set $a_w=0$ for all $w\in W$ of length $\ne k$. Then the above formula turns into 
\begin{equation}\label{eq:decomposition}
[S]= \sum_{w\in W} a_w [X_w]\in \coh_{2k}(\cF).   
\end{equation}

The following two theorems are the main results of our book: 

\begin{introthm}\label{thm:main}Fix a non-elementary $B$-Anosov  subgroup $\Gamma< G$ acting on the flag manifold $\cF=G/B$. 
Also fix a hyperbolic length function $\beta$ on $\Ga$ with its associated critical exponent $s_0$ and Gibbs-measure $\mu_{\cF}$. 
 Then for every $\Lambda(\Ga)$-generic subvariety 
 $S  \subset \cF$, one has 
$$
\lim_{s_j\to s_0+} \frac{1}{{\mathcal P}(s_j)} \sum_{\ga\in \Ga} e^{-s_j\beta(\ga)} \llbracket \ga(S)\rrbracket= \int_{\La}  \sum_{w\in W}  a_w \llbracket\Th_w(\la) \rrbracket d\mu_{\cF}(\la),
$$
where $a_w$ are the coefficients from  \eqref{eq:decomposition} and where 
 the limit is understood as the weak limit in the space $\cD(\cF)$ of currents on $\cF$. 
\end{introthm}

We refer to  equation \eqref{eq:current integration} in Section \ref{sec:Generalities of currents} as well as Section \ref{sec: Anosov Gibbs} 
for the definition of the integral in this theorem. The limit in  Theorem \ref{thm:main}
$$
\int_{\La}  \sum_{w\in W}  a_w \llbracket\Th_w(\la) \rrbracket d\mu_{\cF}(\la), a_w\in \C,$$
will be called the {\em Gibbs currents} on $\cF$ associated with the Anosov subgroups $\Ga< G$ and the Gibbs measures $\mu_{\cF}$.


\medskip
Recall that when $\Gamma$ is Zariski-dense every subvariety is $\Lambda(\Ga)$-generic, see Remark~\ref{rem:genericity}. 
We also obtain an equidistribution result for {\em smooth} currents on $\cF$, where a genericity condition is not needed. The limiting currents are again Gibbs currents. Let $\psi\in \Om^{n-k,n-k}(\cF)$ be a smooth closed 
form on $\cF$. It defines a current $\llbracket\psi\rrbracket$ of {\em bidimension} $(k,k)$  on $\cF$, whose homology 
class is the Poincar\'e dual $PD([\psi])$ of the cohomology 
class $[\psi]\in \coh^{n-k,n-k}(\cF; \mathbb C)$.  We, therefore, set 
\begin{equation}\label{eq_coef_alpha}
PD([\psi])= \sum_{w\in W} a_w [X_w]\in \coh_{2k}(\cF; \mathbb C), a_w\in {\mathbb C},   
\end{equation}
again with the convention that $a_w=0$ when $|w|\ne k$. 

\begin{introthm}\label{thm:main smooth} Fix a non-elementary $B$-Anosov  subgroup $\Gamma< G$ acting on the flag manifold $\cF=G/B$. 
Also fix a hyperbolic length function $\beta$ on $\Ga$, a subsequence $s_j$ converging to the critical exponent $s_0$ and the corresponding Gibbs-measure $\mu_{\cF}$. 
 Then for all $0\leqslant k\leqslant n = \dim \cF$ and for every smooth closed form $\psi$ on $\cF$  of bidegree $(n-k, n-k)$ we have   
$$
\lim_{s_j\to s_0+} \frac{1}{{\mathcal P}(s)} \sum_{\ga\in \Ga} e^{-s_j\beta(\ga)} \llbracket (\ga^{-1})^*\psi\rrbracket= \int_{\La}  \sum_{w\in W}  a_w \llbracket\Th_w(\la) \rrbracket d\mu_{\cF}(\la),
$$
where $a_w$ are the coefficients given in \eqref{eq_coef_alpha} and
where the limit is understood as a weak limit in the space of currents on $\cF$. 
  \end{introthm}

Note that since the  group elements $\gamma$ act as diffeomorphisms, the pushforward of  a current $T$ is $\llbracket \gamma(S) \rrbracket$, if $T = \llbracket S\rrbracket$ is the current of integration on a subvariety $S$ and the pushforward equals $\llbracket(\gamma^{-1})^* \psi\rrbracket$, if $T=\llbracket\psi\rrbracket$ is represented by a smooth closed $(n-k,n-k)$-form $\psi$.

\medskip
{\bf Idea of the proofs of Theorems \ref{thm:main} and \ref{thm:main smooth}.} The overall strategy of the proofs is to consider a 
map
\begin{equation}\label{eq:eta-map}
\eta: \Ga\cup \geo \Ga\to \cD(\cF)
\end{equation}
sending each $\ga\in \Ga$ to the current $\ga_*(\llbracket S\rrbracket)$ for integration currents over algebraic surfaces 
(respectively, $\ga_*(\llbracket \psi \rrbracket)$ for currents corresponding to smooth forms $\psi$) and sending each $\xi\in \geo \Ga$ to the {\em  current}
$$
 \sum_{w\in W}  a_w \llbracket\Th_w(f(\xi))\rrbracket. 
$$
The map $\eta$ is clearly separately continuous on $\Ga$ and on $\geo \Ga$. If the map $\eta$ is overall continuous, then the equidistribution result follows easily, see Proposition \ref{prop:poincare_series_current}. As it turn out, the map $\eta$ is indeed continuous in the case of smooth currents $\llbracket \psi \rrbracket$. The proof of continuity relies upon a quantitative version of a theorem due to Harvey and Lawson, \cite{harvey_lawson_survey,harvey_lawson}, on asymptotics of smooth currents under 
$\mathbb C^*$-actions on algebraic varieties, which we prove in Theorem \ref{thm_extension_harvey_lawson}. Once Theorem \ref{thm_extension_harvey_lawson} is proven, continuity of 
$\eta$ is reasonably straightforward, see the proof of Corollary  \ref{cor_equivariant_continuous}. 

The singular case, of integration currents $\llbracket S\rrbracket$, is harder since the map $\eta$, in general, is not continuous, if $\La_S\ne \emptyset$. 
We prove, however, that since, by the genericity assumption, $\mu(\La_S)=0$, the discontinuity does not prevent convergence in Theorem \ref{thm:main}. Informally speaking, the reason is that the subset $\Ga^c_{P,t}$ of $\Ga$ consisting of {\em bad elements} (those whose inverses are $t$-close to points of $f^{-1}(\La_S)$)  grows (in terms of the length function $\beta$) slowly comparing to the growth of the entire group $\Ga$. This part of the proof is also quite delicate and occupies the bulk of Chapter \ref{sec:Equidistribution of subvarieties} (which, in turn, relies upon results of earlier parts of the book). 

A common theme in the proofs of Theorems \ref{thm:main} and \ref{thm:main smooth} is {\em redressing}: We find a finite subset $E\subset \Ga$ such that for every $\ga\in \Ga$ there exists $g\in E$ satisfying the property that the fixed points of $\tilde\ga= \ga\circ g$ are uniformly separated in $\geo \Ga$. Working with the {\em separated elements} $\tilde\ga$ is simpler than with general elements of $\Ga$ since estimates on the action of $\tilde\ga$ on currents reduce to estimates for elements of uniformly regular 1-parameter subgroups of a fixed maximal torus $T< G$. 

%
%
%
%


\medskip
{\bf Organization of the book.} In Chapter \ref{sec:hyperbolic} we review Gromov-hyperbolic spaces and hyperbolic groups. In Chapter \ref{sec:currents} we discuss 
basics of currents on complex manifolds. In Chapter \ref{sec:counting} we introduce basic counting problems on Gromov-hyperbolic spaces, Poincar\'e series and Gibbs measures. We also define equidistribution problems for currents on smooth manifolds. 
 In   Chapter \ref{section_prelim}  we review symmetric spaces,  flag manifolds and thickenings. In Chapter \ref{sec:Anosov} we discuss Anosov subgroups of semisimple Lie groups. The heart of the book are Chapters \ref{sec:Equidistribution of subvarieties} 
 and \ref{sec:Equidistribution of smooth forms} 
 where we prove the main theorems (Theorems \ref{thm:main} and \ref{thm:main smooth}) on equidistribution of currents on flag manifolds under Anosov group actions.  In Section \ref{sec:examples} of Chapter \ref{sec:Equidistribution of subvarieties}  we describe several classes of examples of Anosov subgroups and subvarieties $S\subset \cF$ satisfying and not satisfying our genericity assumptions. We also give examples of failure of Theorem \ref{thm:main} for subvarieties which are not generic.

The remaining chapters of the book are (Chapters \ref{sec:intersection_Gibbs} 
and \ref{sec:Morse-Smale}), 
strictly speaking, are mostly not needed for the proofs of the equidistribution results. 
However, they  
sheds light on dynamics of Anosov group actions and connect the theory to semigroup dynamics: Iterations of rational maps and 
Morse--Smale flows. Chapter \ref{sec:intersection_Gibbs} is rather exploratory, it gives a criterion for which the wedge product of Gibbs currents can be defined and shows these conditions hold for certain classes of Anosov groups. This construction yields some new measures whose support are invariant under the group. 
In Chapter \ref{sec:Morse-Smale} we interpret the dynamics of Anosov group actions on flag-manifolds $\cF$ in terms of {\em Morse--Smale} flows, via the suspension flow construction combining the geodesic flow on a hyperbolic group with the group action 
on a product of two copies of $\cF$ (see also \cite{Guichard-Wienhard}). It extends Labourie's original definition of Anosov actions in terms of flows: From hyperbolicity of the action on the invariant subset of the suspension flow (corresponding to the limit set) to the {\em Axiom A}-type dynamics on the entire space.  

Finally, in the appendices to the book, we prove some auxiliary results used in the book. In Appendix A (Chapter \ref{sec:appendixA}) we prove a {\em redressing theorem} for hyperbolic groups, taken from a work by Didac Martinez-Granado and Michael Kapovich. In Appendix B (Chapter \ref{sec:appendixB})  we establish several results on location of fixed points of isometries of Gromov-hyperbolic spaces. 
In Appendix C (Chapter \ref{sec:holder}) we prove certain separation result for thickenings of limit points of Anosov subgroups. 
As a consequence, we prove that boundary maps of Anosov subgroups are bi-H\"older. 
Earlier, Konstantinos Tsouvalas proved that boundary maps $f: \geo \Ga\to \cF$ of Anosov subgroups are H\"older. We supplement his result by proving the H\"older property of inverse maps $f^{-1}$, thereby establishing the bi-H\"older property of $f$. 
Results from Appendix \ref{sec:appendix D} are background material in measure and geometric measure theory used throughout the book.

\medskip
{\bf Acknowledgements.} The first author is indebted to Nguyen-Thi, Nguyen-Viet Dang, Thang Nguyen, Anne Moreau, Thomas Gautier for useful comments and suggestions. He has been partly funded by the grants
ANR-24-CE40-1163, ANR-24-CE40-6184 and ANR-25-CE40-0070. The second author is grateful to Subhadip Dey for useful conversations regarding this work. The third author has been supported by the NSF over the years. 
The fourth author is grateful to Jialun Li for some useful discussions. He acknowledges partial support from the grants ANR-24-CE40-3526-01 and ANR-24-CE40-1163. 

This work was initiated during the visit of the first three authors to the Fields institute in Toronto (May 2019). It was continued while the first and fourth authors were Milnor lecturers at the IMS at Stony Brook. 
  It was partly done during the special Advancing Bridges in Complex Dynamics program at MSRI Berkeley (Spring 2022) and during the visit of the third author to Paris and Luminy (July 2025).
 While we were completing  our manuscript, an  independent paper by Bonthonneau-Lefeuvre-Weich \cite{bonthonneau_lefeuvre_weich}
 appeared. It studies, using the thermodynamical formalism,  conformal measures for Anosov actions.  
  
  With great sadness we have to mention that our friend and collaborator Misha Kapovich has passed away on June 16th 2026.

\chapter{Hyperbolic spaces and groups}\label{sec:hyperbolic}

In this chapter we collect some standard definitions and facts about Gromov-hyperbolic spaces and groups. We refer the reader to \cite{BH}, \cite{Buyalo-Schroeder}, \cite{DK18}, \cite{Ghys-book}, \cite{Gromov}, \cite{Vaisala} for details. 


\section{Metric geometry}\label{sec:metric_geometry}


\subsection{Main definitions}


Let $(X,d)$ be a metric space (we will frequently suppress the letter $d$ when referring to a metric space). 
It is called {\em proper} if closed metric balls in $X$ are compact. In this situation, the metric $d$ on $X$ is also called {\em proper}. \index{proper metric} \index{proper metric space}
For a subset $A\subset X$ and a point $x\in X$, define $d(x,A)$, the distance from $x$ to $A$, as \index{$d(x,A)$, distance to a subset $A$} 
$$
\inf \{d(x,a): a\in A\}.
$$ 
Similarly, for two subsets $A, B\subset X$ define \index{$d(A,B)$, minimal distance between subsets $A, B$}
$$
d(A,B):= \inf \{d(a,b): a\in A, b\in B\},
$$
the {\em minimal distance} between $A$ and $B$. 

Given $A\subset X$ 
and a number $r\geqslant 0$, one defines the (closed) 
$r$-neighborhood $N_r(A)$  of $A$ in $X$ as \index{$N_r(A)$, closed $r$-neighborhood of a subset} 
$$
\{x\in X: d(x, A)\leqslant r\}. 
$$
We will use the notation $B(x,r)$ \index{$B(x,r)$, open metric ball} 
for the open metric ball in $X$ centered at $x$ and of the radius $r$ and let $\bar{B}(x,r)$ \index{$\bar{B}(x,r)$, closed metric ball}
denote the closed ball
$$
\bar{B}(x,r)= \{y\in X: d(x,t)\leqslant r\}.
$$  

The {\em Hausdorff distance} $d_H(A_1, A_2)$ between two subsets $A_1, A_2\subset X$, \index{$d_H$, Hausdorff distance} 
is
$$
\inf \{r: A_1\subset N_r(A_2), A_2\subset N_r(A_1)\}. 
$$

An isometric group action of a discrete group $\Ga$ on $X$ is called {\em metrically proper} \index{metrically proper action} 
if for some $x\in X$ (equivalently, every $x\in X$) 
the function $\gamma\mapsto d(x, \gamma x)$ is a proper function on $\Ga$. In other words, for every $R\in \mathbb R_+$,
$$
\{\ga\in \Ga: \ga x \in B(x,R)\}
$$
is a finite subset of $\Ga$.  

A group action $\Gamma\acts X$ is called {\em geometric} if it is isometric, metrically proper and cobounded, i.e. there exists a bounded subset in $X$ whose $\Ga$-orbit is the entire $X$. \index{geometric action}  

A map $c$ from an interval $J\subset {\mathbb R}$ to $X$ is called {\em geodesic} if \index{geodesic in a metric space} 
$$
d(c(s), c(t))=|s-t|
$$
for all $s, t\in J$. In what follows, we will conflate geodesics with their images in $X$. A {\em geodesic segment} \index{$xy$, geodesic segment} 
between points $x, y\in X$, denoted $xy$ or $[x,y]$, is the image of a geodesic segment $c: [a,b]\to X$ such that $c(a)=x, c(b)=y$. One also says that the segment $xy$ {\em connects} $x$ and $y$ and that $x, y$ are the {\em end-points} of $xy$.  In general, geodesic segments connecting $x, y$ are not uniquely determined by $x, y$, so our notation is somewhat ambiguous. A {\em geodesic ray} in $X$ is a geodesic map 
$c: [0,\infty)\to X$. Such a ray is said to {\em emanate} from $c(0)$.

A metric space $(X,d)$ is called {\em geodesic} \index{geodesic metric space} 
if  every two points $x, y\in X$ are connected by a geodesic segment. 

\subsection{Quasigeodesics and quasiisometries} 
\label{sec:quasigeod_quasiiso}

A {\em quasigeodesic} (more precisely, an $(L,A)$-quasigeodesic for some $L\geqslant 1, A\geqslant 0$) \index{quasigeodesic} 
in $(X,d)$ is a map to $c: J\to X$ such that 
$$
 L^{-1}|s-t| -A\leqslant d(c(s), c(t))\leqslant L|s-t| +A 
$$
for all $s, t\in J$, where $J$ is an interval in $\R$. A quasigeodesic is said to be {\em complete} if $J=\R$. \index{complete quasigeodesic} 

A metric space 
$(X,d)$ is called {\em quasigeodesic} if there exist constants $(L,A)$ such that every two points $x, y\in X$ are connected by an \index{quasigeodesic space} 
$(L,A)$-quasigeodesic. A metric space $(X,d)$ is called {\em roughly geodesic} if every two points $x, y\in X$ are connected by a  \index{roughly geodesic space} $(1,A)$-quasigeodesic for some uniform constant $A$. 

More generally, a map $f: (X,d_X)\to (Y,d_Y)$ is called an 
{\em $(L,A)$-quasiisometric embedding} \index{quasiisometric embedding}
if 
 $$
 L^{-1} d_X(p,q) -A\leqslant d_Y(f(p), f(q))\leqslant Ld_X(p,q)+A  
$$
for all $p, q\in X$. A quasiisometric embedding is said to be a {\em quasiisometry} \index{quasiisometry} 
if there exists $D<\infty$ such that 
 $Y=N_D(f(X))$. 
 
 \begin{example}
 Consider the graph $G$ of the function $|x|$ in $\R^2$ with the restriction of the Euclidean distance function. Then $G$ is a quasigeodesic metric space 
 (since the map $x\mapsto (x, |x|)$ is a quasiisometry $\R\to G$) but is not a roughly geodesic space.  
 \end{example}

\subsection{Metrics on graphs and groups}

Let $H$ be a connected graph. By abuse of terminology, we will conflate $H$ and its geometric realization. 
The {\em graph-metric} $d_H(x,y)$ on  $H$ is defined as follows. Let $v, w$ be two vertices in $H$. Then $d_H(v,w)$  
is the minimum of edge-lengths of combinatorial paths in $H$ connecting $v$ and $w$. This metric extends (uniquely) to the rest of $H$ so that each edge is isometric to the unit interval $[0,1]$ in $\R$ with its standard metric.

\medskip
The most important for us example of a graph-metric will be the one of a Cayley graph of a group $\Ga$, inducing the 
{\em word metric} on $\Ga$. Let $\Ga$ be a group with a generating set $S$. This defines the {\em word-metric} $d_S$ on $\Ga$ and the {\em word-length} \index{$d_S$, word-metric} 
$|\ga|=d_S(1_\Ga, \ga)$. The word-metric $d_S$ comes from the graph-metric $d_S$ on the Cayley graph $Cay_\Ga$ \index{$Cay_\Ga$, Cayley graph} 
of $\Ga$ associated with the generating set $S$. A geodesic ray $c$ in $Cay_\Ga$ is said to be {\em normalized} if $c(0)=1_\Ga$. 
\index{normalized ray} 

In the case of finite generating sets, different word-metrics on $\Ga$ are bilipschitz to each other via the identity map. 
We will be only using word-metrics and word-lengths associated with finite generating sets. 


\subsection{Distortion of metrics and maps}\label{sec:Dilatation of metrics and maps}


Given two metrics $\rho, \rho'$, metrizing a topological space $Z$, we define the {\em distortion} or {\em infinitesimal Lipschitz constant} of these metrics at $z\in Z$ as 
\begin{equation}\label{eq:dilatation1}
Lip(\rho',\rho)(z) :=\limsup_{x, y\to z} \frac{\rho'(x, y)}{\rho(x, y)}\in [0,\infty],
\end{equation}
where the limit is taken over pairs of distinct points $x, y\in Z$. More generally, we define the {\em  infinitesimal Lipschitz constant}  of a map $f: (Z,\rho)\to (Z',\rho')$ as 
\begin{equation}\label{eq:dilatation2}
Lip_f(z):= \limsup_{x, y\to z} \frac{\rho'(f x, fy)}{\rho(x, y)}\in [0,\infty]
\end{equation}

This is a metric analogue of the norm of the differential of a smooth map of Riemannian manifolds at a point. Sometimes one can replace $\limsup$ in these definitions 
by the ordinary limit. Then we will use the notation 
$$
\left. \frac{d\rho'}{d\rho}\right\vert_z
$$
for the distortion of metrics and $J_f(z)$ for the  infinitesimal Lipschitz constant of maps. 

Two metrics $\rho, \rho'$ on $Z$ are said to be {\em conformal to each other} if for 
every $z\in Z$ we have 
$$
0< \left. \frac{d\rho'}{d\rho}\right\vert_z= \lim_{(x,y)\to (z,z), x\ne y} \frac{\rho'(x,y)}{\rho(x,y)}<\infty.   
$$ 





A homeomorphism $f: (Z,d)\to (Z',d')$ of two metric spaces is said to be {\em quasisymmetric} \index{quasisymmetric map} if there exists 
an increasing function $\eta: [0,\infty)\to [0,\infty)$ such that for any triple of distinct points $x, y, z\in Z$ we have
$$
\frac{d'(f(x), f(y))}{d'(f(x), f(z))} \le \eta \left( \frac{d(x, y)}{d(x, z)}\right). 
$$
It turns out that if $f$ is quasisymmetric, so is $f^{-1}$, see e.g. \cite[Lemma 5.2.13]{Buyalo-Schroeder}.

\section{Gromov-hyperbolic spaces}\label{sec:Gromov-hyperbolic spaces}

\subsection{General facts}

For points  $x,y,z$ in $X$, the {\em Gromov product} of $y$ and $z$ with respect to $x$ is \index{$(y,z)_x$, Gromov product} 
\[(y,z)_x=\frac{1}{2}\left(d(x,y)+d(x,z)-d(y,z)\right).\]

A {\em geodesic triangle} in a metric space $(X,d)$ is a triple of geodesic segments $xy, yz, zx$. Such triangles are denoted $\Delta xyz$. 
\index{$\Delta xyz$, geodesic triangle} 
The segments $xy, yz, zx$ are called the {\em sides} of  $\Delta xyz$. A geodesic triangle $\Delta xyz$ is called {\em $\delta$-slim} if each point of a side of $\Delta xyz$ is within distance $\delta$ from a point in the union of the two other sides. \index{$\delta$-slim triangle} 

There are two slightly different definitions of hyperbolicity for metric spaces. The first one, which we will adopt as our main definition, is due to Rips:

\begin{defi}\index{$\delta$-hyperbolic metric space} 
A geodesic metric space $(X,d)$ is said to be {\em $\delta$-hyperbolic} (in the sense of Rips) if every geodesic triangle in $X$ is $\delta$-slim. 
\end{defi}

The second definition works for metric spaces which need not be geodesic, it is due to Gromov:

\begin{defi}\index{$\delta$-hyperbolic metric space}
A geodesic metric space $(X,d)$ is said to be {\em $\delta'$-hyperbo\-lic} (in the sense of Gromov) if every  quadruple of points $x_0, x, y, z\in X$ satisfies the inequality:
\begin{equation}\label{eq:Gh}
(x,z)_{x_0}\geqslant \min \left( (x,y)_{x_0}, (y,z)_{x_0}\right) - \delta'. 
\end{equation}
\end{defi}

The two notions of hyperbolicity are essentially the same but one has to modify the hyperbolicity constants:

\begin{enumerate}
\item If $(X,d)$ is $\delta$-hyperbolic in the sense of Rips, then it is $\delta'$-hyperbolic in Gromov's sense for $\delta'=3\delta$. 

\item If $(X,d)$ is a geodesic metric space which is $\delta'$-hyperbolic in Gromov's sense, then it is $\delta$-hyperbolic in the sense of Rips for $\delta=2\delta'$. 
\end{enumerate} 

 A space is called {\em Gromov hyperbolic} 
 if it is $\delta$-hyperbolic for some $\delta<\infty$ in the sense of either definition.  
 Throughout the book, we will refer to geodesic metric spaces $\delta$-hyperbolic in the sense of Rips 
simply as $\delta$-hyperbolic and reserve the symbol $\delta'$ for the constant of 
hyperbolicity in Gromov's sense. Gromov's definition will be used for estimating 
Gromov products and distances between points of the visual compactification of 
Gromov hyperbolic spaces.

Geodesic segments in geodesic $\delta$-hyperbolic spaces are {\em almost} determined by their end-points: For any two segments $c_1, c_2$ with the end-points $x, y$, $d_H(c_1, c_2)\leqslant \delta$. Thus, the ambiguity in the notation for geodesic segments in Gromov-hyperbolic spaces is mostly harmless. 

\medskip 
The {\em quasiconvex hull}, \index{$\hull(A)$, quasiconvex hull} 
$\hull(A)$, of a subset $A\subset X$ is the closure of the union of all geodesic segments in $X$ connecting points of $A$. 
A subset $A\subset X$ is said to be $\eps$-{\em quasiconvex}  if $d_H(A, \hull(A))\leqslant \eps$. A subset of $X$ is said to be quasiconvex if it is $\eps$-quasiconvex for some $\eps<\infty$. \index{quasiconvex subset}

Given a closed $\eps$-quasiconvex subset $A\subset X$ in $X$, we define the {\em nearest-point projection} $P_A: X\to A$ \index{$P_A$, nearest-point projection} 
by sending $x\in X$ to a point 
$\bar{x}\in A$ such that $d(x, \bar{x})=d(x,A)$. Note that the point $\bar{x}$ need not be unique. However, any two nearest points are within distance
$2\eps+ 9\delta$ from each other. More precisely, for any two points $x, y\in X$
\begin{equation}\label{eq:Lip-projection} 
d(P_A(x), P_A(y))\leqslant 2 d(x,y) + 2\eps+ 9\delta,
\end{equation}
see e.g. \cite[Lemma 11.53]{DK18}, i.e. $P_A$ is a {\em coarsely Lipschitz} map.  

\medskip 
Suppose that $X$ is $\delta$-hyperbolic and $c$ is a geodesic connecting $y$ and $z$. Then 
\begin{equation}\label{eq:gromov-inequality} 
(y,z)_x - \delta \leqslant d(x, c)\leqslant (y,z)_x, 
\end{equation} 
see \cite[Lemma 11.22]{DK18}. 

\begin{lem}\label{lem:alt-hyp}
A metric space $(X,d)$ is ${\delta'}$-hyperbolic in Gromov's sense if and only if for every quadruple of points $x, y, z, w$ in $X$ we have
$$
d(w, x) + d(y, z) \le \max\{d(w, y) + d(x, z), d(x, y) + d(w, z)\} + 2{\delta'}, 
$$
see the formula (2.12) in \cite{Vaisala}. 
\end{lem}
\begin{proof} The proof amounts to rewriting Gromov's condition for the points $x, y, z, w$:
$$
(x,y)_w\ge \min \{(x,z)_w, (y,z)_w\} - {\delta'} \iff $$
$$
-d(x,y) + d(x,w) + d(y,w) \ge \min \{ d(x,w) + d(z,w) - d(x,z), d(y,w) + d(z,w) - d(y,z)\} -2{\delta'} \iff
$$
$$
d(x,y) - d(x,w) - d(y,w) \le \max \{ -d(x,w) - d(z,w) + d(x,z), - d(y,w) - d(z,w) + d(y,z)\} + 2{\delta'} \iff 
$$
$$
d(x,y) + d(z,w)\le \max \{ -d(x,w)  + d(x,z), - d(y,w)  + d(y,z)\} + (d(x,w) + d(y,w)) + 2{\delta'}=$$
$$
 \max\{d(w, y) + d(x, z), d(x, y) + d(w, z)\} + 2{\delta'}. 
$$
\end{proof}

One way to visualize the 4-point inequality in the lemma is to consider a formal "tetrahedron" with vertices $x, y, z, w$ in $X$. It has three pairs 
of disjoint "edges" $(xy, zw)$, $(xz, yw)$ and $(xw, yz)$. Each of the three pairs has "total length" equal to the sum of the corresponding distances, resulting in three 
numbers $r, s, t$. Then the inequality becomes 
$$
r\leqslant \max \{s, t\} +2{\delta'}. 
$$

Slimness of geodesic triangles implies slimness of geodesic polygons:

\begin{lem}\label{lem:thinpolygons}
Suppose that $P=x_1x_2...x_n$ is a geodesic $n$-gon in a $\delta$-hyperbolic space in the sense of Rips. Then 
$P$ is $(n-2)\delta$-slim: Every side of $P$ is contained in the $(n-2)\delta$-neighborhood of the union of the other sides. 
\end{lem}
\begin{proof} Let us prove that the side $x_1x_n$ is contained in the $(n-2)\delta$-neighborhood of the union of the other sides. Draw a geodesic $x_1 x_{n-1}$, subdividing $P$ in a triangle $T=x_1x_{n-1}x_n$ and $(n-1)$-gon $P'= x_1....x_{n-1}$. Inductively, $x_1 x_{n-1}$ is contained in the $(n-3)\delta$-neighborhood of the union of the sides 
$x_ix_{i+1}$, $i=1,...,n-2$, of $P'$. The claim now follows from $\delta$-slimness of the triangle $T$. 
\end{proof}

\subsection{Ideal boundaries}\label{sec:ideal boundaries}

We next define the {\em ideal} (or, {\em visual}) boundary $\geo X$ \index{$\geo X$, visual boundary} 
of a Gromov--hyperbolic space $X$. Fix a base-point $x_0\in X$. 
A sequence $(x_n)$ in $X$ is called a {\em Gromov--sequence} if \index{Gromov--sequence} 
$$
\liminf_{m, n\to\infty} (x_n, x_m)_{x_0}= \infty
$$
Two Gromov--sequences $(x_n), (y_n)$ are {\em equivalent} if the combined sequence 
$(z_n)$ defined by $z_{2n-1}=x_n$, $z_{2n}=y_n$ is also a Gromov--sequence. 
 Elements of $\geo X$ are equivalence classes $[x_n]$ of {\em Gromov--sequences}. 
 One verifies that these notions are independent of the choice of the base-point $x_0$.   Define the set $\ol{X}:= X\cup \geo X$ (later on, we will topologize this set).

 A geodesic ray $c: \R_+\to X$ is said to be {\em asymptotic} to $\xi\in \geo X$ if for some (equivalently, every) diverging sequence $t_n\in \R_+$, the sequence $(x_n)=(c(t_n))$ is a Gromov-sequence representing the point $\xi$, $\xi=[x_n]$. 
 A geodesic ray emanating from a point $x\in X$ and asymptotic to $\xi\in \geo X$ is denoted $x\xi$. \index{$x\xi$, geodesic ray}  Two geodesic rays $c_1, c_2$ are {\em equivalent} if they are within finite Hausdorff distance from each other. If $X$ is a proper geodesic Gromov-hyperbolic space, one can equivalently define the visual boundary of $X$ as the set of equivalence classes of geodesic rays in $X$.      
 
 By analogy with geodesic triangles, one defines {\em ideal geodesic triangles} with vertices in $\ol{X}$ by taking unions 
 of finite geodesic segments, geodesic rays and complete geodesics connecting points of  $\ol{X}$; these geodesics 
 are the {\em sides} of ideal triangles and the points connected by them are the {\em vertices}. (Note that we are not claiming that such infinite geodesics necessarily exist.)  The {\em ideal} vertices of the ideal triangles are those which belong to $\geo X$. By abusing the notation, we will refer to ideal triangles as {\em triangles in $X$} since we are primarily interested in their points which belong to $X$ rather than $\geo X$. We {\em do not} allow two ideal vertices of an ideal triangle to be equal. 
 
\begin{lem}\label{lem:idealthinness}
Suppose that $(X,d)$ is a $\delta$-hyperbolic geodesic space in the sense of Rips and $T=xyz$ is an ideal triangle in $X$. Then:

1. If $T$ has one ideal vertex, then $T$ is $2\delta$-slim: Every side is contained in the $2\delta$-neighborhood of the union of the two other sides. 

2. If $T$ has two ideal vertices, then it is $3\delta$-slim.

3. If $T$ has three ideal vertices then it is $4\delta$-slim.

Furthermore, every geodesic quadrilateral $xy\xi \eta, \xi, \eta\in \geo X$, $x, y\in X$ is $3\delta$-slim.  
\end{lem} 
\begin{proof} We will prove only the third statement since the proofs of the rest are similar. Take sequences 
$x_n, y_n\in xy$, $y'_n, z_n\in yz$, $z_n', x'_n\in zx$ such that
$$
[x_n]=[x'_n]=x, \quad [y_n]=[y'_n]=y, \quad [z_n]=[z'_n]=z. 
$$
We get geodesic hexagons $P_n=x_ny_ny'_nz_nz'_nx'_n$ in $X$. By Lemma \ref{lem:thinpolygons}, each $P_n$ is $4\delta$-slim. Pick a point $p\in xy$. Then $p\in x_ny_n$ for all sufficiently large $n$. Since
$$
\lim_{n\to\infty} (x_n,x'_n)_p= \lim_{n\to\infty} (y_n,y'_n)_p= \lim_{n\to\infty} (z_n,z'_n)_p=\infty,
$$
by the slimness of $P_n$, for all sufficiently large $n$, $p$ lies in the $4\delta$-neighborhood of the union $y'_nz_n \cup z'_n x'_n\subset yz\cup zx$. The claim follows. 
\end{proof}

One extends the notion of quasiconvex hulls to subsets of $\geo X$: Given $A\subset \geo X$, its quasiconvex hull $\hull A$ in $X$ is defined as the union of all biinfinite geodesics in $X$ which are asymptotic (in both directions) to points of $A$.

\medskip 
For $p,q\in \geo X$, the Gromov product is  
\[
(p,q)_x=\inf \liminf_{i,j\to\infty}(y_i,z_j)_x\in [0,\infty],  
\]
 \index{$(y,z)_x$, Gromov product}
where  the infimum is taken over all sequences of points $(y_n)_n, (z_n)_n$ in $X$ 
representing, respectively,  $p, q$. Similarly, one defines 
the Gromov product $(y,q)_x$ for $p\in \geo X$ and $y\in X$:
{\[(y,q)_x=\inf \liminf_{j\to\infty}(y,z_j)_x\]}
where  the infimum is taken over all sequences of points $(z_n)_n$ representing $q$. Note that for $y\in X$ and fixed $x\in X, q\in \geo X$, $(y,q)_x$ is 
a continuous function of $y$. Furthermore, $(p,q)_x=\infty$ if and only if $p=q\in \geo X$. 
 

\begin{rem}
There are several slightly different definitions of Gromov-prod\-ucts of points in the literature; their differences are uniformly bounded by $2\delta'$, 
see \cite{Vaisala} for a comparison of these definitions. Our definition follows the one in \cite{CDP}, \cite{Kapovich-Sardar} and \cite{Vaisala}. 
We will be frequently using the inequality
$$
\sup \liminf_{i,j\to\infty}(y_i,z_j)_x - (p,q)_{x} \leqslant 2\delta', 
$$
where $p=[y_i], q=[z_j]$ (if $p, q\in \geo X$) or $p=y_i, q=[z_j]$ (if $p\in X, q\in \geo X$). 
\end{rem}

With this definition of Gromov product on $\ol{X}$ one still has the inequality \eqref{eq:Gh}:
\begin{equation}
\label{eq:Gh1}
(x,z)_{x_0}\geqslant \min \left( (x,y)_{x_0}, (y,z)_{x_0}\right) - \delta'
\end{equation}
for $x_0\in X$ and $x, y, z\in \ol{X}$. 

The inequality \eqref{eq:gromov-inequality} implies that for every geodesic ray $c=y\zeta$, $\zeta\in \geo X$ we have 
\begin{equation}\label{eq:gromov-inequality1} 
(y,\zeta)_x - \delta \leqslant d(x, c)\leqslant (y,\zeta)_x+2\delta'\leqslant (y,\zeta)_x+6\delta. 
\end{equation}

\begin{lem}\label{lem:0}
$$
(y,\zeta)_x\leqslant d(x,y)
$$
for all $x, y\in X$, $\zeta\in \geo X$. 
\end{lem}
\begin{proof} Take a sequence $z_n\in X$ converging to $\zeta$ such that 
$$
\lim_{n\to\infty} (y,z_n)_x= (y,\zeta)_x.
$$ 
By the triangle inequality, $d(x,y)\geqslant (y,z_n)_x$. Taking the limit as $n\to \infty$, lemma follows. 
\end{proof}

Note also that (even for general metric spaces) we have
$$
|(x,z)_{x_0} - (z,y)_x|\leqslant d(x,y),
$$
see \cite[Lemma 2.8]{Vaisala}. Therefore, for every $\zeta\in \geo X$ and 
$x_0, x, y\in X$, by taking a sequence $(z_n)$ in $X$ such that $\lim_{n\to\infty} (y,z_n)_{x_0}= (y,\zeta)_{x_0}$, 
we get
\begin{equation}\label{eq:Lip}
|(x, \zeta)_{x_0} - (y,\zeta)_{x_0}|\leqslant d(x,y)+2\delta'\leqslant d(x,y)+6\delta,
\end{equation}

The space 
$$
\ol{X}=X\cup \geo X
$$
has a topology such that $X$ is topologically embedded in $\ol{X}$ and a sequence $(x_n)$ in $X$ converges to $z\in \geo X$ if and only if $(x_n)$ is a Gromov--sequence and $[x_n]=z$. If $X$ is a 
proper metric space, then $\ol{X}$ is compact. Note that the Gromov-product on $\ol{X}$ is, in general, not continuous (we will discuss this in more detail in \S \ref{sec:strong hyperbolicity}). 

\begin{figure}[htbp]
   \centering
   \includegraphics{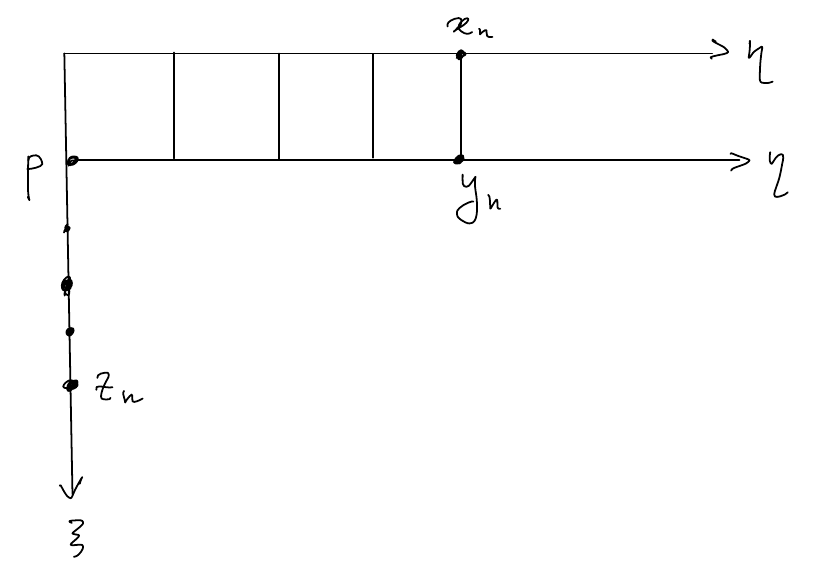} 
   \caption{Discontinuity of the Gromov-product.}
   \label{fig:figurexyz}
\end{figure}

\begin{example}
Let $X$ be the infinite graph as in Figure \ref{fig:figurexyz} with its graph-metric. The visual boundary of $X$ is a 2-point 
set $\{\xi, \eta\}$ with $\xi=[z_n], \eta=[x_n]=[y_n]$, where $(x_n), (y_n), (z_n)$ are Gromov sequences in $X$ such that $d(p,x_n)=n+1, d(p, y_n)=d(p,z_n)=n$. Then 
$$
(x_n,z_n)_p=1, \quad (y_n,z_n)_p=0, \quad (\xi, \eta)_p=0. 
$$
\end{example}

\medskip 
Below is a brief description of the topology of $\ol{X}$. 
For points $x\in X$ we use metric neighborhoods to define the basis of topology.

\medskip
For $\zeta\in \geo X$ and $\alpha, s>0$ we define $U(\zeta,t)$ to be the set of points $q\in \ol{X}$ such that \index{$U(\zeta,t)$, a neighborhood 
of a point at infinity} 
\begin{equation}
\label{eq:U-nbd}
(\zeta,q)_x> - \log t.\end{equation}
Then $\{U(\zeta,t): t >0\}$ is a basis of open neighborhoods of $\zeta$ in $\ol{X}$.  Similarly, for a compact subset $P\subset \geo X$ define 
\begin{equation}
\label{eq:U-nbd1}
U(P,t) = \{ z \in \ol{X} | \sup_{p\in P} (z, p)_{x_0} > - \log(t)  \}= \bigcup_{p\in P} U(p,t).
\end{equation}
Thus, subsets $U(P,t)$ form a basis of neighborhoods of $P$ in $\ol{X}$. 

\begin{rem}
With the above definition of Gromov product for points in $\geo X$, 
the subsets $U(\zeta,t)$ 
are open in $\ol{X}$. This need not be the case 
if we were to use the supremum in the definition of the Gromov product for points at infinity 
instead of the infimum. Furthermore, for fixed $x$ and $z$, the function $y\mapsto (y,z)_x$ is lower semicontinuous. 
\end{rem}

\begin{defi}\label{def:semiproper}
 A geodesic Gromov-hyperbolic space $X$ is {\em semiproper} if: \index{semiproper space}
\begin{enumerate}
\item $\geo X$ is compact. 

\item $X$ satisfies the {\em visibility property}: 
Every two distinct points in $\ol{X}$ are connected by a geodesic in $X$.  
\end{enumerate}
\end{defi}

For instance, every proper geodesic Gromov-hyperbolic space is semiproper, but the converse is, of course, false. 
Semiproper spaces appear naturally when one considers injective hulls of certain Gromov-hyperbolic spaces, see \S \ref{sec:hyp-injective hull}.  

\subsection{Visual metrics}\label{sec:Visual metrics}

Fix a base-point $x_0\in X$.  There exist constants $\alpha,k_1,k_2>0$  and a {\em visual metric} \index{visual metric}
$d_{\alpha}=d_{x_0,\alpha}$ on $\geo X$ metrizing $\geo X$, meaning that for all $p, q\in \geo X$, 
\begin{equation}\label{eq:visual}
k_1\exp(-\alpha (p,q)_{x_0})\leqslant d_{\alpha}(p,q)\leqslant k_2\exp(-\alpha (p,q)_{x_0}),
\end{equation}
see e.g \cite[Section 5]{Vaisala}. The {\em naive} definition of a visual metric is to take 
$$
\rho_\alpha(p,q)= \exp(-\alpha (p,q)_{x_0})$$ 
for some $\alpha>0$. It turns out that for {\em some} Gromov-hyperbolic spaces this does result in a metric 
topologizing $\geo X$ if $\al$ is sufficiently close to $0$: This is the case for  
 {\em strongly hyperbolic spaces} (e.g. $CAT(-1)$ spaces), 
 see Section \ref{sec:strong hyperbolicity}. In general, constructions of metrics $d_\alpha$ are a bit more involved. One such construction is to take a sufficiently small $\al$ 
 and then define  $d_{\alpha}$ as
 $$
 d_{\alpha}(p,q)= \inf \sum_{i=0}^n \rho_\al(p_i, p_{i+1})
 $$
where the infimum is taken over all finite sequences $p_0=p, p_1,...,p_n, p_{n+1}=q$ in $\geo X$ (the number $n$ is not fixed).

\medskip 
In what follows, we will fix $\alpha$ and a visual metric $d_\alpha$ (does not matter which one) and 
use the notation $d_\infty$ for a metric on $\ol{X}$ (topologizing $\ol{X}$) extending the visual metric $d_\alpha$. Such an extension 
always exists due to the Hausdorff extension theorem, see e.g. \cite{h-extension}. 

In view of the interpretation of the Gromov product $(p,q)_{x_0}$ in terms of the distance from $x_0$ to the geodesic $pq$ in $X$ connecting 
distinct points $p$ and $q$ in $\geo X$, we obtain:
\begin{equation*} 
k_1\exp(- \alpha d(x_0, pq)) \leqslant d_\infty(p,q) \leqslant  k_2 \exp (- \alpha d(x_0, pq) +\delta).  
\end{equation*}
If $X$ is proper, a similar property holds for all points $p, q\in \ol{X}$:
\begin{equation}\label{eq:dinfty} 
k'_1\exp(- \alpha d(x_0, pq)) \leqslant d_\infty(p,q) \leqslant  k'_2 \exp (- \alpha d(x_0, pq)),  
\end{equation}
provided that $p, q$ belong to a sufficiently small neighborhood of $\geo X$ in $\ol{X}$.  


\subsection{Busemann functions}\label{sec:Busemann functions}

Given a base-point $p\in X$ and a point $\zeta\in \geo X$ one defines the {\em Busemann function} $b_\zeta(y), y\in X$, as
\begin{equation}\label{eq:busemann} 
b_\zeta(y):= \lim\inf_{z\to \zeta} (d(y, z)- d(p,z))
\end{equation}
where the limit is taken with $z\in X$ converging to $\zeta$ in $\ol{X}$. There exists $C=C(\delta)$ such that
$$
|b_\zeta(y)- \lim\sup_{z\to \zeta} (d(y, z)- d(p,z))|\le C. 
$$
\index{Busemann function} 

\begin{example}\label{ex:Busemann}
Let $\Ga$ be the free group on a free generating set $S$. We let $A=S\cup S^{-1}$ denote the symmetrization of $S$. Let $X$ be the corresponding Cayley tree. A normalized ray in $X$ is a geodesic ray emanating from the neutral element $1_\Ga$ in $\Ga$; such rays are represented by reduced half-infinite words
$$
a_1a_2....a_k....
$$
where each $a_i$ is an element of $A$. Let $\xi\in \geo \Ga$ be the ideal point represented by a ray as above. Consider an element $\ga\in \Ga$ represented by a reduced word $a_1a_2...a_mc_1c_2...c_n$, where $m\ge 0$, $n\ge 0$ and each $c_i$ is an element of $A$. Let $b_\xi$ denote the Busemann function on $X$ normalized to vanish at $1_\Ga$. Then for $x=a_1...a_m$
$$
b_\xi(\ga)= d(\ga, x) - |x| = |\ga^{-1}x| - |x|= n-m. 
$$
\end{example}

We will also define the {\em Busemann almost cocycle} $b_\zeta(x,y):= b_\zeta(x)-b_\zeta(y)$, equivalently, 
\begin{equation}\label{eq:busemanncocycle} 
b_\zeta(x,y):= \lim\inf_{z\to \zeta} (d(x, z)- d(y,z)). 
\end{equation}
Note that if $y=p$ is the normalization point of Busemann functions, then $b_\zeta(x,p)=b_\zeta(x)$. 

\begin{lem}\label{lem:almost cocycle}
$\half b_\zeta(x,y) + \half b_\eta(x,y)= (\zeta,\eta)_x- (\zeta,\eta)_y + C$, where $C$ is a uniformly bounded function with bounds depending only on $\delta$. 
\end{lem}
\begin{proof} The claim follows from the definitions of Gromov products for points at infinity and Busemann functions, with the error terms $C$ coming from the difference between upper and lower limits in these definitions. \end{proof}

{\em Closed horoballs} in $X$ are sublevel sets of Busemann functions:
$$
B(\zeta,t)= \{y\in X: b_\zeta(y)\le t\}, t\in \R. 
$$
Such horoball is said to be {\em centered} at the point $\zeta$. \index{horoball} \index{horosphere}  
Similarly, {\em horospheres} in $X$ are level sets of Busemann functions.

\subsection{Boundary extension of quasiisometries}

{\kap 

It is a fundamental fact of hyperbolic geometry that quasiisometries of Gromov-hyperbolic spaces extend to their visual boundaries:

\begin{thm}\label{thm:boundary extension}
Let $h: (X,d)\to (X',d')$ be an $(L,A)$-quasiisometry of roughly geodesic Gromov-hyperbolic spaces. Then there is an extension of $h$ to a map $(\geo X, d_\infty)\to (\geo X', d'_\infty)$ is an $\eta$-quasisymmetric bi-H\"older homeomorphism  between their visual boundaries equipped with visual metrics. Here $\eta$ depends only on the constants $L, A$, hyperbolicity constants of the spaces $(X,d), (X',d')$ and the parameters of visual metrics. The extension is continuous in the sense that if $x_n\in X$ is a Gromov-sequence representing a point $\xi\in \geo X$ then $f(x_n)$ is a Gromov-sequence in $X'$ representing $f(\xi)$. 
\end{thm}

See e.g. \cite[Theorem 5.2.17]{Buyalo-Schroeder} and \cite{Vaisala}. 

\medskip


For isometries of Gromov-hyperbolic spaces one gets a bit stronger result.   The following was proven in \cite[Proposition 3.1]{Coo93}:

\begin{prop}
Let $X$ be a $\delta$-hyperbolic space. Fix a base-point $p\in X$ and a visual metric $d_\al=d_{\al,p}$ 
defined with respect to this base-point. Then for every isometry $\ga$ of $X$ there is a constant $C$ such that 
$$
C^{-1} e^{-\al |b_z(\ga^{-1}(p))|} \le Lip_\ga(z)\le C e^{\al |b_z(\ga^{-1}(p))|}$$
for all $z\in \geo X$. 
\end{prop}

Here $Lip_\ga$ denotes the distortion function of $\ga$ at $z\in \geo X$ with respect to the metric $d_\al$ (see \S \ref{sec:Dilatation of metrics and maps}).

In the setting of {\em strongly hyperbolic spaces} one gets even stronger {\em conformality} property for $\ga$, 
see Lemma \ref{lem:Lip Jacobian} in \S \ref{sec:Conformality properties}.

}

\subsection{Classification of isometries}

\begin{defi}\label{def:axis}
Let $X$ be a geodesic $\delta$-hyperbolic space. 

1. An isometry $\ga$ of $X$ is said to be {\em loxodromic}\index{loxodromic isometry} 
 if it has two distinct fixed points in $\geo X$, one attractive and one repelling 
($\xi^+, \xi^-$ respectively):
$$
\lim_{n\to\infty} \ga^n(x)=\xi^+, \forall x\in \ol{X}\setminus \{\xi^-\}, 
$$
$$
\lim_{n\to\infty} \ga^{-n}(x)=\xi^-, \forall x\in \ol{X}\setminus \{\xi^+\}. 
$$
Equivalently, for some (equivalently, every) $x\in X$ we have 
$$
\liminf_{n\to\infty} \frac{d(x, \ga^n(x))}{n}>0. 
$$

2. An isometry $\ga$ is said to be {\em elliptic} \index{elliptic isometry} 
 if it generates a subgroup of $\Isom(X)$ which has a bounded  orbit. Such an isometry may or may not have a 
fixed point in $\geo X$. 

3. An isometry $\ga$ is said to be {\em parabolic} \index{parabolic isometry} 
if it is not elliptic and has a unique fixed point in $\geo X$. Equivalently, $\ga$ is such that 
$$\lim_{n\to\infty} d(\ga^n(x), x)=\infty$$ but this divergence is sublinear.  

4. An $(L,D)$-{\em quasi-axis} of an isometry $\gamma$ of $X$ is a complete $(L,D)$-quasi-geodesic $\alpha$ in $X$ \index{quasi-axis} 
invariant under the action of $\gamma$, and such that the ideal end-points $\alpha(\pm \infty)$ are fixed by $\gamma$. 
\end{defi}

Note that a quasi-axis $\alpha$ of a loxodromic isometry is necessarily asymptotic to the fixed points $\xi^\pm\in \geo X$ of $\gamma$. 
Given an isometric action $\Ga\acts X$ on a semiproper gromov-hyperbolic space, 
there exist constants $L, D$ such that every loxodromic element $\ga\in \Ga$ has an 
$(L,D)$-quasi-axis. More precisely, let $\xi^\pm$ denote the attractive/repelling fixed points of $\ga$ in $\geo X$. Then any two geodesics in $X$ connecting 
$\xi^\pm$ are within Hausdorff distance $2\delta$ from each other. The union $A_\ga$ of these geodesics is $\ga$-invariant. It is also $3\delta$-quasiconvex. 
One finds a quasi-axis of $\ga$ by taking an arbitrary point $x\in A_\ga$ and taking the biinfinite concatenation of geodesic segments 
$$
... \ga^{-1}(x)x \star x \ga(x) \star \ga(x) \ga^2(x) ... 
$$ 
We will refer to $A_\ga$ as {\em the axis} of $\ga$ in $X$. \index{$A_\ga$, the axis of $\ga$} 

While a quasi-axis of $\gamma$ is not unique, by the Morse Lemma, all $(L,D)$-quasi-axes of $\gamma$ are within Hausdorff distance $H=H(L,D,\delta)$ from each other. In some cases, when any two distinct points in $\geo X$ are connected by a unique geodesic in $X$ (which is the case, for instance, when $X$ is a $CAT(-1)$ space), each hyperbolic element $\gamma\in \Gamma$ has a unique invariant geodesic. Suppose again that $X$ is a general $\delta$-hyperbolic geodesic metric space. There exists  a function $C(L,D,\delta)$ such that whenever $\ga$ is an isometry of $X$ with an $(L,D)$-quasi-axis $\alpha$ and 
for some $x\in \alpha$,
$$
d(x, \gamma(x))\geqslant C(L,D,\delta),
$$
then $\gamma$ is necessarily a loxodromic isometry of $X$. (Cf. Lemma 3.6 in Chapter 9 of \cite{CDP}.) We will discuss properties of isometries of $X$ (mostly, 
loxodromic isometries) in more detail in Section \ref{sec:projection} and Chapter \ref{sec:appendixB}. 

\subsection{Hyperbolic groups}\label{sec:Hyperbolic groups}

A finitely generated group $\Ga$ is said to be hyperbolic (or Gromov--hyper\-bolic) if $(Cay_\Ga, d_S)$ is a \index{hyperbolic group} 
Gromov-hyperbolic metric space. A group $\Ga$ is said to be $\delta$-hyperbolic if $(Cay_\Ga, d_S)$ is $\delta$-hyperbolic. 
A group is Gromov-hyperbolic if and only if it admits a geometric action on a proper Gromov-hyperbolic 
(geodesic) metric space. We will use the notation $\geo\Ga$ for the Gromov-boundary of $Cay_\Ga$. If $\Ga$ acts geometrically on a 
metric space $X$, then there is an equivariant homeomorphism between $\geo X$ and $\geo \Ga$. A hyperbolic group is said to be {\em elementary} if $\geo \Ga$ consists of at most two points (it turns out that $\geo \Ga$ cannot be a singleton). A group $\Ga$ is said to be \index{elementary hyperbolic group}
{\em nonelementary} otherwise (i.e. if $\geo \Ga$ consists of at least three points, in which case it has the cardinality of continuum). 

Algebraically, the elementary/nonelementary dichotomy can be described as follows:

\begin{itemize}

\item A hyperbolic group $\Ga$ is elementary if and only if it is virtually cyclic, i.e. contains a finite index cyclic subgroup (this subgroup can be finite). Equivalently, $\Ga$ is amenable. 

\item A hyperbolic group is nonelementary if and only if it contains a nonabelian free subgroup. 
\end{itemize}

\begin{example}
\begin{enumerate}
\item Virtually cyclic groups are hyperbolic. 

\item Finitely generated free groups are hyperbolic. 

\item Fundamental groups of compact surfaces of negative Euler characteristic are hyperbolic. 

\item More generally, fundamental groups of compact connected Riemannian manifolds of negative curvature are hyperbolic. 

\item Fundamental groups of compact connected 3-dimensional manifolds (possibly with boundary) are hyperbolic if and only if they do not contain $\bZ^2$. 

\item Free products of finitely many hyperbolic groups are hyperbolic. 

\item If a group $\Ga_1$  contains a finite index subgroup $\Ga_2$ then $\Ga_1$ is hyperbolic if and only if $\Ga_2$ is hyperbolic. 

\item If a group $\Ga$ contains a subgroup isomorphic to $\bZ^2$, then it is not hyperbolic. 

\item More generally, if $\Ga$ contains an amenable subgroup which is not virtually cyclic, then $\Ga$ is not hyperbolic. 
\end{enumerate}
\end{example}

\subsection{Convergence actions and conical points}
\label{sec:convdy}

\def\Homeo{\operatorname{Homeo}}

Let $Z$ be a compact metrizable space with at least three points.

A sequence $(h_n)$ in $\Homeo(Z)$ is called {\em contracting}
if there exist points $z_{\pm}\in Z$ such that 
\begin{equation}
\label{eq:contrtaucv}
h_n|_{Z-\{z_-\}}\to z_+
\end{equation} 
uniformly on compacts as $n\to+\infty$.
This condition is clearly symmetric, i.e.\ (\ref{eq:contrtaucv})
is equivalent to the dual condition 
that
\begin{equation}
\label{eq:contrtaucvdual}
h^{-1}_n|_{Z-\{z_+\}}\to z_-
\end{equation} 
uniformly on compacts as $n\to+\infty$.
The points $z_{\pm}$ are uniquely determined, since $|Z|\geq3$.

A sequence $(h_n)$ in $\Homeo(Z)$
is said to {\em converge} to a point $z\in Z$,
\begin{equation}
\label{eq:convseqhm}
h_n\to z
\end{equation} 
if every subsequence contains a contracting subsequence 
which, outside its exceptional point, converges to the constant map $\equiv z$.

One considers the following stronger form of convergence:

\begin{defi}[Conical convergence]
A converging sequence $h_n\to z$
converges {\em conically},
\begin{equation}
\label{eq:convseqhmcon}
h_n\stackrel{con}{\to} z
\end{equation} 
if for some relatively compact sequence $(\hat z_n)$ in $Z\setminus \{z\}$,
the sequence of pairs of distinct points
$h_n^{-1}(\hat z_n,z)$ is relatively compact in 
$$D(Z):= Z\times Z\setminus \diag(Z\times Z).$$
\end{defi}


\begin{lem}
If $h_n\stackrel{con}{\to} z$,
then the condition in the definition of conical convergence 
holds for all relatively compact sequences $(\hat z_n)$ in $Z\setminus \{z\}$. 
\end{lem}
\proof
Let $(\hat z_n)$ be a relatively compact sequence in $Z\setminus\{z\}$. 
For every contracting subsequence $(h_{n_k})$
there exists a point $\hat z\in Z$ such that 
$$h^{-1}_{n_k}|_{Z-\{z\}}\to\hat z$$ uniformly on compacts.
In particular,
$h^{-1}_{n_k}\hat z_{n_k}\to\hat z$ 
and the relative compactness of $(h_{n_k}^{-1}(\hat z_{n_k},z))$ in $D(Z)$ 
becomes equivalent to the condition 
that the sequence $(h^{-1}_{n_k}z)$ does not accumulate at $\hat z$.
The latter condition is independent of the sequence $(\hat z_n)$.
\qed

\medskip
The following criterion for being a conical limit point of a subsequence is immediate:\footnote{Here
it suffices that $|Z|\geq2$.}
\begin{lem}
\label{lem:recogncnlimcv}
A sequence $(h_n)$ in $\Homeo(Z)$
has a subsequence conically converging to $z\in Z$
iff there exists a subsequence $(h_{n_k})$ and a point $z_-\in Z$
such that the following conditions are satisfied:

(i) 
$h_{n_k}^{-1}|_{Z-\{z\}}\to z_-$ uniformly on compacts.

(ii)
$(h_{n_k}^{-1}z)$ converges to a point different from $z_-$.
\end{lem}

\medskip
We now consider group actions. 
A continuous action $\Ga\acts Z$ of a 
discrete group $\Ga$ is a {\em convergence action} \index{convergence action}\index{convergence property}
if every sequence $(\ga_n)$ of pairwise distinct elements in $\Ga$ 
contains a subsequence converging to a point, equivalently,
a contracting subsequence. An action is a convergence action if and only if the diagonal action of $\Ga$ on
$$
T(Z)=\{(z_1, z_2, z_3)\in Z^3: \#\{z_1, z_2, z_3\}= 3\}
$$
is properly discontinuous, see \cite{Bowditch}. 

The kernel of a convergence action is finite,
and we will identify $\Ga$ with its image in $\Homeo(Z)$
which we will call a {\em convergence group}. 

The {\em limit set} $\La\subset Z$ of a convergence group $\Ga<\Homeo(Z)$ is the subset of all points 
which occur as limits $z_+$ as in (\ref{eq:contrtaucv}),
equivalently,
as limits $z$ as in (\ref{eq:convseqhm})
for sequences $\ga_n\to\infty$ in $\Ga$.
The limit set is $\Ga$-invariant and compact.
A limit point $\la\in\La$ is {\em conical} if it occurs as the limit of a conically converging sequence.
A convergence group is said to have {\em conical limit set} if all limits points are conical,
and to be {\em non-elementary} if $|\La|\geq 3$.
Tukia \cite[Thm.\ 2S]{Tukia} has shown that in the non-ele\-men\-ta\-ry case the limit set is perfect 
and the $\Ga$-action on it is minimal. For nonelementary group actions, the condition that the limit set is conical is equivalent to  the condition that the action of $\Ga$ on $T(\La)$ is cocompact, see \cite{Bowditch}. 

If the limit set is conical, then $\Ga$ and its action on $\La$ are very special:

\begin{thm}[B.~Bowditch, {\cite{Bowditch}}]
\label{thm:bowdchr}
Suppose that $\Ga<\Homeo(Z)$ is a non-elementary convergence group with conical limit set $\La$.
Then $\Ga$ is hyperbolic and $\La\cong\geo\Ga$ equivariantly.
\end{thm}

The converse is easier:

\begin{thm}[{\cite{Gromov,Tukia,Freden}}]
\label{thm:hypgpbdac}
The natural action of a nonelementary hyperbolic group on its Gromov boundary 
is a minimal conical convergence action.
\end{thm} 

More generally, if $X$ is a proper Gromov hyperbo\-lic space and $\Ga< \Isom(X)$ is 
a properly discontinuous subgroup, then the action $\Ga\acts \geo X$  
the action on $\geo X$ of is a convergence action. See \cite[Theorem 3A]{Tukia} and \cite{DK18}. 

Let $G\times Z\to Z$ be a convergence action. An element $g\in G$ is said to be {\em loxodromic} (with respect to this action) if it has two distinct fixed points $z_\pm \in Z$, one of which is attractive and the other one repelling:
$$ 
\lim_{n\to\pm \infty} g^n(z)=z_\pm
$$
for all $z\in Z\setminus \{z_+, z_-\}$. If $Z$ is the ideal boundary of a proper Gromov-hyperbolic space $X$ and $G< \Isom(X)$, then an isometry $g\in G$ acts loxodromically on $X$ (in the sense of Definition \ref{def:axis}) if and only if it acts loxodromically on $Z$. The next theorem was proven by Tukia in \cite[Theorem 2R]{Tukia}: 

\begin{thm}\label{thm:tukia}
Suppose that $G\times Z\to Z$ is a nonelementary convergence group action. Then the set of pairs of fixed points of loxodromic elements of $G$ is dense in 
$\La\times \La$, where $\La\subset Z$ is the limit set of $G$. 
\end{thm}

 \section{Injective hulls of metric spaces}

Let $(X,d)$ be a nonempty metric space. Let $\Delta(X)$ denote the space of all functions $f: X\to \R$ such that
$$
f(x)+f(y)\geqslant d(x,y)
$$
for all $x, y\in X$. In particular, every $f\in \Delta(X)$ is a nonnegative function. We let $\Delta_1(X)$ denote the subspace of $1$-Lipschitz functions 
in $\Delta(X)$.

We equip $\Delta(X)$ with the natural partial order:
$$
f\leqslant g\iff f(x)\leqslant g(x), \forall x\in X. 
$$

Let $E(X)$ denote the set of all {\em minimal elements} of $\Delta(X)$ with respect to this order.  \index{injective hull $E(X)$}

\begin{lem}\label{lem:extreme}
$f\in \Delta(X)$ belongs to $E(X)$ if and only if 
\begin{equation}\label{eq:extreme}
f(x)=\sup_{y\in X} (d(x,y) -f(y))= \sup_{y\in X} (d_x(y) -f(y))
\end{equation}
for every $x\in X$. 
\end{lem} 
\begin{proof} 1. First, note that $f\in \Delta(X)$ if and only if 
for every $f\in \Delta$, $x\in X$, we have
$$
f(x)\geqslant \bar{f}(x):=\sup_{y\in X} (d(x,y) -f(y)). 
$$
 Suppose now that $f\in E(X)$ but  \eqref{eq:extreme} fails, i.e. there exists $x\in X$ for which $f(x)> \bar{f}(x)$. We then define a new function $g: X\to \R$ which equals 
 $f$ on $X\setminus \{x\}$ and $g(x)=\bar{f}(x)$. Clearly, $f> g$. Let's verify that $g\in \Delta(X)$. It suffices to  check each pair of points $x$, $z\ne x$, $z\in X$. We have
 $$
g(x)= \bar{f}(x)\geqslant d(x,z)- f(z)
 $$
(by the definition of $\bar{f}$). Adding $f(z)=g(z)$ to both sides, we get $g(x)+g(z)\geqslant d(x,z)$, as required. Thus, $f$ is not a minimal element of $\Delta(X)$, which is a contradiction. 

\medskip 
2. Conversely, suppose that $f\in \Delta(X)$ satisfies  \eqref{eq:extreme}. Suppose that $g: X\to \R$ is such that $g\leqslant f$ and $g(x)< f(x)$ for some $x\in X$. Thus,
$$
g(x)< d(x,y) -f(y)
$$ 
for some $y\in X$. We also have $f(y)\geqslant g(y)$, hence, $g(x)< d(x,y) - g(y)$, $g(x)+ g(y)< d(x,y)$, which means that $g$ cannot be in $\Delta(X)$. Hence, $f$ is minimal.   
\end{proof}

\begin{cor}\label{cor:C1-extreme}
1. For every $f\in E(X)$, $x\in X$ and $\eps>0$ 
there is a point $y\in X$  such that
$$
d(x,y) \leqslant f(x)+ f(y)\leqslant d(x,y)+\eps. 
$$

2. Each $f\in E(X)$ is $1$-Lipschitz, i.e. $E(X)\subset \Delta_1(X)$. 

3.   Every distance function $d_x$ is an element of $E(X)$.
\end{cor}
\begin{proof} 1. This part follows immediately from the implication $\Rightarrow$ of the lemma. 

2. In view of \eqref{eq:extreme} and the triangle inequality, for all $x, x'\in X$ we have 
$$
f(x)\leqslant \sup_y (d(x,x') + d(x',y) - f(y))= d(x,x') +  \sup_y (d(x',y) - f(y)). 
$$
Applying \eqref{eq:extreme} again, we get
$$
d(x,x') +  \sup_y (d(x',y) - f(y)) = d(x,x') + f(x').
$$
Hence, $f(x)-f(x')\leqslant d(x,x')$. The inequality $f(x')-f(x)\leqslant d(x,x')$ follows by switching the roles of $x, x'$. Hence, $|f(x)-f(x')|\leqslant d(x,x')$. 

\medskip 
3. Triangle inequality implies that $d_x\in \Delta(X)$ for every $x\in X$. Suppose that $d_x< \bar{d}_x$, i.e. there exists $p\in X$ such that 
$$
d_x(p)=d(x,p)< \inf_{y\in X} (d(p,y) - d_x(y))= \inf_{y\in X} (d(p,y) - d(y,x)).
$$
Thus, there is $y\in X$ such that $d(p,x)< d(p,y)- d(x,y)$, i.e. $d(p, x) + d(x,y)< d(p,y)$, contradicting the triangle inequality. 
\end{proof}

In particular, we obtain an embedding $\iota: X\to E(X), x\mapsto d_x$; by abuse of notation, we will retain the name $x$ for $\iota(x)$.

\begin{lem}
For all functions $f, g\in \Delta_1(X)$ we have $||f-g||<\infty$. 
\end{lem}
\begin{proof} Fix $y\in X$. Every  $f\in \Delta_1(X)$ satisfies  
$$
|f(x)-d(x,y)|\leqslant f(y), x\in X, 
$$
i.e. $||f-d_y||\leqslant f(y)$. Now, finiteness of $||f-g||$ follows from the triangle inequality. 
\end{proof}

We, thus, equip $E(X)$ with the metric given by the supremum-norm: 
$$
d_E(f, g)=||f-g||. 
$$

\begin{lem}\label{lem:norm-estimate} 
For every pair of functions $f, g\in E(X)$ and $\eps>0$ there are points $x, y\in X$ such that 
$$
||f-g||\leqslant d(x,y) - f(y)- g(x) +2\eps. 
$$
\end{lem}
\begin{proof} 
After swapping the roles of $f$ and $g$ if necessary, we can assume that  $||f-g||=\sup_X (f-g)$. 
Then there exists $x\in X$ such that $||f-g||\leqslant f(x) - g(x) +\eps$. By Corollary \ref{cor:C1-extreme}(1), 
there is $y\in X$ satisfying the inequality $f(x)\leqslant d(x,y) - f(y) +\eps$. 
By combining the two inequalities, we obtain
$$
||f-g||\leqslant d(x,y) - f(y) - g(x) +2\eps. 
$$
\end{proof}

\begin{lem}
For every $x\in X, f\in E(X)$, $d_E(x,f)=f(x)$. In particular, $\iota$ is an isometric embedding. 
\end{lem}
\begin{proof} 
Indeed,
$$
d_E(x,f)= \sup_{y\in Y} |d_x(y) -f(y)|.
$$
By \eqref{eq:extreme},
$$
f(x)=\sup_{y\in Y} (d_x(y)-f(y)).
$$
Since $f\in \Delta(X)$, $d_x(y)-f(y)=|d_x(y)-f(y)|$, implying $d(x,f)=f(x)$. \end{proof}

\begin{lem}\label{lem:group actions on hulls}
Suppose that $\Ga\acts X$ is an isometric group action on $X$. Then $\Ga$ acts on $E(X)$ isometrically via $f\mapsto f\circ \ga^{-1}$, $\ga\in \Ga, f\in E(X)$.  
\end{lem}
\begin{proof} Clearly, the action of $\Ga$ on the space of functions preserves $\Delta(X)$ and, thus, $E(X)$. We need to verify that the action on $E(X)$ is isometric. 
Take functions $f, g\in E(X)$ and $\ga\in \Ga$. Then for every $x\in X$ we have
$$
|f(\ga x)-g(\ga x)|\leqslant ||f-g||.
$$ 
Hence, $||f\circ \ga - g\circ \ga||\leqslant ||f-g||$. The opposite inequality follows by considering $\ga^{-1}$. 
\end{proof}

We refer to \cite{Lang} for a proof of the following theorem: 

\begin{thm}\label{thm:geodesic and complete}
For every metric space $(X,d)$, the space $(E(X),d_E)$ is geodesic and metrically complete. 
\end{thm}

In fact, much more is true and $E(X)$ is the {\em injective hull} of $(X,d)$.  This, among other things, implies that if $(A,d_A)$ is a metric space, $B\subset A$ is a subset 
and $h: (B,d_A)\to (E, d_E)$ is a 1-Lipschitz map, then $h$ extends to a 1-Lipschitz map $(A,d_A)\to (E,d_E)$. 

Furthermore, $E=E(X)$ satisfies a certain form of nonpositive curvature in the sense of Busemann. In order to state this property we will allow for constant speed reparameterization of geodesics.

\begin{thm}\label{thm:comb}
There is a family $\cC$ of constant speed geodesics $c: [0,1]\to E$ 
in $E$ such that any pair of points in $E$ is connected by a geodesic in $\cC$ and such that for every pair of geodesics $c_i\in \cC, i=1,2,$ for every $t\in [0,1]$ we have:
\begin{equation}\label{eq:B-convex}
d_E(c_1(t), c_2(t))\leqslant (1-t)d_E(c_1(0), c_2(0)) + t d_E(c_1(1), c_2(1)),\end{equation}
i.e. the distance function $d_E$ is convex along the paths $c\in \cC$. 
\end{thm}

\medskip 
 We refer the reader to \cite{Lang} for details.

\subsection{Injective hulls of Gromov-hyperbolic spaces}\label{sec:hyp-injective hull}

Urs Lang proved in \cite[Proposition 1.3]{Lang} that if  $(X,d)$ is a $\delta$-hyperbolic space (in Gromov's sense), then $E(X)$ is also  $\delta$-hyperbolic. 

\begin{defi}\index{almost geodesic metric space}\label{def:almost geodesic metric space} 
A metric space $(X,d)$ is {\em $C$-almost geodesic} if for every pair of points $x, y\in X$ and a number $D\in [0, d(x,y)]$ there exists $z\in X$ such that
$$
d(x,y)\leqslant d(x,z)+ d(z, y) \leqslant d(x,y) +C,
$$
and
$$
|d(x, z)- D|\leqslant C. 
$$
A metric space is {\em almost geodesic} if it is $C$-almost geodesic for some $C$. 
\end{defi}

Recall that a metric space $(X,d)$ is called {\em roughly geodesic} if there exists a constant $A$ such that any two points in $X$ can be 
connected by a $(1,A)$-quasigeodesic. 

\begin{example}
1. A geodesic metric space is $0$-almost geodesic. 

2. Every path-metric space is $\eps$-almost geodesic for every $\eps>0$. 

3. Every group $G$ with word metric $d_S$ is $\frac{1}{2}$-almost geodesic. 

4. Suppose that $(Y,d_Y)$ is a geodesic metric space, $X\subset Y$ is a subset $C$-Hausdorff close to $Y$. Then $(X,d_Y)$ is $2C$-almost geodesic. 

5. Every roughly geodesic metric space  is almost geodesic. 
\end{example} 

The following is a ``coarse version'' of Part 2 of Proposition 1.3 proven by Urs Lang in \cite{Lang} and our proof is a ``coarsification'' of Lang's argument: 

\begin{prop}\label{prop:almost}
Suppose that $(X,d)$ is $C$-almost geodesic $\delta$-hyperbolic in Gromov's sense. 
Then the Hausdorff distance between $E(X)$ and $\iota(X)$ is $\leqslant C':=2(C+\delta)$. 
\end{prop}
\begin{proof} Consider a function $f\in E(X)\setminus X$. According to Corollary \ref{cor:C1-extreme}, for $x\in X$ and every $\eps>0$ there is a point $x, y\in X$ such that 
\begin{equation}\label{eq:0}
d(x,y)\leqslant f(x)+f(y)\leqslant d(x, y) +\eps. 
\end{equation}
Since $f(x)>0$, without loss of generality we may assume that $0<\eps<f(x)$. Thus, 
$D:=f(x)-\eps\in [0, d(x,y)]$. Since $X$ is $C$-almost geodesic, there exists a point $v\in X$ such that 
\begin{equation}\label{eq:1}
d(x,y)\leqslant d(x,v)+ d(v, y) \leqslant d(x,y) +C,
\end{equation}
and
\begin{equation}\label{eq:2}
|d(x,v)-f(x)+\eps|=|d(x, v)- D|\leqslant C. 
\end{equation}
In particular, we have
\begin{equation}\label{eq:3}
d(x,v)\leqslant f(x)+C.  
\end{equation}
We also need to bound from above $d(v,y)-f(y)$.  Combining \eqref{eq:0}, \eqref{eq:1} and \eqref{eq:2} we obtain:
$$
f(y)\geqslant d(x,y) - f(x)\geqslant d(v,y) + d(x,v)  - C - f(x)\geqslant -2C -\eps + d(v,y) 
$$
and, therefore, 
\begin{equation}\label{eq:4}
d(v,y)\leqslant f(y) + 2C +\eps. 
\end{equation}

Now, we apply $\delta$-hyperbolicity of $E$ to the quadruple $f, d_v, d_x, d_y$:
$$
d(f, v) + d(x,y)\leqslant \max\{d(f, x)+ d(v,y), d(f,y) + d(v,x)\} +2\delta,
$$
see Lemma \ref{lem:alt-hyp}. Equivalently, 
\begin{equation}\label{eq:5}
f (v) + d(x, y) \leqslant \max \{f (x) + d(v, y), f (y) + d(v, x)\} + 2\delta.
\end{equation}
As we observed in inequalities \eqref{eq:3} and \eqref{eq:4},
$$
d(v,y)\leqslant f(y) + C, \quad d(v,y)\leqslant f(y) + 2C +\eps,
$$
thus,
$$
\max \{f (x) + d(v, y), f (y) + d(v, x)\} \leqslant f(x) + f(y) +2C+\eps. 
$$
Therefore,
$$
f(v) + d(x,y)\leqslant  f(x) + f(y) +2C+\eps +2\delta. 
$$
But \eqref{eq:0} implies that
$$
f(x)+f(y) + 2C+\eps+\delta\leqslant d(x,y)+2\eps + 2C+2\delta 
$$
and, therefore, in view of \eqref{eq:5}, 
$$
f(v)\leqslant 2(\eps + C)+2\delta. 
$$
In other words, $d(f, X)\leqslant 2(C+\delta)$. 
\end{proof}

In general, there is no reason to expect injective hulls of proper Gromov-hyperbolic spaces to be again proper. Nevertheless, one has:


\begin{prop}\label{prop:hull-visibility}
Suppose that $(X,d)$ satisfies the assumptions of Proposition \ref{prop:almost}. Then $E=E(X,d)$ 
satisfies the {\em visibility} property: Given any pair 
of distinct points $x, y\in \ol{E}$, there exists a geodesic $xy$ in $E$ connecting $x$ to $y$. 
\end{prop}
\begin{proof} If $x, y\in E$, the claim follows from the fact that $E$ is a geodesic space. Assume that $x\in E, y\in \geo E$. Consider a Gromov-sequence $(y_n)$ in $E$ representing $y$, $T_n:= d_E(x,y_n)$.  Take a sequence of unit speed geodesic segments $c_n: [0, T_n]\to E$ which are reparamaterizations of constant speed 
geodesics $[0,1]\to E$ which belong to the combing $\cC$ as in Theorem \ref{thm:comb}, connecting $x$ to $y_n$. Then, in view of hyperbolicity of $E$, there exists a constant $D=D(\delta)$ such that for every $R>0$ and $n$, the images $c_n([0,R])$ are $D$-Hausdorff-close to each other, as long as $T_n\geqslant R$. Then the inequality \eqref{eq:B-convex} implies that for every $t\in [0,\infty)$ the sequence $c_n(t)$ is Cauchy (where $T_n>t$). Metric completeness of $(E,d_E)$ then implies that the sequence of geodesic 
segments $c_n$ converges to a geodesic ray $c: [0,\infty)\to E$ connecting $x$ to $y$. The proof in the case $x, y\in \geo E$ is similar and we leave it to the reader.  
\end{proof}

Thus, $E(X)$ is semiproper provided that $X$ is proper, see Definition \ref{def:semiproper}.


\section{Strongly hyperbolic metric spaces} \label{sec:strong hyperbolicity}

\subsection{Definitions and main properties}\label{sec:Definitions and main properties}

 Strongly hyperbolic metrics (without this name) first appeared in the work of Mineyev \cite{Mineyev} on geodesic flows for hyperbolic groups. 
Our discussion of such  metrics mostly follows \cite{Nica-Spakula}.

\begin{definition}
A metric space $(X,d)$ is {\em $\epsilon$-strongly hyperbolic} \index{strongly hyperbolic space}
(for $\eps>0$) if for all points $x, y, z, p\in X$ and the Gromov product $(\cdot, \cdot)_p$ associated with the metric $d$ we have 
$$
\exp(-\epsilon (x, y)_p)\leqslant \exp(-\epsilon (x, z)_p) + \exp(-\epsilon (z, y)_p). 
$$
Equivalently, every quadruple $x_1, x_2, x_3, x_4$ in $X$ satisfies the inequality 
\begin{equation}\label{eq:Pto}
\exp\left(\frac{\eps}{2}(d(x_1,x_3) + d(x_2,x_4))\right) \leqslant \exp\left(\frac{\eps}{2}(d(x_1,x_2) + d(x_3,x_4))\right) + \exp\left(\frac{\eps}{2}(d(x_1,x_4) + d(x_2,x_3))\right). 
\end{equation}
\end{definition}

\begin{prop}\label{prop:hulls of strongly hyperbolic spaces}
If $(X,d)$ is $\eps$-strongly hyperbolic, then so is its injective hull $(E,d_E)$, $E=E(X)$. 
\end{prop}
\begin{proof} Our proof is almost a copy of the proof that taking the injective hull preserves $\delta$-hyperbolicity given in \cite[Proposition 1.3]{Lang}. 
Let $e, f, g, h$ be a quadruple of functions in $E(X)$. Then Lemma \ref{lem:norm-estimate} (applied separately to the pairs of functions $e, f$ and $g,h$) implies that 
for every $\eta>0$ there are points $x, y, z, w\in X$ such that
$$
d_E(e,f)\leqslant d(x,w) - e(w) -f(x) +2\eta, \quad d_E(g,h)\leqslant d(y,z) - g(y) -h(z) +2\eta. 
$$
Set $\Sigma:= e(w)+ f(x) + g(y) + h(z)$. By $\eps$-strong hyperbolicity of $(X,d)$ we have 
$$
\exp\left(\frac{\eps}{2}(d(x,w) + d(y,z))\right)\leqslant \exp\left(\frac{\eps}{2}(d(w,y) + d(x,z))\right) + \exp\left(\frac{\eps}{2}(d(x,z) + d(y,w))\right). $$
By combining the three inequalities we obtain
$$
\exp\left(\frac{\eps}{2}(d_E(e,f) + d_E(g,h))\right) \leqslant \exp(2\eps\eta) \exp\left(-\frac{\eps}{2} \Sigma\right) \exp\left(\frac{\eps}{2}(d(x,w) + d(y,z))\right) \le
$$
$$
\exp(2\eps\eta)\exp(-\frac{\eps}{2} \Sigma) \left( \exp(\frac{\eps}{2}(d(w,y) + d(x,z))) + \exp(\frac{\eps}{2}(d(x,z) + d(y,w))) \right). 
$$

Since $e, f, g, h$ are in $\Delta(X)$, we obtain the inequalities
$$
d(w, y) + d(x, z)- \Sigma \leqslant e(y) + f(z)- g(y)- h(z),$$ 
$$
d(x, z) + d(y, w)- \Sigma \leqslant - e(w) - f(x)+ g(x)+ h(w).  
$$
Next, the inequalities 
\begin{align*}
d(w, y) + d(x, z)- \Sigma \leqslant e(y) + f(z)- g(y)- h(z) \leqslant  \sup_{X} ( e(t)- g(t)) + \sup_{X} ( f(s)- h(s))= \\
d_E(e,g) + d_E(f,h) 
\end{align*}
and, similar, 
$$
d(x, z) + d(y, w)- \Sigma \leqslant - e(w) - f(x)+ g(x)+ h(w) \leqslant  d_E(f,g) + d_E(e,h)  
$$
imply
\begin{align*}
\exp(2\eps\eta)\exp(-\frac{\eps}{2} \Sigma) \left( \exp(\frac{\eps}{2}(d(w,y) + d(x,z))) + \exp(\frac{\eps}{2}(d(x,z) + d(y,w))) \right)\leqslant \\
\exp(2\eps\eta)  \left(  \exp(\frac{\eps}{2}( d_E(e,g) + d_E(f,h) ) + \exp(\frac{\eps}{2}(  d(x,z) + d(y,w))) \right). 
\end{align*}
Therefore,
$$
\exp(\frac{\eps}{2}(d_E(e,f) + d_E(g,h))) \leqslant \exp(2\eps\eta)  \left(  \exp(\frac{\eps}{2}( d_E(e,g) + d_E(f,h) ) + \exp(\frac{\eps}{2}(  d(x,z) + d(y,w))) \right). 
$$
Since $\eta>0$ was arbitrary, we obtain the desired strong hyperbolicity inequality for $(E,d_E)$. 
\end{proof}

The definition of strong hyperbolicity is pretty much designed so that if the Gromov product $(x,y)_p$ extends continuously 
to points $x, y\in \geo X$ (allowing for the infinite value in the case $x=y$), then the {\em naive formula} $\rho_{\eps,p}:=\exp(-\epsilon (x, y)_p)$ (from 
\S \ref{sec:Visual metrics}) results in a metric on $\geo X$. Here, as before, the convention is that $\exp(-\infty)=0$. 

The next definition is taken from \cite{Buyalo-Schroeder}: 

\begin{defi}
A Gromov-hyperbolic space $X$ is said to be {\em boundary continuous} if the 
Gromov product continuously extends to $\geo X$, $(x,y)_p$, where $p\in X, x, y\in \ol{X}$. \index{Boundary continuous hyperbolic space} 
\end{defi}

\begin{lem}
If $X$ is boundary continuous, then Busemann functions on $X$ can be defined by taking the limit in \eqref{eq:busemann}. 
Namely, for all $p, y\in X, \zeta\in \geo X$, the limit (in the definition of the Busemann function normalized to vanish at $p$) 
$$
b_\zeta(y):= \lim_{z\to \zeta} (d(y, z)- d(p,z))
$$
exists. Busemann functions are continuous as functions on $X\times \geo X$. 
\end{lem}
\begin{proof} Observe that $b_z(x,y):=d(x,z)-d(y,z), x, y, z\in X$, satisfies  
$$ 
b_z(x,y)=2(y,z)_x - d(x,y).
$$
Since $(y,z)_x$ is a continuous function of the triple $(x,y,z)\in \ol{X}^3$ (allowing for the infinite value if $y=z$), 
we obtain continuity of the extension $b_\zeta(x,y)$ ($\zeta\in \geo X$) of the function $b_z(x,y)$ and the fact that 
$$
b_z(x,y)= \lim_{i\to\infty} (d(y,z_i) - d(x,z_i))= \lim_{z_i\to z} (2(y,z_i)_x - d(x,y)), 
$$
where $(z_i)$ is a Gromov sequence representing $z$.  
\end{proof}

\begin{thm}
[see \cite{Nica-Spakula}] \label{thm:strongly}
Every $\epsilon$-strongly hyperbolic metric space $(X,d)$ satisfies:

\begin{enumerate}
\item $(X,d)$ is $\frac{1}{\eps}\log 2$-hyperbolic. 

\item $(X,d)$ is boundary continuous. 
 The boundary extension of the Gromov product still satisfies 
the $\eps$-strong hyperbolicity inequality. 

\item For all points $p\in X$, $x, y\in \geo X$, the formula 
$$
d_{\eps,p}(x,y):= \exp(-\eps(x,y)_p)
$$
defines a metric on $\geo X$. 

\item Busemann functions on $X$ can be defined by taking the limit in \eqref{eq:busemann}. 

\item The metric $d_\eps=d_{\eps,p}$ on $\geo X$ is {\em Ptolemaic}: 
$$
d_\eps(x_1, x_3) d_\eps(x_2, x_4)\leqslant d_\eps(x_1, x_2) d_\eps(x_3, x_4)+ d_\eps(x_1,x_4) d_\eps(x_2, x_3)
$$
for all quadruples $x_1, x_2, x_3, x_4\in \geo X$. 
\end{enumerate}
\end{thm}
\begin{proof} We will verify only part 2: Parts 3 and 5 then immediately follow 
(from continuity, the definition of strong hyperbolicity and equation \eqref{eq:Pto}).  Once Part 2 is proven, Part 4 follows from the previous lemma. 
We refer the reader to \cite{Nica-Spakula} for a proof of Part 1. 

2. Fix $p$ and let $(x_i), (y_i)$ be Gromov-sequences in $X$ representing points $x, y\in \geo X$. Thus,
$$
\lim_{i,j\to\infty} (x_i,x_j)_p=\infty, \quad \lim_{i,j\to\infty} (y_i,y_j)_p=\infty. 
$$
By applying the definition of strong hyperbolicity twice, we obtain
$$
|\exp(-\epsilon (x_i, y_i)_p) - \exp(-\epsilon (x_j, y_j)_p) | \leqslant \exp(-\epsilon (x_i, x_j)_p) + \exp(-\epsilon (y_i, y_j)_p). 
$$
Since the right-hand side of the inequality converges to $0$ as $i, j\to\infty$, we obtain that 
$$
\lim_{i,j\to\infty} |\exp(-\epsilon (x_i, y_i)_p) - \exp(-\epsilon (x_j, y_j)_p) | =0. 
$$
It follows that 
$$
\lim_{i,j\to\infty} |(x_i, y_i)_p -  (x_j, y_j)_p) | =0 
$$
(as long as $x\ne y$) and 
$$
\lim_{i\to\infty} (x_i, y_i)_p = \infty \iff \lim_{i\to\infty} (x_j, y_j)_p) | =\infty.   
$$ 
Thus, the extension of the Gromov product to points $x, y\in \geo X$ given by
$$
(x,y)_p:= \lim_{i\to\infty} (x_i, y_i)_p
$$
is well-defined. 
Continuity of the extension for fixed $p$ follows from the diagonal argument and 
the fact that a sequence $(z_i)$ in $X$ converges to $z\in \geo X$ if and only if $(z_i)$ is a Gromov sequence 
representing $z$. Continuity in $(x, y, p)\in \ol{X}\times \ol{X}\times X$  is left to the reader. Strong hyperbolicity for the extension is an immediate consequence 
of its continuity. 
\end{proof}


\begin{lem}\label{lem:linear}
Suppose that $c: \R\to X$ is a complete geodesic forward asymptotic to $\xi\in \geo \Ga$. Then $b_\xi$ is linear of slope $-1$ along $c$. 
\end{lem}
\begin{proof} We have $\lim_{t\to\infty}c(t)=\xi$. We will normalize $b_\xi$ to vanish at $p=c(0)$. Then 
$$
b_\xi(c(t))= \lim_{s\to \infty} (d(c(t), c(s))- d(p, c(s)))= -t. 
$$
\end{proof}

\medskip 
Earlier in this section we defined {\em Busemann almost cocycles}. In the context of boundary continuous hyperbolic spaces we obtain  the Busemann cocycle \index{Busemann cocycle $b_\zeta(x,y)$} 
$$
b_\zeta(x,y):= b_\zeta(x)- b_\zeta(y)= \lim_{z\to \zeta} ( d(x, z)- d(y,z) )  
$$
by taking the ordinary limit. In particular,
$$
b_\zeta(x,p)= b_\zeta(x). 
$$
The geometric meaning of $b_\zeta(x,y)$ is that it measures the signed distance between the horospheres $H_x, H_y$ in $X$ centered at $\zeta$ and passing through, respectively, the points $x$ and $y$. If $X$ is a complete negatively curved Riemannian manifold then this distance is constant: Assuming that $b_\zeta(x,y)>0$, 
for $p\in H_x$ one takes the unique geodesic ray $p\zeta$. This ray intersects $H_y$ at a unique point $q$ within distance $b_\zeta(x,y)$ from $x$, this point $q$ is the nearest point to $p$ on the horosphere $H_y$. 

By analogy with Lemma \ref{lem:almost cocycle} for boundary continuous hyperbolic spaces we obtain:
\begin{equation}\label{eq:Buscocycle}
\half b_\zeta(x,y) + \half b_\eta(x,y)= (\zeta,\eta)_x- (\zeta,\eta)_y\end{equation}
for all pairs of points $x, y\in X$ and distinct points $\zeta, \eta\in \geo X$. 

\medskip 
Unlike Busemann functions which are only invariant under isometries of $X$ up to a constant, Busemann cocycles are invariant:
$$
b_{g\zeta}(gx, gy)= b_\zeta(x, y)
$$
for all isometries $g$ of $X$. Since we are working with boundary continuous hyperbolic spaces, Busemann cocycles  also satisfy the {\em cocycle condition}
$$
b_\zeta(x,y) + b_\zeta(y,z)=b_\zeta(x,z),
$$
unlike the case of general $\delta$-hyperbolic spaces, where only the {\em almost cocycle} condition holds:
$$
b_\zeta(x,y) + b_\zeta(y,z)=b_\zeta(x,z) + O(1). 
$$
In the setting of isometric group actions $G\times X\to X$, we also define
$$
\sigma_{g}(\zeta):= b_{\zeta}(g^{-1}p,p)=-b_\zeta(p, g^{-1} p).
$$
The cocycle condition for $b_\zeta$ translates into a group cocycle condition for $G$:
\begin{equation}\label{eq:cocycle}
\sigma_{gh}(\zeta)= \sigma_g(h\zeta) + \sigma_h(\zeta)
\end{equation}
(in the case of boundary continuous hyperbolic spaces) and 
$$
\sigma_{gh}(\zeta)= \sigma_g(h\zeta) + \sigma_h(\zeta) + O(1)
$$
for general hyperbolic spaces. Accordingly, for 
$$
J_g(\zeta):=\exp(-\eps\sigma_g(\zeta))$$ 
we get the {\em Chain Rule}: 
\begin{equation}\label{eq:ChainRule}
J_{gh}(\zeta)= J_g(h\xi)\cdot J_h(\zeta) 
\end{equation}
for strongly hyperbolic spaces. (Technically speaking, $J_g$ depends also on $\eps$ and the choice of the base-point $p\in X$.) 
We will discuss the function $J_g$ in more detail below. 

\subsection{Conformality properties} \label{sec:Conformality properties}

\begin{lem}
[Conformality of visual metrics] \label{lem:conformality of visual metrics}
Suppose that $\rho=d_{\eps,p}, \rho'=d_{\eps,p'}$ are visual metrics on the ideal boundary of an 
$\eps$-strongly hyperbolic space $(X,d)$. Then for every point $z\in \geo X$ we have
$$
\left. \frac{d\rho'}{d\rho}\right\vert_z=\exp(-\eps b_z(p',p)). 
$$
\end{lem}
\begin{proof} In view of continuity of the Gromov-product on $\ol{X}$ (for strongly hyperbolic spaces) it 
suffices to prove that for all sequences $x_n, y_n$ in $X$ converging to $z$,
$$
\lim_{n\to\infty} \exp(\eps [(x_n,y_n)_{p'} - (x_n, y_n)_{p}])=  \exp(\eps b_z(p',p)), 
$$ 
equivalently,
$$
\lim_{n\to\infty} [(x_n,y_n)_{p'} - (x_n, y_n)_{p}]= b_z(p',p).
$$
By the definition of the Gromov product and Busemann functions, we get:
$$
(x_n,y_n)_{p'} - (x_n, y_n)_{p}= \frac{1}{2}\left(d(x_n,p)- d(x_n,p')+ d(y_n,p)- d(y_n,p')\right), 
$$
$$
\lim_{n\to\infty} \frac{1}{2}\left(d(x_n,p)- d(x_n,p')+ d(y_n,p)- d(y_n,p')\right)= \frac{1}{2}\left(b_z(p',p) + b_z(p',p)\right)= b_z(p',p),
$$
as required. 
\end{proof}

In the next corollary, $b_z$ is the Busemann function normalized to vanish at the base-point $p\in X$. 

\begin{cor}
[Isometries of strongly hyperbolic spaces are conformal at infinity] 
\label{cor:conformality of visual metrics}
If $X$ is $\eps$-strongly hyperbolic, then for every isometry $\ga$ of $X$ and every point $z\in \geo X$,  
$$
J_\ga(z)=\lim_{\xi, \eta\to z} \frac{d_{\eps,p}(\ga(\xi), \ga(\eta))}{d_{\eps,p}(\xi,\eta)}=  
\exp(-\eps b_z(\ga^{-1}(p)))
$$
i.e. the induced map $\ga: (\geo X, d_{\eps,p})\to (\geo X, d_{\eps,p})$ is conformal. 
\end{cor}
\begin{proof}  Set $p':= \ga^{-1} p$ and $\rho:= d_{\eps,p}, \rho':= d_{\eps,p'}$. Since $\ga: X\to X$ is an isometry, we have 
$$
\rho'= \ga^*\rho, \quad \rho'(\xi,\eta)= \rho(\ga x, \ga y).  
$$
By the definition of $J_\ga$ and Lemma \ref{lem:conformality of visual metrics}, 
$$
J_\ga(z)=\left. \frac{d\rho'}{d\rho}\right\vert_z= \exp(-\eps b_z(p',p))=  \exp(-\eps b_z(\ga^{-1}(p))). 
$$
\end{proof}

\begin{lem}\label{lem:Lip Jacobian}
For every isometry $\ga$ of an $\eps$-strongly hyperbolic space $X$, the function $\log J_\ga$ is $L$-Lipschitz on $(\geo X, d_{\eps,p})$ with $L=\eps^{-1}e^{\eps D}, D=d(p,\ga p)$. 
\end{lem}
\begin{proof} Set $x=\ga^{-1}(p)$. Strong hyperbolicity gives us: 
$$
|J^{-\eps}_\ga(\xi)-J^{-\eps}_\ga(\eta)|=|e^{\eps b_\xi(x)} - e^{\eps b_\eta(x)}|= |e^{-\eps(x,\xi)_p} - e^{-\eps(x,\eta)_p}|\leqslant e^{-\eps (\xi,\eta)_p}=d_\eps(\xi,\eta). 
$$
Hence, $J^{-\eps}_\ga(\xi)= e^{\eps b_\xi(x)}$ is $1$-Lipschitz on $(\geo X,d_\eps)$. From this, we get the claimed Lipschitz property of $J_\ga$. 
\end{proof}

\begin{example}\label{ex:free-jacobian}
1. Let $\ga$ be a translation by $D$ along a geodesic line $L\subset X$ containing $p$ and asymptotic to $z$. We assume that $\ga(p)$ lies in the ray $pz$. Then
$$
J_\ga(z)= e^{-b_z(\ga^{-1}(p))}= e^{-D}<1, 
$$
as expected, since $z$ is the attracting fixed point of $\ga$. 

2. Suppose that $X$ is the Cayley tree of a free group $\Ga$. As usual, we take $p=1_\Ga\in X$. 
As in Example \ref{ex:Busemann} we consider a geodesic ray 
$$
\rho=a_1a_2....a_k....
$$
in $X$ emanating from $1_\Ga$. Here and below, $a_i, c_j$ belong to the symmetrized free generating set of $\Ga$. 
Suppose that $\ga\in\Ga$ is such that
$$
\ga^{-1}=a_1a_2...a_mc_1c_2...c_n,
$$
Then 
$$
J_\ga(z)= e^{m-n}. 
$$
For instance, if $\ga=(a_1...a_m)^{-1}$ and $\rho$ is asymptotic to the repelling fixed point $z$ of $\ga$ then $J_\ga(z)= e^{m}$. 
\end{example}

\subsection{Existence of strongly hyperbolic metrics} 
\label{sec:green_metric}

\begin{thm}
[Theorem 5.1 in \cite{Nica-Spakula}] All $CAT(-1)$-spaces are $1$--strongly hyperbolic.
\end{thm}
 
Another source of strongly hyperbolic space is   the following theorem proven in \cite[Theorem 5]{Mineyev} and, later, 
with a different proof in \cite[Theorem 6.1]{Nica-Spakula} (the latter used Green metrics on hyperbolic groups): 

\begin{thm}\label{thm:strongly hyperbolic} 
For every hyperbolic group $\Ga$ there exists a $\Ga$-invariant metric $d$ on $\Ga$ such that the following holds:

1. $(\Ga, d)$ is roughly geodesic.

2. The identity map from $(\Ga, d)$ to $(\Ga, d_S)$ (where $d_S$ is the word-metric on $\Ga$ for a finite generating set $S$) is a quasi-isometry. 

3. $(\Ga, d)$ is $\epsilon$-strongly hyperbolic for some $\eps>0$. 
\end{thm}

\medskip
Below is a brief description of one construction of a strongly hyperbolic {\em Green metric} on a hyperbolic group $\Ga$, following 
\cite{Nica-Spakula}. 

Let $m$ be a  probability measure on $\Ga$. We will assume that $m$ is {\em symmetric}, i.e. 
$m(\ga)=m(\ga^{-1})$ for all $\ga\in \Ga$ and that the support set of $m$ equals a finite (necessarily symmetric) generating set $S$ of $\Ga$. 
Consider the random walk on $\Ga$ defined by the transition probabilities 
$$
p(\ga_1,\ga_2)= m(\ga_1^{-1}\ga_2). 
$$
Let $F(\ga_1,\ga_2)$ denote the probability that the random walk starting at $\ga_1$ hits $\ga_2$ at some time. 

\begin{defi}\label{def:G-metric}\index{Green metric}
The 
{\em Green metric} $d_G(\ga_1, \ga_2)$ on $\Ga$ is defined by
$$
d_G(\ga_1, \ga_2)= -\log F(\ga_1, \ga_2). 
$$
\end{defi}

Equivalently (and this justifies the name of the metric), if $G(\ga_1, \ga_2)$ denotes the Green function on $\Ga$ associated with the measure $m$, then  
$$
d_G(\ga_1, \ga_2)= \log G(1_\Ga, 1_\Ga) - \log G(\ga_1, \ga_2). 
$$
Green metrics were introduced by Blach\`ere  and Brofferio in \cite{MR2278460}, where they proved that for non-amenable groups 
(such as nonelementary hyperbolic groups), the identity map $(\Ga, d_G)\to (\Ga, d_S)$ is a quasiisometry, see also \cite[Corollary 1.2]{BHM} where it was proven that $d_G$ is roughly geodesic. It was further proven by Nica and \v{S}pakula in \cite{Nica-Spakula} that for hyperbolic groups this metric satisfies the properties stated in Theorem \ref{thm:strongly hyperbolic}. 

\medskip
The space $X=(\Ga, d)$ in Theorem \ref{thm:strongly hyperbolic} need not be a geodesic metric space. However, we can {\em thicken} $X$ to get a 
geodesic metric space with a $\Ga$-action, using the injective hull construction described in Section \ref{sec:hyp-injective hull}: 

\begin{thm}\label{thm:hull}
There exists a metric space $E(X)$, namely the {\em injective hull} of $(X,d)$, with the following properties:

1. $E(X)$ is a strongly hyperbolic complete geodesic metric space. 

2. There exists an isometric embedding $\iota: X\to E(X)$ whose image is at a finite Hausdorff distance from $X$. 

3. The action of $\Ga$ on $X$ (via left multiplication) extends to a metrically proper cobounded $\Ga$-action on $E(X)$. 
(Here we identify $X$ and $\iota(X)$.)  

4. 
$E(X)$ is semiproper. 
\end{thm}
\begin{proof} To prove the theorem we just need to collect various results established earlier. 

1. The fact that $E(X)$ is complete and geodesic was proven by Urs Lang, see Theorem \ref{thm:geodesic and complete}. Strong hyperbolicity of $E(X)$ is proven in Proposition \ref{prop:hulls of strongly hyperbolic spaces}. 

2. This part is proven in Proposition \ref{prop:almost} (one only needs hyperbolicity for this part, strong hyperbolicity is not needed). 

3. The existence of extension of the $\Ga$-action was proven in Lemma \ref{lem:group actions on hulls}. Metric properness of the action on $E(X)$ follows from the metric properness of the action on $X$. Coboundedness follows from Part 2 of the theorem.

4. Since $X$ is proper, its visual boundary $\geo X$ is compact. Since $X$ and $E(X)$ are Hausdorff-close to each  other, the inclusion 
$\iota: X\to E(X)$ induces a homeomorphism of their visual boundaries. Visibility property of $E(X)$ was proven in Proposition \ref{prop:hull-visibility}. 
\end{proof}

Note that, a priori, the space $E(X)$ is not  a proper metric space even if $X$ is. 



In Section \ref{sec:uniqueness}, we will be using this theorem to get $\Ga$-conformal densities on $\geo \Ga$ using strongly hyperbolic metrics on $\Ga$.

\subsection{Horospherical metrics on the punctured visual boundary} 

Lastly, we describe a version of the visual metric on the {\em punctured} ideal boundary of a strongly hyperbolic space. This construction is due to Hamenst\"adt, 
\cite{Hamenstadt}, in the setting of visual boundaries of Riemannian manifolds of negative curvature and Mineyev, \cite{Mineyev}, in the context of special metrics he defined on 
hyperbolic complexes. Our description in strongly hyperbolic setting follows the paper by 
Foertsch and Schroeder, \cite{Foertsch-Schroeder}. This metric is a generalization of the Euclidean metric on the ideal boundary (minus the point $\infty$) 
of the hyperbolic space $\bH^n$ in the upper half-space model. Pick a Busemann function $b$ on $X$ (associated with 
an ideal boundary point $\omega\in \geo X$), normalized to vanish at a point $p\in X$. Define the {\em Gromov product on $\ol{X}$ relative to $b$}:
$$
(x,y)_b:= (x,y)_p - (x,\omega)_p - (y,\omega)_p. 
$$   
Then for $x, y\in \geo X \setminus \{\omega\}$ set
\begin{equation}\label{eq:boundary metric}
d_{\eps,b}(x,y):= \exp(-\eps (x,y)_b). 
\end{equation}
The Ptolemaic inequality in Part 5 of Theorem \ref{thm:strongly} implies that $d_{\eps,b}$ is a metric on $\geo X \setminus \{\omega\}$, see \cite{Mineyev} or 
\cite{Foertsch-Schroeder}. 

\begin{lem}
[See Lemma 17 in \cite{Mineyev}.]\label{lem:horo-rescale}
Let $b, b'$ be two Busemann functions with respect to the same point $\omega\in \geo X$; hence, $b'-b\equiv C$, a constant. Then
$$
d_{\eps,b}(x,y)= e^{\eps C} d_{\eps,b'}(x,y).
$$
\end{lem}
\begin{proof} We let $p, p'$ be, respectively, the normalization points in $X$ used in the definitions of $b, b'$. 
First, consider the case when in the definition of  $d_{\eps,b}$ instead of $\omega\in \geo X$ we use $z\in X$. Then $(x,y)_b$ corresponds 
to
$$
(x,y)_p - (x,z)_p - (y,z)_p= (x,y)_z - d(z,p) 
$$
and, thus, $d_{\eps,b}(x,y)$ corresponds to
$$
\exp(\eps d(z,p)  -\eps  (x,y)_z).  
$$
Similarly,  $d_{\eps,b'}(x,y)$ corresponds to
$$
\exp(\eps d(z,p')  -\eps  (x,y)_z).  
$$
The term $e^{\eps C}$ in the lemma corresponds to $\exp(\eps d(p,z) - \eps d(p',z))$. Combining the terms, we obtain the equality  
$$
\exp(\eps d(z,p)  -\eps  (x,y)_z)= \exp(\eps d(p,z) - \eps d(p',z)) \cdot \exp(\eps d(z,p')  -\eps  (x,y)_z).
$$
Now, the claim of the lemma follows by continuity of the Gromov-product, taking limits as $z\to \omega$.  
\end{proof}

\begin{rem}
In particular, $d_{\eps,b}$ depends only on $b$ and not on the base-point $p$. The property in the lemma is a generalization of the fact that the Euclidean 
metrics on the ideal boundary of $\bH^n$ in the upper half-space model defined via identification with different concentric horospheres, differ by a constant conformal factor. 
Informally, one can think of such metrics as follows. Each horosphere in $\bH^n$ has induced Riemannian metric and, hence, intrinsic (Euclidean) distance function. Pushing 
such metrics to the ideal boundary (punctured at the center $\omega$ of the horosphere) defines a metric on $\geo \bH^n \setminus \{\omega\}$. While a direct generalization 
of this Riemannian construction is impossible for Gromov-hyperbolic spaces, the metrics $d_{\eps,b}$ serve as a replacement and enjoy properties similar to the above 
Euclidean metrics. 
\end{rem}

\begin{lem}\label{lem:horospherical-bilip1}
Consider the horospherical metric $d_{\eps,b}$ on $\geo X\setminus \{\omega\}$ determined by 
the Busemann function $b=b_\omega(p,\cdot)$. The metric $d_{\eps,b}$ is conformal to the visual metric 
$d_{\eps,p}$. Furthermore, the two metrics are $\max(e^2, R^{-2})$-bilipschitz to each other on the complement to the ball $B(\omega,R)\subset \geo X$ with respect to the visual metric  $d_{\eps,p}$. 
\end{lem}
\begin{proof} Fix a point $z\in \geo X\setminus \{\omega\}$. For $x, y\in \geo X\setminus \{\omega\}$ we have 
\begin{equation}\label{eq:two visuals}
\frac{d_{\eps,b}(x,y)}{d_{\eps,p}(x,y)}= \exp(\eps (x,\omega)_p+ \eps (y,\omega)_p).
\end{equation}
Taking the limit as $x, y\to z$, we obtain
$$
\lim_{x,y\to z} \frac{d_{\eps,p}(x,y)}{d_{\eps,b}(x,y)}= \exp(2\eps(z,\omega)_p). 
$$
This proves the conformality claim. Suppose now that 
$$
d_{\eps,p}(x,\omega)\geqslant R, \quad d_{\eps,p}(y,\omega)\geqslant R, 
$$
i.e. $\exp(-\eps (x,\omega)_p)\geqslant R$, $\exp(-\eps (y,\omega)_p)\geqslant R$.  Then \eqref{eq:two visuals} yields
$$
\frac{d_{\eps,b}(x,y)}{d_{\eps,p}(x,y)}\leqslant R^{-2}. 
$$
On the other hand, $(x,\omega)_p\geqslant 0, (y,\omega)_p$ for all points $x\in \geo X$. It follows that 
$\exp(\eps (x,\omega)_p+ \eps (y,\omega)_p)\geqslant e^2$. We, thus, obtain the claimed bilipschitz estimate. 
\end{proof}

\begin{lem}\label{lem:horospherical-bilip2}
Consider two horospherical metrics $\rho=d_{\eps,b}$, $\rho'=d_{\eps,b'}$, where $b=b_\omega(p,\cdot)$ and $b'=b_{\omega'}(p',\cdot)$. We assume that
 $(\omega,\zeta)_p=(\omega',\zeta)_{p'}$. 
Then 
$$
\left. \frac{d\rho}{d\rho'}\right\vert_\zeta= \left. Lip(\rho,\rho')\right\vert_\zeta= \exp(-\eps b_\zeta(p,p')).
$$
In particular, if $b_\zeta(p,p')=0$ (i.e. $p, p'$ belong to the same horosphere centered at $\zeta$) then 
$$
\left. \frac{d\rho}{d\rho'}\right\vert_\zeta= \left. Lip(\rho,\rho')\right\vert_\zeta= 1.
$$
\end{lem}
\begin{proof} In view of Lemma \ref{lem:horo-rescale}, it suffices to consider the case $b_\zeta(p,p')=0$. We have
$$
\rho(x,y)= \exp(-\eps (x,y)_p + \eps (x,\omega)_p + \eps(y,\omega)_p),  
$$
$$
\rho'(x,y)= \exp(-\eps (x,y)_{p'} + \eps (x,\omega')_{p'} + \eps(y,\omega')_{p'}).  
$$
Thus, the ratio $\frac{\rho(x,y)}{\rho'(x,y)}$ equals
$$
\exp(-\eps (x,y)_p + \eps (x,y)_{p'})\cdot \exp( \eps (x,\omega)_p + \eps(y,\omega)_p - \eps (x,\omega')_{p'} - \eps(y,\omega')_{p'}).  
$$
We then consider the limit of this expression as $x,y\to \zeta$. By Lemma 
\ref{lem:conformality of visual metrics}, the limit of the first term equals  $\exp(\eps b_\zeta(p,p'))=1$. The second term depends 
continuously on $x, y$ an has limit equal to 
$$
\exp( 2\eps (\zeta,\omega)_p - 2\eps (\zeta,\omega')_{p'})=1 
$$
by the assumption of the lemma. 
 \end{proof}

\begin{cor}\label{cor:Hausdorff-equal}
1. Metric spaces $(\geo X\setminus \{\xi\}, d_{\eps,b}))$ have the same Hausdorff dimension $\al$ as $(\geo X, d_{\eps,p})$.  In particular, all horospherical metrics have the same Hausdorff dimension. 

2. Assume that $b_\zeta(p,p')=0$ and $(\omega,\zeta)_p=(\omega',\zeta)_{p'}$. Let $\cH, \cH'$ denote the $\al$-dimensional Hausdorff measures defined by the metrics $\rho, \rho'$. Then 
$$
\left. \frac{d\cH}{d\cH'}\right\vert_\zeta=1. 
$$
\end{cor}
\begin{proof} 1. Lemma \ref{lem:horospherical-bilip1} implies that away from singular points $\omega$ of horospherical metrics, the 
visual metric is locally Lipschitz-equivalent to the horospherical metrics.   Part 1 follows.  

2. The second statement follows from Lemma \ref{lem:horospherical-bilip2}. 
\end{proof}


\begin{rem}
Recall that all Hausdorff measures $\cH^\al_x$ are in the same measure class, are positive and finite on $\geo \Ga$. Lemma \ref{lem:horospherical-bilip1} implies that  
the same is true for the $\al$-Hausdorff measures of the metrics $d_{\eps,b}$ on $\geo X\setminus \{\xi\}$, when $\xi$ is fixed, except that finiteness of the measure 
holds for compacts in $\geo X\setminus \{\xi\}$. 
\end{rem}

\section{Length functions on groups}\label{sec:length functions}


We let $\Gamma$ be a countable group equipped with discrete topology. 

\begin{defi} A function $\beta: \Gamma\to [0,\infty)$ is called  a {\em length function}, if it is \index{length function $\beta$} 
a proper function such that the following hold:

1. $\beta(\ga)=0\iff \ga =1\in \Ga$.

2.  $\beta(\gamma)=\beta(\gamma^{-1})$ for all $\gamma\in\Gamma$. 
 
3. $\beta(\gamma_1 \gamma_2)\leqslant \beta(\gamma_1) + \beta(\gamma_2)$  for all $\gamma_1, \gamma_2\in\Gamma$. 
\end{defi}

\begin{example}
Suppose that $\Ga\acts (X,d)$ is an isometric metrically proper action and $x_0\in X$ is a point with trivial $\Ga$-stabilizer. Then the function
$$
\beta(\ga)=d(x_0, \ga x_0)
$$
is a length function on $\Ga$. 
\end{example}

The next lemma shows that all length functions come from this example:

\begin{lem}\label{lem:dbeta}
Suppose that $\beta$ is a length function on $\Ga$. Then the function $d=d_\beta: \Ga\times \Ga\to \R_+$ given by 
 $$d(\gamma_1, \gamma_2):= \beta(\gamma_1^{-1}\gamma_2)$$ 
 is a $\Ga$-invariant proper metric on $\Ga$. (Here $\Ga$ acts on itself via the left multiplication.) 
 In particular,  $\beta(\ga)=d(1, \ga)$. 
\end{lem}
\begin{proof} Let us check the triangle inequality:
$$
d(\gamma_1, \ga_2) + d(\ga_2, \ga_3)= \beta(\gamma_1^{-1}\gamma_2) + \beta(\gamma_2^{-1}\gamma_3) \geqslant 
$$
$$
\beta(\gamma_1^{-1}\gamma_2\gamma_2^{-1}\gamma_3)= \beta(\gamma_1^{-1}\gamma_3)= d(\ga_1, \ga_3). 
$$
The other two axioms of a metric are clear as well. Furthermore, $\Ga$-invariance of $d$ is clear from the definition and properness of $d$ follows from properness of the function $\beta$. 
\end{proof} 

\begin{defi}\index{hyperbolic length function}
\label{def:hyperbolic length function}
A length function $\beta$ on a group $\Ga$ is said to be {\em hyperbolic}  if two conditions are met: 

1. $(\Ga, d_\beta)$ is a Gromov-hyperbolic metric space.   

2. The metric space $(\Ga, d_\beta)$ is almost geodesic. Equivalently, there exists a constant $C$ such that 
for every $x\in \Ga$ and a number $D\in [0, \beta(x)]$ there is $y\in \Ga$ such that
$$
\beta(x)\leqslant \beta(y)+ \beta(y^{-1}x) \leqslant \beta(x) +C,
$$
and
$$
|\beta(y)- D|\leqslant C. 
$$
\end{defi}


\begin{example}\label{ex:hyperbolic length}
1. Word-length function associated with finite generating set of a hyperbolic group is a hyperbolic length function. 

2. Suppose that $\Ga$ is a group acting geometrically on a geodesic metric space $(X,d)$. Then for every $x\in X$, 
the function $\beta(\ga)=d(x, \ga x)$ is a hyperbolic length function. 
\end{example}

Here is an example of a non-hyperbolic length function on a hyperbolic group:

\begin{example}
[Snowflaking] Let $\beta$ be the word-length function on an infinite hyperbolic group $\Ga$. Take a real number $\al\in (0,1)$. 
Then the function $\beta^\al: \ga\mapsto (\beta(\ga))^{\al}$ defines a $\Ga$-invariant Gromov-hyperbolic metric on $\Ga$, but 
is not a hyperbolic length function since it is not almost geodesic. 
\end{example} 


\begin{lem}\label{lem:regular}
Suppose that $\beta$ is a hyperbolic length function on $\Ga$. Then:

1. $\beta$ comes from  a cobounded action of $\Ga$ on a semiproper hyperbolic metric space as in Example \ref{ex:hyperbolic length}. 

2. In particular, the identity map $(\Ga, d_\beta)\to (\Ga, d_S)$ is a quasi-isometry, where $d_S$ is a word metric on $\Ga$, and  $(\Ga,d_S)$ is hyperbolic. 
\end{lem}
\begin{proof} 
To prove the first claim, recall that $(\Ga, d_\beta)$ is an almost geodesic hyperbolic metric space. The action of $\Ga$ on itself via left multiplication 
is isometric and transitive. Let $(X,d_X)$ denote the injective hull $E(\Ga, d_\beta)=(E,d_E)$ of the metric space $(\Ga, d_\beta)$. Then the inclusion map 
$(\Ga, d_\beta)\to (X,d_X)$ is an isometric embedding, $(X,d_X)$ is a Gromov-hyperbolic geodesic complete metric space and the $\Ga$-action on 
$X$ is cobounded  (since the Hausdorff distance between $\Ga$ and its injective hull is bounded, see Proposition \ref{prop:almost}). Semiproperness of $(E,d_E)$ is proven in Proposition \ref{prop:hull-visibility}. 
The second claim of the lemma is a consequence of the Milnor--Schwartz Lemma and quasi-isometry invariance of hyperbolicity for geodesic metric spaces. 
 \end{proof}

\medskip 
 
 We refer the reader to the recent paper by Oreg\'on-Reyes, \cite{OR}, for an in-depth discussion of the {\em space} of invariant metrics on hyperbolic groups coming from geometric actions on geodesic spaces, equivalently, hyperbolic length functions.

\section{Estimates on Gromov products and fixed points of loxodromic isometries}\label{sec:projection}

In this section we prove several technical results that will be used in Section \ref{sec:redressing}.

\medskip 
Given a metric space $(X,d)$, a subset $A\subset X$ and $\eps\geqslant 0$, 
define the {\em $\eps$-nearest point projection} $P^\eps_A: X\to 2^A$ (the power-set of $A$)  
\index{$P^\eps_A$, the $\eps$-nearest point projection}
with $P^\eps_A(x)$ equal to the subset of $A$ 
consisting of all $\bar{x}\in A$ such that for every $a\in A$
$$
d(x, a)\geqslant d(x, \bar{x}) - \eps. 
$$

From now on, we assume that $(X,d)$ is a $\delta$-hyperbolic geodesic metric space. 

\begin{lem}\label{lem:1}
If $l$ is a (finite or infinite) geodesic in $X$, $p\in X$, $\bar{p}\in P^\eps_l(p)$, then for every 
$y\in l$ we have
\begin{equation}\label{eq:projection0}
 d(\bar{p}, py)\leqslant \eps':=\eps+ 3\delta, 
\end{equation}
\begin{equation}\label{eq:projection}
d(p, y) \leqslant d(p, \bar{p}) + d(\bar{p}, y)\leqslant d(p, y) +2\eps', 
\end{equation}
\end{lem} 
\begin{proof} 1. This is proven for $\eps=0$ in \cite[Lemma 1.105]{Kapovich-Sardar}. The proof for general $\eps\geqslant 0$ is similar: 
Take $y\in l$. Since the triangle $\Delta p \bar{p} y$ is $\delta$-slim, there exists $m\in py$ and points $m_1\in 
p \bar{p}$, $m_2\in \bar{p}y$ such that
$$
\max( d(m, m_1), d(m, m_2))\leqslant \delta. 
$$
Thus,
$$
d(p, m_2)\leqslant d(p, m_1)+ 2\delta.$$
Since $\bar{p}\in P^\eps_l(p)$, 
$$
d(p, m_1) + d(m_1, \bar{p}) -\eps= d(p, \bar{p})-\eps \leqslant d(p, m_2).$$ 
Combining the two inequalities, we get
$$
d(p, m_1) + d(m_1, \bar{p}) -\eps= d(p, \bar{p})-\eps \leqslant d(p, m_1)+ 2\delta.$$ 
Thus, $d(m_1, \bar{p})\leqslant \eps+2\delta$ and, therefore, $d(m, \bar{p})\leqslant \eps+3\delta$, as claimed. 

2. The second part is a consequence of the first part and the triangle inequality
\end{proof}

We continue with the notation from the previous lemma: 

\begin{lem}\label{lem:2}
Suppose that $l$ is a complete geodesic connecting points $\xi^\pm\in \geo X$. 
Take $p, q\in X$, $\bar{p}\in P^\eps_l(p)$, $\bar{q}\in P^\eps_l(q)$, and assume that 
$\bar{q}$ belongs to the subray $\bar{p}\xi^+\subset l$. Then for every 
$y\in \bar{p}\xi^+$ we have
$$
(q,y)_p\geqslant d(p,q) - d(q, \bar{q}) -\eps'. 
$$
\end{lem}
\begin{proof} Recall that 
$$
(q,y)_p= \half(d(p,q) + d(p,y)  - d(q,y)). 
$$
Estimating the difference $d(p,y)  - d(q,y)$ using Lemma \ref{lem:1}, we obtain:
\begin{align*}
d(p,y)  - d(q,y) \geqslant (d(p, \bar{p}) + d(\bar{p}y) -2\eps') - (d(q,\bar{q}) + d(\bar{q}, y))\ge\\
d(p, \bar{p}) + d(\bar{p}, \bar{q}) - d(q, \bar{q}) -2\eps'. 
\end{align*}
Combining this inequality with $d(p, \bar{p}) + d(\bar{p}, \bar{q}) + d(q, \bar{q})\geqslant d(p,q)$, we get:
$$
(q,y)_p\geqslant \half \left( 2 d(p,q) -2 d(q, \bar{q}) - 2\eps'\right) = d(p,q) - d(q, \bar{q}) - \eps'. 
$$
\end{proof}

Taking the limit $y\to \xi^+$ and taking into account the fact that
$$
(q,\xi^+)_p\geqslant \lim\inf_{y\to \xi^+} (q,y)_p -2\delta',
$$
as a corollary of the lemma we obtain:

\begin{cor}\label{cor:1}
$$
(q,\xi^+)_p\geqslant d(p,q) - d(q, \bar{q}) - \eps' - 2\delta'. 
$$
\end{cor}

We now consider the setting in the previous lemma and Corollary \ref{cor:1} when $\bar{q}\notin \bar{p}\xi^+$.
 By repeating the proof of Lemma \ref{lem:2} and its corollary we then get:

\begin{lem}\label{lem:3}
Suppose that $l$ is a complete geodesic in $X$ connecting points $\xi^\pm\in \geo X$ and $\bar{q}$ does not belong to the subray $\bar{p}\xi^+\subset l$. Then for every 
$y\in \bar{p}\xi^+$ we have:

1. $(q,y)_p\geqslant d(p,q) - \left(d(q, \bar{q}) + d(\bar{p}, \bar{q}) + \eps'\right)$. 

2. $(q,\xi^+)_p\geqslant d(p,q) - \left(d(q, \bar{q}) + d(\bar{p}, \bar{q}) + \eps' + 2\delta'\right)$. 
\end{lem}

By combining Lemma \ref{lem:2}, Corollary \ref{cor:1} and Lemma \ref{lem:3} we get: 

\begin{prop}
Assume that $d(q, \bar{q})\leqslant R$ and either $\bar{q}\in \bar{p}\xi^+$ or $d(\bar{p}, \bar{q})\leqslant D_0$. Then for every $\zeta\in \geo X$ we have
$$
(q,\xi^+)_p \geqslant (q,\zeta)_p - (R+D_0+\eps' + 2\delta').
$$
\end{prop}

As an immediate consequence of the proposition and the inequality \eqref{eq:Gh1} applied to the points $p=x_0, z=\zeta, x=\xi^+, y=q$, 
we obtain: 

\begin{cor}\label{cor: product inequality}
Under the assumptions of the proposition, 
$$
(\zeta,\xi^+)_p\geqslant (\zeta,q)_p - (R+D_0+\eps' + 3\delta'). 
$$
\end{cor}

Let $\gamma$ be a loxodromic isometry of a $\delta$-hyperbolic geodesic metric space $X$ with the attractive/repelling fixed points $\xi^\pm\in \geo X$ and the axis $A=A_\ga$. We let $P_A: X\to A$ denote a $\ga$-equivariant nearest-point projection (it always can be found using the Axiom of Choice). For a 
base-point $x_0\in X$ let $x:= \ga(x_0)$, $\bar{x}_0:=P_A(x_0)$. Then
$P_A(x)=\ga \bar{x}_0$. 

Let $l\subset A$ denote a complete geodesic in $X$ passing through $\bar{x}_0$ and $l_\pm$ the subrays $\bar{x}_0\xi^\pm$ in $l$. We let $c: \R \to l$ denote the arc-length parameterization of $l$ such that
 $c(0)= \bar{x}_0$ and 
$$
\lim_{t\to \pm \infty} c(t)= \xi^\pm. 
$$
We let $\bar{x}\in l$ denote $P_c(P_A(x))$. Here and below, $P_c$ is the restriction to $A$ of the nearest-point projection to $l$. Recall that the Hausdorff distance between any two complete geodesics in $A$ is $\leqslant 2\delta$. 
Then $d(P_A(x), \bar{x})\leqslant 2\delta$, which implies that 
$\bar{x}\in P^\eps_l(x)$ for $\eps=2\delta$, where $P^\eps_l: X\to 2^l$ is the $\eps$-projection.  Set $D:= d(\bar{x}_0, \bar{x})$.   


 \begin{lem}\label{lem:right move}
Under the above assumptions, $\bar{x}\in l_+$, provided that $D\geqslant D_0=14\delta$.  
\end{lem}
\begin{proof} 

Suppose that $\bar{x}\in l_-$, i.e. $\bar{x}=c(-D)$. Define $U:= P_c^{-1}(l_-)$. 
Our goal is to prove that $\ga(U)\subset U$. Recall that $d(\gamma(x_0), \bar{x})\leqslant 2\delta$. 

Take $s<0$ and consider the quadruple $c(s), \ga(c(s)), \ga(\bar{x}_0), \bar{x}_0$. Noting that the complete geodesics $l$ and $\gamma(l)$ are within distance $2\delta$ and applying triangle inequality we obtain: 
$$
|d(\bar{x}_0, \ga(\bar{x}_0))  - d(c(s), \ga(c(s)))|\leqslant 2\eta=4\delta. 
$$
Therefore, 
\begin{equation}\label{eq:pro0}
P_c\left( \gamma c(s)\right)= c(s'),  |s'- (s-D)|\leqslant 8\delta
\end{equation}

We also have (by applying the triangle inequality and taking into account the fact that 
$d(\gamma c(s), c(s'))\leqslant 2\delta$) 
\begin{equation}\label{eq:pro}
P_c( B( \gamma c(s), 2\delta)\cap A)\subset B( c(s'), 6\delta).  
\end{equation}

Take $y\in U$, i.e. $P_c(y)=c(s)$ for some $s<0$. Then $d(y, c(s))\leqslant 2\delta$, 
$d(\ga y, \ga c(s))\leqslant 2\delta$, hence, by \eqref{eq:pro}, 
$$
d(P_c(\ga y), c(s'))\leqslant 6\delta,
$$
where $c(s')=P_c( \ga c(s))$. Thus, $P_c(\ga y)=c(t)$ where $|t-s'|\leqslant 6\delta$. Applying 
the inequality \eqref{eq:pro0}, we get:
$$
|t- (s-D)|\leqslant 8\delta + 6\delta=14\delta. 
$$
Since $s<0$, it follows that 
$$
t< - D + 14\delta\leqslant 0 
$$
by the assumption $D\geqslant D_0=14\delta$ of the lemma. 
 
Thus, we proved that  $\ga(U)\subset U$ and, therefore, for every $n\geqslant 0$, 
$\ga^n(U)\subset U$. It follows that $\lim_{n\to\infty} \ga^n(\bar{x}_0)\ne \xi^+$,  
contradicting the assumption that $\xi^+$ is the attractive fixed point of $\ga$.  This contradiction proves that $\bar{x}\in l_+$, as claimed. 
\end{proof} 

We are now ready for the main theorem of this section. Its generalization (which does not require an upper bound on $d(x_0, A_\ga)$)  
will be proven in Appendix 2.

\begin{thm}\label{thm:fixed-point-estimate}
Assume that $\ga$ is a loxodromic isometry of a $\delta$-hyperbolic geodesic metric space $X$, 
with the attractive  fixed point $\xi^+\in \geo X$ and the axis $A_\delta$. Suppose also 
that $x_0\in X$ is the base-point such that $d(x_0, A_\ga)\leqslant R$. Set
\begin{align*}
D_0:=14\delta, \eps=2\delta, \eps':=\eps+ 3\delta, \eps'=5\delta, \\ 
R'=  R+D_0+\eps' + 2\delta' = 
R+ 25\delta. 
\end{align*}
Then for every $\zeta\in \geo X$, the point $x=\ga(x_0)$ satisfies 
$$
(x,\xi^+)_{x_0} \geqslant (x,\zeta)_{x_0} - R'.
$$
\end{thm}
\begin{proof} The theorem is an immediate consequence of the fact (noted earlier) that $\bar{x}\in P_l^\eps(x)$, Corollary \ref{cor: product inequality} and 
Lemma \ref{lem:right move}: Lemma \ref{lem:right move} proves that the points $p=x_0, \bar{p}=\bar{x}_0, q=x, \bar{q}=\bar{x}$ satisfy the assumptions of  
Corollary \ref{cor: product inequality}. 
\end{proof}

Informally speaking, if $\ga(x_0)$ is $t$-close to some $\zeta\in \geo X$, then $\xi^+$ is $O(t)$-close to 
$\zeta$. Restating this theorem in terms of neighborhoods $U(\zeta, \cdot)$ in $\ol{X}$ (defined relative to the base-point $x_0\in X$), 
we get:

\begin{cor}\label{cor:fixed-point-estimate}
Under the assumptions of the theorem, if $\ga(x_0)\in U(\zeta, t)$, then $\xi^+\in U(\zeta, e^{R'}t)$. 
\end{cor}

\section{Quasiconvex-cocompact and convex-cobounded groups} \label{sec:quasiconvex}


We refer the reader to \cite{Coo93} for a  treatment of quasiconvex-cocompact groups of isometries of Gromov-hyperbolic spaces and to {\cite[Definition 12.2.5]{DSU}} for a more general treatment of convex-cobounded groups of isometries { when the space is not assumed to be proper}.  Note that all this discussion in this section goes through in the context of $(1,A)$-roughly geodesic Gromov-hyperbolic spaces, with geodesics replaced by $(1,A)$-rough geodesics and geodesic rays replaced by $(1,A)$-rough geodesic rays. 
More precisely, following \S~\ref{sec:green_metric}, any hyperbolic group can be endowed with a $\Gamma$-invariant metric $d$ so that $(\Gamma,d)$ becomes roughly geodesic and $\epsilon$-strongly hyperbolic. In this specific situation, the notions of quasiconvex-cocompact groups will also make sense.

\subsection{Limit points}

Let $X$ be a 
geodesic $\delta$-hyperbolic space and let $\Ga< \Isom(X)$ be a discrete subgroup acting metrically properly 
on $X$. The {\em limit set} $\Lambda(\Gamma)$ \index{$\La(\Ga)$, the limit set} 
of the action of $\Ga$  is the accumulation set in $\geo X$ of an orbit $\Gamma x\subset X$ (the limit set does not depend on the choice of $x$). Thus, the limit set is a closed $\Ga$-invariant subset of $\geo X$. 
The \index{$\Omega(\Ga)$, discontinuity domain} 
{\em discontinuity domain} $\Omega(\Ga)$ of $\Ga$ is the complement $\geo X \setminus \La(\Ga)$. 
The action of $\Ga$ on $X\cup \Omega(\Ga)$ is properly discontinuous. 
Note that if the action of $\Ga$ on $X$ is cobounded, 
then $\La(\Ga)=\geo X$. 
Furthermore, if $(\ga_n)$ is a diverging sequence in $\Ga$ such that for some $x\in X\cup \Omega(\Ga)$, 
$\lim_{n\to\infty} \ga_n(x)=\la\in \La(\Ga)$, then the same holds for all points $y\in X\cup \Omega(\Ga)$:
$$
\lim_{n\to\infty} \ga_n(y)=\la. 
$$
A discrete subgroup acting metrically properly on $X$ is said to be elementary if its limit set is finite, equivalently, 
if it consists of at most two points. 

\medskip 
Below is a brief discussion of taxonomy of limit points. Suppose that $X$ satisfies the visibility property. 
A limit point $\xi\in \La(\Ga)$ is said to be {\em uniformly conical} if for one (equivalently, every) $x\in X$ there 
exists $D<\infty$ such that for 
$$
x\xi \subset N_D(\Ga x). 
$$
A limit point $\xi$ is said to be {\em conical} if, there exists $D<\infty$ and 
a sequence $x_n\in \Ga x$ and  such that $x_n\in N_D(x\xi)$ and converges to $\xi$. It turns out (although we will not need this) that this notion of conical limit points is equivalent to the one for convergence group actions, as defined 
in \S \ref{sec:convdy}.  A limit point 
$\xi\in \La(\Ga)$ is called a {\em horospherical limit point} if every  horoball $B(\xi, t)$ centered at $\xi$ has nonempty intersection with every $\Ga$-orbit $\Ga x\subset X$.  In other words, there is a sequence $x_n\in \Ga x$ (necessarily converging to $\xi$) such that $b_\xi(x_n)\to-\infty$.   \index{horospherical limit point} One defines the {\em conical limit set} $\La^c(\Ga)$ and {\em horospherical limit set} $\La^h(\Ga)$ as subsets consisting of, respectively, conical and horospherical limit points of $\Ga$. It is easy to see that $\La^c(\Ga)\subset \La^h(\Ga)$. In general, this inclusion is proper. 

\subsection{Equivalent definitions}

Quasiconvex-cocompact groups of isometries $\Ga< \Isom(X)$ behave similarly to groups acting geometrically on $X$, but they, typically, are not cocompact. 

\begin{defi}\index{quasiconvex-cocompact subgroup}
1. Let $X$ be a proper geodesic Gromov-hyperbolic metric space. A properly discontinuous subgroup $\Ga< \Isom(X)$ is said to be {\em quasiconvex-cocompact} if its action on the closure of the quasiconvex hull of one (equivalently, every) 
$\Ga$-orbit in $X$ is cocompact. 

2. Suppose that $X$ is merely a geodesic Gromov-hyperbolic metric space and 
$\Ga< \Isom(X)$ is a subgroup acting metrically properly on $X$. Then $\Ga$ is said to be {\em convex-cobounded} if 
the $\Ga$-action on the quasiconvex hull of one (equivalently, every) point in $X$ is cobounded. 

3. A subgroup $\Gamma_0< \Ga$ of a hyperbolic group $\Ga$ is said to be {\em quasiconvex} if its action on $Cay_\Ga$ is quasiconvex-cocompact.  
\index{quasiconvex subgroup} \index{convex-cobounded group} \index{quasiconvex-cocompact group}
\end{defi}


\begin{rem}
1. Clearly, if $X$ is proper then quasiconvex-cocompactness is equivalent to convex-coboundedness for subgroups 
$\Ga< \Isom(X)$. 

2. Suppose that $x, x'\in X$ are points within distance $D$ from each other, $\ga\in \Ga$, $y=\ga x, y'=\ga x'$. 
Then slimness of geodesic quadrilaterals in $X$ implies that each geodesic $x'y'$ is contained in the $D+2\delta$-neighborhood of $xy$. Thus, the hulls $\hull \Ga x, \hull \Ga x'$ are within Hausdorff distance $D+2\delta$ from each other.

3. It follows that  a metrically proper discrete subgroup $\Ga<\Isom(X)$ is  
convex-cobounded if and only if one (equivalently, every) $\Ga$-orbit in $X$ is a quasiconvex subset of $X$.
\end{rem}


\begin{lem}
Suppose that $X$ is semiproper and $\Ga$ is nonelementary. Then $\Ga$ is convex-cobounded if and only if its action on the quasiconvex hull of its limit set is cobounded. 
\end{lem}
\begin{proof} 1. Suppose that $\Ga$ is convex-cobounded: The $\Ga$-action on $\hull \Ga x$ is cobounded, i.e. there 
is $R<\infty$ such that for any pair of points $y=\ga x, y'=\ga'x \in \Ga x$, the segment $yy'$ is contained in the $R$-neighborhood $N_R(\Ga x)$ of $\Ga x$. Consider two distinct points $\xi, \xi'\in \La=\La(\Ga)$ and a geodesic $\xi \xi'$ in $X$ connecting them. We then have two sequences $y_n=\ga_n x, y'_n=\ga'_n x$ in $\Ga x$ converging, respectively, to $\xi, \xi'$. Since the ideal geodesic quasdrilateral with the vertices  $\xi, y_n, y'_n, \xi'$ is $3\delta$-slim according to Lemma \ref{lem:idealthinness}(4), we obtain that every point $p\in \xi \xi'$ belongs to a $3\delta$-neighborhood of the segment $y_ny'_n$ if $n$ is sufficiently large. Thus, $p$ belongs to the $(R+3\delta)$-neighborhood of the orbit $\Ga x$. It follows that the $\Ga$-action on $\hull \La$ is cobounded. 

\medskip 
2.  Suppose that $x\in \hull \La$ and $R<\infty$ are such that $\hull \La$ is contained in $N_R(\Ga x)$. Consider a geodesic segment $xx'$ in $X$, where $x'\in \Ga x\subset \hull \La$. Thus, we find two geodesics $\xi\eta, \xi'\eta'$ in $\hull \La$ such that $x\in \xi\eta, x'\in \xi'\eta'$. We will consider the ``generic'' case when $\xi, \xi'$ are all different and leave the equality case to the reader. We have an ideal quadrilateral $Q=\xi x x'\xi'$. Since $Q$ is $3\delta$-slim, the segment $xx'$ is contained in the $3\delta$-neighborhood of the union of the infinite sides of $Q$. But these infinite sides are contained in $\hull \La$. Thus, 
$$
\hull \Ga x\subset N_{3\delta} \hull \La\subset N_{3\delta +R}(\Ga x),
$$
as required. 
\end{proof}

\begin{lem}
A metrically proper  discrete subgroup $\Ga< \Isom(X)$ is convex-cobounded if and only if $\Ga$ is finitely generated 
and for one (equivalently, every) point $x\in X$ the orbit map $o_x: \Ga\to \Ga x\subset X$ is a  quasiisometric embedding. 
\end{lem}
\begin{proof} 1. Suppose first that $\Ga$ is convex-cobounded. The action of $\Ga$ on $X_0:=\hull \Ga x\subset X$ (equipped with the restriction of the metric of $X$) is metrically proper and cobounded. Since $X_0$ is roughly geodesic, Milnor--Schwartz lemma applies and, hence, $\Ga$ is finitely generated and $o_x:  \Ga\to X_0$ is a quasiisometry. 

2. Conversely, suppose that $\Ga$ is finitely generated and $o_x: \Ga\to X$ is a  quasiisometric embedding. We extend 
$o_x$ to a quasiisometric embedding $f: Cay_\Ga\to X$ of the Cayley graph of $\Ga$. Consider points $x, y=\ga x$ in $\Ga x$. The map $f$ sends a geodesic in $Cay_\Ga$ connecting $1$ to $\ga$ to a uniform quasigeodesic $q$ in $\Ga x$ connecting $x$ to $y$.  Morse Lemma implies that every geodesic $xy\subset X$ is uniformly close to $q$. We conclude that $\Ga x$ is quasiconvex, implying that $\Ga$ is convex-cobounded. 
\end{proof}

\begin{cor}\label{cor:boundary map}
If $\Ga$ is convex-cobounded, then it is Gromov-hyperbolic and there exists an equivariant homeomorphism
$$
h: \geo \Ga\to \La=\La(\Ga)\subset \geo X. 
$$
\end{cor}
\begin{proof} By the lemma we get an equivariant quasiisometry $f: Cay_\Ga\to X_0= \hull \Ga x$. Since $X_0$ is roughly geodesic and $X$ is Gromov-hyperbolic, it follows that $X_0$ is also Gromov-hyperbolic. The space $X_0$, of course, is not geodesic but its injective hull $E(X_0)$ is and $f: Cay_\Ga\to E(X_0)$ is a quasiisometry. Hyperbolicity of $E(X_0)$ then implies hyperbolicity of $Cay_\Ga$. Thus, $\Ga$ is a hyperbolic group. The visual boundary of 
$X_0$ is $\La$.  Hence, the quasiisometry $f$ extends to an equivariant homeomorphism $h: \geo \Ga\to \La$ by Theorem \ref{thm:boundary extension}. 
\end{proof}

Note that a converse to this corollary holds as well, but we will not need this.

\begin{thm}\label{thm:qccc}
Suppose that $X$ satisfies the visibility and condition and $\Ga< \Isom(X)$ is a discrete subgroup acting metrically properly on $X$. Then the following are equivalent:

1. A discrete subgroup $\Ga< \Isom(X)$ is convex-cobounded.  

2. $\La$ is compact and every limit point of $\Ga$ is uniformly conical. 

3. $\La$ is compact and every limit point of $\Ga$ is conical. 

4. $\La$ is compact and every limit point of $\Ga$ is horospherical. 






\end{thm}

We refer to \cite[Theorem 12.2.7]{DSU} for a proof. Note that, in particular, this theorem says that convex-cobounded subgroups $\Ga$ have compact limit sets. Thus, each orbit $X_0=\Ga x$ with metric induced from $X$ has semiproper injective hull $E(X)$.  


\section{Redressing in hyperbolic groups}\label{sec:redressing}

The material of this section is taken from a joint work of Michael Kapovich and Didac Martinez-Granado, \cite{KM}. 
The theorem we will need is:

\begin{thm}
\label{thm:redressing}
Assume that $(X,d)$ is a $\delta$-hyperbolic, proper, geodesic metric space, equipped with a geometric action $\Gamma\times X\to X$   
of a nonelementary hyperbolic group $\Gamma$. Then for every $y_0\in X$ there exists a finite subset $E\subset \Gamma$ and constants $L, K>0$ satisfying the following.  For every $\ga\in \Gamma$, there exists $g\in E$ such that $\{g^{-1}(y_0), y_0, \gamma (y_0)\}$ lies on 
an $(L,K)$-quasi-axis $\rho$ of a loxodromic isometry $\tilde{\gamma}= \gamma\circ g$.   
\end{thm}

\index{$\tilde\gamma=\gamma\circ g$, a redressing of $\gamma$} 
The element $\tilde{\ga}$ is said to be {\em a redressing} of $\ga$ with respect to $y_0$. See Appendix A, Chapter \ref{sec:appendixA}, for a proof of the theorem. 

\medskip
Recall that (by the Morse Lemma) each $(L,K)$-quasigeodesic in $X$ is uniformly close to a geodesic. Thus, there exists $R$ depending only on $L, K$ and $\delta$ (hence, depending only on $y_0$ and $\Ga\acts X$) such that for $A:= A_{\tilde\ga}$ and a  complete geodesic $l\subset A$ we have 
\begin{equation}\label{defi:R}
\{g^{-1}(y_0), y_0, \gamma (y_0)\}\subset \rho \subset N_R(l) \subset N_R(A).  
\end{equation}

Let $m_0$ denote the maximum of distances 
$d(y_0, g(y_0)), g\in E$. Below we will be using neighborhoods $U(\zeta, t)\subset \ol{X}$ of points $\zeta\in \geo X$ defined in 
\eqref{eq:U-nbd}.   

The following is   an addendum to the Redressing Theorem:

\begin{thm}\label{thm:redressing2}
There exists a constant $\kappa$ depending only on $X$ and the action $\Ga\acts X$, such that for all $\zeta\in \geo X$ and $t>0$, if  
$\ga (y_0)\in U(\zeta,t)$, then  the attracting fixed point $\xi^+$ of $\tilde{\gamma}$ belongs to $U(\zeta, \kappa t)$. 
\end{thm}
\begin{proof} Set $y:= \ga(y_0)$, $x:=\tilde\ga(y_0)$. Then $d(x, y)\leqslant m_0$. We also have 
$$
d(y_0, A)\leqslant R. 
$$
By Theorem \ref{thm:fixed-point-estimate} applied to the loxodromic isometry $\tilde\ga$, for the point $x=\tilde \ga(y_0)$ we get
$$
(x,\xi^+)_{y_0} \geqslant (x,\zeta)_{y_0} - R',
$$
with $R'=R+ 25\delta$. However, what we need is a lower bound in terms of $(y,\zeta)_{y_0}$. To get one, we use the inequality 
\eqref{eq:Lip}, which yields
\begin{equation}\label{eq:xyzeta_estimate}
|(x,\zeta)_{y_0}- (y,\zeta)_{y_0}|\leqslant d(x,y) + 2\delta'\leqslant d(x,y) + 6\delta\leqslant m_0 + 6\delta. 
\end{equation}
Combining this with the inequality \eqref{eq:Gh1} (applied to the points $\zeta, x, \xi^+$), we get 
\begin{align} 
\begin{aligned}\label{eq:zeta-estimate}
(\xi^+, \zeta)_{y_0}\geqslant \min \left( (\xi^+, x)_{y_0} , (x, \zeta)_{y_0}\right)\geqslant (x,\zeta)_{y_0} - R' \geqslant \\   (y,\zeta)_{y_0} - R' - m_0 - 6\delta=
(y,\zeta)_{y_0} - R -m_0 -31\delta. 
\end{aligned}
\end{align}
We then take 
$$
 \kappa= e^{R +m_0 +31\delta}.
$$
By the inequality \eqref{eq:zeta-estimate} we get 
$$
y=\ga(y_0)\in U(\zeta, t)\Rightarrow \xi^+\in U(\zeta, \kappa t). 
$$
\end{proof} 

A sort of converse to Theorem \ref{thm:redressing2} is also true. Let $R$ be as in \eqref{defi:R} 
(it depends only  only on $X$ and the action $\Ga\acts X$).

\begin{thm}\label{thm:redressing3}
There exist 
 a constant $\kappa'>0$, both depending only on $X$ and the action $\Ga\acts X$, such that for all $h\in \Ga$, 
$\zeta\in \geo X$ and $t > 0$, if $d(y_0, h(y_0))> - \log(t)$, $d(y_0, A_h)\le R$ and if the attracting fixed point $\xi^+$ of $h$ belongs to $U(\zeta, t)$, then  
$x:=h(y_0)$ is in $U(\zeta, \kappa' t)$. 
\end{thm}
\begin{proof} 
Let $l$ be a complete geodesic in $A=A_\ga$ such that $d(x, l)\le R$. 
By Corollary \ref{cor:1}, applied to $\eps=0$ and the points $p=y_0$, $q=x$ and 
$\bar{q}=P_l(x)$, we have 
$$
(x,\xi^+)_{y_0}\geqslant d(y_0, x) - R - 9\delta. 
$$
Combining with  the inequality \eqref{eq:Gh1}, we obtain, 
\begin{align}\label{eq:yzeta_estimate}
\begin{split}
(x,\zeta)_{y_0}  \geqslant \min( (x, \xi^+)_{y_0}, (\zeta,\xi^+)_{y_0}) -\delta' \geqslant \\
\min( d(y_0, x) - R - 9\delta, -\log(t))-3\delta\geqslant \\
 \min( d(y_0, x) - R - 9\delta, -\log(t))-3\delta. 
 \end{split}
\end{align}

Therefore, 
we get that $d(y_0, x)\geqslant -\log(t)$ implies 
$$
(x,\zeta)_{y_0}  \geqslant   -\log(t) - R - 9\delta \geqslant  -\log(\kappa' t), 
$$
where $\kappa'= e^{9\delta + R}$.  Thus, $x\in U(\zeta, \kappa' t)$. 
\end{proof}

\section{Geodesic flows of hyperbolic groups}\label{sec:Geodesic flows of hyperbolic groups}


\subsection{Definitions}

Let $\Gamma$ be a non-elemen\-tary hyperbolic group with a fixed word metric $d_\Ga$ and Cayley graph $Cay_\Ga$.  Consider a {\em geodesic flow} $\widehat\Ga$ of $\Ga$;  such a  flow was originally constructed by Gromov \cite{Gromov} and then improved by Champetier \cite{Champetier} and Mineyev \cite{Mineyev2005}, resulting in definitions with different properties. We note that the exponential convergence of asymptotic geodesic rays will not be used in our discussion; 
as we will see, it is also irrelevant whether the trajectories of the geodesic flow are geodesics or uniform quasigeodesics in $\widehat\Ga$. 
In particular, it will be mostly irrelevant for us which definition of $\widehat\Ga$ is used. (Apart from sections \ref{sec:Measures on geodesic flows of hyperbolic groups} and  \ref{sec:Mineyev's flow}.) 
The following properties of $\widehat\Ga$ will be used in the sequel:

1. $\widehat\Ga$ is a proper metric space. 

2. There exists a 
properly discontinuous isometric action 
$\Ga\acts\widehat\Ga$. 

3. There exists a $\Ga$-equivariant quasi-isometry $\Pi: \widehat\Ga\to \Ga$; in particular, the fibers of $\Pi$ are relatively compact.

4. There exists a continuous action $\R\acts\widehat\Ga$, denoted $\hat\phi_t$ and called the {\em geodesic flow}, \index{$(\widehat\Ga,\hat\phi_t)$, 
geodesic flow of a hyperbolic group} 
whose trajectories are uniform quasigeodesics in $\widehat\Ga$, i.e. for each $\hat m\in \widehat\Ga$ the 
{\em flow line} 
$$
t\to \hat m_t := \hat\phi_t(\hat m)
$$ 
is a uniform quasi-isometric embedding $\R\to \widehat\Ga$. 

5. The flow $\hat\phi_t$ commutes with the action of $\Ga$. 

6. Each $\hat m\in\hat\Ga$ defines a uniform  
quasigeodesic 
 in $\Ga$ by the formula: 
$$
t\mapsto \Pi(\hat m_t), t\in \R.  
$$
Following the notation in \S \ref{sec:conical}, we let $\geo^2\Ga$ denote the subset of $\geo\Ga\times\geo\Ga$ 
consisting of pairs of distinct points. 
The natural map 
$$
 \re=(\re_+, \re_-): \widehat\Ga\to \geo^2\Ga 
$$
assigns to $\hat m\in  \widehat\Ga$ the pair $(\xi_+, \xi_-)$, with points $\xi_\pm\in \geo \Ga$ defined as 
$$
\xi_\pm=  \lim_{t\to\pm\infty}  \Pi(\hat m_t)\in \geo \Ga. 
$$
The map $\mathrm{e}$ is continuous and surjective. (Recall that $ \geo^2\Ga$ is the complement to the diagonal in $\geo \Ga\times \geo \Ga$.)

In particular,
every uniform quasigeodesic in $\widehat\Ga$ 
is uniformly Hausdorff close to a flow line.

\medskip 
The reader can think of the elements of $\widehat\Ga$ {\em as parameterized geodesics in the Cayley graph 
$Cay_\Ga$, so that $\hat\phi_t$ acts on geodesics via reparameterization}. This was Gromov's original viewpoint, 
although not the one in \cite{Mineyev}.  In Mineyev's work, the space $\hat\Ga$, as a topological space, is homeomorphic to 
$\geo^2 \Ga \times \R$. 

We let $\cG:=\cG(\Ga)$ denote the quotient space $\widehat{\Ga}/\Ga$ and let $m$ denote the projection of 
$\hat{m}\in \widehat{\Ga}$ to $\cG$. Since the flow $\hat\phi_t$ commutes with the $\Ga$-action, it descends to a flow 
$\phi_t$ on $\cG$. \index{$(\cG,\phi_t)$, geodesic flow of a hyperbolic group}

We say that $\hat m\in \widehat\Ga$ is {\em normalized} if $\Pi(\hat m)=1_\Ga\in\Ga$. \index{normalized element of $\widehat\Ga$} 
It is clear that every $\hat m\in \widehat\Ga$ can be sent to a normalized element of $\widehat\Ga$ 
via the action of $m_0^{-1}\in \Ga$.

Since trajectories of $\hat\phi$ are uniform quasigeodesics, for 
each normalized $\hat m\in \widehat\Ga$ we have 
\begin{equation}\label{eq:distanceF}
C_1^{-1} t - C_2\leqslant d_\Ga(1, \Pi(\hat m_t)) \leqslant C_1 t + C_2
\end{equation}
for some positive constants $C_1, C_2$. 

\medskip

\subsection{Examples}

\subsubsection{Example: Free groups} 

Suppose that $\Ga$ is the free group of rank $r$ on generators $a_1,..., a_r$. As $\widehat\Ga$ we will take the set of geodesics in the corresponding Cayley graph $Cay_\Ga$: These are encoded by biinfinite sequences $(g_n)_{n\in \bZ}$ of elements $g_n\in \Ga$ such that $d_S(g_n, g_{n+1})=1$ for all $n\in \bZ$. 
Every such sequence $(g_n)$ yields a geodesic $c=\hat m$ in $Cay_\Ga$. The $\R$-action on geodesics $c$ restricts to the $\bZ$-action on the set of sequences $(g_n)$, which is the shift:
$$
(k, (g_n))\mapsto (g_{n+k}), k, n\in \bZ. 
$$
The projection $\Pi: \widehat \Ga\to \Ga$ corresponds to the map $(g_n)\mapsto g_0$ on sequences and, more generally, sends a geodesic 
$c: \R\to Cay_\Ga$ to $\Pi(c)=c(0)$ if $c(0)\in \Ga$, and if $c(0)\in (g, h)$, then $\Pi(c)=g$. Here $(g, h)$ is the edge in 
$Cay_\Ga$ connecting $g, h\in \Ga$, with the vertices $g, h$ removed. Thus, normalized geodesics are those which have $c(0)=g_0=1$. 
They correspond to biinfinite reduced words in the alphabet $\{a_i^{\pm 1}: i=1,...,r\}$.  Topologically, $\widehat\Ga$ is just 
$\geo^2\Ga\times \R$, where $\geo^2\Ga$ is noncompact, locally compact and totally disconnected.

\subsubsection{Example: Fundamental groups of compact negatively curved Riemannian maniifolds}

Let $M$ be a compact connected negatively curved Riemannian manifold, $UTM$ is the unit tangent bundle of $M$ and $\phi$ is the geodesic flow on $UTM$. The unit tangent bundle lifts to the unit tangent bundle $UT\tilde{M}\to \tilde M$ of the universal cover of $M$. The fundamental group $\Ga=\pi_1(M)$  acts properly discontinuously on $UT\tilde{M}$ with the quotient $UTM$. Using a Levi-Civita connection on $\tilde M$ one lifts the Riemannian metric from $\tilde M$ to a $\Ga$-invariant Riemannian metric on $UT\tilde M$. With this metric, $UT\tilde M$ is a Gromov-hyperbolic geodesic metric space and, thus,  
is naturally quasiisometric to $\Gamma$. Topologically, $UT\tilde M$ is the product $\geo^2\tilde M\times \R\cong \geo^2\Ga\times \R$. The geodesic flow 
on $\tilde M$ defines an $\R$-action $\hat\phi$ on  $UT\tilde M$, which projects to the geodesic flow $\phi$ on $UTM$. The quasiisometric projection 
$\Pi: UT\tilde M\to \Ga$ is the composition of the projection  $UT\tilde M\to \tilde M$ with the quasiisometry $q: \tilde M\to \Ga$ given by a quasi-inverse to 
the orbit map $\Ga \to \Ga x\subset \tilde M$, where $x$ is a base-point in $\tilde M$. For instance, we can take a bounded fundamental subset $F\subset \tilde M$ for the $\Ga$-action so that $x\in F$ and define 
$$
q(y)= \ga, y\in \ga F, y\in \tilde M, \ga\in \Ga. 
$$
We then take $\cG=UTM$, $\widehat\Ga=  UT\tilde M$. 
Then normalized geodesics $\hat m: \R\to UT\tilde M$ are lifts $(c(t), c'(t)), t\in \R$, of unit speed geodesics $c: \R\to \tilde M$ such that $c(0)\in F$.

\subsection{Mineyev's version of the geodesic flow} \label{sec:Mineyev's flow}

The version of a geodesic flow defined by Mineyev in \cite{Mineyev2005} has sharper properties than the ones listed earlier. Below are some details 
(none of which will be used elsewhere in the book). 
Mineyev equips the space $\widehat\Ga= \geo^2\Ga \times \R$ with a roughly geodesic $\Ga$-invariant metric $d$ such that 
the flow-lines are geodesic (not just uniformly quasigeodesic). Note that Mineyev uses the notation $d_*$ for this metric. 
Mineyev further defines certain continuous functions 
$\beta^\times_\zeta(x,y)$ (satisfying the cocycle condition) on 
$\geo \Ga\times \widehat\Ga\times \widehat\Ga$ which play the role of  
Busemann cocycles (Theorem 55 in  \cite{Mineyev2005}). (We will not attempt to define $d$ and $\beta^\times$, and refer the reader to \cite{Mineyev2005}.)  
While these functions {\em are not} Busemann cocycles for the metric $d$, this collection of functions is $\Ga$-invariant:
$$
\beta^\times_{\ga \zeta}(\ga x, \ga y)
$$
and  every function $\beta^\times_\zeta(x,\cdot)$  is linear with slope $-1$ on every trajectory of the flow 
forward-asymptotic to $\zeta$. 

To define a $\Ga$-action on $\widehat\Ga= \geo^2\Ga \times \R$ which projects to the given $\Ga$-action on $\geo^2\Ga$ and is by translations on the lines $\{z\}\times \R$, one needs a continuous cocycle 
$$
\tau: \Ga \times \geo^2\Ga\to \R,
$$
$$
\tau(\ga_2\circ \ga_1, z)= \tau(\ga_2, \ga_1 z) + \tau(\ga_1, z), \ga_1, \ga_2\in \Ga, z\in  \geo^2\Ga. 
$$
Then the action would be given by 
$$
\ga \cdot (z, t)= (\ga z, t- \tau(\ga, z)), z=(\xi,\eta)\in \geo^2\Ga. 
$$
For this, one uses
$$
\tau(\ga, \xi, \eta)= \half(b_\xi(\ga^{-1},1) - b_\eta(\ga^{-1}),1))=   (\ga^{-1},\eta)_1 - (\ga^{-1},\xi)_1. 
$$

One defines {\em horospherical foliations} on $\widehat\Ga$ as follows. 
Given $\zeta\in \geo \Ga=\geo \widehat\Ga$ we set
$$
\widehat{W}^-_\zeta:= \{\hat m\in \widehat\Ga: e_+(\hat m)=\zeta\}, 
$$
a leaf of the {\em stable lamination} on $\widehat\Ga$. Similarly,  one defines 
$$
\widehat W^+_\zeta:= \{\hat m\in \widehat\Ga: e_-(\hat m)=\zeta\}, 
$$
a leaf of the {\em unstable} foliation. These leaves are clearly closed and pairwise disjoint. Moreover, since $e_\pm$ are restrictions 
of projections of $\geo X \times \geo X\times \R$ to the first/second factor, the partition of $\widehat\Ga$ into stable/unstable leaves are foliations, called 
{\em stable/unstable foliations} of $\widehat\Ga$. Since the maps $e_\pm$ are $\Ga$-equivariant, both foliations are $\Ga$-invariant. 
Next, one defines {\em strongly stable and strongly unstable leaves} $\widehat W^{s}_{\hat m}, \widehat W^{u}_{\hat m}$ through 
points $\hat m\in \widehat\Ga$. Given a point $\hat m\in \widehat\Ga$ take the level set 
$$
\widehat W^{s}_{\hat m}:=\{x\in \widehat{W}^-_{\xi_+}: \beta^\times_{\xi_+}(x)=\beta^\times_{\xi_+}(\hat m)\},
$$
where $\xi_+=e_+(\hat m)$; this is the  
 {\em strongly stable} leaf through $\hat m$. Similarly, for $\xi_-=e_-(\hat m)$, 
$$
\widehat W^{u}_{\hat m}:=\{x\in \widehat{W}^-_{\xi_-}: \beta^\times_{\xi_-}(x)=\beta^\times_{\xi_-}(\hat m)\},
$$
is a {\em strongly unstable}  leaf. 

Since every function $\beta^\times_\zeta$ is linear with slope $-1$ on every trajectory of the flow 
forward-asymptotic to $\zeta$, every such trajectory intersects every strongly stable leaf contained in $\widehat W^-_\zeta$. 
The same applies to strongly unstable leaves. In particular, for every function $b=\beta^\times_\zeta$ and $t\in \mathbb R$ we obtain a homeomorphism
\begin{equation}\label{eq:bdry-to-horosphere}
h_b: \geo \Ga=\geo \widehat \Ga\setminus \{\zeta\} \to \widehat W^{s}_{\hat m},
\end{equation}
(where $b(\hat m)=t$ and $\zeta=e_+(\hat m)$), sending $\xi$ to the unique intersection of $\widehat W^{s}_{\hat m}$ with the unique flow-line connecting 
$\xi$ to $\zeta$.  Thus, the partitions of $\widehat\Ga$ into strongly stable/unstable leaves again form foliations of $\widehat\Ga$. 
The (strongly) stable/unstable foliations of $\widehat \Ga$ are invariant under the 
flow and the $\Ga$-action, hence, project to flow-invariant partitions of $\cG$. These partitions of $\cG$ 
are not quite foliations in the traditional sense since the $\Ga$-action is not free: They are foliations in the case when $\Ga$ is torsion-free.

\medskip
Lastly and, most importantly, the flow defined by Mineyev has a ceryain {\em exponential contraction property} which we describe below. 
Fix a number $R$ such that the $\Ga$-orbit of every $R$-ball in $\widehat\Ga$ equals the entire $\widehat\Ga$. Suppose that 
$\zeta\in\geo \Ga$ and $x, y$ are points in $\widehat\Ga$ such that:

\begin{enumerate}
\item $d(x,y)\le R$.

\item 
$y\in W^s_x$. 

\end{enumerate}

Then for every $t\ge 0$,
\begin{equation}\label{eq:contract}
d(\hat\phi_t(x), \hat\phi_t(y))\le Ae^{-ct}
\end{equation}
 for some uniform constants $A>0$ and $c>0$ which depend only on $\widehat\Ga$ and the choice of $R$. See Theorem 
 60(h) in \cite{Mineyev2005} for details.


\chapter{Currents in complex geometry}\label{sec:currents}


In this section we review currents on  K\"ahler manifolds; we refer the reader to \cite[Chapter 4]{morgan} and \cite[Chapter III]{deRham} 
in the case of real manifolds and to \cite[I.2 and III.2]{demailly} in the complex case. We also 
 prove some basic results about currents used later in the book. 

\section{Generalities of currents}\label{sec:Generalities of currents}

\subsection{Definitions}

Let $M$ be a compact $n$-dimensional K\"ahler manifold with the K\"ahler form $\omega$. We will define and use {\em only order $0$} currents on $M$; we refer to \cite{deRham} and \cite[I.2 and III.2]{demailly} for the treatment of general currents.  Let $\Om^*(M)$ denote the topological vector space of smooth complex differential forms on $M$ (i.e. smooth sections of the complexified cotangent bundle), on which the topology is defined by the $C^0$-norm, 
$$
||\varphi||_{\infty}=\sup_{x\in M} \sup_{\nu} \varphi(\nu), \varphi\in \Om^k(M),$$
where the second supremum is taken over all orthonormal $k$-tuples $\nu$ in $T_xM$. 
The space of (order $0$) {\em $k$-dimensional currents} on $M$, ${\cD}_k(M)$, \index{${\cD}_k$, the space of $k$-dimensional currents} is the topological dual space 
of $\Omega^k(M)$. We also refer to $2n-k$ as the {\em degree} of a $k$-dimensional current and, accordingly, use the notation $\cD^{2n-k}(M)$ for this space. \index{${\cD}^d$, the space of degree $d$ currents}
The pairing between currents $T$ and forms $\varphi$ 
is denoted  $\langle T , \varphi\rangle$. \index{$\langle T , \varphi\rangle$, pairing of a current and a form} 
Let $\Om^{p,q}(M)$ denote the space of  smooth forms  of bidegree $(p,q)$ on $M$; its dual is ${\cD}_{p,q}(M)={\cD}^{n-p,n-q}(M)$, 
\index{${\cD}_{p,q}$, the space of  $(p,q)$-bidimensional currents} 
the space of currents of {\em bidimension $(p,q)$ and {\em bidegree $(n-p,n-q)$. 
\index{${\cD}^{p,q}$, the space of currents of bidegree $(p,q)$} 
Then
$$
{\cD}_k(M)= \bigoplus_{p+q=k} {\cD}_{p,q}(M). 
$$

We will be primarily using the {\em topology of weak convergence} on currents: 
$T_i\to T$ (weakly) if for every $\varphi\in \Omega^k(M)$, $\langle T_i, \varphi\rangle \to \langle T, \varphi \rangle$. 
Given a current $T\in {\cD}_m(M)$ and a form $\varphi\in \Omega^k(M)$ with $k\leqslant m$, one defines 
the {\em wedge product} $T\wedge \varphi$ as a current of dimension  $m-k$ such that for $\psi\in \Omega^{m-k}(M)$ 
$$
\langle T\wedge \varphi, \psi\rangle = \langle T , \varphi\wedge \psi \rangle. 
$$
Like measures, currents can be {\em pushed forward} by appropriate smooth maps. Namely, if $f \colon M \to N$ is a smooth map of compact manifolds 
and $T\in \cD_k(M)$, then  the current $f_* T\in \cD_k(N)$ is determined by the formula: 
\[ \langle f_* T , \psi\rangle = \langle T , f^* \psi \rangle\]
for every $\psi\in \Om^k(N)$.  When $f$ is a diffeomorphism, one also defines the {\em pullback} of $T$, $f^* T$, as the 
pushforward by the inverse map: $f^* T:=(f^{-1})_* T$.  A current $T$ is said to be {\em supported} 
on an open subset $U\subset M$ if for every form $\varphi\in \Om^k(M)$ whose support is disjoint from $U$, we have $\langle T, \varphi\rangle =0$. Accordingly, the {\em support} of $T$, denoted $\supp(T)$, is defined as the complement to the union of open subsets $U\subset M$ such that $\langle T, \varphi\rangle=0$ for all forms $\varphi$ supported 
in $U$.

\subsection{Exterior derivatives and homology classes} 

Given a degree $k$ 
current $T$, one defines its  {\em exterior derivative}, denoted $dT$, 
as the adjoint of the exterior derivative of differential forms, namely: \index{$dT$, the exterior derivative of a current}
{ 
\begin{equation}
\langle dT, \varphi \rangle = {(-1)^{k+1}}\langle T , d\varphi \rangle,
\end{equation}
}
for $\varphi\in \Om^{k-1}(M)$. (Note that this sign convention is different from the one for the {\em boundary} of a current, $\partial T$.) A current is said to be {\em closed} if $dT=0$. Currents on $M$ equipped with the boundary operator form a chain complex, \index{closed current}\index{$[T]$, homology class of a current}
whose homology is the dual of the de Rham cohomology of $M$, i.e. is naturally isomorphic to $\coh_*(M; \C)$. 
Accordingly, we will regard homology classes $[T]$ of closed $k$-dimensional currents $T$ on $M$ as elements of $\coh_k(M; \C)$. Homology classes of currents of bidimension $(p,q)$ belong to $\coh_{p,q}(M; \C)$. 

\subsection{Mass of a current} 

The {\em mass} $||T||$ (also frequently denoted ${\mathbf M}(T)$) of a (order $0$) current $T$ 
\index{${\mathbf{M}(T)}$,  the mass of a current} is defined by dualizing the norm $||\cdot ||_{\infty}$ on forms. Let $\{U_\alpha: \alpha\in A\}$ be a finite coordinate cover of $M$. Every current { $T\in \cD^{p,q}(M)$} can be written as $\sum_{\alpha\in A} T_\alpha$ where each $T_\alpha$ is supported in $U_\alpha$. In local coordinates, identifying $U_\alpha$'s with open subsets of $\C^{n}$,  
one can write each $T_\alpha$ as
$$
\sum_{I,J} T_{\alpha,I,J} dz_I\wedge \ol{dz}_J, \ dz_I= dz_{i_1} \wedge \ldots \wedge  dz_{i_{p}}, \ 
 \ol{dz}_J= \ol{dz}_{j_1}\wedge \ldots \wedge  \ol{dz}_{j_{q}} 
$$
where each $T_{\alpha,I,J}$ is a finite complex measure on $U_\alpha$, see \cite[Ch I, Prop. 2.9]{demailly}. 


The mass of $T$ is sandwiched by the total variations of $T_{\alpha,I,J}$:
$$
C^{-1} \int_M \sum_{\alpha,I,J} |T_{\alpha,I,J}|\leqslant  ||T||\leqslant C \int_M \sum_{\alpha,I,J} |T_{\alpha,I,J}|
$$
for some uniform constant $C\geqslant 1$ independent of $T$. See 
\cite[Remark 1.15 p.134]{demailly}. 
  
Just as distributions, currents satisfy a compactness theorem with respect to the topology of weak convergence: 
The subset of $\cD_k(M)$ consisting of currents of mass $\leqslant m<\infty$, 
is weakly compact (see e.g. \cite[Lemma 26.13]{Simon}). 

In the case where the currents $T$ is represented by smooth forms, we will use the following lemma later on.

\begin{lem} \label{lem_sup_coef} Let $\psi$ be a smooth form supported on a compact set $K $ contained in an affine coordinate chart $U_\alpha \subset M$. Consider  the decomposition of $\psi$ written in (local) coordinates on $U_\alpha$  as: 
\begin{equation*}
\psi = \sum_{I,J} f_{\alpha I J} dz_I \wedge \overline{dz_J}, 
\end{equation*}
where $f_{\alpha I J}$ are smooth  functions whose support is contained in $K$. Then there exists a constant $C> 0$ such that for each multi-index $I,J$: 
\begin{equation*}
\sup_{K } |f_{\alpha I J}| \leqslant C || \psi ||.
\end{equation*}
\end{lem}

The proof of this lemma is straightforward and we omit it. 

\subsection{Integration of currents} 

Let $Z$ be a topological space equipped with a finite Borelean measure $\mu$ and $f: Z\to  
\cD_m(M), f(z)=T_z$, be an $L^\infty$-map, i.e. a measurable map such that 
$$
\int_Z ||T_z||d\mu<\infty. 
$$ 
Then one defines the current \index{$\int_Z T_zd\mu$, the integral of a current-valued map} 
$$
\int_Z T_zd\mu\in \cD_m(M)
$$
by 
\begin{equation}\label{eq:current integration} 
\left< \int_Z T_zd\mu, \varphi\right>=  \int_Z \left<  T_z, \varphi \right> d\mu. 
\end{equation}

The main examples of currents of concerns in our book are built from smooth forms and from \textit{currents of integration} along complex subvarieties. If $Z'$ is a Borel subset of a compact analytic subset $Z$ of complex dimension $k$ in $M$, one defines a current  of bidimension $(k,k)$, denoted $\llbracket Z'\rrbracket$, 
\index{$\llbracket Z\rrbracket$, current of integration} 
whose pairing $\langle \llbracket Z'\rrbracket, \varphi \rangle$ with a smooth $(k,k)$-form $\varphi$ on $M$, 
is given by: 

\begin{equation}
\langle \llbracket Z'\rrbracket, \varphi \rangle = \int_{Z^{reg}\cap Z'} \varphi,
\end{equation}
where the integration takes place on the smooth locus $Z^{reg}$ of $Z$. 


\begin{rem}
The notation for the currents of integration in the literature is rather inconsistent and they are frequently denoted $[Z]$. 
We decided against using this notation in order to avoid clash with the notation for the homology class $[Z]$ of $Z$. 
Our notation $\llbracket Z\rrbracket$ follows \cite{Simon}. 
\end{rem}

\subsection{Positive currents} \label{section_positive_currents}

Note that $\llbracket Z'\rrbracket$ is a closed current if and only if $Z'$ is a closed subset of $Z$. Currents of the form $\llbracket Z'\rrbracket$ are specific in the sense that they are currents of order zero, and they are also {\em positive} in the sense of \cite[Definition 1.13]{demailly}. More precisely, a current $T$ is positive if  \index{positive current}
$\langle T, \varphi \rangle \geqslant 0$ for any  smooth {\em strongly positive} form $\varphi\in \Omega^{p,p}(M)$, 
i.e. a form which can be written locally (in holomorphic coordinates $(z_1 , \ldots, z_n)$) as 
\begin{equation*}
\varphi =  \lambda(z_1, \ldots, z_n)  idz_1 \wedge d \bar z_1 \wedge \ldots \wedge idz_p\wedge d\bar z_p, 
\end{equation*} 
where $\lambda \geq 0$ is a smooth real non-negative function. The space of positive currents of bidimension $(k,k)$ (resp. bidegree $(n-k,n-k)$) 
will be denoted $\cD^+_{k,k}(M)$ (resp. $\cD^{n-k,n-k}_+(M)$). 

\subsection{Estimates of masses of currents}


 \begin{lem} \label{lem_mass_estimate} Fix a
 K\"ahler form $\omega$ on $M$.
 There exists a constant $C>0$ such that for every $T\in \cD^{p,p}_+(M)$, 
 \begin{equation}
 || T || \leqslant C \langle T, \omega^{n-p} \rangle.
\end{equation}  
In particular, for any smooth form $\varphi$ of bidegree $(n-p,n-p)$, one has: 
\begin{equation*}
| \langle T, \varphi \rangle| \leqslant C ||\varphi||_\infty \langle T,\omega^{n-p} \rangle.
\end{equation*}
 
\end{lem} 
  
  
  \begin{proof} Note that the second inequality follows directly from the first one.

  Using a partition of unity, we reduce to the case where $\varphi$ is supported in a coordinate neighborhood of a point $x \in M$. Choose some local holomorphic coordinates on this neighborhood and write the current  $T$ locally as: 
  \begin{equation*}
  T = \sum_{|I|=|J|=p} T_{I,J} dz_I \wedge d\bar z_J,
\end{equation*}   
  where $T_{I,J}$ are complex measures and where the coefficients $T_{I,I}$ are positive measures. 
Furthermore, we decompose the K\"ahler form locally near $x $ as $\sum f_{i,i} dz_i \wedge d\bar{z}_{i}$ where $f_{i,i}$ is a positive smooth function.
Since each $f_{i,i}$ is smooth, it is uniformly bounded above and below by $L>0$ and $L^{-1}$ on the coordinate neighborhood of the point $x$ given by the partition of unity.

  By \cite[Proposition 1.14]{demailly} applied to $\lambda_k = 1$ for all $k$, for all multi-indices $I,J$ of length $p$, 
  the total variations of these complex measures satisfy: 
  \begin{equation*}
  |T_{I,J}| \leqslant 2^p \sum_{K}  T_{K,K},  
\end{equation*}   
where the sum is taken over all multi-indices $K$ containing $I\cap J$ and contained in $I\cup J$.
Observe that $\omega^{n-p}$ is locally of the form $ \omega^{n-p} = \sum_{\substack{|I'| = n-p}} g_{I'} dz_{I'} \wedge d\bar{z}_{I'} $, 
where $g_{I'} = \prod_{i\in I'} f_{i,i}$.
We get 
\begin{multline*}
|\langle T , \omega^{n-p} \rangle| = \int_M \sum_{\substack{I,J, I'\\ |I|= |J|= p \\
|I'| = n-p}} |\langle  T_{I,J} , g_{I'}\rangle| dz_{I} \wedge dz_{I'} \wedge d\bar{z}_{J} \wedge d\bar{z}_{I'} \\ = \int_M \sum_{|K|=p} | \langle T_{K,K}, g_{K^c} \rangle | dz_1 \wedge \ldots \wedge  dz_{n} \wedge d\bar{z}_1 \wedge \ldots \wedge d\bar{z}_n 
\geqslant L^{-1} \sum_{|K|=p} \int_M T_{K,K}. 
\end{multline*}
Using the previous inequality, we obtain: 
\begin{multline*}
||T||=\int_M \sum_{ |I| = |J|=p} |T_{I,J}| \leqslant \int_M  \sum_{ |I|= |J|=p } 2^p \sum_{\substack{I\cap J \subset K \subset I\cup J \\ |K| = p} } T_{K,K} \leqslant 2^p \int_M \sum_{\substack{|I|=|J|=p}}  \sum_{|K|=p} T_{K,K} \\ = 2^p {n \choose p}^2 \left ( \sum_{|K|=p} \int_{M} T_{K,K} \right ) \leqslant  2^p {n \choose p}^2  L \langle T, \omega^{n-p} \rangle, 
\end{multline*}
as required.
  \end{proof}

We shall apply  the above mass estimate as follows:
 
\begin{cor} \label{cor_difference_current}
Let $V$ be a (closed) subvariety of codimension $p$ in $M$ for which there is a partition $V = W\sqcup R$, 
where $W$ is an open subset in $V$. Consider a closed positive current $T$ of bidegree $(p,p)$, a smooth form $\varphi$ of bidegree $(n-p,n-p)$ and a constant $\epsilon>0$  satisfying the following conditions:
\begin{enumerate}
\item[(i)] $[T] = [V]$ in $\coh_{p,p}(\mathcal{F})$; and 
\item[(ii)] 
 for all $\Psi  \in \{ \varphi , \omega^{n-p}\}$  where $\omega$ is a K\"ahler form on $M$, one has: 
$$ \left | \langle \llbracket W\rrbracket, \Psi \rangle - \langle T, \Psi \rangle \right | \leqslant  \epsilon \max (||\Psi||_\infty, ||d\Psi||_\infty ). $$
\end{enumerate} 
Then there exists a constant $C> 0$ (independent of $\varphi$) such that:
\begin{equation*}
| \langle T -   \llbracket V\rrbracket, \varphi  \rangle | \leqslant C\epsilon \max (\lVert\varphi\rVert_\infty, \lVert d\varphi\rVert_\infty ). 
\end{equation*}
\end{cor}

\begin{rem} While Corollary \ref{cor_difference_current} is stated for a given test form $\varphi$, it will be used for arbitrary test forms as long as the inequality in the condition (ii) of the corollary  is satisfied. 
\end{rem}

\begin{proof}
Observe that since the currents $T, \llbracket V\rrbracket$ are closed and positive, the first condition implies that $ \langle T, \omega^{n-p} \rangle = \langle  \llbracket V\rrbracket, \omega^{n-p} \rangle $.
Since $ \llbracket V\rrbracket =  \llbracket W\rrbracket +  \llbracket R\rrbracket$ is a sum of positive currents, the second condition gives: $$ 0 \leqslant \langle  \llbracket R\rrbracket , \omega^{n-p} \rangle = \langle  \llbracket V\rrbracket -  \llbracket W\rrbracket , \omega^{n-p} \rangle = \langle T, \omega^{n-p} \rangle - \langle  \llbracket W\rrbracket, \omega^{n-p} \rangle \leqslant C_1 \epsilon,  $$
where $C_1>0$ depends only on $\omega$. 
Applying Lemma \ref{lem_mass_estimate} to the current $ \llbracket R\rrbracket$ and to $T - \intcur{W}$, we see that there is a uniform constant $C> 0$ such that $| \langle\intcur{R},\varphi \rangle| \leqslant C C_1 ||\varphi||_\infty \epsilon $. Using the decomposition  $\intcur{V} = \intcur{W} + \intcur{R}$ together with the triangle inequality, one obtains 
\begin{align*}
 |\langle T -  \llbracket V\rrbracket, \varphi \rangle| &\leqslant |\langle T-  \llbracket W\rrbracket ,\varphi \rangle | + |\langle \intcur{R} , \varphi \rangle| \\
   & \leqslant  \epsilon \max(||\varphi ||_\infty, ||d\varphi ||_\infty) + C C_1 \epsilon ||\varphi||_\infty,
\end{align*}
as required.
\end{proof} 

\medskip

\section{Positive currents trapped in a cylinder}\label{sec:Positive currents trapped in a cylinder}



In this section, we provide some estimates on the difference of positive currents whose support is contained in a cylinder-like region. 
To explain the principle we shall first give the statement for measures, which are currents of dimension zero. 

\begin{lem} Let $S\subset M$ be a finite subset.
Let $U \subset M$ be an open subset containing $S$, whose diameter is smaller than $\alpha >0$. Then for any smooth function $\varphi$ on $M$ and any point $p\in U$, one has: 
\begin{equation}
|\langle \delta_S -  (\#S)\delta_p, \varphi \rangle | \leqslant \# S\alpha   ||d\varphi||_\infty.
\end{equation}
Here, $\delta_S$ denotes the sum of Dirac masses at points of $S$ and $\delta_p$ is the Dirac mass at $p$. 
\end{lem}

\begin{proof}
For each point $s\in S$, take $\sigma_s = [p,s]$, a minimal geodesic segment in $M$ joining $p$ to $s$. 
Applying Stokes theorem to the restriction 
$\varphi_s$ of $\varphi$ to $\sigma_s$, we obtain, by triangle inequality:
\begin{align*}
|\langle \delta_S -  (\#S)\delta_p ,\varphi \rangle | = 
| \sum_{s \in S} \langle \delta_s - \delta_p, \varphi \rangle| &\leqslant \left | \sum_{s \in S} \int_{\sigma_s} d \varphi_s \right | \\
 & \leqslant  \# S ||d\varphi||_\infty \alpha, 
\end{align*}
as required.
\end{proof}

In order to control the behavior of higher dimensional currents, we shall apply similar ideas using the flat norm on currents. For a current $T$ on $\mathcal{F}$ (of order zero) the {\em flat norm} $N_F(T)$ of $T$ (see \cite[p. 41]{morgan}) is defined by: 
\index{$N_F(T)$, the flat norm of a current $T$}
\begin{equation} \label{eq_def_flat_norm}
N_F(T) = \inf_{T = A + d B}  \{ ||A|| + ||B|| \}, 
\end{equation}
where the infimum is taken over all currents $A, B$ such that $T = A + d B$. 
As an immediate consequence of the definition,  
for any smooth form $\varphi$ on $M$ one has 
\begin{equation} \label{eq_apply_flat_norm}
|\langle T , \varphi\rangle| \leqslant  N_F(T) \max  (||\varphi||_\infty , ||d\varphi||_\infty ).  
\end{equation}

We state some basic feature of the flat norm. 
\begin{lem} \label{prop_pushforward_flat}
Let $\pi \colon M \to N$ be a holomorphic map between compact K\"ahler manifolds with $\dim M = n$. { Then there is a constant $C<\infty$ such that every 
current $T$ on $M$ of finite flat norm satisfies}
\begin{equation*}
N_F(\pi_*T) \leqslant {\kap C}N_F(T). 
\end{equation*}  	
\end{lem}
\begin{proof} It follows from the definition of the pushforward that for any current $T= A + d B$ on $M$ we have 
\begin{align*}
	\pi_* T& = \pi_* A + \pi_* d  B \\ 
	& = \pi_* A + d \pi_*B.
\end{align*}

On the other hand, compactness of $M$ implies that there exists a constant $C<\infty$ depending only on $M, N$ and $\pi$ such that  $|| \pi_*A|| \leqslant C ||A||$ and 
$|| \pi_* B || \leqslant C ||B||$. For  $\eps>0$ we take currents $A, B$ such that $T= A + d B$ and 
$$
|| A || + ||B|| \leqslant N_F(A)+\eps. 
$$
Combining the above inequalities, we obtain 
\begin{align*}
	N_F(\pi_* T) &\leqslant || \pi_* A || + || \pi_* B|| \\ 
	& \leqslant C (|| A || + ||B||)\\ 
	& \leqslant C (N_F(T) + \epsilon).
\end{align*} 
Since this holds for arbitrary $\eps>0$, lemma follows. \end{proof}

Our use of the flat norm in this book comes from the following example.

\begin{example}\label{example:cylinder}
Let $D\subset \C$ be a unit disk and  $D_1$ and $D_2$ be two discs in $\C^2$ given by $D \times \{ \alpha\}$ and $D\times \{\beta \}$, 
where $\alpha, \beta \in \C$ satisfy $|\alpha - \beta| \leqslant \epsilon$ for some $\epsilon>0$.  
Then the discs $D_1, D_2$ together with the cylindrical shell 
$$C =\partial D \times \{ t \alpha + (1-t)\beta | t\in [0,1] \} $$ 
define the boundary of the solid cylinder $F$, which is the convex hull of $D_1, D_2$. 
In terms of currents, this means that: 
\begin{equation}
\llbracket D_1\rrbracket - \llbracket D_2\rrbracket + \llbracket C\rrbracket = d F. 
\end{equation}
In particular, this shows that $\llbracket D_1\rrbracket - \llbracket D_2\rrbracket = -\llbracket C\rrbracket + d F$. 
Furthermore, the mass of  the cylindrical shell $\llbracket C\rrbracket$ is equal to the surface area 
$2 \pi   \epsilon$, and the mass of $F$ is equal to the volume, which is $\pi \epsilon$. 
Overall, this shows that $N_F(\llbracket D_1\rrbracket - \llbracket D_2\rrbracket) \leqslant 3\pi  \epsilon$.  
  \end{example}

The next lemma is immediate:  

\begin{lem} \label{lem_interpolation} Let $n > 0$ and $k \in [0,n]$ be two integers. Fix a complex-linear projection 
$\pi_h\colon \C^n= \C^{k}\times \C^{n-k} \to \C^k $ and a subset 
an 
 $S\subset \C^n$. 
Then   
\begin{equation}
I(S) = \{ (1-t)x + t \pi_h(x) \ | \ x\in S, t\in (0,1)   \},
\end{equation} 
 is the image of the map $H: S \times (0,1)\to \C^n$, 
 $$(y,t) \mapsto (1-t)y + t \pi_h(y).$$
\end{lem}

We will refer to $I(S)$ as an {\em interpolating subset} between $S$ and $\pi_h(S)$. 


 In what follows, $U\subset \C^k$ will be a bounded domain whose closure 
is a real semi-algebraic subset. We will use the notation $\vol_m$ for the $m$-dimensional Hausdorff measure (which, up to a uniform multiplicative constant,  is comparable to the $m$-dimensional Riemannian volume). 

For a subset of $\bC^k$, such as $U$, we will denote by $\partial U$ its boundary in $\bC^k$ (not in $\bC^n$).
%



 \begin{prop}\label{prop_subvariety_smaller_cylinder}
Let $n\dyl$ and $k\in [0,n-1]$ be two integers. Let $\pi_h$ be the orthogonal projection of $\bC^n=\bC^k\times \bC^{n-k}$ onto $\bC^k$ where $\C^n$ is endowed with the euclidean metric. 
Fix $D\dyl$, a compact set $K\subset \C^n$ and consider $U=B_{\C^k}(0,D)$, the open Euclidean ball of radius $D$ in $\bC^k$. 
Let $S'$ be a complex-algebraic subvariety of $\bC^n$. 
Assume that 
\begin{enumerate}
\item[(i)] the restriction of $\pi_h$ to $S_U':= S'\cap \pi_h^{-1}(U)$ induces a ramified covering of degree $a$ onto $U$; 
\item[(ii)]
there is a number $\fd \leqslant  1$  such that for all $x\in S_U'$ we have:
\begin{align*}
\begin{aligned}
||x- \pi_h(x)||& \leqslant \fd, 
\end{aligned}
\end{align*}

\item[(iii)] $U \times \{0 \}\cup S_U'\subset K$.
\end{enumerate}

Then there is a constant $C> 0$ which depends on $K$ and the degree $\deg(S')$, such that   for 
every smooth $(k,k)$-form $\varphi$ compactly supported in $K$, 
we have
\[
\left\lvert\langle \llbracket S_U' \rrbracket - a\llbracket U  \times \{0\} \rrbracket, \varphi \rangle\right\rvert \leqslant {C\sqrt{\fd}} \max(  \lVert \varphi \rVert_{\infty}, \lVert d\varphi\rVert_\infty).
\]
\end{prop}

\begin{proof} 

Consider the interpolating set $A=I({S_U'})$ and denote by $\phi : S_U' \times (0,1) \to \C^n$ the homotopy map given in Lemma \ref{lem_interpolation}. 
The map $\phi$ is smooth outside of the points $(p,t)$ where $p$ is a singular point of $S_U'$. 
We consider the current of integration on $S_U' \times (0,1)$ (which exists by e.g \cite{hartd_cie}) and we have: 
\begin{equation}
d \intcur{S_U' \times (0,1)} = \intcur{S_U' \times \{1 \}} -  \intcur{S_U' \times \{0 \}} + \intcur{\partial S_U' \times (0,1)}. 
\end{equation}  
Pushing forward by $\phi$ and using the fact that $\pi_h$ restricted to $S_U' $ has degree $a$, one obtains:
\begin{multline*}
d \intcur{A} = d \phi_* \intcur{S_U'\times (0,1)} = \phi_* d \intcur{S_U' \times (0,1)} = \phi_* \intcur{S_U' \times \{1 \}} - \phi_* \intcur{S_U' \times \{0 \}} + \phi_*\intcur{\partial S_U' \times (0,1)} \\= \intcur{S_U'} - a \intcur{U\times 0}  + \phi_* \intcur{\partial S_U' \times (0,1)}.  
\end{multline*}
Now, the image of $\partial S_U' \times (0,1)$ is exactly $\overline{A} \cap \pi_h^{-1} (\partial U) $, so that we have:
\begin{equation*}
d \intcur{A} = \intcur{S_U'} - a \intcur{U\times \{ 0 \}} + \intcur{\overline{A} \cap \pi_h^{-1}(\partial U) }. 
\end{equation*}

%

By definition of the flat norm, one obtains:
\begin{equation} \label{eq_flat_norm_main}
N_F(\intcur{S_U'} - a \intcur{U\times \{0\}}) \leqslant   \vol_{2k+1} (A) + \vol_{2k}(\overline{A} \cap \pi_h^{-1}(\partial U)).
\end{equation} 

In the following Lemma, we will use the notion of degree depending on the choice coordinates given by $\C^n$, see \S~\ref{appendix_volume_degree} for a discussion of various notions of degrees. 

\begin{lem} \label{lem_two_bound} There exists a constant $C_1> 0$ depending only on the degrees of 
 ${S'}$, of ${S'} \cap \pi_h^{-1}(\partial U)$, and on the diameter of $K$  such that 
one has $\vol_{2k+1}(A) \leqslant C_1 \sqrt{\fd} $ and $\vol_{2k}(\pi_h^{-1}(\partial U) \cap \overline{A}) \leqslant C_1 \sqrt{\fd}$.
\end{lem}

\begin{proof}
As our computation is in a compact set, we endow $\C^n$ with the euclidean metric. The two estimates for $A$ and $\pi_h^{-1}(\partial U) \cap \overline{A}$ are very similar. 
Denote by ${S_{U,reg}'}$ the regular locus of ${S_U'}$. The complement ${S_U' \setminus S_{U,reg}'}$ is a real algebraic set of dimension at most $2k-2$. Since the homotopy  $H$ is a { real} algebraic map, the image $H({S_U' \setminus S_{U,reg}'} \times (0,1)) = \phi(S_U' \setminus S_{U,reg}' \times (0,1))$ is also a real semi-algebraic set of dimension at most $2k-1$ in $\C^n$. 
In particular, this gives: 
\begin{equation*}
\vol_{2k+1}(A)= \vol_{2k+1} (\phi(S_U' \times (0,1)) = \vol_{2k+1} ( \phi (S_{U,reg}' \times (0,1))). 
\end{equation*}
Let $u \colon S_{U,reg}' \to \C^n$ be the embedding of $S_{U,reg}'$ into $\C^n$ so that $\phi$ restricted to $S_{U,reg}'$ is given by: 
\begin{equation*}
\phi \colon (x,t) \in S_{U,reg}' \times (0,1) \mapsto (1-t) u(x) + t \pi_h (u(x)) \in I(S_{U,reg}').  
\end{equation*}
The differential of $\phi $ at $(x,t) \in S_{U,reg}'\times (0,1)$ is of the form
\begin{equation*}
D \phi_{(x,t)}(v, w) = (1-t) D_x u(v) + t \pi_h D_x u(v) - w ( u(x) - \pi_h (u(x))), 
\end{equation*}
for any $(v,w) \in T_{(x,t)} S_{U,reg}' \times (0,1)$.  
 Take $v_1, \ldots, v_{2k}, v_{2k+1}$ a basis of $T_{(x,t)} S_{U,reg}' \times (0,1)$ where the first $2k$ vectors form an orthonormal basis of $T_x S_{U,reg}'$ for the euclidean metric induced by the embedding in $\C^n$. 
The jacobian of $\phi$ at $(x,t) \in S_{U,reg}' \times (0,1)$  is by \cite[Assertion (2) of Lemma 5.1.4]{krantz_parks} the square root of the determinant of the matrix
\begin{equation}
M(x,t)= (\langle D\phi_{(x,t)}(v_p) , D\phi_{(x,t)}(v_q) \rangle)_{p,q \leqslant 2k+1}.
\end{equation}
By assumption (ii), we have $||u(x) - \pi(u(x))|| \leqslant \fd$, implying that 
$$|| D\phi_{(x,t)} (v_{2k+1})|| = || u(x) - \pi_h (u (x))|| \leqslant \fd.$$  
Using the Cauchy--Schwartz inequality and the fact that $\pi_h$ is an orthogonal projection, we see that entries of the last column of the matrix $M(x,t)$ are bounded by a multiple of $\fd$, i.e. for any $p \leqslant 2k$ 
 \begin{multline*}
|\langle D_{(x,t)}\phi(v_p) , D_{(x,t)}\phi(v_{2k+1}) \rangle |  = | \langle D_{(x,t)} \phi(v_p) , u(x) - \pi_h (u(x)) \rangle| \leqslant ||D_{(x,t)} \phi(v_p)|| \fd  \\ \leqslant ||D_x u(v_p) || \fd \leqslant \fd,   
 \end{multline*}
 where we have used the fact that the first $2k$ vectors form an orthonormal basis for the euclidean metric. 
 Similarly, we see that the components in the first $k$ columns are bounded by $1$. Hence we deduce that $\sqrt{\det M(x,t) } \leqslant \sqrt{\fd} $ for all $(x,t) \in S_{U,reg}' \times (0,1)$. 
By the area formula \cite[Theorem 5.1.9]{krantz_parks}, we have: 
\begin{multline*}
\vol_{2k+1}(A)  \leqslant \int_{\C^n} \# \{ S_{U,reg}' \times (0,1) \cap \phi^{-1}(y)\} d\vol_{2k+1}(y)\\ =  \int_{S_{U,reg}' \times (0,1)} \sqrt{\det M(x,t)} d\vol(x,t)  
\leqslant  \sqrt{\fd} \vol_{2k+1}(S_{U,reg}' \times (0,1)) = \sqrt{\fd} \vol_{2k}(S_U').  
\end{multline*}

Using the fact that $\pi_h^{-1}(\partial U) \cap \overline{A} = I(\pi_h^{-1}(\partial U) \cap S')$, the proof of the second bound follows from the same method and yields $\vol_{2k}(\pi_h^{-1}(\partial U) \cap \overline{A}) = \vol_{2k } I(\pi_h^{-1}(\partial U) \cap S') \leqslant  \sqrt{\fd} \vol_{2k-1}(\pi_h^{-1}(\partial U) \cap S')  $. 
Taking $C_2 = \max ( \vol_{2k}(S_U'), \vol_{2k-1}(S'\cap \pi_h^{-1}(\partial U)))$ yields the required inequalities. 
 We then conclude using Proposition \ref{prop_volume_degree} that $C_2$ is bounded above by a constant which depends on the degrees of $U, \partial U$, $S'$ and on the diameter of $K$.
\end{proof}

We now conclude the proof of the result using Lemma \ref{lem_two_bound} together with \eqref{eq_flat_norm_main} and \eqref{eq_apply_flat_norm}.
\end{proof}

We will use the previous proposition in a global context as follows. 
\begin{cor} \label{cor_subvariety_smaller_cylinder}
Let $S$ be a complex subvariety of (complex) dimension $k$ contained in the full flag manifold $\cF$. Fix $O$  a Zariski open subset in $\cF$ which is isomorphic to $\C^n$ via $ \Psi \colon O \to \C^n$  where $n = \dim \cF$.  Consider the  map $\tilde{\pi_h} \colon O \mapsto  O\cap V$ conjugate via the isomorphism $\psi$ to an orthogonal projection, where $V$ is a complex subvariety of dimension $k$ given by the Zariski closure of $\Psi^{-1}(\C^k \times \{ 0\})$. 
Fix a compact set $\tilde{K}$ in $O$, a real number $D>0$ and consider $\tilde{U} = V \cap \Psi^{-1} (B(0,D))$.  Assume that: 
\begin{enumerate}
\item[(i)] The restriction of $\tilde{\pi_h}$ to $S_{\tilde{U}} := S \cap \tilde{\pi_h}^{-1} (\tilde{U})$ induces a ramified covering of degree $a$ onto $\tilde{U}$;
\item[(ii)] there is a number $\fd \leqslant 1$ such that for all $x \in S_{\tilde{U}}$, we have: 
\begin{equation*}
 || x - \tilde{\pi_h}(x) || \leqslant \fd, 
\end{equation*}
\item[(iii)] $\tilde{U} \cup S_{\tilde{U}} \subset K$.
\end{enumerate}
 Then there exists a constant $C>0$ which depends on $K$ and the degree $\deg_{\cF}(S)$, such that for every smooth $(k,k)$-form $\varphi$ compactly supported in $\tilde{K}$, we have: 
 \begin{equation*}
 | \langle \intcur{S_{\tilde{U}}} - a \intcur{\tilde{U}} ,\varphi \rangle  | \leqslant C \fd \max(||\varphi||_\infty, || d\varphi||_{\infty}).
 \end{equation*}
\end{cor}

\begin{proof}
Set $S' =\Psi(S\cap O)$, $U = \Psi(\tilde{U})$, $K = \Psi(\tilde{K})$, $\pi_h =\Psi \circ  \tilde{\pi_h} \circ \Psi^{-1}$ and  $S_U' = S' \cap \pi_h^{-1}(U) $. 
Then $S' $, $\pi_h$ and $U$ satisfy the hypothesis of Proposition \ref{prop_subvariety_smaller_cylinder}. 
 Given a smooth  $(k,k)$-form $\varphi$ supported in $\tilde{K}$, the form $(\Psi^{-1})^* \varphi $ is supported in $K$ so the previous proposition gives: 
 \begin{multline*}
 | \langle \intcur{S_{\tilde{U}}} - a \intcur{\tilde{U}} ,\varphi \rangle| = | \langle \intcur{S_{U}'}  - a \intcur{U \times \{0 \}} , (\Psi^{-1})^* \varphi \rangle| \leqslant C \sqrt{\fd} \max (||(\Psi^{-1})^* \varphi||_\infty , ||d(\Psi^{-1})^* \varphi||_\infty) \\
 \leqslant C \sqrt{\fd} \max(||\varphi||_\infty, || d\varphi||_{\infty}). 
 \end{multline*}
 where $C>0$ is a constant which only depends on the degree of $\deg_{\R}(S')$ (viewed as a real algebraic subset) and $K$. 
 Fix a K\"ahler form $\Omega$ on $\cF$, then by Proposition \ref{prop_degree_real_vs_complex}, there exists a constant $C'>0$ such that $\deg_{\R}(S') \leqslant C' \deg_{\Omega}(S)$ where $\deg_{\Omega}(S)$ is the degree with respect to the K\"ahler form $\Omega$.  
 Overall, this shows that $C$ only depends on $\deg_{\Omega}(S)$ and on $\tilde{K}$, as required.
\end{proof}

\section{Wedge products of currents and Bedford-Taylor's method}
\label{section_bedford}

Currents of order zero can be viewed as differential forms whose coefficients are measures. In particular, if one wants to define the  wedge product of two currents, one has to give a sense to the multiplication of two measures. 
In some instances, there is a notion of wedge product, due to Bedford and Taylor when the currents have bidegree $(1,1)$.

Recall that a function $u \colon \Omega \to \R \cup \{ - \infty \} $ on an open subset $\Omega \subset \C^n$ is called  \textit{ plurisubharmonic} if it is upper-semicontinuous and 
for every $a \in \Omega$, $\xi \in \C^n$, one has: 
\begin{equation*}
u(a) \leqslant \dfrac{1}{2\pi} \int_{[0,2\pi]} u(a + e^{i\theta}\xi) d\theta
\end{equation*}
A global function $u\colon \cF \to \R \cup \{ -\infty\} $ is plurisubharmonic on $\cF$ if its restriction to any complex chart is plurisubharmonic. We denote by $\PSH(\cF)$ the space of plurisubharmonic functions on $\cF$. 
 For these particular functions, their  $dd^c = \dfrac{i}{2\pi} \partial \bar \partial$  yields a particular $(1,1)$ current, denoted $dd^c u$ which is closed and positive.    
  For example, if $ Z_1 , Z_2$ are the subvarieties which realize the zeros and pole locus of a meromorphic function $f$ on $\cF$, then $\log |f|$ is plurisubharmonic and the Poincar\'e Lelong shows that the $dd^c$ yields a difference of integration currents. 
  $$dd^c \log |f| = \intcur{Z_1} - \intcur{Z_2}. $$
  The above formula generalizes if $\sigma$ is a meromorphic section of a vector bundle $E$ endowed with a metric $||\cdot||$ as follows:
  \begin{equation*}
 i dd^c \log || \sigma|| = \intcur{\sigma = 0} - \intcur{\sigma = \infty} - i c_1(E),
  \end{equation*}
  where $c_1(E)$ is the Chern curvature of the bundle $E$.
\smallskip   
  
  Given $u \in \PSH(\cF)$ and a closed positive $T$ current of bidegree $(p,p)$,  Bedford and Taylor define the wedge product between $dd^c u$ and $T$ as
  \begin{equation} \label{eq_wedge} 
  dd^c u \wedge T:= dd^c(uT),
  \end{equation}
  provided the function $u$ is locally bounded on every compact subset, see e.g \cite[p. 145]{demailly}\index{$dd^c u \wedge T$, Bedford-Taylor wedge product between two currents}.  
  Equation \eqref{eq_wedge} is a definition of the wedge product, however it does not show what the wedge product is in concrete examples.
In certain favorable situations, the intersection of currents is not only defined but also coincides with the geometric intersection:  

\begin{prop} [Proposition 4.12, p.156, in \cite{demailly}]  
\label{prop_demailly_algebraic}
Assume $V$ is an algebraic subvariety in $\cF$ of dimension $k \geqslant 1$ and $D_1, \ldots, D_m$ are effective divisors (more precisely, codimension $1$ projective subvarieties in $\cF$) such that every irreducible component of 
$D_1 \cap \ldots \cap D_m \cap V$ has dimension $k-m$. Then the Bedford-Taylor wedge product 
$$\llbracket D_1\rrbracket \wedge \ldots \wedge \llbracket D_m\rrbracket \wedge \llbracket V\rrbracket $$ 
is a well-defined 
current 
which equals
\begin{equation} \label{eq_geometric}
\llbracket D_1 \rrbracket \wedge  \ldots \wedge \llbracket D_m \rrbracket \wedge \intcur{V} = \sum_{C \in D_1 \cap \ldots \cap D_m \cap V } m_C \intcur{C},    
\end{equation} 
where $C$ runs through the irreducible component of the intersection $D_1 \cap \ldots \cap D_m \cap V$ and where $m_C$ is an integer corresponding to the multiplicity of $C$ in the intersection. 
\end{prop} 

Equation \eqref{eq_geometric} means that the intersection  of currents defined via the wedge-product  is geometric in the sense the intersection coincides with the  (weighted) geometric intersection. 

  We now recall some properties of so-called woven currents. 
  
  \begin{defi} \label{defi_woven} A current $T$ on a compact K\"ahler manifold is uniformly woven if locally on any open set $\Omega$, there exists a family of subvarieties $A_\alpha$ with uniformly bounded volume, and a positive measure $\mu$ such that $T = \int \intcur{A_\alpha} d\mu(\alpha)$.   
  \end{defi}
When the family of subvarieties form a lamination, then the current is called uniformly laminar.   
  The following result is due to Dujardin concerning the wedge product of bidegree $(1,1)$ woven currents.
  
  \begin{thm} \cite[Theorem 3.1]{dujardin_support_bifurcation} Let $\cF$  be a complex projective variety. For $i=1, \ldots, k$, let $T_i = \int \intcur{A_{i,\alpha}} d\mu_i(\alpha) $ be a current with locally bounded potential where $\mu_i$ are positive measures and where $A_{i,\alpha}$ are some families of hypersurfaces on $\cF$ and let $S$ be an irreducible analytic subset of $\cF$, then the intersection 
  $T_1 \wedge \ldots \wedge T_k \wedge \intcur{S}$ exists and is geometric in the sense that: 
  \begin{equation}
  T_1 \wedge \ldots \wedge T_k \wedge \intcur{S} = \int  \intcur{A_{1, \alpha_1}} \wedge \ldots \wedge \intcur{A_{k,\alpha_k}} \wedge \intcur{S}  d\mu_1(\alpha_1) \ldots d\mu_k(\alpha_k),
\end{equation}     
where one sets $\intcur{A_{1, \alpha_1}} \wedge \ldots \wedge \intcur{A_{k,\alpha_k}} \wedge \intcur{S} =0$ whenever the intersection is not proper. 
  \end{thm} 
  
  As a corollary of the previous result, one can integrate over the subvarieties $S$. 
  
  \begin{cor}\label{cor_intersection_woven}Let $\cF$  be a complex projective variety.  Let $T = \int \intcur{S_\alpha} d\mu_0(\alpha)$ where $\mu_0$ is a positive measure and where $S_\alpha$ is a family of codimension $k$ subvarieties on $\cF$. For $i=1, \ldots, p$, consider $T_i = \int \intcur{A_{i,\alpha}} d\mu_i(\alpha)$, some uniformly woven currents of bidegree $(1,1)$ with  bounded potentials; then the wedge product $T_1 \wedge \ldots \wedge T_p \wedge T$ exists and one has: 
  \begin{equation}
  T_1 \wedge \ldots \wedge T_p \wedge T = \int \intcur{A_{1,\alpha_1}}\wedge \ldots \wedge  \intcur{A_{p,\alpha_p}} \wedge\intcur{S_\alpha} d\mu_0(\alpha) d\mu_1(\alpha_1) \ldots d\mu_p(\alpha_p).
  \end{equation}
  \end{cor}
  
\begin{proof} The proof of this result is very similar to Dujardin's proof, and proceeds by induction on 
  $0 \leqslant p \leqslant  k $.
  We give the inductive step between $p=0$ and $p=1$ as the next steps are similar. 
Since the result is local, we can reduce to the case where each hypersurface $A_{1,\alpha_1}$ is defined locally as $\{ w = h_{\alpha_1}(z)\}$ for some local coordinates $(w,z)$ where $|h_{\alpha_1} (z)|< 1$.
The local potential of $T_1$ is then $u_1 = \int \log |w - h_{\alpha_1}(z)| d\mu_1(\alpha_1)$ which is negative.  
Since $T$ is a positive current and since $ u_1 T$ has finite mass, one has: 
\begin{equation*}
u_1 T = \int \log |w - h_{\alpha_1}(z)| T d\mu_1(\alpha_1)
\end{equation*}
Applying $dd^c$ we get: 
\begin{equation*}
dd^c (u_1 T) = dd^c \int \log |w - h_{\alpha_1}(z)| T d\mu_1(\alpha_1).  
\end{equation*}
Similarly, for $\mu_1$ almost every point $\alpha_1$, the current $ \log |w-h_{\alpha_1}(z)| T$ has finite mass, hence we have: 
\begin{equation*}
\log |w- h_{\alpha_1}(z)| T = \int \log |w- h_{\alpha_1}(z)| \intcur{S_\alpha} d\mu_0(\alpha).
\end{equation*}
 This gives: 
 \begin{align*}
 dd^c (u_1 T) &=  dd^c \int \log |w - h_{\alpha_1}(z)| \intcur{S_\alpha}  d\mu_1(\alpha_1) d\mu_0(\alpha) \\
 & = \int dd^c (\log |w - h_{\alpha_1}(z)| \intcur{S_\alpha}) d\mu_1(\alpha_1) d\mu_0(\alpha)\\
 & =   \int \intcur{ A_{1,\alpha_1}} \wedge \intcur{S_{\alpha}} d\mu_1(\alpha_1) d\mu_0(\alpha).
 \end{align*}
\end{proof}  
  


\chapter{Gibbs measures and counting on Gromov-hyperbolic spaces}\label{sec:counting}


\section{Patterson--Sullivan measures for real-hyperbolic spaces} 

Patterson--Sullivan measures were first introduced by Patterson for {\em Fuchsian groups} (discrete isometry groups of the hyperbolic plane) and then generalized by Sullivan to higher-dimensional real-hyperbolic spaces, see \cite{Patterson, sullivan_density}. The book by Nicholls, \cite{nicholls}, is a comprehensive treatment of this theory in the setting of real-hyperbolic spaces. In this section we will review some basic facts and relate two different versions of the construction (one using hyperbolic geometry and the second using derivatives).

\subsection{Geometry of real-hyperbolic spaces} 

In what follows, we will be using the unit ball model $\mathbb B^n$ of $\bH^n$ and let $||\cdot||$ denote the Euclidean norm.  
The hyperbolic Riemannian metric on $\mathbb B^n$ is given by the metric tensor
$$
ds^2= 4\frac{dx_1^2+...+dx_n^2}{1-||x||^2}. 
$$
We will not need a general formula for the hyperbolic distance between points in $\bH^n$, it will suffice to know that 
for a point $p\in \mathbb B^n$ the hyperbolic distance $d(0, p)$ is given by  
\begin{equation}\label{eq:distance} 
d(0, p)= \log \frac{1+||p||}{1-||p||}.
\end{equation}
See e.g. \cite[p. 14]{nicholls}.

We first discuss the relation between Busemann functions and Poisson kernels for $\bH^n$.  The {\em Poisson kernel} $k(w,\zeta)$ on $\bB^n$ 
equals 
$$
\frac{1-||w||^2}{||\zeta-w||^2},  \quad w\in \mathbb B^n, \quad \zeta\in \partial \mathbb B^n. 
$$
In particular, $k(0,\zeta)=1$. Take $x, x'\in \mathbb B^n$. Then, according to  Lemma 3.2.1 in \cite{nicholls}, the Busemann cocycle $b_\zeta(x,x')$ is related to the 
Poisson kernel by 
$$
b_\zeta(x,x')= \log k(x',\zeta) - \log k(x,\zeta). 
$$
Setting $x'=0$ and normalizing Busemann functions to vanish at that point, we get
$$
b_\zeta(x)=b_\zeta(x,0)= - \log k(x,\zeta). 
$$
Note that in the upper half-space model, for $\zeta=\infty$ we have
$$
k(x,\zeta)= x_n, 
$$ 
where we use the point $(0,...,0,1)$ for our normalization of Busemann functions.

\medskip 
For a Moebius transformation $g$ of the closed unit ball we let $J_g(z)$ denote the Euclidean operator norm of the derivative of $g$ at the point $z$. 
For every $z\in \mathbb B^n$ we have 
\begin{equation}\label{eq:Jacobian} 
J_g(z)= \frac{1-||g(z)||^2}{1-||z||^2}, 
\end{equation}
see Theorem 1.3.3 in \cite{nicholls}. 
If $w=g^{-1}(0)$ then for $\zeta\in \partial \mathbb B^n$ we also have
\begin{equation}\label{eq:Jac0}
J_g(\zeta)= \frac{1-||w||^2}{||\zeta-w||^2}=k(w,\zeta),
\end{equation}
see Theorem 1.3.4  in \cite{nicholls}.

Thus, for $w=g^{-1}(0)$ the expression $-\log J_g(\zeta)= -\log k(w,\zeta)$ 
equals the value $b_\zeta(w)$ of the Busemann function on $\bH^n$ normalized at $0$ (cf. Corollary \ref{cor:conformality of visual metrics}),  
\begin{equation}\label{eq:J_g}
J_g(\zeta)= e^{-b_\zeta(g^{-1}(0))}= e^{b_\zeta(0,g^{-1}(0))}= e^{b_{g\zeta}(g(0),0)}= e^{b_{g\zeta}(g(0))}. 
\end{equation}

\subsection{Poincar\'e series}\label{sec:Poincare1}

For a discrete group $\Ga$ of Moebius transformations of $\ol{\mathbb B^n}$ (a {\em Kleinian group}) we will define several closely related series associated with certain points $z\in \ol{\mathbb B^n}$. 

\medskip
{\bf The Jacobian Poincar\'e series:}  
\begin{equation}\label{eq:PoincareSeries} 
P_z(s):=\sum_{\ga\in \Ga} (J_\ga(z))^s\equiv \sum_{\ga\in \Ga} J^s_\ga(z). 
\end{equation}

{\bf The metric  Poincar\'e series}: \index{Poincar\'e series $\cP(s)$}
\begin{equation}\label{eq:Poincare0}
\cP(s) = \cP_{x,y}(s)=\sum_{\ga \in \Ga} e^{-s d(x, \gamma y)}.
\end{equation}

For a positive function $a: \Ga\to \R_+, a(\ga)=a_\ga$ and a general series of the form
$$
\sum_{\ga\in \Ga} a_\ga^s, 
$$
we define the {\em exponent of convergence} as the infimum of those values of $s$ for which this series converges. We will see 
in Section \ref{sec:reviewPS} that the exponent of convergence of the metric Poincar\'e series does not depend on the choice of the 
points $x, y$ (for now we just take it for granted), 
it is called the {\em critical exponent} $s_{0,\Ga}$ of $\Ga$. Arguably, it is the single most important numerical invariant of a Kleinian group.

\begin{lem}
[Cf. Lemma 4.1 in \cite{FMS}.] \label{lem:2CE}
Let $\Ga< \Isom(\bH^n)$ be a discrete subgroup with domain of discontinuity $\Omega\subset \partial \mathbb B^n$. 
Then the following numbers are equal:

(1) The critical exponent $s_{0,\Ga}$ of $\Ga$. 

(2) The exponent of convergence of the series
$$
\sum_{\ga\in \Ga} (1-||\ga(z)||)^s  
$$
for every $z\in \bB^n$. 

(3) The exponent of convergence of the series \eqref{eq:PoincareSeries} $P_z(s)$ for $z\in \Omega \cup \mathbb B^n$. 

In particular, the latter exponent is independent of $z$. 
\end{lem}
\begin{proof} Recall  that for a point $p\in \mathbb B^n$,   
$$
d(0, p)= \log \frac{1+||p||}{1-||p||},
$$
see \eqref{eq:distance}. Therefore, for every $z\in \mathbb B^n$, 
$$
\exp(-s d(0, \gamma(z)))= \left(\frac{1+||\gamma(z)||}{1-||\gamma(z)||} \right)^{-s}= \left(\frac{1-||\gamma(z)||}{1+||\gamma(z)||} \right)^{s}. 
$$
Since $1\le 1+||\gamma(z)||\leqslant 2$, we obtain
$$
2^{-s} (1-||\gamma(z)||)^s \leqslant \exp(-s d(0, \gamma(z))) \leqslant  (1-||\gamma(z)||)^s. 
$$
This proves the equality of (1) and (2). 

Consider $z\in \mathbb B^n$. Since $r=||z||<1$, the equation \eqref{eq:Jacobian} implies the inequality 
$$
1-||g(z)|| \leqslant J_g(z)\leqslant \frac{2}{1-r^2}(1-||g(z)||). 
$$
This proves the equality of (2) and (3) for $z\in \mathbb B^n$. 

Now, consider a point $\zeta\in \Omega$. Then $||\zeta-\gamma^{-1}(0)||^2$ is bounded away from zero by a positive 
constant $C$ for all $\ga\in \Ga$. Thus, \eqref{eq:Jac0} implies that for $w=\ga^{-1}(0)$ we have 
$$
2(1-||w||) \le J_g(\zeta)= \frac{1-||w||^2}{||\zeta-w||^2}\leqslant \frac{2}{C}(1-||w||). 
$$
This implies the equality of (2) and (3) for $z=\zeta\in \Omega$. 
\end{proof}

\begin{cor}
For $s> s_{0,\Ga}$, $P_z(s)$ is a continuous function of $z\in \mathbb B^n \cup \Om$.  
\end{cor} 

\begin{defi}
A Kleinian group $\Ga$ is said to be of {\em divergence type} if $P(s)$ (equivalently, $\cP(s)$) diverges at $s=s_{0,\Ga}$; $\Ga$ is said to be of {\em convergence type} otherwise.  
\end{defi}

\subsection{Conformal measures and densities}\label{sec:Conformal measures}

Perhaps the most natural way to formulate the notion of a {\em conformal density} is to start with the setting of measures on Riemannian manifolds. Let $M$ be a manifold equipped with a Riemannian metric. Then for every diffeomorphism 
$\gamma: M\to M$ we define a smooth function 
$$
J_\ga: M\to [0,\infty), J_\ga(x)=||d\ga_x||,
$$
where $d\ga_x$ is the derivative of $\ga$ at $x$, $d\ga_x: T_xM\to T_{\ga(x)}M$ and $||\cdot||$ is the operator norm defined via the norms on the tangent spaces $T_xM, T_{\ga(x)}M$. 
\begin{defi}\label{def:conformal_density_1}
A measure $\mu$ on $M$ is said to be {\em conformal of dimension $s$} with respect to a group of diffeomorphisms $\Ga< Diff(M)$ if 
\begin{equation}\label{eq:conformal1}
\gamma^*(\mu)(A)=\mu(\gamma(A))= \int_A J^s_\gamma(x) d\mu(x) 
\end{equation}
for every measurable subset $A\subset M$ and $\gamma\in\Gamma$.  In other words,
$$
\frac{d\ga^*\mu}{d\mu}= J^s_\gamma, \gamma\in \Gamma. 
$$
\end{defi}

In terms of functions $f$ on $M$ one can formulate the conformality equation as
$$
\int_M f\circ \ga^{-1}(x)d\mu(x)= \int_M f(x)J^s_\ga(x)d\mu(x). 
$$  

In order to make sense of the equation \eqref{eq:conformal1}, suppose that $\ga$ is a  diffeomorphism between two open subsets of $\mathbb R^n$ (equipped with the Euclidean metric). Then for two measurable subsets $A, B$ of the domain and codomain of $\ga$, $B=\gamma(A)$, the change of variables formula for the volume reads:
$$
Vol(B)=\int_A |\det(D_\ga(x))|dx_1....dx_n,
$$
where $D_\ga(x)$ is the Jacobian matrix of $\ga$ at $x\in A$. Then  we have 
$$
Vol(B)=\int_A J^s_\ga(x)dx_1....dx_n 
$$
if and only if $s=n$ and $D_\ga(x)$ is a conformal matrix a.e. on $A$, i.e. $D_\ga(x)$ is the product of an orthogonal matrix and a scalar. If $A$ is open, this is equivalent to conformality of $\ga$ on $A$. Similarly, for a general Riemannian manifold $M$ and the Riemannian measure $\mu$ on $M$, the measure $\mu$ is conformal with respect to a diffeomorphism $\ga: M\to M$ if and only if $\ga$ is a conformal map (i.e. preserves angles between nonzero tangent vectors).

The example we are actually interested in is that of the unit sphere $S^{n-1}$ with the standard spherical Riemannian metric $g_0$ and a group $\Ga$ of Moebius transformations of $S^{n-1}$. We will regard $S^{n-1}$ as the boundary sphere of the unit ball (centered at $0$) equipped with hyperbolic metric, so that $\Ga$ acts isometrically. Then for every $z\in S^{n-1}$ and a Moebius transformation $\ga$ we have (see \eqref{eq:J_g})
\begin{equation}\label{eq:J_g1}
J_\ga(z)= e^{b_{\ga z}(\ga(0))}=e^{-b_z(\ga^{-1}(0))}, 
\end{equation}
where we normalize Busemann functions to vanish at $0$. 

More generally, since $\bH^n$ is homogeneous, instead of the standard spherical Riemannian metric $g_0$ we could use the Riemannian metric 
$g_x$ on $S^{n-1}$  (with $x\in \bH^n$), defined as follows. Let $\ga$ be an isometry of $\bH^n$ sending $0$ to $x$. Then take the metric 
$$
g_x:= (d\ga^{-1})^* g_0,
$$
the pull-back of $g_0$ via $\ga^{-1}$. Another way to think of $g_x$ is to consider the exponential map $\exp_x: T_x\bH^n\to \bH^n$. 
This map defines a diffeomorphism $\log_x: S^{n-1}\to T^1_x\bH^n$ between the boundary sphere of $\bH^n$ and the unit sphere $T^1_x\bH^n$ in the tangent space $T_x\bH^n$. The latter is equipped with the inner product (given by the hyperbolic metric). Restricting this flat Riemannian metric on $T_x\bH^n$ to the unit sphere $T^1_x\bH^n$ we obtain a Riemannian metric on the latter. Now, pull-back this Riemannian metric to $S^{n-1}$ via $\log_x$. The result is the Riemannian metric $g_x$.

Consider an arbitrary isometry $\ga: \bH^n\to \bH^n$, $\ga(x)=y\in \bH^n$. The norm of the derivative of $\ga$ at a point 
$z\in S^{n-1}$ with respect to the Riemannian metric $g_y$ equals
$e^{-b_z(x,y)}$.  

We can now define the notion of an $s$-conformal density on $S^{n-1}$ with respect to a group of isometries $\Ga$ of $\bH^n$:
It is a family of probability measures $\mu_x$ on $S^{n-1}$ 
indexed by points $x\in \bH^n$ (i.e. a map from $\bH^n$ to the space of probability measures $P(S^{n-1})$) 
satisfying two conditions:

\begin{itemize}
\item Equivariance: $\ga^*\mu_x=\mu_{\ga x}$ for all $x\in \bH^n$ and $\ga\in \Ga$. 

\item Conformality:
$$
\left. \frac{d\mu_y}{d\mu_x}\right\vert_z= e^{-s b_z(y,x)}.
$$ 
\end{itemize}

\begin{example}
Take $\Ga=\Isom(\bH^n)$, the full isometry group. For each $x\in \bH^n$ let $\mu_x$ denote the unique  probability measure on $S^{n-1}$ invariant under the stabilizer of $x$ in the isometry group of $\bH^n$. For instance, in the case $x=0$ we get the standard rotationally invariant probability measure on $S^{n-1}$. Then the family of measures $\mu_x$ is an $(n-1)$-conformal density on $S^{n-1}$ with respect to $\Ga$. 
\end{example}

\medskip 


\subsection{Patterson--Sullivan measures on limit sets of Kleinian groups}

We now return to the setting of Kleinian groups acting on the unit ball $\bB^n$ (we defer until 
Section \ref{sec:reviewPS} the discussion of Patterson--Sullivan measures for general Gromov-hyperbolic spaces).  

We will define two versions of Patterson--Sullivan measures on $S^{n-1}$. As before, we let $\Ga$ denote a discrete subgroup of the isometry group of $\bH^n$ with the limit set $\La$ and discontinuity domain $\Omega\subset S^{n-1}$. 

The first version of Patterson--Sullivan measures is {\em Jacobian}, our treatment follows \cite{FMS}. 

Take $z\in \bH^n\cup \Omega$, $s> s_0=s_{0,\Gamma}$ and define probability measures on $\ol{\mathbb B^n}$:
$$
\nu_{s,z}:=\frac{1}{P_s(z)}\sum_{\ga\in \Ga} J^s_\gamma(z) \delta_{\gamma z}. 
$$
These measures depend (weakly) continuously on $z\in \bH^n\cup \Omega$.  The primary interest of \cite{FMS} is 
in the case of groups $\Gamma$ of convergence type, while we are primarily interested in groups of divergence type (specifically, convex-cocompact subgroups). 
For $\Ga$ of divergence type one then takes weak limits (with the caveat of taking a suitable 
sequence $s_k\to s_0+$), 
$$
\nu_{s_0,z}=\nu_z:=\lim_{s\to s_0+} \nu_{s,z}. 
$$

This limiting measure is the {\em Jacobian Patterson--Sullivan measure} of $\Ga$. It, a priori, depends on the choice of a point $z$ (actually, it does not) and, more importantly, on the fact that we are using the standard Euclidean metric to define $J_\ga$. 

\medskip
The other version of this construction is the {\em metric Patterson--Sullivan measure} of $\Ga$. Namely, given points 
$x, y\in \bH^n$ one considers a family of probability measures
$$
\mu_{x,y}(s)=\frac{1}{\cP_{x,y}(s)}\sum_{\ga\in\Ga} e^{-sd(x, \ga y)}\delta_{\ga y}, s> s_0. 
$$
Then one takes a weak limit
$$
\mu_{x,y}:= \lim_{s\to s_0+} \mu_{x,y}(s),
$$
again with the caveat of first passing to a subsequence $(s_k)$. As we will see in Section \ref{sec:uniqueness} in greater generality, for Kleinian groups different  subsequences yield the same limit, justifying the notation $\mu_{x,y}$. Furthermore, the limit does not depend on the choice 
of the point $y$ (but does depend on $x$). 

Since $\Ga$ is of the divergence type, both measures $\mu_{x,y}$ and $\nu_z$ are supported in the limit set $\La$ of $\Ga$.

It is proven in \cite{FMS} that the measures $\nu=\nu_z$  satisfy the $s$-conformality property for $s=s_0$:

\begin{equation}\label{eq:conformal density}
\nu(\gamma(A))= \int_A J^s_\gamma(x) d\nu(x) 
\end{equation}
for every Borel subset $A\subset \ol{\mathbb B^n}$ and $\ga\in \Ga$. Equivalently, since 
$$
J^s_\gamma(\eta)=e^{-sb_\eta(0,\ga(0))},
$$
the conformality condition reads as 
\begin{align*}
(\ga^{-1})_*(\nu)(A)=\gamma^*(\nu)(A)=\nu(\gamma(A))= \int_A J^s_\gamma(\eta) d\nu(\eta)= \\
\int_A e^{sb_\eta(0,\gamma^{-1}(0))} d\nu(\eta).  
\end{align*}

We include a proof of conformality for the sake of completeness:

\begin{lem}\label{lem:J-conformal}
For every function $f$ continuous on $\ol{\mathbb B^n}$, every $z\in  \bH^n\cup \Omega$, $s\geqslant s_0$, and $\gamma\in \Ga$, the measure $\nu=\nu_{s,z}$ satisfies 
$$
\int f(\ga^{-1}(w))d\nu(w)= \int f(w)J^s_\ga(w)d\nu(w). 
$$  
\end{lem}
\begin{proof} By continuity, it suffices to prove the claim for $s> s_0$. We have
$$
\int f(\ga^{-1}(w))d\nu(w)= \frac{1}{P_s(z)}\sum_{g\in \Ga} f(\ga^{-1}g(z))J^s_g(z)=
$$
$$
 \frac{1}{P_z(s)}\sum_{g\in \Ga} f(g(z))J^s_{\ga g}(z)=   
\frac{1}{P_z(s)}\sum_{g\in \Ga} f(g(z)) J^s_\ga(g(z)) J^s_g(z)=
$$
$$
\int f(w)J^s_\ga(w)d\nu(w). 
$$
\end{proof}

In Section \ref{sec:uniqueness} we will verify conformality of the measures $\mu_{x,y}$ as well (in the more general setting of strongly hyperbolic metric spaces). It turns out that $s_0$-conformal measures in $\La$ for groups of divergence type are unique (if the base-point $x$ is fixed), see 
\cite[Theorem 4.2.1]{nicholls} and also Section \ref{sec:uniqueness}. Therefore, we conclude: 

\begin{thm}\label{thm:alt-PS}
Suppose that a Kleinian group $\Ga< \Isom(\bH^n)$ is of the divergence type. Then for every $z\in \bH^n\cup \Om$, the measure $\nu_z$ is nothing but the Patterson-Sullivan measure $\mu_{0,0}$ of the group $\Ga$. In particular, the measure $\nu_z$  does not depend on the choice of the point $z$. 
\end{thm}


\section{Patterson--Sullivan measures for Gromov-hyperbolic spaces}\label{sec:reviewPS}

Let $(X,d)$ be a Gromov hyperbolic metric space and let $\Gamma$ be a discrete group of 
isometries of $X$ acting metrically properly on  $X$. While these assumptions suffice for all definitions and many results, further assumptions on $(X,d)$ will be needed for a satisfactory theory. Coornaert in \cite{Coo93} assumed that $X$ is 
proper and geodesic. The properness assumption was significantly weakened in \cite{DSU}. For our purposes, it will suffice to assume that $(X,d)$ is {\em semiproper} (see Definition \ref{def:semiproper}): This is much stronger than the assumption in \cite{DSU} but weaker than the assumptions in \cite{Coo93}.

In this section we discuss definitions and properties of the 
Patterson--Sullivan measures on $\geo X$ associated with the action of $\Ga$. We refer the reader to \cite{Coo93, Calegari, BF, DSU} for details. 

\subsection{Poincar\'e series and critical exponent}

For $x\in X$ we have 
the {\em Poincar\'e series} defined analogously to the one in Section \ref{sec:Poincare1}: \index{Poincar\'e series $\cP(s)$}
\begin{equation}\label{eq:Poincare}
\cP(s) = \cP_{x,y}(s)=\sum_{\gamma \in \Gamma} e^{-s d(x, \gamma y)}.
\end{equation}
Let $s_0=s_{0,\Ga}\in \mathbb R_+$ be the {\em critical exponent} \index{critical exponent} 
of $\cP$, i.e. 
$$
s_0=\inf \{s: \cP(s)<\infty\}. 
$$ 
A discrete subgroup $\Ga$ is said to be of {\em divergence type} if $s_0<\infty$ and its Poincar\'e series diverges at $s=s_0$. \index{group of divergence type} 

\begin{thm}
Suppose that $X$ is hyperbolic, geodesic and satisfies the visibility property.  

1. Every nonelementary quasiconvex-cocompact group $\Ga$ has divergence type, see \cite{Coo93}.

2. More generally, every  convex-cobounded group $\Ga$ has divergence type, see \cite[Corollary 17.14]{DSU}.
\end{thm}

Triangle inequalities imply that the critical exponent $s_0$ does not depend on the choice of $x$ and $y$. Indeed, for $z\in X$,
$$
|d(x, \ga y) - d(z, \ga y)|\leqslant D=d(x,z),
$$
$$
e^{-sD}\leqslant \frac{ e^{-s d(x, \ga z)}}{ e^{-s d(y, \ga z)}} \leqslant e^{sD}. 
$$
Hence, by the comparison test, 
$$
\cP_{x,y}(s)<\infty \iff \cP_{z,y}(s)<\infty. 
$$
Since
$$
\cP_{x,y}(s)= \cP_{y,x}(s),
$$
it follows that 
$$
\cP_{x,y}(s)<\infty \iff \cP_{x,z}(s)<\infty. 
$$

The geometric meaning of the critical exponent is that 
it determines the exponent of the exponential growth of $\Ga$-orbits in $X$:
$$
s_0= \lim_{R\to \infty} \frac{1}{R} \log \# (\Ga x \cap B(x,R)), 
$$
see \cite{Coo93}. When we want to stress dependence of the critical exponent on $X$ and the subgroup $\Ga< \Isom(X)$, we will use the notation $s_{X,\Ga}$ for $s_0$. 

\subsection{Partial Poincar\'e series} 

 Suppose that $X$ is proper, $\Ga< \Isom(X)$ is a nonelementary discrete subgroup and that $U\subset \ol{X}$ is  a neighborhood 
 (in $\ol{X}$) of a limit point $\la$ of $\Ga$. Define the {\em partial Poincar\'e series} 
\begin{equation}\label{eq:pps}
\cP_{U}(s)=\cP_{x,U}(s):= \sum_{\ga\in \Ga, \ga(x)\in U} e^{-sd(x, \ga x)}. 
\end{equation}
Note that $\cP_{x,U}(s)$ converges if and only if $\cP_{y,U}(s)$ converges for another $y\in X$. 

\begin{lem}\label{lem:cover}
There exists a finite subset $\Phi\subset \Ga$ such that $\ol{X}= \bigcup_{g\in \Phi} g(U)$. 
\end{lem}
\begin{proof}
Since $\Ga$ is nonelementary, for every $z\in \ol{X}$ the accumulation set of the orbit $\Ga z$ contains the point $\la$. Thus, there exists $\ga\in \Ga$ such that $\ga(z)\in U$, equivalently, $z\in \ga^{-1}U$. In other words, the collection
$$
\{\ga U: \ga\in \Ga\}
$$
is an open cover of $\ol{X}$. By compactness of the latter, this cover contains a finite subcover 
$$
\{\ga U: \ga\in \Phi\}
$$
where $\Phi$ is a finite subset of $\Ga$. \end{proof}

\begin{lem}\label{lem:partialPS}
If $\cP(s)$ diverges at the critical exponent $s_0$, so does $\cP_{x,U}(s)$. 
\end{lem}
\begin{proof} According to Lemma \ref{lem:cover}, there exists a finite subset $\Phi\subset \Ga$ such that
$$\ol{X}= \bigcup_{g\in \Phi} g(U). $$
Suppose that $\cP_{x,U}(s_0)<\infty$. Then the same is true for the partial Poincar\'e series $\cP_{x,gU}(s_0), g\in \Phi$ since 
$$\cP_{x,gU}(s)= \cP_{y, gU}(s)= \sum_{\ga\in \Ga, \ga(x)\in U} e^{-sd(y, \ga y)}$$
where $y=g^{-1}(x)$. Then  
$$\cP(s)\leqslant  \sum_{g\in \Phi} \cP_{x,gU}(s)<\infty. $$
A contradiction. \end{proof}

\subsection{Conformal and quasiconformal densities}


In this section we describe the notions of  $s$-conformal and quasiconformal 
densities on metric spaces suitable for working with visual boundaries of hyperbolic groups and more general than the ones defined in \S \ref{sec:Conformal measures} for real-hyperbolic spaces. 
We will be using definitions and notation for distortion, conformality and quasiconformality of metrics and maps introduced in \S \ref{sec:Dilatation of metrics and maps}. 

Let $Z$ denote a locally compact metrizable topological space, 
$Met(Z)$ the space of all metrics metrizing $Z$. We let $Met_c(Z)$ denote the space of all metrics on $Z$ in a fixed conformal class $c$. Let $\cM_a(Z)$ denote the space of probability measures on $Z$ in a fixed measure class $a$. Then for $0<\al<\infty$, an $\al$-conformal density on $Z$ is a partially defined map
$$
\mu: Met_c(Z)\to \cM_a(Z)
$$
sending a metric $\rho$ to a Borel probability measure $\mu_\rho$ such that for $a$-almost every $z\in Z$ we have
$$
\left. \frac{d\mu_{\rho'}}{d\mu_{\rho}}\right\vert_z= \left. \left(\frac{d\rho'}{d\rho} \right)^\al \right\vert_z. 
$$

This definition has a natural {\em quasiconformal} generalization. We let $Met_q(Z)$ denote the space of all metrics on $Z$ in a fixed Lipschitz class $q$. 
Then for $0<\al<\infty$, an $\al$-quasiconformal density on $Z$ is a partially defined map
$$
\mu: Met_q(Z)\to \cM_a(Z)
$$
sending a metric $\rho$ to a Borel probability measure $\mu_\rho$ such that for $a$-almost every $z\in Z$ we have
$$
K^{-1}  \left(Lip(\rho',\rho)(z) \right)^\al \le \left. \frac{d\mu_{\rho'}}{d\mu_{\rho}}\right\vert_z\le K  \left(Lip(\rho',\rho)(z) \right)^\al, 
$$
where $K\in [1,\infty)$ is a constant independent of $\rho, \rho'$. 

\medskip 
We now proceed to the setting of proper $\delta$-Gromov-hyperbolic spaces $X$ and isometric group actions on $X$. 





\medskip
We are now ready to extend the definition of (quasi)conformal densities to general Gromov-hyperbolic spaces. We start with the easier case of $\eps$-strongly hyperbolic 
spaces $X$. Then we obtain visual metrics $\rho_p=d_{\eps,p}$ on $\geo X$ as  in Theorem \ref{thm:strongly}. These metrics belong to the same conformal class of metrics 
$Met_c(\geo X)$ (see Lemma \ref{lem:conformality of visual metrics}). Then each map $X\to \cM(\geo X)$, $x\mapsto \mu_x$, defines a partially defined map 
$Met_c(\geo X)\to \cM(\geo X)$, $\rho_x\mapsto \mu_x$. This will be our setup for conformal densities.   

\begin{defi}\index{conformal density}\label{def:conformal density}
Let $X$ be an $\eps$-strongly hyperbolic space and $\Ga$ be a group of isometries of $X$. Then a map $X\to \cM(\geo X), x\to \mu_x$, from $X$ to the space of probability measures on $\geo X$ is said to be an {\em $s$-conformal density} with respect to $\Ga$ if the following properties hold:
\begin{itemize}
\item Equivariance: $\ga^*\mu_x=\mu_{\ga x}$ for all $x\in X$ and $\ga\in \Ga$. 

\item $s$-conformality:
$$
\left. \frac{d\mu_y}{d\mu_x}\right\vert_z= e^{-s b_z(y,x)}, x, y\in X, z\in \geo X.
$$ 
\end{itemize}

The conformality condition for points $x, y$ in the same $\Ga$-orbit can be restated as: For every $\ga\in \Ga$, $x\in X$ and a measurable subset $A\subset \geo X$,
$$
\mu_x(\ga A)=\int_A J^{s/\eps}_\ga(z)d\mu_x(z),  
$$   
where the metric Jacobian $J_\ga$ is defined with respect to the metric $\rho_{\eps,x}$ on $\geo X$. 
Equivalently, for every measurable  function $f$ on $\geo X$ we have 
 $$
 \int_{\geo X} f(\ga^{-1}z)d\mu_x(z)= \int_{\geo X} J_\ga^{s/\eps}(z)f(z)d\mu_x(z). 
 $$
 We will also say that {\em the measure $\mu_x$ is $s$-conformal with respect to the action of $\Ga$}. When the base-point $x\in X$ is fixed (say, the center of the unit ball 
 in the ball model of $\bH^n$) then we will say that a measure $\mu$ on $\geo X$ is {\em $s$-conformal with respect to $\Ga$} if for every $\ga\in \Ga$, and a.e. 
$ z\in \geo X$, 
 $$
 \left. \frac{d\ga^*\mu}{d\mu}\right\vert_z= e^{-sb_z(\ga^{-1} x)},  
 $$
 where Busemann functions are normalized to vanish at $x$. 
\end{defi}

Note that all $CAT(-1)$ spaces are $1$-strongly hyperbolic and, hence, one can take $\eps=1$. For instance, in the case of the hyperbolic space $\bH^n$ in the unit ball model, 
 the metric Jacobian with respect to the visual metric $\rho_{1,p}$ is the same as the Jacobian with respect to the Riemannian metric $g_x$ on $S^{n-1}$ discussed earlier.  
 Thus, while the visual metrics $\rho_{1,x}$ on $S^{n-1}$ are definitely non-Riemannian, their metric distortions  $Lip(\rho_x,\rho_y)$ are the same as the conformal factors for the pairs of Riemannian metrics $g_x, g_y$, namely,  $e^{-b_z(y,x)}$.

For general $\delta$-hyperbolic spaces, the visual metrics $d_{p,\al}$ (with the same parameter $\al$) 
are only bilipschitz to each other and one has to settle for a weaker version of  
{\em quasiconformal densities}, defined by Coornaert in \cite{Coo93}:   

\begin{defi}\index{quasiconformal density}
Let $X$ be a $\delta$-hyperbolic space and $\Ga$ be a group of isometries of $X$. Then a map $X\to \cM(\geo X), x\to \mu_x$, 
from $X$ to the space of probability measures on $\geo X$, is said to be a {\em $(K,s)$-quasiconformal density} if the following properties hold:
\begin{itemize}
\item Equivariance: $\ga^*\mu_x=\mu_{\ga x}$ for all $x\in X$ and $\ga\in \Ga$. 

\item $(K,s)$-quasiconformality:
$$
K^{-1} e^{-s |b_z(y,x)|} \le \left. \frac{d\mu_y}{d\mu_x}\right\vert_z\le K e^{s |b_z(y,x)|}, x, y\in X, z\in \geo X.
$$ 
\end{itemize}
 We will also say that each measure $\mu_x$ is $(K,s)$-quasiconformal with respect to the action of $\Ga$. 
\end{defi}

Since metric distortions of pairs of visual metrics on $\geo X$ are comparable to  $e^{-b_z(y,x)}$ (see \S \ref{sec:Conformality properties}), this definition is consistent 
with the general notion of a quasiconformal density defined above.


\begin{lem}\label{lem:generators}
Suppose that $\Ga$ is a transitive group of isometries of an $\eps$-strongly hyperbolic space $X$, 
generated by elements $\ga_i$, $i\in I$. If for every $\ga\in \{\ga^{\pm 1}_i: i\in I\}$,
$$
\mu(\ga A)=\int_A J^{s/\eps}_\ga(z)d\mu(z),
$$  
then the measure $\mu$ is $s$-conformal with respect to the group $\Ga$. 
\end{lem}
\begin{proof} Suppose that conformality equation holds for elements $\al, \beta\in \Ga$. Take a measurable subset $A\subset \geo X$ and let $B=\al(A)$. Then 
$$
\mu(\beta \circ \al(A))= \mu(\beta A)= \int_B J_{\beta}^{s/\eps}(w)d\mu(w). 
$$  
Define the function $f(z)=J^{s/\eps}_\beta(\al z)$. Then the integral becomes
$$
\int_B f(\al^{-1} w)d\mu(w). 
$$
By $s$-conformality of $\mu$ with respect to $\al$, we have
$$
\int_B f(\al^{-1} w)d\mu(w)= \int_{A} f(z)J^{s/\eps}_\al(z)d\mu(z)= \int_{A} J^{s/\eps}_{\beta}(\al z) J^{s/\eps}_\al(z)d\mu(z)= \int_A J^{s/\eps}_{\beta\al}(z)d\mu(z),
$$
by the Chain Rule for the Jacobian (see \eqref{eq:ChainRule}). Hence,
$$
\mu(\beta \al(A))= \int_A J^{s/\eps}_{\beta\al}(z)d\mu(z)
$$
and we obtain $s$-conformality of $\mu$ with respect to the composition $\beta\circ \al$. 
\end{proof}

\subsection{Patterson--Sullivan construction}

For every $s>s_0$ we define a family of probability measures on $X$,
 \begin{equation}
\mu_{x,y}(s) = \frac{1}{\cP(s)}\sum_{\gamma\in \Gamma}e^{-s d(x, \gamma y)}\delta_{\gamma y}, \quad s> s_0, x, y\in X.
 \end{equation}
This family of measures is $\Ga$-invariant in the sense that for every $\ga\in \Ga$ we have
$$
\ga^* \mu_{x,y}(s)= \mu_{\ga x, \ga y}(s). 
$$
Assume that the group $\Ga$ is of divergence type. 

Then Patterson--Sullivan measures are obtained as weak limits of converging sequences of measures 
$\mu_{x,y}(s_n)$, $s_n>s_0$, $\lim s_n=s_0$. In general, there is no reason to expect uniqueness of the limits. We will discuss 
the question of uniqueness in Section \ref{sec:uniqueness}. For now, we only note that all limiting measures 
belong to the same measure class provided that the $\Ga$-action is quasiconvex-cocompact, see \cite{Coo93}. By abuse of the terminology, 
the limiting measures will be denoted 
$\mu_{X,PS}$, $\mu_{x,y}(s_0)$ or simply $\mu_x$. \index{Patterson--Sullivan measure $\mu_{PS}$} Note that the points $x$ and $y$ 
are frequently taken to be the same. 
Sequences $(s_n)$ converging to $s_0$ from the right for which the limiting measure exists will be called {\em PS sequences}. \index{PS sequence}
Each family of Patterson--Sullivan measures $\mu_x$ is still $\Ga$-invariant. 
Furthermore, each measure $\mu_x$ is supported on the limit set $\La(\Ga)\subset \geo X$ of the group $\Ga$.  

\begin{thm}\label{thm:ergodic}
Suppose that $\Ga$ is a properly discontinuous nonelementary  group of isometries of a proper geodesic hyperbolic metric space $X$ of divergence type\footnote{e.g. a convex-cocompact group}. Let $\mu=\mu_{X,PS}$. Then:

1. The measure $\mu$ has no atoms. 

2. The $\Ga$-action on $(\La(\Ga), \mu)$ is ergodic. 

3. The diagonal $\Ga$-action  on 
$$(\La(\Ga)\times \La(\Ga), \mu\times \mu)$$ 
is ergodic. 

4. The nonconical limit set $\La^{nc}(\Ga)$ has zero measure.  
\end{thm}

The first two statement were proven by Coornaert in \cite{Coo93}. Part 3 (which is a strengthening of Part 2) was proven by Bader and Furman in \cite{BF} (in the same paper they proved a number of other interesting ergodic results for actions of isometry groups of hyperbolic spaces on their limit sets). Part 4 is proven by Coulon, \cite{COULON2024}, and Yang, \cite{yang2023conformaldynamicsinfinitygroups},  in greater generality, although in the context of Gromov-hyperbolic spaces the result might have been known earlier. 

The following was proven by Coornaert in \cite[Corollary 7.5]{Coo93} (see also Calegari's paper \cite[Corollary 2.5.10]{Calegari}):

\begin{thm}\label{thm:Hausdorff-qc}
Suppose that $\Ga< \Isom(X)$ is a nonelementary quasicon\-vex-cocompact subgroup with critical exponent $s=s_\Ga$. Then every $s_\Ga$-quasiconformal density $\{\mu_x: x\in X\}$
on $\La(\Ga)$ is in the same measure class as the $s$-dimensional Hausdorff measure ${\mathcal H}_x^{s}$ on $\La(\Ga)$, where the Hausdorff dimension is computed with respect to any $\al$-visual metric $d_\al$ on $\geo X$ with the base-point $x$. More precisely, 
there exists a constant $C\geqslant 1$ such that 
for every $x\in X$  the measure $\mu_x$ satisfies  
$$
C^{-1}{\mathcal H}_x^{s}(A)\leqslant \mu_x(A)\leqslant C{\mathcal H}_x^{s}(A)
$$
for every Borel subset $A\subset \La(\Ga)$.  
\end{thm}





\section{Uniqueness of Patterson--Sullivan measures} \label{sec:uniqueness}

In this section we address the question of independence of the Patterson--Sullivan density  of a PS sequence $(s_n)$. 
We let $\Ga$ be a group acting isometrically and properly discontinuously on a Gromov-hyperbolic space $X$, let $\mu_{x,y}(s), s>s_0$, 
denote the associated family of measures as in the previous section. Independence on $y$ and 
existence of the limit $\mu_{X,PS}=\lim_{s\to s_0+} \mu_{x}(s)$ (and, hence, uniqueness of a Patterson--Sullivan density) 
is known in the following situations:

\begin{enumerate}

\item $X$ is the real-hyperbolic space $X=\bH^n$ (see \cite{nicholls}) and, more generally  the  universal cover of a compact negatively curved manifold, see \cite{Yue} as well as \cite{quint_note}. The argument in \cite{Yue} and \cite{quint_note} goes through in greater generality, when $X$ is a $CAT(-1)$-space and the $\Ga$-action on $X$ is quasiconvex-cocom\-pact, see \cite[Corollary 1.8]{Roblin}. See also Theorem \ref{thm:uniqueness} below for a more general result. 



\item The group $\Ga$ is a $P$-Anosov subgroup of a semisimple Lie group $G$ with the symmetric space $\cX$. Unless $\cX$ has rank one, the group $\Ga$ does not act cocompactly on $\cX$. Nevertheless, one can equip $\cX$ with a certain $G$-invariant polyhedral Finsler metric $d_F$ (adopted to the parabolic subgroup $P$); the group $\Ga$ acts cocompactly on closed $R$-neighborhoods (taken with respect to $d_F$) of each $\Ga$-orbit in $\cX$.  Taking $R$ sufficiently large, we can assume that this neighborhood $X$ is path-connected. We obtain a $\Ga$-invariant geodesic metric $d$ on $X$ from the metric $d_F$ and then can define Patterson--Sullivan measures on $\geo X$ using this metric. It is proven in \cite{Dey-Kapovich} that the limit in the definition of this measure again exists provided that $R$ is large enough. 


\end{enumerate} 

There are instances for which the Patterson--Sullivan densities are not unique, this happens when the group is not geometrically finite, see e.g \cite[p. 480]{mc_mullen_conformal_I}. 


\begin{question}
Is there a geodesic hyperbolic metric space $X$ with a geometric action $\Ga\acts X$ such that the Patterson--Sullivan density on $\geo X$ is not unique? 
\end{question}

Below, we will prove that every (nonelementary) hyperbolic group has a geometric action on a hyperbolic metric space $X$ 
such that the Patterson--Sullivan density $\mu_{X,PS}$ for this action is unique. 


\begin{prop}\label{prop:conformal}
Suppose that $(X,d)$  is  a {\em boundary continuous proper hyperbolic metric space}, e.g. 
a strongly hyperbolic proper metric space. 
Let $\Ga\acts X$ be a properly discontinuous isometric action of divergence type with the critical exponent $s_0$. 
Suppose that $(s_n)$ is a PS sequence for the $\Ga$-action  on  $(X,d)$. Define measures
$$
\mu_{x,p}:=\mu_{x,p}(s_0)=\lim_{s_n\to s_0} \mu_{x,p}(s_n), \quad x, p\in X. 
$$
Then:

\begin{enumerate}
\item $$
\frac{d\mu_{x,p}}{d\mu_{y,q}} (\zeta) = \exp(-s_0 b_\zeta(x,y)),
$$
for all $x, y, p, q\in X$ and $ \zeta\in \geo X$. 

\item $\mu_{x,p}=\mu_{x,q}$ for all $x, p, q\in X$. 

\item The family of measures $x\mapsto \mu_x:=\mu_{x,x}$ is an $s_0$-conformal density on $\geo X$ with respect to the action of the group $\Ga$. 
\end{enumerate}
\end{prop}
\begin{proof} The proof mostly follows the one in \cite[Th\'eor\`eme 5.4]{Coo93}.  

1. Since the function $\exp(-s_0 b_\zeta(x,y))$ is continuous (as a function 
of $\zeta\in \geo X$), in order to prove that 
$$
\frac{d\mu_{x,p}}{d\mu_{y,q}} (\zeta) = \exp(-s_0 b_\zeta(x,y)), \zeta\in \geo X,
$$
it suffices to show that for a countable neighborhood basis $\{W_i: i\in \bN\}$ of $\zeta$ in $\ol{X}$,
$$
\lim_{i\to\infty} \frac{\mu_{x,p}(W_i)}{\mu_{y,q}(W_i)}=\exp(-s_0 b_\zeta(x,y)), 
$$ 
see e.g \cite[\S 14.13]{Williams}. For instance, one can take $W_i=B(\zeta, 1/i)$, 
metric balls with respect to a visual metric $d_\infty$ on $\geo X$.  
Take some $U=W_i$. Our goal is to estimate the ratio
$$
 R_{i,s}=\frac{\mu_{x,p}(s)(U)}{\mu_{y,q}(s)(U)}
$$ 
for $s>s_0$. We will show that
$$
\lim_{i\to\infty} \frac{R_{i,s}}{\exp(-s_0 b_\zeta(x,y))}= 1
$$
and the convergence is uniform in $s$. Fix $x, y, p, q\in X$ and let $\eps=\eps_i$ be such that for all $z\in X\cap U$,
\begin{equation}\label{eq:Busemann}
|b_z(x,y)- b_\zeta(x,y)|\leqslant \eps. 
\end{equation}
By the continuity of the extension of the cocycle $b_z, z\in X$ to $\geo X$, $\eps_i\to 0$ as $i\to\infty$.  
For $U=W_i$ set $U_p:= U\cap \Ga p, U_q:= U\cap \Ga q$, 
$$
 R_{i,s}=\frac{\mu_{x,p}(s)(U)}{\mu_{y,q}(s)(U)}= \frac{\sum_{z\in U_p} \exp(-s d(x, z)) }{\sum_{z\in U_q} \exp(-s d(y, z))}. 
$$ 
We have
\begin{align*}
\sum_{z\in U_p} \exp(-s d(y, z))= \sum_{z\in U_p} \exp(-s d(y, z))=\\
 \sum_{z\in U_p} \exp(-s (d(y, z)- d(x,z)) -sd(x,z) )= \\
   \sum_{z\in U_p} \exp(-s (d(y, z)- d(x,z)) ) \exp(-sd(x,z) ).
\end{align*}
The same applies to the denominator where we sum over $z\in U_q$. 
In view of \eqref{eq:Busemann},
\begin{align*}
e^{-s\eps} \exp(-sb_\zeta(y, x)) \leqslant \exp(-s (d(y, z)- d(x,z)) ) = \\
\exp(-s b_z(y,x))\leqslant e^{s\eps} \exp(-sb_\zeta(y, x) ).
\end{align*}
Applying this to $z\in U_p$ and $z\in U_q$, we obtain: 
$$
e^{-s\eps_i} \leqslant \frac{R_{i,s}}{\exp(-sb_\zeta(x, y)) }\leqslant e^{s\eps_i}. 
$$
Taking the limit $s\to s_0$ we obtain
$$
e^{-s_0\eps_i} \leqslant \frac{R_{i,s_0}}{\exp(-s_0b_\zeta(x, y)) }\leqslant e^{s_0\eps_i}. 
$$
Lastly, we take the limit as $i\to\infty$ and, since $\eps_i\to 0$, we get
$$
\lim_{i\to\infty}  \frac{R_{i,s_0}}{\exp(-s_0b_\zeta(x, y)) }=1. 
$$
This proves Part 1. 

To prove Part 2, we take $x=y$ and observe that Part 1 implies that
$$
\left. \frac{d\mu_{x,p}}{d\mu_{x,q}} \right\vert_\zeta = 1, \forall \zeta\in \geo X.
$$
Hence, $\mu_{x,p}=\mu_{x,q}$, proving Part 2. 

Part 3 follows from the fact that for every $s> s_0$, 
$$
\ga^* \mu_{x}(s)= \mu_{\ga x}(s)
$$
which yields equivariance of the family $x\mapsto \mu_x$; conformality is proven in Part 1. 
Proposition follows. \end{proof}

\begin{rem}\label{rem:quasiconformal} 
{\kap We note that without the boundary continuity assumption, Patter\-son--Sullivan \index{quasiconformal measure} \index{quasiconformal density}
measures are only {\em quasiconformal} (see Section \ref{sec:Conformal measures}). Coornaert proved in 
 \cite[Th\'eor\`eme 5.4]{Coo93} that there exists a constant $C=C(\delta, s_0)$ such that 
$$
 \exp(-s_0  b_\zeta(x,y)  - C) \leqslant \left.  \frac{d\mu_{x,p}(s_0)}{d\mu_{y,p}(s_0)}\right\vert_\zeta \leqslant  \exp(-s_0  b_\zeta(x,y)  + C), \zeta\in \La(\Ga).
$$ 
Quasiconformality of a measure suffices for the uniqueness of the measure class of 
Patterson--Sullivan measures but does not seem to be enough for the uniqueness of the actual measures. }
\end{rem}

\hrule


\begin{thm}\label{thm:uniqueness} 
1. Suppose that $X$ is a  boundary continuous  hyperbolic metric space with compact visual boundary, $\Ga\acts X$ is a metrically proper cobounded 
action. Then for every 
$x\in X$ 
there exists a limit in the definition of the Patterson--Sullivan measure 
$$
\mu_{x}(s_0)=\lim_{s\to s_0+} \mu_{x}(s).$$

2. Every hyperbolic group $\Ga$ admits an action $\Ga\acts X$ such that (1) holds for all $x$ in a $\Ga$-orbit in $X$. 
\end{thm}
\begin{proof} The results is a combination of several known results and standard arguments, but we could not find the theorem in the literature. 

(1) Take a PS sequence $(s_n)$ converging to the critical exponent $s_0$ of the action $\Ga\acts X$. Fix a point 
$x\in X$ and consider the limiting measures 
$\mu_{x}(s_0)=\lim_{s_n\to s_0+} \mu_{x}(s_n)$. Set $X_0:=\Ga x$; then the measure $\mu_x$ is supported in $\geo X_0=\geo X$, the limit set of $\Ga$ in $\geo X$.  
By Proposition \ref{prop:conformal}, the measure $\mu_x$ is $s_0$-conformal with respect to the $\Ga$-action. 
The action of $\Ga$ on $\geo X_0$ is also ergodic with respect to each $\mu_{x}(s_0)$ (see e.g. \cite{Coo93} or \cite{BF}: boundary continuity is not needed for this part of the argument). On the other hand, there exists a  unique $s_0$-conformal probability measure on $\geo X_0$ with respect to which the $\Ga$-action is ergodic, 
see \cite[Proposition 4.2.1]{nicholls} for $\bH^n$, 
\cite[Proposition 3.3.1]{Yue} (in the Riemannian case) and \cite[Proposition 15.3.6]{DSU} (in general). More precisely, if $\nu$ is a probability measure on 
$\geo X_0$ 
satisfying 
$$
\left. \frac{d(\ga^* \nu)}{d\nu}\right\vert_\zeta= \exp(-s_0  b_{\zeta}(\ga(x),x) ), \ga\in \Ga, 
$$
then $\nu=\mu_{x}(s_0)$.  We include a proof since it is quite short. Set $\mu:= \mu_{x}(s_0)$, 
 and define the probability measure $\sigma=\frac{1}{2}(\mu+\nu)$. (This step is needed since we do not know in general if $\mu, \nu$ 
 are absolutely continuous with respect to each other, although this is true for measures defined via different PS sequences.)  
 Then the measures $\mu, \nu, \sigma$ are all $s_0$-conformal with respect to the action of $\Ga$ 
 and $\sigma \gg \mu, \sigma \gg\nu$. Thus, the Radon-Nykodym derivatives 
 $$
 \frac{d\mu}{d\sigma},  \frac{d\nu}{d\sigma}
 $$ 
 are measurable functions on $\geo X$. By conformality, we obtain 
 $$
\left. \frac{d\mu}{d\sigma}\right\vert_{\ga\zeta}= \left. \frac{d\mu}{d\sigma}\right\vert_{\zeta}, \left. \frac{d\nu}{d\sigma}\right\vert_{\ga\zeta}= \left. \frac{d\nu}{d\sigma}\right\vert_{\zeta} 
$$
for $\sigma$-a.e. $\zeta\in \geo X$ (conformal factors in numerator and denominator cancel each other).   
Hence the measurable functions given by these Radon-Nykodim derivative are $\Ga$-invariant; 
by ergodicity, these functions are a.e. constant. The constant has to be $1$ since $\mu, \nu, \sigma$ 
are probability measures. Thus, $\mu=\sigma=\nu$. 

Therefore, for every $x\in X$, the measure $\mu_{x}(s_0)$ does not depend on the choice of a PS sequence. 

\medskip 
(2) Part 2 follows from taking $d$ to be a strongly hyperbolic metric on $\Ga$ guaranteed by Theorem \ref{thm:strongly hyperbolic}. 
Then set $X_0:=(\Ga, d)$ and $X:=E(X_0)$ as in Theorem \ref{thm:hull}.  
 \end{proof}

\begin{rem}
{\kap Suppose that $\Ga$ is a nonelementary hyperbolic group, $\mu$ is a symmetric measure on $\Ga$ whose support is a finite generating set of $\Ga$. 
We then associate with $m$ the corresponding Green metric $d_G$ on $\Ga$ as in Definition \ref{def:G-metric}. One can then define the Patterson--Sullivan density 
$\mu_\ga$ on $\geo \Ga$ corresponding to the metric $d_G$. It is proven in \cite{BHM} that $\mu_\ga$ is nothing but the harmonic measure $\nu_\ga$ on $\geo \Ga$ associated with the random walk starting at $\ga$ and corresponding to the measure $m$ on $\Ga$: For a measurable subset $A\subset \geo \Ga$, $\nu_\ga(A)$ is the probability for the random walk from $\ga$ to converge to a point in $A$. } 
 \end{rem}



\section{Construction and basic properties of Gibbs measures} \label{section_Gibbs}

\subsection{Poincar\'e series, critical exponents and Gibbs measures}\label{sec:Gibbs}

Given a length function $\beta: \Ga\to [0,\infty)$ on a group $\Ga$, we define the {\em Poincar\'e series}: \index{Poincar\'e series $\cP(s)$}
\begin{equation}
\cP(s) = \sum_{\gamma \in \Gamma} e^{-s \beta(\gamma)}.
\end{equation}

This series clearly diverges when $s=0$ (provided that $\Ga$ is infinite). We will assume that the Poincar\'e series converges when $s$ is sufficiently large (which will be the case in all examples of interest to us, i.e. when $\beta$ is regular). Denote by $s_0$ the critical exponent of $\cP$, i.e. the (unique) real number such that $\cP(s)<\infty$ for all $s> s_0$ and $\cP(s)=\infty$ for all $s<s_0$. \index{critical exponent $s_0$} 
As we noted in Section \ref{sec:length functions}, under the assumption that $\beta$ is a hyperbolic length function on a 
(hyperbolic) group $\Ga$, it  
 comes from a geometric action of $\Ga$ on a semiproper hyperbolic metric space $(X_0,d)$ with $\geo X_0=\geo \Gamma$. Assuming further that $\Ga$ is {\em nonelementary}  (i.e. $\geo X_0=\geo \Ga$ contains more than two points), the group $\Gamma$ has the {\em divergence type} with respect to the length function $\beta$, i.e. $\cP(s_0)=\infty$, see \cite[Corollaire 7.3]{Coo93}. Hence, in view of Lemma \ref{lem:regular}, every nonelementary hyperbolic group $\Ga$ has divergence type with respect to all hyperbolic length functions on $\Ga$. 

\begin{example}\label{ex:Free}
Suppose that $\Ga$ is a free group on $r=q+1$ generators $\al_1,...,\al_r$, and let $X$ be the Cayley tree of $\Ga$ corresponding to this generating set $S$. Let $\beta=|\ga|$ we the word-length function on $\Ga$ corresponding to the word-metric $d_S$. Then  the corresponding   Poincar\'e series equals 
\begin{equation}
\label{eq:Pser}
\sum_{\ga\in\Ga} e^{-s|\ga|}= 1+ \sum_{n=1}^\infty (q^{n}+q^{n-1})e^{-sn}. 
\end{equation}
\end{example}

We now consider another proper Gromov-hyperbolic space $X$ (possibly, the same as $X_0$ above) and a properly discontinuous isometric action $\Ga\acts X$ (which need not be cocompact). 
Pick $p\in \overline{X}= X\cup \geo X$ and write $\delta_p$ for the probability measure on $\overline{X}$ supported at $p$. 
We define a family of probability measures on $\overline{X}$ by the  formula
 \begin{equation}
\mu_{p}(s) = \frac{1}{\cP(s)}\sum_{\gamma\in \Gamma}e^{-s\beta(\gamma)}\delta_{\gamma(p)}, \quad s> s_0. 
 \end{equation}

 By the weak compactness of the space of probability measures, there are sequences $s_n\rightarrow s_0^+$ such that the sequence of measures  $\mu_{p}(s_n)$ converges weakly to some probability measure on $\overline{X}$. Due to the divergence assumption, such a limit measure is supported on the limit set $\La(\Ga)\subset \geo X$ (cf. \cite{Coo93}), and the limit measure depends, a priori, on the choice of the sequence $(s_n)$. We will show in Lemma \ref{lem_Gibbs_independant_basepoint} that the limit measures do not depend on the choice of the point $p$ as long as $p\in X$.  Later on, in \S \ref{sec:CT}, we will prove the same for all $p\in \ol{X}$, provided that the action $\Ga\acts X$ is nonelementary and admits a 
 {\em Cannon--Thurston map}, e.g. is quasiconvex-cocompact. 
 
 \begin{defi}\index{Gibbs measure $\mu$}
\label{defi_Gibbs} 
The {\em Gibbs measure class} is the class of measures on $\La(\Ga)$ obtained as weak limits (by letting $s_n\rightarrow s_0^+$) of the measures  $\mu_p(s)$, for $p\in X$. We denote by $\mu=\mu_{X}$ any measure in this family of measures when the context is clear. 
\end{defi}
 
 Assume for a moment that $\beta(\gamma)=d(x,\gamma(x))$ for some $x\in X$. Any limit measure as above is a quasiconformal 
 measure (see Remark \ref{rem:quasiconformal}). If we assume further that the action of  $\Gamma$ on $X$ is quasiconvex-cocompact, then it is proved in \cite[Th\'eor\`eme 7.7 assertion (3)]{Coo93} that all quasiconformal measures on the boundary are absolutely continuous with respect to each other. Therefore, under this extra assumption, 
 the Gibbs measure class is really a class of equivalent measures. 

\begin{example}
Suppose that $X=X_0$, fix a point $p=x_0\in X_0$ and let 
$\beta(\ga)=d(x_0, \ga x_0)$. Then $\mu$ is nothing but the Patterson--Sullivan measure 
$\mu_{x_0}=\mu_{X_0,PS}$ on $\geo X_0$.  
\end{example}

\begin{rem}
Our definition of a Gibbs measure seems at first sight to be very close to that of the  Patterson--Sullivan measure. 
However, the point $p$ in our definition is unrelated to the choice of the length function $\beta$, while the Patterson--Sullivan measure 
$\mu_{X,PS}$ is determined by  
\[\mu_{x}= \lim_{s\to s_0+} \frac{1}{\cP(s)}\sum_{\gamma\in \Gamma}e^{-s \beta(\gamma)}\delta_{\gamma(x)},\]
provided that $\beta(\gamma)=d(x, \gamma x)$. There is, to our knowledge, no treatment of this type of construction (as in the definition of Gibbs measures) in the literature. 
\end{rem}


The following follows immediately from the corresponding facts about Patterson--Sullivan measures:

\begin{thm}\label{thm:Gibbs list}
Let $\Ga$ be a nonelementary hyperbolic group, $\beta: \Ga\to \R_+$ a hyperbolic length function with the critical exponent $s$ and $\mu=\mu_\beta$ the corresponding Gibbs measure on $\geo \Ga$. Then $\mu$ satisfies the following properties:

\begin{enumerate}
\item The measure $\mu$ has no atoms. 

\item The diagonal $\Ga$-action  on 
$$(\La(\Ga)\times \La(\Ga), \mu\times \mu)$$ 
is ergodic.

\item There exists a visual metric $d_\infty$ on $\geo \Ga$ of Hausdorff dimension $s$ and a constant $C\geqslant 1$ 
such that there exists a constant $C\geqslant 1$ such that the the $s$-dimensional Hausdorff measure ${\mathcal H}^s$ 
of $d_\al$ on $\geo \Ga$ satisfies the inequality 
$$
C^{-1}{\mathcal H}^{s}(A)\leqslant \mu^\al(A)\leqslant C{\mathcal H}^{s}(A)
$$
for every Borel subset $A\subset \geo\Ga$. 

\item The density $\{\ga^* \mu: \ga\in \Ga\}$ on $\geo \Ga$ is quasiconformal. 
\end{enumerate}

\end{thm}

\medskip 
The proof of the next lemma is a preview of some counting arguments appearing in the next section. We again turn to general length functions $\beta$ on the group $\Gamma$ and general properly discontinuous actions $\Ga\acts X$ on proper geodesic hyperbolic spaces.   

\begin{lem} \label{lem_Gibbs_independant_basepoint}
The Gibbs measure $\mu_{X}$ is independent of the choice of $p=x\in X$. More precisely, assuming the existence of a weak limit  
$\lim_{n\to\infty}\mu_x(s_n)$ for some $x\in X$ and a sequence $(s_n)$ converging to $s_0$ from 
the right, we have   
$$\lim_{n\to\infty}\mu_x(s_n)= \lim_{n\to\infty}\mu_y(s_n)$$ 
for all $y\in X$. 
\end{lem}
\begin{proof} We fix a metric $d_\infty$ on $\overline X$ and points $x, y\in X$. Observe that for every diverging sequence $(\gamma_i)$ in $\Gamma$, $d_\infty(\gamma_i x, \gamma_i y)\to 0$ as $i\to\infty$. 
Pick a function $\varphi\in C(\overline{X})$ with the modulus of continuity $\eta(\epsilon)$. Then (in view of properness of the function $\beta$ and proper discontinuity of the action of $\Gamma$ on $X$) 
for every $\epsilon>0$, there exists $N=N_\epsilon$ such that 
$
d_\infty( \gamma x, \gamma y)\leqslant \epsilon$ whenever  $\beta(\gamma)\geqslant N$.
Hence, for all $\gamma\in \Gamma$ satisfying $\beta(\gamma)\geqslant N_\epsilon$, 
$
|\varphi(\gamma x)- \varphi(\gamma y)|\leqslant \eta(\epsilon), 
$
Consequently, setting  $\Delta = |\< \varphi, \mu_x(s)\> - \<\varphi,  \mu_y(s)\>|$, we have:
\begin{align*}
\Delta &\leqslant 
\frac{1}{\cP(s)} \sum_{\beta(\gamma)<N} e^{-s\beta(\gamma)} |\varphi(\gamma x)- \varphi(\gamma y)| + \frac{1}{\cP(s)} \sum_{\beta(\gamma)\geqslant N} e^{-s\beta(\gamma)} \eta(\epsilon) \\
&\leqslant 
 \frac{C(N)}{\cP(s)} + \eta(\epsilon),
\end{align*} 
where $C(N)$ depends only on $N=N_\epsilon$, on $x, y$ and on $\Gamma$. Moreover, since $\Gamma$ is of the divergence type with respect to $\beta$,  
$$
\lim_{s\to s_0+}    \frac{C(N)}{\cP(s)}=0.
$$
 At the same time, $\lim_{\epsilon\to 0+} \eta(\epsilon)=0$. It follows that  
$$
\lim_{s\to s_0+}   |\< \varphi, \mu_x(s)\> - \<\varphi,  \mu_y(s)\>|=0
$$
concluding the proof of the lemma. 
\end{proof}

\begin{rem}

This {\em universality} feature of Gibbs measures will manifest itself many times in the book in the form of various equidistribution results, where the same Gibbs measure will appear for points, subvarieties and currents. 
\end{rem}

We will eventually generalize this lemma to the case of points $x\in \geo X$, when the action $\Ga\acts X$ is quasiconvex-cocompact and $\beta$ is regular; see Theorem \ref{thm:rank1-equi}. In particular, in this situation, the existence of a (weak) limit of measures $\lim_{n\to\infty}\mu_x(s_n)$ is independent of a space $X$ and a specific action $\Ga\acts X$. 
This justifies the following definition:

\begin{defi} \index{PS-sequence}
Given a fixed length function $\beta$ and the corresponding critical exponent $s_0$, a sequence $s_n>s_0$ converging to $s_0$ is called a {\em PS sequence} if for some (equivalently, every) proper hyperbolic space $X$, a  quasiconvex-cocompact action $\Ga\acts X$ and some (equivalently, every) point $p\in X$,
there exists a weak limit of the sequence of measures  
$\lim_{n\to\infty}\mu_p(s_n)$. 
\end{defi}

For most of this work we will be working with a fixed length function $\beta$ and a fixed PS sequence $(s_n)$. Then the corresponding 
limit measure $\mu=\lim_{n\to\infty}\mu_p(s_n)$ on $\geo X$ will be called {\em the} Gibbs measure of the group $\Gamma$. This 
is justified by the fact that limit sets of  quasiconvex-cocompact actions $\Ga\acts X$ are $\Gamma$-equivariantly homeomorphic to 
$\geo \Ga$. 

\begin{rem}
Theorem \ref{thm:uniqueness} shows that every hyperbolic group $\Ga$ admits a geometric action $\Ga\acts (X,d)$ on a geodesic 
metric space such that for the length function $\beta(\ga)=d(x, \ga x)$, the Gibbs measure $\mu$ is independent of a PS sequence $(s_n)$. 
\end{rem}

\subsection{Pushforward of Gibbs measures via equivariant continuous maps}\label{sec:CT}


We assume now that $\beta$ is regular, coming from a geometric action $\Ga\acts X_0$ on a geodesic hyperbolic metric space $X_0$, 
$\beta(\ga)= d(x_0, \ga x_0)$. Since $\beta(\gamma)>0$ for all $\gamma\in \Gamma \setminus \{1\}$, this means that the point 
$x_0\in X_0$ has trivial stabilizer.  Consider also, as before, a properly discontinuous isometric action $\Ga\acts X$ of $\Ga$ on another proper 
hyperbolic metric space $X$. 
We have a $\Ga$-equivariant map $f: \Gamma x_0\to X\cup \Omega(\Ga)$, sending $x_0\in X_0$ to a point 
$x\in X\cup \Omega(\Ga)$ (here we are using the assumption that $x_0$ has trivial stabilizer). Below, we will be working under the hypothesis that the map $f$ admits a continuous extension 
$$
f: \Gamma x_0 \cup \geo X_0\to \overline{X}. 
$$\index{Cannon--Thurston extension} \index{Cannon--Thurston map}
Then, necessarily, the extension is equivariant and $f(\geo X_0)= \Lambda(\Ga)$. Such extensions are called {\em Cannon--Thurston extensions} 
and their restrictions $\geo X_0\to \La(\Ga)$ are called  {\em Cannon--Thurston maps}. 
Such maps  
 do not always exist, although examples are by no means easy, see references in \cite{Mj}. Cannon--Thurston  extensions exist in the following special cases:

\begin{itemize}
\item The action of $\Ga$ on $X$ is  quasiconvex-cocompact. In this case, $f: \geo X_0\to \Lambda(\Ga)$ is a homeomorphism. 

\item {Relatively hyperbolic groups $(\Ga; P_1,...,P_n)$, where $\Ga$ is a hyperbolic group and $P_i< \Ga$ are quasiconvex subgroups of $\Ga$, 
see Theorem 7.11 in \cite{Bowditch-RH}. The Cannon--Thurston map in this setting is from the Gromov boundary of $\Ga$ to the Bowditch boundary 
of $(\Ga; P_1,...,P_n)$.}

\item The space $X$ is the real-hyperbolic 3-space $\bH^3$. (This is a very hard theorem proven primarily by M.~Mj, see \cite{Mj} for a survey.) 

\item $X$ is the Cayley graph of another hyperbolic group $\Gamma'$ which is obtained as the fundamental group of a graph ${\mathcal G}'$ of hyperbolic groups where edge-groups are quasi-isometrically embedded in the incident vertex-groups. The subgroup $\Gamma$ is the fundamental groups of a subgraph of groups ${\mathcal G}\subset {\mathcal G}'$. (See \cite{Kapovich-Sardar}.) 
\end{itemize}

\noindent Under the above extension assumption, we obtain the push-forward measure $f_*(\mu_{X_0,PS})$, whose support (by equivariance) is the entire  $\Lambda(\Gamma)\subset \geo X$. 

\begin{lem}\label{lem:CT}
Assuming that the orbit map $f: \Ga x_0\to \Ga x$ has a Cannon--Thurston extension, 
for every $x\in X\cup \Omega(\Ga)$, 
$$
\lim_{s_n\to s_0+} \mu_x(s_n)= f_*(\mu_{X_0,PS}).$$ 
\end{lem}
\begin{proof} The claim follows from the observation that the continuous map $f: \Gamma x_0 \cup \geo X_0\to \overline{X}$ induces continuous push-forward map $f_*$ of the spaces of probability measures 
$$
P(\Gamma x_0 \cup \geo X_0)\to P(\overline{X}),$$
where continuity is understood with respect to the weak *-topology. \end{proof}  

\begin{cor}\label{cor:CT}
Under the same assumptions, the measure $\mu_{X}$ equals the push-forward measure 
$f_*(\mu_{X_0,PS})$. In particular, $\mu_{X}$ has no atoms provided that $\Gamma$ is a nonelementary  quasiconvex-cocompact group. 
\end{cor}


Note that under very mild assumptions on $X$, there is always an equivariant measurable map $\geo X_0\to \La(\Ga)\subset \geo X$, see \cite{Monod-Shalom}.

\subsection{Dependence of the Gibbs measure on the length function}

Different hyperbolic length functions on a hyperbolic group $\Ga$ typically give rise to different Gibbs measures. Frequently, these measures are not even absolutely continuous with respect to each other. For instance, consider the fundamental group $\Ga$ of a compact oriented hyperbolic surface $S$. Take two representations 
$\rho_i: \Ga\to \Isom(\bH^2)$ corresponding to different points in the Teichm\"uller space of $S$. They give rise to two hyperbolic length functions $\beta_i: \Ga\to \bR$. 
The corresponding Gibbs measures are just the Patterson--Sullivan measures $\mu_{\Ga_i}$ for the Fuchsian groups $\Ga_i=\rho_i(\Ga), i=1,2$. Let $f: S^1\to S^1$ denote the unique homeomorphism inducing the isomorphism $\phi: \Ga_1\to \Ga_2, \phi=\rho_2\circ \rho_1^{-1}$. Then $f$ is not absolutely continuous and, hence, violates Luzin's Property N. Therefore, $f_*(\mu_{\Ga_1})$ is singular with respect to $\mu_{\Ga_2}$. 

More generally, Furman in \cite[Theorem 2]{Furman},  proves that the following are equivalent for two hyperbolic length 
functions $\beta_1, \beta_2$ on  a hyperbolic group $\Ga$:

\begin{enumerate}
\item The measures $\mu_{\beta_1}, \mu_{\beta_2}$ are not mutually singular with respect to each other. 

\item There exists $a>0$ such that the length functions $\beta_1$ and $a\beta_2$ satisfy
$$
||\beta_1- a\beta_2||_{\infty}<\infty.$$
\end{enumerate}

On the other hand, (in view of Theorem \ref{thm:Hausdorff-qc}), if  $A\subset \geo \Ga$ has zero Hausdorff dimension with respect to a visual metric on $\geo \Ga$, then $\mu_{\beta}(A)=0$ for every hyperbolic length function $\beta$ on $\Ga$. 

\subsection{Conical support set of Gibbs measures}\label{sec:conical}

We continue with the setup of \S \ref{sec:CT}: We have a hyperbolic group $\Ga_0$ and a properly discontinuous isometric action 
$\rho: \Ga_0\acts X$ on a proper geodesic metric space $X$; set $\Ga=\rho(\Ga_0)< \Isom(X)$.  
Assume that $\rho$ admits a Cannon--Thurston map $f: \geo \Ga_0\to \La(\Ga)\subset \geo X$. Accordingly, we have a hyperbolic length function $\beta$ on $\Ga_0$ and the corresponding Gibbs measures $\mu=\mu_{\Ga_0}$ on $\geo \Ga_0$ and $\mu_{\Ga}=f_*(\mu)$ on $\La(\Ga)$. Assume also that the group $\Ga$ is nonelementary.  

\begin{thm} \index{nonconical limit set $\La^{nc}(\Ga)$}
The Gibbs measure $\mu_\Ga=f_*(\mu)$ vanishes on the subset $\La^{nc}(\Ga)\subset \La(\Ga)$ consisting of nonconical limit points. 
\end{thm} 
\begin{proof} 
The interesting case is when the map $f$ is noninjective. The {\em Cannon--Thurston lamination} $\La(\rho)$ 
is the subset of 
$$
\geo^2 \Ga_0:= \geo \Ga_0\times \geo \Ga_0 \setminus \diag(\geo \Ga_0\times \geo \Ga_0)$$ 
consisting of pairs $(\xi,\eta)$ such that $f(\xi)=f(\eta)$ 
(cf. \cite{Mitra}). Since $\Ga$ is assumed to be nonelementary, $\La(\rho)$ is a closed proper $\Ga_0\times \bZ_2$-invariant subset of $\geo^2 \Ga_0$. In particular, $\La(\rho)$ does not have full measure with respect to the product measure $\mu\times \mu$ on $\geo^2 \Ga_0$. Therefore, in view of the double ergodicity theorem (Theorem \ref{thm:ergodic}), $\La(\rho)$ has zero measure. It follows that the projection $A$ of $\La(\rho)$ to $\geo \Ga_0$ has zero measure as well.

\medskip 
Clearly, the restriction of $f$ to $\geo \Ga_0\setminus A$ is 1-1. 
 Next, we discuss points $\zeta\in \geo \Ga_0$ which are mapped via $f$ to nonconical limit points of $\Ga$. 
 It turns out that there are two types of such points. The first type consists of points where $f$ fails to be injective: 
 $f(A)\subset \La^{nc}(\Ga)$ (see \cite[\S 8]{Kapovich1995} and \cite[\S 7]{JKLO}). 
 The second type is harder to describe. We will only describe a {\em necessary condition} 
 for a point to be in $f^{-1}(\La^{nc}(\Ga))\setminus A$ following \cite{JKLO}.

\begin{defi}\index{controlled concentration point}
Let $G$ be a nonelementary hyperbolic group. A point $\la_+\in \geo G$ is  a {\em  controlled concentration point} \index{controlled concentration point}
 if there exists 
$\la_-\in \geo G\setminus \{\la_+\}$ and a sequence of distinct elements $g_i\in G$ such that 
$$
\lim_{i\to\infty} g_i(\la_+, \la_-)= (\la_+, \la_-) 
$$
and, moreover, the sequence $g_i$ converges to $\la_+$ uniformly on compacts in $\geo G\setminus \{\la_-\}$. Pairs of points $(\la_+, \la_-)$ appearing this way 
are {\em concentration pairs} in $\geo^2 G$. \index{concentration pairs}
\end{defi}

Note that if $(\la_+, \la_-)$  is concentration pair, it does not necessarily follow that $(\la_-, \la_+)$  is also a concentration pair. By its definition, the set of concentration pairs is a countable intersection of open subsets, hence, is measurable. 
A pair \index{recurrent point}
$(\la_+, \la_-)\in \geo^2 G$ is a {\em recurrent point} if there exists a sequence of distinct elements $g_i\in G$ such that 
$$
\lim_{i\to\infty} g_i(\la_+, \la_-)= (\la_+, \la_-).  
$$
The convergence property of the action of $G$ on $\geo G$ (see \S \ref{sec:convdy}) 
implies that for every recurrent pair  $(\la_+, \la_-)$ either  $(\la_+, \la_-)$ or 
$(\la_-, \la_+)$  is a concentration pair. 

Since the diagonal action $G\acts \geo^2 G$ is ergodic with respect to the Gibbs measure $\mu$ on $\geo G$ 
corresponding to any hyperbolic length function on $G$ (see Theorem \ref{thm:ergodic}), it follows that 
$\mu$-a.e. pair $(\la_1, \la_2)\in \geo^2 G$ is recurrent under the action of $G$. Applying 
ergodicity again, we see that $\mu$-a.e. pair $(\la_1, \la_2)\in \geo^2 G$ is a concentration pair. 
In particular, $\mu$-a.e. point $\la\in \geo G$ is a controlled concentration point. 

We will need the following proposition proven in \cite[Proposition 7.5]{JKLO}:

\begin{prop}
Suppose that $G$ is a hyperbolic group and $\rho: G\acts X$ is a nonelementary properly discontinuous isometric action on a proper geodesic hyperbolic space for which there is a Cannon-Thurston map $f: \geo G\to \geo X$. Let $A\subset \geo G$ be the subset consisting of points such that $f^{-1} f(\la)\neq\{\la\}$. If $\la_+\in \geo G\setminus A$ is a concentration point, then $\eta_+=f(\la_+)$ is a conical limit point of $\rho(G)$. 
\end{prop}
\begin{proof} Let $\la_-\in \geo G$ and $(g_i)$ be a sequence in $G$ as in the definition of a controlled concentration point, i.e. 
$$
\lim_{i\to\infty} g_i(\la_+, \la_-)= (\la_+, \la_-) 
$$
and $\lim_{i\to\infty} g_i(\xi)=\la_-$ for every $\xi\ne \la_+$. 
Set $\eta_-:= f(\la_-), \eta=f(\xi)$. Since $\la_+\notin A$, $\eta_-\ne \eta_+$. By continuity of the map $f$,   
$$
\lim_{i\to\infty} f(g_i\la_+, g_i\xi) = \lim_{i\to\infty} (g_i\la_+, g_i\eta) = (\eta_+, \eta_-). 
$$
Since $\eta_+\ne \eta_-$, the point $\eta_+$ is a conical limit point of $\rho(G)$. 
\end{proof}

In particular, if $\mu$ is a Gibbs measure on $\geo G$ corresponding to a hyperbolic length function, then 
$\mu(f^{-1}(\La^{nc}(\rho(G))) \setminus A)=0$.   

\medskip 
We can now finish the proof of the theorem. As we noted earlier, 
$\mu(A)=0$. We also just proved that $\mu(f^{-1}(\La^{nc}(\Ga)) \setminus A)=0$. 
Therefore, $\mu_\Ga(\La^{nc}(\Ga))=0$. 
\end{proof}

\begin{question} 
Suppose that $\Ga_0$ is a hyperbolic group, $\rho: \Ga_0\acts X$ is a properly discontinuous nonelementary isometric 
action on a proper geodesic Gromov hyperbolic space $X$ for which there exists a Cannon--Thurston map 
$f: \geo \Ga_0\to \geo X$. Is it true that $\La(\rho)$ has zero Hausdorff dimension 
with respect to any visual metric on $\geo \Ga_0$? 
\end{question}

This is known to be true in the case of surface groups $\Ga_0$ and discrete and faithful representations 
$\rho: \Ga_0\to \Isom(\bH^3)$, cf. \cite{Birman-Series}. As a consequence, it is also true for all convex-cocompact 
subgroups $\Ga_0< \Isom(\bH^3)$ with connected limit sets and discrete and faithful representations $\rho: \Ga_0\to \Isom(\bH^3)$. 

Lastly, we note that Lecuire and Mj proved in \cite{MR3801448} that for every geometrically finite discrete subgroup $\Ga_0< \Isom(\bH^3)$ and 
every isomorphism $\rho: \Ga_0\to \Ga< \Isom(\bH^3)$ to a discrete subgroup $\Ga$, the image of $\La(\rho)$ under the Cannon--Thurston map 
$\La(\Ga_0)\to \La(\Ga)$ is exactly the set of non-horospherical limit points of $\Ga$.

\subsection{Measures on geodesic flows of hyperbolic groups}\label{sec:Measures on geodesic flows of hyperbolic groups}

In this section we consider a non-elementary hyperbolic group $\Gamma$ with a fixed word metric $d_\Ga$ and we will assume that the total space of the geodesic flow $\widehat\Ga$ is naturally homeomorphic to $\geo^2\Ga \times \R$ (see Section \ref{sec:Geodesic flows of hyperbolic groups}). Below, we describe {\em Bowen--Margulis--Sullivan (BMS) measures} on $\geo^2\Ga$ and $\widehat\Ga$, and (flow-invariant) {\em Bowen--Margulis (BM)  
measures} on $\cG$, following a work of Bader and Furman, \cite{BF}.

\begin{rem}
In this setting one cannot assume that the $\Ga$-action on $\widehat\Ga$ is free. If $\Ga$ contains a nontrivial element $\ga$ of finite order with (at least) two distinct fixed points $\xi_\pm$ in $\geo \Ga$, then $\ga$ must have nonempty fixed-point set in $\widehat\Ga$. Indeed, the flow-line $\{\xi_-\}\times \{\xi_+\}\times \R$ is fixed by $\ga$ pointwise. Such examples are easy to find: For instance, a uniform lattice $\Ga< PSL(2,\C)$  which has a nontrivial element of finite order, or one can take 
$\Ga=\bZ * (\bZ_2\times \bZ)$.  
\end{rem}

\begin{defi}
Let $\mu=\mu_{\beta}$ be the Gibbs measure on $\geo\Ga$ with critical exponent $s$, associated with a hyperbolic length function $\beta$ on $\Ga$. 
A $\Ga$-invariant Radon measure $\nu$ in the 
measure class of $\mu\times \mu$ on $\geo^2 \Ga$ is called a {\em BMS measure} if 
$$d\nu(\xi, \eta) = e^{F(\xi,\eta)} d\mu(\xi) d\mu(\eta),$$
where $F(\xi,\eta)$ is a measurable function on $\geo^2\Ga$ such that the difference $|F(\xi,\eta)- 2s(\xi,\eta)_{1_\Ga}|$ is a bounded function on 
$\geo^2\Ga$. \index{$\nu^{BMS}$, BMS measure} 
\end{defi}

\begin{rem}
The factor of $2$ is mistakenly omitted in \cite{BF}. 
\end{rem}

\begin{thm}
[Bader and Furman, \cite{BF}; Tanaka, \cite{Tanaka}] 1. There exists a BMS measure on $\geo^2\Ga$. 

2. Every BMS measure $\nu$ on $\cG$ is ergodic. 
\end{thm}

Given a BMS measure $\nu$ one defines a flow-invariant and $\Ga$-invariant measure $\nu_{\R}$ on $\widehat \Ga=\geo^2\Ga\times \R$ 
as the product of $\nu$ and the Lebesgue measure on $\R$. 

\begin{thm}
[Bader and Furman, \cite{BF}; Tanaka, \cite{Tanaka}]  
1. There exists a flow-invariant measure $\nu^{BM}=q_\bullet(\nu_{\R})$ 
on $\cG$, the projection of  the measure  $\nu_{\R}$ under the 
quotient map $q: \widehat\Ga\to \cG$, see \S \ref{sec:measures}. 


2. Every measure $\nu^{BM}$ as in (1) is ergodic. 
\end{thm}

Every measure $\nu^{BM}$ as in this theorem is called a {\em Bowen--Margulis} (BM) measure on $\cG$. \index{$\nu^{BM}$, BM measure}

\begin{rem}
The projection $\nu^{BM}$ is defined in \cite{BF} using a measurable fundamental domain for the $\Ga$-action on $\widehat \Ga$. Construction of such domains is a bit delicate since the $\Ga$-action is not free. It helps that (assuming that $\Ga$ is nonelementary) 
each $\ga$ either fixes the entire $\widehat \Ga$ (which happens precisely when $\ga\in \Ga$ belongs to the maximal finite normal subgroup in $\Ga$) 
or its fixed-point set has zero measure. The union of proper fixed-point sets is a closed measure zero subset in  $\widehat \Ga$ and removing it results in a space with a 
 covering $\Ga$-action. The projection construction in \S \ref{sec:measures} circumvents all these issues. To see that the definition of the projection measure given in 
 \S \ref{sec:measures} agrees with the one given by measurable fundamental domains, note that every a relatively compact measurable fundamental domain $\Phi$ is 
 precisely-invariant under the $\Ga$-action (and has trivial stabilizer). Hence, the measure of its image in $\cG$ given by the projection construction equals its 
 $\mu_{\R}$-measure. 
\end{rem}


Here is a construction of a BMS measure in the setting of boundary continuous length functions $\beta$ (i.e. boundary continuous metrics) on $\Ga$, see e.g. \cite{Tanaka}. 
Let $(\cdot, \cdot):= (\xi,\eta)_{1_\Ga}$ be the Gromov product on $\geo \Ga$ defined by the continuous extension of the  Gromov product on $\Ga$ itself. 
 Define $\nu$ by 
 $$
 d\nu(\xi,\eta):= e^{2s(\xi,\eta)} d\mu(\xi)d\mu(\eta), \quad \nu=F\mu\times \mu,
 $$
 where $F(\xi,\eta)=\exp(2s(\xi,\eta))$.

 \begin{lem}
 The measure $\nu$ is $\Ga$-invariant. 
 \end{lem}
 \begin{proof} For $\ga\in \Ga$ we have, in view of conformal invariance of the measure $\mu$,
$$
d\ga^*\mu(\zeta)= e^{-sb_\zeta(\ga^{-1},1)} d\mu(\zeta), \zeta\in \geo \Ga. 
$$
Thus,
\begin{align}\label{eq:invariance}
\begin{split}
(d\ga^*\nu)(\xi,\eta)= &F(\ga\xi,\ga\eta)d\ga^*\mu(\xi)\cdot d\ga^*\mu(\eta) \\
& =F(\ga \xi, \ga \eta) e^{-sb_{\xi}(\ga^{-1},1)- sb_{\eta}(\ga^{-1},1)} d\mu(\xi)\cdot d\mu(\eta). 
\end{split}
\end{align}
According to \eqref{eq:Buscocycle}, we have
$$
b_{\xi}(\ga^{-1},1)+ b_\eta(\ga^{-1},1)= 2 ((\xi,\eta)_{\ga^{-1}} - (\xi,\eta))= 2 ((\ga\xi,\ga\eta) - (\xi,\eta)). 
$$
Thus, the right hand side of the equation \eqref{eq:invariance} becomes
$$
F(\ga \xi, \ga \eta) \exp(-2s (\ga\xi,\ga\eta) ) \exp(2s(\xi,\eta))d\mu(\xi)\cdot d\mu(\eta)= F(\xi,\eta)d\mu(\xi)\cdot d\mu(\eta)= d\nu(\xi,\eta). 
$$
Invariance of the measure $\nu$ follows. \end{proof}

\medskip
Given the measure $\nu$ as above, one takes the product measure $\nu_{\R}$ on $\geo^2\Ga \times \R$ of $\nu$ and the Lebesgue measure on $\R$. 
As we noted above, since this measure is clearly both $\R$-invariant and $\Ga$-invariant, it projects to a (flow-invariant) Bowen--Margulis measure 
$\nu^{BM}$ on $\cG$. Ergodicity of the measure $\nu^{BM}$ is equivalent to ergodicity of the measure $\nu$ on $\geo^2\Ga$ with respect to the $\Ga$-action. 
Since $F$ is a measurable function on  $\geo^2\Ga$, ergodicity of $\nu$ is equivalent to the ergodicity of the product measure $\mu\times \mu$ on $\geo^2\Ga$.


\section{Counting}

The goal of this section is to prove some technical results needed for the generalization of the equidistribution lemma (Lemma \ref{lem:CT}) from the case of $\Ga$-orbits  in a hyperbolic space $X$  to orbits in the limit set of $\Ga$. These results are also critical for proving equidistribution theorems for currents on flag-manifolds under Anosov group actions. 
One of the issues that we have to address is that (unlike points in $X$), different points $p$ in the limit set of $\Ga$ behave differently  under the action of (uniformly quasigeodesic) sequences in $\Ga$. This gives rise to partitions of $\Ga$ into {\em good} and {\em bad} subsets with respect to the point $p$ which appear below.

\subsection{Conical counting} \label{section_conical}



\subsubsection{Setup}\label{sec:counting_setting}

In this section, $(X,d)$ is a $\delta$-hyperbolic {\mini proper} geodesic space and $\Ga\acts X$ is a  convex-cobounded action. 
We will fix a base-point $x_0\in X$ which, we will assume, has trivial $\Ga$-stabilizer (so that $\beta(\ga)=d(x_0, \ga(x_0))$ is a hyperbolic length function on $\Ga$). 
Then the group $\Ga$ is necessarily Gromov-hyperbolic and there is an equivariant homeomorphism between the 
Gromov-boundaries of $X$ and $\Ga$. 


\subsubsection{Bad sets}\label{sec:Counting}

Consider a point  $p\in \geo X$. 
Let $\rho=x_0p$ be a geodesic ray in $X$ joining $x_0$ and $p$. For each integer $N$, we define 
$$
\Gamma_{N, bad} = \{ \gamma \in \Gamma | d(x_0, \gamma(\rho)) \leqslant N \}$$   
and \index{$\Gamma_{N, bad}$} \index{$\Gamma_{N, good}$}
$$
\Gamma_{N, good} = \Gamma \setminus \Gamma_{N,bad}.$$ 

\begin{rem}
If $p\in \Omega(\Ga)$ then for every $N$, $\Gamma_{N,bad}$ is finite. In contrast, if $p\in \La(\Ga)$, then  
$\Gamma_{N,bad}$ is infinite for all sufficiently large $N$. This will be the most interesting case. 
\end{rem}

The group $\Ga$ is divided into two disjoint parts: 
\begin{equation}\label{eq:good-bad-partition}
\Gamma =  \Gamma_{N,good} \sqcup \Gamma_{N,bad}. 
\end{equation}

We also define $\Ga_N:=\{\ga\in \Ga: \beta(\ga) \leqslant N\}$.  
Let $A_{n_1}^{n_2}$ be the annulus 
\[
\{ \gamma \in \Gamma | n_1 <  \beta(\gamma) \leqslant n_2 \},
\]
and let $A_n$ denote $A_n^{n+1}$. \index{$A_n$, annulus in a group}

Let $\Sh^\infty_{x_0}(y, r)\subset \geo X$ be the {\em shadow at infinity} from $x_0$ of the ball of radius $r$ centered at 
\index{ $\Sh^\infty_{x_0}(y, r)$, shadow at infinity}
$y \in X$, i.e. it is the set of points 
$\xi\in \geo X$  such that each geodesic $x_0\xi$ joining $x_0$ and $\xi$ intersects the  closed ball 
$\bar{B}(y,r)$. Similarly, one defines the {\em shadow} $\Sh_{x_0}(y, r)\subset X$ as  \index{$\Sh_{x_0}(y, r)$, shadow}
the set of points $z\in \ol{X}$  such that each geodesic $x_0z$ joining $x_0$ and $z$ intersects $\bar{B}(y,r)$. 
Then for $p\in \geo X$, the shadows $\Sh_{x_0}(y,r)$, $y \in x_0p$, $r > 0$, form a neighborhood basis of $p$ in $\ol{X}$. 

We will use an estimate that appears in the proof of  \cite[Theorem 4.11]{quint_note}:  

\begin{prop} \label{prop_conical_counting} For all $N \in \mathbb{N}$ and all $n \geqslant N$ one has:
\begin{equation*}
\# (A_n \cap \Gamma_{N,bad}) \leqslant \# \Gamma_{1+4N}. 
\end{equation*}
Thus, the number of bad points grows linearly with the hyperbolic distance. 
\end{prop}

\begin{proof}


Let us take any two elements $\gamma, \theta \in A_n$ such that
 $$p\in \Sh_{x}( \gamma^{-1}(x_0), N) \cap \Sh_x(\theta^{-1}(x_0) , N) .$$
   Then 
the geodesic ray $\rho=x_0p$ contains a point $y \in \bar{B}(\gamma^{-1}(x_0),N)$. Since $n \leqslant d(x_0, \gamma^{-1}(x_0)) \leqslant n+1$, the triangle inequality yields (assuming $n > N$): 
\begin{equation*}
n - N\leqslant d(x_0,y) \leqslant n+1 +N . 
\end{equation*} 
Similarly, there is a point $z $ on the same geodesic $\rho$, which belongs to the ball $\bar{B}(\theta^{-1}(x_0), N)$ such that 
$n- N \leqslant d(x_0, z) \leqslant n+1 + N$. Since $y, z$ lie on the same geodesic ray $\rho$ emanating from $x_0$, we get 
$d(y,z) \leqslant 1 + 2N$. The triangle inequality yields:
 $$d(x_0, \gamma \theta^{-1}(x_0)) =d(\gamma^{-1}(x_0),\theta^{-1}(x_0)) \leqslant 1 + 2N + 2N= 1 + 4N. $$ 

This shows that for $n >N$ and for any two elements $\gamma, \theta \in A_n \cap \Gamma_{N,bad}$, one has $\gamma\theta^{-1} \in \Ga_{1+4N}$, 
which gives: 
\begin{equation} \label{eq_bad_annuli_count}
\#  \left ( A_n \cap \Gamma_{N,bad} \right ) \leqslant \# \Gamma_{1+4N}.   
\end{equation} 
\end{proof}



\subsection{Shadow counting} \label{section_shadow_counting}

We continue with the setup from the previous section (see \S \ref{sec:counting_setting}). 
We fix a base-point $x_0\in X$ with trivial $\Ga$-stabilizer. 
Recall that in \eqref{eq:U-nbd} and \eqref{eq:U-nbd1} we defined neighborhood bases $U(p,t)$, $U(P,t)$ of points $p\in \geo X$ and compact subsets $P\subset \geo X$ using Gromov products. 
We will use $x_0$ as the base-point in that definition, i.e. 
\begin{align} 
\begin{aligned}\label{eq_Upt}
U(p,t) = \{ z \in \ol{X} |  (z, p)_{x_0} > - \log(t)  \}\subset \ol{X},\\
U(P,t) = \{ z \in \ol{X} | \sup_{p\in P} (z, p)_{x_0} > - \log(t)  \}\subset \ol{X}
\end{aligned}
\end{align}


For $t >0$, a finite subset $E\subset \Gamma$ and a compact subset $P\subset \geo X$ we consider the subsets 
\begin{align} 
\begin{aligned}
\label{eq_Gpt}\index{$\Gamma_{P,t}$}
&\Gamma_{P,t} = \left \{ \gamma \in \Gamma \ | \ \gamma^{-1}(x_0) \notin U(P,t) \right \}\subset \Ga, 
\\
\Gamma_{P,E,t} &= \left \{ \gamma \in \Gamma \ | \ \forall g\in E, g^{-1}\circ\gamma^{-1}(x_0) \notin U(P,t) \right \} = 
\bigcap_{g\in E} \Gamma_{P,t}g^{-1}. 
\end{aligned}\index{$\Gamma_{P,E,t}$}
\end{align}

The complement of $\Gamma_{P,t}$ in $\Gamma$ will be denoted by $\Gamma_{P,t}^c$, similarly for $\Gamma_{P,E,t}^c$. Then we have
\[
\Gamma_{P,E,t}^c=\bigcup_{g\in E}\Gamma_{P,t}^cg^{-1}.
\] 
If $P=\{p\}$ is a singleton then we use the notation $\Gamma_{p,E,t},\Gamma_{p,E,t}^c$.
\index{$\Gamma_{P,E,t}^c$}

Recall that the annuli 
$A_{n_1}^{n_2}$ and $A_n=A_n^{n+1}$ are defined in the previous sections. 
The inversion $\gamma\mapsto \gamma^{-1}$ induces a bijection
$$
B(x_0, R)\cap \Ga x_0 \to B(x_0, R)\cap \Ga x_0. 
$$
and restricts to a bijection $A_n\to A_n$ for all $n\in \mathbb N$. Then we deduce directly from the definition of $\Gamma_{P,E,t}$ the following lemma:

\begin{lem}\label{lem:shadowcounting_inverse_singlepoint}
For any $t>0$ and any compact subset $P\subset \geo X$ we have a bijection:
\[
A_n \cap \Gamma_{P,t}^c \rightarrow A_n(x_0) \cap \Uupt,\quad \gamma\mapsto \gamma^{-1}(x_0)
\]
\end{lem}
As $\Gamma_{P,E,t}^c=\bigcup_{g\in E}\Gamma_{P,t}^c \,g^{-1}$, we then get:

\begin{lem}\label{lem:shadowcounting_inverse}
Let $E$ be a finite subset of $\Gamma$. Set $m_0=\ceil{\max_{g\in E}\beta(g)}$. 
Then for any $t>0$ and any compact subset $P\subset \geo X$ we have:
\[
\# A_n \cap \Gamma_{P,E,t}^c \leqslant (\#E)\# \left(A_{n-m_0}^{n+m_0}(x_0) \cap \Uupt\right).
\]
\end{lem}
\begin{proof}
We have  
\[
A_n\cap\Gamma_{P,E,t}^c=\bigcup_{g\in E}A_n\cap \left(\Gamma_{P,t}^c \, g^{-1}\right).
\] 
If $\gamma\circ g^{-1}\in A_n$, then $n-m_0<\beta(\gamma)\leqslant n+m_0$ by triangle inequality. Therefore, for every $g\in E$, 
the map $\gamma\mapsto \gamma\circ g$ is an injection from $A_n\cap \left(\Gamma_{P,t}^cg^{-1}\right)$ to $A_{n-m_0}^{n+m_0}\cap \Gamma_{P,t}^c$.
As in Lemma \ref{lem:shadowcounting_inverse_singlepoint}, we have a bijection
\[
A_{n-m_0}^{n+m_0}\cap \Gamma_{P,t}^c 
\rightarrow A_{n-m_0}^{n+m_0}(x_0) \cap \Uupt, \quad \gamma\mapsto \gamma^{-1}(x_0).
\]
\end{proof}

Since the $\Ga$-action on $X$ is convex-cobounded, there exists $D<\infty$ such that every geodesic ray $x_0\xi, \xi\in \La(\Ga)$, is contained in the $D$-neighborhood of the orbit $\Ga x_0$, see Section \ref{sec:quasiconvex}. 

\begin{prop} \label{prop_shadow_counting} 
For $n\in \mathbb N$ and positive $t$ satisfying $\log(t) < -2\delta -1$, define $r(t):=  2 + 3D+  4\delta + \log t$. 
Then $A_n \cap \Gamma_{p,t}^c$ is empty if $n+r(t)\leqslant 0$. Whenever $n+r(t)>0$, for every $p\in \La(\Ga)$ 
we have  
\begin{equation*}
\# (A_n(x_0) \cap U(p,t) )=\# (A_n \cap \Gamma_{p,t}^c ) \leqslant \# \Gamma_{n+r(t)}.
\end{equation*}
\end{prop}

\begin{proof}

Take a point $y$ on a geodesic ray $\rho=x_0p$ from $x_0$ to $p$, such that 
$$d(x_0, y)= - \log t -2\delta - 1.$$ 
There exists  $\gamma_0 \in \Gamma$ such that $\gamma_0(x_0)$ is within distance $D$ from $y$. (Here we are using the assumption that $p$ is a limit point of $\Ga$.) Thus,
\begin{equation}\label{eq:yx_0}
d(x_0, \ga_0(x_0))\geqslant - \log t -2\delta -1 - D. 
\end{equation}

Note that for every $x\in U(p,t)\cap X$ we have  
$$
d(x_0, xp)\geqslant (x,p)_{x_0}\geqslant - \log(t) >   d(x_0, y) +2\delta. 
$$
It follows that 
$$
d(y, xp)> \delta. 
$$
By the $\delta$-slimness of geodesic triangle $\Delta x_0 x p$ we, therefore, obtain that every geodesic $x_0x$ intersects the closed ball 
$\bar{B}(y,\delta)$. In other words, $\bac{x}\in \Sh_{x_0}(y,\delta)$, i.e., 
$$
U(p,t)\subset  \Sh_{x_0}(y,\delta). 
$$
Hence,
$$
U(p,t)\subset  \Sh_{x_0}(\ga_0(x_0),\delta +D).  
$$


Suppose now that 
$\gamma \in A_n\subset \Ga$ and  $x:=\ga(x_0)\in U(p,t)$. Then $x\in  \Sh_{x_0}(\ga_0(x_0),\delta +D)$, i.e. there exists 
$$
z\in x_0 x\cap \bar{B}(\ga_0(x_0),\delta +D). 
$$
 
Set $\mathcal{D}=\delta+D$. Combining the inequality $d(x_0, \ga (x_0)) \leqslant n+1$ with the triangle inequality and the inequality \eqref{eq:yx_0}, we obtain: 
\begin{align*}
d(\gamma (x_0), \gamma_0 (x_0)) &\leqslant  d(\gamma(x_0), z) + d(z, \gamma_0(x_0))  \\
& \leqslant n+1 + \mathcal{D} - d(x_0, z) \\
 & \leqslant n+1 + 2\mathcal{D} - d(x_0, \gamma_0(x_0))\\
 & \leqslant n+1 + 2\mathcal{D}  + \log t + 2\delta +1 + D=  n+r(t). 
\end{align*}
This shows that $\gamma_0^{-1}\gamma \in \Gamma_{n+r(t)}$ and we obtain that
\begin{equation*}
\# \left(A_n(x_0) \cap U(p,t)\right) \leqslant\# \Gamma_{n+r(t)}.
\end{equation*}
Then the conclusion follows from Lemma \ref{lem:shadowcounting_inverse_singlepoint}.  
\end{proof}

\subsection{Equidistribution on the ideal boundaries of Gromov-hyperbo\-lic spaces}
\label{sec:Equidistribution on ideal boundaries of Gromov-hyperbolic spaces}


As before, we work in the setting of \S \ref{sec:counting_setting}. We will be defining the Gibbs measures $\mu_{X}$ using this length function. We continue with the notation introduced in the 
previous section and Section \ref{sec:Gibbs}. 

Take a point $p \in \geo  X$ and let  $\rho=x_0p$ be a geodesic ray in $X$ joining $x_0$ and $p$. 
Recall from Section \ref{section_conical} that given a choice of such points $x_0\in X, p\in \geo X$, we have partitioned 
the group $\Gamma$ as $\Gamma = \Gamma_{N,bad} \sqcup \Gamma_{N,good}$ according to how far the group elements 
move $x_0$ from the ray $\rho$.    

According to this partition of the group, one obtains the following decomposition of the measures $\mu_{x_0}(s), \mu_p(s)$, $s>s_0$, as: 
\begin{equation}
\mu_{x_0}(s) =   \mu_{x_0,  N, good}(s) + \mu_{x_0,N, bad}(s). 
\end{equation}
and 
\begin{equation}
\mu_p(s) =   \mu_{p,  N, good}(s) + \mu_{p,N, bad}(s). 
\end{equation}

\begin{prop} \label{prop_mass_bad_orbits_PS} For every $s > s_0$ and $N\in \mathbb N$, 
the mass $|\mu_{x_0,N,bad}|$ is controlled by: 
\begin{equation}
|\mu_{p,N,bad}(s)| =|\mu_{x_0,N,bad}(s)| \leqslant \dfrac{C}{\cP(s)} \#  \Gamma_{1+4N}  
\end{equation}
where $C>0$ is a constant which does not depend on $s, N$.
\end{prop}

\begin{proof}

We estimate the mass of the measure $\mu_{x_0,N,bad}(s)$ using Proposition \ref{prop_conical_counting}.

\begin{align*}
|\mu_{x_0, N,bad}(s)| &= \dfrac{1}{\cP(s)} \sum_{\gamma\in \Gamma_{N,bad}} e^{-s d(x_0, \gamma(x_0))} |\delta_{\gamma(x_0)}| \\
 & = \dfrac{1}{\cP(s)} \sum_{\gamma\in \Gamma_{N,bad}} e^{-s d(x_0, \gamma(x_0))} \\
  &  \leqslant  \dfrac{\# \Gamma_N}{\cP(s)} + \dfrac{1}{\cP(s)} \sum_{n \geqslant N} \sum_{\gamma\in A_n \cap \Gamma_{N,bad}} e^{-s d(x_0, \gamma(x_0))} 
   \\
  & \leqslant \dfrac{\# \Gamma_N}{\cP(s)} + \dfrac{1}{\cP(s)}  \sum_{n \geqslant N} e^{-n s} \# (A_n \cap \Gamma_{N,bad}) 
   \\
  & \leqslant  \dfrac{\# \Gamma_N}{\cP(s)} + \dfrac{1}{\cP(s)} \# \Gamma_{1+4N}  \dfrac{e^{-Ns}}{1 -e^{-s}}   \\
  & \leqslant \dfrac{C}{\cP(s)} \# \Gamma_{1+4N}, 
\end{align*}
for all $s > s_0$, where $C= (1 -e^{-s_0})^{-1}$ does not depend on $s$ and $N$. 
The proof for $ \mu_{p, N,bad}$ is the same. 
\end{proof}

For the next statement, we will exploit the properties of the elements $\gamma \in  \Ga_{N,good}$ for which the distance from $x_0$ to $\gamma(\rho)$ is  $\geqslant N$. 


\begin{prop} \label{prop_good} Given $x_0 \in X$ and $p \in \ol{X}$, 
for any continuous function $\varphi$ on $\ol{X}$ and for any $\epsilon>0$, there exists an 
integer $N=N (x_0,\epsilon, \varphi) $ such that: 
\begin{equation*}
| \langle \mu_{x_0, N, good}(s) - \mu_{p,N,good}(s) , \varphi\rangle | \leqslant \epsilon. 
\end{equation*} 
\end{prop}

\begin{proof}
Recall (see Section \ref{sec:Gromov-hyperbolic spaces}) that there exist positive numbers $a$ and $\alpha$ such that $d_{\infty}(\gamma(x_0),\gamma(p)) \leqslant a e^{-\alpha N}$ if the distance from $x_0$ to the geodesic joining $\gamma(x_0)$ and $\gamma(p)$ is larger than $N$. (One can take $a=k'_1$ from the inequality  \eqref{eq:dinfty}.)  
This implies that for any $\gamma \in \Gamma_{N,good}$, $d_{\infty}(\gamma(x_0), \gamma(p)) \leqslant ae^{-\alpha N}$. 

Take $\eta>0$ to be the modulus of continuity of $\varphi$ with respect to $\epsilon>0$. We choose an integer $N$ large enough satisfying $ae^{-\alpha N} \leqslant \eta$;  we then have for any $\gamma \in \Gamma_{N,good}$: 
\begin{equation*}
d_{\infty}(\gamma(x_0), \gamma(p)) \leqslant ae^{-\alpha N} \leqslant \eta.
\end{equation*}
This gives $|\varphi(\gamma(x_0)) - \varphi( \gamma(p))| \leqslant \epsilon $ for any $\gamma \in \Gamma_{N,good}$.

We conclude that: 
\begin{align*}
| \langle \mu_{x_0, N, good}(s) - \mu_{p,N,good}(s) , \varphi\rangle |  \\
 \leqslant \dfrac{1}{\cP(s)} \sum_{\gamma \in \Gamma_{N,good} } e^{-s\beta(\gamma)} |\varphi(\gamma(x_0)) - \varphi(\gamma(p))| \\
\leqslant \left  (\dfrac{\sum_{\gamma \in \Gamma_{N,good} } e^{-s\beta(\gamma)}}{\cP(s)} \right )  \epsilon \leqslant \epsilon,
\end{align*}
as required. 
\end{proof}

\begin{thm}\label{thm:rank1-equi} 
Suppose that $\Gamma$ is a non-elementary quasiconvex-cocom\-pact group of isometries of a proper geodesic 
Gromov hyperbolic space $X$. 
Fix $x_0 \in X$.  For every $p \in \geo  X$, and every sequence $s_n\rightarrow s_0^+$ such that the sequence of measures   
$ \mu_{x_0}(s_n)$ converges weakly, the limit of $\mu_p(s_n)$ also exists and  we have  
\begin{equation*}
\lim_{n\rightarrow +\infty} \mu_p(s_n) = \lim_{n\rightarrow +\infty} \mu_x(s_n)
\end{equation*}
\end{thm}
  
\begin{proof}
The existence of such a limiting measure $\mu_p$ for a subsequence in $(\mu_p(s_n))$ 
comes by extracting a subsequence in the original PS sequence $(s_n)$, 
converging weakly to a probability measure.   
Let us show that this measure does not depend on the choice of a subsequence and does not depend on the choice of $p$. 

Fix $\epsilon >0$ and $\varphi$, a continuous function on $\bar X$. 
By Proposition \ref{prop_good}, we get an integer  $N = N(\epsilon,\varphi)$  such that for any $s > s_0$ :
\begin{equation} \label{eq_good}
| \langle \mu_{x, N , good}(s) - \mu_{p, N , good}(s) , \varphi \rangle | \leqslant \epsilon.
\end{equation}

By Proposition \ref{prop_mass_bad_orbits_PS}, there exists a constant $C> 0$ such that we can control the masses of ``bad'' parts of 
the measures:  
\begin{equation*} 
|\mu_{p, N , bad}(s)| = |\mu_{x_0, N , bad}(s)| \leqslant \dfrac{C}{\cP(s)} \# \Gamma_{1+4N} .
\end{equation*}

Since the Poincar\'e series $\cP(s)$ diverges for $s=s_0$, for fixed $N$ we can choose $s_1=s_1(N) > s_0$ 
such that for any $s \in (s_0,s_1] $, one has: 
\begin{equation*}
\dfrac{C}{\cP(s)} \# \Gamma_{1+4N} \leqslant \epsilon.
\end{equation*}

In particular, one obtains: 
\begin{equation}
\label{eq_bad}
|\mu_{p, N , bad}(s)| = |\mu_{x, N , bad}(s)| \leqslant \dfrac{C}{\cP(s)} \# \Gamma_{1+4N} \leqslant \epsilon. 
\end{equation}

Overall \eqref{eq_good} and \eqref{eq_bad}, yield for all $s_0< s \leqslant s_1$: 
\begin{align*}
| \langle \mu_{p}(s) - \mu_{x}(s) , \varphi \rangle |& \leqslant  |\varphi|_\infty (  |\mu_{p, N,bad}(s)| + |\mu_{x_0,N,bad}(s)|  )  \\
& \quad + | \langle \mu_{x_0, N , good}(s) - \mu_{p, N , good}(s) , \varphi \rangle | \\
& \leqslant \epsilon \left ( 1 + 2 |\varphi|_\infty   \right ).  
\end{align*}
Here $|\varphi|_\infty$ is the supremum-norm of $\varphi$ on $\ol{X}$. 

We, thus, have thus shown that 
$$
\lim_{n\rightarrow \infty}\mu_{p}(s_n)=\lim_{n\rightarrow \infty} \mu_{x_0}(s_n).
$$
\end{proof}

 \begin{rem}
 The situation is slightly different in the case of infinite elementary groups $\Ga$. In this case, $\La(\Ga)$ consists of two points $\la_\pm$ and either they are both fixed by $\Ga$ or they are swapped by some element  of $\Ga$. If $p\in \La(\Ga)\subset \geo X$ is fixed by $\Ga$, then, of course, the family of measures $\mu_{p}(s)$ is constant and converges to $\delta_p$. In all other cases (i.e. either 
 $p\notin \La(\Ga)$ or 
 $p\in \La(\Ga)$ and $\Ga$ acts nontrivially on $\La(\Ga)$), the family of measures $\mu_{p}(s)$   converges to the probability measure that has equal mass $1/2$ at each point $\la_\pm$ of $\La(\Ga)$. 
 \end{rem}

 As an immediate corollary of the theorem, we obtain:

 \begin{cor}\label{cor:Gibbs}
 Let $\Ga$ be a nonelementary hyperbolic group, $G\times\La\to \La$ be a continuous action of $G$ on a topological space $\La$ and $\xi: \geo \Ga\to \La$ a $\Ga$-equivariant homeomorphism. Let $\beta$ be a length function on $\Ga$ as above and $\cP(s)$ be the associated Poincar\'e series which diverges at the critical exponent $s_0$. 
 Then for each $\la\in \La$, the family of probability measures 
 $$
 \mu_\la(s):=  \frac{1}{\cP(s)}\sum_{\gamma\in \Gamma}e^{-s\beta(\gamma)}\delta_{\gamma(\la)}.
 $$ 
satisfies 
 $$
  \lim_{s_n\to s_0+} \mu_\la(s_n)= \xi_*(\mu_{\Ga}).  
 $$
 \end{cor}
 
 For instance, if $\Ga$ is a nonelementary $P$-Anosov subgroup of a Lie group $G$ and $\La\subset G/P$ is the flag-limit set of $\Ga$, then the above corollary applies. 

\subsection{Convergence of sequences of partial measures to Gibbs measures} \label{sec:partial measures}

We assume the setup of Section \ref{sec:counting_setting}.


We fix $\alpha\dyl$ and an $\alpha$-visual metric $d_\alpha$ on $\geo X$. We also fix constants $k_1,k_2,t_0\dyl$ and a metric $d_\infty$ on $\ol{X}$ that extends $d_\alpha$ and satisfies the following inequality (see Section \ref{sec:Visual metrics}): For all $x,y\in U(\geo X,t_0)$ we have 
\begin{equation}\label{eq:dinfty_gromovproduct_bounded}
k_1\exp(-\alpha (x,y)_{x_0})\leqslant d_\infty(x,y)\leqslant k_2\exp(-\alpha (x,y)_{x_0}).
\end{equation}

\subsubsection{Some equivalent quantities}
\begin{lem}\label{lem:equivalent_poincare_series}
For $s>s_0$, let $\cQ(s)$ be the series:
\[
\cQ(s)=\sum_{n=0}^\infty \#\Gamma_n e^{-sn}.
\]
Then, $\cQ(s) \asymp \cP(s)$ 
as $s\rightarrow s_0$.

\end{lem}
\begin{proof}
We have
\begin{align*}
\cP(s)=\sum_{n\geqslant 0}\sum_{\gamma \in A_n}e^{-s\beta(\gamma)}.
\end{align*}
Any $\gamma \in A_n$ satisfies $e^{-s(n+1)}\leqslant e^{-s\beta(\gamma)}\leqslant e^{-sn}$. Thus, denoting $\cR(s)=\sum_{n\geqslant 0}\# A_n e^{-sn}$, we get 
\[
e^{-s}\cR(s)\leqslant \cP(s) \leqslant \cR(s).
\]
Now the conclusion follows because $\cQ(s)$ is a multiple of $\cR(s)$. Indeed, since $\#A_n=\#\Gamma_{n+1}-\#\Gamma_n$, we have $\cQ(s)=\frac{1}{e^s-1}\cR(s)$. 
\end{proof}

For $\fm\dyl$, an \emph{$\fm$-covering} of a subset $A\subset \geo X$ is a finite covering of $A$ by sets of the form $U(p,\fm)$ (or $U(p,\fm)\cap \geo X$) with $p\in \geo X$. Any $A\subset \geo X$ admits such coverings because $\geo X$ is compact. Then we define
\begin{equation}
\fM_\fm(A):= \fm^{\alpha s_0}\min_{(U_i)_{i\in I}}\# I 
\end{equation}
where the minimum is taken over all $\fm$-covering $(U_i)_{i\in I}$ of $A$, and
\begin{align}
\begin{aligned}
\ol{\fM}(A):=\limsup_{\fm\rightarrow 0}\fM_\fm(A),\\
\underline{\fM}(A):=\liminf_{\fm\rightarrow 0}\fM_\fm(A).
\end{aligned}
\end{align}
The definition of $\ol{\fM},\underline{\fM}$ is similar to that of the Hausdorff (outer) measure but, a priori, $\ol{\fM},\underline{\fM}$ do not define measures. 

\begin{rem}\label{rem:boxcounting_upt}
The inequality \eqref{eq:dinfty_gromovproduct_bounded} implies that in the definition of $\ol{\fM},\underline{\fM}$, when we consider coverings, we can use $d_\al$-balls instead of subsets of the form $U(p,\fm)$.
\end{rem}

It is known that the Gibbs-measure is bounded above and below by constant multiples of the $s_0$-dimensional Hausdorff measure 
(see Theorem \ref{thm:Hausdorff-qc} as well as 
\cite{Coo93}, \cite{Calegari}). A similar sandwich property is true for $\ol{\fM},\underline{\fM}$:
\begin{prop}\label{prop:box_counting_measure}
There exist two constants $k_5,k_6\dyl$ such that: 
\begin{enumerate}
\item for any Borel subset $A\subset \geo X$, we have $\mu(A)\leqslant k_5 \underline{\fM}(A)$;
\item for any compact subset $P\subset \geo X$, we have $
k_6\ol{\fM}(P)\leqslant \mu(P)$.
\end{enumerate}
\end{prop}
\begin{proof}
The key property is the following: by \eqref{eq:dinfty_gromovproduct_bounded} and \cite[Proposition 7.4]{Coo93}, there exist constants $k_7,k_8\dyl$ such that for any $p\in \geo X$, for any $\fm\dyl$ sufficiently small, 
\[
k_7 \fm^{\alpha s_0}\leqslant \mu\left(U(p,\fm)\right)\leqslant k_8 \fm^{\alpha s_0}.
\]

Let $A\subset \geo X$ be a Borel subset. Consider an $\fm$-covering $(U_i)_{i\in I}$ of $A$. Then 
\begin{align*}
\mu(A)\leqslant \sum_{i\in I} \mu(U_i) \leqslant k_8  \fm^{\alpha s_0}\# I.
\end{align*}
Therefore $\mu(A)\leqslant k_8\fM_\fm(A)$. Letting $\fm\rightarrow 0$, we get $\mu(A)\leqslant k_8 \underline{\fM}(A)$. This argument is the same as in \cite{sullivan_density}.

The proof of the reverse inequality is essentially the same as \cite[Corollary 2.5.10]{Calegari}. Let $P\subset \geo X$ be a compact subset. For $\delta_1\dyl$ we take an open subset $U\subset \geo X$ containing $P$ such that $\mu(U-P)\leqslant \delta_1$. By compactness of $P$, there exists a $\delta_2\dyl$ depending on $\delta_1$ such that any subset of the form $U(p,\delta_2)\cap \geo X$ with $p\in P$, is contained in $U$. For 
$\delta\leqslant \delta_2$ we can inductively construct a $\delta$-covering $(U(p_i,\delta))_{1\leqslant i\leqslant l}$ of $P$ such that each $p_j$ belongs to $P$ but does not belong to any $U(p_i,\delta)$ with $i<j$. Then, there exists $k_9\dyl$ such that the subsets $U(p_i,k_9\delta),1\leqslant i \leqslant l$, are pairwise disjoint. Hence
\begin{align*}
l\delta^{\alpha s_0}& =\frac{l}{k_9^{\alpha s_0}}(k_9\delta)^{\alpha s_0}
\leqslant \frac{1}{k_7k_9^{\alpha s_0}}\sum_{i=1}^l  \mu\left(U(p_i,k_9\delta)\right)\\
&\leqslant \frac{1}{k_7k_9^{{\kap \alpha}s_0}}\mu(U).
\end{align*}
Therefore, for some constant $k_6\dyl$, we get $\fM_\delta(P)\leqslant k_6\mu(U)$ and, thus,  
$$
\ol{\fM}(P)\leqslant k_6\mu(U) \le  k_6\mu(P) + k_6\delta_1. 
$$ 
Since this holds for all $\delta_1>0$, 
we get $\ol{\fM}(P)\leqslant k_6\mu(P)$. 
\end{proof}

\subsubsection{Convergence of partial measures}

Recall that at the beginning of Section \ref{section_shadow_counting} we have defined the subsets $\Gamma_{P,t},\Gamma_{P,E,t}$ of $\Gamma$ for any $t > 0$, any compact $P\subset \geo X$ and any finite subset $E\subset \Gamma$. Fix $P,E$ and let $\mu_{s,t}$ be the measure given by:

\begin{equation}\label{eq:def_mu_st}
\mu_{s,t}:= \dfrac{1}{\cP(s)} \sum_{\gamma \in \Gamma_{P,E,t}} e^{-s \beta(\gamma)} \delta_{\gamma(x_0)}. 
\end{equation}

This is a {\em partial measure}, supported on $\Ga_{P,E,t} \cdot x_0\subset \Ga \cdot x_0$. Note that the measure $\mu_{s,t}$ depends on the choice of $x_0$, $P$ and $E$ but we omit the corresponding subscripts to simplify the notation. As before, we let $s_n\to s_0+$ denote a PS-sequence for the family of measures $(\mu_s)$. We will be interested in limits of families of measures $\mu_{s_n,t}$. The situation here is similar to the one in \S 
\ref{sec:Equidistribution on ideal boundaries of Gromov-hyperbolic spaces},  but is more complicated. In \S 
\ref{sec:Equidistribution on ideal boundaries of Gromov-hyperbolic spaces} we were working with sequences of measures 
$\mu_{x_0, N, good}(s_n)$ which have larger support than the measures  $\mu_{s_n,t}$: The  complement of the good set (the bad set $\Ga_{N,bad}$) 
 grows linearly in terms of $n$ when intersected with $\Ga_n$. In contrast, the subsets  excluded in this section,  
$\Gamma_{P,E,t}^c$, have exponential growth with respect to the length function $\beta$. Nevertheless, we will see that this growth has slow rate so that the partial measures $\mu_{s_n,t}$ still limit (in appropriate sense) to $\mu_{X}$, as before.

\begin{prop} \label{prop:control_must_singlepoint} Let  $E\subset \Gamma$ be a finite set and $p\in \geo X$ be a point. 
The measures $\mu_{s,t}$ associated with $E$ and the singleton $\{p\}$ satisfy that for any $t>0$ sufficiently small, we have:
$$
\limsup_{s \rightarrow s_0}\vert \mu_{s,t}-\mu_{s}\vert \leqslant \fC_1 t^{s_0},  
$$
where $\fC_1$ is a constant that does not depend on $t$.
\end{prop}
\begin{proof}
We first bound the mass $|\mu_s - \mu_{s,t}|$ as follows: 
\begin{align*}
|\mu_s - \mu_{s,t}| & =  \dfrac{1}{\cP(s)} \sum_{n=0}^\infty \sum_{\gamma \in A_n \cap \Gamma_{p,E,t}^c} e^{-s\beta(\gamma)} \\
   & \leqslant \dfrac{1}{\cP(s)} \sum_{n=0}^\infty  \# (A_n \cap \Gamma_{p,E,t}^c) e^{-sn} 
\end{align*}
Recall that in Proposition \ref{prop_shadow_counting}, for $t > 0$ we have defined the function $r(t)= \log(t)+2+3D+4\delta$. 
Without loss of generality, we assume that $r(t)\leqslant 0$. Note that by Proposition \ref{prop_shadow_counting}, we have 
$A_n \cap \Gamma_{p,E,t}^c=\emptyset$ when $0>n+r(t)$. Recall that $m_0=\ceil{\max_{g\in E}\beta(g)}$. 
Using Proposition \ref{prop_shadow_counting} and Lemma \ref{lem:shadowcounting_inverse}, we get, 
\begin{align*}
|\mu_s - \mu_{s,t}| 	&\leqslant \dfrac{\#E }{\cP(s)} \sum_{n\geqslant 0}\# \left(A_{n-m_0}^{n+m_0}(x_0)\cap U(p,t) \right) e^{-s n}\\
&= \dfrac{\#E }{\cP(s)} \sum_{n\geqslant 0}\sum_{m=-m_0}^{m_0-1}\# \left(A_{n+m}(x_0)\cap U(p,t) \right) e^{-s n}\\
&= \dfrac{\#E }{\cP(s)} \sum_{n\geqslant 0}\sum_{m=-m_0}^{m_0-1}e^{sm}\# \left(A_{n+m}(x_0)\cap U(p,t) \right) e^{-s (n+m)}\\
&= \dfrac{\#E\sum_{m=-m_0}^{m_0-1}e^{sm} }{\cP(s)} \sum_{n\geqslant 0}\# \left(A_n(x_0)\cap U(p,t) \right) e^{-s n}\\
&\leqslant \dfrac{\#E\sum_{m=-m_0}^{m_0-1}e^{sm}}{\cP(s)} \sum_{n\geqslant -r(t)} \# \Gamma_{n+r(t)} e^{-s n}\\
&= e^{sr(t)}\dfrac{\#E\sum_{m=-m_0}^{m_0-1}e^{sm}}{\cP(s)} \sum_{n\geqslant 0} \# \Gamma_{n} e^{-s n}
\end{align*}
We have  
$$
\lim_{s\to s_0} e^{sr(t)}= \lim_{s\to s_0} \exp(s(\log(t)+2+3D+4\delta))= e^{s_0(2+3D+4\delta)} t^{s_0}$$
and
$$
\lim_{s\to s_0}  \#E\sum_{m=-m_0}^{m_0-1}e^{sm} \leqslant 2m_0e^{s_0 m_0} \#E .
$$
By 
Lemma \ref{lem:equivalent_poincare_series}, we also obtain
$$
\limsup_{s\to s_0} \frac{1}{\cP(s)}   \sum_{n\geqslant 0} \# \Gamma_{n} e^{-s n}\leqslant k_4.  
$$
Proposition follows. 
\end{proof}

\begin{prop}\label{prop:estimate_mus_must}
For any finite set $E\subset \Gamma$, any compact set $P\subset \geo X$, the corresponding measures $\mu_{s,t}$ satisfy:
\begin{equation}\label{eq:double_limit_mus_must}
\limsup_{t\rightarrow 0}\limsup_{n\rightarrow \infty}\vert \mu_{s_n}-\mu_{s_n,t}\vert \leqslant \mathfrak{C}\mu\left(P\right)^{1/\al}.
\end{equation} 
where $\fC$ is a constant. 
\end{prop}

\begin{proof}
In this proof $\fC$, with or without a subscript, will denote a positive constant that does not depend on $s,t$.

For $t\dyl$ sufficiently small, we take a $t$-covering 
\[
\{U(p_i,t)\cap \geo X: p_i\in \geo X, i\in I_t\}
\]
of $P$ such that 
\begin{equation}\label{eq:must_box1}
t^{\alpha s_0}\#I_t = \fM_{t}(P).
\end{equation}

By the sandwich inequality \eqref{eq:dinfty_gromovproduct_bounded}, subsets of the form $U(p,t)$ behave like $d_\infty$-balls. In particular, there exists a constant $\fC_2\dyl$ that does not depend on $t$ 
such that
\[
U(P,t)\subset \bigcup_{i\in I_t}U(p_i,\fC_2t).
\]
Therefore, by Lemma \ref{lem:shadowcounting_inverse}, we have, for $s>s_0$ sufficiently close to $s_0$:
\begin{align*}
|\mu_s - \mu_{s,t}| 
   & \leqslant \dfrac{1}{\cP(s)} \sum_{n=0}^\infty  \# (A_n \cap \Gamma_{p,E,t}^c) e^{-sn} \\
   & \leqslant \dfrac{1}{\cP(s)} \sum_{n=0}^\infty  \# (A_{n-m_0}^{n+m_0}(x_0) \cap U(P,t) ) e^{-sn}\\
   & = \dfrac{1}{\cP(s)} \sum_{n=0}^\infty \sum_{m=-m_0}^{m_0} \# (A_{n+m}(x_0) \cap U(P,t) ) e^{-sn}\\
   & =\dfrac{1}{\cP(s)} \sum_{n=0}^\infty \sum_{m=-m_0}^{m_0}e^{sm} \# (A_{n+m}(x_0) \cap U(P,t) ) e^{-s(n+m)}\\
   &  \leqslant \dfrac{\fC_3 }{\cP(s)}\sum_{i\in I_t}\sum_{n=0}^\infty  \# (A_n(x_0) \cap U(p_i,\fC_2t) ) e^{-sn},
\end{align*}
where $\fC_3$ does not depend on $t$. 
The terms corresponding to a single $p_i$ are estimated in Proposition \ref{prop:control_must_singlepoint}. Therefore:
\begin{align*}
\limsup_{j\rightarrow\infty}|\mu_{s_j} - \mu_{s_j,t}| 
   & \leqslant \fC_4\#I_t t^{s_0},
\end{align*}
for a constant $\fC_4$ which does not depend on $t$. 
Now, by \eqref{eq:must_box1} we get 
\[\limsup_{t\rightarrow 0}\#I_t t^{\alpha s_0}=\ol{\fM}(P).\]
Therefore we have:
\begin{align*}
\limsup_{t\rightarrow 0}\limsup_{j\rightarrow\infty}|\mu_{s_j} - \mu_{s_j,t}| 
   & \leqslant \fC_5\left(\ol{\fM}(P)\right)^{1/\al}.
\end{align*}
Now the conclusion follows from the second part of Proposition \ref{prop:box_counting_measure}.
\end{proof}

For a fixed $t$, let $\cM_t$ be the cluster set of the sequence of measures $(\mu_{s_n,t})_n$ (with respect to the topology of weak convergence), i.e.\
\begin{equation}\label{eq:def_Mt}
\cM_t=\{\lim_{n\rightarrow \infty}\mu_{s'_n,t}: (s'_n) \ \text{subsequence of} \ (s_n)\}.
\end{equation} 
\begin{rem}\label{rem:support_mut}
\begin{enumerate}
\item The set $\cM_t$ is weakly compact. 
\item Like $\mu$, any $\nu$ in $\cM_t$ is supported on $\geo X$ because $\cP(s_0)=\infty$. 
\end{enumerate}
\end{rem}

\begin{prop} \label{prop:double_limit_ordered}
Fix a PS-sequence $(s_n)$, the corresponding Gibbs-measure $\mu$ on $\geo X$, a finite subset $E\subset \Gamma$ and a compact 
subset $P\subset \geo X$ such that $\mu(P)=0$. For any $t>0$ take an element $\nu_t$ in $\cM_t$. Then one has the weak convergence: 
\begin{equation*}
\lim_{t \rightarrow  0^+}\nu_t = \mu. 
\end{equation*}
\end{prop}

\begin{proof}
Pick a subsequence $(s'_n)$ of $(s_n)$ such that $\mu_{s'_n,t}$ converges to $\nu_t$. Taking $\psi$ to be a continuous test function on $\ol{X}$, one decomposes:
\begin{align*}
|\langle \mu_{s'_n,t} - \mu , \psi \rangle| &\leqslant |\mu_{s'_n,t} - \mu_{s'_n}|  |\psi|_\infty + |\langle \mu_{s'_n} - \mu , \psi \rangle|. 
\end{align*}
When $n\rightarrow\infty$, the second term on the right-hand side converges to zero because $\mu_{s'_n}$ converges weakly to $\mu$. Thus, by Proposition \ref{prop:estimate_mus_must}, we have:
\begin{align*}
\limsup_{t\rightarrow 0}|\langle \nu_t - \mu , \psi \rangle| &\leqslant \limsup_{t\rightarrow 0} \limsup_{n'\rightarrow \infty} |\mu_{s'_n,t} - \mu_{s'_n}|  |\psi|_\infty\\
& 
\leqslant \mathfrak{C}|\psi|_\infty\left(\mu(P)\right)^{1/\al}.
\end{align*}
The assumption that $\mu(P)=0$ allows us to conclude.
\end{proof}

%
%
%

\begin{thm} \label{thm:integration_partial_measures}
Assume the setup of Section \ref{sec:counting_setting}. Let $(s_m)_m$ be a PS-sequence and $\mu$ be the corresponding Gibbs-measure on $\geo X$. Let $E\subset \Gamma$ be a finite set and $P\subset \geo X$ a compact subset such that $\mu(P)=0$. Assume that $\psi:\ol{X}\rightarrow \bC$ is a function such that: 
\begin{enumerate}
 \item[(i)] the restriction of $\psi$ to 
 $\geo X$ 
 is continuous;   
 \item[(ii)] for every $t >0$, the restriction of $\psi$ to $\overline{\Gamma_{P,E,t}x_0}\subset \overline{X}$ is continuous.
\end{enumerate} 
Then one has
\begin{equation}\label{eq:main_convergence}
\lim_{m \rightarrow +\infty} \dfrac{1}{\cP(s_m)}\sum_{\gamma \in \Gamma} e^{-s_m \beta(\gamma)} \psi(\gamma x_0) = \int_{\geo X} 
\psi(p) d\mu(p).
\end{equation}
\end{thm}

\begin{proof}
The convergence \eqref{eq:main_convergence} is equivalent to the convergence:
\begin{equation}\label{eq:main_convergence1}
\lim_{m \rightarrow +\infty} \int_{\overline{X}}\psi(x)d\mu_{s_m}(x) = \int_{\geo X} \psi(p)d\mu(p).
\end{equation}

Recall that $\cM_t$ was defined in \eqref{eq:def_Mt}. Set
\[
\cM_t(\psi)=\{\int \psi d\nu: \nu\in \cM_t\}.
\]
Since $\psi$ is continous on $\overline{\Gamma_{P,E,t}x_0}$, the set $\cM_t(\psi)$ is compact and coincides with the cluster set of the numerical sequence $(\int_X \psi  d\mu_{s_m,t})_m$.

By the triangle inequality we have:
\begin{align}
\begin{aligned}\label{eq:three_terms_estimate1}
\left\vert\langle \mu_{s_m} - \mu, \psi\rangle\right\vert
\leqslant & \vert \mu_{s_m} - \mu_{s_m,t}\vert \vert \psi\vert_\infty 
+\min_{\nu\in \cM_t}\vert \int \psi d\mu_{s_m,t}-\int\psi  d\nu  \vert\\
+&\max_{\nu\in \cM_t}\vert \int \psi d\nu-\int\psi  d\mu \vert.
\end{aligned}
\end{align}
Note that the above $\min,\max$ are attained because $\cM_t(\psi)$ is compact. When $m\rightarrow\infty$, the second term on the right-hand side of \eqref{eq:three_terms_estimate1} goes to $0$ by definition of $\cM_t$. We let $m\rightarrow\infty$ and get:
\begin{align}
\begin{aligned}\label{eq:three_terms_estimate2}
\limsup_{m\rightarrow\infty}\vert\langle \mu_{s_m} - \mu, \psi\rangle\vert
&\leqslant  \limsup_{m\rightarrow\infty}\vert \mu_{s_m} - \mu_{s_m,t}\vert \vert \psi\vert_\infty \\
& + \max_{\nu\in \cM_t}\vert \int \psi d\nu-\int\psi  d\mu \vert.
\end{aligned}
\end{align}

First we deal with the second term on the right-hand side of \eqref{eq:three_terms_estimate2}. Let $\nu_t\in \cM_t$ be such that 
\[
\max_{\nu\in \cM_t}\vert \int \psi d\nu-\int\psi  d\mu \vert=\vert\int\psi d\mu-\int\psi d\nu_t\vert.
\]
Recall that support of $\nu_t$ is contained in $\geo X$ (see Remark \ref{rem:support_mut}) and that $\psi$ is continuous on $\geo X$ by the assumption. Proposition \ref{prop:double_limit_ordered} says that $\nu_t$ converges weakly to $\mu$ as $t\rightarrow 0$. Therefore,  
\begin{equation}\label{eq:three_terms_estimate3}
\lim_{t\rightarrow 0}\max_{\nu\in \cM_t}\vert \int \psi d\nu-\int\psi d\mu \vert=0.
\end{equation}

Secondly, as $\mu(P)=0$, we obtain, by Proposition \ref{prop:estimate_mus_must}, that
\[\lim_{t\rightarrow 0}\limsup_{m\rightarrow\infty} \vert\vert \mu_{s_m} - \mu_{s_m,t}\vert\vert=0.\]

Finally, since the left-hand side of \eqref{eq:three_terms_estimate2} does not depend on $t$, we obtain by \eqref{eq:three_terms_estimate2} that
\begin{align*}
\lim_{m\rightarrow\infty}\vert\langle \mu_{s_m} - \mu, \psi\rangle\vert=0.
\end{align*}
Thus, \eqref{eq:main_convergence1} holds and theorem follows. 
\end{proof}

\chapter{Flag manifolds and group actions}
\label{section_prelim}

\section{Basics of symmetric spaces}


We refer the reader to \cite{BGS} and \cite{Eberlein} for detailed treatment of symmetric spaces and their ideal boundaries. 

\medskip \index{$\cX$, a symmetric space}
Throughout the book, $\cX$ will denote a symmetric space of noncompact type associated with a semisimple Lie group $G$ (complex or real) 
commensurable to the isometry group of $\cX$ (i.e. $G$ acts on $\cX$ isometrically with finite kernel and the image of $G$ in the isometry group of $\cX$ is a finite index subgroup). Let $K$ denote the stabilizer of a point $o$ in $\cX$. Then $K$ is a maximal compact subgroup of $G$ and $\cX$ is $G$-equivariantly diffeomorphic to the quotient $G/K$. The defining feature of a symmetric space is the existence of {\em Cartan involutions}\index{Cartan involution} $s_x$ about any given point $x\in \cX$: The involution 
$s_x$ is an isometry of $\cX$ fixing $x$ and inducing the antipodal map on the tangent space $T_x\cX$.

\index{flat} \index{rank of a symmetric space}
A {\em flat} in $\cX$ is a totally geodesic subspace  isometric to a Euclidean space. The {\em rank} $r$ of $\cX$ is the maximal dimension of  maximal flats in $\cX$. The group $G$ acts transitively on the set of maximal flats in $\cX$. We fix a   {\em model maximal flat} $\Fmod\subset \cX$ containing $o$. \index{$\Fmod$, 
model maximal flat}

A {\em transvection}\index{transvection} in $G$ is the composition of two Cartan involutions $s_x, s_y$ of $\cX$ (fixing points $x, y\in \cX$ respectively). 
Geometrically speaking, these are translations along geodesics in $\cX$ 
(called {\em axes} of the transvection) with trivial {\em rotational component} around the axis. In the case of a nontrivial transvection  
$t=s_xs_y$ the unique complete geodesic through $x, y$ is an axis of $t$. If  
$\cX$ is a real-hyperbolic space $\bH^n$, transvections are exactly hyperbolic elements of the isometry group. Unless $\cX$ has rank $1$, transvections never have 
unique axis, the union of axes has dimension $\geqslant r$ and contains a maximal flat in $\cX$.   
The group  $A$ of transvections preserving $\Fmod$ 
is abelian, isomorphic to the $\bR^r$; it is the connected component of identity in the maximal real split torus in $G$.

\subsection{Visual boundary: Two topologies}

\index{$\geo \cX$, visual boundary of a symmetric space}
The {\em visual boundary} $\geo \cX$ (or the {\em ideal boundary}) of $\cX$ consists of equivalence classes of  $[\rho]$ of geodesic rays $\rho: \R_+\to \cX$ in $\cX$, where two rays are equivalent if and only if their images are at a finite Hausdorff distance from each other. One says that every ray $\rho$ representing $\xi=[\rho]$ is {\em asymptotic to $\rho$.} 
The visual boundary of $\cX$ has two natural topologies. The first one is the {\em visual topology}: \index{visual topology}
Every $\xi\in \geo \cX$ is represented by a unique unit speed geodesic ray emanating from $o$. Thus, there is a natural bijection between $\geo \cX$ and the unit sphere in the tangent space $T_o\cX$. The visual topology on $\geo \cX$ is the one making this bijection a homeomorphism. The natural $G$-action on $\geo \cX$ is continuous with respect to this topology. This topology extends to a {\em visual compactification} 
$\overline{\cX}=\cX\cup \geo \cX$: A sequence $(x_n)$ in $\cX$ converges to $\xi=[\rho]\in \geo \cX$ if 
 $$
 \lim_{n\to\infty} \angle_o(x_n, \rho(1))=0 \quad \hbox{and}\quad \lim_{n\to\infty} d(o,x_n)=\infty, 
 $$ 
 where $\rho(0)=o$. 

 The {\em Tits angle} $\tangle(\xi_1, \xi_2)$ between points $\xi_1=[\rho_1], \xi_2=[\rho_2]$ is defined as \index{$\tangle$, Tits angle}
$$
\sup_{x\in Y}  \angle_x(\rho_1(t), \rho_2(t)),
$$ 
 where the supremum is taken over all pairs of rays $\rho_1, \rho_2$ representing $\xi_1, \xi_2$ such that $\rho_1(0)=\rho_2(0)=x$. 
 Since $\cX$ is a symmetric space, there exists a flat $F\subset \cX$ such that $\xi_1, \xi_2$ are represented by rays whose images are contained in $F$. The supremum in the definition of $\tangle(\xi_1, \xi_2)$ is realized by   pairs of such rays.  The Tits angle defines the {\em Tits metric} on $\geo \cX$. This metric is invariant under the natural $G$-action on $\geo \cX$. 

 The second, {\em Tits topology}, \index{Tits topology}
  is the one defined by the Tits metric on $\geo \cX$. With respect to this topology,  $\geo \cX$ has the structure of a certain simplicial complex, the {\em thick spherical (Tits) building} $\tits \cX$, invariant under the action of $G$. Ideal boundaries of maximal flats in $\cX$ are called {\em apartments} in the building $\tits X$. \index{apartment}
  Dimension of the complex $\tits \cX$ is $r-1$; each apartment in $\tits \cX$ is an isometrically embedded copy of the unit sphere of dimension $r-1$.   Top-dimensional simplices in $\tits \cX$ ({\em facets} in the simplicial terminology)  are called {\em (spherical) chambers}. Codimension one faces of $\tits \cX$ are called  {\em panels}. 
  The group $G$ acts transitively on the set of chambers and the set of apartments.   \index{chamber}  
  Thickness of the building means that each  panel \index{panel}
  in $\tits \cX$ is shared by at least three (in fact, infinitely many) chambers.

\begin{assumption}
We will assume that the $G$-stabilizer of each facet in $\tits \cX$ is trivial. 
\end{assumption}  
  
\begin{rem}
1. This is where commensurability of $G$ and the full isometry of $\cX$ comes into play: The group of self-isometries of each facet is finite, isomorphic to the group of automorphisms of the Dynkin diagram of $\cX$. This finite group is sometimes nontrivial; for instance, for the Lie group $\SL(n+1, \C)$, this group is the dihedral group of order $2n$. 

2. All connected semisimple Lie groups $G$ satisfy this assumption. 
\end{rem}

  Below is an illustration of the above concepts in the case $G=\SL(3,\R)$, $K=\SO(3)$, $\cX$ can be (and will be) identified with the space of positive definite symmetric matrices of determinant $1$ (hence, is $5$-dimensional). The space $\cX$ has rank $2$. As a model flat $\Fmod$ one usually takes the subspace of diagonal matrices (with positive diagonal entries, of course). To see the Riemannian metric on $\Fmod$, take the logarithms of the corresponding diagonal matrices: The result is the hyperplane
$$
\{(t_1, t_2, t_3): t_1+t_2+t_3=0\}\subset \R^3. 
$$
The Riemannian metric on $\cX$ then becomes the restriction of the Euclidean metric.

The subgroup $A< G$ stabilizing the model flat is the multiplicative 
group of diagonal matrices with positive diagonal entries and unit determinant. 
Transvections in $\cX$ correspond to elements of $G$ which are diagonalizable matrices with positive eigenvalues.

The visual boundary of $\cX$ is homeomorphic to the $4$-dimensional sphere $S^4$. The spherical building $\tits X$ of $\cX$ is the (bipartite) incidence graph of real projective plane $\R P^2$. It has two types of vertices: One corresponding to points in $\R P^2$ and the other corresponds to lines. A point-vertex is connected by an edge to a line-vertex if and only if the point is incident to the line. The edges of this graph are chambers of   $\tits X$; panels are vertices. Each edge has length $\pi/3$ in Tits' metric. The building is thick: Every vertex is incident to continuum of lines and every line is incident to continuum of points. Cartan involutions do not belong to $G$; they swap point-vertices and line-vertices of $\tits \cX$.  An apartment in  $\tits X$ is a circuit of combinatorial length $6$ (and metric length $2\pi$). Geometrically, this corresponds to a triangle in the projective plane: Vertices of the triangle are the point-vertices of the apartment, while projective lines extending the edges of the triangle and the line-vertices.

One can think of edges in $\tits X$ as {\em projective flags} $(p,l)$, where $p$ is a point and $l$ is a line incident to the point $p$.

\subsection{Weyl chamber, Weyl group and simple roots}  
  
 We will fix a {\em model chamber} of $\geo \cX$, i.e. a  facet $\simod$ 
 of this spherical  building contained in the ideal boundary $\amod$ \index{$\amod$, model apartment}
 of the  model maximal flat $\Fmod\subset \cX$. \index{$\simod$, model chamber}

  The {\em (Euclidean) Weyl chamber} $\Delta$ of $\cX$ is the cone in $\Fmod$ with the tip $o$ over  \index{$\Delta$, the Weyl chamber} 
 $\simod$ (the union of geodesic rays emanating from $o$ and asymptotic to the points of $\simod$). The {\em Weyl group} $W$ of $\cX$ 
 \index{$W$, the Weyl group} 
 is the image of $K\cap \Stab_G(\Fmod)$ in the isometry group of the flat $\Fmod$. Then $\Delta$ is a fundamental domain of the $W$-action on $\Fmod$. 
 The element of $W$ sending $\Delta$ to $-\Delta$ is denoted $w_0$, this is the {\em longest element} of $W$. \index{$w_0$, the longest element of $W$} 
  \index{$\iota$, the opposition involution}
  Identifying $\Fmod$ with $\R^r$ (where $r$ is the rank of $\cX$), we get the {\em opposition involution} $\iota=-w_0$ preserving $\simod$.

The stabilizer of a facet $\sigma$ in $\tits \cX$ is a {\em minimal parabolic subgroup} of $G$; in the case of interest, when $G$ is a complex semisimple Lie group, these are exactly the {\em Borel subgroups} of $G$. \index{$B$, Borel subgroup} \index{$P_{\simod}$, minimal parabolic subgroup} 
Thus, the set of facets in $\tits \cX$ can be then identified with  $G/P_{\simod}$, where $P_{\simod}$ is the $G$-stabilizer of the model facet $\simod$. When equipped with the quotient topology, $G/P_{\simod}$ is the {\em full flag-manifold} $\cF$ of $G$. \index{$\cF$, flag-manifold}
Thus, we will be regarding facets (chambers) in $\tits \cX$ as (maximal) {\em flags}. We will fix a K\"ahler metric with the K\"ahler form $\omega$ on $\cF$ and let $\dF$ denote the Riemannian distance function associated with this K\"ahler metric. \index{$\omega$, K\"ahler form} 

We will use the notation $\taumod$ to denote  faces (possibly equal to the entire $\simod$) of the model chamber $\simod$. We will also use the notation 
$\pm\taumod$ for the faces $\taumod=+\taumod$ and $\iota\taumod=-\taumod$. Accordingly, we define parabolic subgroups $P_{\pm\taumod}$ as the 
$G$-stabilizers of the faces $\pm\taumod$. 


Recall that the $G$-stabilizer of the flat $\Fmod$ contains an abelian Lie subgroup of transvections $A$ isomorphic to 
$\mathbb R^r$; $\mathfrak{a}$ denotes its Lie algebra. We will be identifying $\mathfrak a$ and $\Fmod$. We will use the notation ${\mathfrak a}^+$ for the {\em interior} of the Weyl chamber $\Delta$. 
\index{$\mathfrak{a}$} \index{$\mathfrak{a}^+$}
The maximal torus $T$ in $G$ containing the abelian subgroup $A$ has the Lie algebra \index{${\mathfrak{h}}$} 
$$
{\mathfrak h}= {\mathfrak a}\otimes \C. 
$$
We set  \index{${\mathfrak{h}}^+$}
$$
{\mathfrak h}^+:= {\mathfrak a}^+\oplus i  {\mathfrak a}\subset {\mathfrak h}. 
$$
Then $T^+:= \exp({\mathfrak h}^+)$ is a subsemigroup in $T$.  \index{$T^+$} 
Geometrically speaking, $T^+$ consists of  those regular loxodromic elements of $G$ which move the base-point $o\in \Fmod$ to 
points in $\Delta\subset \Fmod$. 

Let $U< G$ denote the maximal unipotent subgroup of $G$ normalized by $T$ and corresponding to the positive chamber 
${\mathfrak a}^+\subset {\mathfrak a}$, i.e. $U< P_{\simod}$. 
Let $U_{\alpha_i}, i=1,...,r$, denote the  {\em root subgroups} of $U$ corresponding to the simple roots $\alpha_i$ and 
${\mathfrak u}_{\alpha_i}$ their Lie algebras. 

\medskip 
Let $\Phi$ denote the root system of ${\mathfrak g}$ relative to the maximal torus $T$; 
let $\Phi^+$ and $\Phi^-$ denote, respectively, the sets of positive and negative roots in $\Phi$; $\Phi=\Phi^+\sqcup \Phi^-$. We let 
$\Pi^+=\{\alpha_1,...,\alpha_r\}\subset \Phi^+$ denote the set of {\em simple} roots of $G$, the ones defining the chamber $\Delta$, 
\index{$\Phi^+$, the set of positive roots} 
\index{simple roots}  
\index{$\Phi^-$, the set of negative roots} 
$$
\Delta=\{y\in {\mathfrak a}: \alpha_i(y)\geqslant 0, i=1,...,r\}. 
$$
The set of negatives of the simple roots is denoted $\Pi^-$.

Simple roots are in bijective correspondence with generating reflections of the Weyl group $W$; they are called {\em simple reflections}, 
\index{simple reflection}
reflections in the facets of $\Delta$. The length $|w|$ of an element $w\in W$ is its word length with respect to the generating set. The 
element $w_0\in W$ is the unique longest element of $W$ with respect to this word metric. \index{$\vert w\vert$, the word-length}

\medskip 
Below is an illustration of the above concepts in the case $G=\SL(3,\R)$. As before, we will identify the model flat $\Fmod$ with the hyparplane 
$$
\{(t_1, t_2, t_3): t_1+t_2+t_3=0\}\subset \R^3. 
$$
As the positive chamber $\Delta$ one usually takes the one given by the inequalities 
$t_1\geqslant t_2\geqslant t_3, i=1, 2, 3$. The roots are the linear functions $t_i-t_j$, $i\ne j$.  The simple roots are $\al_1=t_1-t_2$ and $\al_2=t_2-t_3$. 
The Weyl group $W$ is the permutation group $S_3$ generated by linear orthogonal reflections in the hyperplanes $t_1=t_2$ and $t_2=t_3$. The longest element $w_0$ is the reflection in the hyperplane $t_1=t_3$. The standard minimal parabolic subgroup $P_{\simod}=B$ 
is Borel, it is the group of (nonstrictly) upper triangular matrices in $G$. The unipotent radical $U< B$ is the subgroup of strictly upper-triangular matrices. The model 
spherical chamber fixed by $B$ has the point-vertex $\Span {\mathbf e}_1$ and the line-vertex is $\Span \{ {\mathbf e}_1, {\mathbf e}_2\}$.

\subsection{Opposite flags}

Two points $\xi, \eta$ in $\geo X$ are called {\em opposite}\index{opposite points}
 if $\tangle(\xi, \eta)=\pi$, equivalently, if there exists a geodesic $c$ in $X$ whose opposite subrays are asymptotic to $\xi$ and $\eta$ respectively. Equivalently, there exists a Cartan involution of $X$ swapping $\xi$ and $\eta$. Two simplices $\tau, \hat\tau$ in $\tits X$ are 
 \index{opposite simplices}\index{antipodal simplices}
 {\em opposite} (or, {\em antipodal}) if and only if they contain opposite {\em generic} points in $\geo X$. (A point in a simplex $\tau$ is {\em generic} if it does not belong to any proper face of $\tau$.) Two simplices in $\tits \cX$ are opposite  if and only if they are swapped by a Cartan involution of $\cX$. 

Given a flag $\tau\in \mathcal F$, we let $\Opp(\tau)$ denote the corresponding {\em open Schubert cell} in $\mathcal F$, i.e. 
the subset of $\cF$ consisting of flags opposite (antipodal) to $\tau$.  \index{$\Opp(\tau)$, the open Schubert cell}
 




To illustrate these concepts in the case of $G=\SL(3,\R)$, if $f_i=(p_i,l_i), i=1,2$, are projective flags representing chambers $\sigma_i, i=1, 2$, then $\sigma_1$ is opposite to $\sigma_2$ if and only if the flags $f_1, f_2$ are {\em transversal}, which in this example means that $p_i\notin l_{3-i}, i=1, 2$. Let's identify $\R P^2$ with the projective compactification of the affine plane $\R^2$. Consider the flag $f=(p,l)$, where $l$ is the line at infinity and $p$ corresponds to the vertical direction in the plane  
$\R^2$. Then the open Schubert cell $\Opp(f)$ consists of all flags $(q,m)$, such that $q$ is a point in the affine plane $\R^2$ and $m$ is an incident line which is not vertical. To parameterize $\Opp(f)$ via $\R^3$, we just use the location of $q$ in $\R^2$ and the slope of the line $m$ (the slope is finite since the line is not vertical).

\subsection{Cartan projection and Lyapunov spectrum}\label{sec:Cartan and Lyapunov}

 For each point $x\in \cX$ one defines the {\em $\Delta$-valued distance} 
\index{$d_\Delta$, the $\Delta$-valued distance}
 $d_\Delta(o,x)$ as the unique point of intersection $Kx\cap \Delta$. 
This definition extends to general pairs of points in $X$ by $G$-invariance. Then the Riemannian distance $d(x,y)$ between 
two points in $\cX$ is the Euclidean norm of the vector $d_\Delta(x,y)$. 

We have the {\em Cartan projection} $\mu: G\to \Delta$ defined by \index{$\mu(g)$, the Cartan projection} 
$$
\mu(g)= d_\Delta(o, go) 
$$
and its components $\mu_i= \alpha_i\circ \mu$. Suppose that $g$ preserves a maximal flat $F\subset X$ and acts on it as a translation. We will refer to such elements of $G$ as {\em loxodromic}. \index{loxodromic isometry of a symmetric space}
The {\em Lyapunov spectrum} $\lambda(g)$ of a loxodromic element $g\in G$ is defined as the vector \index{$\lambda(g)$, Lyapunov spectrum}
$$
\lambda(g)= d_\Delta(x, gx), x\in F.  
$$
Accordingly, we set $\lambda_i:= \alpha_i\circ \lambda$, the {\em $i$-th component of the Lyapunov spectrum}. One defines the Lyapunov spectrum for general elements $g\in G$ using the {\em Jordan decomposition} of $g$, but we do not need this. Note that, unlike the Cartan projection, the Lyapunov spectrum is invariant under conjugation. 

For instance, if $G=\SL(n,\C)$, $\mu(g)$ can be defined as follows. Let $\sigma_1(g)\geqslant ...\geqslant \sigma_n(g)$ 
denote the singular values of a matrix $g\in G$, arranged in the descending order. Then $\mu_i(g)= \log \sigma_i(g)- \log \sigma_{i+1}(g)$, $i=1,...,n-1$. Similarly, $\la_i(g)$ is the difference between logarithms of absolute values of 
the eigenvalues of $g$ (again arranged in the descending order).

Here, as usual, we ignore a uniform normalization constant in the definition of a Riemannian metric on the symmetric space of $\SL(n,\C)$.

We now return to the discussion of general semisimple Lie groups. A loxodromic element $g\in G$ is called {\em regular} if the vector $\lambda(g)$ is regular, \index{regular loxodromic isometry}
i.e. belongs to the interior of $\Delta$. Each regular loxodromic element $g$ has one attractive and one repelling fixed point in $\cF$ (respectively, $x^+, x^-$). We have uniform convergence on compacts 
\begin{align}\label{eq:polar action} 
\lim_{m\to\infty} g^m(x) = x^+, x\in \Opp(x^-), \\ \lim_{m\to-\infty} g^m(x) = x^-, x\in \Opp(x^+). 
\end{align}

Then every loxodromic element $g\in G$ stabilizing $\Fmod$ (i.e. an element of the torus $T$) 
normalizes each  $U_{\alpha_i}$ and acts on ${\mathfrak u}_{\alpha_i}$ 
(via the adjoint representation) with the eigenvalue $\exp(\la_i(g))$.  Thus, regularity of a loxodromic element means that 
$\exp(\la_i(g))\ne 1$ for all $i=1,...,r$. 

Geometrically speaking, regular loxodromic elements of $G$ can be described as follows: Each loxodromic element preserves a {\em regular geodesic} in $\cX$ (a geodesic contained in a unique maximal flat), acts on this geodesic as a nontrivial translation, but also is allowed to ``rotate'' around this geodesic (unlike a transvection which cannot rotate). This generalizes the familiar behavior of loxodromic elements of $\SL(2,\C)$.

We will discuss the dynamics of regular loxodromic elements of $G$ in more detail in Section \ref{summary_torus_action}.



\subsection{Relative position} \label{section_relative_position_bruhat}

Given two points $x, y$ in $\cF$, we think of them as chambers in $\tits \cX$. Consider an apartment $a$ of $\tits \cX$ containing both $x, y$. Since $G$ acts transitively on the set of apartments in $\cX$, there exists $g\in G$ sending $a$ to $\amod$ and sending chambers $x, y$ to 
chambers $x', y'$ in $\amod$. Since the Weyl group $W$ acts transitively on the set of chambers contained in one apartment, there is a unique $w\in W$ such that $w(x')=y'$. 
We define the {\em relative position} of $y$ with respect to $x$ as \index{$\pos(x,y)$, the relative position}
\[\operatorname{pos}(x,y)=w.\]
The relative position is independent of the chosen apartment and of $g$ (\cite{KLP18}) and we, thus, obtain 
the {\em relative position function} 
$$
\operatorname{pos} : \cF^2\to W. 
$$
The relative position is $G$-invariant, i.e.\ $\operatorname{pos}(x,x')=\operatorname{pos}(gx,gx')$ for any $g\in G$; furthermore, it is a complete $G$-invariant 
of pairs: $(x, y)$ and $(x', y')$ belong to the same $G$-orbit if and only if $\operatorname{pos}(x,y)=\operatorname{pos}(x',y')$.

Relative position is not (in general) symmetric, instead it is {\em antisymmetric}: 
$$
\operatorname{pos}(x,y)= \left(\operatorname{pos}(y,x)\right)^{-1}. 
$$
Nevertheless, the relative position function can be regarded as a {\em $W$-valued distance} on $\cF$ as it satisfies a weak form of the triangle inequality with respect to the Bruhat order, see \cite{Ronan}.  It becomes a genuine distance function if we compose it with the word-length on $W$:
$$
D(x,y)= |\operatorname{pos}(x,y)|.  
$$
\index{$D(x,y)$, distance in the Weyl group}
The diameter of $\cF$ with respect to $D$ equals $|w_0|$, the diameter of $W$ with respect to the word-metric.  Geometrically, $D$ is the length of the shortest {\em gallery} in $\tits \cX$ connecting $x$ and $y$, provided that $\tits \cX$ is connected, i.e. $\cX$ has rank $\geqslant 2$. (If $\cX$ has rank one then $\tits X=\cF$ is discrete and $D$ is the discrete metric on $\cF$.) Such shortest gallery is always contained in an apartment $a\subset \tits \cX$ (the apartment $a$ is unique if $x, y$ are antipodal.) The function $D$ satisfies the triangle inequality because of the existence of a retraction $\rho: \tits \cX\to a$ to any apartment $a\subset \cX$, which is a simplicial map fixing $a$ pointwise. In particular, if $z\in \cF$ is a chamber such that 
$$
D(x,z)+ D(z,y)= D(y,x),
$$
then $z$ belongs to an apartment containing $x$ and $y$. We refer the reader to \cite[Chapter IV]{Brown} or \cite[Chapter 3]{Ronan} for details.

\section{Bruhat order and thickenings}\label{sec:Thickenings}

\subsection{Thickenings in the Weyl group}

The Weyl group $W$ has a certain partial order, called {\em (strong) Bruhat order}.  \index{Bruhat order $\le$}
Let $S$ denote the set of simple reflections generating $W$. A word in these generators is their product. 
(The empty word is allowed, it corresponds to $1\in W$.) A word 
$w=s_{1}...s_{l}$, $s_i\in S, i=1,...,l$, is said to be {\em reduced} if it is the shortest word representing $w$. 
Given a reduced word $w=s_{1}...s_{l}$, a {\em subword} of $w$ is a product of the form
$$
w'=s_{i_1}...s_{i_k}, i_1< i_2< ... < i_k. 
$$
Then the Bruhat order is defined by: $u\le v$ if some (equivalently, every) reduced word $s_{1}...s_{l}$ 
representing $v$ contains subword representing $u$. 
Multiplication by $w_0$ reverses the order. Furthermore, $u\le v$ if and only if $u^{-1}\le v^{-1}$. 

We refer the reader to \cite{BB} for a more detailed discussion and alternative definitions of the Bruhat order. 

For $w\in W$ define its {\em dual} $w^\vee:=w_0w\in W$. \index{$w^\vee$, 
dual element of the Weyl group} Then $|w^\vee|=|w_0|-|w|$.

A {\em thickening} $\Th$ in $W$ \index{$\Th$, a thickening} 
is a subset $\Th\subset W$ such that if $v\in \Th$ then $u\in \Th$ for any $u\leq v$. If we view $(W, \le)$ as a poset, then 
thickenings in $W$ are {\em ideals} in this poset:

Given a poset $(\P, \le)$, an ideal in $\mathbb P$ is a subset $I\subset \P$ such that if $p\in I$ and $q\le p$, then $q\in I$. 
 
Therefore for any $w\in W$ the subset $\Th_w$ defined by \index{$\Th_w$, the $w$-thickening} 
\[\Th_w=\{u\in W: u\leq w\}\]
is a thickening, a {\em principal ideal} in the poset terminology. A general ideal in $(W,\le)$ is a union of principal ideals. 
Note that for every $u\in \Th_w$ we have $|u|<|w|$ unless $u=w$. 
Given $k\geqslant 0$ we let $\Th^k$ denote the thickening \index{$\Th^k$, $k$-dimensional thickening}
$$
\bigcup_{w\in W, \vert w\vert\leq k} \Th_w. 
$$

For each thickening $\Th$ we define the {\em dual thickening} $\Th^\vee$ by: 
$$
\Th^\vee_w= \Th_{w^\vee}, 
$$
and if $\Th$ is the union of thickenings $\Th_w$ for $w$ belonging to some subset $V\subset W$, then 
$$
\Th^\vee= \bigcup_{w\in V} \Th^\vee_w. 
$$
The {\em inverse thickening} of a thickening $\Th$ is $\widehat\Th=\{w^{-1}: w\in \Th\}$. The fact that inversion preserves the Bruhat order implies that 
$\widehat\Th$ is again a thickening. For example,
$$
\hat\Th_{w}= \Th_{w^{-1}}. 
$$

A thickening $\Th$ is called {\em slim} if $\Th\cap w_0 \Th=\emptyset$ and {\em fat} if 
$\Th\cup w_0 \Th=W$. A thickening is called {\em balanced} if it is both slim and fat. \index{balanced thickening}

\begin{lem}\label{lem:slim}
1. Set $n:=|w_0|$. If $k=|w|< n/2$, then $\Th_w$ is slim.

2. For every $w\in W$ either $\Th_w$ or $\Th_{w^\vee}$ is slim. 
\end{lem}
\begin{proof} 1. Suppose that $u\in \Th_w \cap w_0 \Th_w$. Then the distances (in the word metric) from $u$ to $1$ and to $w_0$ are $\le k$. By the triangle inequality 
the distance from $1$ to $w_0$ is at most $2k<n$. This contradicts the fact that $n=|w_0|$. 

2. Setting again $k=|w|$, we get $|w^\vee|=n-k$. If $\min(k, n-k)< n/2$, we are done by Part 1. Otherwise, $k=n-k=n/2$. By the same argument as in Part 1, the only 
element that can belong to $\Th_w\cap w_0 \Th_w$ is $w$: The length of all other elements of $\Th_w$ is $<k= n/2$. But then $w=w_0w$, i.e. $w_0=1$, 
which is a contradiction. 
\end{proof}

\subsection{Thickenings in the flag manifold}\label{sec:Thickenings in the flag manifold}

From now on, $G$ is a complex semisimple Lie group with Borel subgroup $B$ and $\cF=G/B$ is the full flag-manifold of $G$ and $n$ is the complex dimension of $\cF$. 

Let $x\in \cF$, we think of $x$ as a chamber in $\partial_{Tits}\cX$. For $w\in W$, the set 
\[S_w(x)=\{x'\in \cF: \operatorname{pos}(x,x')=w\}\] 
\index{$S_w(x)$, Schubert cell} 
is a Zariski open subset of $\cF$ known as a Schubert cell. For example, when $w=w_0$ is the longest element, 
then $S_{w_0}(x)=\Opp(x)$ is the set of all flags opposite to $x$. 

 We define the \textit{thickening} $\Th_w(x)$ to be
\[\Th_w(x)=\{x'\in \cF: \operatorname{pos}(x,x')\leq w\}.\]
We list here the main features on the thickenings needed in our book (see \cite{brion_lectures, KLP18} for the proofs).
\begin{prop} \label{prop_listing_thickenings}
	For any point $x \in \cF$, the following properties are satisfied. 
	\begin{enumerate}
\item[(i)] For any $w\in W$,  the Zariski closure of $S_w(x)$ in $\cF$ equals $\Th_w(x)$. 
\item[(ii)] For any $w , w' \in W$ with $w\leqslant w'$, one has $\Th_w(x) \subset \Th_{w'}(x)$. 
\item[(iii)] Each $\Th_w(x)$ is an irreducible subvariety of $\cF$ of complex dimension $|w|$.  
	\end{enumerate}
\end{prop}


\begin{rem}\label{rem:local Schubert} 
Much is known about geometry of Schubert varieties $X_w=\Th_w(x)$, but their detailed description is quite complicated. 
These varieties are singular if and only if they are singular at the point $x$. 
For some $w$'s the varieties are indeed singular, while for other $w$'s they are smooth. 
The precise description of Zariski tangent spaces, tangent cones and local geometry at singular points is possible but is complicated. See for instance \cite{BL}, \cite{CK}, \cite{Ulfarsson-Woo}. In the case when $G=\SL(m,\C)$, Lakshmibai and Sandhya (in \cite{MR1051089}) proved that $X_w$ is 
smooth if and only if the permutation $w$ avoids the patterns  $3412$ and $4231$. General smoothness criteria can be found in 
\cite{BL} and \cite{CK}; one can find there also descriptions of tangent spaces to Schubert varieties. For instance, in the case when $G=\SL(m,\C)$, 
$G=\Sp(2m,\C)$, the tangent space $T_{x_+}\Th_w(x_+)$ has a particularly simple description: It is a subspace of ${\mathfrak u}^-$ (the Lie algebra of negative unipotent subgroup $U^-< G$) spanned by the root subspaces ${\mathfrak g}_{-\alpha}$, where $\alpha\in \Phi^+$ are positive  roots such that the corresponding reflections $s_\alpha\in W$ satisfy the inequality $s_\alpha\le w$; see \cite[Theorem 5.3.1]{BL}.  More recently, Ulfarsson and Woo in \cite{Ulfarsson-Woo} 
proved that $X_w$ is a local complete intersection if and only if $w$ avoids the patterns 
$53241$, $52341$, $52431$, $35142$, $42513$, and $426153$. There are several constructions of $T$-equivariant resolutions of singularities of 
Schubert varieties, the most common is the {\em Bott--Samelson resolution}, see e.g. \cite{BL}. 
\end{rem}

\begin{example}
Suppose that $w=s_\alpha$ is a simple reflection. Then the Schubert variety $X_w$ is smooth, isomorphic to the projective line $\P^1$.  Let $U_{\pm\alpha}=\exp({\mathfrak g}_{\pm\alpha})$ denote the root subgroups  of $G$ corresponding to the roots $\pm\alpha$ (where $\alpha$ is a simple root). These are 1-dimensional unipotent subgroups of $G$. They fix, respectively, points  $x_+:=\simod$ and $x_w=w(x_{+})$, and we have 
$$
S_w(x_+)= U_{\alpha}(x_w), {S_w(x_w)}= U_{-\alpha}(x_+), 
$$
open Schubert cells isomorphic to $\mathbb C$. Intersection $S_w(x_+)\cap S_w(x_w)$ is the complement 
$$
\Th_w(x_+)\setminus \{x_+, x_w\}. 
$$
Thus, $X_w=\Th_w(x_+)$ is smooth at $x_w$, hence, is smooth. It is isomorphic to $\P^1$ since it contains a dense subset isomorphic to the affine line $\mathbb C$. Lastly, note that the subalgebra of $\mathfrak{g}$ generated by $\mathfrak{g}_{\pm\alpha}$ is isomorphic to $sl(2,\mathbb C)$; the 
maximal torus of $\exp(sl(2,\mathbb C))$ fixing $x_+, x_w$ acts transitively on $S_w(x_+)\cap S_w(\textcolor{red}{x_w})\cong \mathbb C^*$. 
(The action might have nontrivial kernel $\cong \mathbb Z_2$.) 
\end{example}

More generally, if $\Th\subset W$ is a thickening, then we define \index{$\Th(x)$, thickening of $x\in \cF$} 
\[\Th(x)=\{x'\in \cF: \operatorname{pos}(x,x')\in \Th\}\]
and for a subset $\Lambda\subset \cF$ we define
\[\Th(\Lambda)=\cup_{x\in \Lambda}\Th(x).\]

For $m\in \bN$, we define 
\[\Th^m(x)=\{x'\in \cF: D(x,x')= |\operatorname{pos}(x,x')|\leq m\}\] and similarly for $\Th^m(\Lambda)$. 
Note that the thickenings $\Th(x),\Th^m(x)$ are unions of thickenings of the form $\Th_w(x)$.

We will denote by $ \Th_{<w}(x)$ the complement $\Th_w(x)\setminus S_w(x)$. 
It is a thickening of $x$ itself. For example, 
when $w=w_0$ is the longest element, $\Th_{w_0}(x)=\cF$ and $\Th_{<w_0}(x)=\Th^{n-1}(x)$. In general, we have
\begin{equation}
\label{thickenings_boundary_inclusion}
 \Th_{<w}(x)\subset \Th^{\vert w \vert -1}(x).
\end{equation}

Since the relative position function is a complete $G$-invariant for pairs in $\cF$, 
the $G$-stabilizer $B_x=\Stab_G(x)$ acts transitively on $S_w(x)$. The action is not simply-transitive, the stabilizer $\Stab_{B_x}(y)$ of $y\in S_w$ in $B_x$ 
is the intersection of two Borel subgroups, $B_x\cap B_y$. If $w=w_0$, i.e. $y$ is opposite to $x$, the intersection is the maximal torus $T_a$ \index{$T_a$, the maximal torus stabilizing the apartment $a$} 
in $G$ fixing pointwise the unique apartment $a$ in $\cF$ containing $x$ and $y$. Hence, the quotient group $B_x/T_a$ acts simply transitively on 
$\Opp(x)$. The group $B_x$  splits as a semidirect product of $T_a$ and the unipotent radical in $B_x$; the latter acts simply-transitively on  $\Opp(x)$.  
However, the unipotent radical in $B_x$ acts transitively but not simply-transitively on other Schubert cells $S_w, |w|<n$: The stabilizer of $y\in S_w$ 
is the unipotent radical of the solvable group $B_x\cap B_y$. We will discuss this in more detail in Section \ref{section_lie_alg}. 

\medskip 
Observe that  for $g\in G, x\in \cF$, $g \Th(x)=\Th(gx)$, for every thickening $\Th\subset W$.

The following lemma was proven in \cite[Lemma 8.18]{Kapovich-Leeb-finsler} 
in greater generality; we include a proof for the sake of completeness:

\begin{lem}\label{lem:compactness}
1. Suppose that $(x_i)$ is a sequence in $\cF$ converging to $x$. Then for every thickening $\Th\subset W$, the sequence of thickenings $\Th(x_i)$ converges 
to $\Th(x)$ with respect to the Hausdorff metric. 

2. If $C\subset \cF$ is a compact subset and $\Th\subset W$ is a thickening, then $\Th(C)$ is again compact in $\cF$. 
\end{lem}
\begin{proof} 1. There exists a sequence $g_i\in K$ such that $g_i(x_1)=x_i$ for all $i\in \bN$ (since the maximal compact subgroup $K< G$ 
acts transitively on $\cF$). In view of compactness of the subgroup $K< G$, without loss of generality, we may assume that the sequence $g_i$ converges to some $g_0\in K$. Hence, the sequence $g_i(\Th(x_1))=\Th(x_i)$ converges to $\Th(x)=g_0(\Th(x_1))$. 

2.  Take a sequence $x_k\in C$ converging to $x\in C$. Then, by Part 1, $\Th(x_k)\to \Th(x)\subset \Th(C)$. 
It follows that $\Th(C)\subset \cF$ is closed; compactness of $\cF$ implies compactness of 
$\Th(C)$. 
\end{proof}

\begin{cor}
The function $\pos$ is lower semicontinuous on $\cF$. 
\end{cor}
\begin{proof} Suppose that $x_i, y_i\in \cF$ are such that $\pos(x_i, y_i)=w\in W$ and $\lim_{i\to\infty} x_i=x$,  $\lim_{i\to\infty} y_i=y$. Hence, $y_i\in \Th_w(x_i)$ and, by the lemma, $y\in \Th_w(x)$, i.e. $\pos(x,y)\le w$. 
\end{proof}

\subsection{Richardson varieties}\label{sec:Richardson} 

Let $x_\pm$ be opposite points in $\cF$, $\bx=(x_+,x_-)$, $v, w\in W$ are such that $|v|+|w|\geqslant n$. Intersections 
$\Th_{vw}(\bx):=\Th_v(x_+)\cap \Th_w(x_-)$ are called 
{\em Richardson varieties} in $\cF$ and the intersections of corresponding Schubert cells $S_{vw}(\bx):= S_v(x_+)\cap S_w(x_-)$ are called {\em open 
 Richardson varieties}. All open Richardson varieties are smooth and have dimension $|v|+|w|-n$. On the other hand, Richardson varieties could be singular. 
 For instance, if $|v|+|w|=n+1$, then $S_{vw}(\bx)$ is isomorphic to $\bC^*$ and $\Th_{vw}(\bx)$ is isomorphic  to $\P^1$. Up to the $G$-action, 
 the varieties $S_{vw}(\bx)$, $\Th_{vw}(\bx)$ depend only on $(v,w)$ and, for this reason, are usually denoted $\mathring{X}_{vw}, X_{vw}$ respectively. 


\begin{lem}  \label{lem:intersection of opposite thickenings} 
Pick $v, w \in W $, two elements of the Weyl group and let $x, y$ be two antipodal points in $\cF$. 
 Then $X_{vw}= \Th_v(x) \cap \Th_w(y)$ is nonempty if and only if $ w_0 w \le v$ for the Bruhat order. 
 Moreover, when nonempty, the intersection is proper and is  a projective  subvariety of dimension $|w| + |v| - n$.  
\end{lem} 
\begin{proof}  

1. We give a building-theoretic proof similar to the one in Proposition \ref{prop:singleton}. 
Assume first that $ \Th_v(x) \cap \Th_w(y)$ is nonempty. 
Let $a\subset \tits \cX$ be the unique apartment containing $x$ and $y$. Suppose that 
there is $z\in \Th_v(x) \cap \Th_w(y)$. Let $\rho: \tits \cX\to a$ denote the simplicial retraction which restricts to an isometry on an apartment 
in $\cX$ containing $y$ and $z$ (see \cite{Ronan} for the existence of such retractions). 
Set $\bar{z}:= \rho(z)$. Then $\operatorname{pos}(y, z)= \operatorname{pos}(y, \bar z)$. The retraction $\rho$ maps the minimal gallery in $\tits X$ to 
a {\em folded gallery} in $a$ connecting $x$ and $\bar z$. Thus, 
$$
\operatorname{pos}(x, \bar z)\le \operatorname{pos}(x,z),
$$
because the {\em folding order} coincides with the {\em Bruhat order} on $W$ (see \cite{KLP18A}). Since $y$ is opposite to $x$, 
$$
\operatorname{pos}(y, \bar z)= w_0 \operatorname{pos}(x, \bar z)$$
We also have 
$$
 \operatorname{pos}(x,z)\le v, \quad \operatorname{pos}(y,z)\le w.
$$
Thus, by transitivity of the Bruhat order:
$$
w_0 \operatorname{pos}(y,z)\le \operatorname{pos}(x,z) \le v.
$$
Since $\operatorname{pos}(y,z)\le w$ and the multiplication by $w_0$ reverses the Bruhat order,
$$
w_0 w\le w_0 \operatorname{pos}(y,z). 
$$
Hence, as claimed, $ w_0 w \le v$.  The proof of the converse follows the same chain of inequalities. 
Assume that $ w_0 w \le v$. Take an apartment $a$ as above and a chamber $\bar{z}\in a$ on a minimal gallery between $x$ and $y$, 
such that $\operatorname{pos}(y,\bar z)=w$. Then $\operatorname{pos}(x, \bar z)=w_0w\le v$. Then, by the definition of thickenings $\Th_v, \Th_w$,
$$
\bar z\in \Th_v(y), \quad \bar z\in \Th_v(x).
$$ 
Hence $\Th_v(y)\cap \Th_v(x)\ne \emptyset$. 

\medskip 
2. The group $G$ acts transitively on $\cF$. Therefore, by Kleiman Transversality Theorem (see \cite{Kleiman74}), for generic $g\in G$ the intersection 
$$
\Th_v(x) \cap g \Th_w(y) = \Th_v(x) \cap \Th_w(gy)$$
is transversal and, hence, has codimension which is the sum of codimensions of $\Th_v(x), g \Th_w(y)$. 
For generic $g$, $gy\in \Opp(x)$, since $\Opp(x)$ is Zariski open and dense in $\cF$. The parabolic subgroup $P_x< G$ (the stabilizer of $x$ in $G$) 
acts transitively on $\Opp(x)$. Hence, there is $h\in P_x$ such that $h(x)=x, hg(y)=y$. Since $h$ carries the pair $(\Th_v(x), \Th_w(gy))$ to 
$(\Th_v(x), \Th_w(y))$, the intersection $\Th_v(x) \cap \Th_w(y)$ is also transversal. 
\end{proof}

We refer the reader to  \cite{speyer2024richardsonvarietiesprojectedrichardson} for  a detailed discussion of Richardson varieties.

\subsection{Homology}\label{sec:homology}


As $g(\Th(x))=\Th(x')$ for $x'=gx, g\in G$, any 
$\Th_w(x')$ has the same homology class as $\Th_w(x)$ in 
$\coh_*(\cF; \bZ)= \coh_*(\cF)$. We will denote it $[\Th_w]=[X_w]$. It is known that the classes $[\Th_w], w\in W$, form a free basis of $\operatorname{H}_*(\cF)$ (see \cite{BGG73}).
Furthermore, for any $w$, $[\Th_w]$ is the Poincar\'e dual of $[\Th_{w^\vee}]$. More precisely, if 
$PD([\Th_w])\in \coh^{2n-2k}(\cF)$ denotes the Poincar\'e dual of $[\Th_w]\in \coh_{2k}(\cF)$, $|w|=k$, 
then 
$$
\langle PD([\Th_w]) \cup PD([\Th_v]), [\cF]\rangle = \delta_{w, v^\vee}, w, v\in W,
$$
i.e. 
$$
\langle PD([\Th_w]) \cup PD([\Th_v]), [\cF]\rangle = 0, \hbox{\ if\ } v\ne  w^\vee,
$$
and 
$$
\langle PD([\Th_w]) \cup PD([\Th_{w^\vee}]), [\cF]\rangle = 1. 
$$

\section{Unipotent orbits and product decompositions of thickenings}
\label{section_lie_alg}

In this section we describe in detail thickenings of the type $\Th_w(x)$ as closures of unipotent orbits. Given a pair $\bx=(x^+,x^-)$ of  
opposite points in $\cF$ we consider the apartment $a=a(x^+,x^-)=a_{\bx}\subset \tits \cX$ 
containing both $x^\pm$. For $w\in W$ we will also define opposite chambers $x_{w}, x_{w^\vee}$ in this apartment. 
Given the unipotent radical $U_{x_w}$ stabilizing the point $x_{w^\vee}$ we define maximal 
unipotent subgroups $U_{x_w}^\pm$ in the stabilizers of $\{x_{w^\vee}, x^\pm\}$ 
and get the {\em knit product decomposition}
$$
U_{x_w}= U_{x_w}^+ U_{x_w}^-. 
$$
We will identify the Schubert cells $S_w(x^+), S_{w^\vee}(x^-)$ as unipotent orbits $U_{x_w}^+ x_w, U_{x_w}^- x_w$. 
Accordingly, we obtain a decomposition of the open Schubert cell $\Opp(x_{w^\vee})$ as the direct 
product of the cells:  $S_w(x^+)\times S_{w^\vee}(x^-)$. This decomposition is invariant under the action of the maximal torus $T_a$.  
 In the following section (\S \ref{section_torus_action}) 
we define invariant coordinates on 
$\Opp(x_{w^\vee})$; the action of each element $g$ in this torus is diagonal with respect to these coordinates 
with diagonal entries corresponding to exponentials of values of certain roots on $\log(g)$. Real parts of these coordinates are the Lyapunov exponents of $g$. We will see that the cells $S_w(x^+), S_{w^\vee}(x^-)$ are stable/unstable manifolds of the action of $g$ in $\Opp(x_{w^\vee})$ at the fixed point $x_w$.

\subsection{Structure of thickenings of chambers in the model apartment}\label{sec:Structure of thickenings}

Recall from Section \ref{section_relative_position_bruhat} that the position function $\pos$ 
is $G$-invariant. Therefore, in this section we work with thickenings of 
chambers in $\cF$ which belong to the model apartment $\amod$;  thickenings at other points can be 
described accordingly via the $G$-action (this will be done in the next section).

Recall that the Weyl group $W$ is the quotient $N_K(T)/(K\cap T)\cong N_G(T)/T$, the quotient of the normalizer of $T$ by the maximal torus $T< G$. 
 Accordingly, for each element $w \in W$, we choose its preimage $k_w\in N_K(T)$. 
 Note that, in general, one cannot find a homomorphism $W\to N_K(T)$ of the form $w\mapsto k_w$ (this can be seen already in the case of $G=\SL(2,\C)$). Furthermore, 
 the objects involving $k_w$, that appear later in this section,  will not depend on the choices of $k_w$. 
 
 \begin{rem}
 J.~Tits, \cite{Tits66}, proved that there is a subgroup $\tilde{W}< N_K(T)$ which projects onto $W$ and is naturally isomorphic to a central extension
 $$
1\to \bZ_2^r \to \tilde{W}\to W\to 1
 $$
 \end{rem}

Recall that we identify $\cF$ with the set of chambers in $\partial_{Tits}\cX$. The model apartment $\amod$ in $\partial_{Tits}\cX$ is 
determined by $T$ (and determines $T$) and its model chamber $\simod$ corresponds to the trivial coset $B$ in $\cF=G/B$. 
Then the cosets $k_wB, w\in W$ correspond to other chambers in the model apartment. 

We denote the coset $k_wB$ by $\sigma_w$ when we think of it as a point on the manifold $\cF$; when $w=e$ or $w=w_0$ we use the special notation $B=\sigma^+=\simod$ and $k_{w_0}B=\sigma^-$. We will be using the notation $B^+=B$ (the positive Borel subgroup) and $B^-= k_{w_0}B k^{-1}_{w_0}$ (the negative Borel subgroup). Thus, $\sigma_w= k_w \sigma^+$. 
 
 
 The fixed-point set of the $T$-action on $\cF$ is the set of chambers in the model apartment $\amod$, i.e. 
 $$
 W \simod = \{\sigma_w: w\in W\}. 
 $$
The $G$-stabilizer of the $\sigma_w$ is the conjugate Borel subgroup $B_{w}=k_{w} B k_{w}^{-1}$.  
 

The positive Borel subgroup $B$ gives rise to the standard root decomposition of the 
Lie algebra $\mathfrak{g}$ of $G$:
\[\mathfrak{g}=\mathfrak{h}\oplus \left(\bigoplus_{\alpha\in \Phi^+}{\mathfrak g}_{\alpha}\right)\oplus \left(\bigoplus_{\alpha\in \Phi^-}{\mathfrak g}_{\alpha} \right)\]
where 
and the ${\mathfrak g}_{\alpha}$'s are {\em root subspaces}, Lie algebras of 1-dimensional 
unipotent subgroups of $G$ called {\em root subgroups}. 


Let $\mathfrak{b}$ denote the Lie algebra of $B$; more generally, $\mathfrak{b}_w$ will denote the Lie algebra 
of $B_w$. Then 
$$
\mathfrak{b}=\mathfrak{h}\oplus\left(\bigoplus_{\alpha\in \Phi^+}{\mathfrak g}_{\alpha}\right)$$ 
and, since $B_w$ is the conjugate to $B$ via $k_w$: 
\[
\mathfrak{b}_w=\Ad(k_w)(\mathfrak{b})=\mathfrak{h}\oplus\left(\bigoplus_{\alpha\in w(\Phi^+)}{\mathfrak g}_{\alpha}\right).
\]




In what follows we will be using the action of the Weyl group on the set of roots $\Phi$. 
Each $\alpha\in \Phi$ is an element of ${\mathfrak a}^*$ and, accordingly, $W$ 
(a group of automorphisms of ${\mathfrak a}$) 
acts on $\alpha$ by  the standard action on the dual space: $w(\alpha)= \alpha\circ w^{-1}$.

Define  
$$
{\mathfrak u}^\pm= 
\bigoplus_{\alpha\in \Phi^\pm}{\mathfrak g}_{\alpha},$$
the Lie algebra of the unipotent radical $U^\pm$ of the Borel subgroups $B^\pm< G$,  where $U^+=U< B^+=B$. We also define 
\[
\mathfrak{u}_w=\Ad(k_w)(\mathfrak{u}^-)=\bigoplus_{\alpha\in w(\Phi^-)}{\mathfrak g}_{\alpha}.
\]
Then the unipotent subgroup $U_w=\exp(\mathfrak{u}_w) < B_{w^\vee}$ is the maximal unipotent subgroup of the Borel subgroup 
$B_{w^\vee}$ opposite to $B_w$. Recall (\S \ref{sec:Thickenings in the flag manifold}) that $B_w\cap U_w=\{e\}$ 
and the group $U_{w}$ acts simply transitively on $\Opp(\sigma_{w^\vee})$. Hence, the orbit map  
\[
\pi_w: U_w\rightarrow \Opp(\sigma^\vee)\subset \cF, \quad u\mapsto u\cdot \sigma_w
\]
is biholomorphic. 
In particular, the tangent space of $\cF$ at $\sigma_w$ can be identified (via $\log\circ \pi_w$) with the Lie algebra 
$\mathfrak{u}_w$. 
Our next goal is to describe a certain product decomposition of $\Opp(\sigma_{w^\vee})$ and of the tangent space $T_{\sigma_w} \cF$. We will eventually prove that this decomposition will equal to a decomposition of $\Opp(\sigma^\vee)$ as a product of stable and unstable manifolds of $\gamma=\exp(a), a\in {\mathfrak a}^+$, at its fixed point $\sigma_w$ in $\cF$.  

\medskip 
For $w\in W$ we set
\begin{align}
\begin{aligned}
\label{eq_roots_w}
\Phi_w^+ = \{ \beta \in \Phi^+ \ | \exists \alpha \in \Phi^-, w(\alpha)=\beta\},\\
\Phi_w^- = \{ \beta \in \Phi^- \ | \exists \alpha \in \Phi^-, w(\alpha)=\beta \}
\end{aligned}
\end{align}

We also define 
\begin{equation} \label{eq_lie_alg_w}
\mathfrak{u}_w^+=\bigoplus_{\alpha\in\Phi_w^+}{\mathfrak g}_{\alpha}, \quad \mathfrak{u}_w^-=\bigoplus_{\alpha\in\Phi_w^-}{\mathfrak g}_{\alpha}
\end{equation}

According to the definition of $\Phi_w^\pm$, we get: 

\begin{obs}\label{obs:signs}
1. $\Phi_w^+$ consists of positive roots $\alpha$ which are negative on the Weyl chamber $w(\Delta)\subset {\mathfrak a}$, 
while  $\Phi_w^-$ consists of negative roots $\beta$ which are also negative on the Weyl chamber $w(\Delta)$. 

2. Therefore, for every $a\in {\mathfrak a}^+$, the linear operator $ad(a)$ acts with positive eigenvalues on 
$\mathfrak{u}_w^+$ and with negative eigenvalues on $\mathfrak{u}_w^-$. These eigenvalues are the values of the roots: 
$\al(a), \al\in \Phi^+_w$ and $\beta(a), \beta\in \Phi^-_w$, respectively. 
\end{obs}

Since 
$$
\Phi= \Phi_w^+ \sqcup \Phi_w^-, 
$$
we have a direct sum decomposition of $\mathfrak{u}_w$ as a vector space 
\[
\mathfrak{u}_w=\mathfrak{u}_w^+\oplus\mathfrak{u}_w^-.
\]
Recall that $U=U^+=U_{w_0}$ is the unipotent radical in the positive Borel subgroup $B=B^+$, while $U^-=U_e$ is the unipotent radical in the 
negative Borel subgroup $B^-=k_{w_0} B k_{w_0}^{-1}$. 
According to the above interpretation of $\Phi_w^\pm$ in terms of positivity/negativity on various chambers (Observation \ref{obs:signs}), we see that 
$\mathfrak{u}_w^+$ is the Lie algebra of the intersection $U^+\cap U_w$, while  
$\mathfrak{u}_w^-$ is the Lie algebra of the intersection $U^-\cap U_w$. The group 
$U^+\cap U_w$ is the unipotent radical 
in the intersection of Borel subgroups $B^+\cap B_{w^\vee}$, 
while the group $U^-\cap U_w$ is the unipotent radical 
in the intersection of Borel subgroups $B^-\cap B_{w^\vee}$. Thus, we have  
\begin{equation} \label{eq_subgroups_w}
U_w^\pm:=\exp(\mathfrak{u}_w^\pm)= U^\pm \cap U_w< B^\pm \cap B_{w^\vee}. 
\end{equation}

Note that the (real) dimension of $\mathfrak{u}_w^+$ is $2\vert w\vert$ while the real dimension of $\mathfrak{u}_w^-$ is 
$\dim(\cF)-2\vert w\vert= 2n- {2} \vert w\vert$ 
because the adjoint action of $k_w$ changes $\vert w\vert$ positive roots to negative ones and vice-versa 
(see \cite[the proof of Proposition 2.70]{knapp} or Remark \ref{rem:gallery} below). 
Thus $U_w^\pm$ are subvarieties of $U_w$ isomorphic (as varieties) respectively to 
$\bC^{\vert w\vert}$ and $\bC^{n-\vert w\vert}$ and we have an isomorphism of varieties 
$U_w\cong U_w^+\times U_w^-$ (since the exponential map of $U_w$ is  
biholomorphic). In general, this is not a direct (or even semidirect) product of Lie groups. Below we will see this as  
a {\em knit product}\footnote{Also known as the 
{\em Zappa-Sz\'ep product}.} decomposition $U_w= U_w^+ U_w^-$, meaning 
that every $u\in U_w$ splits uniquely as a product $u=u^+ u^-$, where $u^\pm\in U_w^\pm$. 

\medskip 
The next proposition is Theorem 6.15 in \cite{Ronan} (the group $U_w$ in Ronan's book is our group $U_w^+$): 

\begin{prop}\label{unipotentdescriptionofschubertcells1}
The group $U_w^+$ acts simply transitively on $S_w(\sigma^+)$. In particular, $\pi_w(U_w^+)=S_w(\sigma^+)$. 
\end{prop}

\medskip 
Switching the roles of $\sigma^+$ and $\sigma^-$ we also obtain: 

\begin{prop}\label{unipotentdescriptionofschubertcells2}
The group $U_w^-$ acts simply transitively on $S_{w^\vee}(\sigma^-)$. In particular, 
$\pi_w(U^-_w)= S_{w^\vee}(\sigma^-)$. 
\end{prop}

\begin{cor} \label{cor:knit}
We have the knit product decomposition 
$U_w= U_w^+ U_w^-$ or $U_w= U_w^- U_w^+$.
\end{cor}
\begin{proof} We will prove that $U_w= U_w^+ U_w^-$, since the second decomposition is obtained by switching the roles of $\sigma^+$ and $\sigma^-$.

1. Suppose $g=u_1^+ u_1^-=u_2^+u_2^-$, where $u_i^\pm\in U_w^\pm$, $i=1,2$. Then 
$$
u:=(u_2^+)^{-1}u_1^+ = u_2^- (u_1^-)^{-1}, \quad u \in U^+_w\cap U_w^-. 
$$
But $U^+_w\cap U_w^-=\{1\}$ since this subgroup is contained in $U^+\cap U^-$. 
Thus, we obtain uniqueness of decomposition of elements of $U_w$ as 
products of elements of $U_w^\pm$. 

2. Take $g\in U_w$ and apply it to $\sigma^-$. Since 
$$
w=\pos(\sigma_{w^\vee}, \sigma^-)= \pos(\sigma_{w^\vee}, g(\sigma^-)),
$$ 
by applying Proposition \ref{unipotentdescriptionofschubertcells2}, we find  
$u^{-1}\in U_w^+$ such that $u^{-1}(g(\sigma^-))=\sigma^-$. Hence, 
$h=u^{-1}\circ g\in U_w$ fixes $\sigma^-$ and, thus, belongs to $U_w^-$. Therefore, 
$$
g=u \circ h,  u\in U_w^+, h\in U_w^-. 
$$
\end{proof}

\begin{rem}\label{rem:gallery}
Ronan (in \cite[Theorem 6.15]{Ronan}) actually proves more: 
Take $w\in W$ of length $k$ and represent it as a shortest word as a product of simple reflections 
$$
w= s_{i_1}....s_{i_k}.
$$
For every $j=1,...,k$ define the subword
$$
w_j= s_{i_1}....s_{i_j}.
$$
We let $\beta_j$ denote the unique positive root which is positive on the interior of $w_{j-1}(\Delta)$ and negative on the interior of $w_j(\Delta)$. Thus, 
every root in $\Phi^+_w$ appears in this sequence exactly once. 
Then every $g\in U_w^+$ can be represented as a product of elements of the corresponding root subgroups 
$$
g=u_{1}... u_{k}, \quad u_{j}\in U_{\beta_{j}}=\exp({\mathfrak g}_{\beta_j}), 
$$
and this product representation is unique. Such product decompositions are known as {\em Malcev coordinates} on the unipotent group $U_w^+$. 
In other words, $U_w^+$ is the knit product
$$
U_w^+= U_{\beta_1}\ldots U_{\beta_k}.
$$
\end{rem}

Recall that the Zariski closures of $S_w(\sigma^+)$ and $S_{w^\vee}(\sigma^-)$ are respectively $\Th_w(\sigma^+)$ and $\Th_{w^\vee}(\sigma^-)$ whose homology classes are Poincar\'e dual of each other. 

\begin{prop}\label{prop:singleton}
We have $\{\sigma_w\}=\Th_w(\sigma^+)\cap\Th_{w^\vee}(\sigma^-)$.
\end{prop}
\begin{proof} There are several ways to prove this statement, we will give a proof using facts about the distance function $D$ on $\cF$ defined in 
\S \ref{section_relative_position_bruhat}. Note that $\sigma^\pm$ are opposite chambers, thus, $D(\sigma^+, \sigma^-)=|w_0|=n$. At the same time, 
$$
D(\sigma^+, x)\le |w|, \forall x\in \Th_w(\sigma^+), D(\sigma^-, x)\le |w^\vee|, \forall x\in \Th_w(\sigma^-),
$$
while $|w|+|w^\vee|=n$. Thus, for every $x\in \Th_w(\sigma^+)\cap\Th_{w^\vee}(\sigma^-)$ the triangle inequality  for the points $\sigma^+, x, \sigma^-$ is the equality. 
Hence, $x$ lies in the unique apartment $\amod$ containing $\sigma^\pm$. It remains to note that $\sigma_w$ is the unique chamber in $\amod$ which belongs to the intersection 
$\Th_w(\sigma^+)\cap\Th_{w^\vee}(\sigma^-)$. \end{proof}

\begin{prop}\label{nonopposite_thickening}
Let $k\leq n-1$. For any $w\in W$ such that $\vert w \vert =k$, we have $\Th^{n-k-1}(\sigma^-)\cap\Opp(\sigma^\vee)=\emptyset$.
\end{prop}
\begin{proof} Again there are several ways to argue, we use an argument similar to the proof of the previous proposition: 
$$
D(w^\vee, \sigma^-)=k, D(\sigma^-, x)\le n-k-1, \forall x\in \Th^{n-k-1}(\sigma^-),
$$
hence, $D(\sigma^\vee, x)\le k+ (n-k-1)= n-1$. However, for every $y\in \cF$ opposite to $\sigma^\vee$, $D(\sigma^\vee,y)=n$. Hence, $x\in \Th^{n-k-1}(\sigma^-)$ 
cannot be opposite to $\sigma^\vee$. 
\end{proof}

\subsection{Thickenings of general pairs of antipodal chambers}
\label{summary_torus_action}

We introduce the following notation: Let $\cY=\Opp(\cF \times \cF)\subset  \cF \times \cF$ 
denote the subset consisting of pairs of opposite points. We will use the notation  
$\bx=(x_+, x_-)$ or $\bx=(x^+, x^-)$ denote elements of $\cY$. \index{$\cY=\Opp(\cF \times \cF)$, set of pairs of opposite points in $\cF$}

For every $\bx=(x_+, x_-)\in \cY$ we define the following: We let ${\mathfrak a}_{\bx}$ denote the Lie algebra 
of the group of all transvections $A_{\bx}< G$  stabilizing the chambers $x_\pm$. We let ${\mathfrak a}^+_{\bx}$ 
denote the strictly positive Weyl chamber in 
${\mathfrak a}_{\bx}$ corresponding to the spherical chamber $x_+$. We set 
$$
A^+_{\bx}:= \exp({\mathfrak a}^+_{\bx})\subset 
A_{\bx}.$$

For general pairs of antipodal chambers $x^\pm\in \cF$ 
the descriptions of thickenings and product decomposition given in \S \ref{sec:Structure of thickenings} 
can be transported via $G$-action: Take $g\in G$ such that $g(x^\pm)=\sigma^\pm$ and apply $g^{-1}$ to obtain 
a decomposition of $\Opp(x_{w^\vee})$.  We briefly summarize this as follows. 


Set $x_w:= g^{-1}(\sigma_w)$, $w\in W$. We have 
$$
\{x_w\}=\Th_{w}(x^+)\cap \Th_{w^\vee}(x^-),$$ 
$$
\{x_{w^\vee}\}=\Th_{w^\vee}(x^+)\cap \Th_{w}(x^-),$$ 
opposite points in the apartment $a(x^+,x^-)\subset \tits \cX$ containing $x^+, x^-$. 
Similarly to the previous section, we define the maximal unipotent subgroup $U_{x_w}$ 
stabilizing $x_{w^\vee}$ and its subgroups $U_{x_w}^\pm$:
$$
U_{x_w}= g^{-1} U_{\sigma_w} g, \quad U_{x_w}^\pm= g^{-1} U^\pm_{\sigma_w} g. 
$$
As before, we have the knit product decomposition  $U_{x_w}=U_{x_w}^+ U^-_{x_w}$. The 
subgroups $U_{x_w}^+$, $U_{x_w}^-$ act simply transitively on Schubert cells $S_w(x^+)$, $S_{w^\vee}(x^-)$ respectively. Closures of these cells are the thickenings 
$\Th_{w}(x^+)$, $\Th_{w^\vee}(x^-)$ respectively. 

\medskip 
We also introduce a definition for later use. If $(x^+,x^-)$ is a pair of opposite flags, we say that this \index{$\varepsilon$-separated pair}
pair is \emph{$\varepsilon$-separated} if for any $w\in W$ we have 
\begin{equation}\label{e_separation_inequality}
d_{\cF}(x_w,\Th^{n-1}(x_{w^\vee}))\geqslant \varepsilon
\end{equation} 
with respect to the distance function $d_{\cF}$ on $\cF$ defined via the fixed K\"ahler metric on $\cF$. 
If $\gamma\in G$ is a regular loxodromic element, we say that $\gamma$ is $\varepsilon$-separated 
\index{separated loxodromic element} 
if its pair of attractive/repelling fixed chambers $(x^+,x^-)$ is.


\medskip 
We now relate the inequality \eqref{e_separation_inequality} to the more intuitive metric notion of separation for points in $\cF$ (we will be applying this lemma to limit sets $\La=\La(\Ga)$ of Anosov subgroups $\Ga< G$): 

\begin{lem}\label{lem:separate} 
Fix a closed antipodal subset $\La\subset \cF$. Then  
for every $\eps_0>0$ the following hold: 

\begin{enumerate}
\item[(i)] $d_{\cF}(x^+, \Th^{n-1}(x^-))\geqslant \eps_0\Rightarrow  d_{\cF}(x^+, x^-)\geqslant \eps_0$ for all $x^\pm \in \cF$. 
\item[(ii)] The set of $\eps_0$-separated pairs $(x^+,x^-)$ is a compact subset of $\cF^2$.
\item[(iii)]The subset $\{(x^+, x^-)\in \La^2: \dF(x^+, x^-)\geqslant \eps_0\}\subset \La^2\setminus \diag(\La^2)$ is compact.   
\item[(iv)] There exists $\epsilon>0$ such that for all $x^\pm \in \La$ with $d_{\cF}(x^+, x^-)\geqslant \eps_0$,  then $d_{\cF}(x^+, \Th^{n-1}(x^-))\geqslant \epsilon$.
\end{enumerate}
\end{lem}

\begin{proof} The first assertion follows directly from  the fact that  $x^-\in \Th^{n-1}(x^-)$ and the third from the continuity of the distance function together with the compactness of $\Lambda^2$.

For assertion (ii), we first observe that for every thickening $\Th$, the function $x\mapsto \Th(x)$, from $\cF$ to $Cl(\cF)$ (the set  
of closed subsets in $\cF$)  is continuous where $Cl(\cF)$ is equipped with the Hausdorff distance. 

Indeed fix a base point $x_0\in \cF$, by transitivity, any $x \in \cF$ is of the form $g \cdot x_0$ and the map $g \in G \mapsto g \cdot \Th(x_0) = \Th(g x_0) \in Cl(\cF)$ is continuous by Lemma \ref{lem_continuity_hausdorff}. The continuity follows from the fact that a neighborhood of the $g$ in $G$ maps onto a ball around $x \in \cF$. 
In particular, the function $(x,y) \in \cF^2 \mapsto d_{\cF}(x, \Th^{n-1}(y))$ is continuous. Moreover, since for each $w\in W$ the points $x_w, x_{w^\vee}$ depend continuously on the pair of opposite flags $(x_+,x_-)$, we deduce that for each $w \in W$, $(x_+, x_-) \mapsto d_{\cF}(x_w , \Th^{n-1}(x_{w^\vee}))$ is continuous. 
We then conclude that the set of $\epsilon_0$ separated flags is closed, hence compact in $\cF \times \cF$.
%
\smallskip 

Let us prove the last assertion. Consider 
$$
\theta: \La^2\setminus \diag(\La^2) \to \R_+
$$
by
$$
\theta(x^+, x^-):=\min_{w \in W} \{d_{\cF}(x_w, \Th(x_{w^\vee}))\}.
$$
This function is positive and continuous (since the points $x_w, x_{w^\vee}$ depend continuously on $(x^+, x^-)$). 
Hence, its minimum value $\eps = \theta(x_0^+, x_0^-) >0$ is reached on 
the compact $\{(x^+, x^-)\in \La^2: \dF(x^+, x^-)\geqslant \eps_0\}$ at some specific pair $(x_0^+ , x_0^-) \in \cF^2$, as required. \end{proof}

\section{Thickenings of subvarieties}


As in the previous section, we fix two opposite points $x_\pm\in \cF$ and the unique apartment $a=a(x^+, x^-)\subset \tits \cX$ containing $x^\pm$. 

For every $w\in W$ we have 
the chamber $x_w=w(x^+)$ in the apartment $a$. Thus,  
$x_w\in \Opp(x_{w^\vee})$. 

\begin{lem}\label{lem:opp-cover}
The Zariski open subsets $\Opp(x_w), w\in W$, form a finite cover of $\cF$. 
\end{lem}
\begin{proof}  The simplicial complex $\tits \cX$ has dimension $r-1$ and, therefore, every apartment in this complex 
defines a nontrivial homology class in $\tits \cX$. We equip 
$\tits X$ with the Tits metric; this metric is known to be $CAT(1)$, see \cite[Lemma 4.12]{BGS}. 
Suppose that  
$$
U=\bigcup_{w\in W} \Opp(x_w) \ne \cF.  
$$
Let  $\zeta\in \tits \cX$ be the barycenter of a chamber $\sigma$ in $\cF$ which does not belong to $U$. 
Then every point $\xi$ in the apartment $a$ is {\em not opposite} to $\zeta$, hence, is 
connected to $\zeta$ by unique geodesic in $\tits \cX$. Moreover,  this geodesic  depends continuously (with respect to the Tits topology) 
on the point $\xi$. (This is a consequence of the fact that $\tits \cX$ is a $CAT(1)$-space.) 
Therefore, using these geodesics, we see that the apartment $a$ is null-homotopic in $\tits X$. See \cite{Solomon} for 
details. This is a contradiction. 
\end{proof}

\begin{prop}\label{prop:projective}
For every thickening $\Th\subset W$ and a projective subvariety $S\subset \cF$, the thickening $\Th(S)$ is again a projective subvariety in $\cF$. 
\end{prop}
\begin{proof} It suffices to prove the claim in the case of thickenings of the form $\Th=\Th_w, w\in W$. 
For every $x\in \cF$, $\Th(x)$ is a projective subvariety in $\cF$ (a Schubert subvariety). As in the lemma, 
we cover $\cF$ by the open Schubert cells $\Opp(x_{v^\vee}), v\in W$, each of which is the $U_v$-orbit of $x_v\in  \Opp(x_{v^\vee})$. The intersection $\Sigma_v:=S\cap \Opp(x_{v^\vee})$ is an open 
affine subvariety in $S$. 
Since the biholomorphic map $\pi_v: U_v\to \Opp(x_{v^\vee}), u\mapsto u(x_v)$, is 
algebraic, the preimage 
$$
\tilde\Sigma_v:= \pi_v^{-1}(\Sigma_v)
$$
is a quasiprojective subvariety in $U_v$. For every $u\in \tilde\Sigma_v$, if $x=\pi_v(u)=u(x_v)$, then 
$$
\Th_w(x)= \Th_w(u(x_v))= u \Th_w(x_v). 
$$
Taking the union over all $u\in \tilde\Sigma_v$ we obtain
$$
\bigcup_{u\in  \tilde\Sigma_v} u \Th_w(x_v)= \Th_w(\Sigma_v),
$$
is the image of the algebraic map 
$$
\tilde\Sigma_v\times \Th(x_v)\to \cF, \quad (u, z)\mapsto u(z). 
$$
Therefore, the image $\Th(\Sigma_v)$ is quasiprojective. The union of these images over $v\in W$, is 
projective as well (since it is compact in the analytic topology on $\cF$).    
\end{proof}

\section{Torus action and hyperbolicity of type $(p,q)$}

\label{section_torus_action}

Given a diffeomorphism $f$ of a (complex or real) manifold and a fixed point $x$ of $f$, we say that $x$ is a \textit{hyperbolic point of type $(p,q)$} if the tangent space has a $df_x$-invariant splitting: 
$$T_x \cF  = E_s \oplus E_u,$$ 
where $\dim E_s = p $, $\dim E_u=q$ and such that the differential $df_x$ is (strictly) contracting in $E_s$ and (strictly) expanding in $E_u$. 
In this situation, the stable manifold theorem shows that there are $f$-invariant immersed submanifolds $\cW^s(x), \cW^u(x)$ through $x$ 
with $T_x \cW^s(x)= E_s$ and $T_x \cW^u(x)= E_u$, where the dynamics of $f$ is respectively contracting and expanding.
For example, a fixed point $x$ is a hyperbolic point of type $(0, d)$ iff it is a repelling fixed point, and $x$ is attracting iff 
it is of type $(d, 0)$.

\begin{rem}
When the manifold is complex and $f$ is holomorphic, all above data are holomorphic. In this book we mainly work i this situation and unless specified otherwise we use $p,q$ to denote complex dimensions when we use the notion of $(p,q)$-hyperbolicity.
\end{rem}   
 
\smallskip 
We shall verify that for a regular semisimple element $\gamma \in G$ acting by left multiplication $L_\gamma$ on 
$\cF=G/B$, its fixed points are hyperbolic points for which the types are determined by some Lie-theoretic data. In this context, the stable manifolds and unstable manifolds are appropriate Schubert cells. Such hyperbolic features play a crucial role in the dynamics that we are describing in this book. Actually,  in the situation of this section much more is true: We will see that the action of $\gamma$ is linear in some coordinates.

We will be considering regular loxodromic elements $\ga$ of the maximal torus $T$; 
for general regular loxodromic elements we get the same conclusion by conjugation (cf. \S \ref{summary_torus_action}). 


Consider a regular loxodromic element $\ga\in T$. We will assume that $\sigma^+$ and $\sigma^-$ are, respectively, the attractive and repelling fixed points of $\gamma$, i.e. $\ga\in T^+$, $\ga=\exp(a), a\in {\mathfrak h}^+$.  

Then the fixed points of $\gamma$ in $\cF$ are $\sigma_w, w\in W$. 
In particular, $\gamma$ preserves all Schubert cells and thickenings centered at these points.
Recall that we have defined subsets of roots $\Phi_w^{\pm}\subset \Phi$ (see \eqref{eq_roots_w}), 
nilpotent Lie algebras $\mathfrak{u}_w^{\pm}\subset {\mathfrak g}$ (see \eqref{eq_lie_alg_w}) and the 
corresponding unipotent subgroups $U_w^{\pm}< U_w$ (see \eqref{eq_subgroups_w}). 

For each $\alpha\in \Phi_w^\pm$, we have $U_\alpha$, the one-parameter subgroup equal to 
$\exp({\mathfrak g}_{\alpha})< U_w$. 
Recall from the previous section that we have 
$$
U_w=\prod_{\alpha \in \Phi_w^\pm}U_\alpha$$ 
(as varieties and in terms of the knit product of groups) and the biholomorphic orbit map 
$\pi_w: U_w\rightarrow \Opp(\sigma_{w^\vee})$, $g\mapsto g\sigma_w$.

\begin{prop}\label{prop_hyperbolic_torus_action_linear}
Take $\ga\in T$ as above. Then the action of $\gamma$ on $\Opp(\sigma_{w^\vee})$ (the restriction of the left multiplication $L_\ga$ on $G/B$), the conjugation map $Inn_\gamma: U_w\rightarrow U_w, g\mapsto \gamma g \gamma^{-1}$,  
and the linear map $e^{ad(a)}$ on $\mathfrak{u}_w$ are conjugate via $\pi_w$ and $\exp$. 
More precisely, the following diagram is commutative:
 \begin{equation}\label{torusactiondiagram}
 \xymatrix{
\mathfrak{u}_w \ar[d]^{\exp} \ar[r]^{e^{ad(a)}} & \mathfrak{u}_w \ar[d]^{\exp} \\
 U_w \ar[d]^{\pi_w}\ar[r]^{Inn_\gamma}  & U_w \ar[d]^{\pi_w}\\
 \Opp(\sigma_{w^\vee}) \ar[r]^{L_\gamma} & \Opp(\sigma_{w^\vee}).
 }
 \end{equation}
 Here $e^{ad(a)}$ is the exponential of the linear operator $ad(a)$ on the Lie algebra ${\mathfrak u}_w$.
\end{prop}
\begin{proof}
As $\gamma$ belongs to the maximal torus $T$, the conjugation $Inn_\gamma$ (also known as the {\em inner automorphism of $G$ induced by $\gamma$}) 
\index{$Inn_\gamma$, inner automorphism induced by $\gamma$} 
preserves $U_w$. Since $\gamma$ fixes $\sigma_w$, we have $\pi_w \circ Inn_\gamma  = L_\gamma \circ \pi_w $. Moreover, for all $X \in \mathfrak{u}_w$, by  \cite[Equation 1.88]{knapp}, one has $Inn_\gamma(\exp(X)) = \exp(Ad(\gamma)X)$. As $\gamma = \exp(a)$, by \cite[Equation 1.92]{knapp} 
$Ad(\exp(a)) = e^{ad(a)}$, and we obtain $Inn_\gamma(\exp(X)) = \exp ( e^{ad(a)}X )$. \end{proof}

As a consequence, we obtain the following:

\begin{prop} \label{prop_hyperbolic_torus_action} Let $a \in \mathfrak{h}^+$ and $\gamma = \exp(a)$. 
Then for any $w \in W$, the action of $\gamma$ on $\cF$ is hyperbolic of type 
$(2n - 2\vert w \vert , 2\vert w \vert )$ at the fixed point $\sigma_w$ and has the following properties: 
\begin{enumerate}
\item[(i)] Under the linear isomorphism given by $d\pi_w  : \mathfrak{u}_w \to T_{\sigma_w} \cF$, the decomposition $T_{\sigma_w} \cF =  \mathfrak{u}_w^- \oplus \mathfrak{u}_w^+ $ is preserved by $dL_\gamma$. Furthermore, under the chart on $\Opp(\sigma_{w^\vee})$ given by $\log\circ \pi_w^{-1}$, the map 
$\ga$ is linear, diagonalizable with absolute values of the diagonal entries given by $\exp(\al(Re(a))), \al\in \Phi_w$.  
\item[(ii)] $||(d L_\gamma)_{|\mathfrak{u}_w^- }|| < 1$ and $||(d L_\gamma)_{|\mathfrak{u}_w^+ }|| > 1$.
\item[(iii)] The stable and unstable manifolds $ \mathcal{W}^s(\sigma_w), \mathcal{W}^u(\sigma_w)$ 
of $\ga$ through the fixed point $\sigma_w$ are the Schubert cells: 
$$
S_{w^\vee}(\sigma^-)=\pi_w(U_w^-) = \mathcal{W}^s(\sigma_w),$$
$$ 
S_{w}(\sigma^+)=\pi_w(U_w^+) = \mathcal{W}^u(\sigma_w).$$
\end{enumerate} 
Here, as always, $n$ is the dimension of $\cF$. 
\end{prop}
 
 \begin{proof}
 By the definition, the adjoint action $ad(a)$ is diagonalizable and, by the Observation \ref{obs:signs}, 
 the real parts of the eigenvalues are negative on $\mathfrak{u}_w^-$ and positive on $\mathfrak{u}_w^+$ respectively. 
 As a consequence, taking the exponential shows that in the chart given by the inverse of  $\pi_w \circ \exp$, the action of $\gamma$ is linear, diagonal, preserves the decomposition $\mathfrak{u}_w^-\oplus \mathfrak{u}_w^+$ 
  and has absolute values of all eigenvalues $<1$ on $\mathfrak{u}_w^-$ and $>1$ on $\mathfrak{u}_w^+$. This proves assertions (i) and (ii). 
 Finally, the last assertion follows from the definition of $U_w^\pm = \exp(\mathfrak{u}_w^\pm)$ and 
 Propositions \ref{unipotentdescriptionofschubertcells1} and \ref{unipotentdescriptionofschubertcells2}.   
 \end{proof}

\begin{rem}
The same conclusion holds if instead of the chart on $\Opp(\sigma_{w^\vee})$ given by $\log\circ \pi_w^{-1}$ we use Malcev coordinates on $U_w$.  
\end{rem}

\section{Expansion factors of unipotent elements}  \label{sec:expansion_unipotent}
 
In this section we will prove a uniform estimate on {\em expansion factors} of unipotent element $g\in G$ acting on 
$\cF$. This estimate will be used in \S \ref{sec:Contraction estimates}.  
 
For an linear map $T: V_1\to V_2$ between two normed vector spaces, 
the {\em expansion factor} $\boldsymbol{\varepsilon}(T)$  is defined as:
\begin{equation}
\label{eq:expfac}
\boldsymbol{\varepsilon}(T) = \inf_{v\in V_1-\{0\}} \frac{\|Tv\|}{\|v\|}.  
\end{equation}
If $T$ is invertible,
$$
\boldsymbol{\varepsilon}(T) = \| T^{-1} \|^{-1}.
$$

We let $\bx=(x, \hat{x})\in \cF$ be a pair of opposite points; let $T=T_{\bx}< G$ be the maximal torus fixing both $x, \hat{x}$. 
Let $U=U_x$ denote the unipotent radical of the stabilizer of $\hat{x}$ in $G$.

We will identify $T_x\cF$ with the Lie algebra $\mathfrak{u}$ of $U$. 
We equip $\mathfrak{g}$ with an inner product where the root subalgebras ${\mathfrak g}_\al$ (with respect to the torus $T$) are pairwise orthogonal. 
We let $\mathcal{g}$ be a $K$-invariant Riemannian metric on $\cF$ whose restriction to $T_x\cF$ is the restriction 
 of the above inner product on 
$\mathfrak{g}$ to the subalgebra $\mathfrak{u}$. Let $||\cdot||$ denote the corresponding norm on tangent bundle $T\cF$.

At the same time, by compactness, there exists $\eps_0>0$ such that for all 
$k\in K$, $\boldsymbol{\varepsilon}(Ad (k))\geqslant \eps_0$.

\begin{lem}\label{lem:unipotent bounds}
For every unipotent element $g\in U$ we have
$$
\boldsymbol{\varepsilon}(dg_x: T_x\cF\to T_{gx}\cF)\geqslant \eps_0,
$$
where we compute the expansion factor with respect to the Riemannian metric $\mathcal{g}$. 
\end{lem}
\begin{proof} With respect to the above norm on $\mathfrak{g}$, 
$$
\boldsymbol{\varepsilon}(Ad(g)|_{\mathfrak{u}})\geqslant 1
$$ 
since all eigenvalues of $Ad(g)$ are equal to $1$. 

For $y:=g(x)$ we take $k\in K$ such that $k(y)=x$. Then for $h:= k\circ g$ we have 
$\boldsymbol{\varepsilon}(dg_x)= \boldsymbol{\varepsilon}(dh_x)$. 
Since $\boldsymbol{\varepsilon}(Ad(g)|_{\mathfrak{u}})\geqslant 1$ and 
$\boldsymbol{\varepsilon}(Ad(k))\geqslant\eps_0$, 
 we have
$$
\boldsymbol{\varepsilon}(Ad(h)|_{\mathfrak{u}})\geqslant\eps_0. 
$$
On the other hand,  $\boldsymbol{\varepsilon}(Ad(h)|_{\mathfrak{u}}) =  \boldsymbol{\varepsilon}(dh_x)$.  
Lemma follows. 
\end{proof}

%

\section{Examples}

\subsection{Products of $\bP^1$'s}\label{thickeningsofproductsofprojectivelines}

When $G=\PGL(2,\bC)^n$ is the product of $n$ copies of $\PGL(2,\bC)$, the Borel subgroup $B$ of $G$ can be taken to be 
$B_0^n$, where $B_0$ is the Borel subgroup of $\PGL(2,\bC)$ (say, the projection of the subgroup of 
upper triangular matrices in $\GL(2,\bC)$). The flag manifold $\cF=G/B$ is isomorphic to the product of $n$ copies of $\bP^1$. 
The Weyl group $W$ is isomorphic to the product of $n$ copies of $\bZ/2\bZ$. The set of generators $\Sigma$ of $W$ corresponding to simple positive roots is given by the non-trivial elements of each copy of $\bZ/2\bZ$.

Any element of $W$ can be written as a $k$-multi-index $J=(j_1,\cdots,j_k)$ where $\{j_1,\cdots,j_k\}$ is a subset of $\{1,\cdots,n\}$. If $J$ is a $k$-multi-index, then we denote by $\hat{J}$ the complementary $(n-k)$-multi-index. For any $k$-multi-index $J$ with $0< k <n$, we denote by $\pi_J$ the projection of $\cF$ onto $\cF_J=\bP^1_{j_1}\times \cdots\times\bP^1_{j_k}$. For any point $x\in \cF$, we have
\[\Th_J(x)=\pi_{\hat{J}}^{-1}\left(\pi_{\hat{J}}(x)\right).\]

We can choose an affine coordinate on $\bP^1$ so that an element of the torus $T$ acts on $\cF=(\bP^1)^n$ as
\[(x_1,\cdots,x_n)\mapsto (\lambda x_1,\cdots,\lambda x_n), \quad \lambda\in \bC,\]
while an element of $T^+$ has the form
\[(x_1,\cdots,x_n)\mapsto (\lambda x_1,\cdots,\lambda x_n), \quad \lvert\lambda\rvert >1.\]
Fixed points of $T^+$ have the form $(x_1,\cdots,x_n)$ with $x_i$ equals $0$ or $\infty$. Such a fixed point $(x_1,\cdots,x_n)$ is $(p,q)$-hyperbolic for $T^+$ if $p$ is the number of $\infty$ among the $x_i$ and $q$ is the number of $0$.

\subsection{The classical full flag manifold}\label{homologyflagmanifolds}

When $G=\PGL(n+1,\bC)$, we take the Borel subgroup $B$ corresponding to upper triangular matrices. 
The unipotent radical $N$ of $B$ is isomorphic to the group of upper triangular matrices with $1$'s on the diagonal. 
A full (projective) flag (an element of $G/B$) is a sequence of projective subspaces in $\bP^n$
\[ E_1\subset \cdots\subset E_{n}\]
with $\dim E_i=i-1, i=1,...,n$. The flag manifold $\cF=G/B$ is a smooth complex projective variety parametrizing all full flags in $\bP^n$. The Weyl group $W$ is isomorphic to the symmetric group $\mathfrak{S}_{n+1}$. The set of generators $\Sigma$ of $W$ corresponding to simple positive roots is given by the transpositions $(12),(23),\dots,(n(n+1))$ in $\mathfrak{S}_{n+1}$. 

Let $x=(E_1, E_2, \cdots, E_{n})$ be a full flag in $\bP^n$ and let $w\in \mathfrak{S}_{n+1}$ be a permutation. 
The thickening $\Th_w(x)$ is the closure in $\cF$ of the following set of flags:
\begin{equation}
\left \{(F_1,\cdots,F_{n})\in \cF\vert \operatorname{dim}(F_i\cap E_j)=\#\{k\leq i, w(k)\leq j\}-1\right \}
\label{eq:schubertcelldescribedbypermutation}
\end{equation}
where $\operatorname{dim}(F_i\cap E_j)=-1$ means that $F_i, E_j$ do not intersect each other.

\begin{example}\label{dimthreeflagexample1}
The space of lines in $\bC^3$ is $\bP^2=G/P=G/_{P_{\taumod}}$, where $\taumod$ is the vertex of $\simod$ of the {\em point-type}; 
the space of hyperplanes in $\fk^3$ is $\bP^{2\vee}=G/\Pc= G/_{P_{\iota\taumod}}$. 
The full flag manifold $\cF$ parametrizes pairs $(p,l)$ with $p\in \bP^2, l\in \bP^{2\vee}$ and $p\in l$. Denote $\bP^2$ and $\bP^{2\vee}$ by $\bP^2_1,\bP^2_2$. The two projections $\pi_1: (p,\lambda)\mapsto p$ and $\pi_2: (p,l)\mapsto l$ realize $\cF$ in two ways as $\bP^1$-bundle and induce an embedding of the three dimensional $\bC$-variety $\cF$ into $\bP^2_1\times \bP^{2}_2$. The image by $\pi_1$ of a fiber of $\pi_2$ is a line in $\bP^2_1$ and vice versa. For $i=1,2$, let $F_i$ be a fiber of $\pi_i$ and let $\omega_i$ be a $(1,1)$-form on $\cF$ pulled back from the Fubini-Study form on $\bP^2_i$. Then the Poincaré dual classes of $F_1,F_2$ generate $\coh^{4}(\cF,\ZZ)$ while $[\omega_1],[\omega_2]$ generate $\coh^{2}(\cF,\RR)=\coh^{1,1}(\cF,\RR)$. We can normalize $\omega_1,\omega_2$ so that $[\omega_i^2]$ is the current associated with $F_i$, i.e.\ for a $(1,1)$ form $\alpha$ on $\cF$ we have 
\[\int_{F_i}\alpha=\left\langle \alpha,\omega_i^2\right\rangle=\int_{\cF}\omega_i^2\wedge\alpha.\]

The Weyl group of $G=\PGL(3,\bC)$ is $\mathfrak{S}_3$. 


Let $x\in \cF$ be the flag $(p,l)$. It follows from the description \eqref{eq:schubertcelldescribedbypermutation} that 
\[\Th_{(132)}(x)=\{(q,m)\in \cF \vert q=p\}, \Th_{(213)}(x)=\{(q,m)\in \cF \vert m=l\}.\]
In other words these two thickenings are respectively the fibers of $\pi_1,\pi_2$ containing $x$.  

The thickening $\Th_{(231)}(x)$ is the $2$-dimensional surface $\pi_1^{-1}(l)$ where $l$ is viewed as a line in $\bP^2$. The projection $\pi_2$ realizes $\Th_{(231)}(x)$ as a one point blow-up of $\bP^{2\vee}$ at $l$. Similarly, $p$ determines a line $p^{\vee}$ in $\bP^{2\vee}$, the pencil of lines in $\bP^2$ passing through $p$, and $\Th_{(312)}(x)=\pi_2^{-1}(p^{\vee})$. The intersection $\Th_{(231)}(x)\cap \Th_{(312)}(x)$ is the union of $\Th_{(132)}(x)$ and $\Th_{(213)}(x)$.

The torus $T$ is the projection of the set of diagonal $3\times 3$ matrices. An element $\ga$ of $T^+$ can be represented by a matrix of the form
\[
\begin{pmatrix}
\lambda_1&0&0\\ 0&\lambda_2&0\\ 0&0&1
\end{pmatrix},\quad \lvert \lambda_1\rvert>\lvert \lambda_2\rvert>1.
\]

\begin{figure}[htbp] 
   \centering
   \includegraphics[width=3in]{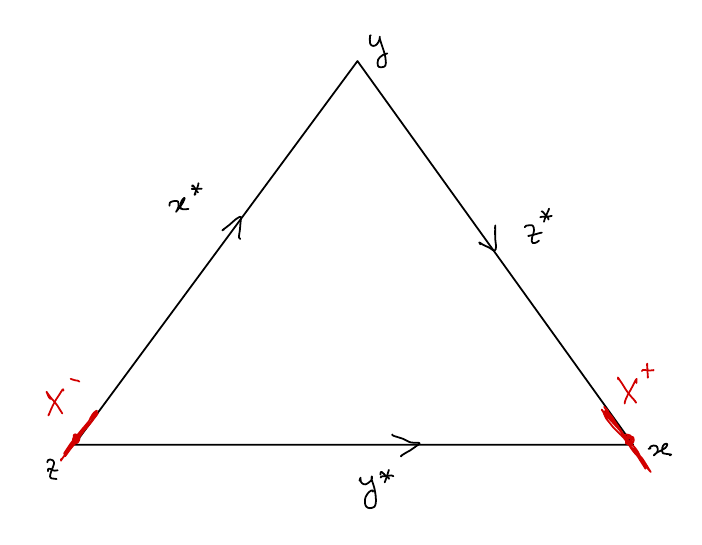} 
   \caption{Action of $\ga$ on $\P^2$.}
   \label{fig:triangle}
\end{figure}

The three coordinate axis of $\bC^3$ give rise to three fixed points $x,y,z\in \bP^2$ for the action of $\ga$ on $\bP^2$; here $x$ is attracting, $y$ is hyperbolic and $z$ is repelling. The three lines $x^{*}=\overline{yz},y^{*}=\overline{xz},z^{*}=\overline{xy}$ in $\bP^2$ are $\ga$-invariant and give rise to three fixed points $x^{*},y^{*},z^{*}\in \bP^{2\vee}$; here $z^{*}$ is attracting, $y^{*}$ is hyperbolic and $x^{*}$ is repelling. The actions of $\ga$ on $x^{*},y^{*},z^{*}$ (considered as lines in $\bP^2$) are respectively conjugate to multiplication by $\lambda_2,\lambda_1,\frac{\lambda_1}{\lambda_2}$. Therefore $\ga$ has six fixed points in $\cF$:
\[(x,y^{*}),(x,z^{*}),(y,x^{*}),(y,z^{*}),(z,x^{*}),(z,y^{*}).\]
\end{example}
They have respectively type
\[
(2,1),(3,0),(1,2),(2,1),(0,3),(1,2).
\]
Consider the repelling fixed flag $X^-=(z,x^{*})$. Its full two-dimensional thickening $\Th^2(X^-)$ contains all $\ga$-fixed points except the attracting fixed point $X^+=(x,z^{*})$. The complement $\cF\setminus \Th^2(X^-)$ is the set of flags opposite to $X^-$, i.e.\ the set of flags $(p,l)$ such that $p\notin x^{*}$ and $z\notin l$. As a line in $\bP^2$ passing through $p$ but not through $z$ is determined by its intersection point with $x^{*}$, the complement $\cF\setminus \Th^2(X^-)$ is isomorphic, as an affine variety, to the product of $\bP^2\setminus x^{*}$ with $x^{*}\setminus \{z\}$, thus to $\bC^3$. It is $\ga$-invariant and the $\ga$-action on $\cF\setminus \Th^2(X^-)$ is conjugate, via this isomorphism, to that of the contracting diagonal matrix 
\[
\begin{pmatrix}
1/\lambda_1 & 0& 0\\ 0& \lambda_2/\lambda_1 &0\\ 0&0&1/\lambda_2
\end{pmatrix}.
\]

We have similar descriptions if we begin with other fixed points. For example $\cF\setminus\Th^2\left((y,x^{*})\right)$ is isomorphic to $\bC^3$, with $(x,y^{*})$ sent to $0$, and the $\ga$ action on this complement is conjugate to 
\[
\begin{pmatrix}
1/\lambda_1 & 0& 0\\ 0&\lambda_2/\lambda_1 &0\\ 0&0& \lambda_2
\end{pmatrix}.
\]
The stable manifold at $(x,y^{*})$ is the Schubert cell whose closure is the two dimensional thickening $\Th_{(312)}\left((x,y^{*})\right)$, i.e.\ the preimage under $\pi_2$ of the pencil of lines through $x$; the unstable manifold at $(x,y^{*})$ has closure $\Th_{(132)}\left((x,y^{*})\right)$, i.e.\ $\pi_1^{-1}(x)$. 





\chapter{Anosov subgroups}\label{sec:Anosov}

\section{Discrete subgroups of higher rank Lie groups}\label{sec:def anosov}

Ano\-sov representations and Anosov subgroups of Lie groups were defined in \cite{Labourie} in the case of fundamental groups of compact negatively curved manifolds and in \cite{Guichard-Wienhard} for general hyperbolic groups; alternative characterizations were subsequently given  in \cite{KLP14}, \cite{KLP17}, \cite{KLP18}, \cite{Kapovich-Leeb-finsler}. We refer the reader to these papers as well as \cite{manicures} for the detailed treatment of the theory. Below is just an overview. 

\subsection{Regularity}

The class of Anosov subgroups of Lie groups is meant to generalize convex-cocompact subgroups of rank 1 Lie groups and, more generally, 
quasiconvex-cocompact subgroups of isometry groups of proper Gro\-mov hyperbolic spaces. We start with a geometric motivation. Discreteness of a subgroup $\Ga$ of a semisimple Lie group $G$ means (geometrically speaking) that for every base-point $x_0\in \cX=G/K$, points in the $\Ga$-orbit $\Ga x_0$ diverge to infinity, i.e.
$$
\lim_{i\to\infty} d(x_0, \ga_i x_0)=\infty,
$$
for every sequence $(\ga_i)$ of distinct elements of $\Ga$.

However, $\cX$ also has a more refined notion of distance, given by the vector-valued distance function $d_\Delta$. A sequence of vectors $(v_i)$ in $\Delta$ can diverge to infinity in different ways with respect to the faces of the chamber $\Delta$. Faces of $\Delta$ are labelled by faces $\taumod$ 
of the spherical chamber $\simod$ (just like $\Delta$ is a cone over $\simod$, the face $V_{\taumod}\subset \Delta$ is a cone over $\taumod$).  
\index{regular sequence} 
A sequence $(v_i)$ in $\Delta$  is said to be {\em $\taumod$-regular} if 
its distance to every face disjoint from the interior of $V_{\taumod}$ diverges to infinity. For instance, if $V_{\taumod}=V_{\simod}=\Delta$ then $(v_i)$ is $\simod$-regular if its distance to the entire boundary  
of $\Delta$ diverges to infinity. One can state this divergence condition in terms of simple roots of the root system of $G$ corresponding to the chamber $\Delta$. 
Let $\alpha_1,...,\alpha_m$ denote those simple roots which do not vanish identically on $V_{\taumod}$. Then $\taumod$-regularity of  $(v_i)$ means that 
$$
\lim_{i\to\infty} \alpha_j(v_i)=\infty, j=1,...,m. 
$$
Accordingly, $\simod$-regularity means that for every simple root $\al$, 
\begin{equation}\label{eq:regular}
\lim_{i\to\infty} \alpha(v_i)=\infty. 
\end{equation}
This divergence to infinity of values of roots is at most linear in terms of $||v_i||$. A sequence $(v_i)$ is said to be {\em uniformly} $\taumod$-regular if \index{uniformly regular sequence}
$$
\liminf_{i\to\infty} \frac{\alpha_j(v_i)}{||v_i||}>0, j=1,...,m. 
$$
In other words, we require divergence in \eqref{eq:regular} to be linear. 
These notions of regularity extend to sequences in $\cX$ and sequences in $G$. A sequence $(x_i)$ in $\cX$ is said to be {\em $\taumod$-regular} (resp. {\em uniformly 
$\taumod$-regular}) if the sequence of vectors $d_\Delta(x_0, x_i)$ is. A sequence $g_i\in G$ is said to be 
$\taumod$-regular (resp. uniformly $\taumod$-regular) if for some (equivalently, every) $x\in \cX$, the sequence $x_i=g_i(x)$ is. For instance, if $g\in G$ is 
a {\em regular loxodromic element} (see Section \ref{sec:Cartan and Lyapunov}), then the sequence of iterates $(g^i)$ is $\simod$-regular.   
A representation \index{regular representation}
$\rho: \Ga\to G$ is said to be $\taumod$-regular (resp. uniformly $\taumod$-regular) if for every sequence of distinct elements $(\ga_i)$ in $\Ga$, the sequence 
of images $\rho(\ga_i)$ is $\taumod$-regular (resp. uniformly $\taumod$-regular). Accordingly, a subgroup $\Ga< G$ is \index{regular subgroup}
$\taumod$-regular (resp. uniformly $\taumod$-regular) if the identity embedding $\rho: \Ga\to G$ is.  

In particular, every $\taumod$-regular representation has finite kernel and discrete image. If $G$ has rank 1 (equivalently, $\cX$ has rank 1) then 
$d_\Delta =d$ and we just recover the notion of discrete subgroups and representations with finite kernel and discrete image. But in higher rank we get something more interesting. For instance, suppose that $G=G_1\times G_2$, $G_1=G_2=\SL(2,\C)$. Then $\Delta=\R_+\times \R_+$. 
A sequence $g_i=(g_i', g_i'')\in G$ is $\simod$-regular if and only if each sequence $(g_i'), (g_i'')$ in $G_1=G_2$ diverges to infinity. To get an example of a $\simod$-regular but not uniformly regular sequence, take $g=(g', g'')$, where $g'$ is loxodromic and $g''$ is parabolic. Then the sequence of iterates $(g^i)_{i\in \bN}$ is $\simod$-regular but not uniformly regular. As another example, consider $G=\SL(n,\C)$. Then $\simod$-regularity of a sequence $(g_i)$ in $G$ means the following. Let $\sigma_1(g)\ge ...\ge \sigma_n(g)$ denote the singular values of a matrix $g\in G$. Then $\simod$-regularity of a sequence $(g_i)$ in $G$ is equivalent to 
$$
\lim_{i\to\infty} \frac{\sigma_k(g_i)}{\sigma_{k+1}(g_i)}=\infty,  k=1,...,n-1, 
$$
cf.  \S \ref{sec:Cartan and Lyapunov}. 

\subsection{Limit sets}
\label{sec:flag_limit_set}

Let $P_{\taumod}$ denote the parabolic subgroup of $G$ corresponding to the simplex $\taumod$, the stabilizer of $\taumod$ in $G$. 
Regularity of a subgroup $\Ga< G$ translates to the existence of certain compactifications of orbits $\Ga x\subset \cX$ by adding a certain subset of $G/P_{\taumod}$. We will explain this 
in the case when $\taumod=\simod$ since the definition is simpler in this case. Let $(x_i)$ be a $\simod$-regular sequence in $\cX$. \index{flag-convergence}
We say that $(x_i)$ {\em flag-converges} to $\sigma\in G/P_{\simod}$ if the following holds: Due to regularity, for every sufficiently large $i$ there exists a unique chamber $\sigma_i\in  
G/P_{\simod}$ such that $x_i\in V(\sigma_i)$, the cone from $x_0$ over $\sigma_i$. Then we say that
$$
\lim_{i\to\infty} x_i=\sigma
$$
if $\lim_{i\to\infty} \sigma_i=\sigma$ where this convergence is understood in the manifold topology of the flag-manifold $G/P_{\simod}$. In the rank one case we 
simply recover the usual convergence of sequences in $\cX$ to points of the visual boundary of $\cX$. Let's illustrate this notion in the same product example as above. Take a regular sequence $x_i=(y_i,z_i)\in \bH^3\times \bH^3$. Chambers $\sigma$ in $\geo \cX$ corresponds to pairs $(\eta, \zeta)\in S^2\times S^2$. Then 
$$
\lim_{i\to\infty} x_i=\sigma \iff \lim_{i\to\infty} y_i=\eta \hbox{~and~}  \lim_{i\to\infty} z_i=\zeta. 
$$
We now return to the general case. Let $\Ga< G$ be a $\simod$-regular subgroup. One defines the {\em flag-limit set} 
$\La_{\simod}(\Ga)\subset \cF= G/P_{\simod}$ as the accumulation set of one (equivalently, every) $\Ga$-orbit in $\Ga x\subset \cX$ in the flag-manifold 
$\cF$. Then we obtain a compactification of $\Ga x$:
$$
\Ga x\cup \La_{\simod}(\Ga). 
$$
The same works for proper faces $\taumod\subset \simod$ and one obtains a flag-limit set $\La_{\taumod}(\Ga)\subset G/P_{\taumod}$. 
Flag-limit sets share some basic properties with limit sets in the rank 1 setting (and, more generally, limit sets of discrete group actions on 
Gromov-hyperbolic spaces): For instance, $\La_{\taumod}(\Ga)$ is closed and $\Ga$-invariant. However, it is no longer true that for (nonelementary $\Ga$) the flag-limit set contains every closed $\Ga$-invariant proper subset of the flag-manifold. 


In the case $\taumod=\simod$ (which is the setting for most of this book) 
one can interpret {\em uniform regularity} and the limit set as follows. Consider the accumulation set 
$L(\Ga)$ of $\Ga x\subset \cX$ in the visual boundary $\geo \cX$ of $\cX$. Then a discrete subgroup $\Ga< G$ is uniformly $\simod$-regular if and only if 
$L(\Ga)$ is disjoint from the codimension 1 skeleton in the simplicial complex $\tits \cX$. The flag-limit set 
$\La_{\simod}(\Ga)$ consists of those chambers $\sigma\in \cF$ which have nonempty intersection with $L(\Ga)$.

The class of regular subgroups of $G$ is quite interesting but is definitely smaller than the class 
of all discrete subgroups (unless $G$ has rank 1). In particular, if $\Ga< G$ is a lattice in a semisimple 
Lie group of rank $\ge 2$, then $\Ga$ is not $\taumod$-regular for any $\taumod$. One advantage of $\taumod$-regular subgroups (and, more generally, sequences in $G$) is that they satisfy a version of the {\em convergence property} familiar in the case of $SL(2,\C)$ and, more generally, rank 1 Lie groups:

Every $\taumod$-regular sequence $g_i$ in $G$ contains a subsequence $(g_{i_j})$ such that there exist two simplices, $\tau_+\in G/P_{\taumod}$ and 
$\tau_-\in G/P_{\iota\taumod}$ such that $(g_{i_j})$ converges to $\tau_+$ uniformly on compacts in the open (and dense) subset 
$\Opp(\tau_-)\subset G/P_{\taumod}$ (cf. \eqref{eq:polar action}).  
Conversely, if $(g_{i})$ is a sequence in $G$ 
whose every subsequence contains a further subsequence  converging to some 
$\tau_+\in G/P_{\taumod}$ uniformly on compacts in $\Opp(\tau_-)$, then $(g_i)$ is $\taumod$-regular. 

Replacing $g_i$ with $g_i^{-1}$ switches the roles 
of $\tau_+$ and $\tau_-$. Since these simplices belong to (a priori) different flag-manifolds, it is convenient to assume (which we will do from now on) 
that $\iota\taumod=\taumod$. (This is automatically satisfied if $\taumod=\simod$, the case we are 
mainly interested in.) It turns out that if $\Ga< G$ is a 
$\taumod$-regular subgroup then its flag-limit set $\La_{\taumod}(\Ga)$ is exactly the set of limiting simplices \index{flag-limit set}
$\tau_+$ for sequences in $\Ga$. Of course, the usual convergence property would mean convergence on the complement to  
a point $\tau_-\in  G/P_{\taumod}$. That would be too much to ask. One can ask, however, if the restriction of the $\Ga$-action to 
$\La_{\taumod}(\Ga)$ is a {\em convergence action}, as defined in \S \ref{sec:convdy}. 
This is not true in general but holds under the following assumption:

$\La_{\taumod}(\Ga)$ is {\em antipodal}, i.e. any two distinct limit points are opposite to each other. \index{antipodal limit set}

The class of {\em $\taumod$-regular antipodal subgroups} appears to be the true higher rank generalization of the notion of Kleinian groups.

By analogy with rank 1 dynamics one can ask if the $\Ga$-action is properly discontinuous on the complement to 
the flag-limit set in $G/P_{\taumod}$. This is again too much to ask (unless $G$ has rank 1) but turns out to be true if 
$\La_{\taumod}(\Ga)$ is antipodal and one removes from $G/P_{\taumod}$ not just the limit set but also its {\em fat thickening}, see 
\cite{KLP17}.  Even for Anosov subgroups, the $\Ga$-dynamics on the complement to $\La_{\taumod}(\Ga)$ is quite complicated 
and is partially chaotic because of presence of {\em hyperbolic points}. Some partial analysis of this action at points where proper 
discontinuity (even non-wandering) fails is the main topic of this book.

\subsection{Definition and properties of Anosov subgroups}

We are now ready for the definition of $\taumod$-Anosov subgroups and representations. 
The next definition was not the original one (the reader can find it in \S \ref{sec:Hyperbolicity of flows}, Definition \ref{defn:original Anosov}) but is proven to be equivalent to it in \cite{KLP18}:

\begin{definition}\index{Anosov representation} \index{Anosov subgroup}
Let $\Ga$ be a finitely generated group and $\rho: \Ga\to G$ a representation. Then $\rho$ is $\taumod$-Anosov 
(or $P_{\taumod}$-Anosov) if it satisfies the  $\taumod$-URU property:

1. $\rho$ is undistorted, i.e. for every $x\in \cX$ the orbit map $o_x: \ga\mapsto \rho(\ga)(x)$   
is a quasi-isometric embedding $o_x: \Ga\to \cX$. 

2. $\rho$ is uniformly $\taumod$-regular. 

A (finitely generated) subgroup $\Ga< G$ is  $\taumod$-Anosov if the identity embedding $\rho: \Ga\to G$ is. 
\end{definition}

It is useful to rewrite this definition in terms of simple roots $\al_1,...,\al_m$ as above. In view of equivariance of the orbit map 
$o_x$, for every hyperbolic length function $\beta$ on $\Ga$ there are constants $C_0, A_0$ such that
\begin{equation}\label{eq:A0}
d(x, \rho(\gamma)x)\leqslant C_0\beta(\gamma) +A_0
\end{equation}
(i.e. the orbit map $o_x$ is coarse Lipschitz). In particular,  since the simple roots 
$\alpha_i\in {\mathfrak a}^*, i=1,...,m,$ are linear, there are constants $C_1, A_1$ such that for each $i=1,....,m$, we have 
\begin{equation}\label{eq:A3}
\mu_i (\rho(\gamma)) \leqslant C_0\beta(\gamma) +A_0. 
\end{equation}
Since $\rho$ is {\em undistorted} (i.e. $o_x$  is a quasiisometric embedding) we also have the reverse inequalities 
\begin{equation}\label{eq:A2}
k_1' \beta(\gamma) - a_1'  \leqslant ||\mu(\rho(\gamma))||= d(x, \rho(\ga)(x)), 
\end{equation}
for some uniform constants $k_1'\geqslant 1$ and $a_1'$.  (Here we are using regularity of the length function $\beta$ on $\Ga$.) 
Uniform $\taumod$-regularity of $\rho$ combined with \eqref{eq:A0} and \eqref{eq:A3} 
implies that for every $i=1,...,m$, 
\begin{equation}\label{eq:A1}
k^{-1} \beta(\gamma) -a \leqslant \mu_i (\rho(\gamma))= \al_i(\mu(\rho(\ga))) \leqslant  k^{-1} \beta(\gamma) -a
\end{equation}
for some uniform constants $k$ and $a$. It is easy to see that the inequality \eqref{eq:A1} implies that $\rho$ satisfies the $\taumod$-URU condition. 

\begin{rem}
This is a geometric version of the {\em domination property} in dynamics. 
\end{rem}

\medskip 
Here are some properties of Anosov subgroups $\Ga< G$:

\begin{enumerate}
\item Each Anosov subgroup is Gromov-hyperbolic. 

\item There exists an equivariant homeomorphism 
$$f: \geo \Ga \to \La_{\taumod}(\Ga)$$
from the Gromov boundary of $\Ga$ to the flag-limit set 
$\La_{\taumod}(\Ga)\subset G/P_{\taumod}$. This homeomorphism is unique if $\Ga$ is nonelementary, i.e. $\geo \Ga$ contains at least three points.   The map $f$ is called the {\em boundary map} of the Anosov subgroup $\Gamma$. 

\item One can say a bit more in the context of (2). Let $x\in \cX$ be a base-point with the orbit map 
$o_x: \ga \mapsto \ga(x), \ga\in \Ga$. Then the boundary map $f$ is a continuous extension of $o_x$ with respect to the topology of flag-convergence. 

\item The boundary map $f: \geo \Ga\to \La_{\taumod}(\Ga)$ is H\"older with respect to any visual metric on $\geo \Ga$ and any Riemannian metric on $G/P_{\taumod}$, see \cite[Theorem 1.6]{Tsouvalas}. In Chapter \ref{sec:holder} we prove a stronger form of the {\em bi-H\"older} property for the map $f$ in the case of Borel-Anosov representations.

\item The limit set $\La_{\taumod}(\Ga)$ is antipodal (i.e. any two distinct points are opposite to each other), see \cite{KLP17}. 

\item Every limit point $\la\in \La_{\taumod}(\Ga)$ is {\em conical}, see \cite{KLP17}. We will define what it means only in the case $\taumod=\simod$ and refer to \cite{KLP17} for the case 
of general $\taumod$. Then $\la$ is an element of the full flag-manifold $\cF=G/P_{\simod}$ and can be regarded as a chamber in the spherical building $\tits \cX$. Let $V$ denote a Euclidean Weyl chamber in $\cX$ whose ideal boundary is $\la$. (This is an analogue of a geodesic ray asymptotic to $\la$ if $\cX$ has rank one, e.g. is the hyperbolic plane.) Then for every orbit $\Ga x\subset \cX$ there exists a constant $C$ and a sequence $x_n\in \Ga x$ converging to $\la$ and contained in the $C$-neighborhood of $V$ in $\cX$.

\item Except for finitely many conjugacy classes (consisting of finite order elements of $\Ga$), every element $\ga$ of $\Ga$ is regular loxodromic. 
The attractive/repelling fixed points of loxodromic elements $\ga\in \Ga$ belong to $\La_{\taumod}(\Ga)$. 
\end{enumerate}

\begin{rem}
Since in this book we will be primarily interested in the case $\taumod=\simod$ (i.e. $P_{\simod}$ is a minimal parabolic subgroup of $G$), 
we will use the notation $\La(\Ga)$ for the flag-limit set $\La_{\simod}(\Ga)$. Accordingly, we will refer to $P_{\simod}$-Anosov subgroups 
simply as Anosov subgroups. Note also that in the case of complex semisimple Lie groups, $P_{\simod}$ is a Borel subgroup of $G$. 
Accordingly, $\simod$-Anosov subgroups are often called {\em Borel-Anosov subgroups} and the corresponding representations are called {\em Borel-Anosov}. 
\end{rem}

\section{Examples of Anosov subgroups}

Anosov subgroups of higher rank Lie groups generalize the notion of convex-cocompact subgroups $\Ga< G_1$ of rank one Lie groups. Since the associated symmetric space $X_1=G_1/K_1$ of a rank one Lie group is a rank one symmetric space of noncompact type, it is a negatively curved symmetric space. In particular, $X_1$ is Gromov-hyperbolic. Thus, one can define convex-cocompact subgroups $\Ga< G_1$ as quasiconvex-cocompact groups of isometries of $X_1$. This, of course, was not the original definition, which is that there exists a closed convex nonempty $\Ga$-invariant subset $C\subset X_1$ such that $C/\Ga$ is compact. However, the two notions are equivalent in the setting of groups of isometries of negatively curved symmetric spaces, see \cite{Bowditch1995}. 


Suppose that $G$ is a semisimple Lie group and $\rho: G_1\to G$ is a representation with finite kernel (where $G_1$ is, as before a rank one Lie group). Then 
$\rho(\Ga)< G$ is a ${\taumod}$-Anosov subgroup for a face $\taumod\subset \simod$ of the model spherical chamber of the symmetric space of $G$ (typically, $\taumod\ne \simod$). Moreover, if $\La_1(\Ga)$ denotes the  limit set of $\Ga$ in $\geo X_1$, then there exists a $\Ga$-equivariant homeomorphism $\La_1(\Ga)\to \La(\Ga)\subset G/P_{\taumod}$ to the flag-limit set 
of $\rho(\Ga)$ in the flag-manifold $G/P_{\taumod}$. Below are several examples of this in the case when $G$ is a complex semisimple Lie group and $P_{\taumod}=P_{\simod}=B$, a Borel subgroup of $G$. 

\begin{example}\label{ex:Anosov1} 
Suppose that $G_1=SL(2,\C), G=SL(3,\C)$, $B<G$ is a Borel subgroup (e.g. the subgroup of upper triangular matrices), 
and $\rho: G_1\to G$ is the unique, up to conjugation, irreducible representation (thinking of $\bC^3$ as the space of homogeneous polynomials of degree two in two variables and letting $\rho$ acting on the variables). There is a $\rho$-equivariant holomorphic embedding $f: \P^1\to \P^2$, the {\em Veronese embedding}:
\[f:[x_0;x_1]\mapsto [x_0^2;x_0x_1;x_1^2].\]
The map $f$ lifts canonically to an equivariant holomorphic embedding 
$\tilde f: \P^1\to \cF= G/B$. Namely, if $f(z)= p\in \P^2$, take the unique projective line $\ell_p$ in $\P^2$ tangent to the curve $V=f(\P^1)$ at $p$. 
Then $(p, \ell_p)$ is an element of $\cF$ and we get 
$$
\tilde f(z)= (p, \ell_p), p=f(z). 
$$
If $\Ga< G_1$ is a convex-cocompact subgroup, then the subgroup $\rho(\Ga)< G$ is $B$-Anosov and its flag-limit set 
$\La=\La(\Ga)$ in $\cF$ is the image under $\tilde f$ of the limit set $\La_1(\Ga)\subset \P^1$. By abusing the notation, we will identify $\Ga$ and its image under $\rho$ in $G$.

Let $\Th\subset W$ be the unique balanced thickening, see \S \ref{sec:Thickenings}. Then the $\Ga$-action on $\Omega_{\Th}(\Ga)=\cF\setminus \Th \La$ is properly discontinuous and cocompact. Explicitly, one can describe $ \Th \La$ 
as follows. Pick a limit flag $\la=(p,\ell)\in \La(\Ga)\subset \cF$, i.e. $p\in V=f(\P^1)$ and $\ell$ is a line tangent to $V$ at $p$. Now, consider all flags of the form $(p, \ell')$ and $(q, \ell)$. The subset of such flags is the thickening $\Th(\la)$. 
It is a union of two smooth rational curves in $\cF$ meeting at $\la$. 
Then $\Th \La$ is the union of thickenings $\Th(\la)$, $\la\in \La$. Alternatively, let $w_1, w_2\in W$ be simple reflections. For each $\la\in \cF$ we have corresponding 1-dimensional Schubert varieties $X_{w_i}(\la)$, $i=1,2$. Then 
$$
\Th\La= \bigcup_{\la\in \La} (X_{w_1}(\la) \cup X_{w_2}(\la)). 
$$
If $\La_1\ne \P^1$ (i.e. $\Ga_1< SL(2,\C)$ is not a lattice), then  $\Th \La$ will have topological dimension $2$ 
(if $\Ga_1$ is a Schottky group) or $3$ (otherwise); it is (non-equivariantly) homeomorphic to the product $\La\times \Th(\la)$, where $\la$ is a chosen point in $\cF$. If $\Ga< SL(2,\C)$ is a lattice (i.e. $\bH^3/\Ga$ is compact) then $\La=\tilde V=\tilde f(\P^1)$. The space $\Th \La$ is a singular 
projective variety, the union of two smooth projective varieties which are (isomorphic) nontrivial 
$\P^1$-bundles over $\P^1$. 

The dynamics of the $\Ga$-action on $\Th\La$ is quite complicated. This $\Ga$-action  is not proper, the nonwandering set $NW(\Ga)$ for the $\Ga$-action on $\cF$ is contained in $\Th\La$ but is strictly smaller, unless 
$\Ga< G_1$ is a lattice. The nonwandering set can be described as follows. For every pair of distinct limit points $\la_\pm\in \La=\La(\Ga)$ there is a unique apartment $a(\la_+, \la_-)$ in $\tits \cX$ containing $\la_\pm$. Let $Ch(\la_+,\la_-)\subset \cF$ denote  the set of chambers in this apartment: There are six such chambers, two of which are $\la_\pm$. Then
$$
NW(\Ga)= \bigcup Ch(\la_+,\la_-)
$$
where the union is taken over all pairs of distinct elements $\la_\pm \in \La$. Equivalently, $NW(\Ga)$ is the closure of the set of points in $\cF$ fixed by elements of infinite order in $\Ga$. The $\Ga$-action on each ``component'' 
$$
\bigcup_{\la\in \La} X_{w_i}(\la), i=1, 2,$$
of $NW(\Ga)$ is ergodic. Equivalently, in this example, the nonwandering set is the closure of the set of fixed points of elements of infinite order of $\Ga$.  
\end{example} 

\begin{example} [Hitchin representations of hyperbolic surface groups to $SL(3,\R)$] \label{ex:hitchin}
{\kap This example, in a sense, lies at the origin of the theory of Anosov subgroups. An open nonempty subset $\Omega$ in the real projective plane $\R P^2$ is said to be 
{\em properly convex} 
if there exists a line $L\subset \R P^2$ disjoint from $\Omega$ such that $\Omega$ is a proper convex subset of the affine plane $\R P^2 \setminus L$. It's not hard to see that convexity does not depend on the choice of a line $L$. One similarly defines {\em strictly convex} subsets. 
Suppose that $\Ga$ is a discrete subgroup of $SL(3,\mathbb R)$ which preserves a strictly convex subset $\Omega\subset \R P^2$ and acts cocompactly on it. It turns out that such $\Ga$ is not a virtually abelian (due to the strict convexity assumption) and, thus, is nonelementary hyperbolic. Up to passing to a finite index subgroup, $\Ga$ is isomorphic to the fundamental group of a compact hyperbolic surface. Furthermore, the loop 
$C=\partial \Omega$ is $C^1$-smooth, hence,  admits a lift $\tilde C$ to the full flag-manifold of the group $G=SL(3,\R)$: For each $p\in C$ we take the unique projective line $\ell$ tangent to $C$ at $p$, hence, lift $p$ to the flag $(p,\ell)$. 
The subgroup $\Ga< G$ is Anosov and its flag-limit set is $\tilde C$. One way to obtain such subgroups $\Ga$ is via {\em bending}. Start with a discrete surface subgroup $\Ga_0< SO(2,1)< SL(3,\R)$, $\Ga\cong \pi_1(\Sigma)$, where $\Sigma=\mathbb H^2/\Gamma$. Pick a simple loop $c$ in $\Sigma$ separating $\Sigma$ in two subsurfaces $\Sigma_1, \Sigma_2$ of genus $\ge 1$. Then $\Ga_0=\pi_1(\Sigma_1)*_{\mathbb Z} \pi_1(\Sigma_2)$, where the amalgamated subgroup $\bZ$ is generated by an element $\ga\in \Ga$ representing the loop $c$. The centralizer $Z_G(\gamma)$ of $\gamma$ in $G$ is 2-dimensional, with identity component $H$ isomorphic to $\R^2$ (commensurable to a maximal torus in $G$). The subgroup $H$ contains a 1-parameter subgroup $H_1$ contained in $SO(2,1)$ (commensurable to a maximal torus in $SO(2,1)$). Take an element $\alpha\in H\setminus H_1$ and define a representation
$$
\rho_\al: \Ga_0\to \rho_\al(\Ga_0)=\Ga< G 
$$
so that the restriction to $\pi_1(\Sigma_1)$ is the identity and the restriction to $\pi_1(\Sigma_1)$ is given by conjugation by $\al$. It turns out that $\Ga$ is a subgroup preserving a strictly convex subset $\Omega\subset \R P^2$ as above and, moreover, its boundary is not a conic, i.e. $\Ga$ is not conjugate to a subgroup of $SO(2,1)$. For more details see e.g. \cite{Labourie-book}.} 
\end{example}

\begin{example}\label{ex:Anosov2} 
We embed $G_1=SL(2,\C)$ in $G=SL(3,\C)$ as the subgroup of block-diagonal matrices with $2\times 2$ and $1\times 1$ blocks. 
The group $SL(2,\C)$ preserves unique projective line $L$  in $\P^2$ and a unique point $q\in \P^2$. We regard $L$ as $\P^1=\geo X_1$.  The holomorphic embedding $f: \P^1\to L\subset \P^2$ also has a canonical holomorphic lift $\tilde f: \P^1\to \cF=G/B$. Namely, for 
each $p\in L$ we take $\ell_p$, the unique projective line through $p, q$. This gives us a flag $(p, \ell_p)\in \cF$ and we set  
$$
\tilde f(z)= (p, \ell_p), p=f(z). 
$$
If $\Ga< G_1$ is a convex-cocompact subgroup, then its image under the embedding  $G_1\to G$ is Borel-Anosov and its flag-limit set in $\cF$ is the image under $\tilde f$ of the limit set $\La_1(\Ga)\subset L=\P^1$. Suppose that $\Ga$ is a uniform lattice in $G_1$; equivalently, $\Ga< G_1$ is convex-cocompact and $\La_1(\Ga)= L=\P^1$. Then the nonwandering limit set of $\Ga$ in $\cF$ is the balanced thickening of $\tilde f(L)$. Thus, in this example, the nonwandering set is not the closure of the set of fixed points of infinite order elements of $\Ga$.    
\end{example}

 \begin{example}\label{ex:Anosov3} 
 Let $\Gamma< G_1=\SU(2,1)< \SL(3,\C)$, where  $\Gamma< \SU(2,1)$ is a convex-cocompact subgroup. Now, $\Ga$ does not preserve any complex curves 
 in $\P^2$. However,  $G_1$ preserves an open ball $\mathbb B\subset \P^2$, a model of the complex-hyperbolic 2-dimensional space, the symmetric space of $G_1$. The boundary of this ball is, of course, also $G_1$-invariant. For every $p\in \partial \mathbb B$ we let $\ell_p$ denote the unique complex-projective  line 
 tangent to  $ \partial \mathbb B$ at $p$. The set of flags
 $$
 (p, \ell_p)\in \cF, p\in  \partial \mathbb B
 $$
 is a submanifold $\Sigma$ in $\cF$ diffeomorphic to the 3-dimensional sphere $S^3=  \partial \mathbb B$. This also defines an equivariant topological embedding 
 $\tilde f: \La_1(\Ga)\to \Sigma$. The image of $\Ga$ under the embedding $\SU(2,1)\to \SL(3,\C)$ is a Borel--Anosov subgroup whose flag-limit set is 
 $\tilde f(\La_1(\Ga))$. 
 \end{example}

\begin{figure}[htbp]
   \centering
   \includegraphics[scale=0.7]{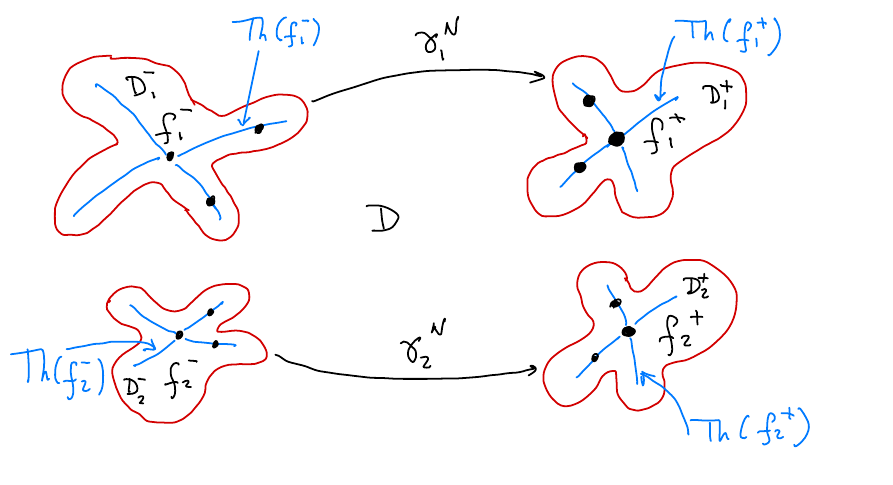} 
   \caption{Anosov-Schottky subgroup of $\SL(3,\C)$. Dots denote fixed flags of the elements $\ga_1, \ga_2$.}
   \label{fig:F2}
\end{figure}

\begin{example}
{\kap 

1. {\bf Anosov--Schottky groups.}  Suppose that $\ga_1,...,\ga_r\in G=\SL(n,\C)$ are regular loxodromic elements 
in ``general position:'' If $E_{\la_j}, E_{\la_j}$ are 
eigenspaces of $\ga_i, \ga_j$ of complementary dimensions, then they span $\C^n$. Then there exists $k_0\in \bN$ such that for all $k\ge k_0$ the elements 
$$
\ga_1^k,...,\ga_r^k
$$
generate a free subgroup of rank $r$ which is an Anosov subgroup  $\Ga< G$. See e.g. \cite{KLP14, Dey-Kapovich-Leeb, KLP25}. 

2. {\bf Anosov--Schottky subgroups of $G=\SL(3,\C)$.} Let $\ga_1, \ga_2$ be two regular loxodromic elements {\em in general position}: Their respective attractive/repelling fixed 
projective flags $f^\pm_1=(p^\pm_1,\ell^\pm_1), (p^\pm_2,\ell^\pm_2)$ are such that $p^+_1\notin \ell^\pm_2, p^+_2\notin \ell^\pm_1, p_1^-\notin \ell^\pm_2, p_2^-\notin \ell^\pm_1$. Let $\Th$ be the unique balanced thickening in the Weyl group $W=S_3$ of $G$: It consists of 
the trivial element and two simple reflections. The general position assumption on the attractive/flags implies that the thickenings
$$
\Th(f_1^\pm), \Th(f_2^\pm)
$$
are pairwise disjoint subvarieties of the full flag-manifold $\cF$ of $G$. Thus, take disjoint 
smooth compact submanifolds with boundary $D^-_k, i=1,2$ containing in their interiors, respectively, the thickenings 
$\Th(f_1^-), \Th(f_2^-)$, and disjoint from  $\Th(f_1^+)\cup \Th(f_2^+)$. Then there exists $N>0$ such that the complements 
$$
D_1^+:=\cF\setminus \ga_1^N(\inte D_1^-), D_2^+:=\cF\setminus \ga_2^N(\inte D_2^-)
$$ 
are disjoint and disjoint from $D_1^-\cup D_2^-$. See Figure \ref{fig:F2}.

Then the subgroup $\Ga<G$ generated by $\ga_1^N, \ga_2^N$ is free of rank two, Anosov and 
$$
D:= \cF \setminus \left(D_1^-\cup D_2^- \cup D_1^+\cup D_2^+ \right)
$$
is the fundamental domain for its action on the domain of discontinuity
$$
\Omega_{\Th}(\Ga)= \cF \setminus \Th\La(\Ga) 
$$
of the $\Ga$-action on $\cF$. 

3. {\bf Nori--Schottky subgroups of $G=\SL(4,\C)$.} These examples were first introduced by M.~Nori in \cite{Nori} 
who was interested in constructing complicated examples 
of compact non-K\"ahler complex 3-dimensional manifolds. 
To make things very concrete, take $\ga_1, \ga_2\in \SL(4,\C)$ as above and consider their action not on the full 
flag-manifold (which is hard to visualize) but on $\P^3$. Let $C^\pm_j\subset  \P^3, j=1,2$, be the (pairwise 
disjoint) projective lines corresponding to the spans of the two highest (for $C_j^+$) and lowest (for $C_j^-$) eigenspaces of $\ga_j$. Then
for every compact $C\subset \P^3\setminus C_j^-$  the sequence of subsets $(\ga_j^k(C))$ accumulates at $C_j^+$. Thus, there exist 
$D_j^\pm, j=1,2$,  pairwise disjoint closed tubular  neighborhoods of  $C_j^\pm$ such that (for some large $k$) $\ga_j^k(\mathrm{ext} D_j^-)= \inte D_j^+$, $j=1,2$. 
The complement to the union of interiors of $D_j^\pm, j=1,2$ is a {\em Schottky fundamental domain} $D$ of the group $\Ga< G$, generated by $\tilde{\ga}_1:=\ga_1^k, \tilde{\gamma}_2:=\ga_2^k$. 
Considering the  $\Ga$-orbit of $D$, we see that the union  
$$
\bigcup_{\ga\in\Ga} \ga D
$$
is an open subset $\Om\subset \P^3$ on which $\Ga$ acts properly discontinuously and cocompactly. 
The complement $\P^3\setminus \Om$ is the {\em Kulkarni-limit set} $L$ for the action of $\Ga$ on $\P^3$. The limit set $L$, topologically, is the product
$$
\P^1\times \geo F_2
$$ 
of the projective line and the Cantor set which is the Gromov-boundary of the free group $F_2$. But the action of $\Ga$ on $L$ is not a product action. 
See \cite{CNS} for further discussion of these interesting examples. 
It would be interesting to study the dynamics of the Nielsen maps $N$ of such groups $\Ga$, 
in particular, equilibrium measures.  
}
\end{example}

\begin{rem}
Unlike the usual Schottky subgroups of $\PSL(2,\C)$, the boundary of the fundamental domain $D$ for a 
Nori--Schottky group  is not a disjoint union of topological spheres: Each boundary component is a nontrivial 
$S^3$-bundle over $S^2$. But even in the case of free convex-cocompact subgroups of isometries of $\bH^4$, one should not expect  a fundamental domain in $S^3$ bounded by topological spheres, see \cite{BC, Matsumoto}. One gets a similar fundamental domain $\cD$ for the action of $\Ga$ on the full flag-manifold $\cF=G/B$. Note that 
the full flag-manifold fibers ($G$-equivariantly) over $\P^3$. Fibers of this fibration are copies of the 3-dimensional flag-manifold for the group $\SL(3,\C)$ 
Accordingly,  $\cD$ fibers over $D$,  via the restriction of this fibration $\cF\to \P^3$. See a schematic picture (over $\R$) in Figure \ref{fig:Nori-Schottky}. 
\end{rem}

\begin{figure}[htbp] 
   \centering
   \includegraphics[width=3in]{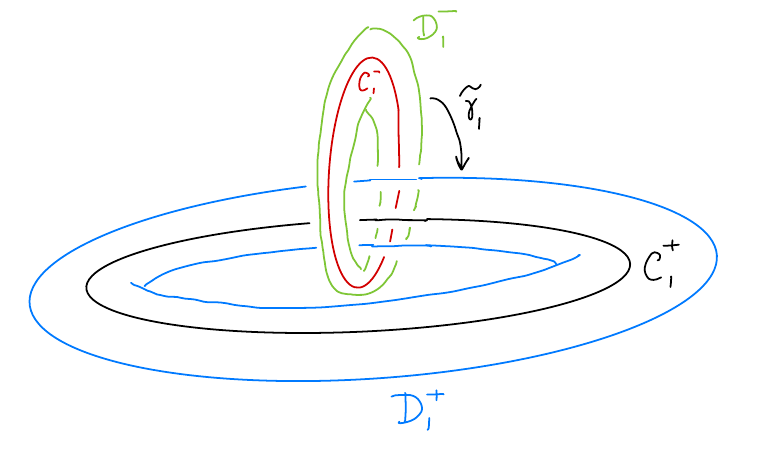} 
   \caption{Nori--Schottky fundamental domain.}
   \label{fig:Nori-Schottky}
\end{figure}

\begin{example}\label{exAnosov4}
Suppose that $G_1,...,G_n$ are rank 1 Lie groups, $\rho_i: \Gamma\to G_i$ are Anosov representations of a nonelementary hyperbolic group $\Gamma$ (i.e. representations with finite kernels and convex-cocompact images). For each $G_i$ let $P_i<G_i$ denote a (maximal) parabolic subgroup. Define
$$
G:= G_1\times ...\times G_n, P:= P_1\times ....\times P_n, \cF:= \prod_{i=1}^n G_i/P_i
$$
and 
$$
\rho=\rho_1\times ...\times \rho_n: \Gamma\to G,
$$
the product representation. Then $\rho$ is $P$-Anosov; the subgroup $P< G$ is minimal parabolic. To describe the flag-limit set of $\rho(\Ga)$ in $\cF$, 
note that for every $i$ we have an equivariant homeomorphism $f_i: \geo \Ga\to \La(\rho_i(\Ga))\subset G_i/P_i$. Set
$$
f= (f_1,...,f_n): \geo \Ga\to \cF. 
$$ 
Then the image of $f$ is the flag-limit set of $\rho(\Ga)$.  
\end{example}

More examples of Anosov subgroups of semisimple Lie groups arise as perturbations of 
Anosov subgroups, see e.g. \cite{Guichard-Wienhard} or \cite{KLP14, KLP25}: 

\begin{thm}
[Stability Theorem] \label{thm:stable} 
Suppose that $\Ga< G$ is a $\taumod$-Anosov subgroup of a semisimple Lie group $G$ and $\rho_i: \Ga\to G$ is a sequence of representations converging to the identity representation, i.e. the sequence of subgroups $\Ga_i=\rho_i(\Ga)$ converges {\em algebraically} to $\Ga$. Then for all sufficiently large $i$ the following hold:

1. Every $\rho_i:  \Ga\to \Ga_i< G$ is a $\taumod$-Anosov representation. 

2. There exists an equivariant homeomorphism $f_i: \La(\Ga)\to \La(\Ga_i)$ between flag-limit sets in $G/P_{\taumod}$. 

3. The sequence of continuous maps $f_i$ converges to the inclusion map $\La(\Ga)\to G/P_{\taumod}$.  
\end{thm}

\begin{example}
Let $\Ga_0< \SL(2,\C)$ be a convex-cocompact torsion-free subgroup with {\em nonempty} domain of discontinuity. To make things interesting, we assume that the limit set of $\Ga$ in $S^2$ is not contained in a circle. The quotient 
$M=\bH^3/\Ga_0$ admits a natural compactification as a compact 3-dimensional manifold with nonempty boundary. Let $\chi$ denote the Euler characteristic of $M$. Thus, $\chi<0$. One then considers the representation variety 
$X(\Ga_0,G):=Hom(\Ga_0, G)//G$ parameterizing equivalence classes of representations $\Ga_0\to G=\SL(3,\C)$. This variety contains a distinguished point $[\rho_0]$ corresponding to the inclusion $\rho_0: \Ga_0\to \SL(2,\C)$ followed by the unique (up to conjugation) irreducible representation $\SL(2,\C)\to \SL(3,\C)$. It turns out that the variety $X(\Ga_0,G)$ is smooth at $[\rho_0]$, of local dimension $-8\chi$. Some of the equivalence classes of representations near $\rho_0$ come from representations $\Ga_0\to \SL(2,\C)$, but this subset forms a smooth subvariety of dimension $-3\chi$, which is much less than $-8\chi$. The representations near $\rho_0$  are Anosov. One can think of these representations as complex analogues of Hitchin representations discussed in Example \ref{ex:hitchin}. Note, however, that unlike in the Hitchin case, it is critical here to assume that $\bH^3/\Ga_0$ is noncompact, otherwise one cannot nontrivially deform $\rho_0$. 
\end{example}


\medskip 
Further examples of Anosov subgroups are obtained via Combination Theorems, i.e. as free products, amalgamated free products and HNN extensions of Anosov subgroups, see \cite{KLP14}, \cite{Dey-Kapovich-Leeb}, \cite{Dey-Kapovich1}, \cite{Dey-Kapovich2}. For 
other constructions of Anosov subgroups we refer the reader to the recent work \cite{DFWZ} and references therein. 

Below is an interesting example of a {\em non-Anosov} subgroup illustrating the necessity of various conditions in the characterizations of Anosov subgroups; the example is taken from \cite[\S 8.11]{Kapovich-Sardar}. 
Let $\Ga$ be a group isomorphic to the fundamental group of a compact hyperbolic surface $S$. Take two representations
$$
\rho_i: \Ga\to G_i\cong \SL(2,\C)
$$
which are isomorphisms of $\Ga$ to {\em singly-degenerate} Kleinian subgroups $\Ga_i< G_i$ such that:

a. Both hyperbolic manifolds $\bH^3/\Ga_i$ have injectivity radius bounded away from zero. 

b. Ending laminations $L_i\subset S$ of the groups $\Ga_i$ are transverse to each other. 

For instance, we can take $L_i$ to be stable/unstable laminations of a pseudo-Anosov homeomorphism of $S$. Now, take 
$$
\rho: \Ga\to G=G_1\times G_2, \rho(\ga)=(\rho_1(\ga), \rho_2(\ga)). 
$$
Then each orbit map of the $\rho$-action on the symmetric space of $G$, $\cX=\bH^3\times \bH^3$, is a quasi-isometric embedding $\Ga\to \cX$ and  
every nontrivial element of $\Ga$ maps to a loxodromic isometry of $\cX$ via $\rho$. In particular, the locally-symmetric space $\cX/\rho(\Ga)$ has injectivity radius bounded from below. Furthermore, the representation $\rho$ is regular with respect to all parabolic subgroups of $G$. However, $\rho$ is not $\taumod$-Anosov for any simplex $\taumod$: Uniform regularity fails.  

Let $f_i: \geo \Ga=S^1\to \La_i=\La(\Ga_i)\subset S^2$ be Cannon-Thurston maps of the groups $\Ga_i, i=1,2$. Then the orbit map $\Ga\to \cX$ 
has a continuous extension to a topological embedding 
$$
f: S^1\to \La_1\times \La_2\subset S^2\times S^2, f(\xi)=(f_1(\zeta), f_2(\zeta)), \zeta\in S^1.
$$
While injective, the map $f$ is not antipodal: Whenever $\zeta, \zeta'\in S^1$ are distinct points such that, say, $f_1(\zeta)=f_1(\zeta')$, the points
$f(\zeta), f(\zeta')\in S^2\times S^2$ are distinct but not antipodal to each other, since their first coordinates coincide. This happens precisely when $\zeta, \zeta'$ 
are either ideal end-points of a common leaf of $\tilde L_1$, the lift of $L_1$ to the hyperbolic plane, or $\zeta, \zeta'$ belong to the ideal boundary 
of a connected component of $\bH^2\setminus L_1$. The subset of $S^1$ consisting of such points $\zeta$ has zero Hausdorff dimension 
(see \cite{Birman-Series}). The same, of course, holds for $f_2$ and $L_2$. Therefore, $f$ is antipodal on a co-null subset of $S^1$ (more precisely, 
its complement has zero Hausdorff dimension). In other words, $f$ is a.e. antipodal.  Due to continuity of the map
$$
\Ga \cup \geo \Ga\to \cX \cup \cF,
$$
the image of $f$ is the flag-limit set of $\rho(\Ga)$ in $\cF$.

\section{Gibbs measures and currents associated with Anosov subgroups} \label{sec: Anosov Gibbs}

\begin{defi}\label{def:Anosov Gibbs}\index{Gibbs measure of Anosov subgroup}
Let $\rho: \Ga\to G$ be a $\taumod$-Anosov representation of a hyperbolic group $\Ga$ to a semisimple Lie group $G$ and let 
$f: \geo \Ga\to \cF= G/P_{\taumod}$ denote the boundary map of $\rho$. Suppose that $\beta$ is a hyperbolic length function on $\Ga$ and $\mu_{\Ga}$ 
the corresponding Gibbs measure on $\geo \Ga$. 
Then the push-forward measure $f_*(\mu_{\Ga})$ is the Gibbs measure $\mu_{\cF}$ of $\rho(\Ga)$. 
\end{defi}

We also have the following alternative description of this Gibbs measure. 
Let $(s_j)$ denote the PS sequence of the Gibbs measure $\mu_{\Ga}$ and consider 
 a $\Ga$-orbit $\Ga x\subset \cX$. For every $z\in Z=\Ga x\cup \La_{\taumod}(\Ga)$ 
 and $s>s_0$ we define a measure 
$$
\mu_{z}(s) = \frac{1}{\cP(s)}\sum_{\gamma\in \Gamma}e^{-s\beta(\gamma)}\delta_{\rho\gamma(z)} 
$$
on $Z$. Since the orbit map $o_x: {\Ga}\to \Ga x$ extends continuously to the boundary homeomorphism 
 $f: \geo \Ga\to \La_{\taumod}(\Ga)$, we can apply according to Corollaries \ref{cor:CT} and \ref{cor:Gibbs} and conclude 
 that  for every choice of $z\in Z$ there exists a limit measure
$$
\lim_{s_j\to s_0} \mu_{z}(s)
$$
on $\La_{\taumod}(\Ga)$ and this measure is exactly  the Gibbs measure $\mu_{\cF}$ of $\rho(\Ga)$. 
 
The equidistribution theorem, Theorem \ref{thm:main}, is a far-reaching generalization of this basic observation.  

Recall that $\cD(\cF)$ denotes the space of {\em currents} on $\cF$, see Section \ref{sec:currents}, and it is endowed with the weak topology. 

\begin{lem}\label{lem:weakly continuous}
Suppose that $G$ is a complex semisimple Lie group, consider the flag manifold $\cF= G/P_{\taumod}$. 
Let $\Th\subset W$ be a thickening. Then the map $\theta: \cF \to {\mathcal D}(\cF)$ sending $x$ to 
$\llbracket \Th(x)\rrbracket$ is continuous.
\end{lem}
\begin{proof} Pick a point $\tau_+\in \cF$ and its opposite simplex $\tau_-$. The unipotent radical $U< P_{\tau_-}$ 
acts simply transitively on an open neighborhood $\Opp(\tau_-)$ of $\tau_+$ in $\cF$. 
Therefore in order to verify continuity of $\theta$ at $\tau_+$, 
it suffices to check weak continuity of the map 
$$
U\to {\mathcal D}(\cF), g\mapsto \llbracket \Th(g\tau_+)\rrbracket =g_* \llbracket \Th(\tau_+)\rrbracket . 
$$
Let $k$ denote the complex dimension of $\Th(\tau_+)$; take a form $\varphi\in \Om^{k,k}(\cF)$. Then the map 
$$
g\mapsto g^*\varphi
$$
is continuous (where, as before, we equip $ \Om^{k,k}(\cF)$ with $C^0$-topology). Hence, we obtain continuity of the map 
$$
g\mapsto \int_{\Th(\tau_+)} g^*\varphi. 
$$
\end{proof}

We will use this lemma to define {\em Gibbs currents} on $\cF$. Pick a thickening $\Th$ in the Weyl group $W$. Let $\La=\La(\Ga)$ be a 
flag-limit set of an Anosov subgroup $\Ga< G$. \index{Gibbs current}
Then, by the lemma, 
we get a continuous map
$$
\la\mapsto \llbracket \Th(\la)\rrbracket, \La \to  \cD(\cF). 
$$
Integrating against the Gibbs measure $\mu$ on $\La$,
 we obtain the definition of  Gibbs currents on $\cF$. 
 
 \begin{defi} \label{defi_Gibbs_current} A bidimension $(k,k)$ {\em Gibbs current} on $\cF$ is a linear combination of currents of the form: 
 $$T=\int_\La \llbracket \Th(\la) \rrbracket d\mu(\la),$$
 where $\mu$ is a Gibbs measure on $\Lambda$ and where $\Th$ is a thickening in the Weyl group corresponding to a $k$-dimensional subvariety in $\cF$. 
 In the special case when $\Th=\Th_w$, we will use the notation $T_w$ for the corresponding Gibbs current $T$. Thus, every Gibbs current is a linear combination of the 
 Gibbs currents $T_w$, $|w|=k$. 
 \end{defi}

The current $T$ is a bit more geometric in the case of {\em slim} thickenings $\Th$. 
(For instance, this is the case when $\Th=\Th^k, k< n/2$.)
In this case $\Th(\la_1)\cap \Th(\la_2)=\emptyset$ 
whenever $\la_1\ne \la_2\in \La$ and, thus, the current $T$ is a {\em measured lamination} supported on $\Th(\La)$. The latter is a bundle over $\La$ with 
typical fiber $\Th(x)$. 

\section{Redressing}

\subsection{Redressing in Lie groups}

For the following discussion we recall that $\lambda$ and $\mu$ denote, respectively, the Lyapunov spectrum and the Cartan projection in a semisimple Lie group $G$ defined in Section \ref{sec:Cartan and Lyapunov}.  
Theorem 2.34 in \cite{GGKW} attributes the following theorem to Y.~Benoist \cite{Benoist} whose proof depends critically on a theorem of Abels, Margulis and Soifer \cite{AMS} on proximal elements of Zariski dense subgroups of reductive groups:

\begin{thm}\label{thm:GGKW}
Let $\Ga< G$ be a finitely generated subgroup of an algebraic semisimple Lie group $G$ with reductive Zariski closure in $G$. Then there exists a constant $C$ (depending on $\Ga$) and a finite subset $E\subset \Gamma$ such that for every $\gamma\in \Gamma$ there exists $g\in E$ satisfying 
$$
||\lambda(\gamma\circ g)- \mu(\gamma\circ g)||\leqslant C.
$$
\end{thm}



\subsection{Redressing in Anosov subgroups}\label{sec_redressing_anosov}

Combining inequalities \eqref{eq:A2} with Theorem \ref{thm:GGKW}, we obtain:

\begin{prop}\label{Mishaseigenvaluelemma}
Suppose that $\Gamma< G$ is an Anosov subgroup and $\beta$ a hyperbolic length function on $\Gamma$. 
Then there are constants  $C\geqslant 1, A\geqslant 0$ and a finite subset $E\subset \Gamma$, depending only on $\Gamma$ 
(regarded as a subgroup of $G$) such that for each $\gamma\in \Gamma$ there exists $g\in E$ such that $\tilde\ga=\ga\circ g$ satisfies 
$$
C^{-1} \beta(\tilde\gamma)- A \leqslant 
\la_i(\tilde\gamma)\leqslant C \beta(\tilde\gamma)+ A 
$$
for all $i=1,...,r$. Here, as always, $r$ is the rank of the group $G$. 
\end{prop}
\proof 1. Since $\Ga$ is Anosov and $\beta$ is a hyperbolic length function on $\Ga$, there are constants $C_1, A_1$ such that
$$
C^{-1}_1 \beta(\gamma) - A_1 \leqslant \mu_i(\gamma)\leqslant C_1 \beta(\gamma) + A_1
$$
for all $i=1,...,r$ and $\gamma\in \Gamma$ (see \eqref{eq:A2}). Therefore, if $\rho(\Ga)$ has reductive closure, 
 Theorem \ref{thm:GGKW} implies existence of a finite subset $E\subset \Ga$ and constants $C, A$ satisfying the assertion of the proposition. 
The only thing that has to be explained is how to eliminate the assumption on reductivity of the Zariski closure coming from Theorem \ref{thm:GGKW}. Let $\Gamma^{red}$ denote the Levi projection of $\Gamma$ to the reductive part of the Zariski closure of $\Gamma< G$. Then  $\Gamma^{red}$ is again Anosov (see \cite{GGKW}) and the projection preserves the Lyapunov spectra of elements of $\Gamma$ and is 1-1. By combining this 
with Theorem \ref{thm:GGKW} applied to the subgroup $\Gamma^{red}$, 
one obtains the existence of the required constants $C, A$ and a 
finite subset $E^{red}\subset \Gamma^{red}$. 
Taking the preimage $E\subset \Ga$ of  $E^{red}$ 
and using the same constants $C, A$, 
we obtain the proposition.

\medskip 
2. We also give a direct proof of the proposition using Theorem \ref{thm:redressing} on redressing in hyperbolic groups and not relying upon Theorem \ref{thm:GGKW}. According to Theorem \ref{thm:redressing}, since $\Ga$ is a hyperbolic group with Gromov boundary equivariantly homeomorphic to the flag-limit set $\La\subset {\mathcal F}=G/P_{min}$, there exists a finite subset $E\subset \Ga$ and $\eps_0>0$, such that for every $\gamma\in \Gamma$  there exists $g\in E$ such that $\tilde\gamma = \gamma\circ g$ has 
distinct fixed points $x^\pm$ in $\La$, one of which is attractive and the other repelling, such that 
\begin{equation}\label{eq:separation} 
d_{\cF}(x^+, x^-)\geqslant \eps_0.\end{equation}
 In particular, $\tilde\gamma$ is a loxodromic element of $G$ preserving the unique maximal flat $F$ in $\cX$ asymptotic to both $x^\pm$. In view of the inequality \eqref{eq:separation},  the flat $F$ passes within uniformly bounded distance from the base-point $o\in \cX$.  Thus, both points $o, \tilde\gamma(o)$ are uniformly close to $F$. The triangle inequality for $d_\Delta$ now implies that 
$$
||\la(\tilde\ga)- \mu(\tilde\ga)||\leqslant C_2
$$ 
for a uniform constant $C_2$ depending only on $\epsilon_0$. The rest of the proof is the same as in 1. \qed

\subsection{Application to separation properties of elements of Anosov subgroups} 

In view of Lemma \ref{lem:separate}, Proposition \ref{Mishaseigenvaluelemma} can also be formulated in terms of $\eps$-separation in the sense of 
 \eqref{e_separation_inequality}:

\begin{thm}\label{Mishaseigenvaluelemma1}
Suppose that $\Gamma< G$ is an Anosov subgroup, $\beta$ is a hyperbolic length function on $\Ga$. 
Then there are constants $\varepsilon>0, C>0, A>0$ and a finite subset $E\subset \Gamma$, depending only on $\Gamma$ such that for each $\gamma\in \Gamma$ there exists $g\in E$ satisfying:

1.  $\tilde \gamma= \ga\circ g$ is $\varepsilon$-separated.

2. For all $i=1,...,r$, 
$$
C^{-1} \beta(\tilde\gamma)- A \leqslant \la_i(\tilde\gamma)\leqslant C \beta(\tilde\gamma)+ A.  
$$
\end{thm}
\begin{proof} It remains to note that the inequality \eqref{eq:separation} implies $\varepsilon$-separation of the attractive/repelling 
fixed points $x^\pm\in \cF$ of $\tilde\ga$ 
in the sense of \eqref{e_separation_inequality}, where $\eps_0$ is as in the 2nd part of the proof of Proposition \ref{Mishaseigenvaluelemma}, 
and $\varepsilon$ depends only on $\eps_0$, see Lemma \ref{lem:separate}. \end{proof} 

\medskip 
Just as in Section \ref{sec:redressing}, we get an addendum to Theorem \ref{Mishaseigenvaluelemma1}. Let $\Ga\acts X$ be a geometric 
action of $\Ga$ on a $\delta$-hyperbolic geodesic metric space $X$. Let $x_0\in X$  be a base-point and $U(\zeta, t)$ be the neighborhoods in $\ol{X}$ 
of points $\zeta\in \geo X$, as defined in \eqref{eq:U-nbd}. Then, in addition to the conclusion of Theorem \ref{Mishaseigenvaluelemma1}, we have:

\begin{lem}\label{Mishalimitpointslemma4}
In the setting of the previous theorem, there exists a constant $\kappa$ depending only on $X$ and the $\Ga$-action on $X$ such that 
 if $\gamma(x_0)\in U(\zeta,t)$ for some $\zeta \in \partial \Gamma$, then the attracting fixed point of $\tilde{\gamma}$ in $\geo X$ 
 is in $U(p, \kappa t)$. 
\end{lem}
\begin{proof} This is an immediate consequence of Theorem \ref{Mishaseigenvaluelemma1} and Theorem \ref{thm:redressing2}. 
\end{proof}

\subsection{Application to currents on flag manifolds}

The redressing results above allow one to transfer the properties of elements of a fixed maximal torus $T< G$ 
to general elements of the subgroup $\Ga< G$. The following corollary illustrates 
the primary way in which the redressing is utilized in our book.

\begin{cor} \label{cor_equivariant_continuous} Fix an Anosov subgroup $\Gamma< G$ of a complex semisimple Lie group $G$, a hyperbolic length function $\beta$ on $\Gamma$, a finite subset $E\subset \Gamma$ as in Theorem \ref{Mishaseigenvaluelemma1} and a boundary map $\xi \colon \partial_\infty \Gamma \to \Lambda$ onto the limit set $\Lambda$ in the flag manifold $\cF=G/B$ of dimension $n$. Let $\Gamma^*$ be a subset of $\Gamma$. Let $\psi$ be a closed bidegree $(n-k,n-k)$-current on $\cF$ whose homology class equals 
$[\psi] = \sum_{w\in W} a_w(\psi) [\Th_w]$, where $a_w \in \C$. 
    Assume that for any $\epsilon > 0$ and any smooth $(k,k)$-form $\varphi$ on $X$, there exists $M>0$ such that for any $g\in E$, for any $\tga\in\Gamma^*g$ whose Lyapunov spectrum satisfies $\lambda_i(\tga) \geqslant M $ for all $i$, one has: 
\begin{equation*}
	|\langle  \tga_* \left((g^{-1})_*\psi\right) -  \sum_{w\in W} a_w(\psi)  \intcur{\Th_w(x^+)} , \varphi  \rangle    | \leqslant \epsilon ,
\end{equation*}
where $x^\pm\in \cF$ are the attracting and repelling flags of $\tga$. 
Then the function $\eta \colon \Gamma \cup \partial_\infty \Gamma \to \mathcal{D}^{k,k}(\cF)$ 
defined by: 
\begin{equation}
	\eta (\ga ) = \left \lbrace  \begin{array}{ll}
		\ga_* \psi & \text{if } \ga\in \Gamma \\ 
		\sum_{w\in W } a_{w}(\psi) \intcur{\Th_w(\xi(\ga))} & \text{if } \ga\in \geo\Gamma
	\end{array}\right. 
\end{equation}
is weakly continuous when restricted to the closure of $\Gamma^*$ in $\ol{\Ga}=\Gamma\cup\geo \Gamma$. 
\end{cor}

\begin{rem} 
Later when we apply this corollary to equidistribution of currents, working with subsets $\Ga^*\subset \Ga$ (instead of the entire $\Ga$) 
is necessary in the case of singular currents. However, for proving equidistribution of smooth forms we use $\Ga^*=\Ga$.
\end{rem}

\begin{proof}
The continuity of the restriction of $\eta$ on $\Gamma$ is immediate, while the continuity of the restriction  on $\geo \Gamma$ follows from Lemma \ref{lem:weakly continuous}. Consider 
$p\in \geo \Gamma\cap \overline{\Gamma^*}$. We need to prove that $\eta$ is continuous at $p$ 
when restricted to $\overline{\Gamma^*}$. 

Fix a closed current $\psi$ (as in the corollary), a smooth $(k,k)$-form  $\varphi$ on $\cF$, the corresponding constant $M$ (given by the hypothesis of the corollary) and an auxiliary number $m> 0$. Set $a_w:= a_w(\psi)$, $w\in W$. 
By Theorem \ref{Mishaseigenvaluelemma1}, for every $\gamma \in \Gamma$, there exists an element $g\in E$ such that the repelling and attracting fixed points  
of $\tilde{\gamma} = \gamma \circ g$ are $\epsilon$-separated and   
\begin{equation*}
 \lambda_i (\tilde{\gamma}) \geqslant C^{-1} \beta ({\tilde{\gamma}}) - A, i=1,...,r, 
\end{equation*}
for some uniform constants $A,C>0$.

Note that, if $M_2:=C M  + \max_{g\in E} \beta(g) + A $ and 
$\beta(\gamma) \geqslant M_2$,  then 
\begin{align*}
\lambda_i(\tilde{\gamma}) &\geqslant C^{-1} \beta(\tilde{\gamma}) - A 
 \geqslant C^{-1} \beta(\gamma) - C^{-1} \beta(g) -A  \geqslant M, i=1,...,r. 
\end{align*}

We will choose three real numbers $\delta_1,\kappa,\delta_2$ as follows. 
First, since the length function $\beta$ is proper, there exits $\delta_1>0$ such that if 
$\gamma \in U(p,\delta_1)$, then $\beta(\gamma)\geqslant M_2$. 
Next, by Lemma \ref{Mishalimitpointslemma4}, if $\gamma \in U(p,\delta_3)$ for some 
$\delta_3\dyl$, then the attracting fixed flag $x^+$ of $\tga$ satisfies 
$$
x^+\in \xi(U(p,\kappa \delta_3)),$$ 
where the constant $\kappa$ depends only on $\Gamma$.

By Lemma \ref{lem:weakly continuous}, 
there is a $\delta_2>0$ such that if $x\in \xi\big(U(p,\kappa \delta_2)\cap \geo\Gamma\big)$, then 
\begin{equation}\label{close_thickenings_inequality}
\big\lvert\langle  \sum_{w\in W} a_w\llbracket\Th_{w}(x)\rrbracket - \sum_{w\in W} a_w \llbracket\Th_{w}(\xi(p))\rrbracket , \varphi \rangle \big\rvert \leqslant \frac{m}{2}.
\end{equation}

Now take $\delta_3:= \min(\delta_1,\delta_2)$ and consider $\gamma\in U(p,\delta_3)\cap \Gamma^*$ and its redressing $\tilde{\gamma} = \gamma \circ g$ as before. Note first that $\tga\in \Gamma^*g$, $\tga_*\left((g^{-1}_*(\psi))\right)=\ga_*\psi$ and $\lambda_i(\tga)\geqslant M$ (for all $i=1,...,r$) because $\delta_3\geqslant\delta_1$. Hence, 
 from the hypothesis of Corollary  \ref{cor_equivariant_continuous} we obtain 
  the inequality:
   
\begin{equation}
\label{stupid_inequality_number_one}
\big\lvert\langle \ga_*\psi - \eta(x^+), \varphi \rangle \big\rvert \leqslant \frac{m}{2}
\end{equation}
where $x^+\in \cF$ is the attracting fixed flag of $\tga$. Secondly, we have $x^+\in \xi (U(p,\kappa \delta_3))$ and by \eqref{close_thickenings_inequality} we obtain that 
\begin{equation}\label{stupid_inequality_number_two}
\lvert\langle  \eta(\xi^{-1}(x^+)) -\eta(p) , \varphi \rangle \rvert \leqslant \frac{m}{2}.
\end{equation} 
Finally, the inequalities \eqref{stupid_inequality_number_one} and \eqref{stupid_inequality_number_two} imply that for every $\gamma\in U(p,\delta_3)\cap \Gamma^*$, we have
\[
\lvert\langle   \eta(\gamma) -\eta(p) , \varphi \rangle \rvert \leqslant m.
\]
The weak continuity at $p$ follows.
\end{proof}

\chapter{Equidistribution of subvarieties on flag manifolds}\label{sec:Equidistribution of subvarieties}

\section{The toy case}\label{sec:toy}

In this and the next section we illustrate the main ideas of the proof of Theorem \ref{thm:main} under certain simplifying assumptions. 
The main difficulty that we will have to deal with,in general, is the lack of continuity of the map $\eta$ defined in \eqref{eq:eta-map}. 
In this section, we will simply assume that  $\eta$ is continuous.  In fact, there is nothing special here about flag manifolds and the proof works for actions on general projective varieties. 


Let $Y$ be a smooth $n$-dimensional com\-plex-projective variety. We let $\Gamma\acts Y$ be an action of a hyperbolic group $\Ga$ by algebraic automorphisms of $Y$. 
In the case when $Y=\cF$, a flag-manifold, the current $T$ in the next proposition plays the role of a Gibbs current defined in Section \ref{sec: Anosov Gibbs}.

\begin{prop} \label{prop:poincare_series_current}
 Fix an integer $k$, $1 \leqslant k \leqslant n$. Let $\beta$ be a hyperbolic length function on $\Ga$ with critical exponent $s_0$. We fix a PS sequence $(s_m)$ converging from above to $s_0$, and the corresponding Gibbs measure $\mu=\mu_{\Ga}$ on $\geo \Ga$. 
 Then for every continuous map $\theta \colon \ol{\Ga}= \Gamma \cup \geo \Gamma \to \mathcal{D}_{k,k}(Y)$ we have 
\begin{equation}
\lim_{m \rightarrow \infty} \dfrac{1}{\cP(s_m)}\sum_{\gamma \in \Gamma} e^{-s_m \beta(\gamma)} \theta(\gamma) = 
\int_{\geo \Ga} \theta(p) d\mu(p),  
\end{equation}
where convergence is understood as a weak convergence  of currents in $Y$. 
\end{prop}
\begin{proof}
Let $\mu_s=\mu_{1_\Ga,s}, s>s_0$, be the family of probability measures on $\Ga$ as in Section \ref{section_Gibbs} such that the sequence $(\mu_{s_m})$ converges to $\mu$. (Here we use the distance function $d_\beta$ on $\Ga$ defined in Lemma \ref{lem:dbeta} and 
$(\Ga,d_\beta)$ as the hyperbolic metric space in the construction of a Gibbs measure.) 
For every $s>s_0$ consider the current $T(s)\in \mathcal{D}_{k,k}(Y)$ on $Y$:
\begin{equation}
T(s) = \int_{\Gamma} \theta (\ga) d\mu_s(\ga)= \dfrac{1}{\cP(s)}\sum_{\gamma \in \Gamma} e^{-s \beta(\gamma)} \theta(\gamma), 
\end{equation}
and the $(k,k)$-current $T$ on $Y$ given by 
\begin{equation}
T = \int_{\geo \Gamma} \theta(p) d\mu(p). 
\end{equation}

\medskip
For any smooth $(k,k)$-form $\varphi$ on $Y$ and for any $s > s_0$, by integrating on $\ol\Gamma$, we obtain:  
$$
\left | \langle T(s) - T , \varphi \rangle \right |  
= \left | \int_{\ol\Gamma}\langle\theta(x),\varphi \rangle d\mu_{s}(x) -  \int_{\ol\Gamma}\langle\theta(x),\varphi \rangle d\mu(x)  \right |. 
$$
By continuity of the map $\theta$, the function $\psi(x)= \langle\theta(x),\varphi\rangle$ on $\ol{\Ga}$ is continuous. Therefore, by the weak convergence of measures $\mu_{s_m}\to \mu$, we obtain
$$
\lim_{m\to\infty} \left  | \langle T(s_m) - T , \varphi \rangle \right | = \lim_{m\to\infty} \left| 
\int_{\ol\Gamma} \psi(x) d\mu_{s_m}(x) - \int_{\ol\Gamma} \psi(x) d\mu_{s_m}(x) \right |=0. 
$$
This shows weak convergence of currents 
 \begin{equation*}
\lim_{m\rightarrow \infty} T(s_{m}) = T.  
\end{equation*}  
 \end{proof}

\section{Illustrative example: diagonal action on $\bP^1\times \bP^1$}


In this section, to illustrate the ideas of our work without using too much preliminaries from previous chapters, we describe the situation of the flag manifold $\bP^1\times \bP^1$ associated with the Lie group $G=\PGL(2,\bC)\times \PGL(2,\bC)$. We  outline the proof of our main equidistribution theorem when the Anosov representation is the diagonal representation of a convex-cocompact Kleinian group in $\PGL(2,\bC)\times \PGL(2,\bC)$. The strategy of the proof for a general complex Anosov representation will be essentially the same while many more technical details will be needed in some steps. 

\subsection{Action}

Let $\Gamma_1< \PGL(2,\bC)$ be a convex-cocompact Kleinian group and consider its action on $\bP^1\times \bP^1$ defined by 
\[\gamma\cdot(x,y)=(\gamma\cdot x, \gamma\cdot y),\quad \gamma\in \Gamma_1.\]
We will use the letter $\Ga$ to denote the image of $\Gamma_1$ in $G$ given by this action. 
Note that projections of $\bP^1\times \bP^1$ onto both factors are $\Gamma$-equivariant. When we use metrics on $\bP^1$ or on $\bP^1\times \bP^1$, we will be referring to the Fubini--Study metric or the product of two such metrics. 

\subsubsection{Dynamics of a single loxodromic element}

Let $\gamma\in PSL(2,\C)$ be a loxodromic element, i.e. (up to conjugation) an element acting on the extended complex plane as 
$z\mapsto \rho z$, 
$\rho=\rho(\ga)\in \C^*$, $|\rho|>1$. We have 
\[\rho(\ga^n)=\big(\rho(\ga)\big)^n,\quad n\in \bN.\]
When convenient, we will identify $\P^1$ and the extended complex plane. 

The element $\gamma$ has the attracting fixed point $x^+(\ga)=\infty$ and the repelling fixed point $x^-(\ga)=0$. The action of $\ga$ on $\bP^1\times \bP^1$ is given by 
\[(x,y)\mapsto (\rho x,\rho y).\]
Thus, it has  four fixed points in $\bP^1\times \bP^1$:
\[(x^+,x^+), (x^-,x^-), (x^+,x^-), (x^-,x^+).\]
The first two are respectively attracting and repelling, while the last two are hyperbolic fixed points where the differential of $\ga$ corresponds to the matrix 
\[\begin{pmatrix}
\rho&0\\0&\frac{1}{\rho}
\end{pmatrix}.\]

There are also four $\gamma$-invariant rational curves in $\bP^1\times \bP^1$:
\[ \bP^1\times \{x^+\}, \bP^1\times \{x^-\},\{x^+\}\times \bP^1, \{x^-\}\times \bP^1.\]
Each invariant curve contains two fixed points and the action of $\gamma$ on each invariant curve is conjugate to $z\mapsto \rho z$. 


The action of a generic element of a general Anosov subgroup on the corresponding flag manifold is described in Section \ref{section_prelim} (see in particular Section \ref{section_torus_action}). It has finitely many fixed points and  the dynamics is of Morse--Smale type. There are subvarieties passing through some of the fixed points, called Schubert cycles, that play the role of the above invariant curves (see Section \ref{sec:Thickenings in the flag manifold}). What plays the role of $\rho$ are the exponentials of the Lyapunov spectrum (see Section \ref{sec:Cartan and Lyapunov}); here we ignore the argument of $\rho$ as it is ultimately irrelevant for the dynamical purposes. Thus, the more useful parameter 
is $\la=\log |\rho|$,  the only component of the Lyapunov spectrum of $\ga$. 


\subsubsection{Group dynamics}

Let $\Lambda_1\subset \bP^1$ be the limit set of the Kleinian group $\Gamma_1$, i.e.\ an invariant closed subset such that every $\Gamma_1$-orbit in $\Lambda_1$ is dense, while the action of $\Gamma_1$ is properly discontinuous on $\Omega:=\bP^1\setminus \Lambda_1$. By convex-cocompactness, 
the group $\Gamma$ is Gromov-hyperbolic (see Section \ref{sec:Hyperbolic groups}) 
and there exists a $\Gamma$-equivariant homeomorphism 
\[
h:\partial_\infty \Gamma\rightarrow \Lambda_1
\]
where $\geo \Gamma$ is the Gromov boundary of $\Gamma$. 
The quotient $\Omega/\Gamma_1$, if nonempty, is a finite union of compact Riemann surfaces.

The set $\Lambda=\{(x,x)\in \bP^1\times \bP^1\mid x\in \Lambda_1\}$ is the flag-limit set for the action of 
$\Gamma$ on $\bP^1\times\bP^1$. 
It is homeomorphic to $\Lambda_1$, thus to $\geo \Gamma$. We let $f: \geo \Ga \to \La$, be the homeomorphism $f(z)=(h(z), h(z))$. 
The action of $\Gamma$ is properly discontinuous on $\Omega\times \bP^1$ and 
on $\bP^1\times \Omega$, but not on the complement $\bP^1\times \bP^1\setminus \Lambda$. 

The quotient $\big(\Omega\times \bP^1\big)/\Gamma$ is a holomorphic $\bP^1$-bundle over $\Omega/\Gamma$. Here to get the discontinuity domains $\Omega\times \bP^1$ or $\bP^1\times \Omega$, it is crucial to remove a thickened limit set $\Lambda\times \bP^1$ or $\bP^1\times \Lambda$. See Section \ref{sec:Thickenings} for the construction of thickened limit sets in the general Anosov case and see \cite{KLP18A} for more on discontinuity domains of Anosov subgroups.

\medskip
A consequence of the convex-cocompactness of the Kleinian group $\Gamma_1$, which is crucial for us, is the following:

\begin{lem}\label{ill_statement_one}
There exist $\varepsilon,A,C\dyl$ and a finite subset $E\subset \Gamma$ such that for any $\gamma\in \Gamma$, there exists $g\in E$ such that 
$\tga=\gamma\circ g$ is a loxodromic element satisfying the following:
\begin{enumerate}
\item the distance between the two fixed points $x^+(\tga)$ and $x^-(\tga)$ is larger than $\varepsilon$;
\item $\la(\tilde\ga)= \log \lvert\rho(\tga)\rvert\geqslant A^{-1}\lvert \tga \rvert -C$ where $\lvert \tga\rvert$ is the word length of $\tga$, equivalently, 
$$
\lambda(\tga)\geqslant A^{-1} d(\tga_1 x_0, x_0)  -C,
$$
where $x_0$ is a base-point in $\bH^3$, with trivial $\Ga_1$-stabilizer, $d$ is the hyperbolic distance on $\bH^3$, and 
$\tga_1\in \Ga_1$ is the element corresponding to $\tga\in \Ga$. 
\end{enumerate}
\end{lem}
See Section \ref{sec_redressing_anosov} for a proof in the general case. 

\subsection{Gibbs measures}

We need a construction of measures on $\Lambda_1$, thus also on $\geo \Gamma$ and on $\Lambda$ via the above homeomorphism $\La_1\to\La$, that we call {\em Gibbs measures}, which is slightly different from the classical construction of Patterson--Sullivan measures. These measures are treated in greater generality in Chapter \ref{sec:counting}. Here we collect some descriptions for our illustrative example without proofs. For $\ga\in\Ga_1$ let $\beta(\ga)$ be a length function on $\Ga_1$ (see Chapter \ref{sec:counting} for definition; here the reader can take for example $\beta$ to be the word length of $\ga$ in a fixed Cayley graph or the distance function $\beta(\gamma)=d(x_0, \gamma x_0)$ coming from the action on $\bH^3$). Consider the Poincar\'e series:
\[\cP(s)=\sum_{\ga\in \Gamma_1}e^{-s \beta(\ga)}, \quad s>0.\]
There is a critical parameter $s_0\dyl$ such that $\cP(s)$ is convergent for $s>s_0$ and divergent for $s\leqslant s_0$. For $s>s_0$ and $x\in \bP^1$, we define a probability measure on $\bP^1$:
\[\mu_s=\frac{1}{\cP(s)}\sum_{\ga\in \Gamma_1}e^{-s\beta(\ga)}\delta_{\ga x}\]
where $\delta$ is the Dirac mass. For some sequences $(s_n)_n, s_n> s_0$, such that $\lim s_n=s_0$, the measures $\mu_{s_n}$ converge to a limit probability measure $\mu_{Gibbs}$ supported on $\Lambda_1$ which only depends on $(s_n)$ but not on the base point $x$ used in the construction (see Lemma \ref{lem_Gibbs_independant_basepoint} and Theorem \ref{thm:rank1-equi}). Note that, in the setting of length functions $\beta$ coming from the distance function on $\bH^3$, the limit does not depend on the sequence $(s_n)$ either  (see Section \ref{sec:uniqueness}). (For general length functions $\beta$, the limit measures obtained from different sequences are absolutely continuous with respect to each other. In the sequel we will talk about $\mu$ by fixing one sequence $(s_n)$.)

{ For the ease of notation, we will identify $\geo \Ga$ with $\La_1$, so that the map $f$ becomes the diagonal embedding $z\mapsto (z,z)$. 
We will use the notation $\mu_{\cF}$ for the $f$-pushforward of the Gibbs measure $f_*(\mu)$ on $\La$.} 

\subsection{Bad locus and current-valued map}

\subsubsection{Current-valued map}\label{sec:illustration:Current-valued map}

From now on, we consider an irreducible projective curve $S\subset \bP^1\times \bP^1$, which is different from the diagonal $\{(x,x)\mid x\in \bP^1\}$. As 
$\Lambda$ is contained in the diagonal, the set $S\cap \Lambda$ is contained in the intersection locus of two curves, thus, is a finite set. 

For any point $x\in\bP^1$, the homology classes of $\bP^1\times \{x\}$ and $\{x\}\times \bP^1$ generate $\coh_2(\bP^1\times \bP^1,\bZ)$ as a free abelian group of rank two. We, thus, can express the homology class of $S$ as:
\[
[S]=a_1[\bP^1\times \{x\}]+a_2[\{x\}\times \bP^1]
\]
where $a_1, a_2$ are nonnegative integers. The number of intersection points between $S$ and $\bP^1\times \{x\}$ is exactly $a_2$ for a generic $x$, 
similarly for $a_1$.
 Given a point $p \in \partial_\infty \Ga$, recall that $x=f(p)$ is a point $(p,p)$ in the diagonal of  $\bP^1 \times \bP^1$  
and we denote by $\Th_v(x)$ the vertical line $\{p\} \times  \bP^1$ and by $\Th_h(x)$ the horizontal line $\bP^1 \times \{p\}$.

Our main goal is to prove the following:
\begin{thm}\label{ill_statement_two}
We have the weak convergence of currents:
\begin{equation}\label{eq:ill_eq_one}
\lim_{n\rightarrow +\infty}\frac{1}{\cP(s_n)}\sum_{\ga\in\Gamma}e^{-s_n\beta(\ga)}\llbracket \ga(S) \rrbracket= \int_{\Lambda} a_1 \intcur{\Th_h(x)}+ a_2 \intcur{\Th_v(x)} \big)d\mu_{\cF}(x).
\end{equation}
\end{thm}

To prove Theorem \ref{ill_statement_two} we fix a smooth $2$-form $\varphi$ and define a map $\eta$ from the compact space 
$\overline{\Gamma}=\Gamma\cup \geo \Gamma= \Gamma _1 x_0 \cup \La_1$ to $\bC$ by
\begin{equation*}
\phantom{\eta=}  
    \begin{cases}
      \eta(\gamma)=\int_{\ga(S)}\varphi & \text{if $\ga\in \Gamma$},\\
      \eta(p)=a_1\int_{\Th_h(f(p))}\varphi +a_2\int_{\Th_v(f(p))}\varphi & \text{if $p\in \geo\Gamma$}.
    \end{cases}       
\end{equation*}

\subsubsection{Discontinuity and non-genericity}\label{sec:flawedproof}

One of the main difficulties that we will have to deal with is lack of continuity of the map $\eta$. Note that, after pairing with $\varphi$, the left-hand side (without the limit) of \eqref{eq:ill_eq_one} is $\int_{\overline{\Gamma}}\eta\, d\mu_{s_n}$ while the right-hand side is $\int_{\overline{\Gamma}}\eta \, d\mu$. 
If $\eta$ was a continuous function then Theorem \ref{ill_statement_two} would follow from the weak convergence $\mu_{s_n}\rightarrow \mu$.  

To see that $\eta$ is not necessarily continuous, consider a loxodromic element  $\gamma\in \Gamma$ which, up to coordinate change, can be assumed to act on $\bP^1\times \bP^1$ as $(x,y)\mapsto (\rho x, \rho y)$ with $|\rho|>1$. Assume that $S$ is the embedding of $\bP^1$ into $\bP^1\times \bP^1$ defined by $z\mapsto(z,z^2)$. We have 
\[[S]=[\bP^1\times \{x\}]+2[\{x\}\times \bP^1].\]
For $n\in \bN^*$ the curve $\ga^n(S)$ is 
\[
\{(x,y)\in \bP^1\times \bP^1\mid y=\rho^{-n}x^2\}.
\]
Consider a $2$-form $\varphi$ supported on a small neighborhood of $(0,\infty)\in \bP^1\times \bP^1$ such that 
\begin{equation}\label{eq:ill_eq_two}
\int_{\bP^1\times \{\infty\}}\varphi\neq 0=\int_{\{\infty\}\times \bP^1}\varphi.
\end{equation}
Then for large $n$, we have
\begin{equation}
\label{eq:ill_eq_three}
\int_{\ga^n(S)}\varphi=0.
\end{equation}
As $f^{-1}(\infty)$ is the limit point of $(\ga^n)_n$ in $\overline{\Gamma}$, we deduce from \eqref{eq:ill_eq_two} and \eqref{eq:ill_eq_three} that $\eta$ is not continuous at $f^{-1}(\infty)$.

\begin{rem}
The discontinuity in this example is caused by the fact that $S$ passes through the repelling fixed point of $\ga$. 
This repelling point belongs to $\Lambda = \Th^{2-codim(S)-1}(\Lambda)$.
 As we will see later this is, in some sense, the only bad situation.
\end{rem}

\subsubsection{Approximation}
As we have seen in Section \ref{sec:flawedproof}, the intersection between $S$ and $\Lambda$ gives rise to discontinuity of $\eta$. To remedy this problem we define (here $E$ is from Theorem \ref{ill_statement_one})
\[
P=\{x\in \geo \Gamma\mid \exists g\in E, {f(x)}\in g^{-1}(S)\}.
\]
Since $\Lambda\cap S$ and $E$ are finite sets, the subset $P\subset \geo \Gamma$ is finite. We then define, for $t\dyl$,
\[
\Gamma_t=\{\ga\in\Gamma\mid \forall g\in E, g^{-1}\circ \ga^{-1}(x_0)\notin U(P,t)\}
\]
where $U(P,t)$ is a $t$-neighborhood of $P\subset \geo \Gamma$ in $\overline{\Gamma}$ (See Sections \ref{sec:ideal boundaries} and \ref{section_shadow_counting} for the precise definitions). Intuitively an element $\ga$ of $\Gamma$ belongs to $\Gamma_t$ if the repelling fixed point of each $\ga\circ g$ is $t$-away from $S$. 

Let us digress briefly to explain what $U(P,t)$ looks like. Here is the definition of $U(P,t)$ in the case when $P$ is the singleton 
$\{p\}=\{\infty\}\subset \geo \bH^3$ in the upper half-space model 
$$
\{(x_1,x_2,x_3): x_3>0\},
$$
and the point $x_0$ has coordinates $(0,0,1)$ (the neighborhoods 
$U(p,t)$ are defined relative to this base-point). We have
$$
U(p,t)= \{ y\in \ol{\bH^3}: (p,y)_{x_0}> -\log(t)\}. 
$$
Let $b_p(y)$ denote the value at $y$ of the Busemann function normalized to vanish at $x_0$. Thus, for $y=(y_1,y_2,y_3)$, we have 
$b_p(y)=-\log(y_3)$ and 
$$
(p,y)_{x_0}= \frac{1}{2}(d(x_0,y) - b_p(y)). 
$$
The inequality $(p,y)_{x_0}> -\log(t)$ then reads
$$
\frac{1}{2}(d(x_0,y) + \log(y_3)) > -\log(t). 
$$
See Figure \ref{fig:Fig1} for an approximate sketch of the region $U(p,t)$.

\begin{figure}[htbp]
   \centering
   \includegraphics[scale=0.7]{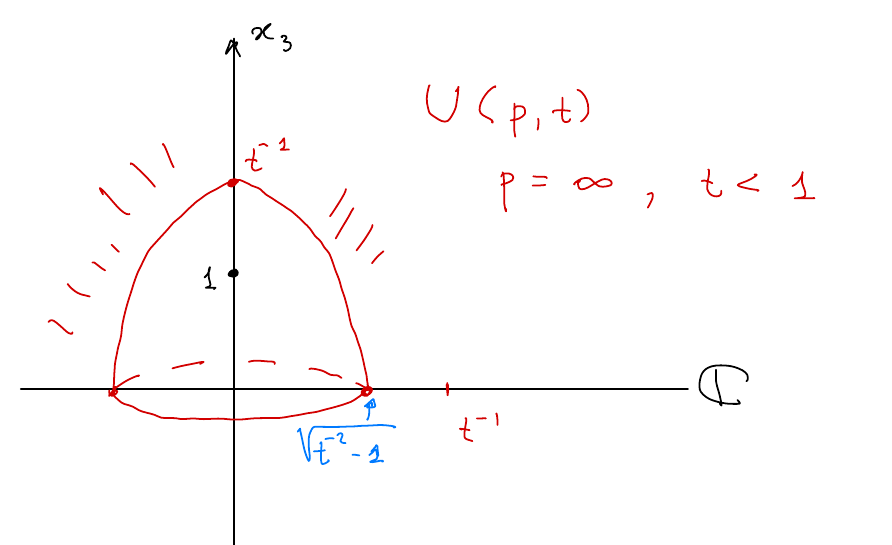} 
   \caption{Neighborhood $U(p,t)$ in the upper half-space model.}
   \label{fig:Fig1}
\end{figure}

One of our main results in Section \ref{section_shadow_counting} implies the following (see Theorem \ref{thm:integration_partial_measures}):

\begin{lem}\label{ill_statement_three}
Theorem \ref{ill_statement_two} holds if $\eta$ restricted to the closure $\overline{\Gamma_t}$ is continuous for all $t\dyl$.
\end{lem}

Informally, one can think of this lemma as follows: The subset $P$ has zero Gibbs measure. This translates into the fact that the total contribution to the sum in left hand side in \eqref{eq:ill_eq_one} of elements $\ga\notin \Ga_t$ 
tends  to zero as $t\to 0$. Thus, for the purpose of proving the theorem one can restrict attention to $\Ga_t$ and, thus, continuity of $\eta$ on $\overline{\Gamma_t}$ suffices.

\subsection{Stretching disks}

By Lemma \ref{ill_statement_three}, it remains to prove that $\eta$ restricted to $\overline{\Gamma_t}$ is continuous. Note that $\eta$ is continuous on $\Gamma$ because $\Gamma$ has discrete topology and that the restriction of $\eta$ to $\geo \Gamma$ is continuous by construction. Fix $t\dyl$ and $p\in \geo \Gamma \cap \overline{\Gamma_t}$. We need to prove that 
\[\int_{\gamma(S)}\varphi\]
is sufficiently close to 
\[
a_1\int_{\Th_v(f(p))}\varphi +a_2\int_{\Th_h(f(p))}\varphi
\]
when $\gamma\in \Gamma_t$ is sufficiently close to $p$. To prove this we note that 
\[\ga(S)=\tga(g^{-1}(S))\]
where $\tga$ is provided by Lemma \ref{ill_statement_one}, and consider the dynamics of $\tga$ instead of that of $\ga$, after replacing the initial curve $S$ with $g^{-1}(S)$. In other words, the point of Lemma \ref{ill_statement_one} is to allow us to deal with only loxodromic elements with the uniform control 
\[\lambda(\tga)\geqslant A\beta(\tga)-C.\]
Note also that by the definition of $\Gamma_t$, the repelling fixed point of $\tga$ is $\epsilon$-bounded away from $g^{-1}(S)$ (with respect to the metric on $\bP^1\times \bP^1$) for some constant $\epsilon\dyl$ depending only on $t$ but not on $\ga$ (see Lemma \ref{uniformboundedawayfromthelimitset1}). Then some analysis on Gromov hyperbolic spaces (see Theorem \ref{thm:redressing2}, Lemma \ref{Mishalimitpointslemma4} and the proof of Proposition \ref{continuousonGammaSt}) reduce  the proof of continuity to the following statement:
\begin{lem}\label{ill_statement_four_bis}
For every  $m\dyl$, there exists $N\dyl$ such that if the element $\tga$ as above satisfies $\lambda(\tga) \geqslant N$, then 
\[
\big\lvert\int_{\ga(S)}\varphi-a_1\int_{\bP^1\times \{x^+(\tga)\}}\varphi-a_2\int_{\{x^+(\tga)\}\times \bP^1}\varphi\big\rvert\leqslant m.
\]
\end{lem}
By using some compactness arguments (see Sections \ref{sec:thickenings_parametrizations} and \ref{sec:cylinders_and_continuity}), one reduces Lemma \ref{ill_statement_four_bis} to the following statement about a single map:
\begin{lem}\label{ill_statement_four}
Assume that $(0,0)$ is $\epsilon$-away from $S$. For every $m\dyl$ there exists $N\dyl$ such that if $\rho\in \bC^*$ satisfies $\lvert\rho\rvert \geqslant N$ then the map $\tau:(x,y)\mapsto (\rho x,\rho y)$ satisfies the inequality
\[
\big\lvert\int_{\tau(S)}\varphi-a_1\int_{\{\infty\}\times \bP^1}\varphi-a_2\int_{\bP^1\times \{\infty\}}\varphi\big\rvert\leqslant m.
\]
\end{lem}

\begin{proof}[Sketch of proof]
This lemma is similar to the classical Inclination Lemma (see for instance, \cite[Proposition 6.2.23]{KHbook}). 
We outline the main ideas behind the proof. We refer to 
Sections \ref{sec:thickenings_parametrizations} and \ref{sec:cylinders_and_continuity} for the general treatment the the case of general flag manifolds. 

For simplicity, we assume that the intersection between $S$ and $\{0\}\times \bP^1$ is transversal and consists of $a_2$ points $v_1,\cdots,v_{a_2}$ while $S\cap (\bP^1\times \{0\})$ is transversal and consists of $a_1$ points $u_1,\cdots,u_{a_1}$. By the assumption, all the points $u_i, v_j$ are bounded away from $(0,0)$. Therefore, there exist two constant $r_1,r_2\dyl$ and disks $B(0,r_1), B(0,r_2)\subset \bC\subset \bP^1$ such that: \\

\begin{enumerate}
\item The intersection of $S$ with $\big(\bP^1\setminus B(0,r_2)\big)\times B(0,r_1)$ consists of $a_2$ disks $V_1,\cdots,V_{a_2}$ centered at $v_1,\cdots,v_{a_2}$ which all project biholomorphically under the first projection onto $B(0,r_1)$.\\

\item The intersection of $S$ with $B(0,r_1)\times \big(\bP^1\setminus B(0,r_2)\big)$ consists of $a_1$ disks $U_1,\cdots,U_{a_1}$ centered at $u_1,\cdots,u_{a_1}$ which all project biholomorphically under the second projection onto $B(0,r_1)$.\\

\item The disks $U_i,V_j$ are pairwise disjoint.
\end{enumerate}

\medskip 
Then each $\tau(V_j)$ is a disk in 
\[B(0,\lvert \rho\rvert r_1)\times \big(\bP^1\setminus B(0,\lvert \rho \rvert r_2)\big)\]
which projects biholomorphically onto $B(0,\lvert \rho\rvert r_1)$, and similarly for the $\tau(U_i)$s, as shown in Figure \ref{fig:Fig2}.

\begin{center}
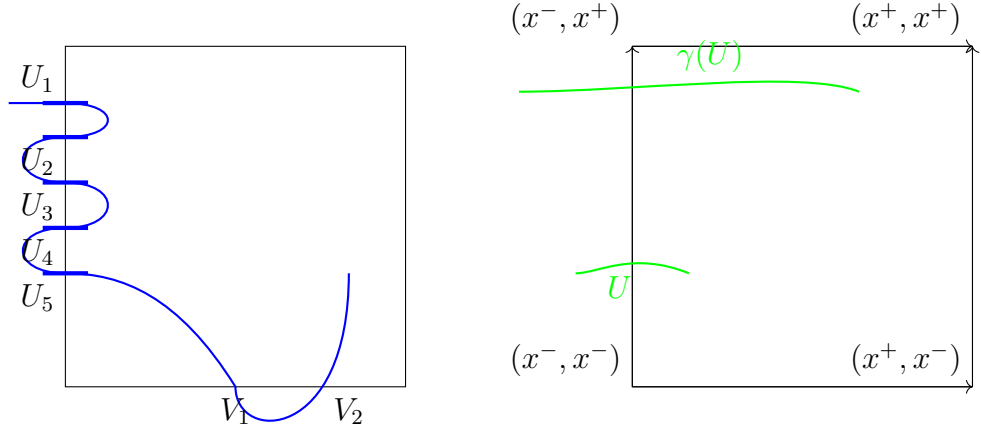
\begin{figure}
\begin{tikzpicture}[scale=1.5]

\draw (0, 0) rectangle (3, 3); 

\draw (5, 0) rectangle (8, 3); 

\draw[->] (5, 0) -- (8, 0); 
\draw[->] (8, 0) -- (8, 3); 
\draw[->] (5, 3) -- (8, 3); 
\draw[->] (5, 0) -- (5, 3); 

\node[above left] at (5, 0) {$(x^-,x^-)$};
\node[above left] at (5, 3) {$(x^-,x^+)$};
\node[above left] at (8, 3) {$(x^+,x^+)$};
\node[above left] at (8, 0) {$(x^+,x^-)$};

\draw[thick, blue, smooth] 
    (-0.5, 2.5) -- (0, 2.5) 
    .. controls (0.5, 2.5) and (0.5, 2.2) .. (0, 2.2) 
    .. controls (-0.5, 2.2) and (-0.5, 1.8) .. (0, 1.8) 
    .. controls (0.5, 1.8) and (0.5, 1.4) .. (0, 1.4) 
    .. controls (-0.5, 1.4) and (-0.5, 1.0) .. (0, 1.0) 
    .. controls (0.5, 1.0) and (1.0, 0.8) .. (1.5, 0) 
    .. controls (1.5, -0.5) and (2.5, -0.5) .. (2.5, 1); 

\draw[ultra thick, blue] 
    (-0.2, 2.5) -- (0.2, 2.5); 
\draw[ultra thick, blue] 
    (-0.2, 2.2) -- (0.2, 2.2); 
\draw[ultra thick, blue] 
    (-0.2, 1.8) -- (0.2, 1.8); 
\draw[ultra thick, blue] 
    (-0.2, 1.4) -- (0.2, 1.4);
\draw[ultra thick, blue] 
    (-0.2, 1) -- (0.2, 1);

\node[above left] at (0, 2.5) {$U_1$}; 
\node[below left] at (0, 2.2) {$U_2$}; 
\node[below left] at (0, 1.8) {$U_3$}; 
\node[below left] at (0, 1.4) {$U_4$};
\node[below left] at (0, 1) {$U_5$};
\node[below] at (1.5, 0) {$V_1$}; 
\node[below] at (2.5, 0) {$V_2$}; 

\draw[thick, green, smooth] 
    (4.5, 1.0) .. controls (4.7, 1.0) and (5.0, 1.2) .. (5.5, 1.0) 
    node[midway, below] {$U$};
\draw[thick, green, smooth] 
    (4, 2.6) .. controls (5, 2.6) and (6.5, 2.8) .. (7, 2.6) 
    node[midway, above] {$\gamma(U)$};

\node at (5, 1.0) {}; 
\node at (5, 2.0) {}; 

\end{tikzpicture}
\caption{stretching disks}\label{fig:Fig2}
\end{figure}
\end{center}


In this case, the main observation follows an idea going back to Ahlfors (see e.g \cite[\S~2]{brunella}), even though the currents of integration along $\tau(U_i)$ are not closed currents, their area increases so that as $\rho \rightarrow +\infty$, the limiting current becomes closed (see for instance, Example \ref{example:cylinder}).  
In other words,  
 when $\lvert \rho \rvert\rightarrow +\infty$, we have for each $i$,
\[
\int_{\tau(U_i)}\varphi\rightarrow 
\int_{\bP^1\times \{\infty\}}\varphi.
\]
Similarly for $\tau(V_j)$'s. Therefore if we denote by $S'$ the subset of $S$ equal to the union of all the disks $U_i,V_j$, then $S'$ satisfies the conclusion of Lemma \ref{ill_statement_four} for $\lvert \rho \rvert$ sufficiently large. To conclude it suffices to note that the mass of the complement $\tau(S\setminus S')$ converges uniformly to zero since the mass of $\tau(S')$ converges to the total mass of $S$.  Thus, 
the Wirtinger inequality (see Corollary \ref{cor_difference_current}) implies that 
the current of integration on $\tau(S\setminus S')$ does not contribute to the limit integral. 
\end{proof}

\section{Setup}\label{sec:flag_setup}

In the rest of this chapter we let $G$ be a complex semisimple Lie group, $B< G$ a Borel subgroup. Let $n$ be the 
complex dimension of the flag manifold $\cF=G/B$. As before, we will fix an arbitrary K\"ahler metric $\mathcal{g}$ on $\cF$ with the distance function $\dF$. 
We will fix, for every root subspace $\mathfrak{g}_\alpha$ of the Lie algebra $\mathfrak{g}$, a real generator vector $\vec{b}_\alpha\in \mathfrak{g}_{\alpha,\bR}$. Then 
the subspace  $\fu^+\oplus \fu^-\subset \mathfrak{g}$ 
is equipped with a Euclidean norm $\vert\vert \cdot \vert\vert$ where the 
generators $\vec{b}_\alpha$ form an orthonormal basis. 

We will also fix an Anosov subgroup $\Gamma< G$ acting on $\cF$ with flag-limit set $\Lambda\subset \cF$ and boundary homeomorphism 
$f: \geo \Gamma \rightarrow \Lambda$. Let $\beta$ be a hyperbolic length function on $\Gamma$. We fix a PS-sequence $(s_j)_j$ 
and the corresponding limiting Gibbs-measure $\mu_{\Ga}$ on $\geo \Gamma$. 
For the ease of notation, we let $\mu$ denote the Gibbs measure $\mu_{\cF}=f_*\mu_{\Ga}$ on $\La$.

Recall that a finite subset $E\subset \Ga$ and a constant $\varepsilon> 0$ associated with $\Gamma$ were introduced in Propositions \ref{Mishaseigenvaluelemma} and \ref{Mishaseigenvaluelemma1}. In this section 
we fix $E$ and $\varepsilon$. If $\gamma\in \Gamma$, then we denote by $\tilde{\gamma}=\gamma\circ g$ "the" redressing of $\gamma$: here $g\in E$ is not necessarily unique and depends on $\gamma$, but we will omit the suffixes for the ease of notation.

\section{The genericity condition}\label{sec:genericity}

Let $S$ be an irreducible projective subvariety of (complex) dimension $k$ in $\cF$. As in the introduction, we define \index{$\La_S$}
$$
\La_S=\{\la\in \La: S\cap \Th^{n-k-1}(\la)\ne \emptyset \}.
$$
Since $\Th^{n-k-1}(\la)$ depends continuously on $\la$, the subset $\La_S\subset \La$ is closed. We say that $S$ is $\Lambda$-generic (or $\Gamma$-generic) if the following condition is satisfied:
\begin{equation}
\label{eq:the_genericity_condition}
\mu(\La_S)=0.
\end{equation}
From now on, we assume that $S$ is $\Lambda$-generic.

	
For example if $S$ is a hypersurface then Condition \eqref{eq:the_genericity_condition} means that $\mu(S\cap \Lambda)=0$, and if $S$ is a singleton $p$, then it means 
for $\mu$-a.e. $\la\in \La$,  $p\notin \Th^{n-1}(\la)$, i.e. $p$  is opposite to $\mu$-a.e.  limit point.



We will need the following recent theorem proven by Kim and Oh in  \cite{Kim-Oh}. They assumed that Anosov subgroups are Zariski dense, but their proof actually shows the following: 

\begin{thm}\label{thm:kim-oh}
Let $\Ga< G$ be a nonelementary Anosov subgroup. Then for every Gibbs measure $\mu$ on $\La$, and every projective subvariety  $V\subset \cF$ 
either $\La\subset V$ or $\mu(V\cap \La)=0$.  
\end{thm}
\begin{proof} We include a proof for the sake of completeness, following the arguments in  \cite{Kim-Oh}. Suppose that  $\mu(V\cap \La)>0$ but 
$Z_1:=V\cap \La$ is a proper subset. Without loss of generality, we may assume that $V$ is the smallest subvariety with this property. 
Then for every $\ga\in \Ga$ either $\mu(\ga Z_1\cap Z_1)=0$ or $\ga Z_1=Z_1$ (otherwise, we replace $V$ with $V\cap \ga V$). Since $Z_1$ is closed and the $\Ga$-action on $\La$ is minimal and satisfies the convergence property (see \S \ref{sec:convdy}), there exists $\ga_0\in \Ga$ such that $\ga_0(Z_1)\cap Z_1=\emptyset$ (for instance, one can take a high power of loxodromic element of $\Ga$ whose fixed points lie outside of $Z$). 
The subset $Z_2:=\ga_0(Z_1)\subset \La$ also has positive measure. Furthermore, $Z_2$ also has the property that for every $\ga\in \Ga$  either $\mu(\ga Z_2\cap Z_2)=0$ or $\ga Z_2=Z_2$. Set $Z:= Z_1\times Z_2\subset \La^2$. Recall that  the diagonal action of $\Ga$ on $(\La\times \La, \mu\times \mu)$ is ergodic (Theorem \ref{thm:Gibbs list}). Hence, according to Lemma \ref{lem:infinitely conservative},  there is an infinite sequence $\gamma_i\in \Ga$ consisting of 
distinct elements, such that for every $i$
$$
 \mu^2(Z\cap \ga_i Z)>0,
$$
i.e.
$$
\mu(Z_1\cap \ga_i Z_1)>0, \mu(Z_2\cap \ga_i Z_2)>0.  
$$
Therefore, as we observed earlier, $\ga_i(Z_1)=Z_1, \ga_i(Z_2)=Z_2$ for every $i$. Applying the convergence property to the sequence $(\ga_i)$, after extraction, we can assume that there exist $\la_\pm\in \La$ such that $\ga_i\to \la_+\in \La$ uniformly on compacts in $\La\setminus \{\la_-\}$. Assume that $\la_+\notin Z_1$.  
Since $Z_1$ is not a singleton (as $\mu$ has no atoms), there exists a point $z\in Z_1\setminus \{\la_-\}$ such that $\lim_{i\to\infty} \ga_i(z)=\la_+$. But then for large $i$, $\ga_i(z)\notin Z_1$,  which is a contradiction. 
Assume finally that $\la_+ 
\in Z_1$, then $\la_+ 
\notin Z_2$ as $Z_1$ and $Z_2$ are disjoint. We thus apply the previous argument, 
switching roles of $Z_1$ and $Z_2$. 
\end{proof}

\begin{prop}\label{prop:genericity}
Suppose that $\Ga< G$ is an Anosov subgroup with the flag-limit set $\La\subset \cF$. Then for every hyperbolic length function $\beta$ on $\Ga$ and every 
algebraic subvariety $S\subset \cF$ either $\La_S=\La$ or $\mu(\La_S)=0$, where $\mu=\mu_\beta$ is a Gibbs measure on  $\cF$ associated with $\beta$ ($\mu$ is, of course, supported on $\La$). 
\end{prop}
\begin{proof}  Observe that for two points $x, y\in \cF$ we have
$$
y\in \Th^m(x) \iff x\in \Th^m(y)\iff D(x,y)\le m,
$$
where $D(x,y)=|\pos(x,y)|$ is the distance function on $\cF$ defined in \S \ref{section_relative_position_bruhat}. Therefore, setting $m=n-1-\dim(S)$, we obtain:
$$
\la\in \La_S \iff \la\in \Th^m(S)\cap \La 
$$
and, thus, $\La_S=\Th^m(S)\cap \La$. Since $S$ is a projective subvariety in $\cF$, so is $V:=\Th^m(S)$ (see Proposition \ref{prop:projective}). We now apply Theorem \ref{thm:kim-oh}. If $\mu(V\cap \La)=0$, then we are done. If $V\cap \La=\La$ then $\La_S=\La$ and we are done as well. 
\end{proof}

\begin{cor}\label{cor:all-subvarieties-generic}
If $\Ga< G$ is Zariski dense then every projective subvariety $S\subsetneq \cF$ is $\La(\Ga)$-generic. 
\end{cor}
\begin{proof} Suppose that $\Th^m(S)$ contains $\La$. Then the intersection
$$
\bigcap_{\ga\in \Ga} \ga \Th^m(S)
$$
is a subvariety in $\cF$ which contains $\La$ (hence, is nonempty) and is $\Ga$-invariant. But this means that $\Ga< G$ is not Zariski dense. 
\end{proof}


For the $\Lambda$-generic subvariety $S$, we introduce the following compact subset of $\Lambda$:
$$
\Lambda_{S,E}:=\bigcup_{g\in E}\{\lambda\in \Lambda: \Th^{n-k-1}(\lambda)\cap g^{-1}(S)\ne \emptyset \}= \bigcup_{g\in E} g^{-1}\La_S
$$
and its counterpart in the visual boundary of the group $\Ga$:
\begin{equation}\label{eq:P1} 
P:=f^{-1}(\Lambda_{S,E})\subset \geo \Ga. 
\end{equation}
Note that $\Lambda_{S,E}$ and, hence, $P$, are compact sets of zero Gibbs-measure because $S$ satisfies Condition \eqref{eq:the_genericity_condition}.


Recall from \eqref{eq_Gpt} that for all $t > 0$, we have defined a subset 
$\Gamma_{P,E,t}\subset \Ga$ associated with a finite subset $E\subset \Ga$ and a closed 
subset $P\subset \geo \Ga$. Note that in the setting that we are considering here, the finite 
set $E$ depends only on $\Gamma$ (regarded as a subgroup of $G$) 
while $P$ also depends on $S$. 
For  $t>0$ we also define 
\begin{equation} \label{eq_defi_gammaSt_flags1}
\Gamma_{S,t} = \Gamma_{P,E,t}\subset \Ga. 
\end{equation}

Let $\beta$ be a hyperbolic length function on $\Ga$, $d_\beta$ the corresponding $\Ga$-invariant metric on $\Ga$. We will be 
using Theorem \ref{thm:redressing3}, where $X=(\Ga, d_\beta)$ and $y_0\in X$ corresponds to the neutral element in $\Ga$.

\begin{lem}\label{uniformboundedawayfromthelimitset1}
Let $S\subset \cF$ be a $\Lambda(\Gamma)$-generic subvariety of dimension $k$. For any $t >0$, there exist $\epsilon_t>0$ and $N_t>0$ such that for any $\gamma \in \Gamma_{S,t}$ with $\beta(\gamma)>N_t$, the repelling fixed point $x^-\in \cF$ of $\tilde{\gamma}$ is in the set 
\[
\Lambda_{\epsilon_t}=\{\lambda\in \Lambda\vert \forall g\in E, d_\cF\left(\Th^{n-k-1}(\lambda),g^{-1}(S)\right)\geq \epsilon_t \}.
\]
\end{lem}

\begin{proof} 
Let $\xi^-\in \geo \Ga$ be the repelling fixed point of $\tilde\ga$, i.e. the attractive fixed point of its inverse $h\in \Ga$.  
We fix $t >0$. We need to bound from below the minimal distance between $\xi^-$ and 
$$
P':=\bigcup_{g\in E} g^{-1} f^{-1}(\Th^{n-k-1}(S)\cap \Lambda). 
$$

According to Theorem \ref{thm:redressing3}, there is 
a constant $\kappa'>0$ such that if $h\in \Ga$, 
is such that the axis of $h$ in $X$ is within distance $R$ from $y_0$, $\beta(h)> -\log(t')$ 
and the attractive fixed point of $h$ 
is in $U(\zeta, t')$ for some $\zeta\in \geo \Gamma$, then $h(y_0)$ 
 belongs to $U(\zeta, \kappa' t')$. We will use $h:= (\tilde\ga)^{-1}$: Its axis in $X$ is the same as the axis of $\tilde\ga$, i.e. within distance 
 $R$ from $y_0$, where $R$ is defined as in \eqref{defi:R}.   The repelling fixed point $\xi^-$ of $\tilde\gamma$ is the attracting fixed point  of 
 $h$. Since, by the assumption, $\gamma \in \Gamma_{S,t}$, we have
 $$
 h(y_0)\notin U(P', t). 
 $$
 Set $t':= (\kappa')^{-1}t$. Then $ h(y_0)\notin U(P',t)= U(P', \kappa' t')$.  Therefore, as we noted above, if $\beta(h)> -\log(t')$, then 
 $\xi^-\notin U(P', t')= U(P', (\kappa')^{-1}t)$. However,
 $$
 \beta(h)=\beta(h^{-1})= \beta(\tilde \ga)\ge \beta(\gamma)-\beta(g) .  
 $$
 Taking $m_0:=\max\{\beta(g): g\in E\}$, we obtain that if 
 $$
 \beta(\gamma)> N_t:= -\log(t)+ m_0,$$  
 then $\beta(h)> -\log(t)$. 
 
 Now, let $\Lambda'$ denote the compact subset 
 $$
 \La \setminus f(U(P', (\kappa')^{-1}t)). 
 $$
 The thickening $\Th^{n-k-1}(\Lambda')$ is a compact subset of $\cF$ disjoint from the compact
 $$
 \bigcup_{g\in E} g^{-1}(S). 
 $$
 Therefore
 $$
 \epsilon_t:= d(\Th^{n-k-1}(\Lambda'), \bigcup_{g\in E} g^{-1}(S))>0. 
 $$
As we observed above, for every $\gamma \in \Gamma_{S,t}$ satisfying $\beta(\gamma)> N_t$, we have $\xi^-\in \La'$, 
where $\xi^-\in \geo\Ga$ is the repelling fixed point of $\tilde\ga$. Hence,  
the repelling fixed point $x^-\in \cF$ of $\tilde{\gamma}$ belongs to $\Lambda_{\epsilon_t}$, as required. 

\end{proof}

\subsection{Thickenings and parametrizations}\label{sec:thickenings_parametrizations}

Recall that $\varepsilon>0$ is the constant from Proposition \ref{Mishaseigenvaluelemma1}. In this section we will be dealing with several objects depending continuously on $\varepsilon$-separated pairs of opposite flags (chambers) $\bx=(x^+,x^-)\in \cF^2$ (see \eqref{e_separation_inequality} for the definition of separation used in what follows). Recall that in Section \ref{section_prelim} we already defined the following objects:  
\begin{itemize}
\item the apartment $a_{\bx}=a(x^+, x^-)$; 
\item the chambers $x_w, w\in W$ in $a(x^+, x^-)$;
\item the Schubert cells $S_{w}(x^+), S_{w^\vee}(x^-), \Opp(x_{w^\vee})$;
\item the thickenings $\Th_{w}(x^+), \Th_{w^\vee}(x^-)$.
\end{itemize}
We set
\begin{align}\label{eq:bundle_separated_pairs}
\begin{aligned}
\mathscr{K}_\cF &:= \{\bx=(x^+,x^-)\in \cF^2, (x^+,x^-)\ \hbox{is}\ \varepsilon\text{-separated}\},\\
\mathscr{K}_\Lambda &:= \mathscr{K}_\cF \cap \Lambda^2,\\
\mathscr{U}_w &:=\{(x^+,x^-,y)\in \cF^3, (x^+,x^-) \  \ \hbox{is}\ \varepsilon\text{-separated}, y\in \Opp(x_w)\}.
\end{aligned}
\end{align}
Here $\La\subset \cF$ is a compact antipodal subset. 

Then 
every $\mathscr{U}_w$
 is a fibration over $\mathscr{K}_\cF$ 
with affine 
spaces as  fibers.

\subsubsection{Parametrization with varying axii} \label{sec:parametrization}

Let $\bx=(x^+,x^-)$ be an $\varepsilon$-separated pair. The $\varepsilon$-separation implies that, for any $w\in W$, the closed ball (defined with respect to the metric $\dF$ on $\cF$)  
$$\bar{\mathbf{B}}(x_w,\frac{\varepsilon}{2})\subset \cF$$
is contained in $\Opp(x_{w^\vee})$. 

Let $w\in W$ and $k=\vert w \vert$. Recall from \S~\ref{sec:Structure of thickenings} that the nilpotent Lie algebra $\fu_w$ is biholomorphic to the cell $\Opp(x_{w^\vee})$ via the map $\pi_w\circ \exp : u \mapsto \exp(u) x_{w}$ . We let 
$$\Psi_w:\Opp(x_{w^\vee})\rightarrow \fu_w$$
denote the inverse map $\left(\pi_w\circ \exp\right)^{-1}$.
We will denote by
$$
\Psiwp: \Opp(x_{w^\vee})\rightarrow \fuwp,\ \Psiwm:\Opp(x_{w^\vee})\rightarrow \fuwm
$$
the compositions of $\Psi_w$ with the projections $\fu_w\rightarrow \fuwp$, $\fu_w\rightarrow\fuwm$. In particular, we have
\[
\Psi_w(x)=\Psiwp(x)+\Psiwm(x), x\in \Opp(x_{w^\vee}).
\]
The restriction of $\Psiwp$ (resp.\ $\Psiwm$) to the Schubert cell $\Swx$ (resp.\ $\Swxd$) is exactly the inverse map of the restriction of $\pi_w\circ \exp$ to $\fuwp$ (resp.\ $\fuwm$). 

\begin{rem} \label{rem_algebraic_projections} The maps $\Psi_w^{\pm}$ are algebraic as the logarithm restricted to a maximal 
unipotent subgroup of $G$ is an algebraic function. 
\end{rem}

Note that all maps $\Psi_w^\pm, \Psi_w$ send the point  $x_w$ to $0$. 
We shall denote these maps $\Psi_{w,\bx}^\pm$, $\Psi_{w,\bx}$, respectively, when $\bx=(x^+,x^-)$ is not fixed. In the next sections, we will often use the following observation. 
\begin{rem}\label{rem_continuity_psi_inverse} Observe that if $g \in G$ is such that $g \bx_0 = \bx$ where  $\bx, \bx_0 \in \cF\times \cF$ are two pairs of antipodal elements. One has for all $u \in \mathfrak{u}_w$,  $\Psi_{w,\bx}^{-1} (u) = \exp(Ad(g)(u)) x_{w}$ for all $w \in W$. In particular, 
 the mapping $\bx \mapsto \Psi_{w,\bx}^{-1}$ are continuous. 
\end{rem}

\subsubsection{Comparison of metrics}

We prove in this subsection some estimates of how the K\"ahler metric on $\cF$ compares to the Euclidean norm on $\fu_w$ via $\Psi_{w,\bx}^{-1}$. All the quantities involved here will be functions involving the distance $d_{\cF}$, the biholomorphic mappings $\Psi_{w,\bx}^{-1}$, which a priori depend on the choice of $\bx$. However, all the quantities that we will consider will be continuous in $\bx$ (see Remark \ref{rem_continuity_psi_inverse}). Hence if we force $\bx$ to be in the set $\mathcal{K}_{\cF}$ which is compact  by assertion (ii) Lemma \ref{lem:separate}, the estimate will be uniform.  To illustrate this, let us give a first statement of this type.

Recall that we have fixed generators $\vec{b}_\al$ of root subspaces ${\mathfrak g}_{\al}$.
\begin{lem}\label{lem:genarators_bounded_length}
There exists constants $l_1,l_2\dyl$ such that for any $\xpair\in \mathscr{K}_\cF$, for any $w\in W$, for any root $\al$ such that $\vec{b}_\al\in \fu_w$, we have 
\[
l_1\leqslant d_\cF\left(x_w, \Psi_{w,\bx}^{-1}(\vec{b}_\al)\right)\leqslant l_2.
\]
\end{lem}
\begin{proof}
The lemma follows from the compactness of $\sK_\cF$ by assertion (ii) of Lemma \ref{lem:separate} and the continuity of the maps $\bx \mapsto\Psi_{w,\bx}^{-1}(\vec{b}_\alpha)$ by Remark \ref{rem_continuity_psi_inverse}, of $(x_+,x_-) \in \mathcal{K}_{\cF} \mapsto x_{w} $ and of the distance function. 
\end{proof}

Recall that for $x\in \cF$ we denote $\Th_{<w}(x)=\Th_w(x)\setminus S_w(x)$. 
Recall also that the generators $\vec{b}_\al$ form an orthonormal basis and, hence, determine a Euclidean metric on every $\fu_w$. The maps $\Psi_w^\pm$ satisfy the following uniform properties:

\begin{lem}\label{comparison_metric_balls}
There is a positive function $\tau_1(R)$ (which depends only on $\varepsilon$) with 
$$\lim_{R\rightarrow +\infty}\tau_1(R)=0$$
such that for any $\xpair\in \mathscr{K}_\cF$ and any $w\in W$, 
\begin{enumerate} 
\item the Hausdorff $d_{\cF}$-distances between $\Th^{n-1}(x_{w^\vee})$ 
and $\Psi_{w,\bx}^{-1}( \mathbf{B}(0,R))$, between $\Th_{<w^\vee}(x^-)$ and 
$\Psi_{w,\bx}^{-1}(\mathfrak{u}^-_w \cap \mathbf{B}(0,R))$, between 
$\Th_{<w}(x^+)$)  
and $\Psi_{w,\bx}^{-1}(\mathfrak{u}^+_w \cap \mathbf{B}(0,R))$, 
are bounded above by $\tau_1(R)$; 
\item the $n$-dimensional (resp.\ $n-\vert w\vert$, $\vert w \vert$-dimensional) volume of the complement of $\Psi_{w,\bx}^{-1}(\mathbf{B}(0,R))$ in $\cF$ (resp.\ $\Th_{w^\vee}(x^-)$, $\Th_w(x^+)$) is bounded above by $\tau_1(R)$. 
\end{enumerate}
\end{lem}
\begin{proof}
As in the proof of the previous statement, for each fixed $\bx\in \mathcal{K}_{\cF}$, the subsets $\Psi_{w,\bx}^{-1}(\mathbf{B}(0,R))$ exhaust the open Schubert cell $\Opp(x_{w^\vee})$ when $R$ tends to infinity. We conclude using the continuity of $\bx \mapsto d_{\cF}(\Psi_{w,\bx}^{-1}(\mathbf{B}(0,R)), \Th^{n-1}(x_{w^\vee}))$ and the compactness of $\mathcal{K}_{\cF}$.    
Similarly, $\Psi_{w,\bx}^{-1}(\mathfrak{u}^\pm_w \cap \mathbf{B}(0,R))$ also exhaust $S_w(x^+)$ and  $S_{w^\vee}(x^-)$ respectively. 
For the second statement, observe that the complement  $\cF \setminus \Psi_{w,\bx}^{-1}(\mathbf{B}(0,R))$ are semi-algebraic sets and as $R \rightarrow +\infty$,  the Hausdorff distance between these sets and $\Th^{n-1}(x_{w^\vee})$ converges to zero. Since $\Th^{n-1}(x_{w^\vee})$ is $n-1$ dimensional, its $n$-dimensional volume is zero. We then conclude (using e.g \cite[Theorem 5.10]{yomdin_comte}) that their $n$-volume converge to zero. We conclude similarly when one restricts to the thickenings $\Th_{w^\vee}(x^-)$ and $\Th_w(x^+)$. 
\end{proof}




As an immediate consequence of the first part of the lemma, we obtain: 

\begin{lem}\label{uniform_r1r2}
For any $\epsilon>0$, there exists $R\dyl$ such that for any $\xpair\in \mathscr{K}_\cF$, for any $w\in W$, if the $d_{\cF}$-distance between a point $z\in \Opp(x_{w^\vee})$ (resp.\ $S_w(x^+)$, $S_{w^\vee}(x^-)$) and $\Th^{n-1}(x_{w^\vee})$ (resp.\ $\Th_{<w}(x^+)$, $\Th_{<w^\vee}(x^-)$) is larger than $\eps$, then $\Psi_w(z)\in \mathbf{B}(0,R)$.
\end{lem}

In the next lemma, we shall relate the difference between a neighborhood of large unstable polydisc with the euclidean neighborhood of its image in the chart $\Psi_{w,\bx}^{-1}$. 
More precisely, recall that $\Psi_{w,\bx}^{-1} (\mathbf{B}(0,R) \cap \mathfrak{u}_w^{+})$ is a subset of $\Th_w(x_+)$ containing $x_w$ and can be see as an unstable plaque of radius $R$ around $x_w$. A nearby point can be taken by multiplying by an element of  $G$ close to identity, but we measure its effect in the chart $\Psi_{w,\bx}$.
Recall also from \S~\ref{sec:expansion_unipotent} that tangent space at $x_w$ can be identified with the Lie subalgebra $\mathfrak{u}_w$ and given a $K$ invariant metric on $\cF$, its restriction yields a metric we denote by $|\cdot |$. 
%
\begin{lem}\label{lem:vertical_BCH_estimate}
There exists a constant $C>0$ 
such that for any $\xpair\in \sK_\cF$, any $w\in W$, any $R_h,R_v\dyl$ such that $R_v\geqslant R_h\geqslant 1$, any $y\in \fuwp$ with $\vert y \vert \leqslant 1$ and any $z\in \fuwm$ with $\vert z\vert \leqslant 1$, we have, for $\mathbf{B}=\mathbf{B}\left(0,C\left(\frac{R_h^{N_0}}{R_v}\right)\right)\subset \fu_w$ where $N_0$ is the order of nilpotence of the Lie algebra $\fu_w$, 
\[
\Psi_{w,\bx}^{-1} \left (R_hy+\frac{1}{R_v}z \right )\in \exp\left(\mathbf{B}\right)\left(\Psi_{w,\bx}^{-1}(R_hy)\right).
\]
\end{lem}
\begin{proof}
Fix $(x^+,x^-,w)$. 
Let 
$y\in \fuwp$ with $\vert y \vert \leqslant 1$, and let $z\in \fuwm$ with $\vert z\vert \leqslant 1$.
For any $t,t'\in \C$, we let $v(t,t') \in \mathfrak{u}_w$ determined by the relation:
\begin{equation}\label{eq:BCH_formula0}
\exp(v(t,t'))=\exp(t y+ t' z)\exp (-t y).
\end{equation}

By Baker-Campbell-Hausdorff Formula (see e.g. \cite[Theorem 1.2.1]{corwin}), setting $y'(t,t') = ty + t' z$, one has:
\begin{multline*}
v(t,t')=t' z-\frac{1}{2}[y'(t,t'), t y]-\frac{1}{12}[y'(t,t'),[y'(t,t'), t y]]+\frac{1}{12}[t y,[y'(t,t'),t y]]+\cdots \\ 
= t' z + \sum_{k=1}^{+\infty} P_k(ad(t y)  , ad( y'(t,t'))) t y ,
\end{multline*}
where $P_k$ are particular (non-commutative) homogeneous polynomials of degree $k$ in two variables.
Since $y'(t,t'), z $ belong to $\fu_w$ which is a nilpotent, the right hand side of the previous formula is a finite sum. Namely, let $N_0$ be the degree of nilpotency of the subalgebra $\fu_w$, then for $k \geqslant N_0+1$, the polynomial term $P_k(ad(ty), ad(y'(t,t')))$ vanish, and we get for all $t,t'  \in \C$: 
\begin{equation*} 
v(t,t') = t' z + \sum_{k=2}^{N_0} P_k(ad(t y)  , ad( y'(t,t'))) t y . 
\end{equation*} 
Observe also that $ad(y'(t,t')) = ad(t y) + ad(t' z)$, so expanding the previous expression shows that there exists a non-commutative polynomial $Q_k$ of degree $k$ such that
\begin{equation} \label{eq_BCH_simplified}
v(t,t') = t' z + \sum_{k=2}^{N_0} Q_k(ad(ty), ad(t'z)) ty.
\end{equation} 
In the above formula $v$ is a polynomial in $t,t'$ whose pure term in $t$ all vanish, because they involve brackets in the span of $y$. This shows that the non-vanishing coefficients correspond to monomials $t^p t'^q$ with $q \geqslant 1$.  
Note that the polynomials $Q_k$ in the above formula do not depend on the choice of $y,z$. 
Taking $t = R_h$ and $t' = 1/R_v$ with $R_h, R_v > 1$ and using the triangular inequality, we see that the norm of $v(R_h, 1/R_v)$ is bounded by $C R_h^{N_0}/R_v$ where $C> 0$ is independent of the choice of $y,z$.   
Let $\mathbf{B}$ be the ball $\mathbf{B}(0, C R_h^{N_0}/ R_v)$,  we have obtained for all $y \in \mathfrak{u}_w^+, z\in \mathfrak{u}_{w}^-$ in the unit ball:
\begin{align}
\label{eq:BCH_bigO}
\exp(R_h y+\frac{1}{R_v}z)=\exp\left(\mathbf{B} \right)\exp(R_hy).
\end{align}
Then the definition $\Psi_w^{-1}=\pi_w\circ \exp$ implies the conclusion of the Lemma.
\end{proof}

For $\bx=\xpair\in \mathscr{K}_\cF$, $w\in W$, $z\in \Opp(x_{w^\vee})$, we define $\cah(\bx,w,z)$ as the length (with respect to the background K\"ahler metric on $\cF$) of the preimage under $\Psi_{w,\bx}$ of the line segment 
\[
\{y_1+cy_2 \mid y_1=\Psi_{w,\bx}^+(z), y_2=\Psi_{w,\bx}^-(z), 0\leqslant c \leqslant 1\}\subset \fu_w.
\]
Then $\cah$ is a continuous function with respect to $x^+,x^-,z$. We will sometimes write $\cah(z)$ for ease of notation. As a corollary of the proof of Lemma \ref{lem:vertical_BCH_estimate}, we get
\begin{cor}\label{cor:pizza_BCH}
For any $\epsilon>0$, there exists $\theta > 0$ such that 
for any $\bx\in \sK_\cF$, any $w\in W$, any $R_h,R_v\dyl$ such that $R_v\geqslant R_h\geqslant 1$, any $z\in  \Opp(x_{w^\vee})$ with $\vert \Psi_{w,\bx}^+(z) \vert \leqslant R_h$, $\vert \Psi_{w,\bx}^-(z)\vert \leqslant 1/R_v$ and such that $R_h^{N_0}/R_v\leqslant \theta$ where $N_0$ is the order of nilpotence of the Lie algebra $\fu_w$, we have 
\begin{align*}
\cah(\bx,w,z)\leqslant \epsilon.
\end{align*}
\end{cor}
\begin{proof}
Fix $\bx,w,R_h,R_v$ as in the statement. Let $z\in  \Opp(x_{w^\vee})$ with $\vert \Psi_w^+(z) \vert \leqslant R_h$ and $\vert \Psi_w^-(z)\vert \leqslant 1/R_v$. 
By Lemma \ref{lem:vertical_BCH_estimate}, any point $\Psi_{w,\bx}^{-1} (R_h y + z/R_v)$ belongs in the neighborhood $\exp(\mathbf{B}) \Psi_{w,\bx}^{-1}(R_h y)$ where $\mathbf{B} = \mathbf{B}(0, C R_h^{N_0}/R_v)$ with $C>0$.  
For any $\bx \in \sK_{\cF}$, any $w \in W$, any $R_v \geqslant R_h \geqslant 1$ and any $y \in \mathfrak{u}_{w}^+$ of norm bounded by $1$, the diameter of $\exp(\mathbf{B}) \Psi_{w,\bx}^{-1}(R_h y) \subset \cF$ for the Riemannian metric depends continuously on $(\bx,y)$ which belong to the compact set $ \sK_{\cF} \times  \{ y \in \mathfrak{u}_{w}^+ \ | \ |y|\leqslant 1 \}$. Moreover, as $R_h^{N_0}/R_v$ tends to zero, this diameter converges to zero. We conclude using the fact that $\cah(\bx, w, z)$ is bounded above by this diameter. 
\end{proof}

\section{Continuity of the current-valued map}\label{sec:cylinders_and_continuity}

Let $S$ be an irreducible $k$-dimensional subvariety  in $\cF$ which is $\Lambda$-generic. 
Recall that $\mathcal{D}^{n-k,n-k}_+(\cF)$ denotes the space of positive currents of bidegree $(n-k,n-k)$ on $\cF$. Let 
\[ \intcur{S}=\sum_{w\in W}a_w \intcur{\Th_w(x)}\]
be the homological decomposition of $S$ (see equation \eqref{eq:decomposition} and Section \ref{sec:homology}). 
The coefficients $a_w$ are nonnegative integers. We define a current-valued map 
\[
\eta_S: \ol{\Ga}= \Gamma\cup \geo \Gamma\rightarrow \mathcal{D}^{n-k,n-k}(\cF)
\] 
as follows:
\begin{equation}\label{eq:eta_flags1}
\eta_S(u) = \left  \lbrace
\begin{array}{ll}
u_*(\llbracket S\rrbracket ) & \text{if } u\in \Gamma,\\
 \sum_{w\in W} a_w \llbracket \Th_{w}(f(u))\rrbracket  &\text{if }u\in \geo \Gamma.
\end{array} \right .
\end{equation}

\begin{example}\label{ex:discontinuous}
Let $\Ga< PSL(2,\C)$ be a (nonelementary) convex-cocompact Kleinian group. Take $S=\{\la_-\}\subset \La$, the repelling fixed point 
of a hyperbolic element $\ga\in\Ga$.  Let  $\la_+$ denote the attractive fixed point of $\ga$; let $\xi_\pm:=f^{-1}(\la_\pm)$. 
Then the sequence $(\ga^k)_{k>0}$ converges to $\xi_+$ in $\ol\Ga$ and $f(\xi_+)=\llbracket \la_+\rrbracket$. 
At the same time, the sequence of currents $\eta(\ga^k)=\ga_k \llbracket S\rrbracket=\llbracket S\rrbracket$ converges to 
$\llbracket S\rrbracket=\llbracket \la_-\rrbracket$. Thus, $\eta$ is discontinuous at $\la_+$. 
\end{example}

The goal of this section is to prove that $\eta_S$ is continuous on the closure of $\Gamma_{S,t}$ in $\ol\Ga$ (see \eqref{eq_defi_gammaSt_flags1} for 
the definition of the subset $\Gamma_{S,t}\subset \Ga$). 

\subsection{Setup}\label{sec:hypothesis_cylinder_stretch}
In the rest of this section we fix $t>0$, as well as the constants $\epsilon=\epsilon_t$ and $N_t$ determined by $t$ from Lemma \ref{uniformboundedawayfromthelimitset1}.
Throughout this section, when we use the letters $\bx,\gamma,g,\tga$, we will be considering an element $w\in W$ with $\vert w\vert=k$, a pair of flags $\bx=\xpair$, and a redressing $\tga=\gamma\circ g$ of an element $\gamma\in \Gamma_{S,t}$ with $\beta(\gamma)>N_t$ such that $x^+$ (resp. $x^-$) is the attractive (resp. repelling) fixed point of $\tga$. Note that $\gamma,g,\tilde{\gamma}=\gamma\circ g$ satisfy the conclusion of Proposition \ref{Mishaseigenvaluelemma1}, in other words, we have $\bx\in \sK_\Lambda$.

We also note that the subvariety $S\subset \cF$ of dimension $k$ satisfies the genericity condition \eqref{eq:the_genericity_condition} and by Lemma \ref{uniformboundedawayfromthelimitset1} we have:
\begin{equation}\label{eq:epsilon_condition_for_single_gamma}
d(g^{-1}S,\Th^{{n-k-1}}(x^-))>\epsilon.
\end{equation}
Remark that $\tga(g^{-1}S)=\ga(S)$ and this subvariety satisfies for any $x\in \cF$
\begin{align*}
[\ga(S)]=\sum_{\vert w\vert=k}a_w[\Th_w(x)]
\end{align*}
because $[S]=\sum_{\vert w\vert=k}a_w[\Th_w(x)]$ and every element of the connected Lie group $G$ acts trivially on the homology.

\subsection{Cylinders}
We introduce some objects that will be used throughout this section. If $\bx=\xpair\in \sK_\Lambda$, then for any $m_1,m_2\dyl$ and $w\in W$, we define two cylinder-like subsets of $\Oppxwd$:
\begin{align}\label{eq:definition_cylinder}
\begin{aligned}
& \Cyli(\bx,w,m_1):=\Psi_{w,\bx}^{-1}\left(\mathbf{B}(0,m_1)\oplus \fuwm \right),\\
& \Cylii(\bx,w,m_1,m_2):=\Psi_{w,\bx}^{-1}\left(\mathbf{B}(0,m_1)\oplus \mathbf{B}(0,m_2)\right),
\end{aligned}
\end{align}
where $\mathbf{B}(0,m_1)\subset \fuwp$ and $\mathbf{B}(0,m_2)\subset \fuwm$ are balls with respect to the Euclidean metric on $\fu_w$. When $\bx,w$ are fixed we simply denote them by $\Cyli(m_1), \Cylii(m_1,m_2)$. Note that $\Cylii(m_1,m_2)$ is relatively compact in $\Oppxwd$ but $\Cyli(m_1)$ is not.

\begin{lem}\label{finiteintersectionlemma1}
The subvariety $g^{-1}(S)$ intersects $S_{w^{\vee}}(x^-)$ at $a_w$ points (counted with multiplicity).
\end{lem}
\begin{proof}
The complement $ \Th_{<w^\vee}(x^-)=\Th_{w^\vee}(x^-)\setminus S_{w^\vee}(x^-)$ is a thickening of dimension $n-k-1$ of $x^-$. By \eqref{eq:epsilon_condition_for_single_gamma}, $g^{-1}S$ does not intersect $ \Th_{<w^\vee}(x^-)$. Therefore the intersection between the two projective subvarieties $S$ and $\Th_{w^\vee}(x^-)$ is contained in $S_{w^\vee}(x^-)$. Since $S_{w^\vee}(x^-)$ is an affine variety while the intersection locus is projective, we claim that this intersection has to be of dimension $0$, 
i.e.\ a finite set of points. The number of intersection points counted with multiplicity is given by the homological coefficient $a_w$.
Indeed assume that $S \cap \Th_{w^\vee}(x^-)$ contains an algebraic curve $C$. Then the closure of $C\cap S_{w^\vee}(x^-)  $ contains a point in $$\Th_{w^\vee}(x^-) \setminus S_{w^\vee}(x^-) =  \Th_{<w^\vee}(x^-),$$ 
which contradicts the inequality \eqref{eq:epsilon_condition_for_single_gamma}.
\end{proof}

\subsection{Cylinders of uniform sizes}

\begin{lem}\label{cylindertwoconstants1}
There exist two constants $r_1,r_2>0$ such that for any $w\in W$ with $\vert w\vert=k$, for any $\bx=\xpair \in \sK_\Lambda$ with $x^-\in \Lambda_\epsilon$, for any $\gamma\in \Gamma$ with redressing $\tga=\gamma\circ g$, we have
$$
g^{-1}(S)\cap \Cyli(\bx,w,r_1)=g^{-1}(S)\cap \Cylii(\bx,w,r_1,r_2),
$$
and that $\Psi_{w,\bx}^+$ induces a ramified covering 
$$
g^{-1}(S)\cap \Cylii(\bx,w,r_1,r_2)\rightarrow \mathbf{B}(0,r_1)\subset \fuwp.
$$
\end{lem} 
\begin{proof}
As $W$ is a finite set, we can assume without loss of generality that $w$ is fixed. 
Let $\bx=\xpair \in \sK_\Lambda$ with $x^-\in \Lambda_\epsilon$. Then by the definition of $\Lambda_\eps$, the $\dF$-distance between $S$ and $\Th^{n-k-1}(\Lambda_\epsilon)$ is bounded from 
below by $\epsilon$. Recall that $\Th_{<w^\vee}(x^-)\subset \Th^{n-k-1}(x^-)$ (see \eqref{thickenings_boundary_inclusion}). We obtain by Lemma \ref{uniform_r1r2} that there exists $r_2'>0$ such that for any $\bx \in \sK_\Lambda$ with $x^-\in \Lambda_\epsilon$, for any $z\in S\cap \Swxd$, we have $\vert \Psi_{w,\bx}^{\bac{-}}(z)\vert<r_2'$ (here we can assume without loss of generality that $\eps$ is small enough so that Lemma \ref{uniform_r1r2} can be applied).

We fix an arbitrary constant $r_2>r_2'$ and for $\bx \in \mathcal{K}_{\Lambda}$, we consider the following set $A_{\bx}$ given by: 
\begin{equation*} 
A_{\bx} = \{r \geqslant 0 \ | \ g^{-1}(S)\cap \mathcal{C}_I(\bx, w , r) = g^{-1}(S) \cap \mathcal{C}_{II}(\bx,w, r, r_2) \}. 
\end{equation*}
Let us show that $A_\bx$ is non empty for each $\bx \in \mathcal{K}_{\Lambda}$. 

Note that, according to Lemma \ref{finiteintersectionlemma1}, 
for a fixed $x^-\in \Lambda_\epsilon$, the intersection $S\cap\Th_{w^\vee}(x^-)$ is finite. 
For each intersection point $p_i \in g^{-1}(S) \cap \Th_{w^{\vee}}(x^-)$ with $i \leqslant a$, we choose an open polydisk $\mathcal{D}_i$ around $p_i$ contained in the horizontal cylinder $\Psi_{w,\bx}^{-1}(\mathfrak{u}_w^+ \oplus \mathbf{B}(0,r_2))$. 
Up to reducing the size of the polydisk, we can ensure that the $\mathcal{D}_i$ are pairwise disjoint.  
 The image of each $\mathcal{D}_i \cap g^{-1}(S)$ by the projection $\Psi_{ w, \bx}^+$ is open  since the restriction of $\Psi_{w,\bx}^+$ to $g^{-1}(S)$ is proper.
Reducing these disks, we can ensure that $\Psi_{w,\bx}^+$ restricted to $\mathcal{D}_i \cap g^{-1}(S)$ is a ramified covering, as the fiber over $0$ is a finite set of points.  
  Since these images all contain the origin in $ \C^n$, these images all contain a polydisc of some radius $r>0$. We thus have: 
\begin{align*}
g^{-1}(S) \cap \mathcal{C}_I(\bx, w,r) = g^{-1}(S) \cap (\Psi_{w,\bx}^{+})^{-1} (\mathbf{B}(0,r)) \subset \\
\coprod_{i=1}^a  (g^{-1}(S) \cap \mathcal{D}_i) \subset g^{-1}(S) \cap (\Psi_{w, \bx}^-)^{-1} (\mathbf{B}(0,r_2)). 
\end{align*} 
 Intersecting with $\mathcal{C}_{I}(\bx,w,r)$ yields the inclusion 
 $ g^{-1}(S) \cap \mathcal{C}_I(\bx,w,r) \subset g^{-1}(S) \cap \mathcal{C}_{II}(\bx,w, r,r_2)$, and since $\mathcal{C}_{II}(\bx,w,r,r_2)  \subset \mathcal{C}_I(\bx, w,r)$ holds by definition, this shows that $g^{-1}(S) \cap \mathcal{C}_I(\bx, w, r) = g^{-1}(S) \cap \mathcal{C}_{II}(\bx,w, r, r_2)$ and $r \in A_{\bx}$.
\smallskip 
  
  Observe furthermore that if $r \in A_{\bx}$, then the whole interval $[0,r]$ is contained in $A_{\bx}$. Let us now consider $r(\bx) = (\sup A_{\bx})/2 > 0 $ for $\bx \in \mathcal{K}_{\Lambda}$.  
We will show that $r_1 =\inf_{\bx \in \mathcal{K}_{\Lambda}} r(\bx) \neq 0$. 
Assume for the sake of a contradiction that $r_1 = 0 $; then there exists a sequence of points 
$\bx_i \in \mathcal{K}_{\Lambda}$ such that $r(\bx_i)< 1/i$ converges to zero. 
By the compactness of $\mathcal{K}_{\Lambda}$, we can assume that the sequence $\bx_i$ converges to a point $\bx$. We can thus construct a sequence of points 
$$
p_i \in g^{-1}(S) \cap \mathcal{C}_I(\bx_i ,w, 1/i) \setminus \mathcal{C}_{II}(\bx_i ,w, 1/i,r_2).
$$ 
Let us first assume that the sequence $(p_i)$ is bounded and remains in the chart $\Opp(x_{w^\vee})$. 
We can  assume (up to taking a subsequence) that $p_i$ converges to $p \in g^{-1}(S) \cap \Opp(x_{w^\vee})$. We claim  that $p$ projects to $0$ under $\Psi_{w, \bx}^+$. Indeed, the family of maps $\bx \mapsto (\Psi_{w,\bx}^+)_{| \Opp(x_{w^\vee})} $ is continuous (even algebraic) and $p_i$ projects under $\Psi_{w,\bx_i}^+$ to a sequence of point converging to $0$, hence, by continuity, the limit $p$ projects to $0$ under $\Psi_{w,\bx}^+$, which implies that $p \in \Th_{w^{\vee}}(x_-)$. 
Similarly, since $p_i \in \mathcal{C}_{I}(\bx_i,w, 1/i) \setminus \mathcal{C}_{II}(\bx_i , w ,1/i , r_2)$, the limit $p$ does not belong to $(\Psi_{w,\bx}^-)^{-1} (B(0,r_2 - \epsilon))$ for all $\epsilon>0$. This shows that $p \in \Th_{w^\vee}(x_-) \cap g^{-1}(S) \setminus \mathcal{C}_{II}(\bx,w , r(\bx), r_2 - \epsilon)$. 
Choosing $\epsilon$ so that $r_2 - \epsilon > r_2'$, we obtain 
a contradiction with the fact that  $g^{-1}(S) $ intersects $\Th_{w^\vee}(x_-)$ at exactly $a$ points counted with multiplicity as we have found $a+1$ intersections points. 

This shows that if $(p_i)$ converges, then it must converge to a point $p$ outside of $\Opp(x_{w^\vee})$. 
Recall from \S~\ref{summary_torus_action} that $\Opp((x_{w^\vee})$ is isomorphic to the image of the map $ U_{x_w} \mapsto U_{x_w} x_w \in \cF $. 
By the knit product decomposition $U_{x_w} = U_{x_w}^+ U_{x_w}^-$, the point $p_i$ can be written as $p_i = u_i^+ u_i^- x_w$, 
where $x_i^{\pm} \in U_{x_w}^\pm$. Similarly, the condition that $p_i \in \mathcal{C}_{I}(\bx_i , w, 1/i)$ shows that $p_i$ can also be written as: 
\begin{equation*}
p_i = u_i^+ u_i^- x_w = v_i^+ v_i^- (x_{i})_w,
\end{equation*}
where $v_i^{\pm} \in U_{(x_i)_w}^{\pm}$. 
As $p_i \in \mathcal{C}_I(\bx_i , w , 1/i)$, this shows that $(v_i^+$ converges to $1$. Using the fact that $(\bx_i)$ converges to $\bx$, we thus deduce that 
\begin{equation*}
p = \lim_{i \rightarrow +\infty} v_i^- x_w.
\end{equation*}
As $i\rightarrow +\infty$, the sequence of subgroups $U_{(x_i)_w}^-$ converges to $U_{x_w}^-$; hence,  we can write 
$U_{(x_i)_w}^- = g_i U_{x_w}^- g_i^{-1}$ where $(g_i)$ converges to the identity. 
Since the sequence $(v_i^-)$ is unbounded, we have $v_i^-  = g_i h_i g_i^{-1}$, where the sequence $(h_i)$ 
is unbounded in $U_{x_w}^-$. This gives: 
\begin{equation}
p = \lim_{i\rightarrow +\infty} g_i h_i g_i^{-1} x_w = \lim_{i\rightarrow +\infty} h_i x_w  \in \overline{S_{w^\vee}(x^-)} = \Th_{w^\vee}(x^-).
\end{equation} 
Since $p \in g^{-1}(S)$, this shows that the intersection $g^{-1}(S) \cap \Th_{w^\vee}(x^-) $ has $a+1$ points counted with multiplicity, which is impossible. 

We have thus shown that $r_1 > 0$, hence $g^{-1}(S) \cap \mathcal{C}_{II}(\bx, w, r_1,r_2) = g^{-1}(S) \cap \mathcal{C}_{I}(\bx, w, r_1)$. 
The fact that the restriction of $\Psi_{w,\bx}^+$ to this subset is ramified onto its image follows from our construction. \end{proof}

%

\begin{lem}\label{disjointcyliders1}
In Lemma \ref{cylindertwoconstants1}, we can choose $r_1$ so that, in addition to the conclusion of Lemma \ref{cylindertwoconstants1}, the following is satisfied: For fixed $\bx\in \sK_\Lambda$ with $x^-\in \Lambda_\eps$, for $g\in E$ and for any two 
elements $w_1 \neq w_2 \in W$ with $|w_1|=|w_2|=k$, the intersection 
$$
g^{-1}(S)\cap \Cylii(w_1,r_1,r_2) \cap \Cylii(w_2,r_1,r_2)
$$
 is empty.
\end{lem}

\begin{proof}
It suffices to prove the claim for $g=1$. By using the compactness of $\sK_\Lambda$ as in the proof of Lemma \ref{cylindertwoconstants1}, we can assume that $\bx$ is fixed. 

When $w$ runs through $\{w\in W: \vert w\vert=k\}$, the dual $w^\vee$ runs through $\{w\in W: \vert w\vert=n-k\}$. We list all $w$ with $\vert w \vert=k$: $w_1,w_2,\cdots,w_j$. The intersections between the thickenings $\Th_{w^\vee}(x^-), \vert w\vert=n-k$ are thickenings of $\{x^-\}$ of dimension $\leq n-k-1$. They do not intersect with $S$ because $x^-\in \Lambda_\eps$. They are uniformly bounded away from $S$ by compactness. By definition of $\Cylii$, for fixed $r_2$, the set $S\cap\Cylii(w,r_1,r_2)$ can be thought of as living in a tubular neighborhood of $\Th_{w^\vee}(x^-)$. Therefore we can find $r_1(w_1)$ small enough so that 
\[
S\cap\Cylii(w_1,r_1(w_1),r_2)\cap \Th_{w'}(x^-)=\emptyset
\]
for any $w'$ such that $\vert w'\vert=n-k$ and $w'\neq w_1^\vee$. Then we can find $r_1(w_2)\leqslant r_1(w_1)$ such that 
\begin{align*}
& S\cap\Cylii(w_2,r_1(w_2),r_2)\cap \Th_{w'}(x^-)=\emptyset
\\
& \Cylii(w_2,r_1(w_2),r_2)\cap \Cylii(w_1,r_1(w_1),r_2)=\emptyset
\end{align*}
for any $w'$ such that $\vert w'\vert=n-k$ and $w'\neq w_1^\vee,w_2^\vee$. We continue and in the end we take $r_1=r_1(w_j)$. 
\end{proof}

\subsection{Stretching cylinders}

We will fix $r_1,r_2$ so that conclusions of Lemma \ref{cylindertwoconstants1} and \ref{disjointcyliders1} hold. Recall that $\la(\ga)$ is the Lyapunov spectrum of an element $\ga\in G$. 
 
\begin{lem} \label{keylemma1}
For any $q>0$, there exists $M \geq 0$ such that for any smooth $(k,k)$-form $\varphi$ on $\cF$ the following assertion holds. If $\gamma,g,\tga,\bx=\xpair,w,k,\eps,S$ satisfy the hypothesis in \ref{sec:hypothesis_cylinder_stretch} and if 
\[
\exp\left(\alpha(\la(\tga))\right)>M
\]
for all positive roots $\alpha$, 
then there exists a subset
\[
\fS\subset \tga\left(g^{-1}(S)\cap \Cylii(\bx,w,r_1,r_2)\right)
\]
such that
\[ 
\left\lvert\langle \llbracket \fS\rrbracket  - a_w\llbracket \Th_{w}(x^+)\rrbracket , \varphi \rangle \right\rvert \leqslant q\max(  \lVert \varphi \rVert_{\infty}, \lVert d\varphi\rVert_\infty).
\]
\end{lem}

\begin{proof}
Fix $q>0$. We are looking for an $M$ that works for arbitrary $\tga,\bx,w$ and $\varphi$ as in the statement. In what follows, for the ease of notation we often omit suffixes depending on $\tga,\bx,w$ and the constants we obtain will always be uniform. Without loss of generality we will assume that $g$ is the neutral element so that $\gamma=\tga$.



Recall that $\omega$ is our fixed K\"ahler form on $\cF$. By Lemma \ref{comparison_metric_balls}, we can take $D\dyl$ so that 
\begin{equation}\label{eq:introduce_the_constant_D}
\left\lvert \langle \llbracket(\Psi_{w,\bx}^+)^{-1}(\mathbf{B}(0,D))\rrbracket - \llbracket \Th_w(x^+)\rrbracket, \omega^k\rangle\right \rvert
\end{equation}
is small enough so that, by Corollary \ref{cor_difference_current} we get
\begin{equation}\label{eq:ball_approximating_thickening}
\left\lvert \langle \llbracket(\Psi_{w,\bx}^+)^{-1}(\mathbf{B}(0,D))\rrbracket - \llbracket \Th_w(x^+)\rrbracket, \varphi\rangle\right \rvert \leqslant \frac{q}{2a_w}\max(  \lVert \varphi \rVert_{\infty}, \lVert d\varphi\rVert_\infty).
\end{equation}

The action of $\gamma$ on $\fu_w$ is conjugate to the linear map $e^{ad(\la(\gamma))}$. In the sequel we assume that 
\[\exp(\alpha(\la(\gamma)))>\frac{D}{r_1}\]
for any positive root $\alpha$, so
\begin{equation}\label{eq:horizontal_big_ball}
\Psi_{w ,\bx}^+\left(\gamma\left(S\cap\Cylii(r_1,r_2)\right)\right)  \supset \mathbf{B}(0,D)\subset \fuwp.
\end{equation}
We then set $U=\mathbf{B}(0,D)\subset \fuwp$ and
\begin{equation*}
\fS=\gamma(S)\cap \Cyli(D)\subset \gamma\big(S\cap \Cylii(r_1,r_2)\big).
\end{equation*}
Note that 
\[\exp(\alpha(\la(\gamma)))<\frac{r_1}{D}\]
for any $\alpha\in \Phi_w^-$ and consequently
\[\fS=(\Psi_w^+)^{-1}(U)\cap \gamma(S)=\gamma(S)\cap\Cylii\big(D,\frac{r_1r_2}{D}\big).\] 
We write $\fS_\bx$ when we need to emphasize on the dependence on $\bx$. 
Choose a compact set $K$ containing both $U$ and $\fS_{\bx}$. 
We fix $\epsilon'$ and choose a smooth bump function $\chi_{\epsilon'}$ which is identically $1$ on  $K$ and $0$ outside an $\epsilon'$-neighborhood of $K$. We then decompose $\varphi$ into $\varphi = \chi_{\epsilon'} \varphi + (1 - \chi_{\epsilon'}) \varphi$.

Since $\chi_\epsilon \varphi$ is supported in an $\epsilon$-neighborhood of $K$,  by Corollary \ref{cor_subvariety_smaller_cylinder} we obtain that 
\begin{multline*}
\left\lvert\langle \llbracket \fS \rrbracket - a_w\llbracket U \rrbracket, \varphi \rangle\right\rvert  = \left \lvert  \langle \intcur{\fS} - a_w \intcur{U} , \chi_{\epsilon'} \varphi \rangle \right \rvert  \leqslant C\sqrt{\fd} \max(  \lVert \chi_{\epsilon'} \varphi \rVert_{\infty}, \lVert d(\chi_{\epsilon'}\varphi)\rVert_\infty) \\  \leqslant C \sqrt{\fd} \max ( \lVert d\chi_{\epsilon'} \rVert_\infty \lVert \varphi \rVert_{\infty}, \lVert d\varphi \rVert_\infty )
\end{multline*}
where ($\cah$ is defined before Corollary \ref{cor:pizza_BCH})
\[
\fd=\max_{\bx\in \sK_\Lambda, z\in \fS_\bx}\cah(\bx,w,z)
\]
and $C$ is a constant depending only   on the degree of $S$. 
Since $\epsilon'$ is fixed, the derivative $d\chi_{\epsilon'}$ is bounded above so the previous equation gives: 
\begin{equation}\label{eq:from_prop_cylinder}
\left\lvert\langle \llbracket \fS \rrbracket - a_w\llbracket U \rrbracket, \varphi \rangle\right\rvert \leqslant C'\sqrt{\fd} \max ( \lVert \varphi \rVert_\infty, \lVert d\varphi \rVert_\infty ),
\end{equation} 
where  $C' = C \lVert d\chi_{\epsilon'}\rVert_\infty$ is a uniform constant, independent of $q$.

Recall that $D$ was introduced in \eqref{eq:introduce_the_constant_D}, \eqref{eq:ball_approximating_thickening}. It depends only on $q$. Assume that some number $M$ satisfies $M>\frac{D}{r_1}$ and 
\[\exp(\alpha(\la(\gamma)))>M\]
for any positive root $\alpha$. Then we have 
\[
\fS=\gamma(S)\cap\Cylii\big(D,\frac{r_2}{M}\big).
\]

Applying Corollary \ref{cor:pizza_BCH} to $R_h =D$, $R_v= M/r_2$ and to $\epsilon = q/(4C'^2)$, if $M$ is chosen so that  $M>>D^{N_0}$ where $N_0$ is the order of nilpotence appearing in Corollary \ref{cor:pizza_BCH},  then 
$
\fd\leqslant 
\dfrac{q^2}{4{C'}^2},
$
which combined with \eqref{eq:from_prop_cylinder} gives 
\begin{equation}\label{eq:prop_cylinder_and_height}
\left\lvert\langle \llbracket \fS \rrbracket - a_w\llbracket U \rrbracket, \varphi \rangle\right\rvert \leqslant \dfrac{q}{2} \max(  \lVert \varphi \rVert_{\infty}, \lVert d\varphi\rVert_\infty).
\end{equation}
Now \eqref{eq:ball_approximating_thickening} and \eqref{eq:prop_cylinder_and_height} imply
\begin{equation*}
\big\lvert\langle \llbracket \fS\rrbracket - a_w\llbracket \Th_{w}(x^+)\rrbracket , \varphi \rangle \big\rvert \leqslant q\max(  \lVert \varphi \rVert_{\infty}, \lVert d\varphi\rVert_\infty)
\end{equation*}
and the proof is finished.
\end{proof}

\begin{lem}\label{lem:geometric_redressing}
For any $m>0$, for any smooth $(k,k)$-form $\varphi$ on $\cF$, there exists $N \geq 0$ so that the following holds. If $\gamma \in \Gamma_{S,t}$ satisfies $\beta(\gamma)>N$ and if $x^+,x^-$ are the attracting/repelling fixed points of $\tga$, then we have
\[
|\langle \llbracket\ga(S)\rrbracket  - \sum_{\vert w\vert=k}a_w\llbracket \Th_{w}(x^+)\rrbracket , \varphi \rangle | \leqslant m.
\]
\end{lem}

\begin{proof}
Recall that we work in the setup of \ref{sec:hypothesis_cylinder_stretch}. In particular $\beta(\gamma) >N_t$, $x^-\in \Lambda_{\epsilon_t}$. Recall also that we have fixed $r_1,r_2$ so that Lemmas \ref{cylindertwoconstants1} and \ref{disjointcyliders1} hold.     

Let us fix $\varphi$ and $m$. Fix $q>0$, to be chosen later on.

Proposition \ref{Mishaseigenvaluelemma1}  says that $\alpha(\la(\tilde{\gamma}))\geqslant C\beta(\tga)$ for some constant $C$ (when $\beta(\tga)$ is large). Therefore by Lemma \ref{keylemma1}, for any $q >0$, there exists $M>0$ such that  for any $\gamma$ with $\beta(\tga)> M$ is large enough, there exists some  subsets
\[
\fS_w\subset\tga\big(g^{-1}(S)\big)
\]
such that, for all $w$ with $\lvert w \rvert=k$ and for all smooth test form $\varphi$ of bidegree $(k,k)$ on $\cF$, 
\begin{equation*}\label{eq:constant_m1}
\big\lvert\langle \llbracket \fS_{w}\rrbracket - a_w\llbracket \Th_{w}(x^+)\rrbracket , \varphi \rangle \big\rvert \leqslant q \max ( \lVert \varphi \rVert_\infty , \lVert d\varphi \rVert_\infty ).
\end{equation*}
Moreover, we have 
\[
\fS_{w}\subset \tga\left(g^{-1}(S)\cap \Cylii(\bx,w,r_1,r_2)\right).
\]
Thus, by Lemma \ref{disjointcyliders1}, the subsets $\fS_w,\lvert w\rvert=k$, are pairwise disjoint and hence their union $\fS_\tga$ is a subset of $\tga(g^{-1}(S))=\ga(S)$ such that for any smooth test $(k,k)$-form $\varphi$, one has 
\begin{equation} \label{eq_flat_term}
\vert\langle \llbracket\fS_\tga \rrbracket- \sum_{\vert w\vert=k}a_w\llbracket \Th_{w}(x^+)\rrbracket, \varphi \rangle \vert\leqslant N q \max ( \lVert \varphi \rVert_\infty , \lVert d\varphi \rVert_\infty ), 
\end{equation}
where $ N = \# \{ w \in W \ | \ |w|=k \}$.

Recall that the homology class of $\tga(g^{-1}S)=\ga(S)$ equals 
\begin{align*}
[\ga(S)]=\sum_{\vert w\vert=k}a_w[\Th_w] \in H_{k,k}(\cF).
\end{align*}
Equation  \eqref{eq_flat_term} together with the above formula show that the hypothesis  of Corollary \ref{cor_difference_current} are satisfied for  $V =\ga(S)$, $T = \sum_{w} a_w \intcur{\Th_w(x_+)}$, $W = \fS_{\tilde{\gamma}}$ and $\epsilon = K q$, hence there exists a constant $C' > 0$ such that for any smooth $(k,k)$-form $\varphi$, 
\[
\vert\langle \llbracket\ga(S)\rrbracket - \sum_{\vert w\vert=k}a_w\llbracket \Th_{w}(x^+)\rrbracket , \varphi \rangle \vert \leqslant C' N q  \max( \lVert \varphi \rVert_\infty , \lVert d\varphi \rVert_\infty  ), 
\]
 For our given $m, \varphi$, we thus set $q = m/(C' N \max (\lVert \varphi \rVert_\infty , \lVert d\varphi \rVert_\infty )) $. Recall that we chose $M$ so that if $\beta(\tilde{\gamma})>M$ then the above inequality holds, and this  concludes the proof.
\end{proof}

We can finally prove the main result of this section:  

\begin{thm}\label{continuousonGammaSt}
For every $t> 0$, the map $\eta_S$ is continuous on the closure $\overline{\Gamma_{S,t}}$.
\end{thm}

\begin{proof}
Remark that the inequality $\beta(\tilde{\gamma}) - \max_{g\in E} \beta(g)\leqslant \beta(\gamma) \leqslant \beta(\tilde{\gamma}) + \max_{g\in E} \beta(g)$ shows that $\beta(\gamma)$ is large enough is equivalent to the condition that $\beta(\tga)$ is large enough.
By Proposition \ref{Mishaseigenvaluelemma}, these conditions are equivalent to the condition that the Lyapunov spectrum of $\tga$ is large enough, i.e $\forall i\leqslant r$,  $\lambda_i(\tga) \geqslant M$ for some $M>0$. 
By Lemma \ref{lem:geometric_redressing}, the hypothesis of Corollary \ref{cor_equivariant_continuous} are satisfied for $\Psi = \intcur{S}$ and for $\Gamma^* = \Gamma_{S,t}$.  Thus the map $\eta_S$ is continuous on the closure $\overline{\Gamma_{S,t}}$ as required.
\end{proof}

\section{Proof of the equidistribution theorem, Theorem \refthmmain}\label{sec:mainproof}

\begin{proof}
Let $S\subset \cF$ be a subvariety satisfying the genericity condition. Without loss of generality 
we can assume that $S$ is irreducible. Consider the current-valued map $\eta_S:\ol{\Ga}\rightarrow\mathcal{D}_{k,k}(\cF)$ defined in \eqref{eq:eta_flags1}. To prove Theorem \ref{thm:main} we need to prove that 
\begin{equation}
\label{eq:main_convergence_current}
\lim_{m\rightarrow +\infty}\int_\Gamma \eta_S d\mu_{s_m}=\int_{\geo \Ga}\eta_S d\mu.
\end{equation}

Let $\varphi$ be a smooth $(k,k)$-form on $\cF$ and let $\psi:\overline{\Gamma}\rightarrow \bC$ be the function defined by
\[
\psi(x)=\langle \eta_S(x),\varphi \rangle.
\]
The weak convergence \eqref{eq:main_convergence_current} is the same as the numerical convergence (for all $\varphi$):
\begin{equation}\label{eq:main_convergence1_current}
\lim_{m \rightarrow +\infty} \int_{\overline{\Ga}}\psi(x)d\mu_{s_m}(x) = \int_{\geo \Ga} \psi(p)d\mu_{\Ga}(p).
\end{equation}

We will be using Theorem \ref{thm:integration_partial_measures}. Let us explain some data needed in Theorem \ref{thm:integration_partial_measures} in our setting. The hyperbolic space $X$ in that theorem will be simply (the Cayley graph of) $\Gamma$ and the base point $x_0 $ will be 
the neutral element.
Since $\Gamma$ was an Anosov subgroup, Theorem \ref{Mishaseigenvaluelemma1} produced a finite subset $E \subset \Gamma$. 
 Following this, we defined the subset $\Gamma_{S,t} = \Gamma_{P,E,t}\subset \Ga$ in \eqref{eq_defi_gammaSt_flags1}.  

By the definition of $\eta_S$, the function $\psi$ is continuous on $\geo \Ga$. It is also continuous on $\Gamma$ because 
$\Gamma$ is discrete. For every $t\dyl$, $\psi$ is continuous on $ \overline{\Gamma_{S,t}} = \ol{\Gamma_{P,E,t}}$ by Theorem \ref{continuousonGammaSt}, showing that the hypothesis of Theorem \ref{thm:integration_partial_measures} hold. Hence, Theorem \ref{thm:integration_partial_measures} yields the convergence \eqref{eq:main_convergence1_current}, and, thus, Theorem \ref{thm:main} follows.
\end{proof}

\section{Examples}\label{sec:examples}

We will start with a general observation about points in flag manifolds which are not $\La(\Ga)$-generic.  

\begin{prop}
Let $G$ be a complex semisimple Lie group with $n$-dimensional flag manifold $\cF=G/B$. Let 
$\Ga< G$ be a nonelementary Anosov subgroup with flag-limit set $\La=\La(\Ga)\subset \cF$. 
Suppose that $x\in \cF$ is such that $\La\subset \Th^{n-1}(x)$, i.e. $x$ is $\La$-nongeneric.  
Then the closure of the orbit $\Ga x\subset \cF$ is disjoint from $\La$. 
\end{prop}
\begin{proof} Suppose that there exists a sequence $\gamma_i\in \Ga$ such that 
$$
\lim_{i\to\infty} \ga_i(x)=y\in \La.$$ 
For every $i$ we have $\La\subset \Th^{n-1}(x_i)$, $x_i=\ga_i(x)$. By Lemma \ref{lem:compactness}, it follows that 
$\La\subset \Th^{n-1}(y)$. Since the limit set $\La$ is an antipodal subset of $\cF$, we have $\{y\}=\La$. This is a contradiction. 
\end{proof}

\begin{cor}
In the setting of the proposition, assume that $\La\subset \Th^{n-1}(x)$. Then 
Theorem \ref{thm:main-points} fails for the point $x$. 
\end{cor}
\begin{proof} Let $\beta$ be any hyperbolic length function on the group $\Ga$, let  
$\cP(s)$ be the corresponding Poincar\'e series and $s_0$ its critical exponent. According to the proposition, the accumulation set of 
$\Ga x$ is disjoint from $\La$. 
Therefore, the support of any limiting measure 
$$
\mu_{x}= \lim_{s_j\to s_0} \frac{1}{\cP(s_j)} \sum_{\ga\in \Ga} e^{-s_j\beta(\ga)} \delta_{\ga x}
$$
is disjoint from $\La$. In particular, $\mu_x$ cannot be the Gibbs measure on $\La(\Ga)$, which means that 
Theorem \ref{thm:main-points} fails for $x$. 
\end{proof}

Below,  we will discuss generic and nongeneric subvarieties and limiting Gibbs currents for some classes of Anosov subgroups $\Ga< G=\SL(3,\C)$  (except for the very last example where the Lie group will be a product of rank 1 Lie groups) and generic subvarieties of flag-manifolds. We let $\cF=G/B$ denote the full flag-manifold of $G=\SL(3,\C)$. Recall that $\Th^2(x)= \cF\setminus \Opp(x)$ is the maximal (2-dimensional) thickening of $x\in \cF$.

Since all Zariski dense (over $\bC$) Anosov subgroups $\Ga< G$ have the property that every subvariety $S\subsetneq \cF$ is $\La(\Ga)$-generic, we consider subgroups which are not Zariski dense. However, we will assume that $\Ga$ is nonelementary.

\begin{example}\label{ex:E1}
We consider Anosov representations coming from an irreducible representation 
$\SL(2,\C)\to \SL(3,\C)$ described in Example \ref{ex:Anosov1}. 
The flag-limit set $\La=\La(\Ga)\subset \cF$ is contained in the rational curve 
$\tilde V=\tilde f(\P^1)\subset \cF$. 

Every hypersurface $S\subset  {\mathcal F}$ intersects $\tilde V$ in a finite subset or contains $\tilde V$. In the former case, $S$ is $\La$-generic and in the latter case, it is not. For instance, one can take $S$ to be the preimage of $V$ in $\cF$ under the projection $\cF\to \P^2$. Such $S$ is $\Ga$-invariant and, clearly, Theorem \ref{thm:main} fails for this surface. 

Since $V$ is a nonsingular quadric, for any three distinct points $\lambda_i\in V$ the projective lines $\ell_{\lambda_i}$ 
(tangent to $V$ at $\lambda_i$) have empty triple intersection.   
Hence, for every Anosov subgroup $\Gamma< \SL(3,\C)$ preserving $V$, every 
$x\in \cF$  is $\Lambda$-generic. 

Suppose that $S$ is a curve in $\cF$. If the projection of $S$ to $\P^2$ or $(\P^2)^\vee$ (the space of lines in $\P^2$) contains $V$ or its dual quadric $V^\vee$, then $S$ is obviously non-generic with respect to $\La$. Moreover, the $\Ga$-orbit of $S$ is contained in the preimage of $V$ (or $V^\vee$) 
in $\cF$ which is $\Ga$-invariant. In the extreme case, $S=\tilde V$, the curve $S$ is $\Ga$-invariant and, thus, Theorem \ref{thm:main} fails for such $S$. 
Suppose that both projections $C, C^\vee$ of $S$ to $\P^2$ and $(\P^2)^\vee$ are different from 
$V$ and $V^\vee$ respectively. Then $\Th^1(S)\cap \tilde V$ is a finite subset and, hence, $S$ is $\La$-generic. 
\end{example}

 \begin{example}\label{ex:E2}
 We take a discrete subgroup $\Gamma< \SL(2,\C) < \SL(3,\C)$, as in Example \ref{ex:Anosov2}. The subgroup $\SL(2,\C)$ preserves a projective line 
 $L\subset \P^2$ and a point $q\in \P^2$ not incident to $L$. The flag-limit set $\La$ of $\Ga$ is contained in a lift to $\cF$ 
 of the limit set $\La_1(\Ga)\subset L$: Each $p\in \La_1(\Ga)$ lifts to the flag $(p, \ell_p)$, where $\ell_p$ is the line incident to both $p$ and $q$.  
 Consider any flag $x=(q, \ell)$. Then $\Th^2(x)$ contains the lift $\tilde L$ of the entire $L$ to $\cF$, hence, $x$ is $\La$-nongeneric. The set of flags of the form  $(q,\ell)$ is a rational curve in $\cF$ invariant under $\SL(2,\C)$ and disjoint from $\tilde L$, hence, Theorem \ref{thm:main-points} 
 fails for the flag $x=(q,\ell)$. The same applies to all flags of the form $(p, L)$. However, for any flag $x=(p,\ell)$ not of the above form, intersection 
 $\Th^2(x)\cap \tilde L$ is a singleton and, hence, $x$ is $\La$-generic. 
 
 Suppose that $S\subset \cF$ is an algebraic curve whose projection to $\P^2$ is not $L$ and whose projection to $(\P^2)^\vee$ is not the set of lines incident to $q$. Then, as in the previous example, $S\cap \Th^1(\tilde L)$ is finite and, hence, $S$ is $\La$-generic. The remaining type of curves consists of two rational curves, both of which are not generic. Both are $\Ga$-invariant, hence, Theorem \ref{thm:main} fails for such curves. 
 
 The analysis of hypersurfaces $S\subset \cF$ is similar to the previous example. 
  \end{example}

 \begin{example}\label{ex:E3}
 We take a discrete subgroup $\Gamma< \SU(2,1)< \SL(3,\C)$, as in Example \ref{ex:Anosov3}. Then the flag-limit set of $\Ga$ in $\cF$ is the image of $\La_1(\Ga)\subset S^3=\partial \mathbb B$ under the lift $S^3\to \Sigma\subset \cF$. (Given a limit point $\la\in \la_1(\Ga)$, take the unique projective line $\ell_\la$ tangent to $S^3$ at $\la$. Then the flag $(\la,\ell_\la)$ is the lift of $\la$ to $\Sigma$.)  The only class of (nonelementary) subgroups $\Ga$ which are not Zariski-dense (over $\bC$) in $\SL(3,\C)$ are the {\em $\C$-Fuchsian subgroups}: These are subgroups fixing a point $q\in \P^2 \setminus 
 \hat{\mathbb B}$ (the closure of $\mathbb B$ in $\P^2$), equivalently, these are the subgroups preserving the projective line $L\subset \P^2$ dual to the point $q$ as above, under the duality induced by the pseudo-hermitian form on $\C^3$ invariant under the $\SU(2,1)$-action. Then the discussion reduces to the one 
 in Example \ref{ex:E2}. 
 \end{example}


\begin{example}
Consider the Anosov representation $\rho: \Ga\to G= G_1\times ...\times G_n$ as in Example \ref{exAnosov4}. We will assume that $\Ga$ 
is nonelementary. Let $\La\subset \cF= \prod_{i=1}^n G_i/P_i$ denote the flag-limit set of $\rho(\Ga)$ as described in Example \ref{exAnosov4}. 
Note that for every $i=1,...,n$, the projection of $\La$ to $G_i/P_i$ is 1-1. 
We claim that every $x\in \cF$ is generic with respect to the limit set of $\rho(\Gamma)$ in $\cF$. 
Indeed, $x=(x_1,...,x_n)\in \cF$. For $\lambda=(\lambda_1,...,\lambda_n)$ in the limit set $\La$, $x\in \Th^{n-1}(\lambda)$ 
if and only if there exists an index $i\in \{1,...,n\}$ such that $x_i= \lambda_i$. Since the projection of the limit set of $\rho(\Ga)$ in 
$\cF$ to each factor $G_i/P_i$ is 1-1, it follows that at most $n$ maximal thickenings of limit points of $\rho(\Gamma)$ can have a 
common intersection.   

Assume that $n=2$ and $\rho_1=\rho_2$. Then every curve 
$S\subset \P^1\times \P^1$ either equals the diagonal $D=\diag(\P^1\times \P^1)$ or has finite intersection with 
$D$. The flag-limit set of $\rho(\Ga)$ is contained in $D$. Thus, $S$ is $\La$-generic unless $S=D$. In the latter case, $S=D$ is 
$\Ga$-invariant and Theorem \ref{thm:main} fails for $S$. 

We will use the notation from \S \ref{thickeningsofproductsofprojectivelines} and \ref{sec:illustration:Current-valued map}: For each factor $\P^1_i$ in the product $\P_1^1\times \P_2^1$ and points 
$z=(x,y)\in \La(\Ga)$ we have the thickenings $\Th_h(z)= \{(x', y): x'\in \P^1\}$ and $\Th_v(z)= \{(x, y'): y'\in \P^1\}$. Let $S\subset \cF$ be a curve of 
{\em bidegree} $(p,q)$, i.e. its projection to $\P^1_1$ is $p$-to-$1$ and its projection to $\P^1_2$ is $q$-to-$1$. Thus, the homology class 
$[S]\in \coh_2(\cF)$ equals $p[\P^1_1] + q [\P^1_2]$. Let $\mu=\mu_\beta$ denote the Gibbs measure on $\La$ corresponding to a hyperbolic length function $\beta$ on 
$\Ga$. For instance, if $\rho_1=\rho_2$ and $\Ga_1=\Ga_2$ is a cocompact subgroup in $G_1=G_2$, then $\La=D$ and 
we can take $\mu$ to be a multiple of the Lebesgue measure on $D$. Then the limiting Gibbs current in Theorem \ref{thm:main} (for the curve $S$) equals
$$
\int_D \left( p\llbracket \Th_v(z)\rrbracket  + q \llbracket \Th_h(z)\rrbracket  \right)d\mu(z). 
$$
\end{example}

\chapter{Equidistribution of smooth forms}\label{sec:Equidistribution of smooth forms}

Fix $m \geqslant 0 $.
In this section, we fix a smooth closed $(n-m, n-m)$-form $\psi$ on $\mathcal{F}$. For $g\in G$, we will use the notation 
$$
g_*\psi=(g^{-1})^*\psi.
$$ 

Recall that $\psi$ defines a current in $  \cD_{k,k}(\cF)$ via the formula
\[
\langle \psi, \varphi\rangle=\int_\cF \psi \wedge \varphi.
\]
For the ease of the notation, we shall view $\psi$ as a current via the above identification.

The aim of this chapter is to prove the equidistribution for smooth forms on flag manifolds. 
The most simple case is the rank one situation. We give a direct proof in this situation in Section \S~\ref{subsection_equidistrib_P1} as this serves as a roadmap for the general situation. The mass estimates and equidistribution speed arguments in the rank one situations are replaced in the general case by an extension of a result of Harvey--Lawson in \S~\ref{section_harvey}.

\section{Equidistribution of volume forms on $\bP^1$} \label{subsection_equidistrib_P1}

To illustrate the proof of Theorem \ref{thm:main smooth} we first prove this result in the setting of actions of convex-cocompact Kleinian groups on complex-projective line $\bP^1$. 

We identify $\bP^1$ with $\C \cup \{\infty\}$, 
and fix a semigroup action $\gamma_\la \colon \bP^1 \to \bP^1$,  $z \mapsto \lambda z$, with $|\lambda| > 1$, where the repelling fixed point is $x_{\omega_0}= 0$ and the attracting point is $x_{e}= \infty$.  
\smallskip

$\bullet$ \textbf{Step 1} (Mass estimates):  
%

For $\eta>0$ which will be chosen later on,
we cover $\bP^1$ as $\{ |z| < 2 \eta \} \cup \{ |z| > \eta \}$. We fix a partition of unity $(\rho_0, \rho_\infty)$ associated with this 
cover (by constructing a bump function that vanishes for all $z\in \C$ with $|z|> 2\eta$ and is identically $1$ on $\{ |z| \leqslant \eta\}$). 
This gives a decomposition:
\begin{align*}
\psi = \rho_0 \psi + \rho_\infty \psi, \\
\varphi = \rho_0 \varphi + \rho_\infty \varphi. 
\end{align*}
Recall that $\psi$ is a $(1,1)$-form, hence it can be viewed as a measure. 
For fixed $\epsilon>0$, we choose $\eta$ so that
\begin{equation*} 
Mass(\rho_0 \psi) \leqslant \epsilon.
\end{equation*}
Then for $\epsilon>0$, there exists $\eta >0$ such that:
\begin{equation} \label{eq_small_mass}
Mass((\gamma^{-1})^* (\rho_0\psi)) = Mass(\rho_0\psi) \leqslant \epsilon.
\end{equation}
This first estimate shows that for any smooth function $\varphi$, one has:
\begin{equation} \label{eq_mass_1}
\left | \langle (\gamma^{-1})^* (\rho_0\psi), \varphi \rangle \right | \leqslant |\varphi|_\infty \epsilon. 
\end{equation}

$\bullet$ \textbf{Step 2} (Support of the limit): Consider now the measures $(\gamma^{-1})^* (\rho_\infty \psi) $. Since $\rho_\infty \psi$ is supported outside of the disc of radius $\eta$ around $z=0$, the pullback measures are supported in a disc near $\infty$ contained in $\bP^1(\C) \setminus \{ |z| \leqslant |\lambda| \eta  \}$. 
 This means that any weak limit of $(\gamma_{\la_n}^{-1})^* (\rho_\infty \psi)$ with with $|\lambda_n| \rightarrow +\infty$  will 
  be supported at the point $z = \infty$. 
Moreover, if we decompose $\rho_\infty \psi = \psi_+ - \psi_-$ into a difference of positive measures, the positive measures $(\gamma_{\la_n}^{-1})^* \psi_{\pm}$ converge to a positive measures supported at $\{\infty\}$, hence they converge to 
$ Mass(\psi_\pm) \delta_{\infty} $.  
Since the masses of $\psi_\pm$ are equal to $ \langle \psi_\pm , \intcur{\bP^1} \rangle $, this shows that the difference $(\gamma_{\la_n}^{-1})^*(\psi_+ - \psi_-)$  converges to  $ \langle  \rho_\infty \psi , \intcur{\bP^1}\rangle \delta_\infty  $.   

\smallskip 

$\bullet$ \textbf{Step 3} (Speed of convergence): 
We estimate the speed of convergence to the limit in terms of the modulus of $|\lambda|$. There exists $M > 0$ such that the support of $ (\gamma_\la^{-1})^* ( \rho_\infty \psi ) $ is contained in a fixed ball $\mathbf{B}$ around the point at infinity for all $|\lambda| \geqslant M$. 
Decompose the action by $\gamma_\la$ into $\gamma_{\lambda/M}\circ \gamma_M$ where $\gamma_M$ corresponds to the multiplication by $M$ and where $\gamma_{\lambda/M}$ corresponds to the multiplication by $\lambda/M$. The form $(\gamma_M^{-1})^* (\rho_\infty \psi )$ is supported in $\mathbf{B}$, while the support $\gamma_{\lambda/M}(\mathbf{B}) \subset \mathbf{B}$ of $(\gamma_\la^{-1})^* (\rho_\infty \psi )$ converges to the point $\infty$. In view of the uniform continuity of any compactly supported test form in $\mathbf{B}$, this means that for any test form $\varphi$, there exists $M' \geqslant M$ such that for any $|\lambda| \geqslant M'$, one has: 
\begin{equation} \label{eq_dirac_estimate}
	\left | \langle   (\gamma_\la^{-1})^* (\rho_\infty \psi ) , \varphi  \rangle - \langle \rho_\infty \psi  , \intcur{\bP^1} \rangle \varphi(\infty)   \right | \leqslant \epsilon,
\end{equation}

Finally using \eqref{eq_dirac_estimate} together with \eqref{eq_mass_1}, one obtains that for any smooth test function $\varphi$ and any $\epsilon > 0$, there exists $M '> 0$ such that for any $|\lambda| \geqslant M' $, one has:
 \begin{align*}
  \left | \langle (\gamma_\la^{-1})^*\psi  , \varphi  \rangle  - \langle \psi  , \intcur{\bP^1} \rangle  \varphi(\infty ) \right |   &  	\leqslant \left | \langle (\gamma_\la^{-1})^* (\rho_0 \psi) , \varphi  \rangle  \right | \\ 
  &+ \left | \langle (\gamma_\la^{-1})^* (\rho_\infty \psi) , \varphi  \rangle   -\langle \psi  , \intcur{\bP^1} \rangle \varphi(\infty) \right | \\
  & \leqslant  |\varphi|_\infty \epsilon + \epsilon. 
\end{align*}
 
 We thus obtain that for any smooth test function $\varphi$ and any $\epsilon> 0$, there exists $M'> 0$ such that for any $|\lambda| \geqslant M'$, one has: 
 \begin{equation} \label{eq_convergence_volume_uniform}
 	\left | \langle (\gamma_\la^{-1})^*\psi  , \varphi  \rangle  - \langle \psi  , \intcur{\bP^1} \rangle  \varphi(\infty ) \right |   \leqslant \epsilon. 
 \end{equation}
 \smallskip 
 
 $\bullet$ \textbf{Step 4} (Equivariant boundary map extension and conclusion):
Given a smooth $(1,1)$-form $\psi $ on $\bP^1$, we consider the  map $\eta \colon \Gamma \cup \partial_\infty  \Gamma \to \mathcal{D}^{1,1}(\bP^1) $ with values in bidegree $(1,1)$ currents given by: 
\begin{equation}
	\eta( g) = \left \lbrace \begin{array}{ll}
		(\gamma^{-1})^* \psi & \text{ if } g\in \Gamma  \\
		 \langle \psi, \intcur{\bP^1}\rangle  \delta_{f(g)} & \text{otherwise,}
	\end{array} \right .  
\end{equation}   
 where $f \colon \partial_\infty \Gamma \to \Lambda$ is the boundary map from the boundary of the group to the limit set $\Lambda$. 
 The map $\eta$ is continuous and this is a consequence of Corollary \ref{cor_equivariant_continuous} together with \eqref{eq_convergence_volume_uniform}. As we saw in Sections \ref{sec:toy} and \ref{sec:mainproof}, continuity of $\eta$ implies the equidistribution result. 
 
\section{Algebraic $(\C^*)^n$-actions and exponential rate of convergence} \label{section_harvey}

The aim of this section is prove a quantitative version of some results due to Harvey and Lawson (see \cite{harvey_lawson}) in a more specific setup, where one considers algebraic $(\C^*)^n$ actions on smooth projective varieties.  
Their result is as follows:

\begin{thm} (see \cite[Theorem 2.3]{harvey_lawson_survey}) \label{thm_harvey_lawson} Let  $\theta_t,t\in \R$ be a smooth flow on a compact Riemannian manifold $X$ induced by a Morse function. Assume that the flow has finite volume in the sense of \cite[Definition 1.1]{harvey_lawson_survey}, that the stable and unstable manifolds for each critical point of the flow have finite volume, and that for any two points $p \neq q$ for which there exists finitely many continuous path for the positive time flow starting at $p$ and ending at $q$, the dimension of the stable manifold at $p$ has dimension strictly smaller than the dimension of the stable manifold at $q$,
then  the pull-back operators $\theta_t^*$ converge weakly to the limit operator: 
\begin{equation} 
\psi \in \Omega(X) \mapsto \lim_{t\in \bR, t\rightarrow +\infty} \theta_t^* \psi = \sum _{V} \langle S_V , \psi \rangle \intcur{U_V} \in \mathcal{D}'(X),
\end{equation}
where the sum is taken over the critical manifolds for the Morse function and where $U_V$ and $S_V$ denote the unstable and stable manifolds respectively.
\end{thm}

\begin{rem} Further results of this nature were proved with the same conclusion but with other types of flow. For example, if for any critical points $p,q$, the stable and unstable manifolds at $p$ and $q$ intersect transversely, and if the flow is linearizable at the points for which the Lyapounov exponenent is $\pm 1$. Then the same conclusion also holds. The linearizability condition can now be dropped using the result in \cite[Theorem 2]{meddane_axiom_A}.    
\end{rem}

In this section we follow closely the main strategy of Harvey--Lawson's paper. Some parts are simplified thanks to the algebraicity of the action (which implies that the flow has finite volume by \cite[Proposition 1.7]{harvey_lawson_survey}) while a difference is that we need to control the behavior of the action of a higher dimensional torus $(\C^*)^r$. 
As a first step, Harvey and Lawson show that the flow can be compactified, which is a consequence of the technical condition of  finite volume.
In our case, a compactification of the action always exists because the action is algebraic. We will only need a particular compactification: 

\begin{prop} Let $X$ be a smooth projective complex variety equipped with a free algebraic action of the multiplicative group $T=(\C^*)^r$.
Then the Zariski-closure in $\bP^r \times X\times X$ of the graph of the action $ \{(\lambda, x, \lambda\cdot  x) | \lambda\in (\C^*)^r , x\in X \}$  is an algebraic variety of complex dimension $r+ \dim_{\C} X$.   
\end{prop}

\begin{proof} 
We will view $T$ as an open subset of the affine patch of $\bP^r(\C)$; hence, the action of $T$ by multiplication on itself extends to an  algebraic action of $T$ 
on  $\bP^r(\C)$. Consider the following algebraic action of $T$ on $\bP^r(\C) \times X \times X$ defined by 
\begin{equation*}
\lambda \cdot (\mu, x, y) = (\lambda \cdot \mu , x, \lambda \cdot y).
\end{equation*} 
Set $\mathbb{I} =\{(\lambda, x, \lambda\cdot  x) | \lambda\in (\C^*)^r , x\in X \} $ and denote by $\overline{\mathbb{I}}^Z$ its Zariski closure in 
$\bP^r(\C) \times X \times X$.
We take $p = \{ \underline{1} \} \times \Delta_X$, where $\Delta_X$ is the diagonal in $X\times X$ and where $\underline{1}$ corresponds to the point $[1 : 1 : \ldots : 1] \in \bP^r(\C)$. 
The set 
$$\mathbb{I} =  \{(\lambda, x, \lambda\cdot  x) | \lambda\in (\C^*)^r , x\in X \} = T \cdot p$$ 
coincides with  the union of $T$-tranlates of $p$ and is a smooth algebraic variety of dimension $r + \dim_{\C} X$  since the map $(\lambda, x , y) \mapsto  (\lambda, x)$  from $\mathbb{I}$ to $T \times X$ is an isomorphism.  

Since $T \cdot p$ is the image of the morphism 
$ (\lambda,x) \in T\times X \mapsto (\lambda,x,\lambda\cdot x) \in \bP^r(\C) \times X \times X$, it is a constructible subset by Chevalley's Theorem (see  e.g.\ \cite[Chapter 2, \S 6.E, Theorem 6]{matsumura}), i.e., it is a finite union of subsets obtained as the intersection of a Zariski open subset with a Zariski closed set. Moreover, as $\mathbb{I}$ is irreducible, its Zariski closure $\overline{\mathbb{I}}^Z$ contains a dense Zariski-open subset { which is also contained in $\mathbb{I}$} by \cite[Lemma 1.4.5]{tauvel_yu}. 
As $T\cdot p = \mathbb{I}$ is also open and dense in $\overline{\mathbb{I}}^Z$, this  subset is exactly $T\cdot p$.  
Consider  the morphism $u \colon T \times X  \mapsto \overline{T \cdot p}^Z$ induced  by $(\lambda, x) \in T \times X \mapsto (\lambda, x, \lambda \cdot x) \in T \times X \times X$. It is a dominant morphism so by \cite[Part (ii) of Corollary 15.5.5]{tauvel_yu}, and we have: 
\begin{equation*}
\dim_{\C} \overline{ T\cdot p}^Z = \dim_{\C} (T\times X) - k,
\end{equation*}
where $k$ is the dimension of any fiber over $u(T\times X)$. Note that $k=0$ by construction, and one obtains 
$\dim_{\C} \overline{ T\cdot p}^Z = r + \dim_{\C} X$, as required. \end{proof}

We will apply the previous result in our specific setting where $X $ is the flag-manifold and where $r$ is the rank of the Lie group $G$. 
In the same spirit as \cite[Definition 2.21]{manicures}, we say  that a cone $\mathcal{C}$ in $\mathfrak{h}^+$ is \textit{uniformly regular} if  $\mathcal{C}$ is a closed convex set and the directions induced by the non-zero elements of $\mathcal{C}$ are uniformly bounded away from the walls of $\Delta$ (i.e., away from  the hyperplanes defined by the simple roots).

\begin{lem} \label{lem_basic_cone} Fix an (algebraic) isomorphism between $(\C^*)^r$ and $\exp(\mathfrak{h})$. Let $\mathcal{C}$ be a polyhedral cone in $\mathfrak{h}^+$, the following properties hold:  
\begin{enumerate}
\item[(i)] The image of $\mathcal{C}$ by the exponential map is contained in a semi-algebraic set in $(\C^*)^r \simeq \exp(\mathfrak{h})$. 
\item[(ii)] There exists a constant $A> 0$ such that for any positive simple root $\alpha_i$ and any $u \in \mathcal{C} \setminus \{0 \}$, one has $\alpha_i(u) \geqslant A || u||$. 
\item[(iii)] The closure of $\exp(\mathcal{C})$ in $\bP^r(\C)$ contains $H_\infty = \{ [x_1 : \ldots : x_r : 0] \ | \ (x_1, \ldots , x_r) \neq 0\}$. 
\item[(iv)] For any point $q \in H_\infty$, there exist a multiplicative subgroup $u_q(t)$ such that $u_q(t) \in \exp(\mathcal{C})$ for all   
$t\in \C^*$ with $|t | $ is sufficiently large. 
\item[(v)] For any point $\lambda = (\lambda_1 , \ldots, \lambda_r) \in \exp(\mathcal{C})$, the segment joining $\lambda$ and the point $[\lambda_1 : \ldots : \lambda_r : 0]$ is contained in $\overline{\exp{\mathcal{C}}} \subset \bP^r(\C)$.  
\end{enumerate}
\end{lem}

\begin{proof} We will fix a particular isomorphism, as the result does not depend on the choice of an isomorphism between $(\C^*)^r$ and $\exp(\mathfrak{h})$.
Take $\alpha_1, \ldots, \alpha_r$ the $r$ simple positive roots. 
Recall that $\mathfrak{h}^+$ is the dual cone of the cone generated by the simple roots $\alpha_i$. 
In particular, we can find a basis $(u_1, \ldots, u_r) $ in $\mathfrak{h}^+$ which is dual to $(\alpha_1, \ldots, \alpha_r)$.
We then fix an isomorphism which is the vertical map in the following diagram:
\begin{equation*}
\xymatrix{
u= \sum_{i=1}^r t_i u_i \in \mathfrak{h} \ar[r] \ar[rd]_{\exp} & (e^{t_1} = e^{\alpha_1(u)}, \ldots , e^{t_r}= e^{\alpha_r(u)}) \in (\C^*)^r \ar[d] \\
& \exp(\sum_{i=1}^r t_i u_i) \in\exp(\mathfrak{h})  
}
\end{equation*}
Let us now prove the first statement. 
As $\mathcal{C} $ is a uniformly regular polyhedral cone, there exists $C > 0$ such that for each $u \in \mathcal{C}\setminus \{ 0\}$, $C^{-1} \leqslant \alpha_i(u)/ \alpha_j(u ) \leqslant C$. Setting $\lambda_i(u) = e^{\alpha_i(u)}$, we choose a rational number $p/ q > C$ with $p$ and $q$ coprime. 
The polyhedrality condition implies that $|\lambda_i(u)|^q \leqslant |\lambda_j(u)|^p$ and $|\lambda_j(u)|^q \leqslant |\lambda_i(u)|^p$.  
As these conditions are algebraic inequalities, the image is  contained in a semi-algebraic set. 
\smallskip

For the second statement, we use the compactness of the unit sphere intersected with $\mathcal{C}$ together with the continuity of the function $\min_{i} \alpha_i(u)$ on this set. The constant $A$ can be chosen as $A = \min_{\substack{||u||=1 \\ u\in \mathcal{C}\\ i}} \alpha_i(u)$.  
\smallskip 

Observe that the third assertion follows from assertion (iv). 
As before, we fix some non-zero integers $p,q$ with $p > q$ such that $u \in \mathcal{C}$, then $\alpha_i(u)/\alpha_j(u) \leqslant p/q$ for all $i,j \leqslant r$. 
 Let us fix $x \in H_\infty$ of the form $[x_1: \ldots : x_r : 0]$ with $x_i \in \C$ not all vanishing. Let us normalize so that  $|x_i| \leqslant 1$ for all $i \leqslant r$. 
We consider 
the one-parameter subgroup $$ u_x \colon  t \in \C^* \mapsto [v_1(t): \ldots : v_r(t) : t^2] \in \bP^n(\C), $$
where 
$$v_i(t) = \left \{ \begin{array}{ll}
t^{2p-1} x_i & \text{if } x_i \neq 0, \\
t^{2q} & \text{otherwise}.
\end{array} \right . $$ 
We check that for each $t$, $u(t) = \exp h(t)$ where $\alpha_i (h(t)) = (2p-1) \log |t| + \log |x_i|$ if $x_i \neq 0$ and $\alpha_j(h(t)) =  2q \log |t|$ if $x_j = 0$. 
The quotient $\alpha_i(h(t))/\alpha_j(h(t))$ converges when $|t| \rightarrow +\infty$ to $1$, $(2p-1)/2q$ or $2q/(2p-1)$. In both cases, $ \lim_{|t|\rightarrow +\infty }\alpha_i(h(t))/ \alpha_j(h(t)) \in ]q/p, p/q[$.    
This ensures that $u_x(t)$ belongs to $\exp(\mathcal{C})$ for sufficiently large $|t|$. 

\smallskip

The last assertion follows from the fact that $\gamma(t) = (\exp(t), \ldots, \exp(t))$ belongs to $\exp(\mathcal{C})$ for all $t \in \R_+$. Since $\exp(\mathcal{C})$ is a multiplicative semigroup,  $\lambda \cdot \gamma(t)$ is in this set and its image is a segment from $\lambda$ to $[\lambda_1 : \ldots : \lambda_r: 0]$, as required.  
\end{proof}

As a consequence of the above result, we can compare the Zariski closure with the closure with the classical topology induced by the Riemannian metric in our specific situation.  
\begin{prop} \label{prop_degen_fiber} Let $\cF$ be the flag-variety
endowed  with the  algebraic action of the multiplicative group $(\C^*)^r $ which is identified with $\exp(\mathfrak{h})$. Consider the graph of the action given by: 
\begin{equation*}
\mathbb{I} = \{ (\lambda , x , \lambda \cdot x) | \lambda \in (\C^*)^r , x\in \cF \}. 
\end{equation*}
 and consider $H_\infty = \{ [q_1 : \ldots : q_n : 0]\  | \ q_i \in \C  ,\forall  i \leqslant  r\}$ the hyperplane at infinity.
Then $$p_1^{-1} (H_\infty) \cap \overline{\mathbb{I}}^Z = p_1^{-1} (H_\infty) \cap \overline{\mathbb{I}}^R = \bigcup_{w \in W} H_\infty \times \Th_{w^\vee}(\sigma_-) \times \Th_w(\sigma_+). ,$$ where $\sigma_+, \sigma_-$ are the attractive and repulsive points respectively for the action of $T^+ = \exp(\mathfrak{h}^+)$, where $p_1$ is the projection onto $\bP^r(\C)$ and where $\overline{\mathbb{I}}^Z$ and $\overline{\mathbb{I}}^R$ denote the closure in $\bP^r(\C) \times \cF\times \cF$ with respect to the Zariski topology and the  manifold topology respectively. 
Moreover, all the components in this decomposition have multiplicity one, i.e for any $q \in H_\infty$, we have the equality: 
\begin{equation*}
	(p_2)_* \intcur{p_1^{-1} (H_\infty) \cap \overline{\mathbb{I}}^Z} = \sum_{w\in W} \intcur{\Th_{w^\vee}(\sigma_-) \times \Th_w(\sigma_+)},
\end{equation*} 
in the sense of currents where $p_2$ is the projection onto $\cF\times \cF$.
\end{prop}

\begin{proof} 
Let us denote by $\overline{\mathbb{I}}^Z$ the closure with respect to the Zariski topology and by $\overline{\mathbb{I}}^R$ the closure with respect to the manifold topology. 
Since the manifold topology is finer than the Zariski topology, we have $\overline{\mathbb{I}}^R \subset \overline{\mathbb{I}}^Z$. 
By the assertion (iii) of Lemma   \ref{lem_basic_cone}, any point $q \in H_\infty$ is in the closure of $\exp(\mathfrak{h}^+)$.

\smallskip
 We claim that $p_1^{-1} (\{ q\}) \cap  \overline{\mathbb{I}}^R$ contains $\cup_{w\in W} \{ q\} \times \Th_{w^\vee}(x_-) \times \Th_w(x_+)  $. 
 By the assertion (iv) of Lemma \ref{lem_basic_cone}, we can find a one-parameter subgroup $u_q \colon \C^* \to  (\C^*)^r$ whose closure contains $q$ and such that for $|t|$ large enough, $u_q(t)$ satisfies the Anosov conditions given in assertion (ii) of the same Lemma. 
We adopt the same notation as in \S~\ref{sec:Structure of thickenings}. Let $\sigma_+$ and $\sigma_-$ be, respectively, the attractive and repulsive points of the action of $T^+$ and for each $w \in W$, we consider the point $\sigma_w$ given by $k_w B$, 
where $k_w \in K$ is a choice of an element projecting to $w$. 
By Proposition \ref{prop_hyperbolic_torus_action}, for any element $w \in W$,  there exists a chart (denoted $\pi_w$) near the point $\sigma_w  \in \cF$ on which $u_q(t)$ acts linearly diagonally with exactly $ 2n - 2|w|$ negative Lyapounov exponents and $2 |w|$ positive ones. Moreover, the third assertion  shows that the stable and unstable manifolds at $\sigma_w$ in that chart are given by $S_{w^\vee}(\sigma^-) , S_w(\sigma^+)$. A direct computation shows then that the limit of the action in that chart must contain   $S_{w^\vee}(\sigma^-) \times S_w(\sigma^+)$. By the assertion (i)  of Proposition \ref{prop_listing_thickenings}, the subvarieties $\Th_{w^\vee}(\sigma^-)$ and $\Th_w(\sigma^+)$ are the closure of $S_{w^\vee}(\sigma^-)$ and $S_{w}(\sigma^+) $ respectively, hence the following  inclusions hold: 
\begin{equation} \label{eq_inclusion}
\cup_{w \in W} \Th_{w^\vee}(\sigma^-) \times \Th_w(\sigma^+) \subset  p_1^{-1}(H_\infty) \cap \bar{\mathbb{I}}^R \subset p_1^{-1}(H_\infty) \cap  \bar{\mathbb{I}}^Z,
\end{equation}
proving the claim. 

\smallskip

We now prove the equality  $p_1^{-1} (H_\infty) \cap \overline{\mathbb{I}}^Z =  H_\infty^* \times \cup_{w \in W} \Th_{w^\vee}(\sigma^-) \times \Th_w(\sigma^+) $.
Observe now that for each $q' \in \bP^n(\C)$, the  fibers $p_1^{-1} (\{ q'\}) \cap \overline{\mathbb{I}}^Z$ are analytic subvarieties whose image $p_2(p_1^{-1} (\{ q'\}) \cap \overline{\mathbb{I}}^Z) $ by $p_2$ are also analytic subvarieties whose class in $H_{2n}(\cF \times \cF)$ are the same.

The inclusion \eqref{eq_inclusion} shows  that for any $q \in H_\infty$, the difference of currents: 
\begin{equation*} 
R =	\intcur{p_2 (p_1^{-1}(\{q\}) \cap \bar{\mathbb{I}}^Z)} - \sum_{w\in W} \intcur{\Th_{w^\vee}(\sigma^-) \times \Th_w(\sigma^+)}
\end{equation*}
is a closed positive current  in the sense of \S~\ref{section_positive_currents}.
In particular, by Lemma \ref{lem_mass_estimate}, the mass $|| R||$ of $R$ is non-negative and controlled: 
$ 0 \leqslant 	||R|| \leqslant C \langle R , \omega^{2n} \rangle,  
$ where $C$ is a uniform constant. 

Since the fiber over the point $\underline{1}= [1: 1 : \ldots :1]$ is the diagonal $\{ [1: \ldots : 1] \} \times \Delta_{\cF}$, we deduce that for any $q\in H_\infty^*$, one has 
$	\{ p_2 (p_1^{-1}(\{q\}) \cap \bar{\mathbb{I}}^Z) \} = \{\Delta_{\cF} \} \in H_{2n}(\cF\times \cF). 
$
Moreover, the equality of classes
$
	\sum_{w\in W} \{  Th_{w^\vee} (\sigma^-) \times \Th_w(\sigma^+)\} = \{ \Delta_{\cF} \} 
$
in $H_{2n}(\cF, \C)$ shows that $\{ R\} = 0$ in $H_{2n}(\cF\times \cF)$. As a result, 
we have $\langle R, \omega^{2n} \rangle = \{ R\} \cup \{ \omega^{2n}\} =0 $, so $R$ is a positive current whose mass is zero, hence it is equal to zero, as required.  
\end{proof}

\begin{thm} \label{thm_extension_harvey_lawson} Consider a flag-manifold $\cF$ of dimension $n$ and let $\mathcal{G}$ be a semi-group defined by the image of a uniformly regular polyhedral cone in $\mathfrak{h}^+$ under the exponential map, with the attracting and repelling points $\sigma^+, \sigma^-$ respectively. 
We denote by  $f_{\lambda}$ the action of $ \lambda\in \mathcal{G}$ on $\cF$. There exists
some constants $A,C>0$ such that for any smooth  $(n-k,n-k)$-form $\psi$, for any smooth $(k,k)$ form $ \varphi$ on $\cF$ and
 for any (normalized) vector  $ u\in \mathfrak{h}^+$ with $\exp(u) \in \mathcal{G}$, one has: 
\begin{equation}\label{eq_diff_smooth}
\left| \langle f_{\exp({-tu})}^* \psi - \sum_{w\in W} a_w(\psi) \intcur{\Th_w(\sigma^+) } , \varphi \rangle \right | \leqslant C e^{- At} \max (||\psi ||, || d\psi || , ||\varphi ||, || d\varphi||),
\end{equation}    
where $a_w (\psi) = \langle \intcur{\Th_{w^\vee}(\sigma^-)} ,\psi \rangle$ for all $w\in W$.
\end{thm}

\begin{rem} 
The above theorem implies that the family of pullback currents converges exponentially fast to the limit. 
\end{rem}

\begin{proof}
We first rewrite the formula.
Set $\lambda = \exp(t u)$ for $u \in \mathfrak{h}^+$ and $t > 0$ with $\exp(u) \in \mathcal{G}$. 
Using the projections $\pi_1,\pi_2$ from $\cF \times \cF$ onto the first and second component, one first rewrites the left-hand side of \eqref{eq_diff_smooth} as: 
\begin{equation} \label{eq_difference_graphs}
\langle f_{\lambda^{-1}}^* \psi - \sum_{w\in W} a_w(\psi)\intcur{\Th_w(\sigma^+)} , \varphi \rangle = \langle \intcur{\mathcal{T}_\lambda}  - \sum_{w\in W}  \intcur{ \Th_{w^\vee}(\sigma^-) \times \Th_w(\sigma^+)} , \pi_1^* \psi \wedge \pi_2^* \varphi  \rangle ,
\end{equation}
where $\mathcal{T}_\lambda = \{(x, \lambda \cdot x) | x\in \cF  \}$ is the graph of $f_\lambda$ in $\cF\times \cF$. 
The previous equality shows that it is enough to show that the difference of current converges to zero uniformly in $\lambda$. 
To obtain a finiteness statement, we use the previous compactification of the action which allows one to interpolate between the two currents. 
Following the notation of Proposition \ref{prop_degen_fiber}, we consider the subset $ \bar{\mathbb{I}}^Z$ defined in that statement.
By Proposition \ref{prop_degen_fiber},  $p_1^{-1}(H_\infty) \cap \bar{\mathbb{I}}^Z$ is an complex-analytic subvariety in $\bP^n(\C)\times \cF \times \cF$ of complex dimension  $2n-1$. 
 Note also that the closure on the hyperplane at infinity $H_\infty \simeq \bP^{n-1}(\C)$ in $\bP^n(\C)$ of the semi-group $\mathcal{G}$ consists of   a convex set of $H_\infty$.  

Let $p_1, p_2$ be the projection of $\bP^n(\C) \times \cF \times \cF$ onto $\bP^n(\C)$ and $\cF \times \cF$ respectively. 
Given $\lambda= (\lambda_1, \ldots, \lambda_n)  \in (\C^*)^n$, we take $q = [\lambda_1 : \ldots : \lambda_n :0] \in H_\infty$ and set $R = \max_{i\leqslant n}(|\lambda_i|)$

By Proposition \ref{prop_degen_fiber}, since $q\in H_\infty$, the subset $\mathcal{T}_q =\bar{\mathbb{I}}^Z \cap p_1^{-1}(\{q \})$ is an algebraic variety of real dimension $2n$ contained in $\bar{\mathbb{I}}^Z$. 

 We next show that 
 the fibers $\intcur{p_1^{-1}(\{\lambda\}) \cap \bar{\mathbb{I}}^Z }$ and $\intcur{\{q\}\times \mathcal{T}_q}$ viewed as currents are close for the flat norm uniformly in $R$.  To that end, we let $\sigma_q$ be the line segment connecting $\lambda$ and $q$. By the assertion (v) of Lemma \ref{lem_basic_cone}, this segment is contained in the closure of $\mathcal{G}$ in $\bP^r(\C)$. 
Now let us consider $
	\mathcal{J}_{q} = \bar{\mathbb{I}}^Z \cap p_1^{-1} (\sigma_q).
$
  It is a semi-algebraic set of real dimension $2n +1$. 
  We claim that its volume is controlled so that there exists a constant $C>0$ which is independent of $\lambda$ such that: 
  \begin{equation} \label{eq_control_volume_slice}
  	\vol_{2n+1}(\mathcal{J}_{q,R}) \leqslant \dfrac{C}{R}.   
  \end{equation}
  Indeed, the fiber over the segment $\sigma_q$ by $p_1$ is a real smooth manifold of dimension $2n+1$ while the fibers over the points $\lambda$ and $q$ are projective varieties (non irreducible or smooth) of real dimension $2n$. By the smooth coarea formula applied to the map $v = (p_1)_{| \mathcal{J}_{q,R}}\colon \mathcal{J}_{q,R} \to \bar \sigma \subset \bP^n(\C)$, we have: 
  \begin{equation*}
  \int_{\mathcal{J}_{q,R}} J_1(p_1) d\vol_{2n+1} = \int_{ \sigma_q} \vol_{2n}(v^{-1}(\{ y\})) d\vol_1(y), 
  \end{equation*}
  where $J_1(p_1)$ is the maximal minor of the differential of $p_1$ with respect to the Riemannian metric. 
   Note the fiber $v^{-1}(y)$ over all points $y \in \sigma$ is a complex projective variety in $\cF\times \cF$, in particular they are closed and their volume with respect to the K\"ahler metric $\omega$ is controlled by the cohomology using \cite[Equation 1.24]{demailly}: 
   \begin{equation*}
   \vol_{2n} ( v^{-1}(\{ y\})) = \dfrac{1}{2^{n} n!} \langle \intcur{v^{-1}(\{ y \})} , \omega^{n} \rangle = \dfrac{1}{2^{n} n!} \{ \Delta \} \cup \{ \omega^{n} \} = D \in \R \sim H^{4n}(\cF\times \cF, \R)   
   \end{equation*}
   where $D> 0$ and where we used the fact that $\{v^{-1}(\{ y\}) \}$ and the diagonal  $\{\Delta \}$ in $\cF\times \cF$ have the same in homology classes.    The previous two equalities then give: 
   $$
  D \vol_1(\sigma_q)  =   \int_{\mathcal{J}_{q,R}} |J_1(p_1)| d\vol_{2n+1}  \geqslant \inf_{\mathbb{I}_{>R}} |J_1(p_1)| \vol_{2n+1}(\mathcal{J}_{q,R}). 
$$
    Observe that the restriction of $p_1$ to  $T \times \{ x_0\}\times \{ x_0\}$ is the identity on $T$, in particular $|J_1(p_1)| \geqslant C_1>0$ at all points of $\bar{\mathbb{I}}^Z$. 
    We thus obtain for all $q \in H_{\infty} \cap \overline{\mathcal{G}}$ and all $R> 0$: 
    \begin{align*}
    \vol_{2n+1} (\mathcal{J}_{q,R}) \leqslant \dfrac{D \vol_1(\sigma_q)}{C_1}.
\end{align*}     
A straighforward computation of the length of $\sigma_q$ for the Fubini-Study metric gives: 
\begin{equation*}
\vol_1(\sigma_q) = \int_{[0,\max |q_i|/R]} \dfrac{dr}{\sum_{i=1}^n |q_i|^2 + r^2} \leqslant \int_{[0,\max |q_i|/R]} \dfrac{dr}{\sum_i |q_i|^2} \leqslant \dfrac{1}{R}.
\end{equation*}
The previous two inequalities yields qquation \eqref{eq_control_volume_slice} where $C = D/C_1$ is a constant independent of $\lambda$, as required. 

\medskip

  Consider the currents of integration $\intcur{\mathcal{J}_{q,R}}$; these currents exist by \cite[\S~3 p.11]{hartd_cie} and we claim that: 
  \begin{equation} \label{eq_claim}
  	d \intcur{\mathcal{J}_{q,R}} = \intcur{\mathbb{I}_{\geqslant  R} \cap p_1^{-1}(\{\lambda\})} - \intcur{\{q \}\times \mathcal{T}_q}.
  \end{equation} 
   Indeed, consider the restriction of the projection $p_1$ to $\mathcal{J}_{q,R}$, by construction it is a semi-algebraic map and the image $p_1(\mathcal{J}_{q,R})$ is the segment $\sigma_q $.
   By \cite[Theorem 3.5]{hartd_cie}, we have 
   \begin{equation} \label{eq_projection}
   (p_1)_* d  \intcur{\mathcal{J}_{q,R}} = d (p_1)_*\intcur{\mathcal{J}_{q,R}} = d \intcur{\sigma_q} = \intcur{\lambda} - \intcur{q}.
   \end{equation}
   Observe that  $p_1^{-1}(\sigma_q) \cap \mathcal{J}_{q,R}$ is a semi-algebraic set of real dimension $1 + 2n $ while the fibers over $\lambda,q$ are of dimension $2n$. 
   Since the boundary  $\partial \mathcal{J}_{q,R}$ is a semi-algebraic set of dimension $2n$ contained in the fibers over $\lambda$ and $q$, we deduce that the boundary consists of exactly the two fibers over $\lambda$ and $q$. 
   We thus deduce that $d \intcur{\mathcal{J}_{q,R}} = \alpha \intcur{p_1^{-1}(\{\lambda \})\cap \mathcal{J}_{q,R}} + \beta\intcur{p_1^{-1}(\{q \})\cap \mathcal{J}_{q,R}} $. 
   We then deduce that $\alpha =1$ and $\beta =-1$ from the projection formula \eqref{eq_projection} and  this proves that \eqref{eq_claim} holds.  
By \eqref{eq_control_volume_slice} , the uniform convergence of $\vol_{2n+1} (\mathcal{J}_{q,R})$ to zero shows that the boundary converges to zero uniformly in $q$ for the flat norm, more precisely we have that for any $\lambda \in \mathcal{G}$, 
\begin{equation*} 
N_F \left  (\intcur{\{\lambda\}\times \mathcal{T}_\lambda} - \sum_{w\in W} \intcur{ \{q \}\times \Th_{w^\vee}(\sigma^-) \times \Th_w(\sigma^+)} \right ) \leqslant \dfrac{C }{\max |\lambda_i|}.
\end{equation*}
By Lemma \ref{prop_pushforward_flat} applied to the map $p_2$, there exists a constant $C'>0$ such that for any current $T$ on $\bP^n(\C) \times \cF \times \cF$, $N_F({p_2}_* T ) \leqslant C' N_F(T)$, which gives: 
\begin{equation*}
	N_F(\intcur{\mathcal{T}_\lambda} - \sum_{w\in W} \intcur{\Th_{w^\vee(\sigma^-)}) \times \Th_{w}(\sigma^+)}) \leqslant \dfrac{C C'}{\max |\lambda_i|}
\end{equation*}
Finally, since  the action on $\cF$ of $u\in \mathfrak{a}_+$ is determined by $\lambda_i = \exp( \alpha_i(u))$ where $\alpha_i$ are the positive simple roots in $\mathfrak{a}^*$,  we find an upper bound when  
 $\lambda(t)= (\exp(\alpha_i(tu)))$ with $ \exp(u) \in \mathcal{G}$. 
 By the assertion (ii) of Lemma \ref{lem_basic_cone}, there exists a constant $A> 0 $ such that $ \alpha_i(u)/ || u ||> A$ for all $u \in \mathfrak{h}$ with $\exp(u) \in \mathcal{G}$. The previous formula then gives: 
 \begin{equation} \label{eq_exp_flatnorm}
 N_F \left ( \intcur{\mathcal{T}_{\exp(tu)}} - \sum_{w\in W} \intcur{\Th_{w^\vee}(\sigma^-) \times \Th_w(\sigma^+)} \right ) \leqslant CC' e^{-A t}.
\end{equation}    
One then concludes the proof using the previous formula together with \eqref{eq_apply_flat_norm} and the identity \eqref{eq_difference_graphs}. 
  \end{proof}

\section{Equidistribution of smooth forms} \label{section_equivariant_smooth}

Let $r\leqslant n $ and fix a smooth closed $(n-r,n-r)$-form $\psi$ on $\cF$ and an Anosov group $\Gamma$ acting on $\cF$ with its equivariant homeomorphism 
to the flag-limit set $f \colon \partial_\infty \Gamma \to \Lambda \subset \cF$. 
We set: 
\begin{equation}
a_w =  \langle \psi, \intcur{\Th_{w^\vee}(x_-)} \rangle
\end{equation}
for all $w\in W$ where we use the convention that $a_w=0$ when the dimensions do not match.  
We define the following equivariant map  
$$\eta_\psi \colon \Gamma \cup \partial_\infty \Gamma \to \mathcal{D}_{(k,k)}(\cF)$$ 
\begin{equation} \label{eq_eta_smooth}
\eta_\psi(\gamma) = \left  \lbrace \begin{array}{ll}
(\gamma^{-1})^* \psi & \text{if } \ga \in \Gamma \\
\sum_{w \in W}  a_w \llbracket \Th_w(f(\gamma))\rrbracket  & \text{otherwise.}
\end{array} \right .  
\end{equation}

The main result of this section is the following theorem:

\begin{thm} \label{thm_continuity_equivariant_smooth} Fix $r \leqslant n$. Let $\Gamma$ be a (non-elementary) Anosov subgroup of $G$ acting on the full flag manifold $\cF$.  For any smooth closed  $(n-r,n-r)$-form $\psi$ on $\cF$, the map $\eta_\psi \colon \Gamma \cup \partial_\infty \Gamma \to \mathcal{D}_{(r,r)}(\cF)$ is continuous.  
\end{thm}

\begin{proof} 
Since $\Gamma$ is Anosov , it is uniformly regular by \cite[Theorem 3.41]{manicures}. This means that except for finitely many elements of the group, directions of the Cartan projections of the elements of the $\Ga$ are at uniformly positive distance away from the walls. 
In particular, the Cartan projection is contained in a uniformly regular polyhedral cone in $\mathfrak{h}^+$.  
Using Theorem \ref{thm_extension_harvey_lawson}, the hypothesis of 
 Corollary \ref{cor_equivariant_continuous} are satisfied. This shows that the map $\eta_\psi$ is continuous, as required.
 \end{proof}

As a consequence of the previous theorem, we obtain an equidistribution result for smooth forms as follows.

 Recall from Section \ref{section_Gibbs} and Definition \ref{def:Anosov Gibbs} 
 that for an Anosov subgroup $\Gamma< G$ with the boundary map $f: \geo \Ga\to \La=\La(\Ga)\subset \cF$ 
 and a hyperbolic length function $\beta \colon \Gamma \to \R_+$, we defijned the Poincar\'e series  $\cP(s) = \sum_{\gamma \in \Gamma} e^{-s \beta(\gamma)}$ and defined the Gibbs-measure $\mu=\mu_{\cF}$ on $\La$  as the pushforward of the Gibbs measure on $\geo \Ga$ via  $f$.  

\begin{thm} \label{thm_equidistribution_smooth} Fix $r \leqslant n$ and a smooth closed $(n-r,n-r)$-form $\psi$ on $\mathcal{F}$. Let $\Gamma$ be an Anosov group acting on  $\mathcal{F}$; fix a hyperbolic length function $\beta$ on $\Ga$ 
and  a  PS sequence $(s_j)$ converging to the critical exponent $s_0$ for its associated Poincar\'e series $\mathcal{P}(s)$. 
Then the following convergence holds in the sense of weak convergence of currents: 
\begin{equation}
\lim_{j\rightarrow +\infty} \dfrac{1}{\cP(s_j)} \sum_{\gamma \in \Gamma} e^{-s \beta(\gamma)} (\gamma^{-1})^* \psi = 
\sum_{w\in W} a_w \int_{\La} \llbracket \Th_w (\la)\rrbracket d\mu_{\cF}(\la),
\end{equation}
 where $W$ denotes the Weyl group and where $a_w = \langle \psi, \intcur{\Th_{w^\vee}(x_-)} \rangle $ for all $w\in W$. 
\end{thm}

\begin{proof} By Theorem \ref{thm_continuity_equivariant_smooth}, the map $\eta_\psi \colon \partial_\infty \Gamma \to \mathcal{D}_{(k,k)}(\cF)$ given in \eqref{eq_eta_smooth} is continuous.  Applying Proposition \ref{prop:poincare_series_current} to $\eta_\psi$, we obtain: 
\begin{equation*}
\lim_{j\rightarrow +\infty} \dfrac{1}{\cP(s_j)} \sum_{\gamma \in \Gamma} e^{-s_j \beta(\gamma)} \eta_\psi(\gamma) = \int_{\partial_\infty \Gamma} \eta_\psi(p) d\mu_{\Ga}(p).
\end{equation*}
We conclude that the result holds using the definition $\eta_\psi$ from \eqref{eq_eta_smooth}. 
\end{proof}

Note that the previous theorem is stated only for bidegree $(n-r,n-r)$-forms, the reason is that for arbitrary $(n-r, n-r')$-forms with $r \neq r'$, the currents converge to zero, as we will see below. 

Given $r \neq r' \leqslant n$ and a smooth $(n-r,n-r')$-form $\psi$, we also consider  the equivariant map $\eta_{\psi} \colon \Gamma \cup \partial_\infty \Gamma $: 
\begin{equation*}
\eta_{\psi} (\gamma) = \left  \lbrace \begin{array}{ll}
(\gamma^{-1})^* \psi & \text{ if } \gamma \in \Gamma, \\
 0  & \text{ otherwise.}
\end{array} \right .
\end{equation*}

\begin{cor} For any $r\neq r' \leqslant n$ and any smooth $(n-r,n-r')$-form $\psi$, the map $\eta_{\psi}$ is continuous.  
\end{cor}  
  
\begin{proof}
We first prove an analog of Theorem \ref{thm_extension_harvey_lawson} and adopt the same notation.  We show that under an exponential flow, the currents $f^*_{\exp(-t u)} \psi$ converge to zero where $u \in \mathfrak{h}^+$. Setting $\lambda = \exp(t u)$.
Fix a smooth $(r,r')$-form $\varphi$, we equivalently prove that  $\langle \intcur{\mathcal{T}_{\lambda}} , \pi_1^* \psi \wedge \pi_2^* \varphi \rangle  $ converges to zero as $t$ tends to infinity. 
By \eqref{eq_exp_flatnorm}, we have that: 
\begin{multline*}
\lim_{t\rightarrow +\infty} \langle \intcur{\mathcal{T}_{\lambda}} , \pi_1^* \psi \wedge \pi_2^* \varphi \rangle = \sum_{w\in W} \langle  \intcur{\Th_{w^\vee}(x_-)\times \Th_w(x_+)}, \pi_1^* \psi \wedge \pi_2^*\varphi \rangle \\= \sum_{w \in W} \langle \intcur{\Th_{w^{\vee}}(x_-)}, \psi \rangle \langle \intcur{\Th_{w}(x_+)}, \varphi \rangle =0, 
\end{multline*}
where we have used that $r \neq r'$ and $n-r \neq n-r'$ so one of the pairing is always zero. 
By Corollary \ref{cor_equivariant_continuous}, we conclude that $\eta_{\psi}$ is continuous. 
\end{proof}  
  
 As an immmediate consequence, we obtain the following result. 
\begin{cor} Fix $r \neq  r'\leqslant n$ and a smooth $(n-r,n-r')$-form $\psi$ on $\cF$. Let $\Gamma$ be an Anosov group acting on $\cF$; fix a hyperbolic length function $\beta$ on $\Gamma$ and a $PS$ sequence $(s_j)$ converging to the critical exponent $s_0$ for its associated Poincar\'e series $\cP(s)$. Then the following convergence holds in the sense of weak convergence of currents: 
\begin{equation*}
\lim_{j\rightarrow +\infty} \dfrac{1}{\cP(s_j)} \sum_{\gamma \in \Gamma} e^{-s \beta(\gamma)} (\gamma^{-1})^* \psi = 0
\end{equation*} .
\end{cor} 
   
   \begin{proof}
   By the previous statement, the map $\eta_{\psi} $ is continuous, so by Proposition \ref{prop:poincare_series_current}, we have
  \begin{equation*}
  \lim_{j\rightarrow +\infty} \dfrac{1}{\cP(s_j)} \sum_{\gamma \in \Ga} e^{-s_j \beta(\gamma)} \eta_{\psi}(\gamma) = \int_{\partial_{\infty} \Gamma} \eta_{\psi}(p) d\mu_{\Gamma}(p) = 0 ,
  \end{equation*}
  as required.
   \end{proof}

\chapter{Intersection of Gibbs currents}
\label{sec:intersection_Gibbs}
In  \S~\ref{sec: Anosov Gibbs}, given an Anosov subgroup $\Gamma$, we have defined some currents, called Gibbs currents whose support is invariant by $\Gamma$.  One can "loosely" view these as some analytical objects whose support is some union  higher dimensional  
 submanifolds. Our goal in this section is to produce a measure out of this construction. The main point is that if we are given $n$ codimension $1$ hypersurfaces in the flag manifold $\cF$, then their wedge product, when it is well-defined gives the Dirac mass measure at the intersection. To do that, we will use the wedge product method due to Bedford-Taylor presented in \S~\ref{section_bedford}. 
%


\section{Wedge product  of Gibbs currents}

Recall from Section \ref{section_bedford} that one defines wedge product of a closed positive bidegree $(1,1)$-current (which has bounded local potential) 
with any other closed positive bidegree $(p,p)$ current. 


%
%
%
%

We now describe the specific situation where one averages some current of integration along some thickening on the flag manifold. 
First, we give some conditions for which  the wedge product of a current with any other closed positive current is well-defined. 
To that end, recall that $d_{\cF}(p,V)$ denotes the minimal distance (with respect to the 
the Riemannian distance function on $\cF$) between a point $p\in \cF$ and a subvariety $V \subset \cF$ (cf. \S \ref{sec:metric_geometry}).

\begin{thm} \label{thm_wedgeable} Fix $w\in W$ with $|w|=n -1$.  Let $\mu $ be a probability measure on $\cF$ with no atoms. 
Assume that for any point $p\in \cF$, there exists a constant $M> 0$ and a neighborhood  $N$ of $p$   such that for any $q \in N$ 
\begin{equation} \label{eq_cond_integrability}
\int_{\cF} \log \left (  d_{\cF}(q, \Th_w(\lambda))\right ) d\mu(\lambda) > -M.
\end{equation} 
Then the wedge product of the bidegree $(1,1)$ current 
$$
  T = \int_{\cF} \llbracket \Th_w(\lambda) \rrbracket d\mu(\lambda) $$
  with any other closed positive current is well-defined and the current $T$ has bounded local potential.  

\end{thm}

\begin{proof} Let us show that the current $T$  admits a global potential.
Fix $w\in W$ of length $n-1$ and $\lambda_0 \in \Lambda$, by \cite[Proposition 1.4.1]{brion_lectures}, each thickening $\Th_w(\lambda_0)$ can be realized as the pullback of a hyperplane section under a projection $\pi \colon \cF \to Gr(d,n)$ onto a Grassmannian. In particular, this means that there exists a holomorphic section $\sigma $ on the hyperplane bundle on $Gr(d,n)$ such that $\Th_w(\lambda_0) $ is the zero locus of $\sigma \circ \pi$. 
This implies that the function $u=\log ||\sigma\circ \pi||$ is plurisubharmonic and by the Poincar\'e--Lelong equation, 
\begin{equation}
 i \omega +i dd^c  u = 2 \pi \intcur{\Th_w(\lambda_0)} ,
\end{equation} 
where $\omega $ is the smooth $(1,1)$ form representing the first Chern class $c_1(\Th_w(\lambda_0))$. 
Now for each $\lambda \in \cF$, by homogeneity of the action, there exists $g_{\lambda} \in G$ such that $g_{\lambda}\cdot \lambda = \lambda_0$.  Hence the function $u\circ g_\lambda$ is a global potential of $\Th_w(\lambda')$. 
We thus set: 
\begin{equation} \label{eq_potential}
\mathcal{U}_T=\int_{\cF} u \circ g_{\lambda} d\mu(\lambda),
\end{equation}
so that one has
\begin{equation*}
\int_{\cF} g_\lambda^* \omega  d\mu(\lambda) + i dd^c \mathcal{U}_T = T.
\end{equation*}  
\medskip

Let us now show that the potential $\mathcal{U}_T$ is bounded. To that end, fix two opposite points $p, p_- \in \cF$. 
 Using a partition of unity, we can reduce to the case where $\mu$ is supported in a small polydisc in some affine chart in $\Opp(x_-)$ for some point $x_- \in \cF$. We thus assume that $\supp \mu \subset \Opp(x_-)$. 
We shall rewrite \eqref{eq_potential} more explicitly. 
Fix a base-point $x_0$ in the support of $\mu$ and take a global potential $u_0$ of $\Th_w(x_0)$ (as the logarithm of the norm of some section). 
Recall that $\Opp(x_-) = U^- x_0$ where $U^-$ is the unipotent radical of the $G$-stabilizer of $x_-$. Consider the biholomorphism 
$$\phi \colon g\in U^- \to g \cdot x_0 \in \Opp(x_-)$$
and let $\mu' = \phi^{-1}_* \mu$. We rewrite \eqref{eq_potential}, using the relation $g^{-1} \cdot \lambda = x_0$, as: 
 \begin{equation*}
 \mathcal{U}_T = \int_{U^-} u_0 \circ g^{-1} d\mu'(g).
 \end{equation*}
Our aim is to show that $\mathcal{U}_T$ is bounded near $p$.

By the construction, $u_0= \log ||\sigma_0||$ is the logarithm of the Hermitian norm of a holomorphic section $\sigma_0$ of a line bundle $\pi \colon( E, ||\cdot||) \to \cF$.  
Choose a polydisc $V\subset \cF$ near $p$ over which the bundle can be trivialized so that the section $\sigma_0$ is of the form $x \mapsto f_0(x)$, where $f_0$ is a local analytic function on $V$. 
This means that we can write for all $x \in V$:
\begin{equation} \label{eq_herm_norm}
|| \sigma_0(x)|| = e^{-\varphi_0(z)} |f_0(x)|,
\end{equation}
where $\varphi_0 \colon V\to \mathbb{R}$ is a smooth function.

By the Lojasiewisz inequality \cite{Lojasiewicz}, there exist $C>0$ and $\alpha > 0$, such that: 
\begin{equation*} 
d_{\cF}(x, \Th_w(x_0))^\alpha \leqslant C |f_0(x)|, 
\end{equation*}  
for any $x\in V $.

Note that since the support of $\mu$ is a small polydisc near $x_0$, the support of $\mu'$ is then a compact set containing the identity. Choose a smaller polydisc $\tilde{V} \subset V$ such that for all $g\in \supp \mu'$ and all $x\in \tilde{V}$, $g^{-1} \cdot x \in V$. 

Since the action by $g\in\supp \mu'$ is continuous, there exists a constant $C_2>0$ such that
for all $g \in  \supp \mu'$ and all $x\in \tilde{V}$, one has:  
 $$d_{\cF}(x, g \Th_w(x_0)) \leqslant C_2 d_{\cF}(g^{-1}(x), \Th_w(x_0)) $$ for some constant $C_2> 0$.

The Lojasiewisz inequality then implies that for any $g \in  \supp \mu'$ and any point $x \in \tilde{V}$: 
  \begin{equation*}
  d_{\cF}(x, \Th_w(g x_0))^\alpha \leqslant C_2^\alpha C |f_0 \circ g^{-1} (x)|.  
\end{equation*}
In particular, using \eqref{eq_herm_norm}, one obtains for all $g\in \supp \mu'$ and all $x\in \tilde{V}$:
\begin{align*}
u_0 \circ g^{-1}(x) &= \log ||\sigma_0 \circ g^{-1}(x)|| \\
 & = - \varphi_0 \circ g^{-1} (x) + \log |f_0 \circ g^{-1}(x)| \\
  & \geqslant 
\alpha \log d_{\cF}(x, \Th_w(g x_0)) - \varphi_0(g^{-1} x) -  \log(C_2^\alpha C) .
\end{align*}
  
Summing over  all $ g\in  \supp \mu'$ yields for all $x\in \tilde{V}$: 
\begin{multline*}
\mathcal{U}_T(x) = \int u_0 \circ g^{-1}(x) d\mu'(g)  \geqslant \alpha \int{ \log} d_{\cF} (x, \Th_w(gx_0)) d\mu'(g) - \sup_{V} \varphi_0 M' - \log (C_2^\alpha C), 
\end{multline*}
where $M'$ is the total mass of $\mu'$.
Using \eqref{eq_cond_integrability}, this shows that 
$\mathcal{U}_T(x) \geqslant -\alpha M  - A > -\infty$  for all $x\in \tilde{V}$, 
where $A= \sup_V \varphi_0 M' + \log (C_2^\alpha C) $. 
We have thus shown that the potential  $\mathcal{U}_T$ is bounded near $p$, as required.
\end{proof}

In order to compute the intersection of Gibbs currents, we will use the properties of the algebraic intersections of specific algebraic varieties in the flag manifold $\cF = G/B$.

\begin{thm} \label{thm_intersection_laminar_flag} Let $\Gamma$ be an Anosov group acting on $\cF$ with length function $\beta$. Let $w, w_1, \ldots$, $w_l \in W $ with $|w_i| = l-1$ for all $i \leqslant l$ and consider $(l+1)$ Gibbs currents $T_{w}, T_{w_1}, \ldots$, $T_{w_l}$ for a Gibbs-measure $\mu$ on $\La(\Ga)\subset \cF$ and assume that the integrability condition \eqref{eq_cond_integrability} holds for each $w_i$ with $i\leqslant l$. Then 
the wedge product $T_{w_1} \wedge \ldots \wedge T_{w_l} \wedge T_{w} $ (in the sense of \S~\ref{section_bedford}) exists, 
 and
\begin{equation} \label{eq_geom_intersection}
T_{w_1} \wedge \ldots \wedge T_{w_l} \wedge T_w = \int_{\cF^{p+1}} \sum_{C} m_C \intcur{C}  d\mu^{\otimes (p+1)}(\lambda, \lambda_1 , \ldots, \lambda_l),
\end{equation}
where the sum is taken over all  irreducible components $C$ of the intersection 
$$
\Th_{w}(\lambda)\cap \bigcap_{i\leqslant l} \Th_{w_i}(\lambda_i)$$
 and where $m_C$ is the multiplicity of the intersection. 
\end{thm}

\begin{proof}
By Theorem \ref{thm_wedgeable}, the integrability condition for each $w_i \in W$ implies that $T_{w_1} \wedge \ldots \wedge T_{w_l} \wedge T_w$ is well-defined and each of the current $T_{w_i}$ has bounded potential. 
Recall from \S~\ref{section_bedford} that we have established when the wedge product of $p+1$ currents can be expressed as an integral, when the currents are of bidegree $(1,1)$, are uniformly woven with bounded potentials.  
By Corollary \ref{cor_intersection_woven}, one has: 
\begin{equation}
T_{w_1} \wedge \ldots \wedge T_{w_l} \wedge T_w =  \int_{\cF^{l+1}} \intcur{\Th_{w}(\lambda)} \wedge \intcur{\Th_{w_1}(\lambda_1)} \wedge 
\ldots \wedge \intcur{\Th_{w_l}(\lambda_l)}  d\mu^{\otimes (l+1)}(\lambda, \lambda_1, \ldots, \lambda_l). 
\end{equation}
Since the sum is taken over all the components which intersect properly, we obtain the required formula using Proposition \ref{prop_demailly_algebraic}.
\end{proof}

We now provide a setting for which  wedge products of a Gibbs-current with other currents are well-defined. 

\begin{thm} \label{thm_bounded_potential_curve} Let $\Gamma$ be an Anosov group acting on the flag-manifold $\cF$ with hyperbolic 
length function $\beta$ and with associated Gibbs measure $\mu_{\cF}$. Assume that the limit set of $\Ga$ is contained in a $G$-invariant algebraic curve in $\cF$ which intersects properly any codimension $1$ thickening $\Th_w(p)$ in $\cF$. Then for any $w\in W$ with $|w|= n-1$, the Gibbs-current
\begin{equation*}
T_w = \int \intcur{\Th_w(\lambda)} d\mu_{\cF}(\lambda)
\end{equation*}
has bounded local potential, hence the wedge product of $T_w$ with any other closed positive current is well-defined. 
Moreover, \eqref{eq_geom_intersection} is satisfied for these currents.
\end{thm}

\begin{proof} Denote by $C$ the invariant algebraic curve containing the limit set $\Lambda$ and let us fix $w \in W$ with $|w|=n-1$. By Theorem \ref{thm_wedgeable}, we just need to check the integrability condition \eqref{eq_cond_integrability}.

Take any point $ p \in \cF$. Observe that $ p \in \Th_w(\lambda) \Leftrightarrow \lambda \in \Th_w(p)$. Moreover, by Theorem \ref{thm:thickening_distance_inequality}, locally near any point $\lambda \in \cF$, there exists a constant $M>0$ such that: 
\begin{equation} \label{eq_distance_thick}
 d_{\cF}(p, \Th_w(\lambda)) \geqslant M d_{\cF}(\Th_{w^\vee}(p), \lambda).  
\end{equation}
Consider the intersection points between $\Th_{w^\vee}(p)$ and $C$. Since there are finitely many of those points, we can just  prove the statement locally near each of these intersection points. 
 Let $q$ be an intersection point of $\Th_{w^\vee}(p) $ and $C$. 
This means that the position $\pos(p,q) \leqslant w^\vee$. Set $w' = \pos(p,q)$, we have $q \in S_{w'}(p)$. By definition, there exists a point $p_-$ opposite to $p$ such that $q = p_{w'}$ and $q = p_{w'}$ is the intersection between $\Th_{w'}(p)$ and $\Th_{w'^\vee}(p_-)$.

  Note that if the point $p$ moves, then the intersection point $C\cap \Th_{w^\vee}(p)$ moves accordingly. More precisely, take $p' = g \cdot p$ for some $g \in G$ close to $1\in G$.
{As $C$ is $G$-invariant, the subvarieties $\Th_{w^\vee}(p')$ and $ C $ meet at $g \cdot q$.}  
   
  Since $\Th_{w^\vee}(p)$ and $C$ intersect properly, there exists an open neighborhood $V$ of $q=p_{w'}$ and a constant $M_1> 0$ such that for any point ${c} \in C \cap V$ and all $g\in G$ close enough to $1\in G$, one has $d_{\cF}(c, {g\cdot} \Th_{w^\vee}(p)) \geqslant M_1 d_{\cF}({c}, {g\cdot} p_{w'})$. Using \eqref{eq_distance_thick} together with the previous inequality applied to $c =\lambda$, one obtains  
  \begin{equation*} 
I(g\cdot p):=  \int_V \log d_{\cF}( g\cdot p , \Th_w(\lambda)) d\mu_{\cF}(\lambda) \geqslant \log (MM_1)  + \int_V \log d_{\cF}(\lambda,g\cdot  p_{w'}) d\mu_{\cF}(\lambda).  
  \end{equation*}
 By Proposition \ref{prop:box_counting_measure},  the Gibbs measure $\mu$ is in the Hausdorff measure class on   
 $\partial_\infty \Gamma$.  Since  the equivariant map $\partial_\infty \Gamma \to \Lambda$ is bi-H{\o}lder to its image by Theorem \ref{thm:biholder}, the measure $\mu_{\cF}$ is also in the Hausdorff measure class $\mathcal{H}^s$ for the metric $d_{\cF}$ 
 for some $s> 0$. Choose an integer $k$ such that $ k s > 1$. 
Take an annulus $A_m = \{ x \in \cF \ | \ m^{-k} \leqslant  d_{\cF}(x,g\cdot p_{w'}) \leqslant
 (m-1)^{-k}  \}$. Observe that  $\Lambda \cap A_n$ has a diameter bounded by $2/(m-1)^k$ so we get   a lower bound 
 \begin{equation*}
 I(g\cdot p) \geqslant \log (M M_1)  - \sum_{m\geqslant 2} k \log (m) \mu_{\cF}(A_m) \geqslant \log (M M_1) - k \sum_{m\geqslant 2} \dfrac{2^{s}\log(m)}{(m-1)^{s k}} > -\infty.    
\end{equation*}
This shows that $I(g\cdot p)$ is uniformly bounded for $g$ in a neighborhood of the identity in $G$. Since this holds for each intersection point in $C \cap \Th_w(p)$, the entire integral is uniformly bounded, as required. 
\end{proof}

\section{Examples of wedge products of Gibbs currents}

Recall from Section \ref{sec: Anosov Gibbs} that we have defined some natural currents $T$ of bidegree $(\dim \cF -k , \dim \cF - k)$ associated to an Anosov group $\Ga$ and its Gibbs measure:  
\begin{equation}
T = \sum_{w \in W} a_w \int \llbracket \Th_w(\lambda ) \rrbracket d\mu_{\cF}(\lambda),
\end{equation}
where 
$a_w \in \C$.

Before stating the main result, we shall use the following version of the shadow lemma (Theorem \ref{thm:shad}). 
The main result of this section is that certain Gibbs currents of  bidegree $(1,1)$ admits a wedge product with any other closed positive current.  

\begin{thm} \label{thm:shad}
Let $\Ga$ be an Anosov group  acting on the flag manifold $\cF= (\bP^1)^n$ and let $w \in W$ an element of the Weyl group with $|w| = \dim \cF-1$. 
Consider a Gibbs measure $\mu_{\cF}$ and a Gibbs current $T$ of the form: 
\begin{equation*}
T:= \int \llbracket \Th_w(\lambda) \rrbracket d\mu_{\cF}(\lambda).
\end{equation*}
Then for any closed positive current $R$ on $\cF$, the wedge product $T\wedge R$ in the sense of Bedford-Taylor exists.  
\end{thm}

\begin{proof} Let $\Lambda$ be the limit set in $\cF$. 
Our aim is to apply  Theorem \ref{thm_wedgeable} to the current $T$, we thus establish the integrability condition \eqref{eq_cond_integrability}. 

Take any point $p\in \cF$, we consider the integral: 
\begin{equation*}
I(p) = \int_{\cF} \log (d_{\cF}(p,\Th_w(\lambda))) d\mu_{\cF}(\lambda) . 
\end{equation*}

Write $p = (p_1,\ldots, p_n)$ where each $p_i$ belong to $\bP^1$. 
We also denote by $\Lambda_i$ the limit set on each component, which are the projection of the limit set $\Lambda$ on the $i$-th component. Fix $w\in W$ with $|w|= n-1$, recall that $\Th_w(\lambda)$ is exactly a fiber over the point $\pi_i(\lambda)$ by $\pi_i$ for some index $i$. 

$\bullet$ Case 1: Assume $p \notin \Th_w(\Lambda)$, this means none of the components $p_i $ belong to $\Lambda_i$. 
In that case, the triangular inequality gives for each $\lambda \in \Lambda$,  $d_{\cF}(p, \Th_w(\lambda)) \geqslant d_{\bP^1}(p_i, \Lambda_i)$. By continuity, there exists a neighborhood $V$ of $p$ on which $d(q_i, \Lambda_i) \geqslant  d_{\bP^1}(p_i, \Lambda_i)/2$ for all $q=(q_1,\ldots, q_n) \in V$.
By integration over the Gibbs measure, one then obtains:
\begin{equation}
I(q) \geqslant \int_{\cF} \log(d_{\cF} (q, \Th_w(\lambda))) d\mu_{\cF}(\lambda) \geqslant \log d_{\bP^1}(p_i,\Lambda_i) - \log(2) , 
\end{equation}
for all $q \in V$.  

\medskip

$\bullet$ Case 2: Assume that $p \in \Th_w(\Lambda)$. This means that $p_i \in \Lambda_i$ some $i\leqslant n$. 
 Consider the pushforward measure $\mu_i= (\pi_i)_* \mu_{\cF}$ by the $i$-th projection. 
 Take a neighborhood $V$ of $p$,
by the triangular inequality, we have for all $q=(q_1,\ldots,q_n)\in V \cap \Th_w(\Lambda)$ and all $\lambda=(\lambda_1,\ldots, \lambda_n) \in \Lambda$, $d_{\cF}(q, \Th_w(\lambda)) \geqslant d_{\bP^1}(q_i, \lambda_i)$, so by integrating over the Gibbs current, one obtains:
\begin{equation*}
I(q) \geqslant \int_{\bP^1} \log d_{\bP^1} (q_i,\lambda_i ) d\mu_i(\lambda_i)
\end{equation*}
Let $\xi_i \colon \partial_\infty \Gamma \to \Lambda_i$ the equivariant boundary map associated to the representation on the $i$-th component. There exists a point $\zeta_q \in \partial_\infty \Gamma$ such that $\xi_i(\zeta_q) = q_i$, and write $\tilde{\mu_i} = (\xi_i^{-1})_* \mu_i$, one has: 
\begin{equation*}
I(q) \geqslant \int_{\partial_\infty \Gamma} \log d_{\bP^1}(\xi_i(\zeta_q), \xi_i(\zeta)) d \tilde{\mu_i}(\zeta)
\end{equation*}

 Theorem \ref{thm:biholder} shows that there exists some uniform constants $C,a > 0$ such that for any point $\zeta, \zeta_q \in \partial_\infty \Ga$, we have: 
\begin{equation*}
d_{\bP^1}(\xi_i(\zeta_q), \xi_i(\zeta)) \geqslant C d_\infty^a(\zeta_q, \zeta).  
\end{equation*}
Summing over the measure $\tilde{\mu_i}$ together with the lower on $I(q)$ gives for all $q\in V\cap \Th_w(\Lambda)$: 
\begin{equation} \label{eq_bound_I_1}
I(q)  \geqslant  \log C + a \int_{\partial_\infty \Ga} \log d_{\infty}(\zeta_q, \zeta) d\tilde{\mu_i}(\zeta).
\end{equation}
Observe that $\mu_i$ is also a Gibbs measure, hence it  
it has a positive critical exponent $s_i$. Consider an integer  $k$ for which $k s_i > 1 $.  
By assertion (3) of Theorem \ref{thm:Gibbs list}, there exists a constant $C'> 0$ such that for any Borel set $A$ on $\partial_\infty \Ga$, one has $C'^{-1} \leqslant \tilde{\mu_i}(A)/\mathcal{H}^{s_i}(A) \leqslant C'$ where $\mathcal{H}^{s_i}$ is the Hausdorff measure with respect to the distance $d_\infty$. 

For $n \geqslant 1$ and $\zeta_q\in \partial_\infty \Ga$, set $$A_{n,\zeta_q} = \left \{ \zeta\in \partial_\infty \Ga \ | \  \dfrac{1}{n^k}\leqslant d_\infty(\zeta_q, \zeta) \leqslant \dfrac{1}{(n-1)^k} \right  \}.$$ 
\eqref{eq_bound_I_1} yields  for all $q \in V \cap \Th_w(\Lambda)$: 
\begin{multline*}
I(q) \geqslant \log C - a k \sum_{n\geqslant 2} \log (n) \tilde{\mu_i}(A_{n,\zeta_q}) 
   \geqslant \log C - a k C' \sum_{n \geqslant 2} \log(n) \mathcal{H}^{s}(A_{n,\zeta_q}) \\
     \geqslant \log C - C_3 \sum_{n \geqslant 2} \dfrac{\log (n)}{(n-1)^{ks_i}} > -\infty,  
\end{multline*}
where $C_3> 0$ is a uniform constant. 
This gives $I(q) > M$ where 
$$
M =\log C -C_3 \sum_{n \geqslant 2} \dfrac{\log (n)}{(n-1)^{ks_i}}$$ 
is finite since $ks_i >1$. Since $I$ is bounded near any point $p$, Theorem \ref{thm_wedgeable} can be applied and the result holds.
%
\end{proof}


\begin{cor} \label{cor_veronese_gen} For any $n\geqslant 2$, consider an  Anosov subgroup $\Gamma$ on $\cF = \GL(n+1,\C)  / B$ induced by a subgroup of the irreducible representation of $\SL_2(\C)$ into $\SL_{n+1}(\C)$ with length function $\beta$. Then for any $w \in W$ with $|w|= \dim \cF -1$, the associated Gibbs-current $T_w = \int \intcur{\Th_w(\lambda)} d\mu_{\cF}(\lambda)$ has bounded local potential, hence the wedge product of $T$ with any other closed positive current is well defined. Moreover, we have for any $w, w_1,\ldots, w_k \in W$ with $|w_i| = \dim \cF -1 $ for all $i\geqslant 1$, one has: 
\begin{equation} \label{eq_wedge_veronese}
T_{w} \wedge \ldots \wedge T_{w_k}  = \int \intcur{\Th_{w}(\lambda)}\wedge \intcur{\Th_{w_1}(\lambda_1)} \ldots \intcur{\Th_{w_k}(\lambda_k)} d\mu_{\cF}(\lambda, \la_1, \ldots, \lambda_k).
\end{equation}
\end{cor}

\begin{proof} By Corollary \ref{cor_intersection_woven}, since the Gibbs-current are automatically woven, we just need to check that the bidegree $(1,1)$ current have bounded local potential.  Using  Theorem \ref{thm_bounded_potential_curve}, we are reduced  to show  that the curve $C \subset \cF$ induced by the Veronese embedding intersects properly  any codimension $1$ thickening. Take $w = w_0 (i, i+1)$ where $w_0 \in W$ is the longest element and where $(i,i+1) \in W$ is the transposition of $i, i+1$. 
Assume by homogeneity that $p $ is given by the standard flag associated to a basis $v_1 ,\ldots , v_{n+1}$.  Take a point $q $ in $\Th_{w}(p) \cap C$. In particular, there exists $w'\leqslant w$ such that $q \in S_{w'}(p) \cap C$.

By \cite[p. 157]{fulton_young}, the point $q$ is a flag induced by the vectors $f_1 ,\ldots, f_{n+1}$ where $f_j$ are of the form: 
\begin{equation} \label{eq_formula_fulton}
f_j = v_{w'(j)} + \sum_{k< w'(j)} \alpha_{j,k} \epsilon_{j,k} v_k \text{ with } \alpha_{j,k} \in \C,  \epsilon_{j,k} = \left \lbrace \begin{array}{ll}
0 & \text{ if } \exists i <j \text{ with } w'(i)=k ,\\
1 & \text{otherwise}. 
\end{array}  \right .
\end{equation} 
This formula gives an affine isomorphism between $\C^{\dim S_{w'}(p)}$ and $S_{w'}(p)$, and under this isomorphism, we obtain that the tangent space to $S_{w'}(p)$ at $q$ is identified with the subspace $V$ of $(n+1)\times (n+1)$ matrices whose non-zero entries $(k,j)$ are determined by the condition $\epsilon_{j,k}=1$ and $k < w'(j)$. 

Consider $i_0 = (w')^{-1}(1)$.
Using the fact that $\mathfrak{n}_+ e_{k} = (n+2 - k) e_{k+1}$ for any $k\leqslant n$ and $\mathfrak{n}_+ e_{n+1}=0$ where $\mathfrak{n}_+ \in \mathfrak{sl}_2$ is the upper-triangular nilpotent generator. 
We consider the tangent vector $X $ induced by $\mathfrak{n}_+$ at $q$. Its flow gives:
\begin{equation*}
\exp(tX) f_{i_0} =  v_1 +  (e^{(n+1)t}-1) v_{2} . 
\end{equation*} 
Recall also that since $q \in C$ and the curve $C$ is invariant by the action of $\SL_2(\C)$, this shows that $\exp(tX) q \in C$. 
The above computation shows that $\exp(tX) f_{i_0}$ is not a linear combination of $v_1$, hence the points $\exp(tX) q $ do not belong in the image $\Psi(\C^{\dim S_{w'}(p)}) = S_{w'}(p)$. 
This shows that the intersection of $C$ with $S_{w'}(p)$ at $q$ is proper. 
Since this argument holds for each intersection point $q \in C\cap \Th_{w}(p)$, we deduce that the intersection is proper, as required.  
\end{proof}

Following the previous results, we can write explicitely the wedge product in these two special cases.

\begin{cor} 
 Let $\Ga$ be an Anosov group acting on $\cF = (\bP^1(\C))^n$.  Consider the Gibbs $(1,1)$-current $T$ obtained as the equidistribution of a generic hypersurface using Theorem \ref{thm:main} whose cohomogy class is equal to $\sum_{w \in W} a_w \{ \Th_w(x_0) \}$. Denote by $\sigma_1, \ldots , \sigma_n$ the standard generators of the Weyl group $W \simeq (\mathbb{Z}_2)^n$.  Then for any $k\leqslant n$, the wedge product $T^{\wedge k}$ in the sense of Bedford-Taylor exists and satisfies: 
\begin{equation*}
T^{\wedge k} = \sum_{ |I| =  k } b_I \int_{(\bP^1)^k} \intcur{\pi_I^{-1}(\{(x_1,\ldots, x_k)\})} d {\pi_{i_1}}_*\mu_{\cF}(x_1) d{\pi_{i_2}}_*\mu_{\cF}(x_2) \ldots d {\pi_{i_k}}_*\mu_{\cF}(x_k),
\end{equation*}
where the sum runs over all multi-indices $I = (i_1, \ldots, i_k)$ with $1 \leqslant i_1 < i_2 < \ldots < i_k \leqslant n$, where $\pi_I$ is the projection of $\cF$ onto the $k$ factors determined by the multi-index $I$ while $\pi_{i_j}$ is the projection onto the $i_j$ factor and where $b_I = \prod_{j=1}^k a_{\sigma_{i_j}}$.  
\end{cor}

\begin{proof}
By Theorem \ref{thm:main}, the Gibbs current is a sum of $n$ laminar currents, each determined by the dual generator $\sigma_i$ of the Weyl group and whose transversal measure is determined by ${\pi_i}_* \mu_{\cF}$. 
So we have: 
\begin{equation*}
T = \sum_{i=1}^n a_{\sigma_i} \int_{\bP^1} \intcur{ \pi_i^{-1} (\{x_i \})} d {\pi_{i}}_* \mu_{\cF}(x_i). 
\end{equation*} 
By Theorem \ref{thm:shad} together with Theorem \ref{thm_intersection_laminar_flag}, the wedge product $T^{\wedge k}$ exists and equals: 
\begin{equation*}
T^{\wedge k} = \sum_{I \subset \{1, \ldots,n \}^k} b_I \int_{(\bP^1)^k} \intcur{\pi_{i_1}^{-1}(\{x_1\})}\wedge \ldots \wedge \intcur{\pi_{i_k}^{-1}(\{x_k\})} d {\pi_{i_1}}_* \mu_{\cF}(x_1) \ldots d {\pi_{i_k}}_* \mu_{\cF}(x_k).  
\end{equation*} 
We separate the sum into two parts, the sum over multi-indices of cardinal $k$ with distinct elements and the ones $R$ where at least two indices coincide. Using the fact that $\intcur{\pi_{i_1}^{-1}(\{x_1\})}\wedge \ldots \wedge \intcur{\pi_{i_k}^{-1}(\{x_k\})} = \intcur{\pi_I^{-1}(\{(x_{1},\ldots, x_{k} ) \})}$, 
we get 
\begin{equation*}
T^{\wedge k} = \sum_{|I|=k} b_I \int_{(\bP^1)^k} \intcur{\pi_I^{-1}(\{(x_1,\ldots, x_k)\})} d {\pi_{i_1}}_*\mu_{\cF}(x_1) d{\pi_{i_2}}_*\mu_{\cF}(x_2) \ldots d {\pi_{i_k}}_*\mu_{\cF}(x_k) + R.
\end{equation*}
We claim that $R$ is equal to zero. Using an appropriate permutation, we can reduce to the case where $I$ is a multi-index with $i_1 = i_2$, then denote by $\Delta$ the diagonal in $\bP^1 \times \bP^1$, the contribution of this term in $R$ is determined by: 
\begin{align*}
R_I &= b_I \int_{(\bP^1)^k} \intcur{\pi_{i_1}^{-1}(\{x_1\})}\wedge \ldots \wedge \intcur{\pi_{i_k}^{-1}(\{x_k\})} d {\pi_{i_1}}_*\mu_{\cF}(x_1) d{\pi_{i_1}}_*\mu_{\cF}(x_2) \ldots d {\pi_{i_k}}_*\mu_{\cF}(x_k) \\ 
& = b_I \int_{(\bP^1 \times \bP^1) \setminus \Delta \times (\bP^1)^{n-2}} \intcur{\pi_{i_1}^{-1}(\{x_1\})}\wedge \intcur{\pi_{i_1}^{-1} (\{x_2\}) } \wedge \ldots \wedge \intcur{\pi_{i_k}^{-1}(\{x_k\})} \prod_{j=1}^n d {\pi_{i_j}}_*\mu_{\cF}(x_j) \\ 
& = 0,
\end{align*}
where we have used Fubini's theorem together with the fact that the transversal measures have no atoms and that the intersection of two distinct fibers of $\pi_{i_1}$ is empty.  
\end{proof}

We now treat the case of Anosov groups in $\SL_3(\C)$ preserving the Veronese curve.

\begin{cor} Let $\Ga$ be an Anosov group obtained via the embedding of $\SL_2(\C) \hookrightarrow \SL_3(\C)$. Let $T_1$ be the Gibbs  $(1,1)$-current obtained as the equidistribution of a generic surface cohomologuous to $ a  \{ \Th_{w_p}(x) \} + b \{ \Th_{w_l}(x) \}$ where $w_l = (231), w_p = (312) \in \mathfrak{S}_3$ and let $T_2$ be the Gibbs $(2,2)$-current obtained as the equidistribution of a generic curve cohomologuous to $ a'  \{ (\pi^{\vee})^{-1}(l_x) \} + b'\{ \pi^{-1}(p_x)\}$ with $x = (p_x, l_x)\in \cF$
Denote by $C_{xy}$ the curve $\{ (p, (p p_x)) \ | \ p \in l_y \}$ for all $x = (p_x, l_x), y = (p_y, l_y)  \in \cF$.  
Then the wedge products in the sense of Bedford-Taylor $T_1 \wedge T_2$, $T_1^{\wedge k}$ with $k\leqslant 3$ exists and satisfy the following relations:  
   \begin{align*}
   T_1 \wedge T_1 & = \int \left (  a^2 \intcur{(\pi^{\vee})^{-1}((p_x p_y))}  + 2 ab \intcur{C_{xy}} + b^2 \intcur{\pi^{-1}(l_x \cap l_y)} \right ) d\mu_{\cF}^{\otimes 2}(x,y), \\
   T_1^{\wedge 3} & = \int 3 a^2 b \delta_{(l_z \cap (p_x p_y), (p_x p_y))} + 3 a b^2 \delta_{(q_{x y}, (p_z q_{x y}))} d\mu_{\cF}^{\otimes 3}(x,y,z) ,\\
   T_1 \wedge T_2  & = \int ab'  \delta_{(p_y, (p_x p_y))} + a' b \delta_{(l_x \cap l_y,l_y)}  d\mu^{\otimes 2}_{\cF}(x,y)
   \end{align*}
   where $q_{xy}$ denotes the intersection $l_x\cap l_y$
\end{cor}

\begin{proof} By Theorem \ref{thm:main}, we have: 
\begin{align*}
T_1 & =  \int a \intcur{\Th_{w_p}(x)} + b \intcur{\Th_{w_l}(x)} d\mu_{\cF}(x), \\ 
T_2 &= \int a'  \intcur{(\pi^{\vee})^{-1}(l_x)} + b' \intcur{\pi^{-1}(p_x)} d\mu_{\cF}(x).
\end{align*}
By Corollary \ref{cor_veronese_gen}, the wedge products $T_1 \wedge T_1 $, $T_1 \wedge T_2$ and $T_1^{\wedge 3}$ exist and Formula \eqref{eq_wedge_veronese} can be applied. 
We shall apply the following relations for $x, y,z$ three pairwise opposite flags:
\begin{align*}
\intcur{\Th_{w_p}(x)} \wedge \intcur{\Th_{w_p}(y)} \wedge \intcur{\Th_{w_p}(z)} = 0 ;  &  \quad & \intcur{\Th_{w_l}(x)} \wedge \intcur{\Th_{w_l}(y)} = \intcur{(\pi)^{-1}(l_x \cap  l_y)} ;  \\
\intcur{\Th_{w_p}(x)} \wedge \intcur{\Th_{w_l}(y)} = \intcur{C_{xy}}  ; & \quad & \intcur{C_{yz}} \wedge \intcur{\Th_{w_p}(x)} = \delta_{l_z \cap (p_x p_y), (p_x p_y)} ;  \\
\intcur{\Th_{w_p}(x)} \wedge \intcur{\Th_{w_p}(y)} = \intcur{(\pi^{\vee})^{-1}((p_x p_y))} ;  & \quad & \intcur{C_{xy}} \wedge \intcur{\Th_{w_l}(z)} = \delta_{(q_{xy}, (p_z q_{xy}))}.   
\end{align*}
Since the three computations are similar, we shall explain the computation for $T_1 \wedge T_2$. 
Observe also that the Gibbs measure $\mu_{\cF}$ has no atoms, so the integrals can be computed away from the union of the diagonals in $\cF \times \cF \times \cF$ (for the computation of $T_1\wedge T_1 \wedge T_1$) and $\cF \times \cF$. 
We have: 
\begin{multline*}
T_1 \wedge T_2 = \int (a \intcur{\Th_{w_p}(x)} + b \intcur{\Th_{w_l}(x)}) \wedge (a'  \intcur{(\pi^{\vee})^{-1}(l_y)} + b' \intcur{\pi^{-1}(p_y)}) d \mu_{\cF}^{\otimes 2}(x,y) \\ 
= \int_{\cF \times \cF \setminus \Delta} (a \intcur{\Th_{w_p}(x)} + b \intcur{\Th_{w_l}(x)}) \wedge (a'  \intcur{(\pi^{\vee})^{-1}(l_y)} + b' \intcur{\pi^{-1}(p_y)}) d \mu_{\cF}^{\otimes 2}(x,y) \\
= \int_{\cF \times \cF \setminus \Delta }  (a a' \cdot 0 + a b' \intcur{\Th_{w_p}(x)} \wedge \intcur{\pi^{-1}(p_y)} + a' b \intcur{\Th_{w_l}(x)} \wedge \intcur{(\pi^{\vee})^{-1}(l_y)}  ) d\mu_{\cF}^{\otimes 2}(x,y) \\
 = \int a b' \delta_{(p_y,(p_x p_y))} + a' b \delta_{(l_x \cap l_y,l_y)} d\mu_{\cF}^{\otimes 2}(x,y),
\end{multline*}
as required.
\end{proof}

%

\chapter[Axiom A aspects]{Axiom A aspects of Anosov group actions}\label{sec:Morse-Smale}
The original Labourie's definition of Anosov representations was in terms of hyperbolicity (along certain sections) of  certain lifts of the geodesic 
flow of a hyperbolic group $\Ga$ to flat bundles over the total space of the flow, where the fibers were products of flag manifolds
$$
G/P_{\taumod} \times G/P_{\iota\taumod}.  
$$
 The goal of this chapter is to show that a similar hyperbolicity phenomenon persists even if we consider the bundle $E\to \cG$ together with a flow $\theta_t$, over the geodesic flow of $\Ga$, whose fibers are diffeomorphic to a single flag manifold of $G$. 
 We will do so only in the case of $B$-Anosov subgroups $\Ga< G$ and  flag manifolds $\cF=G/B$. The general case  is more complicated. 
 We refer to \cite{Sambarino-report} and \cite{Delarue-Monclair-Sanders} 
 for discussion of other {\em metric Anosov} hyperbolic aspects of flows associated with Anosov group actions. 
 
 We will describe compact flow-invariant {\em hyperbolic subsets} $\Sigma_w\subset E$ (images of flow-invariant sections $s_w$) and ambient open 
 flow-invariant subsets $\cO_w\subset E$, indexed by elements of Weyl group $W$. The open subsets $\cO_w$ will be (nonuniformly) transversally 
 foliated by stable/unstable foliations $\cW^\pm_w$  of the flow, with fiberwise leaves diffeomorphic to certain complementary Schubert cells, 
 orbits of unipotent subgroups  of $G$, thus, 
 resulting in a flow satisfying a version of the Axiom A/Morse--Smale  conditions. We will have finitely many  
 flow-invariant compact subsets $\Sigma_w$.  
 The union of these subsets is  precisely the non-wandering subset of 
 the flow which turns out to be the hyperbolic subset. Each $\Sigma_w$ is equivariantly homeomorphic to the geodesic 
 flow $\cG$ of $\Ga$; these are the {\em basic sets} of the flow on $E$.  In particular, the restrictions of the flow to these subsets are topologically mixing, 
 have dense subsets of periodic orbits  and, moreover,  are ergodic with respect to certain flow-invariant ergodic measures 
 (Bowen--Margulis measures on geodesic flows of hyperbolic groups).  
 
 In the case when $\Ga$ is the fundamental group of a compact 
 negatively curved manifold $M$, we literally get an Axiom A flow on $E$ if we take as $\cG$ the unit tangent bundle of $M$ with its geodesic flow. 
 Then we can talk about the tangent spaces of $E$ and they split at $\Sigma_w$ as the direct sums of 1-dimensional neutral subspaces (the flow direction), 
 exponentially contracting and expanding subspaces. In general, due to lack of a smooth structure on $\cG$, 
 we will only talk about stable and unstable foliations in the direction of the fibers of $E$ (over $\cG$). 
 \footnote{Even in this case one can define horospherical foliations $\cG$ using the {\em metric Anosov} flow structure of the latter. 
 Multiplying leaves of these foliations with the leaves of $\cW^\pm_w$ one then obtains metric analogues of stable/unstable foliations of 
 Axiom A flows. We will not do this.} 
 These foliation $\cW^\pm_w$ will form {\em ascending} and {\em descending flags} with respect to the Bruhat order on $W$: 
 On the common domain of the foliations, $\cO_{vw}$, each leaf of $\cW^+_w$ will be contained in a leaf of $\cW^+_v$ provided $v\ge w$, 
 while  each leaf of $\cW^-_v$ will be contained in a leaf of $\cW^-_w$. Leaves 
 of $\cW^+_v$, $\cW^-_w$, $w<v$, will be pairwise transversal (but not uniformly transversal). 
 
 Different invariant subsets $\Sigma_v, \Sigma_w$ will be ``connected'' by {\em separatrices} $\mathring\cL_{wv}$ 
 (open Richardson varieties in flag manifolds) and these connections will be governed by the Bruhat order on $W$:  $\mathring\cL_{wv}$  connects 
  $\Sigma_w$ to $\Sigma_v$ if $w<v$. Each orbit of the semi-flow $(\theta_t, t> 0)$, will be accumulating at the subsets $\Sigma_w$ and while the semiflow will be accumulating locally uniformly at {\em thickenings} of $\Sigma_{w_0}$. 
  On each separatrix $\mathring\cL_{wv}$ the semiflow will locally uniformly accumulate at $\Sigma_v$. 
 
 The ``speed'' of accumulation will be also governed by the Bruhat order on $W$. The forward semi-flow $(\theta, t>0)$ will be exponentially contracting on 
 $\cW^+_w$ and the backward semi-flow $(\theta, t< 0)$ exponentially contracting on $\cW^-_w$, for instance:
 $$
|| d\theta_t(z)||^+ \leqslant C^+_w(z) e^{-\eps t},
 $$  
 where $\eps>0$ depends only on $\Ga< G$ and $C_w^+$ is a continuous function bounded on the semi-oribit $\{\theta_t(z): t> 0\}$.    
 Accordingly, the subsets $\Sigma_w$ will be {\em the components of the hyperbolic set of the flow}: 
At each $\Sigma_w$ the flow will uniformly exponentially contract vectors tangent to  $\cW^+_w$ and uniformly exponentially 
expand vectors tangent to $\cW^-_w$.

\section{Flat bundles and suspension flows}\label{sec:flat bundles}

Let $\Ga$ be a hyperbolic group, $\F$ be a smooth compact manifold, and $\rho: \Ga\to Diff(\F)$ be a smooth 
action of $\Ga$ on $\F$. We let $\widehat\Ga$ denote the geodesic flow for $\Ga$ as discussed in Section \ref{sec:Geodesic flows of hyperbolic groups}. 
The group $\Ga$ acts diagonally on the product $\widehat{E}=\widehat\Ga\times \F$. 
Then the quotient $E=(\widehat\Ga\times \F)/\Ga$ has natural structure of a flat 
fiber bundle $\pi: E\to \cG$ over $\cG=\widehat\Ga/\Ga$, with fibers diffeomorphic to $\F$. Due to lack of smoothness of the base $\cG$, flatness of the 
bundle $E\to \cG$ is understood as a prescription of a  foliation $\cE$ on $E$, whose leaves are projections of the leaves 
$\{\hat{m}\}\times F\subset \widehat E$. However, if $\Ga$ is the fundamental group of a compact negatively curved Riemannian manifold $M$ and we use the unit tangent bundle $UTM$ (with the geodesic flow) as our $\cG$, then $E$ is a smooth fiber bundle and flatness can be understood via a flat Ehresmann connection on $E\to UTM$. 

We let $\widehat{E}_{{\hat m}}$ denote 
the fiber of $\widehat{E}\to \widehat\Ga$ over $\hat{m}$; we will use  the notation $E_{m}$ for $\pi^{-1}(m)$, 
the fiber of $E\to \cG$ over $m\in \cG$.     
We equip $E$ with a fiberwise Riemannian metric continuous along the base 
$\cG$; equivalently, we equip each fiber $\widehat{E}_{{\hat m}}$ 
with a Riemannian metric which varies continuously with respect to ${\hat m} \in \widehat{\Ga}$ 
and is invariant under the $\Ga$-action. 

Sections $s$ of the bundle $E\to \cG$ can be identified with a $\Ga$-equivariant map  
$\hat{s}: \widehat\Ga\to \F$.  A section $s$  is said to be {\em parallel along flow lines} if 
$$
\hat{s}(\hat m)= \hat{s}(\hat m_t)
$$
for all $t\in \R$ and $\hat{m} \in \widehat\Ga$. 

\begin{defi}
A section $s$ of $E$ is said to be 
{\em strongly parallel along flow lines} if for any two flow-lines 
$\hat m,\hat m'$ with the same ideal endpoints, we have
$\hat{s}(\hat m)=\hat{s}(\hat m')$. In order to simplify the terminology, we will simply call such sections {\em flow-invariant}. \index{flow-invariant sections}
\end{defi}

Note that this property is automatic for sections parallel along flow-lines  
for the geodesic flows constructed by Champetier and Mineyev
since (for their flows) 
any two (oriented) flow lines with the same ideal end-points 
are actually equal. 
Flow-invariant sections $s$ define (continuous)  $\Ga$-equivariant maps 
\begin{equation*}
f=f_s: \geo^2\Ga\to \F 
\end{equation*}
 by
\begin{equation}
\label{eq:bdmapsect}
f(\xi_+,\xi_-) = \hat{s}(\hat{m}),
\end{equation}
where $\hat m$ is a flow-line with the ideal endpoints $\hat{m}_\infty=\xi_+, \hat{m}_{-\infty}=\xi_-$.

In the context of our work, flow-invariant sections will play the role of {\em uniformly hyperbolic subsets} of non-uniformly hyperbolic flows.

We define smooth maps 
$$
\hat\theta_{\hat{m},t}: \widehat{E}_{\hat{m}}\to \widehat{E}_{\hat{m}_t}
$$ 
(which is the identity map on the direct factor $\F$), abbreviated as an $\R$-action $\hat\theta_t$ on $\widehat E$,  
and  its projection
$$
{\theta}_{m,t}: {E}_{m}\to {E}_{\phi_t(m)},
$$
abbreviated as $\theta_t$. Note that this projection is well-defined since the flow $\hat\theta_t$ 
commutes with the action of $\Ga$. 
The flow $\theta_t$ on $E$ is the {\em lifted geodesic flow} of 
$\Ga$. \index{$(E,\theta_t)$, lifted geodesic flow of a hyperbolic group}

One can regard $(E,\theta_t)$ as the {\em suspension} of the $\Ga$-action $\rho$ with respect to the geodesic flow $(\cG,\theta_t)$. In the case 
when $\Ga$ is infinite cyclic, generated by $\ga$, $(E,\theta_t)$ amounts to the union of 
suspension flows of the diffeomorphisms $\rho(\ga)$ and $\rho(\ga^{-1})$ of the manifold $\F$. 

Since the subset of $\widehat\Ga$ consisting of normalized elements is relatively 
compact, it follows that for normalized $\hat m, \hat m'\in \widehat\Ga$, the identity map
$$
\widehat{E}_{\hat m} \to \widehat{E}_{\hat m'}
$$
is uniformly bilipschitz (with bilipschitz constant independent of $\hat m, \hat m'$). The same is true for all 
$\hat m, \hat m'$ whose projections to $\Ga$ belong to a fixed finite subset. Even though the maps $\hat\theta_t$ 
are the identity maps in disguise, they (as well as the maps $\theta_{t}$) 
will typically distort the Riemannian metrics of the fibers more and more as $|t|\to\infty$, 
which will be discussed below.

The next lemma allows us to convert the questions about metric distortion by the flow $\theta$ to questions about 
metric distortion on $\F$ by elements of the group $\rho(\Ga)$. 

\begin{lem}\label{lem:from flow to group}
Fix a Riemannian metric $\mathcal{g}_0$ on $\F$ and let $||\cdot||_0$ denote the corresponding norm of tangent vectors. 
Then there is a constant $C\geqslant 1$ depending only on $\mathcal{g}_0$, fibration $E\to \cG$ and the background fiberwise 
Riemannian metric on $E$, such that for every $t\in \R, m\in \cG$ and every nonzero vector 
$\bv\in TE_m$ we have: 
$$
C^{-1} \frac{||d\rho(\ga_t^{-1})(\bu)||_0}{||\bu||_0}\leqslant \frac{||d\theta_t(\bv)||}{||\bv||}\leqslant C 
\frac{||d\rho(\ga_t^{-1})(\bu)||_0}{||\bu||_0}.
$$ 
Here we take $\hat{m}\in \widehat{\Ga}$ to be the normalized element  representing $m$ and  
$\ga_t:= \Pi(\hat{m}_t)\in \Ga$. Lastly, $\tilde\bv$ is  the lift of $\bv$ to $TE_{\hat{m}}$ and 
$\bu\in T\F$ is the vector equal to the projection  of  $\tilde\bv$ under the map $\widehat{\Ga}\times \F\to \F$.    
\end{lem}


\begin{proof}  Consider a normalized element $\hat m\in \widehat\Ga$ and  $t\in\R$. 
Set $\hat m_t:= \hat \phi_t(\hat m)$ and $\ga_t:= \Pi(\hat{m}_t)\in \Ga$.   
Consider the composition
$$
\ga_t^{-1} \circ \hat\theta_{\hat m, t}:  \widehat{E}_{\hat m} \to \widehat{E}_{\ga_t^{-1} \hat m_t}.
$$
Note that $\Pi(\ga_t^{-1}\hat m_t)=\ga_t^{-1} \Pi(\hat{m}_t)=1_\Ga$,
i.e. both $\hat m$ and $\ga_t^{-1}\hat{m}_t$ are normalized.  
Since the group $\Ga$ acts isometrically on the fibers of the bundle $E$, the metric distortion of the above composition 
is exactly the same as the distortion of $\hat\theta_{\hat m,t}$. 
Furthermore, since, as we noted above, the metrics on 
$E_{\hat m}$ and $E_{\ga_t^{-1}\hat m_t}$ and the metric $g_0$ on $\F$ are uniformly bilipschitz to each other (via the ``identity'' map), 
the norm distortion of tangent vectors $\bv \in TE_{\hat m}$ 
by the above composition  (up to a uniform multiplicative error) is the same as the norm distortion 
by the differential of the map 
$$
\rho(\ga_t^{-1}): \F  \to \F
$$
with respect to the metric $g_0$ on $\F$. 
\end{proof}

Similarly to flow-invariant sections, one defines {\em flow-invariant subsets} with respect to $\theta_t$: 
They are projections of subsets $\widehat\cO\subset \widehat E$ such that:

(a) If $\widehat \cO$ contains a point $\hat m$, then it contains the entire trajectory $\hat m_t$. 

(b) If $\mathrm{e}(\hat m)= \mathrm{e}(\hat m')$ then $\hat m\in \widehat \cO$ if and only if $\hat m'\in \widehat \cO$.

\medskip 

In the context of Anosov subgroups, we will have finitely many flow-invariant sections $s_w$ and subsets $\cO_w$, all indexed by elements $w$ of the Weyl group $W$. Moreover, the image of each $s_w$ will be contained in $\cO_w$. 
Note that while images of sections $s_w$ are pairwise disjoint, the subsets 
$\cO_w\subset E$ will be open and dense, hence, will overlap.

We will also construct {\em stable} and {\em unstable} foliations $\cW^+_w, \cW^-_w$ of  $\cO_w$; dimensions of leaves of these foliations will be given by 
$|w|$ and $n-|w|$ respectively. While these foliations will be, clearly, different for different $w$'s, they form a natural {\em flag}: 
On the common domain $\cO_{vw}=\cO_w\cap \cO_v$, the leaves of $\cW^+_w$ will be contained in the leaves of   $\cW^+_v$ 
whenever $v\ge w$ in the Bruhat order on $W$, 
and leaves of $\cW^-_w$ will contain the leaves of $\cW^-_v$. 
The action of $w_0$ will ``swap''  stable and unstable foliations:
$$
\cW^+_w \cap \cO_{ww^\vee}= \cW^-_{w^\vee} \cap \cO_{ww^\vee}. 
$$

\section{Morse property of Anosov subgroups}

One of the fundamental facts of geometry of Gromov-hyperbolic geodesic metric spaces is the {\em Morse Lemma}: Images of quasigeodesics are uniformly close to geodesics. This rigidity property of quasigeodesics fails rather badly in symmetric spaces of rank $\geqslant 2$. For instance, a complete quasigeodesic need not belong to an $r$-neighborhood of any flat for any $r<\infty$. Nevertheless, some forms of the Morse property still hold in higher rank. One of the earliest manifestation of this phenomenon was established in \cite{Kleiner-Leeb} where it was proven that if $\cX$ is a symmetric space of rank $r$ then for any $(L,A)$-quasi-isometric   embedding $q: \R^r\to \cX$, there are maximal flats $F_1,...,F_m$ (with $m$ depending only on $r$) 
and $D=D(L,A)$ such that $q(\R^r)$ is contained in 
$$
N_D(F_1\cup...\cup F_m).
$$
However, what is relevant for us is Morse behavior of {\em uniformly regular quasigeodesics} in $\cX$, proven in \cite{KLP18}. We will state only a week version of the result proven there, which will suffice for our purposes, since we are interested only in {\em complete} 
$\simod$-uniformly regular quasigeodesics  in $\cX$ (where $\cF=G/B$ is the complete flag manifold). 

Before stating the result, we will need several notations and definitions. We first quantify the notion of uniform regularity in $\cX$ introduced in 
\S \ref{sec:def anosov}. We let $\Theta$ denote compact subsets in the interior of the simplex $\simod$, which we will identify with the unit sphere in $\Fmod$. Then a nonzero vector $v\in \Delta$ is $\Theta$-regular if 
$$
\bar{v}=\frac{1}{||v||}v\in \Theta. 
$$  
In other words, the cone $\Delta_\Theta\subset \Delta$ of $\Theta$-regular vectors is the cone with the tip $o$ over the subset $\Theta$.

It is also convenient to think numerically, by fixing a positive number $\eps>0$ and considering the subset $\Delta_\eps\subset \Delta$ consisting of nonzero vectors $v$ such that for every positive root $\alpha\in \Phi^+$,
$$
\alpha(\bar{v})\geqslant \eps. 
$$
We similarly define $\mathfrak{a}_\eps^+\subset \mathfrak{a}^+$. 

An oriented geodesic segment $xy$ in $\cX$ is {\em $\Theta$-regular} if $d_\Delta(x,y)\in \Delta_\Theta$. 
(Alternatively, one says that $xy$ is $\eps$-regular if  $d_\Delta(x,y)\in \Delta_\eps$.) 
An $(L,A)$-quasigeodesic $c: I\subset \R\to \cX$ is {\em $\Theta$-regular} if whenever $t-s> A$, we have $d_\Delta(c(s), c(t))\in \Delta_\Theta$.
A family of quasigeodesics in $\cX$ is said to be {\em uniformly Morse} if there exist $L, A, \Theta$ such that every quasi-geodesic in this family 
is a $\Theta$-regular $(L,A)$-quasi-geodesic, see \cite{KLP14, KLP18}. 
(Equivalently, there exist $\eps> 0, L $ and $A$ such that every quasi-geodesic in this family 
is an $\eps$-regular $(L,A)$-quasi-geodesic.) The terminology {\em uniformly Morse} will be justified by Theorem \ref{thm:Morse} below. 

Then the definition of URU representations $\rho: \Ga\to G$ given in \S \ref{sec:def anosov} can be interpreted as follows:

For every $x\in \cX$ there are consists $L, A$ and a compact subset $\Theta$ in the interior of $\simod$,  
 the orbit map $o_x: \Ga\to \cX, o_x(\ga)=\rho(\ga)x$  send geodesics in a Cayley graph of $\Ga$ to 
$\Theta$-regular $(L,A)$-quasi-geodesics in $\cX$.  In other words, $o_x$ sends geodesics in $Cay_\Ga$ to uniformly Morse quasi-geodesics in $\cX$. 

We are now ready to state a weak form of the {\em Morse Lemma} proven in \cite{KLP18}:

\begin{thm}\label{thm:Morse} 
There exists $D=D(\Theta, L,A)$ with the following properties. 
Let $c: \R\to \cX$ be a $\Theta$-regular $(L,A)$-quasigeodesic. Then the following hold:

1. There exists a pair of antipodal elements $x_\pm\in \cF$ such that $\lim_{t\to\pm \infty} c(t)=x_\pm$, where the limit is understood in the sense of flag-convergence. 

2. For $\bx=(x_+, x_-)$ let $F=F_{\bx}\subset \cX$ be the unique flat in $\cX$ asymptotic to both $x_\pm$ and let 
$V\subset F$ denote a cone over $x_+$. We will think of $F$ as a vector space and the tip of $V$ as the zero vector in $F$. 
Then the image of $c$ is contained in the $D$-neighborhood of $F$. 

3. For all $s\leqslant t$ in $\R$ there exist points $\bar{c}(s), \bar{c}(t)\in F$ within distance $D$ from $c(s), c(t)$, such that the oriented segment 
$\bar{c}(s) \bar{c}(t)$ defines a vector in the cone $V$. In particular, if we take the tip of $V$ to be $\bar{c}(0)$, then $c(\R_+)$ is contained in the $D$-neighborhood of the cone $V$. 
\end{thm}

The {\em Morse property} in this theorem is that uniformly regular uniform quasigeodesics in $\cX$ are uniformly close to flats (more precisely, 
{\em Finsler geodesics} in flats, see \cite{KLP18, Kapovich-Leeb-finsler} for an explanation).

\begin{prop}\label{prop:Morse-prop}
Let $\rho: \Ga\to G$ be an Anosov representation. Then  there exists a compact $C\subset G$ such that the following holds:

Suppose that $\hat m\in \widehat\Ga$ is a normalized element with $\re(\hat m)=(\zeta_+,\zeta_-)$, $\bx=(x_+,x_-)=f(\zeta_+,\zeta_-)\in \La^2$.
 Then for every $t\geqslant 1$ and the projection 
$\ga_t:= \Pi(\hat m_t)\in \Ga$, there exists an element $g_t\in  C$ and a regular transvection $\tau_t\in A^+_{\bx} \subset G$ 
such that:
$$
\rho(\ga^{-1}_t)=g_t \circ \tau_t^{-1}.    
$$
\end{prop}
\begin{proof} The proof is essentially the same as the one in \cite[Theorem 2.41]{KLP17}. 
Since $\hat m$ is normalized, the quasigeodesic $l$ in $\Ga$ corresponding to $\hat m$ has the property that $l(0)=1_\Ga$. The image $c$ of the quasigeodesic $l$ under the orbit map $\Ga\to \Ga x\subset \cX$ is a uniform Morse quasigeodesic forward/backward asymptotic to antipodal 
points $x_\pm\in \cX$. Hence, $c$ is contained in a $D$-neighborhood of the flat $F=F_{\bx}$. In particular, the base-point $o\in \cX$ is in $N_D(F)$. 
Let $\bar{o}\in F$ denote the nearest-point projection of $o$ to $F$. Given $\gamma_t$ as above, we have that $y=c(t)=\rho(\ga_t)o\in N_D(F)$ as well.   
Let $V\subset F$ denote the cone in $F$ with the tip $\bar{o}$ over the chamber $x_+$. 
By the Morse property of the quasigeodesic $c$, there exist a point $\bar{y}\in V$ within distance $D$ from $y$. Let $\tau_t\in A^+_{\bx}$ denote the transvection stabilizing $F$ and sending $\bar{o}$ to $\bar{y}$.  Then 
\begin{equation}\label{eq:2D}
d(o, \tau_t^{-1}\rho(\ga_t) o)\leqslant 2D. 
\end{equation}
We let $\tau'_t$ denote the transvection in $G$ moving $\tau_t^{-1}\rho(\ga_t) o$ back to $o$. By the inequality \eqref{eq:2D}, $\tau'_t$ belongs to a compact subset $C'\subset G$ which depends only on $\Ga$ and $\rho$. The composition
$$
k_t= \tau'_t \circ \tau_t^{-1}\circ \rho(\ga_t) 
$$
belongs to $K$, hence, 
$$
\rho(\ga^{-1}_t)= (k_t^{-1}\tau'_t)  \tau_t^{-1} = g_t\circ \tau_t^{-1}, 
$$
where $g_t$ belongs to the compact subset $C=KC'\subset G$. 
\end{proof}

Continuing with the notation from Proposition \ref{prop:Morse-prop}, we get: 

\begin{cor}\label{cor:from group to transvections}
There exists a constant $C_0\geqslant 1$ depending only on the representation $\rho$ such that the following holds for every normalized element $\hat m\in \widehat\Ga$. For 
every nonzero tangent vector $\bv\in T\cF$, we have 
$$
C_0^{-1} \frac{||d\tau^{-1}_t(\bv)||}{||\bv||}\leqslant \frac{||d\rho(\ga^{-1}_t)(\bv)||}{||\bv||} \leqslant C_0  \frac{||d\tau^{-1}_t(\bv)||}{||\bv||}.
$$ 
\end{cor}

In combination with Lemma \ref{lem:from flow to group}, this corollary will allow us to estimate the contraction rates of the flow 
$\theta_t$ on vectors in $TE_{m}$ 
in terms of contraction rates of appropriate transvections acting on the tangent bundle of the flag manifold $\cF$. 
We note that it is critically important that in Lemma \ref{lem:from flow to group} we get an estimate in terms of $\ga_t^{-1}$ and in the corollary we also work 
with $\ga_t^{-1}$. This transition from flow to transvections would have been impossible if not for the inversion.

\section{Geometry of suspended group actions}

\subsection{Foliations and contraction properties of flows} 

If we have a flat bundle $E\to \cG$, then we can talk about the fiberwise smooth structure and fiberwise Riemannian metrics. We let $E\to \cG$ denote a flat bundle associated with a representation $\rho: \Ga\to Diff(\F)$, where $\F$ is a smooth compact manifold, as described in \S \ref{sec:flat bundles}. We will also equip fibers of the bundle with a background Riemannian metric depending continuously on the points in the base $\cG$.  Let $\theta: E\times \R \to E$ be the lifted geodesic flow. We will consider the following additional structure on $E$:

\begin{enumerate}
\item A finite poset $(\PP, \le)$. (Later on, it will be the  
Weyl group $W$ with its Bruhat order.)  

\item A finite cover of $E$ by flow-invariant open and dense subsets $\cO_p\subset E$, $p\in \PP$. Each subset $\cO_p$ will correspond to a $\rho(\Ga)$-invariant family of open subsets $\cO_{p,\re(\hat{m})} \subset \F$:
$$
\rho(\ga) \cO_{p,\re(\hat{m})}= \cO_{p,\re(\ga\hat{m})}. 
$$
Set $\cO_{pq}:= \cO_p\cap \cO_q$. 

\item For each $p\in \PP$ we will have a pair of flow-invariant transversal foliations $\cW_p^\pm$ of $\cO_p$, 
i.e., for every $m\in \cG$,  and $p\in \PP$, $\cW^\pm_{p,m}$ are smooth 
transversal  foliations of $\cO_{p,m}$: Every leaf of $\cW^+_{p,m}$ intersects each leaf of $\cW^-_{p,m}$ and the intersection is $0$-dimensional. 

We will say that a vector $\bv\in T_zE_{m}, z\in \cO_{p,m}$, is {\em horizontal} (resp. {\em vertical}) if it is tangent to the leaf of  
$\cW^+_{p,m}$, (respectively, $\cW^-_{p,m}$) through $z$. Thus, we have a splitting of the tangent bundle 
$T\cO_p$ as the Whitney sum of horizontal and vertical subbundles, $T^+\cO_p\oplus T^-\cO_p$. 
We will be computing norms of horizontal/vertical tangent vectors with respect to the fiberwise background Riemannian metric on $E$. 
For $p\in \PP$ and $z\in \cO_p$ define 
$$
||d\theta_t||_{p,z}^\pm:= \sup_{\bv \in T^\pm_z \cO_{\pi(z)}\setminus \{0\}} \frac{||d\theta_t(\bv)||}{||\bv||}, 
$$
where the horizontal/vertical subbundles are defined with respect to the foliations $\cW^\pm_p$. 

\item For every $p\in \PP$ we will have a flow-invariant 
section $s_p$ of $\pi: E\to \cG$. 
\index{$s_w$, flow-invariant section}
Accordingly, we will have distinguished leaves 
$\cL^\pm_{p,m}$ of foliations $\cW^+_{p,m}$ containing $s_p(m)$. 
We set 
$$
\cL^\pm_{p}=\coprod_{m\in \cG}\cL^\pm_{p,m}.
$$
The image $\Sigma_p$ of $s_p$ will be disjoint from all $\cO_q, q\ne p$. 
The subsets $\Sigma_p$  are analogues of periodic orbits of Morse--Smale flows. While $\theta_t$ 
is no longer periodic on $\Sigma_p$, it is measure-preserving and ergodic with respect to the push-forward via $s_p$ of the 
Bowen--Margulis measure on $\cG$. 

\item We will have a positive constant $\eps>0$ and positive flow-invariant continuous functions $C_p^\pm: \cL^\pm_p\to \R_+, p\in \PP$. These functions 
extend to continuous functions $C_p^\pm: \cO_p\to \R_+$ such that $C_p^\pm$ is bounded 
along each leaf of $\cW_p^\mp$ and along each semi-orbit $(\theta_t, \pm t\geqslant 0)$.  

\item Foliations $\cW^\pm_{\bullet}$ form two {\em flags} with respect  to the partial order on $\PP$:  
$\cW^+_p\subset \cW^+_q$ for all $p\le q$, in the sense that  each leaf of $\cW^+_{p}\cap \cO_{pq}$ 
is contained in a leaf of $\cW^+_{q}\cap \cO_{pq}$. 
Similarly, if $q\le p$ then each  leaf of $\cW^-_{q}\cap \cO_{pq}$ is 
contained in a leaf of $\cW^-_{p}\cap \cO_{pq}$. 
\end{enumerate}

\begin{rem}
1. In the context of $\taumod$-Anosov representations, the poset $(\PP,\le)$ will be the quotient $W/W_{\taumod}$ of $W$ by its {\em parabolic} Weyl subgroup 
$W_{\taumod}$ (the $W$-stabilizer of the face $\taumod\subset \simod$), and the partial order $\le$ on $\PP$ is the {\em folding order}, the 
projection of the Bruhat order to $\PP$, cf. \cite{manicures, KLP18A}.  

2. More generally, and it is relevant in the case of bundles associated with $\taumod$-Anosov representations, instead of flow-invariant images 
of sections $s_p$ of $E$  one uses flow-invariant compact subbundles $\Sigma_p\subset E$, $\Sigma_p\to \cG$ is the restriction of $\pi: E\to \cG$. 
Fibers of these bundles, in general, will depend on $p$. In the $\taumod$-Anosov setting, 
for $p=[1_W]$ and $p=[w_0]$, the subsets $\Sigma_p$ will be images of flow-invariant sections $s_p$ of $E$. For other $p\in \PP=W/W_{\taumod}$, 
fibers of $\Sigma_p\to \cG$ will be various compact homogeneous spaces  of the Levi subgroup of 
$P_{\taumod}$. 

3. There is nothing special about geodesic flows of hyperbolic groups. It makes sense to consider flows $\phi_t$ on compact topological spaces $\cB=\widehat\cB/\Ga$, where $\widehat\cB= Z\times \R$ and $\Ga$ is a discrete group acting on a compact set $Z$ so that its action 
on $\widehat\cB$ is properly discontinuous. For instance, one can take $\Ga$ to be a relatively hyperbolic group $(\Ga; P_1,...,P_m)$, 
and $Z=D(\geo(\Ga; P_1,...,P_m))$, with $\geo(\Ga; P_1,...,P_m)$ the Bowditch boundary of $(\Ga; P_1,...,P_m)$. This construction is applicable in the setting of {\em relatively Anosov} representations $\rho: \Ga\to G$, cf. \cite{MR4622401, MR4760445}. 

4. The leaf-wise intersections $\cW^-_p\cap \cW^+_q, p\le q$ form foliations $\cW_{pq}$ of $\cO_{pq}$. 
We will see that, in the context of bundles associated with Anosov representations, 
 leaves of $\cW_{wv}, v, w\in W$, (when intersected with fibers of $E\to \cG$) 
are {\em open Richardson varieties} in $\cF$ and their closures will be {\em Richardson varieties} see, \S \ref{sec:Richardson}. 

5. Note that in (3) we do not require {\em uniform transversality} of the leaves of foliations $\cW_w^\pm$: Minimal angles between tangent spaces to leaves can decay to zero as we approach the boundary of $\cO_w$. This happens, for instance, because of singularities of Schubert varieties. However, already in the case of $G=SL(3,\C)$ uniform transversality fails even though all Schubert varieties are smooth. This follows from the description of tangent 
spaces $T_{x_+}\Th_w(x_+)$ to Schubert varieties $\Th_w(x_+)$ mentioned in Remark \ref{rem:local Schubert}, where $|w|=2$. 
\end{rem}

\begin{defi}\index{exponentially contracting flow}
We will say that the flow $\theta_t$ is {\em exponentially forward-contracting} on $\cW^+$ and {\em backward-contracting} on $\cW^-$ with the parameters $C^\pm(z)$ and $\eps$  
if: 

(1)  
$$
||d\theta_{\pm t}||_z^\pm \leqslant C_p^\pm(z) e^{-\eps t} 
$$
for all $t\geqslant 1$ and $p\in \PP$. In other words, for every horizontal  tangent vector $\bv\in T_z^+\cO_p$,  and $t\geqslant 1$, 
$$
||d\theta_t(\bv)||\leqslant C_p^+(z) e^{-\eps t}||\bv||
$$
and for every vertical  tangent vector $\bv\in T_z^-\cO_p$,  
$$
||d\theta_{-t}(\bv)||\leqslant C_p^-(z) e^{-\eps t}||\bv||.
$$

(2) Functions $C_p^+, C_p^-$ are bounded along trajectories of semiflows $(\theta_t, t\geqslant 0)$ and $(\theta_t, t\leqslant 0)$ respectively. 

Such a flow $\theta$ will be called an {\em exponentially contracting Morse--Smale flow} on $E$. 
\end{defi}

Note that, due to compactness of each $\Sigma_p=s_p(\cG)$, the flow $\theta$ will be not only forward/backward exponentially contracting, but also backward/forward uniformly exponentially {\em expanding} at all points $z\in \Sigma_p$ and $t\geqslant 1$:
$$
||d\theta_{\mp t}||_z^\pm \geqslant C_p e^{\eps t}, 
$$
where the constants $C_p$, depend only on $p\in \PP$.

The main goal of the rest of this section is to prove that the lifted geodesic flow $\theta_t$ on the bundle $E\to \cG$ associated with an Anosov representation 
$\rho: \Ga\to G$ is an exponentially contracting Morse--Smale flow, Theorem \ref{thm:main-C}. 

\subsection{Original definition of Anosov representations in terms of flows}

We are now ready to give the original definition of Anosov representations, as introduced for  fundamental groups of compact negatively curved manifolds by Labourie in \cite{Labourie}, and by Guichard and Wienhard for general hyperbolic groups in \cite{Guichard-Wienhard}. 

Our treatment follows  \cite[Appendix]{KLP17}. 
Let $\Ga$ be a hyperbolic group, $\rho: \Ga\to G$ a representation to a semisimple Lie group $G$. Let $\taumod, \iota\taumod$ be faces of the model chamber $\simod$. Let $\cF_\pm$ be the corresponding flag manifolds of $G$. 
Then we obtain a diagonal action of $\Ga$ on the product $\F= \cF_+\times \cF_-$. Recall  
that flow-invariant sections $s$ of the associated bundle 
$\Pi: E\to \cG=\widehat\Ga/\Ga$ amount to equivariant maps $\hat s: \widehat\Ga\to \F$ and corresponding 
maps $f=(f_+,f_-): \geo^2 \Ga\to \F$. 
A section $s$ is {\em antipodal} if for every pair $(\xi_+,\xi_-)\in \geo^2 \Ga$, the simplices 
$f_+(\xi_+), f_-(\xi_-)$ are antipodal. 
The product decomposition of $\F$ defines continuous invariant splittings 
$$
T_{s(m)} E=E^+_{s(m)}\oplus E^-_{s(m)}$$ 
of tangent spaces of fibers $E_{m}$ of $\pi$ at the points $s(m)$.

\begin{defi}\label{defn:original Anosov} \index{Anosov representation}
A representation $\rho: \Ga\to G$ is $(P_{\taumod}, P_{\iota\taumod})$-Anosov if and only if 
there exists an antipodal flow-invariant section $s$ of the bundle $E\to \widehat\Ga/\Ga$ and positive constants $C>0, \eps>0$ 
such that the flow $\theta_{t}$ is 
 uniformly exponentially contracting on $E^\pm$:
 $$
||d\theta_{\pm t}||_z^\pm \leqslant C e^{-\eps t} 
$$ 
\end{defi}

The section $s$ corresponds to a pair of the equivariant maps $f_\pm: \geo \Ga\to \cF_\pm$ whose images are, respectively, the {\em forward} and {\em backward} flag-limit sets $\La_\pm\subset \cF_\pm$. Antipodality of $s$ means that the simplices $f_\pm(\xi_\pm)$ are antipodal to each other, whenever 
$\xi_+\ne \xi_-$. In the situation where $\taumod=\iota\taumod$ 
(which is the case, for instance when $\taumod=\simod$), $\La_+=\La_-$ and $\cF_+=\cF_-$. Equivalence of this definition to other definitions of Anosov representations was established in \cite{KLP17}.

As we will see in Theorem \ref{thm:main-C} below, for $B$-Anosov representations, the 
above definition is a part of a general hyperbolic behavior of the lifted geodesic flow: The exponential contraction in  Definition 
\ref{defn:original Anosov} corresponds to hyperbolicity of $\theta_t, \theta_{-t}$ along just flow-invariant section $s_{1}$ and  
$s_{w_0}$, which is only a small fragment of the overall hyperbolic behavior of $\theta_t$ on $E$, the $\cF$-bundle over $\cG$. 

\section{Foliations associated with lifted geodesic flow}\label{sec:Hyperbolicity of flows}

%

\subsection{Stable/unstable foliations and flow-invariant sections associated with suspended Anosov representations} 
\label{sec:def_flow_section}

We now return to the setting of $B$-Anosov representations to complex semisimple Lie groups (as before, we simply refer to these as Anosov representations). Let $\rho: \Ga\to G$ be an Anosov representation with the boundary map 
$$
\beta: \geo \Ga\to \La\subset \cF.$$
This map gives rise to the map 
$$
\beta^{2}: \geo^2 \Ga\to \La^2\setminus \diag(\La^2)\subset \cY=\Opp(\cF\times \cF).$$

For every $(\zeta_+,\zeta_-)\in \geo^2\Ga$ with $\bx=\beta^2(\zeta_+,\zeta_-)=(\beta(\zeta_+), \beta(\zeta_-))$ 
we define the apartment $a_{\bx}\subset \tits X$ containing the limit points $x_\pm$ and chambers $x_w, w\in W$, in this apartment. The chambers, of course, depend not only on $w$ but also on $\zeta_\pm$. We, thus, obtain equivariant maps 
$$
f_w: \geo^2\Ga\to \cX, \quad f_w(\zeta_+,\zeta_-)= x_w. 
$$ 
These, in turn, define equivariant maps  $\hat{s}_w: \widehat\Ga\to \cF$, $\hat{s}_w(\hat m)= f_w(\mathrm{e}(\hat m))$ and the corresponding 
flow-invariant sections $s_w$  of the associated bundle $E\to \cG$. \index{$s_w$, flow-invariant section}
The image  $s_w(\cG)$ inside $E$ is denoted by $\Sigma_w$. \index{$\Sigma_w$, image of the flow invariant section $s_w$}
We next define flow-invariant open subsets $\cO_w\subset E$. Again, take a pair 
$(\zeta_+,\zeta_-)\in \geo^2\Ga$ with $\bx=\beta^2(\zeta_+,\zeta_-)$ and chambers $x_w$ in the apartment $a_{\bx}$. For every $w\in W$ consider the open and dense subset $\Opp(x_{w^\vee})\subset \cF$. Then for $\hat m\in \widehat\Ga$ with $\mathrm{e}(\hat m)= (\zeta_+,\zeta_-)$ define \index{$\widehat{\cO}_{w}$, flow invariant open subset indexed by $w $ in the Weyl group inside the suspension $\widehat{E}$}
\index{$\widehat{\cO}_w$, flow-invariant open subset index by $w$ in the Weyl group inside the fibration $E$}
$$
\widehat\cO_{w,\hat{m}}:= \{\hat{m}\}\times \Opp(x_{w^\vee}) \times \R  \subset \widehat{E}_{\hat{m}}
$$
and 
$$
\widehat\cO_{w}:= \coprod_{\hat{m}\in \widehat\Ga} \{ \hat{m} \} \times \Opp(x_{w^\vee})\times \R\subset  \widehat{E}. 
$$
Here, as before, $x_v= f_v(\re(\hat{m}))$, $v\in W$. By the construction, $\widehat\cO_{w}$ is $\Ga\times \R$-invariant. 
Lastly, let $\cO_w$ denote the projection of $\widehat\cO_{w}$ to $E$. 

\begin{lem}
$\widehat\cO_{w}$ and $\cO_w$ are flow-invariant, open and dense in $ \widehat{E}$ and $E$ respectively. 
\end{lem}
\begin{proof} Flow-invariance is clear from the construction. Denseness of $\widehat\cO_{w}$ follows from denseness of each $\Opp(x_{w^\vee})$ in $\cF$. To verify that $\widehat\cO_{w}$ is open, observe that its complement is the union 
$$
\coprod_{\hat{m}\in \widehat\Ga} \{ \hat{m} \} \times \Th^{n-1}(x_{w^\vee})\times \R. 
$$
The thickenings $\Th^{n-1}(x_{v})\subset \cF$ depend continuously on $x_v$ (see Lemma \ref{lem:compactness}); 
from this, it follows that the above union is closed in $ \widehat{E}$. 
Hence, $\widehat\cO_{w}$ is open. The statements about $\cO_w$ follow from the fact that the projection $\widehat E\to E$ is a quotient map and 
$\widehat\cO_{w}$ is saturated with respect to this map. 
\end{proof}

\medskip 
Now, we define foliations $\hat{\cW}_w^\pm$ of the subsets $\hat{\cO}_w$. Recall that in \S \ref{section_torus_action}, for every $\bx\in \cY$ we defined 
two transversal foliations of $\Opp(x_{w^\vee})$, namely, $\cW^s_{w,\bx}$ and $\cW^u_{w,\bx}$ using images of affine foliations of the Lie algebra 
$\mathfrak{u}_{x_w}$ under the maps 
$$
\Psi^{-1}_{x_w}:  \mathfrak{u}_{x_w}\to \Opp(x_{w^\vee})$$
which are compositions of the exponential map 
of the Lie group $U_{x_w}$ and the orbit map $U_{x_w}\to U_{x_w} x_w$. Below, we will modify these definitions and define new pairs of transversal foliations 
$\cW^\pm_{w,\bx}$ of $\Opp(x_{w^\vee})$. The new foliations will have better behaved asymptotics and, thus, nicer compactifications in $\cF$. 

\medskip

In the same section \S \ref{section_torus_action}, we defined the maximal torus $T_{\bx}< G$, its split subtorus $A_{\bx}< T_{\bx}$,  both 
stabilizing the apartment $a_{\bx}$;  the Lie algebra of $A_{\bx}$ is 
$\mathfrak{a}_{\bx}$. We, furthermore, defined a maximal unipotent subgroup $U_{x_w}< G$, stabilizing the point $x_{w^\vee}$, 
its nilpotent Lie algebra $\mathfrak{u}_{x_w}\subset \mathfrak{g}$ 
and the decomposition as a direct sum of two subalgebras
$$
\mathfrak{u}_{x_w}= \mathfrak{u}^+_{x_w} \oplus \mathfrak{u}^-_{x_w}. 
$$
Here  $\dim(\mathfrak{u}^+_{x_w})=|w|$, $\dim(\mathfrak{u}^-_{x_w})=n-|w|$ and  
$$
\mathfrak{u}^\pm_{x_w}= \mathfrak{u}^\mp_{x_w^\vee}.
$$ 
The subgroup $A_{\bx}$ normalizes all the unipotent subgroups $U_{x_w}, w\in W$. 
We also defined unipotent subgroups
$$
U^\pm_{x_w}:= \exp(\mathfrak{u}^\pm_{x_w})< U_{x_w}. 
$$
The relations between these subalgebras and subgroups under the change of $w\in W$ is described as follows:

\begin{lem}\label{lem:root-inclusions}
If $v\ge w$ in the Bruhat order on $W$, then  
\begin{equation}\label{eq:plus}
\mathfrak{u}^+_{x_w}\subset \mathfrak{u}^+_{x_v}, \quad U^+_{x_w}\leqslant U^+_{x_v}, 
\end{equation}
\begin{equation}\label{eq:minus}
\mathfrak{u}^-_{x_w}\supset \mathfrak{u}^-_{x_v},  \quad U^-_{x_v}\leqslant U^-_{x_w}. 
\end{equation}
\end{lem}
\begin{proof} Recall that in \S \ref{section_torus_action}  we defined a root system 
$\Phi=\Phi_{\bx}$ in ${\mathfrak a}_{\bx}^*$ and its  subsets $\Phi^\pm_{w}$ for $w\in W$.  
Since $v\ge w$, every root $\alpha\in \Phi$ which is negative on $w {\mathfrak a}_{\bx}^+$ 
 is also negative on $v{\mathfrak a}_{\bx}^+$. Thus,
\begin{equation}\label{eq:root-inclusion} 
\Phi_w^+\subset \Phi_v^+, \Phi_v^-\subset \Phi_w^-
\end{equation}
hence,
$$
{\mathfrak u}_w^+\subset {\mathfrak u}_v^+, {\mathfrak u}_v^-\subset {\mathfrak u}_w^-, 
$$
as claimed. Exponentiating, we get inclusions of the corresponding unipotent subgroups. 
\end{proof}

Given $\bx= \beta^2(\hat{m})\in \cY$ with $\hat{m} \in \partial^2_\infty \Gamma$ and $w\in W$ we define  the foliations $\hat{\cW}^\pm_{w,\bx}$ of $\hat{\cO}_{w,\hat{m}}$ as foliations by orbits of the groups  $U^\pm_{x_w}$ where the group acts trivially on $\hat{\Gamma}$. 
Compared to the foliations $\cW^s_{w, \bx} ,\cW^u_{w, \bx}$ constructed in \S \ref{sec:Structure of thickenings}, this double foliation is 
biholomorphic to the double foliation $\cW^u_{w,\bx}, \cW^s_{w,\bx}$ on $\hat{\cO}_{w,\bx}$. In particular, every leaf of 
$\hat{\cW}^+_{w,\hat{m}}$ intersects every leaf of $\hat{\cW}^-_{w,\hat{m}}$ exactly once. 
Since the subgroups $U^\pm_{x_w}$ are normalized by the maximal 
torus $T_{\bx}$, the action of the torus preserves both foliations $\hat{\cW}^\pm_{w,\hat{m}}$. Leaves of $\cW^s_{w,\bx}, \cW^u_{w,\bx}$ through the point $x_w$ are also leaves of the foliations $\hat{\cW}^\pm_{w,\hat{m}}$, we will use the notation $\hat{\cL}^\pm_{w,\hat{m}}$ for these leaves. Thus, $\hat{\cL}^+_{w,\hat{m}}= S_w(x_+)$ and 
$\hat{\cL}^-_{w,\hat{m}}= S_{w^\vee}(x_-)$, Schubert cells based at $x_\pm$. In particular, intersections 
$$
\hat{\cL}_{vw,\hat{m}}:=\hat{\cL}^-_{w,\hat{m}}\cap \hat{\cL}^+_{v,\hat{m}}
$$
are open Richardson varieties in $\cF$ of dimension $|v|-|w|$ (see \S \ref{sec:Richardson}).  Lastly, 
define 
$$
\hat{\cO}_{vw,\hat{m}}:= \hat{\cO}_{v,\hat{m}}\cap \hat{\cO}_{w,\hat{m}}, 
$$
and similarly:
$$
{\cO}_{vw,m}:= {\cO}_{v,m}\cap {\cO}_{w,m}, 
$$
where $\Pi(\hat{m}) = m \in \hat{\Gamma}/\Gamma= \mathcal{G}$.

\section{Basic properties of vertical/horizontal foliations}

As a consequence of Lemma \ref{lem:root-inclusions} we obtain:

\begin{lem}
Suppose that $v\ge w$ in the Bruhat order. Then each leaf of $\hat{\cW}^+_{w,\hat{m}}\cap \hat{\cO}_{vw,\hat{m}}$ is contained in a leaf of 
$\hat{\cW}^+_{v,\hat{m}}\cap \hat{\cO}_{vw,\hat{m}}$ and each  leaf of $\hat{\cW}^-_{v,\hat{m}}\cap \hat{\cO}_{vw,\hat{m}}$ is 
contained in a leaf of $\hat{\cW}^-_{w,\hat{m}}\cap \hat{\cO}_{vw,\hat{m}}$. 
\end{lem}

We also have:  

\begin{lem}[Non-uniform Morse--Smale Property] \label{lem:Morse--Smale Property}
If $w\le v$ then each leaf of $\hat{\cW}^-_{w,\hat{m}}\cap \hat{\cO}_{vw,\hat{m}}$ is transversal to each leaf of 
$\hat{\cW}^+_{v,\hat{m}}\cap \hat{\cO}_{vw,\hat{m}}$. 
\end{lem}

\begin{rem} This lemma  and the previous one are stated for $\hat{m}$ in the image of $\beta^2$ of a geodesic line in the group $\Gamma$, however it can be generalized to any pair of opposite points $\bx \in \mathcal{Y}$ where one replaces $\hat{\cO}_{w, \hat{m}}$ with the subset $\Opp(x_{w^{\vee}})$  
\end{rem}

\begin{proof} In view of the previous lemma, it suffices to check transversality of the pair of foliations $\hat{\cW}^\pm_{w,\hat{m}}$.  
Let $z$ be an intersection point of two leaves of these foliations. Since $U_{x_w}$ acts simply transitively on 
$\Opp(x_{w^\vee})$, we can as well take $z$ as our base-point for the identification of $\Opp(x_{w^\vee})$ with 
$U_{x_w}$ via the orbit map. Thus, the claim reduces to transversality of the Lie subalgebras $\mathfrak{u}^\pm_{\bx_w}\subset \mathfrak{u}_{\bx_w}$, 
which, in turn, follows  from the fact that 
$$
\Phi_{\bx}= \Phi^+_{x_w} \sqcup \Phi^-_{x_w}. 
$$
 \end{proof}


It appears that the key issue to analyze is the geometry of intersection of Schubert varieties $\Th_w(x_+)$  and $\Th_{v}(x_{w^\vee})$ at the chamber $x_+$ where $v=(w^\vee)^{-1}=  w^{-1}w_0 $; 
more precisely, if tangent spaces of these varieties at $x_+$ have nonzero intersection. 

\medskip 
Observe that the action of $w_0$ swaps $w$ and $v=w^\vee$, $x_+$ and $x_-$, $x_w$ and $x_{v}$ 
and, accordingly, sends  $\Phi^\pm_{x_w}$ to  $\Phi^\mp_{x_v}$. 
Hence, the element $k_{w_0}\in N_G(T_{\hat{m}})$ corresponding 
to $w_0$ conjugates $U^\pm_{x_w}$ to $U^\mp_{x_v}$. It follows that, $k_{w_0}$ sends 
$\hat{\cW}^\pm_{w,\hat{m}}\cap \hat{\cO}_{vw,\hat{m}}$ to $\hat{\cW}^\mp_{v,\hat{m}}\cap \hat{\cO}_{vw,\hat{m}}$, swapping vertical and horizontal foliations. 

In view of  the knit-product decomposition $U_{x_w}= U^+_{x_w} U^-_{x_w}= U^-_{x_w} U^+_{x_w}$ 
(Corollary \ref{cor:knit}) we see that the action of $U_{x_w}$ preserves both foliations 
$\hat{\cW}^\pm_{w,\hat{m}}$ and, furthermore, every leaf of $\hat{\cW}^\pm_{w,\hat{m}}$ intersects $\hat{\cL}^\mp_{w,\hat{m}}$ 
in exactly one point, which allows us to label leaves of $\hat{\cW}^\pm_{w,\hat{m}}$ by points of 
$\hat{\cL}^\mp_{w,\hat{m}}$, equivalently, by elements of $U^\mp_{x_w}$.  We, thus, define maps
\begin{equation}\label{eq:pm-proj} 
\pi_{w}^{\pm}: \hat{\cO}_{w, \hat{m}}\to \hat{\cL}^\pm_{w,\hat{m}},
\end{equation}
by sending each $z\in \hat{\cO}_{w}$ to the unique intersection point of $U^\mp_{x_w} z$ with  $\hat{\cL}^\pm_{w,\hat{m}}$. These maps define a biholomorphic map 
$\hat{\cO}_{w,\hat{m}}\to \cL^+_{w,\hat{m}}\times \cL^-_{w,\hat{m}}$, $z\mapsto (\pi_{w}^+(z), \pi_{w}^-(z))$.

\medskip
With all these preliminaries out of the way, we can finally  define foliations $\cW^\pm_w$ of open subsets $\cO_w\subset E$. 

\medskip 
As before, for every $\hat m\in \widehat\Ga, m=\Pi(\hat{m})\in \cG$ we let $\bx:= \beta^2(\re(\hat m))$ and, accordingly, define the chambers 
$x_w=f_w(\hat m)\in a_{\bx}, w\in W$.  We then
define the foliations $\cW^\pm_{w,m}$
 on $\cO_{w,m}$  as projections of   the foliations 
$$
  \hat{\cW}^\pm_{w,\hat{m}}\subset \widehat{E}_{\hat m}
$$
to $E_m$. Similarly, define $\cL^\pm_{w,m}$ as the projection to $E_m$ of 
$$
 \hat{\cL}^\pm_{w,\hat{m}}\subset \widehat{E}_{\hat m}. 
$$
Lastly, set
$$
\cW^\pm_{w}:= \coprod_{m\in \cG} \cW^\pm_{w,m}  
$$
and
$$
\cL^\pm_{w}:= \coprod_{m\in \cG} \cL^\pm_{w,m} 
$$

Below we collect basic properties of these foliations which follow from the definition and 
our earlier results on foliations $\hat{\cW}^\pm_{w,\hat{m}}$:

\begin{enumerate}
 \item Using the automorphism induced by the multiplication by $k_{w_0} \in N_G(T_{\hat{m}})$ inducing the element $w_0$ in $W$ and using the fact that the conjugate of the groups $U_{w,\bx}^{\pm}$ by $k_{w_0}$ are $U_{w^{\vee},\bx}^{\mp}$, one has  
$$
\cW^\pm_{w}\cap \cO_{ww^\vee} = \cW^\mp_{w^\vee}\cap \cO_{ww^\vee}. 
$$

\item Leaves of $\cW^+_{w,m}, \cW^-_{w,m}$ are properly embedded complex submanifolds of dimensions 
$|w|$ and $n-|w|$ respectively in  $\cO_{w,m}$, each biholomorphic to $\C^{|w|}, \C^{n-|w|}$ respectively. 
If $w=1_W$ then leaves of $\cW^+_{w,m}$ are singletons and there is only one leaf of $\cW^-_{w,m}$, namely 
$\cO_{w,m}$. Similarly, if $w=w_0$ then leaves of $\cW^-_{w,m}$ are singletons and there is only one leaf of $\cW^+_{w,m}$, namely 
$\cO_{w,m}$ .

\item If $w\le v$, then the foliation $\cW^-_{w}$ is transversal to $\cW^+_{v}$ (on $ \cO_{wv}$). 

\item If $w\le v$, then  $\cW^+_{w}\subset \cW^+_{v}$ and  
$\cW^-_{w}\subset \cW^-_{v}$ (on $ \cO_{wv}$), meaning that every leaf of  $\cW^+_{w}\cap \cO_{vw}$ is contained in a leaf of 
$\cW^+_{v}\cap \cO_{vw}$, and every leaf of  $\cW^-_{v}\cap \cO_{vw}$ is contained in a leaf of $\cW^-_{w}$. 
In other words, foliations $\cW^+_{w}$ and $\cW^-_{w}$ 
form (respectively) ascending and descending flags with respect to the Bruhat order on $W$. 

\item 
Foliations $\cW^\pm_{w}$ are flow-invariant, meaning that $\theta_t$ sends leaves to leaves. 

\item Foliations $\cW^\pm_{w,m}$ depend continuously on $m\in \cG$. 

\item Stable/unstable submanifolds $\cL^\pm_{w,m}$ depend  continuously on $m\in \cG$ and $\cL^\pm_{w}$ are flow-invariant.  

\item If $w\le v$, then intersections of leaves of $\cW^-_{w,m}, \cW^+_{v,m}$ with $\cO_{vw,m}$ 
are transversal, nonempty and have dimension $|v|-|w|$ . If $v=w$ then every leaf of 
 $\cW^-_{w,m}$ intersects every leaf of $\cW^+_{w,m}$ in exactly one point. 

\item If $w< v$, then the intersections $\cL_{vw,m}:=\cL^-_{w,m}\cap \cL^+_{v,m}$  
are open Richardson varieties in $E_m$ of dimension $|v|-|w|$. If $v=w$ then this intersection is the singleton $\{s_w(m)\}$.

\item $s_w(m)$ belongs to $\cO_{w,m}$ and does not belong to $\cO_{v,m}$ for any $v\in W\setminus \{w\}$. 

\item Each leaf of $\cW^\pm_{w,m}$ is properly embedded in $\cO_{w,m}$. Moreover, one has ``locally uniform properness'' of these embeddings: 
Consider an open subset $Q\subset \cG$ over which the bundle $E$ trivializes as the product $Q\times \cF\cong \pi^{-1}(Q)$. 
For every compact subset $C\subset  \pi^{-1}(Q) \cap \cO_{w,m}$, there is a constant $D=D_C$ such that 
for every leaf $L$ of $\cW^\pm_{w,m}$ ($m\in Q$), if $x, y\in L\cap C$, then the distance between $x, y$ measured in the induced Riemannian metric on $L$ is bounded by $D$. 

\end{enumerate}

\begin{lem}\label{lem:disjoint}
$\cL^+_{w,m}$ is disjoint from $\bigcup_{v>w} \cO_{v,m}$. 
\end{lem}
\begin{proof} For every $v> w$ we have $|v^\vee|< |w^\vee|< n-|w|$. On the other hand, for every $x\in \hat{\cL}^+_{w,\hat{m}}= S_w(x_+)$ we have 
$D(x_+, x)=D(x_+,x_w)=|w|$. Therefore, $D(x_{v^\vee}, x)< n-|w|+|w|<n$. It follows that $x$ cannot be antipodal to $x_{v^\vee}$, i.e. 
cannot belong to $\cO_{v,m}$. 
\end{proof}

\section{Tubes around thickenings}\label{sec:tubes}

Our next goal is to understand compactifications of leaves of $\hat{\cW}^\pm_{w,\hat{m}}$ in $\cF$ and how unions of these form tubular neighborhoods of thickenings in $\cF$.  

 Consider the thickening $\Th_{w^\vee}(x_-)$. It is the closure of the 
Schubert cell $S_{w^\vee}(x_-)=U^-_{w^\vee}\cdot x_w$. Let $\CC\subset \mathfrak{u}_w^+$ 
be a compact convex subset containing $0$. 

\begin{defi}\label{def:tube}
Given such $\CC$,  
we define the {\em $\CC$-tube} around $S_{w^\vee}(x_-)$, denoted $N_{\CC}(S_{w^\vee}(x_-))$, as
$$
\exp(\CC)U^-_w \cdot x_w. 
$$
Similarly, when $\bx = \beta^2(\hat{m})$ with $\hat{m} \in \hat\Gamma$, we set $N_C(\hat{m} \times S_{w^\vee}(x_-))$ to be the set: 
\begin{equation*}
\exp(\CC) U_w^- \cdot (\hat{m},x_w )= \{ \hat{m} \} \times N_{\CC}(S_{w^\vee}(x_-))  \subset \hat{E}_{\hat{m}}. 
\end{equation*}
\end{defi} 

Thus, the tube $N_{\CC}(S_{w^\vee}(x_-))$ is the union of Schubert cells 
$$
g  S_{w^\vee}(x_-)= S_{w^\vee}(g x_-), \quad g\in \exp({\CC}).
$$ 
Note that for every $\ga^{-1}\in T_{\bx}^+$, its inverse $\ga$ 
maps $N_{\CC}(S_{w^\vee}(x_-))$ into itself since  $Ad(\ga)\CC\subset \CC$.  
If $0$ belongs to the interior of $\CC$, then $N_{\CC}(S_{w^\vee}(x_-))$ is a neighborhood of $S_{w^\vee}(x_-)$ in 
$\cO_{x_w}$. We also define
$$
N_{\CC}(\Th_{w^\vee}(x_-)):= \bigcup_{v\ge w} N_{\CC}(S_{v^\vee}(x_-)). 
$$

For a subset $E$ of a topological space $X$ we define the {\em partial boundary} $\partial' E$ as the complement $\ol{E}\setminus E$. 
For instance, if $S_{w^\vee}(x)$ is a Schubert cell in $\cF$, then its partial boundary in $\cF$ equals 
$$
\partial' S_{w^\vee} (x) = \Th_{< w^\vee}(x) =\bigcup_{v> w} \Th_{v^\vee}(x).
$$
where $ \Th_{<w^\vee}(x)$ was defined earlier  in \S \ref{sec:Thickenings in the flag manifold}.

\begin{lem}\label{lem:compactifications of leaves} 
1. The partial boundary $\partial' N_{\CC}(S_{w^\vee}(x_-))$ of $N_{\CC}(S_{w^\vee}(x_-))$ in $\cF$ equals
$$
\partial' S_{w^\vee}(\exp({\CC})x_-)= \bigcup_{g\in \exp({\CC})} g \partial' S_{w^\vee}(x_-)= \bigcup_{g\in \exp({\CC})}  \partial' S_{w^\vee}(g x_-). 
$$
 In particular, $\partial' N_{\CC}(S_{w^\vee}(x_-))$ is compact and equals the thickening $\Th_{<w^\vee}( \exp({\CC})x_-)$, where 
$$
 \Th_{< u}= \{v\in W: v< u\}. 
$$ 

2. When $\bx = \beta(\hat{m})$ and $m=\Pi(\hat{m})$,
the projection of $\{ \hat{m}\} \times \partial'(N_{\CC}(S_{w^\vee}(x_-)))$ to $E_m$ is contained in $\bigcup_{v>w} \cO_{v,m}$ and, thus,
$$
\ol{N_{\CC}(\{ \hat{m} \} \times S_{w^\vee}(x_-))}= N_{\CC}(\{ \hat{m} \}\times \Th_{w^\vee}(x_-))\subset \bigcup_{v\ge w} \hat{\cO}_{v, \hat{m}}. 
$$
\end{lem}
\begin{proof} 1. The inclusion 
$$
\bigcup_{g\in \exp({\CC})} g \partial' S_w(x_-) \subset \partial' N_{\CC}(S_w(x_-))
$$
is clear. To prove the opposite inclusion, consider a sequence $x_i\in N_{\CC}(S_w(x_-))$ converging to a point $x\in \cF$. Then  
$x_i= \exp(c_i)(y_i)$, for some $c_i\in \CC$, $y_i\in S_w(x_+)$. After extraction, we can assume that  $c_i\to c\in \CC$, $y_i\to y\in \partial' S_w(x_-)$. 
But then $x= \exp(c)(y)$. 

2.  Note that $\CC\subset  \mathfrak{u}_w^+\subset  \mathfrak{u}_v^+$, $v> w$. By Part 1, every $x\in \partial' N_{\CC}(S_{w^\vee}(x_-))$ lies in  
$S_{v^\vee}(g x_-)$ for some $g\in \exp({\CC})$ and $v> w$. But  $\{ \hat{m} \} \times S_{v^\vee}(g x_-)\subset \hat{\cO}_{v,\hat{m}}$ since $g\in U_w^+\subset U_v^+$. 
\end{proof}

\begin{cor} For any $\bx = \beta^2(\hat{m})$ with $\hat{m} \in \partial_\infty^2 \Gamma$, one has:
$$
\partial' N_{\CC}( \{\hat{m}\} \times S_w(x_-))= \bigcup_{v<w} N_{\CC}( \{ \hat{m} \} \times \Th_v(x_-)). 
$$
\end{cor}


Of course, even if $\CC$ has nonempty interior and, thus, $N_{\CC}(S_w(x_-))$ is a neighborhood of $S_w(x_-)$ in $\cF$, 
$\partial' \Th_w(\exp({\CC})x_-)$ will not be a neighborhood of $\partial \Th_w(x_-)$, simply because its dimension is too small.



\section{Proof of exponential contraction} 

\subsection{Contraction estimates for a Finsler metric}\label{sec:Contraction estimates}

For every $w\in W$ define a $U_{x_w}$-invariant measurable Riemannian metric on $\cF$ by  
$$
||\bv||_w=||d\Psi_w(\bv)||_0
$$
if $\bv\in T_{x_w}\cF$ and $||\bv||_w= 0$ if $\bv\notin T \Opp(x_{w^\vee})$.  Here, as before $\Psi_w=\Psi_{x_w}$ 
is the chart on $\Opp_{x_{w^\vee}}$ defined in \S \ref{sec:parametrization}.
Recall that if $\bx = \beta^2(\hat{m})$ with $\hat{m} \in \partial^2_\infty \Ga \setminus \Delta$ and $m = \Pi(\hat{m})$, then there is a natural identification between $\{\hat{m} \} \times \cF $ and $\hat{\cO}_{\hat{m}}$. In particular, the image of the tangent vector $\bv$ via this biholomorphism is a tangent vector in $T \hat{\cO}_{\hat{m}}$. We call these vector vertical as they live on the fiber.  

\begin{lem}\label{lem:diverge} Given $\hat{m} \in \partial_\infty^2 \Gamma \setminus \Delta$, $\bx = \beta^2(\hat{m})$. There exists $\epsilon_0> 0$ such that if 
if  $z_i\in \hat{\cO}_{w,\hat{m}}$ is any sequence which diverges to infinity in $ \hat{\cO}_{w,\hat{m}}$, i.e. accumulates only at $ \{\hat{m}  \} \times \Th^{n-1}(x_w)= \{\hat{m} \} \times 
\cF\setminus \hat{\cO}_{w, \hat{m}}$, if  
$\bv_i\in T_{z_i}\hat{\cO}_{\hat{m}}$ is  any sequence of unit vertical vectors with the respect to the background Riemannian metric on $\{ \hat{m} \} \times \cF$, then 
$$
\limsup_{i\to\infty} ||\bv_i||_w \leqslant  \epsilon_0^{-1} <\infty. 
$$
\end{lem}
\begin{proof} Let $g_i\in U_{x_w}$ be a sequence of unipotent elements such that $z_i=g_i \cdot (\hat{m}, x_{w})$. Then $\bv_i= dg_i(\tilde\bv_i)$ for some vectors 
$\tilde\bv_i\in T_{x_w}\cF$.  
According to Lemma \ref{lem:unipotent bounds}, there exists $\eps_0>0$ (which depends only on the pair $\bx$ and the choices the norm $||\cdot||_0$ and the background Riemannian metric on $\cF$) such that 
$$
||\bv_i||\geqslant  \eps_0 ||\tilde\bv_i||_0= ||\bv_i||_w. 
$$
Since, by the assumption, $||\bv_i||=1$, we get
$$
||\bv_i||_w\leqslant \eps_0^{-1}. 
$$
\end{proof}

\begin{defi}\label{def:w-norm}
Set 
$$
\hat{\cO}_{\hat{m}}^w:= \bigcup_{v\ge w} \hat{\cO}_{v,\hat{m}}
$$
and define the norm $||\bv||^w$ on tangent spaces to $\hat{\cO}_{\hat{m}}^w$ by 
$$
||\bv||^w= \max_{v\ge w} ||\bv||_v. 
$$
\end{defi}

\begin{lem}\label{lem:comparison}
$||\cdot||^w$ is a measurable Finsler metric on $T\hat{\cO}_{\hat{m}}^w$ which 
is locally bilipschitz equivalent to the restriction of the Riemannian metric on $\cF$ in the following sense:  

There exists a positive continuous function $c_w: \hat{\cO}_{\hat{m}}^w\to \R_+$  such that:
 $$
\frac{1}{c_w(z)}||\bv|| \leqslant ||\bv||^w\leqslant c_w(z)||\bv||, 
$$
where  $z\in \hat{\cO}_{\hat{m}}^w$ and where $\bv\in T_z\hat{\cO}_{\hat{m}}^w$ is a vertical vector.  
\end{lem}
\begin{proof} The fact that $||\cdot||^w$ is measurable Finsler is clear. We need to check the bilipschitz property. The key is to verify that for every 
$z\in \hat{\cO}_{\hat{m}}^w$, there exists $v \geqslant w$ such that $z \in \hat{\cO}_{v,\hat{m}}$ and 
\begin{equation}\label{eq:inf-sup} 
0<\liminf_{x\to z} \inf_{\bv\ne 0, \bv\in T_x \cF} \dfrac{||\bv||_v}{|| \bv||} \leqslant \liminf_{x\to z} \inf_{\bv\ne 0, \bv\in T_x \cF} \frac{||\bv||^w}{||\bv||}\leqslant \limsup_{x\to z} \sup_{\bv\ne 0, \bv\in T_x \cF} \frac{||\bv||^w}{||\bv||} \leqslant \epsilon <\infty. 
\end{equation}
It suffices to consider vectors $\bv\in T_x \cF$ with $||\bv||=1$. Suppose that $x_i\ne z$ is a sequence in $\hat{\cO}^w_{\hat{m}}$ 
converging to  $z$.  Lemma \ref{lem:diverge} implies that
$$
\limsup_{i\to\infty}  ||\bv_i||^w \leqslant  \epsilon_0^{-1} <\infty, 
$$
where $\epsilon_0$ depends only on $w$ and $\bx$. 
This establishes one inequality (the upper bound). To prove the lower bound, observe that 
$$
||\bv_i||^w \geqslant ||\bv_i||_w. 
$$ 
Let $v\ge w$ be such that $z\in \hat{\cO}_{v, \hat{m}}$. Then for all large $i$, $x_i\in \hat{\cO}_{v, \hat{m}}$. We also have (by the Definition \ref{def:w-norm}) 
$$
||\bv_i||^w \geqslant ||\bv_i||_v. 
$$
Since the norm $||\cdot||_v$ comes from a Riemannian metric on $\cO_{x_v}$, it follows that 
$$
\liminf_{i\to\infty}  ||\bv_i||_v>0,$$
and that it is continuous on $\hat{\cO}_{v,\hat{m}}$. 
Thus, we proved \eqref{eq:inf-sup}. 

Inequalities \eqref{eq:inf-sup} imply that there exists a positive continuous function $c_w(z)$ on $\cO_{\bx}^w$ 
such that  
$$
\frac{1}{c_w(z)}\leqslant \liminf_{x\to z} \inf_{\bv\ne 0, \bv\in T_x \cF} \frac{||\bv||^w}{||\bv||} \leqslant 
\limsup_{x\to z} \sup_{\bv\ne 0, \bv\in T_x \cF} \frac{||\bv||^w}{||\bv||} \leqslant 
c_w(z).  
$$
Lemma follows. 
 \end{proof} 

\begin{rem}
In fact, Finsler metric $||\cdot||^w$ is not just measurable but is lower semicontinuous. 
\end{rem}

Our next goal is to prove the following exponential contraction theorem for an abelian subgroup of transformations.

\begin{thm}\label{thm:transvection contraction} 
Fix $\bx = \beta^2(\hat{m})$ with $\hat{m} \in \partial_\infty^2 \Ga \setminus \Delta$. 
and  $\eps>0$. Then for every $w\in W$ there exists a pair of continuous functions $C_w^\pm: \hat{\cO}_{w,\hat{m}}\to \R_+$ 
such that the following holds for every unit vector $\ba\in \mathfrak{a}^+_{\eps}$, every $z\in \hat{\cO}_{w,\hat{m}}$ and $t\geqslant 1$: 
$$
||de^{\mp t\ba}||_z^\pm \leqslant C_w^\pm(z) e^{-\eps t}.  
$$
Furthermore, $C_w^+$ is bounded on each tube $N_{\CC}(\Th_w(x_-))$ and $C_w^-$ is bounded on each tube $N_{\CC}(\Th_{w^\vee}(x_+))$. In particular, 
$C_w^\pm$ is bounded along each semi-orbit of $e^{\mp t\ba}$, $t\geqslant 0$. 
\end{thm}
\begin{proof} We will prove the claim for the action $e^{-t\ba}$ on the direction induced by the foliation $\hat{\cW}^+_w$ as the claim for  $e^{t\ba}$ follows by reversing the roles of $x_+, x_-$ and exchanging the roles of $\hat{\cW}^+_w$ with $\hat{\cW}^-_w$. 

Note that the contraction estimate is vacuous in the case $w=1$ since then there are no nonzero $w$-horizontal tangent vectors. 
(Tangent vectors to leaves of the horizontal foliation $\cW^+_{x_w}$ on $\cO_{x_w}$ are $w$-horizontal.)

We equip ${\mathfrak g}$ with a hermitian inner product with respect to which the root subalgebras ${\mathfrak g}_\al, \al\in \Phi$, 
are pairwise orthogonal. We let $||\cdot||_0$ denote the corresponding norm. 

\medskip 

We first prove the contraction estimate for the Finsler metric:
\begin{equation}\label{eq:contract} 
 ||d\ga (\bv)||^w\leqslant  e^{-\eps t} ||\bv||^w 
\end{equation}
for every $w\in W$ and every $w$-horizontal vector $\bv$ in $\hat{\cO}_{w,\hat{m}}$, where $\gamma = e^{-t \ba}$. 
The proof is by induction on $|w^\vee|=n-|w|$. 

\medskip 
{\bf Base of induction.} Suppose that $w=w_0, w^\vee=1$. Then $\hat{\cO}_{w,\hat{m}}=\Opp(x_+)$ and the horizontal foliation $ \hat{\cW}_{w,\hat{m}}^+$ consists of a single leaf, namely, $\hat{\cO}_{w,\hat{m}}$.

 Recall that the condition $\ba\in \mathfrak{a}^+_{\eps}$ means that for every positive root $\al\in \Phi^+$ we have
 $\al(\ba)\geqslant \eps$. But the numbers  $e^{\al(t\ba)}$ are the eigenvalues of $Ad(e^{t\ba})$ acting on ${\mathfrak u}$, $t\in \R$, where ${\mathfrak u}$ is the Lie algebra of the unipotent radical of the stabilizer of $x_+$.  Hence, for the norm $||\cdot||_w$, the contraction estimate at $x_w=x_-$ for a vector $\bv \in T_{(\hat{m}, x_-)} \hat{\cO}_{w, \hat{m}}$ becomes a computation on ${\mathfrak u}$, which is 
 $$
||de^{-t\ba}(\bv)||_0= ||de^{-t\ba}(\bv)||^w \leqslant e^{-\eps t} ||\bv||_0= e^{-\eps t} ||\bv||^w. 
 $$
 
For general points $z\in \hat{\cO}_{w, \hat{m}}$ we get a contraction estimate as 
follows. Let $\bv\in T_z \hat{\cO}_{w, \hat{m}}$, $\ga=e^{-t\ba}$, $\ba$ is a unit vector in $\mathfrak{a}_\eps^+$. Our first goal is to compare the norms 
$$
||\bv||_w, \quad ||d\ga (\bv)||_w. 
$$
Let $g, h\in U_{x_w}$ be such that $g(x_-)=z$, $h(\ga z)=x_-$. 

\begin{lem}
 $h\circ \ga \circ g=\ga$. 
\end{lem} 
\begin{proof} By the choice of $g, h$, the composition $\hat\ga:=h\circ \ga \circ g$ fixes the point $x_-$. Note that $\hat\ga$ belongs to 
$B_{w^\vee}= U_w\rtimes T$ (see \S~\ref{sec:Thickenings in the flag manifold}). Since $T$ fixes $x_w=x_-$ and $U_w$ acts simply transitively on $\cO_{x_w}$, the stabilizer of $x_w=x_-$ in 
$B_{w^\vee}$ equals $T$. Hence, $\hat\ga\in T$. Since $U_w$ is a normal subgroup of $B_{w^\vee}$, we also have 
$$
\hat\ga:=h\circ \ga \circ g= \ga \circ h^\ga \circ g, h^\ga \circ g\in U_w. 
$$
Since $\hat\ga\in T$, it follows that $\hat\ga=\ga$, as claimed. 
\end{proof}

Since the metric 
$||\cdot||_w$ is $U_{x_w}$-invariant, it follows that for $\bu= dg^{-1}(\bv)$, 
$$
||\bu||_0=||\bu ||_w= ||\bv||_w, \quad  ||d\ga (\bv)||_w= ||d\ga (\bu)||_0. 
$$

Hence, we obtain the same exponential contraction estimate 
$$
||de^{-t\ba}(\bv)||^w \leqslant e^{-\epsilon t}  ||\bv||^w 
 $$
 at all points $z\in \hat{\cO}_{w, \hat{m}}$ 	and all $\bv \in T_z \hat{\cO}_{w, \hat{m}}$. This concludes the proof of \eqref{eq:contract} in the case $w=w_0$.

\medskip 
{\bf Inductive step.} We assume that \eqref{eq:contract} holds for all $w\in W$ with $n-|w|<k$. Consider $w$ with $n-|w|=k$.  


Arguing as in the base of induction case, and taking into account the fact that action of $U_{x_w}$ 
preserves the vertical/horizontal foliations of $\hat{\cO}_{w,\hat{m}}$,   the contraction estimate 
\begin{equation}\label{eq:w-estimate} 
||d\ga(\bv)||_w\leqslant e^{-\eps t} ||\bv||_w,
\end{equation}
(for $w$-horizontal vectors $\bv\in T \hat{\cO}_{w,\hat{m}}$)  reduces to the case of the base-point $z= (\hat{m},x_w)$. But then the claim follows from the fact that the eigenvalues 
of $Ad(\ga)$ acting on ${\mathfrak u}^+_w$ are given by $e^{-t\al(\ba)}$ for $\al\in \Phi_w^+\subset \Phi^+$. By the assumption, all these eigenvalues are 
$\leqslant  e^{-\eps t}$, as required. 

In order to get an estimate in terms of $||\cdot||^w$, we have to consider two cases:  

{\bf Case 1.} Assume that $z\in \hat{\cO}_{w,\hat{m}}$ and $\bv\in T_z \hat{\cO}_{w,\hat{m}}$ are such that  
$$
||\bv||_v= ||\bv||^w, \quad ||d\ga(\bv)||_v= ||d\ga(\bv)||^w 
$$
for some $v\in W, v \geqslant w$. If $v = w$, then the inequality \eqref{eq:contract} follows from \eqref{eq:w-estimate}. Otherwise, $v> w$ and as $n- |v| < n-|w|$, the induction assumption can be applied to $v$ and we conclude.
For instance, this applies 
to all points $z\in \{ \hat{m} \} \times S_{w^\vee}(x_-)=\cL_{w,\hat{m}}^-$, since $\cL_{w, \hat{m}}^-$ is disjoint from 
$$
\bigcup_{v <  w} \hat{\cO}_{v, \hat{m}}
$$
(see Lemma \ref{lem:disjoint} with the sign reversed, exchanging $x_+$ with $x_-$ and $\mathcal{L}^+$ with $\cL^-$) and is $T$-invariant. Hence,  for all vectors $\bv$ tangent to $\hat{\cW}_w^+$ at points of $\cL_w^-$ we have
$$
||\bv||_w=||\bv||^w. 
$$

\medskip 
{\bf Case 2.} 
$$
||\bv||^v= ||\bv||_v= ||\bv||^w, \quad ||d\ga(\bv)||^w= ||d\ga(\bv)||_{v'} 
$$
for some $ v \neq  v'\in W$ satisfying $v\ge w, v'\ge w$.
Observe that this implies that $z \in \hat{\cO}_{wv,\hat{m}}\cap \hat{\cO}_{wv',\hat{m}}$ since $T$ preserve each domain $\hat{\cO}_{u}$ for $u \in \{ w,v,v'\}$. 
So that the vector $\bv$ belongs to $\hat{\cW}_{v,\hat{m}}^+ \cap \hat{\cW}_{w,\hat{m}}^+ \cap \hat{\cW}_{v', \hat{m}}^+$. 
If $v' > w$ then by the induction hypothesis and using the fact that $|| \cdot ||^{v'}\leqslant || \cdot ||^{w}$, we have: 
\begin{equation*}
|| d \gamma(\bv)||^w = || d\gamma (\bv)||^{v'} \leqslant e^{-\epsilon t} || v||^{v'} \leqslant e^{-\epsilon t} || v||^{w},  
\end{equation*}
as required. 
Otherwise, we have $v'=w$ and using \eqref{eq:w-estimate} 
we obtain: 
\begin{equation*}
|| d \gamma(\bv)||^{w} = || d\gamma(\bv) ||_w \leqslant e^{-\epsilon t} || \bv ||_w \leqslant e^{-\epsilon t} || \bv||^w, 
\end{equation*}
as required.   
This proves the inequality 
$$
 ||d\ga (\bv)||^w\leqslant  e^{-\eps t} ||\bv||^w  
$$
at all points $z\in \hat{\cO}_{w,\hat{m}}$. 

\medskip 
{\bf Estimating contraction in terms of the Riemannian metric on $\cF$.}  
The restriction of 
the Riemannian norm on $T\cF$ to any compact in $\cO^w_{\bx}$  is uniformly bilipschitz on compacts to the measurable Finsler metric  
$||\cdot||^w$ (cf. Lemma \ref{lem:comparison}). The issue is that even if $z\in \cO^w_{\bx}$
 lies in a given compact subset, we, a priori, do not know where its image under $\ga$ is. Consider 
 $z\in \cO_{x_w}$. Recall that earlier in \eqref{eq:pm-proj} we defined continuous maps $\pi_{w}^{\pm}: \cO_{x_w}\to \cL^\pm_{w,\bx}$. 
 Set $z_\pm:= \pi_{w}^{\pm}(z)$  and define 
 $$
 u_z:= \Psi_w(z_+)\in \mathfrak{u}_w^+\subset 
 \mathfrak{u}_w, \quad R=R(z):= ||u_z||_0. 
 $$
 Thus, $R(z)$ is a continuous function of $z$ and $z$ belongs to the $\CC$-tube $N_{\CC}(S_{w}(x_-))\subset \cO_{x_w}$, where 
$$
\CC=\bar{\mathbf{B}}(0, R)\subset \mathfrak{u}^+_w
$$ 
is the closed $R$-ball with respect to the norm $||\cdot||_0$ on $ \mathfrak{u}_w$. 
(See Definition \ref{def:tube} for the notion of $\CC$-tubes.) 
Clearly, $\CC$ is a compact convex subset of $\mathfrak{u}^+_w$ containing $0$. 

We proved in Lemma \ref{lem:compactifications of leaves} 
that the closure in $\cF$ of $N_{\CC}(S_{w}(x_-))$ equals the compact 

{$$
N_{\CC}(\{\hat{m} \} \times \Th_{w}(x_-))\subset \bigcup_{v\geq w^\vee}\hat{\cO}_{v, \hat{m}}=\hat{\cO}_{\hat{m}}^{w^\vee}.$$}
 
Moreover, in \S \ref{sec:tubes} we also observed that 
for every $\gamma\in \exp(- \mathfrak{a}^+)$, 
$$
\ga(N_{\CC}(\Th_{w}(x_-)))\subset N_{\CC}(\Th_{w}(x_-)).$$
Define the function $C_w^-(z)$ as follows.


Recall that in Lemma \ref{lem:comparison} we proved existence of a positive continuous function $c_{w^\vee}$ 
on $\hat{\cO}^{w^\vee}_{\hat{m}}$ providing bilipschitz comparison between the norms $||\cdot||$ and $||\cdot||^{w^\vee}$ on the tangent bundles $T\hat{\cO}_{\hat{m}}^{w^\vee}$. 
Then we set
$$
C_w^+(z):=\max_{x\in N_{\CC}(\Th_{w}(x_-))} c_{w^\vee}(x). 
$$
Continuity of the function $R(z)$ implies continuity of $C_w^+(z)$.
 
We, thus, obtain (using the inequality \eqref{eq:contract}) that for any $\bv \in T_z \hat{\cO}_{w,\hat{m}} \cap \hat{\cW}^+_{\hat{m}}$, one has:
\begin{equation*}
|| d (e^{-t \ba}) \bv ||_{e^{-t  \ba } \cdot z} \leqslant C_w^+(z) || d (e^{-t \ba}) \bv ||^w \leqslant C_w^+(z) e^{-\epsilon t} ||\bv ||^w \leqslant (C_{w}^+(z))^2 e^{-\epsilon t} ||\bv ||_{z}, 
\end{equation*}
Theorem follows. \end{proof}

We are now ready to prove

\begin{thm}\label{thm:main-C}
Suppose that $\Ga$ is a hyperbolic group, $\rho: \Ga\to G$ is an Anosov representation, $E\to \cG$ is the associated $\cF$-bundle and $\theta_t$ the lifted 
geodesic flow on $E$. Then $\theta_t$ is an exponentially contracting Morse--Smale flow.   
\end{thm} 
\begin{proof} 

Given $\hat{m} \in \partial^2_\infty \Ga \setminus \Delta$, recall from \S~\ref{sec:flat bundles} that the flow is defined by acting by the inverse of $\gamma_t$ where $\gamma_t = \Pi(\hat{m}_t)$. Since the action on $\cF$ is defined using the representation $\rho$, by Corollary \ref{cor:from group to transvections}, the distortion of the metric $|| \cdot ||$ by $\rho(\gamma_t^{-1})$ is equivalent to the distortion of the metric by the inverse of an element in $e^{t \mathfrak{a}_+}$. The previous result then gives:
$$
||d\theta_t||^{\pm}_z\leqslant C_w^{\pm}(z) e^{-\eps t}, t\geqslant 1
$$
This proves the exponential contraction of the flow. 
The Morse--Smale property was established in Lemma 
\ref{lem:Morse--Smale Property}. 
\end{proof}

\subsection{ Accumulation sets of $\theta$-semi-orbits}

\begin{defi}
Let $\{f_t: t\in \R_+\}$, 
be a family of maps of a metric space $(Z,d)$, and $A, B\subset Z$ be two subsets: $A$ is open and $B$ is closed. Then we say that the  
family $\{ f_t|_A: t\in \R_+\}$  {\em locally uniformly accumulates at $B$} if for every $r>0$ and every compact $C\subset A$, there exists $T<\infty$ such that 
$$
f_t(C)\subset N_r(B), \forall t\geqslant T. 
$$
In this case, we say that $\omega(C, f_t)\subset B$. 
\end{defi}

Recall from \S~\ref{sec:def_flow_section} that $\Sigma_w$ is the image of $\cG$ by the flow invariant section $s_w$ for each $w \in W$. 
\begin{lem}
$\theta_t|_{\cO_w}$ locally uniformly accumulates at $\Th_{w^\vee}(\Sigma_{w_0}) = \cup_{m \in \cG} \{ m\} \times \Th_{w^\vee}(s_{w_0}(m))$.  
\end{lem}

\begin{lem}
Given $z\in E$, there exists $w\in W$ such that $z\in \cL^+_w$. Furthermore, the $\omega$-limit set $\omega(z, \theta_t, t> 0)$ 
of  the flow $\theta_t$  is contained in $\Sigma_w$.  
\end{lem}
\begin{proof} Consider $z\in E_m\cong \cF$ for some $m\in \cG$. The flag manifold $\cF$ contains distinguished points $s_w(m), w\in W$, lying in a single apartment $a_{\bx}$, $\bx=(s_1(m), s_{w_0}(m))$. The flag manifold $\cF$ is the disjoint union of Schubert cells $S_w(x^+), w\in W$.  
Thus, $z\in S_w(x^+)\subset \cL^+_w$.  According to Theorem \ref{thm:main-C}, 
$$
\lim_{t\to\infty} d_{E_{m_t}}(\theta_t(x_w), \theta_t(z))=0. 
$$
Since $x_w\in \Sigma_w$ and $\Sigma_w$ is closed and flow-invariant, the $\omega$-limit set of the semi-orbit $(\theta_t(x_w), t >0)$, is contained in $\Sigma_w$. Thus, the same is true for the semi-orbit of $z$.  \end{proof}

For every thickening $\Th\subset W$ and $(x_+, x_-)\in \cY$, we have thickenings $\Th(x_\pm)$. These thickenings depend continuously on $x_\pm$ and, thus, produce laminations $\Th(\Sigma_{1})$, $\Th(\Sigma_{w_0})$ of $E$, whose leaves in each fiber $E_m$ are 
$$
\Th(s_1(m)), \Th(s_{w_0}(m)). 
$$
For every thickening $\Th\subset W$ one defines the {\em complementary thickening} $\Th^c\subset W$ (see \cite{KLP18A}) by 
\begin{equation}\label{eq:complementary} 
\Th^c:= w_0 (W\setminus \Th). 
\end{equation}

It was proven in \cite[Corollary 6.3]{KLP18A} that for every thickening 
$\Th\subset W$ the $\Ga$-orbits in $\cF \setminus \Th(\La)$ uniformly on compacts accumulate at 
$\Th^c(\La)$.

\begin{prop}
For every thickening $\Th\subset W$, the positive semi-flow $(\theta_t, t>0)$ locally uniformly accumulates at 
$\Th^c(\Sigma_{w_0})$ on $E\setminus \Th(\Sigma_1)$. Similarly, the negative semi-flow 
$(\theta_t, t< 0)$ locally uniformly accumulates at 
$\Th^c(\Sigma_{1})$ on $E\setminus \Th(\Sigma_{w_0})$. 
\end{prop}  
\begin{proof} We first analyze an easier case, proving 
locally uniform accumulation in a single fiber $E_m\cong \cF$. It suffices to work with a normalized element 
$\hat m\in \widehat E$ projecting to $m$. Set  
$\xi_\pm:= \re(\hat m)$, $x_\pm=\beta(\xi_\pm)$, where $\beta: \geo \Ga\to \La\subset \cF$ is the boundary map of $\rho: \Ga\to G$. 

Let $C\subset \widehat E_{\hat m}\cong \cF$ be a compact disjoint from  
$\Th(\hat s_{w_0}(\widehat\Ga))$. Pick $\eta>0$. We need to show that for all sufficiently large $t\in \R_+$, 
$$
\hat\theta_t(C)\subset N_\eta(\Th^c \hat s_{1}(\hat m_t) )$$
where the $\eta$-neighborhood is taken with respect to the $\Ga$-invariant fiberwise Riemannian metric on $\widehat E$. 
Since the metric is invariant, it suffices to work with the preimage of $E_{\hat m_t})$ under $\gamma_i^{-1}$, where $\gamma_i= \Pi(\hat m_t)\in \Ga$.


Then the claim amounts to saying that the sequence $\rho(\ga^{-1}_i)$ (restricted to a compact in $\cF \setminus \Th(x_+)$) 
accumulates  at some $\Th^c(y_-), y_-\in \La$ (with possibly varying points $y_-$). But this is what's proven in 
\cite[Corollary 6.3]{KLP18A}. More generally, consider a compact subset $C\subset \widehat E \setminus \Th(\hat{s}_{w_0}(\widehat \Ga))$. 
Then there is a compact family of elements $\{g_j\in G: j\in J\}$ such that for every $\hat m\in \hat\pi(C)\subset \widehat\Ga$, $g_j(f(\re(m))=(x_+,x_-)$ for 
two fixed antipodal points $x_\pm \in \cF$.  Take $h_{ij}:=g_j\circ \rho(\gamma^{-1}_i), j\in J, i\in \bN$, where $\gamma_i$'s are the sequences in $\Ga$ defined 
via the elements $\hat m_t$, $\hat m\in C$, $t\geqslant 0$, as in the   the easy case analyzed earlier in the proof. 
Then $\{h_{ij}: j\in J, i\in \bN\}$ is a uniformly regular family of elements of $G$ with the 
attractive/repelling points $y_-, x_+$ (the points $y_-\in \La$ again can be variable). 
We can now again apply \cite[Corollary 6.3]{KLP18A} to conclude that this family accumulates at  $\Th^c(y_-)$ on $\cF \setminus \Th(x_+)$. 
\end{proof}
 
\section{Examples}

\subsection{Convex-cocompact subgroups of $SL(2,\C)$}

Let $\Ga=\pi_1(\Sigma)$, fundamental group of a compact hyperbolic surface, $\Sigma=\bH^2/\Ga$. Let $UT\Sigma$ denote the unit tangent bundle of $\Sigma$, $UT\bH^2$ the unit tangent bundle of $\bH^2$. We have two commuting actions, $\Ga\acts UT\bH^2$ and 
$\hat\phi: \R\acts UT\bH^2$, the geodesic flow. 
The flow $\hat\phi_t$ projects to the geodesic flow $\phi_t$ on $UT\Sigma$.  
To every (unit speed) geodesic $\hat m=c: \R\to \bH^2$ we assign points $\re(\hat m)=\re(c)=(\xi_+, \xi_-)$, $\xi_\pm=c(\pm \infty)\in S^1=\geo \bH^2$. 

Let $\rho: \Ga\to G=SL(2,\C)$ be a quasifuchsian representation, $\beta: S^1=\geo \bH^2\to S^2$ the unique equivariant topological embedding; its image is the limit set $\La$ of $\rho(\Ga)$. 

Set $\widehat E:= UT\bH^2\times S^2$. The group $\Ga$ acts diagonally on $\widehat E$; set $E:= \widehat E/\Ga$. It is an $S^2$-bundle over $UT\Sigma$:
$$
\pi: E\to UT\Sigma. 
$$ 
It is naturally a {\em flat bundle} with Ehresmann connection $\nabla$: Horizontal vectors for this connection are projections of vectors tangent to $UT\bH^2\times \{p\}$, $p\in S^2$. Horizontal leaves of this connection are projections of subsets of the form 
$UT\bH^2\times \{p\}$; they form a foliation $\cE$ on $E$. 

We have two $\Ga$-equivariant maps  
$$
f_\pm: UT\bH^2\to S^2, f_\pm(c)= \beta(\re(c)). 
$$
These maps define sections $s_\pm$ of the bundle $E\to UT\Sigma$: 
$$
s_\pm(m)= [\hat m, f_\pm(\hat m)]\in E, 
$$
where $m\in UT\Sigma$, $\hat m\in UT\bH^2$ belongs to its preimage under the covering map $UT\bH^2\to \Sigma$.  
One can also label these sections $s_+=s_{1}, s_-=s_{w_0}$, where $1=1_W$, $w_0$ is the unique nontrivial element of the Weyl group $W$ of $S$. 

We have the lift $\hat\theta_t: (\hat m, z)\mapsto (\hat\phi_t(\hat m), z), t\in \R$, of the geodesic flow $\hat\phi_t$. The flow 
$\hat\theta_t$ descends to a flow $\theta_t$ on $E$, the {\em suspension flow} of the action $\rho: \Ga\acts S^2$. 

Sections $s_\pm$ define homeomorphisms 
$$
UT\Sigma\to s_\pm(UT\Sigma)\subset E. 
$$
The images of these sections are invariant under $\theta_t$ (since $f_\pm$ is constant along the trajectories of $\hat\phi_t$). The push-forward of the Bowen--Margulis measure on $UT\Sigma$ via $s_\pm$ defines ergodic flow-invariant measures on $s_\pm(UT\Sigma)$. 

Every fiber $E_m$ of the bundle $\pi: E\to UT\Sigma$ is $S^2$ with two distinguished points, $s_\pm(m)$ (North/South poles). 

Equip $E$ with a background Riemannian metric $g$ (it is {\em not} flow-invariant). Let $\cW^+$ denote the partial foliation of $E$; for every $m\in E$, 
$$
\cW^+_m= \cW^+\cap E_m
$$
is the once-punctured sphere, $E_m\setminus s_-(m)$. This foliation is the {\em stable foliation} of the flow $\theta_t$: The flow is (nonuniformly) exponentially contracting (with respect to the metric $g$) on the leaves of this foliation:
$$
||d\theta_t(\bv)||\leqslant C^+(z)e^{-\eps t}, \bv\in T_z E_m, t\in \R, 
$$ 
for some constant $\eps>0$ and some continuous function $C^+$ on the union $\cO_{w_0}$ of leaves of  $\cW^+$. 
Accordingly, for every $z\in \cO_{w_0}$, its $\theta_{t>0}$-trajectory accumulates at $s_+(UT\Sigma)$. The semiflow $\theta_{t>0}$ is proper on 
$\cO_{w_0}$. 

The semiflow is uniformly exponentially contracting at the invariant submanifold $s_+(UT\Sigma)$ and 
uniformly exponentially expanding at the invariant submanifold $s_+(UT\Sigma)$.

The subset
$$
\mathring\cL:= E\setminus (s_+(UT\Sigma)\cup s_-(UT\Sigma))
$$
is the {\em open separatrix} of the flow $\theta_t$, connecting $s_-(UT\Sigma)$ to $s_+(UT\Sigma)$. It is a $\C^*$-bundle over $UT\Sigma$. Fibers 
of the bundle are open Richardson varieties in $\P^1$. The action $\theta_t$ is proper on $\mathring\cL$. 

%
%
%
%
%

%
%
%
% ...
%
%
%
%
%
%
%
 

\appendix



\chapter[Appendix A]{
Redressing in hyperbolic groups}\label{sec:appendixA}

In this appendix we prove Theorem \ref{thm:redressing}. 

\medskip
The first ingredient that we will need is a {\em local-to-global} characterization of quasigeodesics in geodesic hyperbolic spaces.

\begin{defi}
Fix two nonnegative constants $C, D$. A concatenation $xy\star yz$ of two geodesic segments in a metric space $X$, 
is a {\em $(C,D)$-hinge} if $d(x,y)\geqslant D, d(y,z)\geqslant D$ and $d(y, xz)\leqslant C$. 
\end{defi}

The next proposition is well-known and easy to prove, but we could not locate it in the literature, hence, include a proof for the sake of completeness. 

\begin{prop}
[A local-to-global principle]\label{P1}
Suppose that $X$ is a $\delta$-hyperbolic geodesic metric space. Then for every $C$ there exists $D=D(C,\delta)$ satisfying the following. Consider a broken geodesic (arc-length parameterized) 
path 
$$
p= x_1x_2 \star x_2 x_3 \star ....
$$
in $X$, such that the concatenation of every pair of consecutive segments is a $(C,D)$-hinge. Then $p$ is an $(L, K)$-quasigeodesic for some $L, K$ depending only on $\delta$ and $C$. 
\end{prop}
\proof First, recall that a map $f: [a,b]\to X$ is a $k$-local $(L,A)$-quasigeodesic if the restriction of $f$ to each subsegment of length $\leqslant k$ is an $(L,K)$-quasigeodesic. It is known, see e.g. \cite{CDP}, that for all $(L,K, \delta)$, there exist $k=k(L,K,\delta)$ and $L'=L'(L,K,\delta), K'=A'(L,K,\delta)$, such that if $X$ is a $\delta$-hyperbolic geodesic metric space, then every $k$-local $(L,K)$-quasigeodesic in $X$ is $(L',K')$-quasigeodesic. Thus, we only have to verify that there exists $D=D(C,\delta)$ such that each $(C,D)$-hinge $xy\star yz$ (regarded as an arc-length parameterized path) is a $(1,A)$-quasigeodesic for some uniform $A$ depending only on $C$ and $\delta$. By the assumption, the distance from $y$ to $xz$ is $\leqslant C$. Thus, applying Lemma 1.80 in \cite{Kapovich-Sardar}, we obtain a monotonic (possibly discontinuous) map 
$$
f: xz\to xy\cup yz, 
$$
such that $f(x)=x, f(z)=z$ and $d(f, \id)\leqslant C_1=2(\delta+C)$. It follows that for all $u, v\in xy\cup yz$, the length of 
$xy\star yz$ between $u, v$ is at least
$$
d(u,v)-2C -4C_1= d(u,v) - 6C -8\delta. 
$$ 
Hence, the hinge $xy\star yz$ is $(1, 6C +8\delta)$-quasigeodesic. Taking $D=\frac{1}{2}k(1, 6C +8\delta, \delta)$, we obtain that the path $p$ is a $k$-local $(1, 6C +8\delta)$-quasigeodesic. Proposition follows. \qed


 
  Observe that, by compactness of $\ol{X}$, there exist finitely many shadows that cover the entire ideal boundary $\geo X$. From now on, we fix $x\in X$ and set  $X_R:= X\setminus B(x,R)$. Below we will be using the notion of {\em shadows} introduced in \S \ref{sec:Counting}.

\begin{prop}\label{P2} 
Given $x\in X, r\in [0,\infty)$, there exists a finite subset $\{z_1,...,z_n\}\subset X$ such that:

1. The shadows at infinity $\Sh^\infty(i):= \Sh^\infty_x(z_i,r)$ cover $\geo X$.  

2. Set $\Sh(i):= \Sh_x(z_i,r)$,  $i=1,...,n$.  
There exists $C<\infty$ such that for every $i\in [n]:=\{1,...,n\}$ there exists  
$\hat{i}\in [n]$ such that $\forall y_i\in \Sh(i)$, $\forall y_{\hat{i}}\in \Sh(\hat i)$,   
$$
(y_i, y_{\hat{i}})_x\leqslant C. 
$$
\end{prop}
\proof Note that, if $(z_k)$ is a sequence in $X$ diverging to infinity, the diameters (taken with respect to the metric $d_\infty$ on $\ol{X}$) of the shadows $\Sh_x(z_k,r)$ tend to zero. Suppose that $z_k^\pm$ are two sequences in $X$ converging to distinct points $\xi^\pm \in \geo X$. Then the sequence 
$(z_k^+, z_k^-)_x$ is bounded and, moreover, for any choice $y_k^\pm \in \Sh_x(z_k^\pm, r)$, the Gromov-products   
$(y_k^+, y_k^-)_x$ are also bounded. Now, lemma follows from the assumption that $\geo X$ contains at least two distinct points (since $\Ga$ is nonelementary). \qed 

\medskip
From now on, we fix the points $z_1,...,z_n$ as above and use the notation $\Sh(i), \Sh^\infty(i)$, introduced in the proposition.

\begin{cor}\label{cor:P2}
Let $C$ and $\{z_1,...,z_n\}\subset X$ be as in Proposition \ref{P2},  and let $D=D(C)$ be as in Proposition \ref{P1} 
and such that $X_D\subset \bigcup_{i\in [n]} \Sh(i)$. Then for all $i, j\in [n]$ (allowing for $i=j$) there exist 
$\gamma_{ij}\in \Gamma$  such that

1. $\gamma^{\pm 1}_{ij} x\in X_D$. 

2. $\gamma_{ij} x\in \Sh(\hat i)$, $\gamma^{-1}_{ij} x\in \Sh(\hat j)$. 

3. In particular, for all $y_i\in \Sh(i), y_j\in \Sh(j)$, the hinges  
$$
y_i x \star x\gamma_{ij}x, \quad y_j x \star x\gamma^{-1}_{ij}x
$$
satisfy the $(C,D)$-hinge condition. 
\end{cor}
\proof Recall that pairs of fixed points of loxodromic elements of $\Gamma$ are dense in  
$\geo X \times \geo X$, see Theorem \ref{thm:tukia} (also \cite[Corollary 8.2.G]{Gromov}, \cite[Theorem 2R]{Tukia} or \cite[Proposition 7.4.7]{DSU}). Thus, 
we find $g_{ij}\in \Ga$ such that for the attractive fixed point of $g_{ij}$ is in $\inte \Sh^\infty(\hat i)$ and the repelling fixed point is in $\inte \Sh^\infty(\hat j)$. 
Then take $\ga_{ij}$ to be a sufficiently high power of $g_{ij}$. This proves Parts 1 and 2. Part 3 is a consequence of Parts 1 and 2 combined with  Proposition \ref{P2}.  \qed

\medskip 
Define a finite subset ${\mathcal G}=\{\gamma_{ij}: i, j\in [n]\}$ in $\Gamma$. 

\begin{prop}\label{P3}
Consider $\gamma\in \Gamma$ and $x\in X$. Set $y_\pm=\ga^{\pm1}x$. Suppose that $d(x, y_+)\geqslant D$ (equivalently, $d(x,y_-)\geqslant D$), 
where $D=D(C)$ is as in Proposition \ref{P1}. 
 Then there exists $\gamma_{ij}\in {\mathcal G}$ such that the hinge  
$$
y_- x\star x y_+  
$$
belongs to a $(L,K)$-quasi-axis of $\tilde{\gamma}:= \gamma\circ \gamma_{ij}$ in $X$. 
\end{prop}
\begin{proof} 
Since $X_D\subset \bigcup_{i\in [n]} \Sh(i)$ and $d(x, y_\pm)\ge D$, there exist (possibly equal!) $i, j\in [n]$ such that 
$$
y_+\in \Sh(j), y_-\in \Sh(i). 
$$
Using Part 2 of Corollary \ref{cor:P2},  we find $\gamma_{ij}\in {\mathcal G}$ such that 
$$
x_-:=\gamma^{-1}_{ij} x\in \Sh(\hat j), x_+:=\gamma_{ij}x\in \Sh(\hat i). 
$$
By Part 3 of the same corollary, both concatenations $x_- x\star x y_+$ and $y_-x\star xx_+$ satisfy the $(C,D)$-hinge condition.   
Since $\gamma$ is an isometry, for $x'_+:= \ga x_+$, the concatenation 
$$
x y_+ \star y_+ x'_+= \gamma(y_-x\star xx_+)$$
 also satisfies  the $(C,D)$-hinge condition. See Figure \ref{fig:quasi-axis}.

\begin{figure}[htbp]
   \centering
   \includegraphics[scale=0.7]{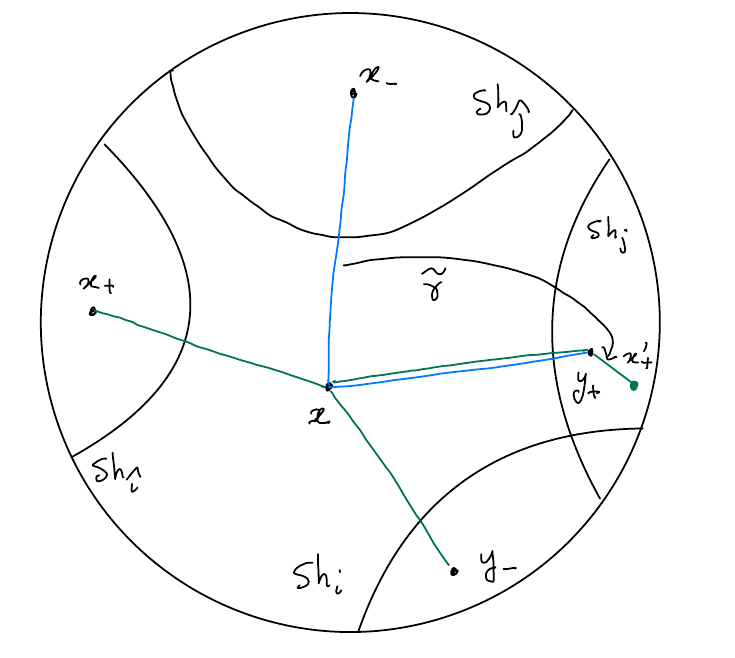} 
   \caption{Constructing a quasi-axis.}
   \label{fig:quasi-axis}
\end{figure}

Note that $\tilde{\gamma}(x_-x)= \gamma(xx_+)= y_+ x'_+$. Thus, every hinge in 
the $\tilde{\gamma}$-periodic piecewise-geodesic path   
$$
\alpha:= ...  \tilde{\gamma}^{-1}(x_- x \star xy_+) \star (x_- x\star x y_+) \star \tilde{\gamma}(x_- x \star xy_+) \star  \tilde{\gamma}^2(x_- x \star xy_+) \star ... 
$$
is congruent either to the hinge $x_- x\star x y_+$ or to the hinge $x y_+ \star y_+ x'_+$. Hence, all the hinges in $\alpha$ satisfy  the $(C,D)$-hinge condition.
By Proposition \ref{P1}, $\alpha$ is an $(L,A)$-quasigeodesic in $X$, hence, an $(L,A)$-quasi-axis of $\tilde{\gamma}$. \qed 

\medskip
We now can complete the proof of the Redressing Theorem.  We define constants $C$, $D=D(C)$ and $L, K$ as in Propositions \ \ref{P2},  \ref{P1} and \ref{P3}. Set $E:={\mathcal G}$, as above. There are only finitely many elements $\gamma\in \Gamma$ such that for $x=y_0$, 
$d(x, \gamma x)<D$, hence, we can safely ignore these. By Proposition \ref{P3}, for every remaining 
$\gamma\in \Gamma$, there exists $g\in E$ such that $\{x, g^{-1}(x), \gamma x\}$ lies on an 
$(L,K)$-quasi-axis $\alpha$ of $\tilde{\gamma}:= \gamma\circ g$. (We will be taking $g=\ga_{ij}$ using the notation of Proposition \ref{P3}.) 
The element $\tilde \gamma$ is loxodromic since it acts as a shift along the quasi-axis $\alpha$.  \end{proof}

\chapter[Appendix B]{
Estimates on fixed points of isometries of Gromov-hyperbolic spaces} \label{sec:appendixB}

In this appendix we prove a generalization of Theorem \ref{thm:fixed-point-estimate} by establishing bounds on fixed points (and approximate fixed points) of isometries of Gromov hyperbolic spaces. Informally speaking, the main result of this section (Theorem \ref{thm:appendixB}) 
 is that if for a loxodromic isometry $\ga$ of $X$ and a base-point $x_0\in X$, the point $\ga(x_0)$ is close to a point $\zeta\in \geo X$, then so is the attractive fixed point of $\ga$ in $\geo X$. the same holds for parabolic/elliptic 
isometries and their fixed points in $\geo X$. Elliptic isometries need not have fixed points in $\geo X$; we prove instead that if 
for an elliptic isometry  $\ga$ the point $\ga(x_0)$ is close to $\zeta$, then {\em approximate fixed points} of $\ga$ in $X$ are also close to $\zeta$.

\begin{defi}
Let $\ga$ be a self-map of a metric space $(X,d)$. A $C$-fixed point of $\ga$ is a point $x\in X$ such that $d(x, \ga(x))\leqslant C$. 
\end{defi}

Such a point $x$ is an {\em approximate} fixed point of $\ga$. Loxodromic and parabolic isometries of a Gromov hyperbolic space $X$ always have fixed points 
in $\geo X$; elliptic isometries need not have fixed points but they always have $C$-fixed points in $X$ with $C$ depending only on $\delta$. More precisely, 
according to \cite[Ch. 8, Proposition 24]{Ghys-book}, 

\begin{prop}
If $\ga$ is an isometry of a $\delta$-hyperbolic metric space $X$ is such that its displacement $d(\ga)$ satisfies 
$$
d(\ga):= \inf_{x\in X} d(x, \ga x)> 26\delta',
$$
then $\ga$ is loxodromic.
\end{prop}

In other words, every elliptic and every parabolic isometry of $X$ has a $C$-fixed point in $X$, with $C=26\delta'$. Of course, loxodromic isometries also can have 
arbitrary small (positive) displacements.  

In what follows, $X$ is a geodesic $\delta$-hyperbolic space with fixed base-point $x_0\in X$. 
The main result of this section is

\begin{thm}\label{thm:appendixB}
1. Suppose that $\ga$ is a loxodromic isometry of $X$ with the attractive fixed point $\xi^+\in \geo X$. Then 
$$
(\zeta, \xi^+)_{x_0}\geqslant \half (\zeta, \ga x_0)_{x_0} - 26\delta. 
$$

2.   Suppose that $\ga$ is an elliptic or parabolic isometry of $X$ fixing a point $\xi\in \geo X$. Then  
$$
(\zeta, \xi)_{x_0}\geqslant \half (\zeta, \ga x_0)_{x_0} - 50\delta. 
$$

3.  Suppose that $\ga$ is an  isometry of $X$ with a $C$-fixed point $y\in X$. Then   
$$
(y,\zeta)_{x_0} \geqslant  \half (\ga x_0, \zeta)_{x_0} - \left(\half C +  3\delta \right). 
$$
\end{thm}

While this theorem is not needed for the main body of the book, it is motivated by the proof of 
Theorem \ref{thm:fixed-point-estimate} in \S \ref{sec:projection} and it does not appear to be in the literature.

\section{Approximate fixed points}

Our first result deals with approximate fixed points in $X$ of all isometries, including loxodromic ones:

\begin{prop}\label{prop:C-fixed-etimate} 
Suppose that $\ga\in \Isom(X)$ is an isometry with a $C$-fixed point ${x}'_0$. 
Then for every $\zeta\in \geo X$,
$$
(x_0',\zeta)_{x_0} \geqslant  \half (\ga x_0, \zeta)_{x_0} - \left(\half C +  \delta'\right). 
$$
\end{prop}
\begin{proof}

\begin{lem}\label{lem:prod-dist} 
Consider a quadruple of points $p, p', q, q'\in X$ such that $d(p', q')\leqslant C$ and 
$d(p, p')= d(q, q')$. Then: 
$$
\min( (q, p')_p, d(p, p'))  \geqslant \half (d(p, q) - C). 
$$
\end{lem}
\begin{proof} First, 
\begin{align*}
(q, p')_p = \half( d(p, q) + d(p, p')  - d(p', q))\geqslant \\
 \half(d(p, q) + d(p, p')  - d(p, p')- C)= \half ( d(p, q) - C). 
\end{align*}
Secondly, by the triangle inequality and the assumptions of the lemma,
\begin{align*}
d(p,q)\leqslant d(p, p') + d(p', q') + d(q', q)\leqslant \\
d(p, p') + C + d(p, p')= 2 d(p, p') +C. 
\end{align*}
Lemma follows. \end{proof}

\medskip
Applying this lemma to a $C$-fixed point $x_0'$ of an isometry $\ga$ of $X$ (with $p=x_0, p'=x_0'$, $q=\ga x_0$, $q'=\ga(x_0')$) we obtain:

\begin{cor}\label{cor:prod-dist}
If $d(x_0', \ga(x_0'))\leqslant C$, then:  

1. $$
\min( (\ga x_0, x_0')_{x_0}, d(x_0, x_0'))  \geqslant \half (d(x_0, \ga x_0) - C). 
$$ 

2. $$
\min( (\ga x_0, x_0')_{x_0}, d(x_0, x_0'))  \geqslant \half ((\ga x_0, \zeta)_{x_0} - C). 
$$ 
\end{cor}
\begin{proof} The first part, as we noted above, is a consequence of  Lemma \ref{lem:prod-dist}. The second part follows from the first 
combined with Lemma  \ref{lem:0}:
$$
 d(\ga x_0,x_0)\geqslant (\ga x_0, \zeta)_{x_0}.  
$$
\end{proof}

We are now ready to conclude the proof of the proposition. 
We have to estimate from below $(x_0', \zeta)_{x_0}$ in terms of $(\ga x_0, \zeta)_{x_0}$ for points $\zeta\in \geo X$. 
Applying \eqref{eq:Gh1} to the points $x_0, x_0', \ga x_0$ and $\zeta$, we get 
$$
(x_0',\zeta)_{x_0} \geqslant \min( (x'_0, \ga x_0)_{x_0}, (\ga x_0, \zeta)_{x_0}) -\delta'.$$
Combining this inequality with Corollary \ref{cor:prod-dist}, we then obtain
\begin{align*}
(x_0',\zeta)_{x_0} \geqslant   \min \left( \half ((\ga x_0, \zeta)_{x_0} - C), (\ga x_0, \zeta)_{x_0}\right) -\delta' =\\ 
\half \left((\ga x_0, \zeta)_{x_0} - C\right) - \delta'. 
\end{align*}
This proves the proposition. \end{proof}

\medskip
For elliptic isometries, in general, this proposition is the best we can do, since they need not have genuine fixed points in $\ol X$. 
Assume now that $\ga$ fixes a point $\xi\in \geo X$ (and still has a $C$-fixed point $x_0'\in X$); this applies to all parabolic isometries with $C=26\delta'$. 
Consider a geodesic ray $x_0'\xi$ parameterized 
by $c: \R_+\to x_0\xi$ so that $c(0)=x_0'$. 

\begin{lem}
Then for every $t > C+2\delta$,
$$
d(c(t), \ga c(t))\leqslant C'= C+4\delta. 
$$
\end{lem}
\begin{proof} Consider the geodesic ray $c_1(t)=\ga c(t), t\in \R_+$. By $2\delta$-slimness of the ideal triangle 
equal to the union of geodesics
$$
c(0) c_1(0) \cup  c(\R_+) \cup c_1(\R_+),
$$
for every $t\geqslant 0$ there exists a point in 
$$
c(0) c_1(0) \cup c_1(\R_+)
$$
within distance $2\delta$ from $c(t)$. Since, by the assumption, $t > C+ 2\delta$, this point has to be some $c_1(t')\in c_1(\R_+)$. 
By the triangle inequality, $|t-t'|\leqslant C+2\delta$. 
Therefore, again by the triangle inequality, 
$$
d(c(t), \ga c(t))= d(c(t), c_1(t)) \leqslant C+4\delta.$$  
 \end{proof}

In other words, $c(t)$ is a $C'$-fixed point (provided that $t> C+2\delta$). 
Taking the lower limit of $(c(t), \zeta)_{x_0}$ as $t\to\infty$ and applying 
Proposition \ref{prop:C-fixed-etimate} to the point $c(t)$ (playing the role of $x_0'$) we obtain:
$$
(\xi, \zeta)_{x_0} + 2\delta' \geqslant \liminf_{t\to\infty} (c(t), \zeta)_{x_0} \geqslant \half (\ga x_0, \zeta)_{x_0} - (\half C' +\delta').
$$
Hence,
$$
(\xi, \zeta)_{x_0} \geqslant  \half (\ga x_0, \zeta)_{x_0} - (\half C + 2\delta + 3\delta')\geqslant  \half (\ga x_0, \zeta)_{x_0} - (\half C + 11\delta). 
$$

We obtain:

\begin{prop}\label{prop:small displacement case} 
Suppose that the displacement $d(\ga)$ is $\leqslant C$. 
Then for every $\xi\in \geo X$ fixed by $\ga$ and every  $\zeta\in \geo X$, we have
$$
(\xi, \zeta)_{x_0} \geqslant   \half (\ga x_0, \zeta)_{x_0} - (\half C + 11\delta). 
$$
\end{prop}

\begin{cor}
Suppose that $\ga$ is an elliptic or parabolic isometry. Then for every $\zeta\in \geo X$ and every $\xi\in \geo X$ fixed by $\ga$, 
$$
(\xi, \zeta)_{x_0} \geqslant   \half (\ga x_0, \zeta)_{x_0} - 50\delta. 
$$
\end{cor}

\begin{cor}
If $d(\ga)\leqslant C$ and $\ga x_0\in U(\zeta, t)$, then the fixed-point set of $\ga$ in $\geo X$ is contained in 
$U(\zeta, e^{\half C + 11\delta}\sqrt{t})$. 
\end{cor}

\section{Loxodromic fixed points}

We now consider loxodromic isometries $\ga$ with ``large'' (i.e. $> 26\delta'$) displacement. Actually, it is convenient to assume that $d(\ga)> 30\delta'$. 

\begin{prop}\label{prop:unbounded case} 
Suppose that $\ga$ is a loxodromic isometry of $X$ with attractive fixed point $\xi^+\in \geo X$ and satisfying $d(\ga)> C=30\delta'$. 
Then for every $\zeta\in \geo X$ we have  
$$
(\xi^+, \zeta)_{x_0} \geqslant \half (\ga x_0, \zeta)_{x_0} -20\delta. 
$$
\end{prop}
\begin{proof} Our argument in this case is similar to the one in \S \ref{sec:projection}. 
We let $A=A_\ga$  denote the axis of $\ga$ in $X$. As in \S \ref{sec:projection}, we define $\bar{x}_0:=P_A(x_0), x:=\ga(x_0)$ and 
$\bar{x}:=P_l(\ga \bar{x}_0)$, 
projection of $\bar{x}_0$ to a complete geodesic $l\subset A$ passing through $\bar{x}_0$. 
Then, by our assumption,
$$
d(\bar{x}_0, \ga \bar{x}_0)> 26\delta'. 
$$
Therefore,
$$
D= d(\bar{x}_0, \bar{x})\geqslant d(\bar{x}_0, \ga \bar{x}_0) - 2\delta > 30\delta'-2\delta \geqslant D_0=14\delta. 
$$
By Lemma \ref{lem:right move}, $\bar{x}\in l_+$, which is the subray $\bar{x}_0\xi^+\subset l$.  Since $d(P_A(x), \bar{x})\leqslant 2\delta$, we obtain
$$
\bar{x}\in P^\eps_l(x), \quad \eps=2\delta. 
$$

Our goal is to get a lower bound on $(x,\xi^+)_p$ in terms of $\half d(x_0, x)$ and a uniform additive constant. 
This will be done with the help of Lemma \ref{lem:pq-ine} below. We will need:

\begin{lem}\label{lem:B1}
If $p, q\in X$, $l\subset X$ is a geodesic, $\bar{p}\in P_l(p)$, $\bar{q}\in P^\eps_l(q)$ and $D=d(\bar{p}, \bar{q})> 4\delta$, then 
$$
d(p, \bar{p}) + d(\bar{p}, \bar{q}) + d(\bar{q}, q)-10\delta-\eps \leqslant d(p,q)\leqslant d(p, \bar{p}) + d(\bar{p}, \bar{q}) + d(\bar{q}, q).
$$
\end{lem}
\begin{proof} The second inequality is just the triangle inequality, we will prove the first one. 
Let $y\in \bar{p} \bar{q}$ denote the midpoint. Then, in view of the assumption that $d(\bar{p}, \bar{q})> 4\delta$ and the fact that
$$
d(p, \bar{p})\geqslant d(p, y), d(q, \bar{q})\geqslant d(q, y),
$$
and the $2\delta$-slimness of the quadrilateral with the vertices $p, \bar{p}, \bar{q}, q$, there exists $z\in pq$ within distance $2\delta$ from $y$.  
By Lemma \ref{lem:1},
$$
d(p, y)\geqslant d(p, \bar{p})+ \half D - 3\delta, d(q, y)\geqslant d(q, \bar{q})+ \half D - 3\delta -\eps,
$$
hence,
\begin{align*}
d(p,q)= d(p,z)+ d(z,q)\geqslant d(p, y) + d(y,q) -4\delta \geqslant \\ 
(d(p, \bar{p})+ \half D - 3\delta) + (d(q, \bar{q})+ \half D - 3\delta-\eps) -4\delta\\ 
= d(p,q)\geqslant d(p, \bar{p}) + d(\bar{p}, \bar{q}) + d(\bar{q}, q) -10\delta-\eps. 
\end{align*}
\end{proof}

We combine this inequality with the estimate from Corollary \ref{cor:1} and obtain, in the setting of Lemma \ref{lem:B1}:

\begin{lem}\label{lem:pq-ine}
Suppose, in addition to the assumptions of Lemma \ref{lem:B1}, 
that $l$ is a complete geodesic connecting points $\xi^\pm\in \geo X$ and 
\begin{equation}\label{eq:pq-proj}
d(p, \bar{p})\leqslant d(q, \bar{q})\leqslant d(p, \bar{p}) + \eps. 
\end{equation}
Then 
$$
(q,\xi^+)_p\geqslant 
 \half d(p,q) - \frac{3}{2}\eps - 8\delta-2\delta'\geqslant   \half d(p,q) - \frac{3}{2}\eps - 14\delta. 
$$
\end{lem}
\begin{proof} 
By Lemma \ref{lem:B1},
$$
d(p, \bar{p}) + d(\bar{p}, \bar{q}) + d(\bar{q}, q) -(\eps+10\delta)\leqslant d(p,q).
$$
Combining this with the inequality $d(p, \bar{p})\geqslant d(q, \bar{q})-  \eps$ (an assumption of the corollary), we get
 $$
2 d(q, \bar{q}) -(\eps+10\delta)\leqslant  d(p,q) 
 $$
 and, thus,
 \begin{equation}\label{eq: half}
 d(q, \bar{q}) \leqslant \half d(p,q) + \half (\eps+10\delta). 
 \end{equation}
According to  Corollary \ref{cor:1},
\begin{equation*}
(q,\xi^+)_p\geqslant d(p,q) - d(q, \bar{q}) - \eps' - 2\delta'. 
\end{equation*}
Combining this with the estimate \eqref{eq: half}, we get
$$
(q,\xi^+)_p\geqslant \half d(p,q) - \half (\eps+10\delta) - (\eps' + 2\delta')
$$
Since $\eps'=\eps+ 3\delta$ (see \eqref{eq:projection0}), lemma follows. 
\end{proof}

We apply this lemma to estimate $(x,\xi^+)_{x_0}$. We will be taking $p=x_0, \bar{p}=\bar{x}, q=x, \bar{q}=\bar{x}$ and 
$l\subset A=A_\ga$ a geodesic through $\bar{x}_0, \bar{x}$ as above. Since 
$\bar{x}\in l_+$, $\bar{x}\in P^\eps_l(x)$ with $\eps=2\delta$, we are in position to apply Lemma \ref{lem:pq-ine} and obtain:
\begin{equation}\label{eq:distance_to_xi}
(x, \xi^+)_p\geqslant  \half d(x_0, x) - 3\delta - 8\delta-2\delta'\geqslant   \half d(x_0, \ga x_0) - 17\delta. 
\end{equation}
Since (according to Lemma  \ref{lem:0}),
$$
d(x_0, \ga x_0)\geqslant (\ga x_0, \zeta)_{x_0},  
$$
we get 
$$
(\ga x_0, \xi^+)_{x_0}\geqslant  \half  (\ga x_0, \zeta)_{x_0} -17\delta.   
$$
The rest of the proof is the same as in Proposition \ref{prop:C-fixed-etimate}: 

Applying \eqref{eq:Gh1} to the points $x_0, \xi^+, \ga x_0$ and $\zeta$, we get 
\begin{align*}
(\xi^+, \zeta)_{x_0} \geqslant \min( (\xi^+, \ga x_0)_{x_0}, (\ga x_0, \zeta)_{x_0}) -\delta'= \\
\half (\ga x_0, \zeta)_{x_0} -20\delta.  
\end{align*}
This concludes the proof of Proposition \ref{prop:unbounded case}.  \end{proof} 

\medskip
We will now finish the proof of the theorem. We already proved the claims for of parabolic and elliptic isometries, as well as loxodromic isometries with 
displacement bounded from above. Suppose that $\ga$ is a general loxodromic isometry. 

\medskip 
{\bf Case 1.} Suppose that $d(\ga)> C=30\delta$. Then, by Proposition \ref{prop:unbounded case}, 
$$
(\xi^+, \zeta)_{x_0} \geqslant \half (\ga x_0, \zeta)_{x_0} -20\delta. 
$$

{\bf Case 2.}  Suppose that $d(\ga)\leqslant C=30\delta$. We then can use Proposition \ref{prop:small displacement case}:
\begin{align*}
(\xi, \zeta)_{x_0} \geqslant   \half (\ga x_0, \zeta)_{x_0} - (\half C + 11\delta)=  \\
  \half (\ga x_0, \zeta)_{x_0} - (15\delta + 11\delta)= 
 \half (\ga x_0, \zeta)_{x_0} - 26\delta. 
\end{align*}

Combining the two cases, we obtain
$$
(\xi, \zeta)_{x_0} \geqslant \half (\ga x_0, \zeta)_{x_0} - 26\delta. 
$$
This concludes the proof of Theorem \ref{thm:appendixB}. 

\begin{cor}
Suppose $\ga$ is a loxodromic isometry with the attractive fixed point $\xi^+$. Then for every $\zeta\in \geo X$ and $t>0$, if $\ga(x_0)\in U(\zeta, t)$, then 
$$
\xi^+\in U(\zeta, e^{26\delta}\sqrt{t}). 
$$  
\end{cor}

\chapter[Appendix C]{
Bi-H\"older property of boundary maps} 
\label{sec:holder}

Let $G$ be a complex semisimple Lie group, $\cX=G/K$ the associated symmetric space of noncompact type, $\cF=G/B$ the corresponding full flag-manifold, $\Ga$ a hyperbolic group, $\rho:  \Ga\to G$ an Anosov representation, $f: \geo \Ga\to \La=\La(\Ga) \subset \cF$ the boundary map of $\rho$. As usual, we equip $\geo \Ga$ with a visual metric $d_\infty$ and $\cF$ with a $K$-invariant K\"ahler metric. As we noted in \S \ref{sec:def anosov}, it was proven by Tsouvalas in \cite{Tsouvalas} that the map $f$ is H\"older. However, in the paper Tsouvalas does not prove that $f^{-1}$ is also H\"older (he proves this result only for some special classes of Anosov representations, namely, ones which are $(1,1,2)$-{\em hyperconvex}). The main goal of this section is to prove a stronger form of the H\"older property of  $f^{-1}$:  

\begin{thm}\label{thm:biholder}
There exist positive constants $a$ and  $C$ (depending on $\rho$) such that for all points $\zeta_\pm\in \geo \Ga$  we have
$$
d_{\cF}(f(\zeta_-), \Th^{n-1}(f(\zeta_+)))\geqslant C d^a_\infty(\zeta_-, \zeta_+). 
$$
More generally, for any thickening $\Th\subset W$ and its complementary thickening (see \eqref{eq:complementary}) $\Th^c$, 
there exist positive constants $a$ and  $C$ (depending on $\rho$) such that for all points $\zeta_\pm\in \geo \Ga$  we have
$$
d_{\cF}(\Th(f(\zeta_-)), \Th^{c}(f(\zeta_+)))\geqslant C d^a_\infty(\zeta_-, \zeta_+).
$$
\end{thm}

Combining this theorem with \cite[Theorem 1.6]{Tsouvalas}, we obtain: 

\begin{cor}\label{cor:biholder}
The map $f$ is bi-H\"older. 
\end{cor}

The key to the proof is the following purely geometric proposition:

\begin{prop}\label{prop:holder}
Let $g\in G$ be a transvection along a geodesic $c$ through the origin $o\in \cX$, with displacement $D=d(o, g(o))$. Then $g: \cF\to \cF$ is 
$L$-bi-Lipschitz with 
$$
L= L_0 e^{a_0 D}, 
$$ 
for some constants $L_0>0, a_0\geqslant 0$ depending only on $\cX$. 
\end{prop}
\begin{proof} We let $a_0$ denote 
$$
\max_{\alpha\in \Phi, ||\ba||=1} \langle \alpha, \ba\rangle, 
$$
i.e. the maximum of norms of the linear functionals $\alpha\in \Phi$ on the normed vector space $\mathfrak{a}$. For an irreducible symmetric space (equivalently, irreducible root system), there is only one or two $W$-orbits in $\Phi$ (simply-laced versus non-simply-laced root system), hence, the maximum is actually a maximum of one or two different values.

Since $g^{-1}$ is also a transvection along the same geodesic $c$ with the same displacement $D$, it suffices to bound the Lipschitz constant of $g$. It suffices to consider the face of transvections which belong to the maximal torus $T$ stabilizing the model flat $F=\Fmod$. Let $x_\pm$ denote the chambers in $\amod=\geo \Fmod$ which are attractive/repelling fixed points of $g$ in $\cF$; set $\bx=(x_+,x_-)$. 
Recall that in \S \ref{sec:Contraction estimates} for every $w\in W$ we defined an open Schubert cell $\cO_{x_w}=\Opp(x_{w^\vee})\subset \cX$; these open subsets cover $\cF$. For every $w\in W$ we defined a Riemannian metric $||\cdot||_w$ on   $\cO_{x_w}$; we also defined measurable Finsler metrics $||\cdot||^w$ 
on
$$
\cO_{\bx}^{w}:= \bigcup_{v\ge w} \cO_{x_v}.
$$
Set $||\cdot||:= ||\cdot||^w$ for $w=1\in W$; this is a measurable Finsler metric on $\cO_{\bx}^{w}=\cF$. According to Lemma \ref{lem:comparison}, this metric is uniformly bi-Lipschitz to the background K\"ahler metric on $\cF$. Thus, it suffices to bound the Lipschitz constant of $g$ with respect to $||\cdot||$. 
The proof follows the one of Theorem \ref{thm:transvection contraction}, more specifically, of the inequality \eqref{eq:contract}. 

Without loss of generality (after changing the $W$-invariant metric on $\mathfrak{a}$ if necessary), 
we may assume that $g=e^{-t\ba}$, $t=D\geqslant 0$ where $\ba\in \Delta$ is a unit vector (this is a bit different from the proof of 
Theorem \ref{thm:transvection contraction} where we had $\ba\in \mathfrak{a}^+$, the interior of $\Delta$). 

\medskip 
{\bf Base of induction.} Suppose that $w=w_0, |w|=n$. Then the same proof as in Theorem \ref{thm:transvection contraction} shows that 
$$
||de^{-t\ba}(\bv)||_w= ||de^{-t\ba}(\bv)||^w \leqslant  ||\bv||^w =  ||\bv||_w
 $$
 for all tangent vectors $\bv$ at all points $z\in \cO_{x_w}$. (The map $g$ restricted to $\cO_{x_w}$ is $1$-Lipschitz with respect to the metric 
 $||\cdot||^w$.)

\medskip 
{\bf Inductive step.} We consider $w\in W$ with $|w|=n-k$ and assume that 
the inequality 
\begin{equation}\label{eq:contract1} 
 ||dg (\bv)||^v\leqslant  e^{a_{0} t} ||\bv||^v 
\end{equation}
holds for all $v\in W$ with $|v|> k$ and vectors $\bv\in T\cO_{\bx}^v$.  As in the proof of Theorem \ref{thm:transvection contraction}, the Lipschitz estimate for 
the norm $||\cdot||_w$ on $ \cO_{x_w}$ reduces to the case of vectors in the tangent space  $T_z\cO_{x_w}$, $z=x_w$, i.e. estimating from above 
the eigenvalues   of $Ad(g)$ acting on the Lie algebra $\mathfrak{u}_w$. These eigenvalues are of the form $e^{-t\al(\ba)}$ for 
$\al\in \Phi$. Since $\ba$ is a unit vector, these eigenvalues are bounded from above by  
$ e^{a_{0} t}$. 

In order to get an estimate in terms of $||\cdot||^w$, we consider two cases:  

\medskip 
{\bf Case 1.} Assume that $z\in \cO_{x_w}$ and $\bv\in T_z\cF$ are such that  
$$
||\bv||_v= ||\bv||^w, \quad ||dg(\bv)||_v= ||dg(\bv)||^w 
$$
for some $v\in W$, $v\ge w$. Then the inequality
\begin{equation}\label{eq:contract2} 
||dg(\bv)||^w= ||dg (\bv)||_v\leqslant  e^{a_{0} t} ||\bv||_v =||\bv||^w
\end{equation}
follows from the induction assumption and the inequality
$$
 ||dg (\bv)||_w\leqslant  e^{a_{0} t} ||\bv||_w. 
$$

\medskip 
{\bf Case 2.} 
$$
||\bv||^v= ||\bv||_v= ||\bv||^w, \quad ||dg(\bv)||^w= ||dg(\bv)||_{v'} 
$$
for some $v, v'\in W$ satisfying $v\ge w, v'\ge w$. The proof of the inequality \eqref{eq:contract2} is exactly the same as in Theorem \ref{thm:transvection contraction} (Case 2), except there is no need to restrict to  horizontal vectors; therefore, we omit the proof. 
\end{proof}

We can now prove Theorem \ref{thm:biholder}. 

\begin{proof}
First of all, recall that the group $K$ acts transitively on the set of maximal flats passing through the base-point $o\in \cX$. Therefore, for every thickening $\Th\subset W$ there exists $\eps_0>0$ 
(depending on $\cX$ and the choice of the $K$-invariant metric on $\cF$) such that for any maximal flat $F\subset \cX$ through $o$ and pair of opposite chambers $\sigma_\pm\subset \geo F$, 
$$
d_{\cF}(\Th(\sigma_-), \Th^{c}(\sigma_+))=\eps_0. 
$$ 
Suppose that $F'\subset \cX$ is a maximal flat containing a point $o'$ at the distance $D$ from $o$ and let $\sigma'_\pm\subset \geo F'$ be a pair of opposite chambers. Let $g\in G$ be a transvection along the geodesic through $o, o'$ sending $o$ to $o'$. Then, according to Proposition \ref{prop:holder}, $g: \cF\to \cF$ is $L= L_0e^{a_0 D}$-bi-Lipschitz.    

Therefore, 
$$
d_\cF(\Th(\sigma'_-), \Th^c(\sigma_+))\geqslant L^{-1}\eps_0. 
$$
Suppose now that $\sigma'_\pm=f(\zeta_\pm)$ for some (necessarily distinct) points $\zeta_\pm\in \geo \Ga$. We assume that the visual metric $d_\infty$ on $\geo \Ga$ comes from a geometric $\Ga$-action on a $\delta$-hyperbolic geodesic metric space $X$ and 
\begin{equation}\label{eq:visual1} 
k_1\exp(-\alpha (\zeta_-,\zeta_+)_{x_0})\leqslant d_{\infty}(\zeta_-,\zeta_+)\leqslant k_2\exp(-\alpha (\zeta_-,\zeta_+)_{x_0}),
\end{equation}
with uniform positive constants $\alpha, k_1, k_2$. Here $x_0$ is  a base-point in $X$. Thus, denoting $c$ a geodesic in $X$ asymptotic to 
$\zeta_\pm$, we get
$$
d(x_0,c) \leqslant  (\zeta_-,\zeta_+)_{x_0}\leqslant  d(x_0, c)+\delta. 
$$
We have an equivariant quasiisometric embedding $q: \Ga x_0\to \rho(\Ga)o\subset X$.  This map sends the geodesic $c$ to a (uniform) Morse quasigeodesic in $\cX$ and sends $x_0$ to $o$. Morse quasigeodesic $q(c)$ is uniformly close to the flat $F'$. Therefore, 
$$
D=d(o, F')\leqslant  k_3 (\zeta_-,\zeta_+)_{x_0} + a_3 
$$ 
for some constants $k_3, a_3$ depending only on $\rho$, $X$, $x_0$ and $o$. Combining this with inequalities \eqref{eq:visual1}, we get:
$$
d(o, F')\leqslant  - k_3 \alpha^{-1} \log (k_1^{-1} d_{\infty}(\zeta_-,\zeta_+)) + a_3. 
$$ 
Thus,
\begin{align*}
d_\cF(\Th(\sigma'_-), \Th^c(\sigma'_+))\geqslant \eps_0 L^{-1}= \eps_0 L^{-1}_0 e^{-a_0 D}\geqslant \\ 
\eps_0 L^{-1}_0 \exp\left(k_3 a_0\alpha^{-1} \log (k_1^{-1} d_{\infty}(\zeta_-,\zeta_+)) - a_0a_3\right)= \\
C \exp(a \log (d_{\infty}(\zeta_-,\zeta_+)) ) 
\end{align*}
for $a:= k_3 a_0 \alpha^{-1}$ and 
$$
C= \eps_0 L^{-1}_0 e^{-a_0(a_3+ k_3\log(k_1)/\alpha)}. 
$$
  Thus,
$$
d_\cF(\Th(\sigma'_-), \Th^c(\sigma'_+))\geqslant C d^a_{\infty}(\zeta_-,\zeta_+),
$$
as claimed. \end{proof}

\begin{rem}
Theorem \ref{thm:biholder} (and its corollary) also holds for $\taumod$-Anosov representations to general semisimple Lie groups 
(with $\Th^{n-1}$ replaced by the maximal proper thickening in $W/W_{\taumod}$). With minor modifications (e.g. replacing maximal flats with parallel 
sets of geodesics asymptotic to a pair of antipodal simplices of type $\taumod$), the same proof as above goes through. 
\end{rem}

%
%
%
%
%


\chapter[Appendix D]{
Measure and geometric measure theory}
\label{sec:appendix D}


\section{Push-forward and pull-back of measures}\label{sec:measures}


Recall that the {\em push-forward} \index{push-forward of a measure}
of a Borelean measure $\mu$ on a topological space $X$ via a measurable map $\ga: X\to X$ is a Borel measure $\nu=\ga_*(\mu)$ on $X$ satisfying 
$$
\nu(A)=\mu(\ga^{-1}(A)). 
$$
Equivalently, for every measurable function $h$ on $X$,
$$
\int_X h d\nu= \int_X h\circ \ga d\mu. 
$$
If $\ga$ is invertible with measurable inverse then the measure $\ga^*\mu=(\ga^{-1}_*)(\mu)$ is the {\em pull-back} \index{pull-back of a measure}
of the measure 
$\mu$. Thus, 
$$
\ga^*\mu(A)= \mu(\ga(A)) 
$$
for every measurable subset $A$. We have
$$
(\ga_1 \ga_2)^*\mu= \ga_2^*(\ga_1^*\mu) 
$$
and, accordingly (the Chain Rule),
$$
\left.\frac{d(\ga_1 \ga_2)^*\mu}{d\mu}\right\vert_z= \left. \frac{d(\ga_2^* \ga_1^*\mu)}{d\ga_1^*\mu}\right\vert_{\ga_1(z)} \cdot 
\left. \frac{d \ga_1^*\mu}{d\mu}\right\vert_{z}. 
$$

More generally, one defines pull-back of measures via finite-to-one maps. 
Namely, suppose that $f: X\to Y$ is a finite-to-one map between two sets. 
The pull-back $f^*\mu$ of a measure $\mu$ on $Y$ is a measure $\nu$ on $X$ such that  
$$
\nu(A)=\int_{Y} \# (f^{-1}(y)\cap A) d\mu(y). 
$$
(The measure $\nu$ is defined on subsets $A$ for which the above integral exists.) If $h$ is a measurable function on $Y$, then the pull-back $\nu=f^*(h\mu)$ satisfies
$$
d\nu(x)= h(f(x))d f^*(\mu)(x). 
$$

This pull-back construction using the counting measure on the fibers of $f$ has a {\em weighted (or orbi) analogue}. Suppose that $X, Y$ are locally compact metrizable 
topological spaces and the measure $\mu$ is a Radon measure on $Y$. Assume that $f$ is the quotient map with respect to a properly discontinuous 
action $\Ga\times X\to X$. For every $x\in X$ we define the {\em weighted counting} measure $m$ on the orbit $\Ga x$, multiplying the counting measure  by $|\Ga_x|$ (which is finite). Then for a subset $Z\subset X$ we define the {\em weighted pull-back} $\nu=f^\bullet \mu$ of $\mu$ under $f|_Z$ by
$$
\nu=\int_Y m(f^{-1}(y)\cap A)d\mu(y), 
$$
$A\subset Z$ is a Borelean subset. 

\medskip
We will also need a construction of a projection of a Radon measure under a properly discontinuous group action; it is a variation on the construction given in \cite[\S 2.6]{PPS}. 
Let $X$ be a locally compact 2nd countable metrizable space equipped with a properly discontinuous action $\Ga\times X\to X$ and a $\Ga$-invariant Radon measure $\mu$. Let $f: X\to Y:=X/\Ga$ denote the quotient map. 

\begin{prop}
There exists a (unique) Radon measure $\nu=f_\bullet \mu$ on $Y$ such that for every relatively compact Borelean subset $K\subset X$, the restriction 
$f|_K: K\to L=f(K)\subset Y$ has the property that $f^\bullet \nu|_L= \mu|_K$. 
\end{prop}
\begin{proof} Recall that for a set $Z$, a subset $\cS$ of $2^Z$ is called a {\em semi-ring}  if it satisfies the following properties:

\begin{enumerate}
\item $\emptyset \in \cS$. 

\item If $A, B\in \cS$ then $A\cap B\in \cS$. 

\item If $A, B\in \cS$ then there are finitely many pairwise disjoint elements $K_0,...,K_n\in \cS$ such that
$$
A\setminus B= \coprod_{i=0}^n K_i. 
$$
\end{enumerate}

A subset $W\subset X$ is said to be {\em precisely invariant} under the $\Ga$-action, if for every $\ga\in \Ga$ either $\ga W\cap W= \emptyset$ or $\ga W=W$. We let 
$\cS$ denote the subset of $2^X$ consisting of all relatively compact precisely invariant Borelean subsets of $X$. Proper discontinuity implies that the 
$\Ga$-stabilizer of every $W$ is finite and  that every point $x\in X$ has a basis of topology consisting of precisely invariant neighborhoods. 

\begin{lem}
$\cS$ is a semi-ring. 
\end{lem}
\begin{proof} Properties 1 and 2 of a semi-ring are clear. Let's prove the last property. Suppose that $A, B\in \cS$. Then $A\setminus B$ is clearly locally 
compact. Let $B_0=B, B_1,...,B_n$ denote the subsets of the form $\ga B, \ga\in \Ga$ which have nonempty intersection with $A$. (There are only finitely many of such 
subsets due to our relative compactness assumption and proper discontinuity of the $\Ga$-action.) Precise invariance of $B$ implies that the subsets $B_i$ are 
pairwise disjoint.  Set $K_i:= A\cap B_i, i=1,...,n$, and let 
$$
K_0:= A\setminus \bigcup_{i=1}^n B_i. 
$$  
Then each $K_i, i=0,...,n$, belongs to $\cS$ and 
$$
A\setminus B= \coprod_{i=0}^n K_i. 
$$
\end{proof}

We leave a proof of the next lemma to the reader:

\begin{lem}
1. The set of projections $\{f(W): W\in \cS\}$ again forms a semi-ring $\cT$ in $2^Y$.  

2. Each point $y\in Y$ has a basis of topology consisting of elements of $\cT$. 
\end{lem}

Part 2 of this lemma implies that the $\sigma$-algebra generated by the semi-ring   $\cT$ is the 
entire $\sigma$-algebra of Borelean subsets of $Y$. We define a real-valued set-function $\nu$ on $\cT$ by 
$$
\nu(f(W)):= \frac{1}{|\Ga_W|}\mu(W), \quad W\in \cS,$$
where $\Ga_W$ is the stabilizer of $W$ in $\Ga$. (Finiteness of $\mu(W)$ follows from the relative compactness of $W$.) 

The next lemma is again straightforward:

\begin{lem}
$\nu$ is a $\sigma$-finite premeasure on $\cT$. 
\end{lem}

Thus, by the Caratheodory Extension Theorem, the   premeasure $\nu$ extends (uniquely) to a $\sigma$-finite measure (again denoted $\nu$) on the $\sigma$-algebra of 
Borelean subsets of $Y$. To check that $\nu$ satisfies the pull-back property stated in the proposition note that this property holds for the premeasure $\nu$ on the 
semi-ring $\cT$ (by the construction of $\nu$). Every compact $C\subset Y$ admits a finite cover by elements of $\cT$, therefore, $\nu(C)<\infty$. 
It follows that the measure $\nu$ is Radon. 
\end{proof}

\medskip 
We will refer to the measure $\nu=f_\bullet \mu$ as the {\em projection} of $\mu$ under the quotient map $f$. \index{Projection of a measure under a quotient map}

\section{Ergodic aspects of group actions}

Let $(X,\mu)$ be a  measure space, where $\mu$ is a probability measure on $X$.  
We will assume that $\mu$ is {\em non-atomic} (has no atoms) i.e. every subset $A\subset X$ of positive measure contains a measurable subset $B\subset A$ such that 
$$
0< \mu(B)< \mu(A). 
$$ 
It follows that $A$ contains subsets of arbitrarily small positive measure. (Even more is true: Sierpinski proved that for every $m\in (0, \mu(A))$ there exists 
a measurable subset $B\subset A$ such that $\mu(B)=m$.)

Let $G$ be a countable group and $G\times X\to X$ be a measurable $G$-action. We will assume that $G$ preserves the measure class of $\mu$. (We will 
not require $G$-invariance of $\mu$.) The $G$-action on $(X,\mu)$  is said to be {\em ergodic} if one of the two equivalent conditions holds:

1. If $X=Y\sqcup Z$ is a measurable $G$-invariant partition of $X$, then $\mu(Y)\mu(Z)=0$. 

2. If $f\in L^\infty(X)$ is a $G$-invariant function, then $f$ is constant almost everywhere. 

The $G$-action on $(X,\mu)$ is said to be {\em dissipative} if there exists a measurable subset $D\subset X$ of positive measure such that for every 
$g\in G\setminus \{1\}$, $\mu(gD\cap D)=0$. Note that an ergodic action cannot be dissipative. Indeed, since $\mu$ has no atoms, there exists a subset $A\subset D$ 
whose measure satisfies $0<\mu(A)<\mu(D)$. Then both $GA\subset X$ and its complement have positive measure and are $G$-invariant, contradicting ergodicity.  
The $G$-action is said to be {\em completely dissipative} 
if the subset $D$ can be chosen so that 
$$
\mu(G D)=\mu(X). 
$$
Such subset $D$ is said to be a {\em measurable fundamental domain} for the $G$-action on $X$. 

A subset $A\subset X$  is said to be 
{\em $G$-recurrent} if there exists $g\in G\setminus \{1\}$ such that $\mu(A\cap gA)>0$. 
A subset $A\subset X$ is said to be {\em infinitely $G$-recurrent} if there are infinitely many elements $g\in G$ such that $\mu(gA\cap A)>0$. 
The $G$-action is said to be {\em (infinitely) conservative} if for every measurable subset $D\subset X$ of positive measure is (infinitely) recurrent. 

For the next theorem we will assume that $(X,\mu)$ is {\em Lebesgue}, i.e. is measure-isomorphic to the interval $[0,1]$ with its Lebesgue measure. 
 All Borel probability measures on compact metrizable spaces satisfy this property.

\begin{thm}
[Hopf's conservative--dissipative dichotomy, see e.g. Theorem 14 in \cite{Kaimanovich}] 
There is a $G$-invariant partition of $X$ in two (disjoint) subsets $Y\sqcup Z$ such that the $G$-action on $Y$ is conservative while the action on $Z$ is totally dissipative. 
\end{thm}


Note that  {\em conservative} in this theorem cannot be replaced with {\em infinitely conservative}. A trivial example is given by a non-effective action 
$G\times X\to X$ with finite kernel $K< G$, such that the induced action on $X$ of $G/K$ is dissipative. In general, the discrepancy between 
{\em conservative} and {\em infinitely conservative} is due to existence of finite order elements $g\in G$ whose 
fixed point sets in $X$ has positive measure, see  Part (iii) of \cite[Theorem 14]{Kaimanovich}. 

\medskip 
What we need for this book is a much weaker statement, which is a weak form of Proposition 10 in \cite{Kaimanovich}:

\begin{lem}\label{lem:infinitely conservative} 
If the $G$-action on a non-atomic probability space $(X,\mu)$ preserving the measure class of $\mu$ is ergodic then it is also infinitely conservative. 
\end{lem}
\begin{proof} We include a proof for the sake of completeness. Let $A\subset X$ be a subset of positive measure. 
Suppose that there are only finitely many elements $g_1,...,g_n\in G$ such that $\mu(g_iA\cap A )> 0$. 
Since $\mu$ has no atoms and the $G$-action preserves the sigma-algebra of  the measure $\mu$, there is a subset $B\subset A$ of positive measure such that  
$$
\mu(g_i B)< \frac{\mu(A)}{n},  i=1,...,n.$$
In particular,
$$
\mu(A)> \mu(g_1B \cup ....\cup g_nB)
$$  
and, thus, 
$$
\mu(X\setminus GB)\geqslant \mu(A\setminus GB)>0. 
$$
Hence, we obtain a $G$-invariant measurable subset $Y=GB\subset X$ of positive measure, whose complement in $X$ 
also has positive measure. This contradicts ergodicity of the action.
\end{proof}

\begin{example}
Let $G$ be a countably infinite group with a bijection $\phi: \bN\to G$. We equip $G$ with the probability measure  
$$
\mu=\frac{6}{\pi^2}\sum_{n=1}^\infty \frac{1}{n^2}\delta_{\phi(n)}. 
$$
The action of $G$ on $(X=G,\mu)$ via left multiplication preserves the measure class of $\mu$. The action is transitive, hence, ergodic. At the same time, 
it is totally dissipative. 
\end{example}

\section{Volumes and degrees of real and complex algebraic sets}
\label{appendix_volume_degree}

In this section, we will show  how to bound the $k$-dimensional volume of a algebraic set by its degree. 
This bound can be obtained using \cite[p. 239, Triangulation lemma]{gromov_entropy_semi_algebraic} but also using the Crofton formula. 
Before stating the main proposition, we recall various  notions of degrees associated to real or complex algebraic subsets. 

Let $AGr(k,n)$ denote the affine grassmannian of affine $k$-dimensional subspaces in $\R^n$. Thus, $AGr(k,n)$ has natural structure of a manifold diffeomorphic 
to the product $\R^n\times Gr(k,n)$, where $Gr(k,n)$ is the Grassmannian of linear $k$-dimensional subspaces in $\R^n$. We will equip $AGr(k,n)$ with a 
Riemannian metric invariant under $\R^n \rtimes O(k)$ and the corresponding Riemannian measure $\nu$. 

A {real semi-algebraic subset} in $\R^n$ is 
the intersection of subsets of $\R^n$ defined by polynomial equalities or inequalities in $n$ variables. 
Similarly, a complex algebraic subset of $\C^n$ is defined by intersecting zero loci of complex polynomials. 
Such sets have a notion of dimension, thanks to the stratification result (see e.g \cite[Proposition 4.4]{yomdin_comte}), where the dimension of a real semi-algebraic subset is the dimension of the stratum of maximal dimension.
Given a real algebraic subset $Y$ in $\R^n$ of dimension $k$, the {\em degree} $\deg(Y)$ of $Y$ is the maximal number of connected components of 
$Y \cap F$ where the maximum is taken over all $F\in AGr(n-k,n)$.   
It turns out that $\deg(Y)$ is finite by \cite[Corollary 4.9]{yomdin_comte} and can be bounded by a 
constant which depends on the degrees of the polynomials defining $Y$. 
%

We start with the following basic result. 

\begin{prop} Fix $k \leqslant n$. Then for any closed ball $K \subset \R^n$ with non-empty interior and any real semi-algebraic set $Y$ of dimension $k $, one has $\vol_i(Y\cap K)=0$ for all $i \geqslant k +1$. 
\end{prop}

\begin{proof}
By the stratification of a semi-algebraic set $Y$, we can decompose $Y$ into $\sqcup_{j} Y_j$ where each $Y_j$ is a smooth submanifold of dimension at most $k$. 
In particular, this gives $\vol_i(Y\cap K) = \sum_{j} \vol_i(Y_j\cap K)$.
 Since $i\geqslant k+1$, all these volumes vanish, hence $\vol_i(Y\cap K) = 0$. 
\end{proof}

More interestingly, when the dimension  matches with the $k$-th dimensional volumes, we use the following result. 

\begin{prop} \label{prop_volume_degree}
Fix $k \leqslant n$ and an integer $d> 0$. Then for any closed ball  $K\subset \R^n$ with non-empty interior, there exists a constant $C> 0$ (depending only on $n,k$ and the diameter of $K$) such for any real semi-algebraic subset  $Y$ of 
dimension $k$ in $K$ and degree $\deg(Y) \leqslant d$, one has $\vol_{k}(Y) \leqslant C d$. 
\end{prop}
\begin{proof}
Since $K$ is compact, it is contained in ball $B_r$ of radius $r$ in $\R^n$. 
By \cite[Corollary 5.2]{yomdin_comte}, we have: 
\begin{equation}
\vol_k(Y) \leqslant C_{n,k} \deg(Y) r^k,  
\end{equation} 
where $C_{n,k}$ is a constant which depends only on $n$ and $k$.
\end{proof}


We now turn to the complex setting. 
Given a projective subvariety $R$ of complex dimension $k$ in a complex projective variety $\cF$ of dimension $n$, fix a K\"ahler form $\Omega$, then the degree of $R$ with respect to $\Omega$ is 
\begin{equation} \label{eq_degree_complex}
\deg R = \deg_{\Omega} R = \langle [ R ] \cup  [ \Omega^{n-k} ] , [\cF] \rangle  = \langle \intcur{R} , \Omega^{n-k} \rangle.
\end{equation} 

This degree depends on the choice of the K\"ahler form $\Omega$. 
Given any other K\"ahler form $\Omega'$, then there exists a constant $C>0$ such that $C \Omega' - \Omega$ is a strongly positive form in the sense of \cite[Definition 1.1]{demailly}. By \cite[Proposition 1.11]{demailly}, the difference $C^k \Omega'^k - \Omega^k$ is also positive and this shows that: 
\begin{equation} \label{eq_comparison_degrees}
 \deg_{\Omega}(V) \leqslant C^k \deg_{\Omega'} V.
\end{equation}
By exchanging the roles of $\Omega$ and $\Omega'$, one obtains that the different notion of degrees are equivalent to one another. 

Furthermore, this degree controls the volume. 

\begin{prop} Let $\cF$ be a complex (smooth) projective variety  of dimension $n$.
Fix a K\"ahler form $\Omega$ on $\cF$, then there exists a constant $C>0$ such that for any subvariety $V$ of dimension $k$, the $k$-dimensional volume with respect to some Riemannian metric satisfies: 
\begin{equation*}
\vol_k(V) \leqslant C \deg(V). 
\end{equation*}  
\end{prop}

\begin{proof} This follows from Wirtinger's inequality. Fix a Riemannian metric induced by a K\"ahler form $\Omega'$. 
 By \cite[Equation 1.24]{demailly}, we have $\vol_k (V) = \dfrac{1}{2^k k!} \langle \intcur{V}, \Omega'^k \rangle = \dfrac{1}{2^k k!} \deg_{\Omega'}(V) $. 
 Using \eqref{eq_comparison_degrees}, we obtain 
 \begin{equation*}
 \vol_k(V) = \dfrac{1}{2^k k!} \deg_{\Omega'}(V) \leqslant \dfrac{C^k}{2^k k!} \deg_{\Omega}(V), 
 \end{equation*}
 as required.
\end{proof}

A priori, the discussion above yield two notions of degrees, one for real semi-algebraic subsets and the other for complex subvarieties. These notions are compatible in the following sense. 

\begin{prop}  \label{prop_degree_real_vs_complex}
Let $\cF$ be a smooth projective variety of dimension $n$ and let $\Omega$ be a K\"ahler form on $\cF$. Fix $O \simeq \C^n $ be an affine chart in $\cF$. 
There exists a constant $C> 0$ such that for 
 any  complex subvariety  $Y$ of (complex) dimension $k$, the  degree $\deg_{\R}(Y\cap O)$ of $Y \cap O$, viewed as a real algebraic subset is controlled by the complex degree $\deg_\Omega(Y)$: 
\begin{equation*}
\deg_{\R} (Y\cap O) \leqslant C \deg_{\Omega}(Y), 
\end{equation*}   
\end{prop}

\begin{proof}
Indeed, intersecting with a complex hyperplane $H$ in $\C^n$ is the same as intersecting with two real hyperplanes in $\R^{2n}\simeq \C^n$. Thus, intersecting with $k$ generic complex hyperplanes is the same as intersecting with $2k$ real hyperplanes. 
 This shows that $\deg_{\R}(Y \cap O) = \# \{ Y \cap H_1 \cap \ldots \cap H_k \cap O \}$ where $H_i$ are generic hyperplanes, where the intersection is counted with multiplicity.  
 We then have $\deg_{\R}(Y \cap O) \leqslant \# \{ \bar Y \cap H_1 \cap \ldots \cap H_k \} = [\bar Y] \cup [H]^k$ where $\bar Y$ is the Zariski closure of $Y$ and where $H$ is the cohomology class of the hyperplanes $H_i$. Choosing a K\"ahler form $\Omega'$ representing the class $[H]$, we have thus shown that $\deg_{\R}(Y\cap O) \leqslant \deg_{\Omega'}(\bar Y)$. The conclusion follows from \eqref{eq_comparison_degrees}. 
\end{proof}

\chapter[Appendix E]{Transformation groups}

In this chapter we will prove a certain general result about displacements of elements of smooth transformation groups of Riemannian manifolds and then apply it to thickenings in flag manifolds.  

\section{Displacements in transformation groups}

In what follows, we fix a compact connected Riemannian manifold $M$. Let $d$ denote the Riemannian distance function on $M$ and   
$\mathfrak{X}(M)$ the Lie algebra of vector fields on $M$. 
We will use the notation $ev_x$ for the {\em evaluation map}  
$$
ev_x: \mathfrak{X}(M)\to T_xM, \quad ev_x(\xi)=\xi(x)\in T_xM.
$$
We equip $\mathfrak{X}(M)$ with the 
supremum-norm $||\xi||$ using the Riemannian metric of $M$. For $\xi\in \mathfrak{X}(M)$, set 
$$
Max(\xi):= \{x\in M: ||\xi(x)||=||\xi||\}.$$ 
Similarly, for a self-diffeomorphism $f: M\to M$ and $x\in M$ define its {\em displacement at $x$,} $Dis_x(f):=d(x, f(x))$ and the {\em maximal displacement} 
$Dis(f):=  \max_{x\in M} Dis_x(f)$. We also set 
$$
Max(f):=\{x\in M: Dis_x(f)=Dis(f)\}.
$$   
 
Below, we consider a Lie group $G$ with Lie algebra $ {\mathfrak g}$ and a smooth transitive 
action  $G\times M\to M$ with discrete kernel. (The action need not be isometric although for our application it 
will suffice to consider the case of compact groups $G$ acting isometrically.) 
This action induces an isomorphism of the Lie algebra ${\mathfrak g}$ to a 
subalgebra $\cL\subset \mathfrak{X}(M)$ of the Lie algebra of vector fields on $M$. 
By abusing the notation, we will identify ${\mathfrak g}$ and $\cL$ via this isomorphism and identify elements $g\in G$ 
with the corresponding self-diffeomorphisms of $M$ given by the $G$-action. 
The main result of this chapter is the following theorem: 

\begin{thm}\label{thm:displacements}
Under the above assumptions,  there exist constants $r > 0, C >0$ such that for all $x, z\in M$ with $d(x,z)<r$, there is 
$g\in G$ such that $g(x)=z$ and $Dis(g)\leqslant C Dis_x(g)$. 
\end{thm}
\begin{proof} 1. We first observe that it suffices to prove the result for a particular point $x\in M$. Indeed, fix $p\in M$ and let $K\subset G$ 
be a relatively compact subset such that the map $K\to M$, $k\mapsto kp$, is a bijection.
Assume that $r=r_p, C=C_p$ (as required by the theorem) exist for this point $p$, i.e. for every $z\in B(p,r_p)$ there is $g\in G$ such that 
$g(p)=z$ and $ Dis(g)\leqslant C_p Dis_p(g)$. Then for any other point $x\in M$ let $k\in K$ be such that $k(p)=x$. Then for $w=k(z)\in k(B(p,r))$ we take 
$g\in G$ such that $g(p)=z$ and set
$$
h= k g k^{-1}.
$$
The element $h$ sends $x$ to $w$. Moreover, since the elements $k$ belong to a relatively compact subset $K\subset G$, the Lipschitz constants of $k$ and $k^{-1}$ 
are uniformly bounded by some $L\geqslant 1$. Thus,
$$
Dis_x(h)\geqslant L^{-2} Dis_p(g), \quad Dis(h)\leqslant L^2 Dis(g)
$$
and, therefore,
$$
Dis(h) \leqslant L^{4}C_p Dis_x(h). 
$$
Of course, we do not claim that $k(B(p,r))= B(x,r)$ (for this, we would have to assume that the action of $K$ is isometric). However, in view of the same 
Lipschitz bound as above,
$$
B(x, L^{-1}r)\subset k(B(p,r)),
$$
where $r=r_p$. Thus, for the theorem, we can take $C:= L^{4}C_p$ and $r:=L^{-1}r_p$. This reduces the proof to the case of the base-point $x=p\in M$.

\medskip 
2. We next prove an infinitesimal version of the theorem. 
We equip the subalgebra $\cL\cong {\mathfrak g}$ in ${\mathfrak X}(M)$ with the restriction of the supremum-norm $||\xi||$. 
Since the $G$-action on $M$ is transitive, for every point $x\in M$ there is a linear subspace  $\cL_x\subset \cL$ 
which maps isomorphically  to $T_xM$ via the evaluation map $ev_x: \cL\to T_x(M)$. In line with Step 1 of the proof, we choose the subspaces $\cL_x$ a bit 
more carefully, using a particular base-point $p\in M$ and a relatively compact subset $K\subset G$ as in Step 1. Then, given 
a subspace $\cL_p$ for the base-point $p$, we define $\cL_{g(p)}$ by $g_* \cL_p, g\in K$. We thus obtain a family of linear subspaces 
$\cL_x\subset \cL$, $x\in M$.

\medskip 
Each subspace $\cL_x$ has the second norm, the one given by the pull-back of the norm $||\cdot||_x$ on $T_xM$ via the map $ev_x$. 
This norm is equivalent to the supremum-norm $||\cdot||$, as these are norms on a finite-dimensional vector space. 
Hence, for every vector $v\in T_xM$ there exists (unique) $\xi\in \cL_x$ such that $ev_x(\xi)=v$ and 
$||\xi||\leqslant C_x ||v||$ for a constant $C_x$ depending only on $x$.  In view of the construction of the subspaces $\cL_x$ via the $K$-action, 
the constants $C_x$ are uniformly bounded away from $0$ as $x$ varies through $M$. We, thus, obtain:

 \begin{lem}\label{lem:E1}
There exists $C_1>0$ such that for every $x\in M$ and $\xi\in \cL_x$, 
$$
||\xi|| \leqslant C_1 ||\xi(x)||=C_1 ||ev_x(\xi)||. 
$$
\end{lem}  

This is an infinitesimal version of Theorem \ref{thm:displacements}. 


3. Our next goal is to relate displacements $Dis_x(g)$ of elements of $g\in G$ and norms at $x$ of vector fields in $\cL$. Here and in what follows, $\exp:  
{\mathfrak g}\to G$ is the exponential map of the Lie group $G$.

\begin{lem}\label{lem:E2}
For every $x\in M$ and every $\xi\in {\mathfrak g}$, the curve $g_t:=  \exp(t\xi)\in G, t\in \mathbb R$, satisfies 
$$
\lim_{t\to 0+} \frac{1}{t}Dis_{x}(g_t)= ||ev_x(\xi)||=||\xi(x)||. 
$$
Moreover, if $||\xi||\leqslant 1$, then convergence is uniform with respect to $x\in M$ and $\xi\in \cL$.  
\end{lem}  
\begin{proof} This result is a standard differential geometry fact, we present a proof for the sake of completeness. 
We let $n$ denote the dimension of $M$. 
Consider the curve $c(t)=g_t(x)$ in $M$.  In the normal coordinates at $x$ 
(where the point $x$ corresponds to the origin in $\R^n$), this curve can be written 
as $x(t)=(x^1(t),...,x^n(t))$. Then, by the Gauss Lemma,  
$$
d(x, c(t))= \sqrt{\sum_{i=1}^n (x^i(t))^2}, 
$$  
since we are using normal coordinates. We obtain: 
$$
\lim_{t\to 0+} \frac{1}{t}Dis_{x}(g_t)=\lim_{t\to 0+} \frac{d(x, c(t))}{t} = \lim_{t\to 0}  \sqrt{\sum_{i=1}^n \left(\frac{x^i(t)}{t}\right)^2}= 
 \sqrt{\sum_{i=1}^n \left(\frac{d}{dt}x^i(0)\right)^2}=||c'(0)||=||\xi(x)||. 
$$ 
To get a uniform bound on convergence in a normal neighborhood of $x$, 
we just need to get a uniform bound on the 2nd derivative of the curve $x(t)$ near $t=0$. But the 2nd derivatives are uniformly bounded since 
the curve $c(t)$ comes from a 1-parameter subgroup in $G$. 
\end{proof}

\begin{lem}\label{lem:E3}
Fix $\xi\in \cL$ which does not vanish at a point $x\in M$ and take, as before, $g_t:= \exp(t\xi)\in G$. 
Then for every $y\in M$ 
$$
\lim_{t\to 0} \frac{Dis_y(g_t)}{Dis_x(g_t)}= \frac{||\xi(y)||}{||\xi(x)||}. 
$$
In particular, this equality holds for every nonzero $\xi\in \cL_x$. Moreover, convergence is uniform (with respect to $x$ and $\xi$) 
for $\xi\in \cL_x$, provided that  $||\xi||= 1$. 
\end{lem}
\begin{proof} The only part which is not immediate is uniformity of convergence. By Lemma \ref{lem:E2}, we have
$$
Dis_y(g_t)= ||\xi(y)|| + R_y(t), \quad Dis_x(g_t)= ||\xi(y)|| + R_x(t),
$$
with $R_x(t)\to 0, R_y(t)\to 0$ as $t\to 0$, and convergence is uniform in $x$, $y$ and $\xi$. Thus,
 \begin{equation}\label{eq:E1}
 \left| \frac{Dis_y(g_t)}{Dis_x(g_t)}- \frac{||\xi(y)||}{||\xi(x)||}\right|= \left|\frac{R_y(t)||\xi(x)|| - R_x(t) ||\xi(y)||}{||\xi(x)|| (||\xi(x)|| + R_x(t))}\right|. 
 \end{equation}
We take $T$ such that $|t|\leqslant T$ implies that $R_x(t)<||\xi(x)||/2$ (this $T$ can be chosen to be independent of $x$ and $\xi\in \cL_x$, as $||\xi||=1$). 
Then, since $||\xi(x)||\leqslant ||\xi||\leqslant 1$ and $||\xi(y)||\leqslant 1$)  the right-hand side of \eqref{eq:E1} is at most
$$
2\frac{|R_x(t)| + |R_y(t)|}{||\xi(x)||^2}\leqslant 2 C_1^2 \left(|R_x(t)| + |R_y(t)| \right), 
$$
provided that $|t|\leqslant T$. (Here we use the inequality $||\xi||\leqslant C_1||\xi(x)||$ estblished in Lemma \ref{lem:E1}.) 
This ensures uniform convergence claim. \end{proof}

\medskip
4. We next relate these results to maximal displacements: 


\begin{lem}\label{lem:E4}
Assuming that $\xi\in \cL_x$ satisfies $||\xi||=1$, we have
$$
\lim_{t\to 0} \frac{Dis(g_t)}{Dis_x(g_t)}= \frac{||\xi||}{||\xi(x)||}.
$$
Moreover, convergence is uniform in $x$ and $\xi$. 
\end{lem}
\begin{proof} Take a sequence $t_i\to 0$ and $y_i\in Max(g_{t_i})$ converging to some $y\in M$. 
Then, by Lemma \ref{lem:E3},
$$
\lim_{t\to 0} \frac{Dis(g_t)}{Dis_x(g_t)}=\lim_{i\to\infty} \frac{Dis(g_{t_i})}{Dis_x(g_{t_i})}= \lim_{i\to\infty} \frac{Dis_{y_i}(g_{t_i})}{Dis_x(g_{t_i})}= 
 \frac{||\xi(y)||}{||\xi(x)||}\leqslant \frac{||\xi||}{||\xi(x)||}
$$
and, moreover, convergence is uniform in $x$ and $\xi$.  To get the opposite (less useful) inequality, take $y\in Max(\xi)$ and observe that 
$$
\lim_{t\to 0} \frac{Dis(g_t)}{Dis_x(g_t)} \geqslant \lim_{t\to 0} \frac{Dis_y(g_t)}{Dis_x(g_t)}=  \frac{||\xi(y)||}{||\xi(x)||}= \frac{||\xi||}{||\xi(x)||}. 
$$
Again, convergence is uniform in $x$ and $\xi$. 
\end{proof}

\medskip
5. In this step of the proof we analyze the image of the {\em action map} $a_p$ at the base-point $p\in M$. Define the {\em action map}  
$$
a_p: 
\cL_p \to M, a_p(\xi)= \exp(\xi)(p)\in M. 
$$
The derivative of this map at $0$ is the evaluation map $ev_p: \xi\mapsto \xi(p)$ (where we identify $T_0\cL_p$ and $\cL_p$). Since, by the definition of $\cL_p$, the 
evaluation map is an isomorphism, the inverse mapping theorem gives us:

\begin{lem}\label{lem:E5}
There is a number $\delta_0>0$ such that the restriction of $a_p$ to the ball 
$$
B(\delta_0):= \{\xi\in \cL_p: ||\xi||<\delta_0\}\subset \cL_p
$$
is a diffeomorphism onto a neighborhood of $p$ in $M$. In particular, there is a function $r(\delta)>0$ defined for $\delta\in (0,\delta_0)$ such that 
the image $a_p(B(\delta))$ contains the geodesic ball  
$B(p, r(\delta))$. 
\end{lem} 

6. We can now finish the proof of the theorem. As we established in the first step of the proof, it suffices to prove the result for $x=p$. We will consider unit vectors 
$\xi\in \cL_p, ||\xi||=1$. According to Lemma \ref{lem:E4}, there exists $\eps>0$ such that for all $t\in (0,\eps]$ we have 
$$
\left| \frac{||\xi||}{||\xi(p)||} - \frac{Dis(g_t)}{Dis_p(g_t)}\right|< C_1, 
$$
where $C_1$ is the constant from Lemma \ref{lem:E1}. Therefore, for such $t$, since $\frac{||\xi||}{||\xi(p)||} \leqslant C_1$, 
$$
\frac{Dis(g_t)}{Dis_p(g_t)} \leqslant 2C_1=C, \quad Dis(g_t) \leqslant C Dis_p(g_t). 
$$
(For $t=0$ the last inequality is trivially true as well.) This is the desired bound on the displacement of $g_t$ at $p$. The set of images $g_t(p)$ is the image of the 
action map, $a_p(B(\eps))$. Without loss of generality, we may assume that $\eps<\delta_0$, where $\delta_0$ is the constant from Lemma \ref{lem:E5}. Therefore, 
by Lemma \ref{lem:E5}, $a_p(B(\eps))$ will contain a ball $B(p,r)$, $r=r_p=r(\eps)$. This concludes the proof of the theorem. 
\end{proof}

\begin{question}
Suppose that $G$ is compact and its action on $M$ is isometric. Can one take $C=1$ in this case? 
\end{question}

\section{Application to thickenings}

 Recall that for every thickening $\Th$, we use the notation $\widehat\Th$ for the inverse thickening in a Weyl group $W$ of a complex semisimple Lie group $G$. 
 Thus, for points $x, y$ in the corresponding flag manifold $\cF=G/B$ we have  $x\in \Th(y)\iff y\in \widehat\Th(x)$. 

The following lemma is elementary and we will apply to specific thickenings.
\begin{lem} \label{lem_continuity_hausdorff}  Let $\cF$ be a Riemannian manifold  on which there exists a Lie group $G$ acting by automorphism. 
Fix a compact subset $A \subset \cF$.  Then the function $g \in G \mapsto gA$ is continuous, where the set  $\{ gA \subset \cF | g \in G\}$ is endowed with the Hausdorff distance.
\end{lem}

\begin{proof}
Denote by $d_{\cF}$ the Riemannian distance on $\cF$. Since the group $G$ acts by automorphism and is a Lie group, we are reduced to proving the continuity near the identity.
For any $g \in G$, recall from \cite[Definition 1.5.1]{krantz_parks} the definition of the Hausdorff distance: 
\begin{equation*}
d_H(A,gA) = \max (\sup_{a \in A} d_{\cF}(a, gA), \sup_{a' \in gA} d_{\cF}(a', A)).
\end{equation*} 
Take $a \in A$ and  $g = \exp(t \xi)$ with $\xi \in \mathfrak{g}$ of norm bounded by $1$. Then by definition and using Lemma \ref{lem:E3}, we have: 
\begin{equation*}
d_{\cF}(a, g A) \leqslant d_{\cF}( a , g a ) \leqslant Dis_a(\exp(t \xi)) \leqslant 2 t  ||\xi(a) ||, 
\end{equation*}
for all $t \leqslant 1$. 
Since $a \in A \mapsto || \xi(a)||$ is a continuous function and $A$ is compact, this is uniformly bounded above by $M$ and we obtain that $d_{\cF}(a, \exp(t \xi )A) \leqslant 2 t M$ for $t \leqslant 1$. 
This shows that $\sup_{a\in A} d_{\cF}(a , \exp(t\xi)A) \leqslant 2 t M$ and by symmetry we also have $\sup_{a'\in \exp(t\xi) A} d_{\cF}(a', A) \leqslant 2 t M$. 
This gives that $d_H(A, \exp(t \xi) A) \leqslant 2 t M$ for all $t \leqslant 1$. Hence $g \in G \mapsto gA$ is continuous.  
\end{proof}

\begin{thm}\label{thm:thickening_distance_inequality}
There exists a constant $C>0$ such that for every thickening $\Th\subset W$ and all $x, y\in \cF$, 
$$
C^{-1} d_{\cF}(x, \Th(y)) \leqslant d_{\cF}(y, \widehat\Th(x))\leqslant C d_{\cF}(x, \Th(y)). 
$$
\end{thm}
\begin{proof} We will use Theorem \ref{thm:displacements} with the flag manifold $\cF=G/B$ playing the role of the manifold $M$ and the $G$ 
(or its maximal compact subgroup) playing the role of the transitive transformation group. 
Let $r, C$ be the constants as in Theorem \ref{thm:displacements}. Let $z\in \Th(y)$ be a point closest to $x$. It suffices to consider the case when $d_{\cF}(x,z)<r$, 
for otherwise, we simply take $C':= r^{-1}\diam(\cF, d_{\cF})$ and obtain 
$$
d_{\cF}(y, \widehat\Th(x))\leqslant C' d_{\cF}(x, \Th(y)). 
$$
Take $g\in G$ as in Theorem \ref{thm:displacements}. Then 
$g(\widehat\Th(x))=\widehat\Th(g(x))$. Since $g(x)\in \Th(y)$,  $y\in \widehat\Th(z)= g(\widehat\Th(x))$. The latter implies that the point $y':=g^{-1}(y)$ is in 
$\widehat\Th(x)$. Lastly,  
$$
d_{\cF}(y, \widehat\Th(x)) \leqslant d_{\cF}(y, y')\leqslant Dis(g)\leqslant C Dis_x(g)=d_{\cF}(x, z)= C d_{\cF} (x,  \Th(y)).$$
The opposite inequality follows from the fact that inversion of thickenings is an involution. 
\end{proof}





\backmatter

\bibliographystyle{amsalpha}

\bibliography{references}

\newpage

\printindex

\end{document}